\documentclass[reqno, 10pt]{amsart}
\usepackage[utf8]{inputenc}
\usepackage{amssymb, amsthm}
\usepackage{mathrsfs}
\usepackage{bm,bbm}
\usepackage{upgreek}
\usepackage{microtype} %
\usepackage{comment} %

\usepackage{tikz}
\usetikzlibrary{shapes.geometric, arrows}
\tikzstyle{choice} = [rectangle, rounded corners, 
minimum width= 100, 
minimum height=30,
text centered, 
draw=black, 
text width = 100 ,
fill=red!30]

\tikzstyle{result} = [rectangle, 
minimum width= 100 , 
minimum height=1 , 
text centered, 
text width=100, 
draw=black, 
fill=orange!30]

\usepackage{accents}
\newlength{\dhatheight}
\newcommand{\hhat}[1]{%
    \settoheight{\dhatheight}{\ensuremath{\hat{#1}}}%
    \addtolength{\dhatheight}{-0.35ex}%
    \hat{\vphantom{\rule{1pt}{\dhatheight}}%
    \smash{\hat{#1}}}}

\usepackage{amsfonts}
\DeclareMathAlphabet{\mathpgoth}{OT1}{pgoth}{m}{n}
\DeclareMathAlphabet{\mathesstixfrak}{U}{esstixfrak}{m}{n}
\DeclareMathAlphabet{\mathboondoxfrak}{U}{BOONDOX-frak}{m}{n}
\usepackage{amsmath,amssymb}
\usepackage[bbgreekl]{mathbbol}
\usepackage{yfonts,mathtools}

\usepackage{tikz-cd}
\usepackage[all,cmtip]{xy}

\numberwithin{equation}{section}
\usepackage{soul}

\usepackage{hyperref}
\hypersetup{colorlinks}
\definecolor{darkred}{rgb}{0.5,0,0}
\definecolor{darkgreen}{rgb}{0,0.5,0}
\definecolor{darkblue}{rgb}{0,0,0.5}

\hypersetup{colorlinks, linkcolor=darkblue, filecolor=darkgreen, urlcolor=darkred, citecolor=darkblue}
\makeatletter %
\@addtoreset{equation}{section}
\makeatother  %

\numberwithin{equation}{section}

\newtheorem{thma}{Theorem}

\newtheorem{situationf}{Situation}

\newtheorem{situationm}{Situation}

\newtheorem{situationc}{Situation}

\newtheorem{situationh}{Situation}

\newtheorem{situationg}{Situation}

\renewcommand{\subset}{\subseteq}

\usepackage{ulem}

\newtheorem{thm}{Theorem}[section]
\newtheorem{cor}[thm]{Corollary}

\newtheorem{prop}[thm]{Proposition}

\newtheorem{lemma}[thm]{Lemma}
\theoremstyle{definition}
\newtheorem{defn}[thm]{Definition}

\theoremstyle{remark}
\newtheorem{rem}[thm]{Remark}

\usepackage{xcolor}
\newtheorem{example}[thm]{Example}
\newtheorem{notation}[thm]{Notation}

\newcommand{\beq}{\begin{equation}}
\newcommand{\eeq}{\end{equation}}
\newcommand{\beqn}{\begin{equation*}}
\newcommand{\eeqn}{\end{equation*}}
\newcommand{\ov}{\overline}
\newcommand{\mb}{\mathbb}
\newcommand{\mc}{\mathcal}
\newcommand{\mf}{\mathfrak}
\newcommand{\ms}{\mathscr}

\newcommand{\bF}{{\mathbb F}}

\newcommand{\dep}{{\rm depth}}

\newcommand{\regpos}{\uds{\bf RegPos}}

\newcommand{\p}{{\sf p}}
\newcommand{\fp}{{\bf F}_{\sf p}}
\newcommand{\kp}{{\bf K}_{\sf p}}
\newcommand{\tkp}{\wt{\bf K}_{\sf p}}
\newcommand{\zp}{{{\mb Z}_{\sf p}}}

\newcommand{\ep}{{{\mathbb E}}}

\renewcommand{\outer}{{}_{\lfloor}}

\newcommand{\tu}{\mathtt{U}}

\newcommand{\pman}{\uds{{\bf \Psi}{\bf Man}}}

\newcommand{\qst}{{\it QSt}_{\p}}
\newcommand{\fst}{{\it FSt}_{\p}}

\newcommand{\psu}{{[\![\tu]\!]}}
\newcommand{\pst}{{[\![\uptheta]\!]}}
\newcommand{\lsu}{{(\!(\tu)\!)}}
\newcommand{\lst}{{(\!(\uptheta)\!)}}

\newcommand{\wt}{\widetilde}

\newcommand{\wh}{\widehat}

\newcommand{\wham}{{\rm H}\wt{\rm am}}

\newcommand{\bsy}[1]{\textcolor{purple}{#1}}

\newcommand{\uds}[1]{\underline{\smash{#1}}}

\newcommand{\ev}{{\rm ev}}

\title[Integral Hamiltonian Floer theory]{Integral Hamiltonian Floer Theory: Foundations}

\author{Shaoyun Bai}
\address{MIT, 77 Massachusetts Avenue Cambridge, MA 02139, USA}
\email{shaoyunb@mit.edu}

\author{Guangbo Xu}
\address{Department of Mathematics, Rutgers University, Hill Center--Busch Campus, 110 Frelinghuysen Road, Piscataway, NJ 08854-8019, USA}
\email{guangbo.xu@rutgers.edu}

\thanks{The first author is supported by NSF DMS-2404843}

\thanks{The second author is supported by NSF DMS-2345030, NSF DMS-2506403, and Simons Foundation's Travel Support for Mathematicians.}

\begin{document}

\begin{abstract}
Continuing our previous work on the integral Arnold conjecture \cite{Bai_Xu_Arnold}, we establish the foundation of Hamiltonian Floer theory over integer coefficients on any compact symplectic manifold. This package produces a well-defined chain homotopy class of Floer complexes, chain-level continuation maps, the Piunikhin--Salamon--Schwarz (PSS) isomorphism, and an associative pair-of-pants product. When the Hamiltonian is the $\p$-th iteration for a prime number $\p$, we also establish the $\zp$-equivariant Floer theory in characteristic $\p$, including the equivariant Floer complex and the $\zp$-equivariant pair-of-pants product. Moreover, we define quantum Steenrod operations for any compact symplectic manifold and compare them with the equivariant pair-of-pants products via the $\zp$-equivariant PSS map. We also obtain chain-level invariants such as spectral numbers and barcodes and prove quantitative properties of the equivariant pair-of-pants product.

The counting theory is based on Fukaya--Ono's normally polynomial perturbation scheme \cite{Fukaya_Ono_integer} and its realization by the authors in \cite{Bai_Xu_2022}, termed FOP perturbations, adapted in the abstract setting of flow categories, flow multimodules, and homotopies of flow bimodules. Along the way, we construct global Kuranishi charts for the relevant moduli spaces following Abouzaid--McLean--Smith's framework and its adaptation to the Hamiltonian Floer setting due to the authors, which may be of independent interest.

\end{abstract}

\maketitle

\setcounter{tocdepth}{1}
\tableofcontents

\section{Introduction}

Floer theory is one of the most influential discoveries in geometry and topology. In symplectic topology, the {\it Hamiltonian Floer homology} \cite{Floer_CMP} provides a powerful framework for studying Hamiltonian dynamics. Most notably, it leads to proofs of various versions of the famous Arnold conjecture \cite{Arnold_conjecture}. Several variations, such as quantitative Floer theory and local Floer theory, played important roles in many other applications in Hamiltonian dynamics. 

The development of Floer homology was accompanied by the advancement of the understanding of pseudoholomorphic curves \cite{Gromov_1985} and regularizations of their moduli spaces. In earlier works \cite{Floer_CMP}\cite{Hofer_Salamon}\cite{Ono_1995} where certain vanishing or positivity conditions were assumed, when defining the Floer complexes, one only needs to consider moduli spaces of virtual dimensions at most one; one can achieve transversality of such moduli spaces by perturbing geometrically and can obtain Floer complexes with coefficients in ${\mb Z}$. When the semipositivity condition was removed, there is no known way to achieve transversality via geometric perturbations without imposing further conditions (cf. \cite{Cieliebak_Mohnke}). The abstract perturbation methods used in \cite{Liu_Tian_Floer}\cite{Fukaya_Ono}\cite{HWZ_book}\cite{Pardon_virtual} require one to regularize the moduli spaces via {\it Kuranishi models} or {\it polyfolds} and to apply virtual perturbations or integrations. One then has to deal with moduli spaces of negative virtual dimensions and to organize all the constructions in a compatible way which can be formalized using Cohen--Jones--Segal's notion of  {\it flow categories} \cite{Cohen_Jones_Segal}. The lack of geometric transversality is closely intertwined with the orbispace structure of the moduli spaces, which naturally leads to fractional counts of   curves and chain-level structures only defined over ${\mb Q}$.

The main purpose of the current work is to extend Floer theory from rational to integer coefficients for general closed symplectic manifolds. In \cite{Bai_Xu_Arnold} we have already established a portion of the integral Floer theory which provides a proof of the Arnold conjecture over ${\mb Z}$. The motivation to go beyond the scope of \cite{Bai_Xu_Arnold}, besides a complete package of the integral Hamiltonian Floer theory, resides in the potential applications in Hamiltonian dynamics. For this purpose, we will construct quantum Steenrod operations for general symplectic manifolds (see \cite{Wilkins_2020}\cite{Seidel_Wilkins_2022}\cite{Seidel_pants}\cite{shelukhin-zhao} under semipositivity or stronger assumptions). In an upcoming joint work with Shelukhin and Wilkins \cite{BSWX}, we will extend many dynamical applications (see \cite{shelukhin20}\cite{shelukhin-21}\cite{shelukhin-22}\cite{CGG}\cite{Atallah_Lou}) of quantum Steenrod operations to greater generality. %

The integral extension is possible due to the recent developments in the virtual technique. The most crucial idea was given by Fukaya--Ono \cite{Fukaya_Ono_integer}. Roughly speaking, the moduli spaces carry natural stable complex structures, which allow a specific type of single-valued perturbations (normally polynomial perturbations) resulting in well-defined integer counts. Fukaya--Ono discussed the original idea for a local Kuranishi model, while B. Parker \cite{BParker_integer} further explored strategies to globalize. In \cite{Bai_Xu_2022} we provide a rigorous treatment of Fukaya--Ono's proposal and define the so-called ``FOP transversality condition'' for normally polynomial perturbations. 

Another crucial concurrent development is the Abouzaid--McLean--Smith construction of global Kuranishi charts (see \cite{AMS, Hirschi_Swaminathan, AMS2}). This new technique helps avoid the complicated inductive arguments involving multiple charts (such as having good coordinate systems \cite{FOOO_2016}). Meanwhile, the smoothing argument of \cite{AMS} replaces the hard analytic approach \cite{FOOO_smooth} by a soft topological approach. These novel ideas helped significantly in the case of Hamiltonian Floer theory. One of the main challenges in Floer theory is to carry out AMS construction for infinitely many moduli spaces compatibly.

We hope that the foundational work laid in this paper can serve as a useful package for applications of Floer theory to Hamiltonian dynamics on general compact symplectic manifolds. Our joint paper with Shelukhin and Wilkins \cite{BSWX} showcases how the technical work done here can efficiently lead to interesting statements for Hamiltonian maps on general compact symplectic manifolds. Meanwhile, we hope the technical procedures can serve as a model to be adopted in other situations such as in Lagrangian Floer theory or contact topology.

Our main results fall into three parts:

\begin{enumerate}

\item Integral Hamiltonian Floer theory at the chain level (Theorems \ref{thma_Floer_complex}---\ref{thma_semipositive}).

\item Equivariant Floer theory and quantum Steenrod operations (Theorems \ref{thma_QST}---\ref{thma_Floer_Steenrod}).

\item Quantitative refinements, including spectral invariants and barcodes (Theorems \ref{thma_Floer_filtered}---\ref{thm_pants_filtered}) as well as results about Floer complexes of homologically nontrivial orbits (Theorem \ref{thm_Floer_noncontractible} and Theorem \ref{thm_noncontractible_pants}).
\end{enumerate}
Below we state these main results in a way that readers can directly cite.

\subsection{Main results on Floer complexes}

\subsubsection{Floer complexes and Floer homology}

First and foremost, by applying the FOP perturbation scheme to the moduli spaces of Floer trajectories, we construct a Floer complex with integer coefficients, which is well-defined up to chain homotopy equivalence. Under certain restrictions of the Hamiltonian, such a  complex was constructed by the authors \cite{Bai_Xu_Arnold} without proving the  homotopy equivalence between different choices. 

We start by recalling the basic setup. Let $(X, \omega)$ be a compact symplectic manifold. Let 
\beq\label{Novikov_group}
\Pi:= \frac{ \pi_2(X) }{{\rm ker} (\omega) \cap {\rm ker} (c_1)},
\eeq
which is a finitely generated free abelian group. Let $\Lambda^\Pi$ be the Novikov ring with integer coefficients:
\beq\label{Novikov_ring}
\Lambda^\Pi = \Big\{ \sum_{i=1}^\infty a_i T^{A_i}\ |\ a_i \in {\mb Z},\ A_i \in \Pi,\ \lim_{i \to \infty} \omega(A_i) = +\infty \Big\}.
\eeq
If $H$ is a nondegenerate 1-periodic Hamiltonian generating a 1-periodic family of Hamiltonian vector fields $X_{H_t}$, a 1-periodic orbit is a map $x: S^1 \to X$ solving
\beqn
x'(t) = X_{H_t}(x(t)).
\eeqn
A {\bf capped 1-periodic orbit} is a pair $p = (x, \tilde x)$ where $x$ is a 1-periodic orbit and $\tilde x$ is an equivalence class of extensions of $x$ to the disk ${\mb D}$, where two extensions are equivalent if their difference, regarded as a map from $S^2$ to $X$, has zero symplectic area and zero Chern number. Each capped 1-periodic orbit $p$ has a well-defined symplectic action
\beq\label{symplectic_action}
{\mc A}_H(p) = - \int_{\mb D} \tilde x^* \omega - \int_{S^1} H_t(x(t)) dt
\eeq
and a well-defined Conley--Zehnder index
\beqn
{\rm deg} p:= {\rm CZ}(p) \in {\mb Z}.
\eeqn
Denote the corresponding cohomological grading of capped periodic orbits by 
\beqn
|p|:= \frac{1}{2} {\rm dim}_{\mb R} X - {\rm CZ}(p).
\eeqn
Let ${\mc O}(H)$ be the set of all contractible 1-periodic orbits and $\tilde {\mc O}(H)$ be the set of all capped 1-periodic orbits. 

\begin{thma}[Floer chain complex]\label{thma_Floer_complex}
Let $(X, \omega)$ be a compact symplectic manifold. Let $H$ be a nondegenerate 1-periodic Hamiltonian. Then for each $1$-periodic family of $\omega$-compatible almost complex structures $J$ and a certain choice $\Xi$, there is a $\Lambda_{\mb Z}^\Pi$-linear differential on the chain group 
\beqn
CF_*(H):= \Big\{ \sum_{i=1}^\infty a_i p_i\ |\ a_i \in {\mb Z},\ p_i \in \tilde {\mc O}(H),\ \lim_{i \to \infty} {\mc A}_H(p_i) = -\infty \Big\}
\eeqn
(which has the structure of a module over $\Lambda^\Pi$), making it a chain complex of finitely generated free $\Lambda^\Pi$-modules
\beqn
CF_*(H, J, \Xi; \Lambda^\Pi)
\eeqn
with a ${\mb Z}$-grading defined by the Conley--Zehnder index. Moreover, the chain homotopy type of $CF_*(H, J, \Xi; \Lambda^\Pi)$ does not depend on the choice $\Xi$.
\end{thma}

Next we establish the continuation map with expected properties.

\begin{thma}[Continuation Map]\label{thma_continuation_map}
Let $(H_1, J_1, \Xi_1)$ and $(H_2, J_2, \Xi_2)$ be two triples such that the Floer chain complexes
\begin{align*}
&\ CF_*(H_1, J_1, \Xi_1; \Lambda^\Pi),\ &\ CF_*(H_2, J_2, \Xi_2; \Lambda^\Pi)
\end{align*}
are defined.
\begin{enumerate}

\item There exists a $\Lambda^\Pi$-linear chain homotopy equivalence, called the continuation map
\beqn
\Phi^{12}(\Theta): CF_*(H_1, J_1, \Xi_1; \Lambda^\Pi) \to CF_*(H_2, J_2, \Xi_2; \Lambda^\Pi)
\eeqn
whose definition further depends on a system of choices $\Theta$.

\item The chain homotopy class of the continuation map does not depend on the system of choices $\Theta$. 

\item For $(H_1, J_1, \Xi_1) = (H_2, J_2, \Xi_2)$, any continuation map is homotopic to the identity map of $CF_*(H_1, J_1, \Xi_1; \Lambda^\Pi)$.\footnote{However we do not know if the identity map on the Floer chain complex can be realized as a continuation map.}

\item If $(H_3, J_3, \Xi_3)$ is another triple for which the complex $CF_*(H_3, J_3, \Xi_3; \Lambda^\Pi)$ is defined, and 
\beqn
\Phi^{23}(\Theta'): CF_*(H_2, J_2,\Xi_2; \Lambda^\Pi) \to CF_*(H_3, J_3, \Xi_3;\Lambda^\Pi)
\eeqn
is a continuation map. Then the composition $\Phi^{23}(\Theta') \circ \Phi^{12}(\Theta)$ is homotopic to a continuation map from $CF_*(H_1, J_1, \Xi_1; \Lambda^\Pi)$ to $CF_*(H_3, J_3, \Xi_3; \Lambda^\Pi)$.

\end{enumerate}
\end{thma}

Theorem \ref{thma_Floer_complex} and Theorem \ref{thma_continuation_map} lead to the definition of a canonical object
\beqn
HF_*(X, \omega; \Lambda^\Pi)
\eeqn
which we call the {\bf integral Floer homology} of $(X, \omega)$.

\subsubsection{Pair-of-pants product}
Next we consider the pair-of-pants product. 

\begin{thma}[Pair-of-pants product] \label{thma_product}
Let $(H_i, J_i, \Xi_i)$, $i = 1, 2, 3$ be triples for which the Floer chain complexes are defined as stated in Theorem \ref{thma_Floer_complex}. Then there is a $\Lambda^\Pi$-linear chain map 
\beq\label{product_chain}
CF_*(H_1, J_1, \Xi_1; \Lambda^\Pi)  \underset{\Lambda^\Pi}{\otimes}  CF_*(H_2, J_2, \Xi_2; \Lambda^\Pi) \to CF_*(H_3, J_3, \Xi_3; \Lambda^\Pi) \footnote{In general the tensor product over the Novikov ring needs to be completed; however, when the modules are finitely generated, the completion is unnecessary.}
\eeq
of degree $- \frac{1}{2} {\rm dim} X$ whose definition depends on a system of choices. Moreover, the chain map satisfies the following properties.
\begin{enumerate}

    \item The homotopy class of \eqref{product_chain} is independent of the choices. 
    
    \item This chain map is compatible with continuation maps. Hence it induces a canonical linear map 
\beqn
HF_*(X, \omega; \Lambda^\Pi) \underset{\Lambda^\Pi}{\otimes}   HF_*(X, \omega; \Lambda^\Pi) \to HF_*(X, \omega; \Lambda^\Pi)
\eeqn
called the pair-of-pants product. 

\item The map \eqref{product_chain} is associative up to homotopy. Hence $HF_*(X, \omega; \Lambda^\Pi)$ is a graded  ring.

\item There exists a multiplicative identity 
\beqn
{\bf 1} \in HF_* (X, \omega; \Lambda^\Pi),
\eeqn
of cohomological degree zero called the Floer fundamental class.\footnote{We do not know if it coincides with the image of the classical fundamental class under the PSS map.}
\end{enumerate}
\end{thma}

\subsubsection{Poincar\'e duality}
We then discuss Poincar\'e-type pairings on Hamiltonian Floer homology.

\begin{thma}[Poincar\'e duality] \label{thma_Poincare} Let $CF_*(H_1, J_1, \Xi_1)$ and $CF_*(H_2, J_2, \Xi_2)$ be Floer chain complexes provided by Theorem \ref{thma_Floer_complex}. Then up to a choice $\Theta$, there is a $\Lambda^\Pi$-bilinear pairing
\beq
\Phi^{\rm PD}(\Theta): CF_*(H_1, J_1, \Xi_1) { \underset{\Lambda^\Pi}{\otimes} } CF_*(H_2, J_2, \Xi_2) \to \Lambda^\Pi
\eeq
called the {\bf chain-level Poincar\'e pairing} which satisfies the following conditions.
\begin{enumerate}

\item The homotopy type of the pairing is well-defined.

\item The pairing is compatible with continuation maps. Namely, suppose $CF_*(H_3, J_3, \Xi_3)$ is another Floer chain complex and $\Phi^{31}: CF_*(H_3, J_3, \Xi_3) \to CF_*(H_1, J_1, \Xi_1)$ is a chain-level continuation map. Then the pairing
\beqn
\Phi^{\rm PD}(  \Phi^{31}(\cdot), \cdot ): CF_*(H_3, J_3, \Xi_3) { \underset{\Lambda^\Pi}{\otimes}} CF_*(H_2, J_2, \Xi_2) \to \Lambda^\Pi
\eeqn
is homotopic to a chain-level Poincar\'e pairing. 

\item In field coefficients, the induced Poincar\'e pairing on homology is nondegenerate. More precisely, let ${\bf K}$ be a field and $\Lambda_{\bf K}^\Pi$ be the Novikov field. Then the induced bilinear pairing
\beqn
\Phi^{\rm PD}: \left( HF_*(X) \widehat{\underset{\Lambda_{\mb Z}^\Pi}{\otimes}} HF_*(X) \right) \underset{\Lambda_{\mb Z}^\Pi}{\otimes} \Lambda_{\bf K}^\Pi \to \Lambda_{\bf K}^\Pi
\eeqn
is nondegenerate.

\item The Poincar\'e duality is compatible with the pair-of-pants product. More precisely,  the following diagram commutes.
\beqn
\xymatrix{ HF_*(X, \omega) \widehat{ \underset{\Lambda^\Pi}{\otimes}} HF_*(X, \omega) \widehat{ \underset{\Lambda^\Pi}{\otimes}} HF_*(X, \omega)  \ar[rr]^-{\Phi^{\rm pants} \otimes {\rm Id}} \ar[d]_{{\rm Id}\otimes \Phi^{\rm pants}}   & &  HF_*(X, \omega) \widehat{\underset{\Lambda^\Pi}{\otimes}} HF_*(X, \omega) \ar[d]^{\Phi^{\rm PD} } \\
HF_*(X, \omega) \widehat{\underset{\Lambda^\Pi}{\otimes}} HF_*(X, \omega) \ar[rr]_-{\Phi^{\rm PD}} & & \Lambda^\Pi }
\eeqn

\end{enumerate}
\end{thma}

\subsubsection{PSS map and isomorphism}
Second, we consider the Piunikhin--Salamon--Schwarz (PSS) map (cf. \cite{PSS}) from the Morse complex to the Floer complex and a map in the reversed direction called the SSP map. 

\begin{thma}[PSS and SSP maps] \label{thma_PSS}
Let $(f, g)$ be a Morse--Smale pair on $X$ with the associated Morse--Smale--Witten complex $CM_*(f, g)$. 
\begin{enumerate}

\item Let $CF_*(H, J, \Xi; \Lambda^\Pi)$ be a Floer chain complex provided by Theorem \ref{thma_Floer_complex}. Then there exists a $\Lambda^\Pi$-linear chain map
\beqn
\Phi^{\rm PSS}: CM_*(f, g) \underset{\mb Z}{\otimes} \Lambda^\Pi \to CF_*(H, J, \Xi; \Lambda^\Pi) 
\eeqn
of degree $-n:= - \frac{1}{2} {\rm dim} X$ 
and a $\Lambda^\Pi$-linear chain map
\beqn
\Phi^{\rm SSP}: CF_*(H, J, \Xi; \Lambda^\Pi) \to CM_*(f, g) \underset{\mb Z}{\otimes} \Lambda^\Pi
\eeqn
of degree $n$ whose definitions depend on various choices. 

\item The chain homotopy classes of $\Phi^{\rm PSS}$ and $\Phi^{\rm SSP}$ are independent of these choices.

\item The PSS and SSP maps are compatible with continuation maps, hence induce well-defined maps
\beqn
\Phi^{\rm PSS}: H_*(X; \Lambda^\Pi) \to HF_{*-n}(X, \omega; \Lambda^\Pi)
\eeqn
\beqn
\Phi^{\rm SSP}: HF_*(X, \omega; \Lambda^\Pi) \to H_{*+n}(X; \Lambda^\Pi).
\eeqn

\end{enumerate}
\end{thma}

\begin{thma}\label{thma_PSS_isomorphism}
The PSS map on homology is an isomorphism of graded $\Lambda^\Pi$-modules.
\end{thma}

\begin{rem}
    We point out that we do not know if the PSS and the SSP maps are inverses to each other up to homotopy. In particular, the proof of Theorem \ref{thma_PSS_isomorphism} involves new moduli spaces beyond those appearing in the constructions of Theorem \ref{thma_PSS}.
\end{rem}

Of course, the Arnold conjecture over integer coefficients now follows from Theorem \ref{thma_PSS_isomorphism} as a corollary, although in \cite{Bai_Xu_Arnold} we proved it by only considering the Floer chain complex for special $H$ and proving only the injectivity part of Theorem \ref{thma_PSS_isomorphism}.

\begin{cor}[Integral Arnold conjecture, cf. \cite{Bai_Xu_Arnold}]
Let $\phi: X \to X$ be a nondegenerate Hamiltonian diffeomorphism. Let $N$ be the minimal Chern number of $(X, \omega)$. Then the number of fixed points of $\phi$ is no less than the minimal number of generators of a ${\mb Z}/2N$-graded complex of free abelian groups whose homology is isomorphic to $H_*(X; {\mb Z})$ viewed as a ${\mb Z}/2N$-graded abelian group.
\end{cor}

If we compose the Floer-theoretic Poincar\'e pairing with the PSS map, then it induces a nondegenerate bilinear pairing on the homology $H_*(X; \Lambda^\Pi)$. 

\begin{defn}
The {\bf FOP Poincar\'e pairing} on $H(X; \Lambda^\Pi)$ is the bilinear pairing
\beqn
\langle \Phi^{\rm PSS} (\cdot), \Phi^{\rm PSS} (\cdot) \rangle: H (X; \Lambda^\Pi) {\underset{\Lambda^\Pi}{\otimes}}  H (X; \Lambda^\Pi) \to \Lambda^\Pi.
\eeqn
\end{defn}

If we use the PSS map defined using the technology of \cite{Fukaya_Ono}\cite{Liu_Tian_Floer}\cite{Pardon_virtual}, then over rational coefficients, such a pairing coincides with the classical Poincar\'e pairing on homology. This is because the orbifold Euler classes satisfy properties needed for $S^1$-localization, so there is no contribution from moduli spaces of curves with nonzero energy. However, the FOP Euler characteristic classes arising in our perturbation scheme do not necessarily satisfy these properties and we do not know if the pairing is classical or not. This phenomenon is similar to the case of quantum $K$-theory (see \cite{Lee-K}) or quantum generalized cohomology developed by Abouzaid--McLean--Smith \cite{AMS2}. The point of introducing the quantum-corrected Poincar\'e pairing is to define an associative product analogous to the quantum product on $H_*(X; \Lambda^\Pi)$ in our context, where we use the pairing to turn three-pointed Gromov--Witten type counts into a bilinear operation on $H_*(X; \Lambda^\Pi)$ to dualize a marked point into an output marked point.

Another corollary is the invertibility of the Seidel map with integer coefficients defined in \cite{BPX}. Recall that for each loop of Hamiltonian diffeomorphisms $\phi: S^1 \to {\rm Ham}(X)$, one has a Hamiltonian fibration $\wt M_\phi \to S^2$ with fibres being $M$. By counting holomorphic sections $u: S^2 \to \wt M_\phi$ (via FOP perturbation), one can define a map
\beqn
S(\phi): H(X; \Lambda^\Pi) \to H(X; \Lambda^\Pi)
\eeqn
which only depends on the homotopy class $[\phi] \in \pi_1({\rm Ham}(X))$.\footnote{One uses the classical Poincar\'e duality in the definition. Hence we do not expect that $S(\phi) S(\psi) = S(\phi \psi)$. However, the main use in \cite{BPX} of the Seidel map is its invertibility.} In the semipositive case, using classical perturbation method, one can prove that $S(\phi)$ is invertible by showing $S(\phi \psi) = S(\phi) S(\psi)$ and $S(1) = {\rm Id}$. However, the invertibility of the Seidel maps can only be shown using the Hamiltonian model for the quantum cohomology.

\begin{cor}[Seidel representation]
The Seidel map $S(\phi)$ is invertible.
\end{cor}

\begin{proof}
Fix an almost complex structure $J$. Choose a nondegenerate Hamiltonian $H_t$. Let $H_t'$ be the pullback by the loop $\phi_t$. On $S^2 = {\mb R} \cup S^1 \cup \{\pm\infty\}$ one can turn on a family of Hamiltonian connections parametrized by $S \in [0, \infty)$ which is $H_t dt$ on $[-S, S] \times S^1$ and which vanishes near $\pm\infty$. We regard this as a family of Hamiltonian connections $\sigma_S$ on $\wt M_\phi \to S^2$ written with respect to the negative trivialization $\psi_-$. Then with respect to the positive trivialization, the connection reads $H' dt$. Then consider the family of moduli spaces of holomorphic sections with respect to this family of Hamiltonian connections. When $S = 0$, one obtains a bimodule defining the Seidel map $S(\phi)$, given as a chain map
\beqn
\tilde S(\phi): CM (f) \to CM(f).
\eeqn
The induced map on homology is the Seidel map $S(\phi)$. When $S = +\infty$, one obtains a triple concatenation of a bimodule for PSS map, a bimodule for continuation map, and a bimodule for SSP map. Therefore, $\tilde S(\phi)$ is chain homotopic to the composition of these three maps. By Theorem \ref{thma_continuation_map} and Theorem \ref{thma_PSS}, on cohomology level $S(\phi)$ is invertible.
\end{proof}

\subsubsection{Comparison with known construction}

When $(X, \omega)$ is weakly monotone (also called semi-positive), then by perturbing either the almost complex structure $J$ or the Hamiltonian $H$, Hofer--Salamon \cite{Hofer_Salamon} and Ono \cite{Ono_1995} defined the Floer chain complex over integers. Such constructions can also be carried out to define continuation maps, pair-of-pants products, and K\"unneth isomorphisms with integer coefficients for weakly monotone $(X, \omega)$. It is natural to wonder if our constructions agree with the classical definitions in the weakly monotone setting. As a proof of principle, we demonstrate that this is the case for the Floer homology.

\begin{thma}\label{thma_semipositive}
When $(X, \omega)$ is weakly monotone, our Floer chain complex $CF_*(H, J, \Xi)$ is chain homotopy equivalent to the one defined using the classical method.
\end{thma}

We expect that the proof methods can be generalized to establish the identification between our Floer-theoretic invariants and the classical versions.

\begin{rem}
On the other hand, for a general compact $(X, \omega)$, Floer chain complex $CF_*^{\mb Q}(H, J)$ over rational coefficients is constructed using the classical virtual technique (see \cite{Liu_Tian_Floer}\cite{Fukaya_Ono}\cite{Pardon_virtual}). Using the setup similar to here, $CF_*^{\mb Q}(H, J)$ is a module over the rational Novikov field $\Lambda_{\mb Q}^\Pi$. Although there is an isomorphism
\beqn
HF_*(X, \omega; \Lambda^\Pi) \underset{\Lambda^\Pi}{\otimes} \Lambda_{\mb Q}^\Pi \cong HF_*^{\mb Q}(X, \omega)
\eeqn
via the PSS maps in the two setups and the universal coefficient theorem for ordinary homology, one does not have an obvious relation on the chain level, between
\beqn
CF_*(H, J, \Xi; \Lambda^\Pi) \underset{\Lambda^\Pi}{\otimes} \Lambda_{\mb Q}^\Pi\ \ \ {\rm and}\ \ \ CF_*^{\mb Q}(H, J)
\eeqn
without going through Morse theory.
\end{rem}

\subsection{Main results on equivariant Floer theory and Steenrod operations}

To discuss Steenrod operations and other aspects of equivariant Floer theory, we switch to the cohomological convention, partially to match with the fact that classical Steenrod operations are cohomological operations. In this way, by considering ascending flows rather than descending flows, one obtains the Floer cochain complex $CF^*(H, J, \Xi)$ associated to a Floer flow category whose cochains are formal infinite linear combinations of capped 1-periodic orbits of $H$ whose actions go to $+\infty$.

Let $\p$ be a prime number. Abbreviate $\zp:= {\mb Z}/\p {\mb Z}$. When we view it as a field serving as the coefficient field of cohomology, we use $\fp$ instead. In characteristic $\p$, on any topological space, there are the Steenrod power operations which can be organized as a single nonlinear map. 

We will primarily use the ``Tate'' version of equivariant theory. Recall that the more familiar ``Borel'' version of the $\zp$-equivariant cohomology of a point is 
\beqn
H^{\rm Borel}_\zp ({\rm pt}; \fp) = \fp[\tu]\langle \uptheta\rangle
\eeqn
over $\fp$ generated by a degree $2$ variable $\tu$ and a degree $1$ variable $\uptheta$ satisfying 
\begin{align*}
&\ \tu\uptheta = \uptheta \tu,\ &\ \uptheta^2 = \left\{ \begin{array}{cc} \tu,\ & \p = 2,\\
                          0,\ & \p > 2. \end{array}\right.
\end{align*}
It is defined as the ordinary cohomology of the classifying space $B\zp$ and admits a nice Morse model. We often need to take the completion by allowing formal power series. 

The Tate version of the equivariant cohomology of a point differs from the Borel one by inverting the variable $\uptheta$ (in characteristic $2$) resp. $\tu$ (in odd characteristics), i.e. (after completion)
\beqn
H_{\zp}^{\rm Tate}({\rm pt}; \fp) = \left\{ \begin{array}{ll} {\bf F}_2 \lst,\ &\ \p = 2,\\
\fp \lsu \langle \uptheta \rangle,\ &\ \p > 2. \end{array}\right.
\eeqn

\begin{rem}
There are two main reasons why we prefer the Tate version. First, all the nice properties of the equivariant Floer theory hold for the Tate version but not the Borel one, for example, localization. Second, the inversion of the variable $\tu$ follows more naturally from our geometric construction. The Borel version, instead, is regarded as a truncation of the Tate version. From now on, all equivariant cochain complexes or cohomology, without further declaration, will be the Tate version. 
\end{rem}

We introduce the following notions. Let 
\beqn
{\bf K}_\p:= \left\{ \begin{array}{ll} {\bf F}_2\lst,\ &\ \p = 2,\\
\fp \lsu,\ &\ \p > 2
\end{array}\right.
\eeqn
The chain-level construction we use will always produce algebraic structures linear over ${\bf K}_\p$. On the other hand, denote
\beqn
\tkp:= \left\{ \begin{array}{ll} {\bf K}_\p,\ &\ \p = 2,\\
{\bf K}_\p \langle \uptheta \rangle ,\ &\ \p> 2 \end{array} \right. .
\eeqn
The resulting cohomology groups (not complexes) will be modules over $\tkp$. However, such module structure is not an algebraic consequence, but needs to be constructed geometrically.

\subsubsection{The quantum Steenrod operations}

The classical Steenrod square resp. power operations are cohomological operations which can be formulated as a single map (for $\p = 2$)
\beqn
St_2: H^*(X; {\bf F}_2) \to H^*(X; {\bf F}_2[\uptheta])\subset  H^*(X; \tkp)
\eeqn
resp. ($\p > 2$)
\beqn
St_\p: H^*(X; \fp) \to H^*(X; \fp[\tu] \langle \uptheta \rangle ) \subset H^*(X; \tkp).
\eeqn
The map $St_\p$ is additive and linear over $\fp$. Moreover, if we regard $\tu$ of degree $2$ and $\uptheta$ of degree $1$, then $St_\p$ multiplies the grading by $\p$. Fukaya proposed the idea of defining quantum Steenrod operation \cite{Fukaya_Morse}. It has been adapted in the Floer setting by Seidel \cite{Seidel_pants} and Shelukhin--Zhao \cite{shelukhin-zhao} and in the quantum cohomology setting by Wilkins \cite{Wilkins_2020} and Seidel--Wilkins \cite{Seidel_Wilkins_2022}. A key point in the construction is the integrality of curve counts, which then admit well-defined mod $\p$ reduction. Hence the existing constructions are all restricted to monotone or semipositive settings. The FOP perturbation method developed in \cite{Bai_Xu_2022} allows one to extend to the general setting. 

\begin{thma}[Quantum Steenrod operation on general symplectic manifold]\label{thma_QST} For any compact symplectic manifold $(X, \omega)$, for each prime $\p$, there are well-defined additive maps
\beqn
\begin{split}
QSt_2: &\ H^k (X; {\bf F}_2) \to H^{2k}( X; \Lambda_{{\bf F}_2 \pst}^\Pi) \subset H^{2k} (X; \Lambda_{\tilde {\bf K}_2}^\Pi),\ (\p = 2)\\
QSt_\p: &\ H^k(X; \fp) \to H^{\p k}(X; \Lambda_{\fp \psu}^\Pi \langle \uptheta \rangle ) \subset H^{\p k}(X; \Lambda_{\tkp}^\Pi)\ (\p > 2)
\end{split}
\eeqn
satisfying the following conditions.
\begin{enumerate}

\item The classical part (i.e., setting the Novikov parameter to be zero) of $QSt_\p$ coincides with the classical Steenrod operation.

\item When $(X, \omega)$ is semipositive, $QSt_\p$ coincides with the quantum Steenrod operation defined by Seidel--Wilkins \cite{Seidel_Wilkins_2022} (which, in $\p=2$ and monotone case, agrees with Wilkins' definition of quantum Steenrod square in \cite{Wilkins_2020}).
\end{enumerate}
\end{thma}

\subsubsection{Equivariant Floer cohomology}

The Floer-theoretic version of the Steenrod operation has played important roles in many recent dynamical applications of Floer theory, see \cite{Seidel_pants} \cite{shelukhin20, shelukhin-21, shelukhin-22} \cite{CGG} etc. In this paper, we also generalize the construction of the Floer-theoretic Steenrod operation to general compact symplectic manifolds. Let $CF^*(H)$ be a Floer chain complex associated to a nondegenerate Hamiltonian $H$, a compatible almost complex structure $J$, and a system of choices $\Xi$; to simplify notation, we suppress $J$ and $\Xi$ from the notation.

\begin{thma}[Equivariant Floer cohomology]\label{thma_equivariant_Floer}
Let $(X,\omega)$ be a compact symplectic manifold and $\p$ be a prime.

\begin{enumerate}

\item Let $H$ be a 1-periodic Hamiltonian such that its $\p$-th iteration $H^{\natural\p}$ is nondegenerate. Then there exists a cochain complex (the Tate $\zp$-equivariant Floer complex)
    \beqn
    CF_{\zp} (H^{\natural \p}, \Xi^\zp)
    \eeqn
    of finitely generated $\Lambda_{\kp}^\Pi$-modules whose definition depends on a choice $\Xi^\zp$.

    \item For Hamiltonians $H_1, H_2$ with both $H_1^{\natural \p}$ and $H_2^{\natural \p}$ nondegenerate and choices $\Xi_1^\zp$ and $\Xi_2^\zp$ to define their equivariant Floer cochain complex, there is a chain homotopy equivalence
    \beqn
    CF_\zp (H_1^{\natural\p}, \Xi_1^\zp) \to CF_{\zp} (H_2^{\natural \p}, \Xi_2^\zp)
    \eeqn
    called the equivariant continuation map.

    \item The equivariant continuation maps are well-defined up to chain homotopy. Moreover, the composition of two equivariant continuation maps is homotopic to an equivariant continuation map.     Hence there is a well-defined symplectic invariant of $(X, \omega)$, denoted by 
    \beqn
    HF_\zp (X)
    \eeqn
    called the (Tate) {\bf $\zp$-equivariant Hamiltonian Floer cohomology}.

    \item There is a natural equivariant PSS map
    \beqn
    \Phi_{\zp}^{\mb{MF}}: H^*(X; \Lambda_{\tkp}^\Pi) \to HF_{\zp}^*(X)
    \eeqn
    and an equivariant SSP map
    \beqn
    \Phi_{\zp}^{\mb{FM}}: HF_{\zp}^*(X) \to H^*(X; \Lambda_{\tkp}^\Pi)
    \eeqn
    which are isomorphisms of $\Lambda_{\kp}^\Pi$-modules.
\end{enumerate}
\end{thma}

\begin{rem}
The equivariant cohomology $HF_{\zp}(X)$ is naturally a module over the field $\Lambda_{\kp}^\Pi$. When $\p=2$, this is the same as the ring $\Lambda_{\tkp}^\Pi$. When $\p>2$, there is no $\Lambda_{\tkp}^\Pi$-module structure on the complex $CF_{\zp}(H^{\natural \p})$ as ``multiplying by $\uptheta$'' is not a chain map. However, one could define a $\Lambda_{\tkp}^\Pi$-module structure on the cohomology $HF_{\zp}(H^{\natural \p})$. However, we do not consider this additional structure in the current paper. 
\end{rem}

There is an algebraically defined $\zp$-equivariant chain complex 
\beqn
C_{\zp}( CF (H)^{\otimes \p})
\eeqn
whose underlying cochain group is a suitable completion of 
\beqn
\Big( CF (H) \underset{\Lambda_{\fp}^\Pi}{\otimes} \cdots \underset{\Lambda_{\fp}^\Pi}{\otimes} CF (H) \Big) \otimes \Lambda_{\tkp}^\Pi.
\eeqn
The differential is linear over $\Lambda_{\kp}^\Pi$. The quasi-Frobenius map is an $\fp$-linear map 
\beqn
qF: HF(H; \Lambda_{\fp}^\Pi) \to H_{\zp} ( CF(H)^{\otimes \p} )
\eeqn
which is, roughly speaking, induced by the chain-level map $x \mapsto x\otimes \cdots \otimes x$. 

\subsubsection{Floer-theoretic Steenrod operation}

\begin{thma}[Floer-theoretic Steenrod operation]\label{thma_Floer_Steenrod}
Let $(X, \omega)$ be a compact symplectic manifold and $J$ be an $\omega$-compatible almost complex structure.
\begin{enumerate}

    \item For any pair of Hamiltonians $H_-$ and $H_+$ with $H_-$ nondegenerate and $H_+^{\natural \p}$ nondegenerate, there exists a cochain map called the {\bf equivariant pair-of-pants product}
    \beqn
    \wt P: C_{\zp}(CF (H_-)^{\otimes \p}) \to CF_\zp (H_+^{\natural\p})
    \eeqn
    It is well-defined up to chain homotopy.
    
    \item The map on cohomology, denoted by 
    \beqn
    P: H_{\zp}(CF(H_-)^{\otimes \p}) \to HF_{\zp}(H_+^{\natural \p})
    \eeqn
    is an isomorphism of $\Lambda_{\kp}^\Pi$-modules.

    \item The composition defining the {\bf Floer-theoretic Steenrod operation} is given by 
    \beqn
    \xymatrix{  HF(X; \Lambda_{\fp}^\Pi) \ar[r] & H(CF(H_-) \otimes \Lambda_{\fp}^\Pi) \ar[r]^{qF} & H_{\zp}(CF(H_-)^{\otimes \p}) \ar[r]^{P} & HF_{\zp}(H_+^{\natural \p}) \ar[r] & HF_{\zp}(X) }
    \eeqn
    denoted by 
    \beqn
    \fst: HF(X; \Lambda_{\fp}^\Pi) \to HF_{\zp}(X)
    \eeqn
    is independent of choices.
    
    \item The following diagram commutes: when $\p=2$; \beq
    \vcenter{ 
    \xymatrix{      H  (X; \Lambda_{{\bf F}_2}^\Pi )  \ar[rr]^{\qst}  \ar[d]_{\Phi^{\mb{MF}}} & & H(X; \Lambda_{\fp \pst}^\Pi) \\
    HF (X; \Lambda_{{\bf F}_2}^\Pi) \ar[rr]_-{{\it FSt}_2} &  & HF{}_\zp (X; \Lambda_{{\bf F}_2 \pst}^\Pi ) \ar[u]_{\Phi_{\zp}^{\mb{FM}}}} }
    \eeq
    when $\p>2$,    \beq\label{qst_fst_commute}
    \vcenter{ 
    \xymatrix{      H  (X; \Lambda_{\fp}^\Pi )  \ar[rr]^{\qst}  \ar[d]_{\Phi^{\mb{MF}}} & & H(X; \Lambda_{\fp \psu}^\Pi \langle \uptheta \rangle) \\
    HF (X; \Lambda_{\fp}^\Pi) \ar[rr]_-{\fst} &  & HF{}_\zp (X; \Lambda_{\fp \psu}^\Pi \langle \uptheta \rangle ) \ar[u]_{\Phi_{\zp}^{\mb{FM}}}} }
    \eeq

\end{enumerate}
\end{thma}

\subsection{Main results on quantitative Floer theory}

Many successful applications of Hamiltonian Floer theory have a quantitative flavor. We also have the quantitative package under our framework. As many applications take advantage of using different coefficients, we phrase the statements using a general Novikov ring/field.

\subsubsection{Action filtration and spectral invariants}

We state the following theorem for the convenience for readers to cite. The proof has no difference from cases proved using other approaches is hence omitted. 

\begin{thma}[Filtered Floer complex and continuation map]\label{thma_Floer_filtered} Let $R$ be a commutative ring with unit and let $CF (H; \Lambda_R^\Pi)$ be a Floer complex provided by Theorem \ref{thma_Floer_complex} tensored with $\Lambda_R^\Pi$.
\begin{enumerate}

\item The differential of the Floer complex preserves the energy filtration. More precisely, for each $\tau \in {\mb R}$, the differential preserves the subgroups
\beqn
CF (H; \Lambda_R^\Pi)^{\leq \tau}:= \Big\{ \sum_{i=1}^\infty a_i p_i \in CF (H; \Lambda_R^\Pi)\ |\ {\mc A}_H(p_i) \leq \tau \Big\}
\eeqn
and 
\beqn
CF (H; \Lambda_R^\Pi)^{<\tau}:= \Big\{\sum_{i=1}^\infty a_i p_i\in CF (H; \Lambda_R^\Pi)\ |\ {\mc A}_H(p_i) < \tau \Big\}.
\eeqn

\item Let $CF(H_1; \Lambda_R^\Pi)$ and $CF (H_2; \Lambda_R^\Pi)$ be two Floer complexes. Denote
\beqn
c:= \int_{S^1} \sup_X (H_{1, t} - H_{2, t}) dt.
\eeqn
Then there exists a continuation map $\Phi^{12}(\Theta)$ such that for all $\tau \in {\mb R}$,
\begin{align*}
&\ \Phi^{12}(\Theta) ( CF(H_1; \Lambda_R^\Pi)^{\leq \tau}) \subset CF (H_2; \Lambda_R^\Pi )^{\leq \tau + c},\ &\ \Phi^{12}(\Theta) (CF(H_1; \Lambda_R^\Pi)^{<\tau}) \subset CF(H_2; \Lambda_R^\Pi)^{<\tau + c}.
\end{align*}
\end{enumerate}
\end{thma}

The filtered version of continuation map implies that for any nondegenerate Hamiltonian $H$, the filtered homology groups
\begin{align*}
&\ HF (H; \Lambda_R^\Pi)^{\leq \tau},\ &\ HF (H; \Lambda_R^\Pi)^{<\tau}
\end{align*}
only depend on $H$. When $R$ is fixed in the local context, we omit the coefficient from the notations. Moreover, there are well-defined maps
\beqn
HF(H)^{<\tau} \to HF(X)
\eeqn
compatible with maps (for $\tau_1 \leq \tau_2$)
\beqn
HF(H)^{<\tau_1} \to HF(H)^{<\tau_2}.
\eeqn
This allows us to define the {\bf spectral invariants} of Floer homology classes.

\begin{defn}
Given a Hamiltonian $H$, the {\bf spectral invariant} of $x \in HF(X; \Lambda_R^\Pi)$ is 
\beqn
c_H(x):= \lim_{i \to \infty} \Big\{ \tau\in {\mb R}\ |\ x \in {\rm Im} \big( HF(H_i)^{<\tau} \to HF (X) \big) \Big\} \in {\mb R} \cup \{-\infty\}.
\eeqn
Here the limit is taken for a sequence of nondegenerate Hamiltonians $H_i$ converging to $H$ in $C^0$-norm.
\end{defn}

\begin{thma}[Spectral invariants]\label{thma_spectral_invariants}
Fix a Noetherian ring $R$. The spectral numbers $c_H (x)$ for $x \in HF(X; \Lambda_R^\Pi)$ satisfy the following properties.
\begin{enumerate}

\item {\bf (Spectrality)} For any $x \in HF(X; \Lambda_R^\Pi)$, the number
\beqn
c_H(x) = -\infty \Longrightarrow x = 0
\eeqn
and for $x \neq 0$, there exists a cycle ${\mf x}$ representing $x$ such that 
\beqn
{\mc A}_H({\mf x}) = c_H(x).
\eeqn

\item {\bf (Shifting Property)} For $x \in HF (X; \Lambda_R^\Pi)$ and $a \in \Lambda_R^\Pi$, one has 
\beqn
c_H(a x) = c_H(x) - {\mf v}(a).
\eeqn

\item {\bf (Lipschitz Continuity)} $c_H(x)$ is Lipschitz continuous with respect to the Hofer metric on the space of 1-periodic Hamiltonians.

\item {\bf (Isotopy Invariance)} $c_H(x)$ only depends on the element $\tilde \phi$ in $\wham (X, \omega)$, the Hamiltonian isotopy induced from $H$ as an element of the universal cover of the Hamiltonian diffeomorphism group. Hence we can denote $c_H(x)$ by $c_{\tilde\phi}(x)$.

\item {\bf (Symplectic Invariance)} For any symplectomorphism $\rho: X \to X$, one has $c_{H\circ \rho}(x)  = c_H(\rho_*(x))$ where $\rho_*$ is the automorphism
\beqn
\rho_*: HF(X) \to HF(X)
\eeqn
induced by $\rho$.

\item {\bf (Triangle Inequality)} For any $x,y \in HF(X)$, $\tilde\phi, \tilde \psi \in \wham (X, \omega)$, one has 
\beqn
c_{\tilde\phi\tilde\psi}(x \star y )  \leq c_{\tilde\phi}(x) + c_{\tilde\psi}(y).
\eeqn
Here $\star$ is the pair-of-pants product on Floer homology defined by Theorem \ref{thma_product}.

\item {\bf (Poincar\'e Duality)} Let $R = {\bf K}$ be a field. Let $\Phi_0^{\rm PD}$ be the composition
\beqn
\xymatrix{ HF(X; \Lambda_{\bf K}^\Pi) \otimes HF(X; \Lambda_{\bf K}^\Pi) \ar[r]^-{\Phi^{\rm PD}}    &  \Lambda_{\bf K}^\Pi \ar[r]  &   {\bf K} }.
\eeqn
Let $\ov{H}$ be the reversal of the Hamiltonian $H$. Then one has
\beqn
c_H(x) = - \inf \Big\{ c_{\ov{H} }( y  )\ |\ y \in HF (X; \Lambda_{\bf K}^\Pi),\  \Phi_0^{\rm PD} ( x, y) \neq 0 \Big\}.
\eeqn

\end{enumerate}
\end{thma}

\subsubsection{Barcodes}

Another quantitative structure widely used in applications is the barcode. We follow the approach of Usher--Zhang \cite{Usher_Zhang_2016} which requires slightly modified setup. 

For a symplectic manifold $X$ with the abelian group $\Pi$ defined by \eqref{Novikov_group}, a (concise) {\bf barcode} (see \cite{Usher_Zhang_2016}) is a finite multiset of elements of $({\mb R}/\omega(\Pi)) \times (0, +\infty]$. 

It is most useful when considering Floer barcodes with field coefficients. Then the complex
\beqn
CF (H; \Lambda_{\bf K}^{\Pi}):= CF(H; \Lambda_{\mb Z}^\Pi) \otimes \Lambda_{\bf K}^{\Pi}
\eeqn
is a (slight generalization of) Floer-type complex in the sense of \cite[Definition 4.1]{Usher_Zhang_2016}. It follows from \cite[Definition 6.3, Theorem 7.1]{Usher_Zhang_2016} that there is a barcode canonically associated to the complex $CF(H; \Lambda_{\bf K}^{\Pi})$. 

\begin{cor}[Floer barcodes] The barcode associated to the Floer complex $CF(H; \Lambda_{{\bf K}}^{\Pi})$ only depends on the Hamiltonian $H$. Moreover, the dependence on $H$ is Lipschitz with respect to the bottleneck distance on barcodes and the Hofer metric on the space of $1$-periodic Hamiltonians.
\end{cor}
\begin{proof}
    It follows from applying Proposition \ref{prop738} and Theorem \ref{thm_Usher_Zhang_barcodes} to the filtered Floer complex constructed in Theorem \ref{thma_Floer_filtered}.
\end{proof}

\subsubsection{Quantitative features of equivariant Floer theory}

The equivariant Floer cochain complex also has an energy filtration. Notice that the equivariant cochain complex $CF_{\zp} (H^{\natural \p})$ is a module over $\Lambda_{\kp}^\Pi$. We can prove that it is a Floer-type complex. On the other hand, the ordinary Floer complex $CF(H; \Lambda_{\fp}^\Pi)$ also induces an equivariant complex coming from its own $\p$-fold tensor power, denoted by 
\beqn
C_{\zp}(CF(H)^{\otimes \p})
\eeqn
which is also a Floer-type complex.

\begin{thma}[Barcode stretching under prime iteration]\label{thm_pants_filtered}
Let $X$ be a compact symplectic manifold, $\p$ be a prime, and $H$ be a Hamiltonian on $X$ such that the $\p$-th iteration $H^{\natural \p}$ is nondegenerate.

\begin{enumerate}

\item The equivariant Floer cochain complex $CF_{\zp} (H^{\natural \p})$ is a Floer-type complex over the Novikov field $\Lambda_{\kp}^\Pi$. 

\item The equivariant pair-of-pants product has a chain-level realization
\beqn
P: C_{\zp}(CF(H)^{\otimes \p}) \to CF_{\zp}(H^{\natural \p})
\eeqn
which is a filtered chain homotopy equivalence for Floer-type complexes.

\item Suppose the bar-length spectrum of $CF(H)$ is $0 < \beta_1(H) \leq \cdots \leq \beta_m(H)$. Then the bar-length spectrum of $CF_{\zp}(H^{\natural \p})$ is 
\beqn
\beta_{i, {\mb Z}_2}(H^{\natural 2}) = 2 \beta_i(H),\ i = 1, \ldots, m
\eeqn
when $\p = 2$ and 
\beqn
\beta_{2i-1, \zp}(H^{\natural \p}) = \beta_{2i, \zp}(H^{\natural \p}) = \p \beta_i(H),\ i = 1, \ldots, m
\eeqn
when $\p > 2$.
\end{enumerate}
\end{thma}

\subsubsection{Quantitative Hamiltonian Floer theory for noncontractible orbits}

We state the main quantitative results about Floer theory for noncontractible orbits. We follow very much the treatment of Sugimoto \cite{Sugimoto_2026}. Since the resulting homology group always vanishes, there are no spectral invariants; hence the barcodes are the main quantitative aspects here. 

Choose a possibly noncontractible loop $\gamma: S^1 \to X$ which represents a nonzero class
\beqn
[\gamma] \neq 0 \in H_1(X; {\mb Z})_{\rm free}.
\eeqn
One can consider noncontractible periodic orbits of a Hamiltonian $H$ which are homologous (over ${\mb Z}$) to $\gamma$. 
\begin{defn}\label{defn_gamma_capping}
Let $H$ be a 1-periodic Hamiltonian on $(X, \omega)$.

\begin{enumerate}
\item Let $x: S^1 \to X$ be a 1-periodic orbit of $H$. A {\bf $\gamma$-capping} of $x$ is an oriented cobordism $u: \Sigma \to X$ from $\gamma$ to $x$ (where $\Sigma$ could be a higher genus surface).

\item Two $\gamma$-cappings are equivalent if their difference, which is a class in $H_2(X; {\mb Z})$, has zero symplectic area. A {\bf $\gamma$-capped orbit} of $H$ is a 1-periodic orbit together with an equivalence class of $\gamma$-cappings. 

\item Let 
\beqn
\tilde {\mc O}_\gamma(H)
\eeqn
denote the set of $\gamma$-capped orbits of $H$, on which there is a free action by the group
\beqn
\Gamma:= H_2(X; {\mb Z})/ {\rm ker} (\omega).
\eeqn

\end{enumerate}
\end{defn}

In our setup we do not keep the ${\mb Z}$-grading but only consider the ${\mb Z}_2$-grading on $\gamma$-capped orbits. Each $\gamma$-capped orbit $p$ has a well-defined symplectic action ${\mc A}_{H,\gamma} (p)\in {\mb R}$. The construction of the filtered Floer complex is summarized below and will be proved in Section \ref{section_noncontractible}.

\begin{thm}[Filtered Floer chain complex for noncontractible periodic orbits]\label{thm_Floer_noncontractible}
Let $(X, \omega)$ be a compact symplectic manifold and let $\gamma:S^1 \to X$ be a smooth loop with $[\gamma] \neq 0 \in H_1(X; {\mb Z})_{\rm free}$. Let $H$ be a nondegenerate $1$-periodic Hamiltonian. Then for each $1$-periodic family of $\omega$-compatible almost complex structures $J$, there is a differential on the ${\mb Z}/2$-graded abelian group
\beqn
CF_\gamma ( H):= \Big\{ \sum_{i=1}^\infty a_i p_i\ |\ a_i \in {\mb Z},\ p_i \in \tilde {\mc O}_\gamma(H),\ \lim_{i \to \infty} {\mc A}_{H,\gamma} (p_i) = -\infty \Big\}
\eeqn
leading to a ${\mb Z}/2$-graded  complex of abelian groups
\beqn
CF_\gamma(H, J, \Xi).
\eeqn
Moreover, for each field ${\bf K}$, the tensor product
\beqn
CF_\gamma(H) \otimes \Lambda_{\bf K}^\Gamma
\eeqn
is a Floer-type complex over $\Lambda_{\bf K}^\Gamma$ whose filtered isomorphism class only depends on $H$.
\end{thm}

Under the assumption that $[\gamma] \neq 0 \in H_1(X)$, we expect that a continuation map argument can prove the vanishing of the resulting homology group. The main use of the noncontractible Floer chain complex is on the quantitative side. The main purpose of developing this theory is to generalize the work of Sugimoto \cite{Sugimoto_2026}, which uses the equivariant pair-of-pants product as a primary tool. We state the main result regarding the quantitative implications, which is also proved in Section \ref{section_noncontractible}.

\begin{thm}\label{thm_noncontractible_pants}
Under the assumption of Theorem \ref{thm_Floer_noncontractible}, let $\p$ be a prime and suppose $H^{\natural \p}$ is nondegenerate. 
\begin{enumerate}
\item There exists a Floer-type complex
\beqn
CF_{\zp}(H^{\natural \p}, \gamma)
\eeqn
over $\Lambda_{\kp}^\Gamma$ which is freely generated by $\p$-periodic orbits homologous to $\gamma$. Its filtered isomorphism class only depends on $H$.

\item There exists a filtered homotopy equivalence
\beqn
P_\gamma: C_{\zp}(CF_\gamma(H)^{\otimes \p}) \to CF_{\zp}(H^{\natural \p}, \gamma^{(\p)})
\eeqn
of Floer-type complexes over $\Lambda_{\kp}^\Gamma$.

\item Let 
\beqn
\beta_1^\gamma (H) \leq \cdots \leq \beta_N^\gamma (H)
\eeqn
be the bar-length spectrum of $CF_\gamma(H) \otimes \Lambda_{\fp}^\Gamma$. Then it only depends on $H$ and the homology class of $\gamma$. Moreover, the bar-length spectrum of $CF_{\zp, \gamma^{(\p)}} (H^{\natural \p})$ is 
\beqn
2\beta_1^\gamma(H), \ldots, 2 \beta_N^{\gamma}(H)
\eeqn
when $\p = 2$ (with the same number of finite bars) and is 
\beqn
\p \beta_1^\gamma(H), \p \beta_1^\gamma(H), \ldots, \p \beta_N^\gamma(H), \p \beta_N^\gamma(H)
\eeqn
when $\p>2$ (with twice the number of finite bars).
\end{enumerate}
\end{thm}

\subsection{What we do ``virtually''}

We would like to explain the general principle of our construction (on the virtual side) along with some historical remarks.

In a very general sense, the Floer complexes and other chain-level structures are constructed via counting solutions to the Floer equation. To count beyond the cases where transversality can be achieved by geometric perturbations (\cite{Floer_CMP}\cite{Hofer_Salamon}\cite{Ono_1995}), the so-called virtual technique was invented (see \cite{Li_Tian}\cite{Fukaya_Ono}\cite{FOOO_Kuranishi}\cite{HWZ_book}\cite{Pardon_virtual} etc.). Such constructions, depending on the extent of applications, usually require building a complete package of abstract theory (such as Kuranishi theory and polyfold theory); the application of these abstract theories can typically be broken down into two main components: (virtual) regularization and (virtual) perturbation/integration.

From now on we only discuss under the finite-dimensional/Kuranishi approach, the one upon which the current paper is based, rather than the infinite-dimensional/polyfold approach. The issues and challenges we will address in the finite-dimensional approach may not be manifested in the infinite-dimensional situation.

A central issue in the finite-dimensional construction is the choice-dependence. Unlike the algebro-geometric situation, one rarely obtains canonical objects such as a perfect obstruction theory (\cite{Behrend_Fantechi}). To obtain choice-independent results, one needs to include various kinds of equivalence relations in the abstract package, or perhaps higher categorical structures. The choice-independence is not yet a significant challenge in cases where one deals with a single or finitely many moduli spaces, such as Gromov--Witten theory. But it causes much more serious difficulties in Floer-type theories. In Floer theory, essentially one obtains virtual fundamental chains rather than cycles; the virtual fundamental chains satisfy expected properties (such as $d^2 = 0$) only if we make the constructions for infinitely many moduli spaces compatibly. It could be a delicate issue when developing the abstract package incorporating such compatibility conditions and when fulfilling these conditions in the concrete regularization procedure.\footnote{The work of Abouzaid \cite{Abouzaid_axiomatic} takes a different approach where the choices are made universally in his notion of ``theories of virtual counts.''} 

We would like to emphasize two important developments which precede the current piece of work. The first one, regarding the perturbation component, starts from \cite{Fukaya_Ono_integer} which proposed the principle leading to integer-valued curve counts even when the orbifold locus of a moduli space is obstructed. After technical exploration by Parker \cite{BParker_integer}, the authors rigorously established this new perturbation method in \cite{Bai_Xu_2022}, which we call the FOP perturbation scheme. The second one, regarding the regularization component, is the construction of global Kuranishi charts starting from \cite{AMS}. It provides a powerful method producing global Kuranishi models regularizing any single moduli space and greatly reduces the technicality when discussing the compatibility among local Kuranishi charts. In \cite{Bai_Xu_Arnold} we produced the global Kuranishi chart construction for Floer complexes and PSS maps, which is being extended to other chain-level structures in the current paper (with certain upgrades coming from \cite{AMS2} and \cite{Hirschi_Swaminathan}).

Now we summarize the technical framework underlying our implementation of the virtual technique. Compared to the regularization component, the perturbation component is to simply apply a blackbox theorem (Theorem \ref{thma_FOP}) which was already established in \cite{Bai_Xu_2022}. However, the delicate prerequisite for the blackbox theorem, along with other requirements, demands careful fine-tunings of the regularization component. Since each stage produces different kinds of objects whose compatibility conditions all look similar to flow category structure, each stage-wise construction is phrased as producing a flow category enriched in a particular ``regular stratification category.'' The objects of such a category are stratified in a way similar to manifolds with corners, with a monoidal structure corresponding to (direct) product and disjoint union. 

To be more concrete, the whole procedure we are describing can be summarized in the following zig-zag of functors
\beqn
\rlap{$\overbrace{\phantom{ \uds{\bf Top} \leftarrow \uds{\bf Kur}_{\rm finite}^{\rm NC} \rightarrow \uds{\bf dOrb}^{\rm NC}} }^{{\rm regularization}}$} \uds{\bf Top} \leftarrow \uds{\bf SKur}^{\rm NC} \rightarrow  \underbrace{\uds{\bf dOrb}^{\rm NC} \leftarrow \uds{\bf dOrb}^{\rm FOP} \rightarrow \pman}_{\rm perturbation}.
\eeqn
The items are all regular stratification categories where
\begin{itemize}

\item $\uds{\bf Top}$ is the category of stratified topological spaces;

\item $\uds{\bf SKur}^{\rm NC}$ is the category of {\it normally complex smooth Kuranishi charts};

\item $\uds{\bf dOrb}^{{\rm NC}}$ is the category of {\it normally complex derived orbifolds};

\item $\pman$ is the category of {\it pseudomanifolds}.
\end{itemize}
In the case of Floer complex, the geometric data of Floer theory (the Hamiltonian and the almost complex structure) provides a flow category enriched in $\uds{\bf Top}$. The regularization component can be summarized as constructing a ``lift'' of this flow category to $\uds{\bf SKur}^{\rm NC}$. This step is basically accomplished by repackaging the {\it AMS construction} \cite{AMS, AMS2} (where the construction of the normal complex structure also borrows ideas from \cite{Abouzaid_Blumberg_2024}), which contains several delicate steps. Once this lift is constructed, it automatically descends to a flow category enriched in $\uds{\bf dOrb}^{\rm NC}$, where the blackbox theorem on FOP perturbations can be directly applied. After the FOP perturbation, the flow category becomes enriched in $\pman$, which automatically produces a chain complex over integers.

To obtain chain-level structures such as continuation maps, pair-of-pants products etc., we also consider two ``higher'' analogues of flow categories, i.e., flow multimodules and flow homotopies. Typically, flow categories induce chain complexes, multimodules induce chain maps, and flow homotopies induce chain homotopies. These chain-level structures are also obtained from compatibly implementing the AMS+FOP procedure, namely, producing multimodules and homotopies enriched in categories appearing in the above zig-zag diagram. The differences from the flow category cases only lie in various sub-technical levels.

\subsubsection{Other AMS-type regularization schemes}

There have been several related works which carry out AMS-type constructions in Floer theory, which aim at constructing not only a single global chart construction, but a compatible system for all involved moduli spaces. These works include  \cite{Rezchikov_Arnold},  \cite{Hirschi_Hugtenburg}, \cite{Chanda_Hirschi_contact}, \cite{Mak_Seyfaddini_Smith_1}, and \cite{Rabah_thesis}, as well as our previous work \cite{Bai_Xu_Arnold}.

One distinct feature of our work is the formal  description of the compatibility condition using the notion of ``lifts'' of flow categories to various ``regular stratification categories.'' Having realized that each step of the AMS construction requires a certain compatibility for the choices made for individual moduli spaces (such as the thickening data, stable smoothing, normal complex structure etc.), the formal notions of regular stratification categories and lifts are more convenient and rigorous, although one needs to carefully define each regular stratification category.

The compatibility condition we require here is also most refined since we need to construct genuine perturbations (essentially the same level of refinedness should be established in \cite{Rezchikov_Arnold} and \cite{Rabah_thesis}, both of which use the FOP perturbation method developed by us). The specific requirement is contained in the definition of various regular stratification categories, most notably $\outer \uds{\bf SKur}_{\rm rig}^{\rm NC}$, the category of {\it collared rigidified normally complex smooth Kuranishi spaces}. Here the collaredness and rigidifiedness are what make the inductive perturbation possible, from low-energy moduli spaces to high-energy ones. Some other works such as \cite{Hirschi_Hugtenburg} and \cite{Chanda_Hirschi_contact} where genuine perturbations are not needed, the requirement of compatibility can be relaxed. On the other hand, the work \cite{Mak_Seyfaddini_Smith_1}, which also relies on constructing compatible (multi-valued) perturbations, uses a different strategy to interpolate not-so-compatible choices over the collar regions.

\subsection{Structure of the paper}

This paper is divided into four main parts. In Part 1, we introduce a series of formal notions, including various categories such as those of pseudomanifolds with corners, Kuranishi spaces, derived orbifolds, and their variants; the notions of flow categories, multimodules, and homotopies; the abstract construction of chain-level objects; and the blackbox of FOP perturbation method. In Part 2, we introduce flow categories, multimodules, and homotopies appearing in (non-equivariant) Floer theory, state the main technical theorems about AMS regularization, and prove most of the main results stated in the introduction. In Part 3, we provide the details of the AMS construction asserted in Part 2. In Part 4, we modify all the previous constructions in the equivariant setting and establish results related to the quantum Steenrod operations.

\subsection{Acknowledgements}

We thank Egor Shelukhin and Nicolas Wilkins for related discussions involving our collaboration \cite{BSWX}. We thank Zihong Chen, Helmut Hofer, Rohil Prasad, Sobhan Seyfaddini for their interest in our work and encouragement.

\newpage

\part{FORMAL NOTIONS}

\section*{Outline of Part 1}

We introduce the abstract framework underlying the formulation of Floer theory in various levels. Several concepts related to Floer theory and the virtual construction are described in more abstract ways, including the abstract notion of flow categories/multimodules/homotopies/concatenations enriched in so-called regular stratification categories and their lifts to different categories. Concrete examples of such regular stratification categories will appear in the geometric construction, including Kuranishi spaces, derived orbifolds, pseudomanifolds, etc. We also formulate the abstract mechanism for deriving chain-level structures from flow categories/multimodules/homotopies. Lastly, we provide the FOP perturbation method as a blackbox theorem which can be directly applied. 

\section{Regular Stratification Categories}

\subsection{Regular posets}

The purpose of this subsection is to introduce notions associated with partially ordered sets which will be used to define the index sets of stratifications of various moduli spaces coming from geometric constructions.

\subsubsection{Partially ordered sets}

The term ``partially ordered set'' is abbreviated as {\bf poset}, which is always assumed to be countable. The partial order relations are denoted by $\leq$, $\preceq$, etc. A poset map is a map of the underlying sets which respects the partial orders.

For a poset $A$, let $A^{\max} = A^{\rm top} \subseteq A$ be the subset of maximal elements. A poset $A$ has a canonical {\bf Alexandrov topology}: a subset $U \subset A$ is Alexandrov open if $\alpha \in U$ and $\alpha \leq \beta$ imply that $\beta \in U$. By definition, poset maps are continuous under the Alexandrov topology. A special class of closed subset indexed by $\alpha \in A$
\beqn
\partial^\alpha A = \{ \beta \in A\ |\ \beta \leq \alpha \}
\eeqn
is used frequently in this paper.

\begin{rem}
Our convention is opposite from that of Abouzaid \cite{Abouzaid_axiomatic}, where he regards top strata as minimal elements of the labelling poset. 
\end{rem}

The product of finitely many posets carries a canonically induced partial order from the lexicographic order. Indeed, if $A_1, \ldots, A_k$ are posets, then the relation 
\beqn
(\alpha_1, \ldots, \alpha_k) \leq ( \beta_1, \ldots, \beta_k )\ \text{if and only if}\ \alpha_i \leq \beta_i\ \forall i = 1, \ldots, k,
\eeqn
is a partial order on $A_1 \times \cdots \times A_k$.

\begin{notation}\label{notation21}
Let $S$ be a countable subset. Let $A^{(S)}$ denote the set of all finite subsets of $S$ whose partial order is $\alpha \leq \beta \Longrightarrow \beta \subseteq \alpha$. So its maximal element is the empty set.
\end{notation}

Now we introduce the kind of posets which can model manifolds with corners.

\begin{defn}\label{defn_homogeneous_poset} (Homogeneous posets)

\begin{enumerate}

\item A poset $A$ is called {\bf homogeneous} if for each $\alpha \in A$, the length of a maximal sequence of elements $\alpha = \alpha_0 < \alpha_1 < \cdots < \alpha_k$ with $\alpha_k \in A^{\max}$, is finite and only depends on $\alpha$ (but not the maximal string and $\alpha_k$). This length is called the {\bf depth} or {\bf codimension} of $\alpha$, denoted by $\dep( \alpha )$. In particular, 
\beqn
\alpha \in A^{\max} \Longleftrightarrow \dep(\alpha) = 0.
\eeqn

\item A poset map between homogeneous posets is called {\bf homogeneous of degree $k$} if the map increases the depth function by $k \geq 0$.

\item Given a homogeneous poset $A$, define
\beqn
\mb{Face}^\gamma := \{ \beta < \gamma \ |\ \dep(\beta) = \dep(\gamma) + 1 \}
\eeqn
and 
\beqn
\mb{Face}_\alpha^\gamma:= \{\beta \in \mb{Face}^\gamma \ |\ \alpha \leq \beta < \gamma \}.
\eeqn
\end{enumerate}
\end{defn}

\begin{defn} \label{defn_regular_poset}
A homogeneous poset $A$ is called {\bf regular} if every $\gamma \in A$ satisfies the following conditions.
\begin{enumerate}

\item For any $\alpha \leq \gamma$, $\mb{Face}_\alpha^\gamma$ is always a finite set (but not necessarily $\mb{Face}^\gamma$).

\item The natural map to the set of subsets of $\mb{Face}^\gamma$
\beqn
f_\gamma: \partial^\gamma A \to A^{(\mb{Face}^\gamma)},\ \alpha \mapsto \mb{Face}_\alpha^\gamma
\eeqn
is injective.

\item For any finite subset $S \subset \mb{Face}^\gamma$
\beqn
\bigcap_{\beta \in S} \partial^\beta A = \left\{ \begin{array}{cc} \emptyset, & S \notin {\rm Im} f_\gamma,\\
                  \partial^\alpha A, & S = f_\gamma(\alpha)\end{array}\right.
\eeqn

\end{enumerate}
\end{defn}

We often need to identify the faces that are only adjacent to one top stratum. 

\begin{defn}\label{defn_normal_poset}
Let $A$ be a regular poset.
\begin{enumerate}

\item A {\bf true boundary stratum} of $A$ is an element $\beta \in A$ of depth $1$ which has only one top stratum above. The {\bf true boundary} of $A$ is the Alexandrov closure of the set of true boundary strata.

\item A regular poset is called {\bf normal} if each depth-one element is either adjacent to one top stratum (true boundary) or adjacent to two top strata (fake boundary).
\end{enumerate}
\end{defn}

\begin{rem}
The poset $A^{(S)}$ in Notation \ref{notation21} is regular with depth function given by $\dep(\alpha) = \# \alpha$. On the other hand, the poset of three elements $\alpha < \beta < \gamma$ is homogeneous but not regular.
\end{rem}

We prove several important properties of regular posets.

\begin{lemma}[Properties of regular posets]\label{lemma_regular_poset}
Let $A$ be a regular poset. Then the following is true.
\begin{enumerate}
\item Any Alexandrov closed subset or Alexandrov open subset is also regular.

\item For $\alpha \leq \gamma$, the natural map 
\beq\label{eqn_regular_poset}
f_{\gamma \alpha}:\{ \beta \in A\ |\alpha \leq \beta \leq \gamma \} \to A^{(\mb{Face}_\alpha^\gamma)},\ \beta \mapsto \mb{Face}_\beta^\gamma
\eeq
is a homogeneous bijection.

\item For any $\alpha \leq \gamma$, 
\beqn
\dep(\alpha) = \dep(\gamma) + \# \mb{ Face}_\alpha^\gamma.
\eeqn

\item When $\alpha \leq \beta \leq \gamma$, there is a canonical bijection
\beqn
\mb{Face}_\alpha^\gamma \setminus \mb{Face}_\beta^\gamma \cong \mb{Face}_\alpha^\beta.
\eeqn
Moreover, when $\alpha\leq\beta\leq\gamma\leq \delta$, the composition $\mb{Face}_\alpha^\beta \to \mb{Face}_\alpha^\gamma \to \mb{Face}_\alpha^\delta$ coincides with the inclusion $\mb{Face}_\alpha^\beta \to \mb{Face}_\alpha^\delta$ defined from the above bijection.

\item If $i: A \to A'$ is a homogeneous poset injection to another regular poset $A'$, then for any $\alpha \leq \gamma$ in $A$, the induced map $\mb{Face}_\alpha^\gamma \to \mb{Face}_{i(\alpha)}^{i(\gamma)}$ is a bijection.

\end{enumerate}
\end{lemma}

\begin{proof}
(1) follows from the definition of the Alexandrov topology and the definition of regular posets. For (2), we show that the $f_{\gamma\alpha}$ is bijective. By the definition of regular posets, this is an injection. For surjectivity, assume on the contrary that $S \subset \mb{Face}_\alpha^\gamma$ is not in the range of $f_{\gamma \alpha}$. Then it is not in the range of $f_\gamma$ either. By the definition of regular posets, we must have 
\beqn
\bigcap_{\beta \in S} \partial^\beta A = \emptyset.
\eeqn
However, the above intersection at least contains $\alpha$, a contradiction. Hence $f_{\gamma\alpha}$ is surjective. Moreover, it is clear that $f_{\gamma\alpha}$ is a poset map. One can verify that $f_{\gamma\alpha}$ is homogeneous by an easy inductive argument. Then (3) follows from the homogeneity condition.

For (4), for any $\lambda \in\mb{Face}_\alpha^\gamma \setminus \mb{Face}_\beta^\gamma$, consider 
\beqn
\partial^\lambda A \cap \partial^\beta A = \partial^\lambda A \cap \bigcap_{\tau \in \mb{Face}_\beta^\gamma } \partial^\tau A.
\eeqn
Then the subset $\{ \lambda \}\cup \mb{Face}_\beta^\gamma \subset \mb{Face}_\alpha^\gamma$, under the inverse of $f_{\gamma\alpha}$, gives an element $\kappa$ between $\alpha$ and $\beta$ whose depth is one bigger than that of $\beta$, hence defines a map from $\mb{Face}_\alpha^\gamma  \setminus \mb{Face}_\beta^\gamma $ to $\mb{Face}_\alpha^\beta $. One can verify that this is a bijection.

For (5), the homogeneity condition implies that $i$ induces an injection from $\mb{Face}_\alpha^\gamma $ to $\mb{Face}_{i(\alpha)}^{i(\gamma)}$. Then by (3), we have
\beqn
\# \mb{Face}_\alpha^\gamma  = \dep(\gamma) - \dep(\alpha) = \dep(i(\gamma)) - \dep(i(\alpha)) = \# \mb{Face}_{i(\alpha)}^{i(\gamma)},
\eeqn
which shows that the induced map is a bijection.
\end{proof}

It is straightforward to check that regular posets admit finite products (Cartesian product) and finite coproducts (disjoint union). 

\begin{defn}
The category of regular posets, denoted by $\regpos$, has objects being countable regular posets and morphisms being injective homogeneous poset maps. 
\end{defn}

\begin{prop}
$\regpos$ is a distributive %
monoidal category with respect to the direct product and disjoint union, whose units are posets with a single element.
\end{prop}
\begin{proof}
This is straightforward from checking the definitions.
\end{proof}

\subsection{Regular stratification categories and presentations}

\begin{defn}\label{defn_regular_stratification_category}
A {\bf regular stratification category} is a distributive monoidal category $\uds{\bf C}$ together with a distributive monoidal functor $\uds{\bf C} \to \regpos$, called the {\bf forgetful functor}, satisfying the following properties. In the following, we say an object $K$ of $\uds{\bf C}$ an $A$-stratified object if the forgetful functor sends $K$ to $A$.
\begin{enumerate}

\item There is an initial object which is $\emptyset$-stratified. 

\item For any morphism $f: A \to B$ in $\regpos$ and any $B$-stratified object $K$, there is a canonical ``pullback,'' which is an $A$-stratified object $f^* K$, together with a canonical morphism $\tilde f: f^* K \to K$. The pullback needs to satisfy the functoriality condition. Namely, the composition of pullbacks is the pullback by the composition.

In particular, for any $A$-stratified object $K$ and $\alpha \in A$ there is a {\bf restriction} $\partial^\alpha K$ for each $\alpha \in A$.
\end{enumerate}
\end{defn}

As we will see below, there are a few examples of regular stratification categories such as those of regularly stratified topological spaces (with morphisms being embeddings), regularly stratified topological/smooth manifolds, regularly stratified Kuranishi charts, and Kuranishi presentations. On the other hand, if $\uds{\bf C}$ is any category, then the product $\uds{\bf RegPos}\times \uds{\bf C}$ is a regular stratification category with respect to the obvious forgetful functor.

\subsection{Stratified spaces, manifolds, and pseudomanifolds}

\subsubsection{Regularly stratified spaces}

Although we can define stratified spaces in general, in this paper all such spaces have strata indexed by regular posets. Therefore, we restrict to this situation. We will specify a regular stratification category $\uds{\bf Top}$ of regularly stratified topological spaces.

\begin{defn}\label{defn_stratified}
Let $A$ be a countable regular poset. An \textbf{$A$-stratified space} ($A$-space for short) is a locally compact, Hausdorff and second countable topological space $X$ endowed with a map
\beqn
s: X \to A
\eeqn
which is continuous with respect to the Alexandrov topology on $A$ whose range is finite. 
\end{defn}

In particular, we can write
\beqn
X = \bigsqcup_{\alpha \in A} X_\alpha,\ {\rm where}\ X_\alpha:= s^{-1}(\alpha),
\eeqn
which satisfies the following conditions. 
\begin{enumerate}
    \item Each $X_\alpha$ (called a {\bf stratum}) is locally closed.
    
    \item All but finitely many strata are empty.
    
    \item For each $\alpha\in A$, the subset
    \beqn
    \partial^\alpha X:= \bigsqcup_{\beta \leq \alpha} X_\beta
    \eeqn
    is a closed set (which may contain the closure of $X_\alpha$ properly\footnote{This condition is weaker than the ``axiom of frontier:'' $X_\alpha \cap \ov{X_\beta} \neq \emptyset \Longrightarrow X_\alpha \subset \ov{X_\beta}$, which is part of the usual definition of a stratification.}). Note that this condition follows from the continuity of $s$.
\end{enumerate}

\begin{rem}
An open subset of an $A$-stratified space is still regarded as $A$-stratified, although the subset may have more empty strata.
\end{rem}

We introduce the following notions for stratified spaces.

\begin{defn}(Category of regularly stratified spaces)
\begin{enumerate}
    \item A map from an $A_1$-space $X_1$ to an $A_2$-space $X_2$ is a commutative diagram
    \beqn
    \xymatrix{ X_1 \ar[r]^f \ar[d]_{s_1} & X_2 \ar[d]^{s_2} \\
               A_1 \ar[r]_i & A_2} 
               \eeqn
    where $i: A_1 \to A_2$ is a morphism in $\regpos$ and $f$ is a continuous map. We usually call such a map a {\bf stratified map} to emphasize that it respects the stratifications.

    \item A stratified map $f: X_1 \to X_2$ as above is called an {\bf embedding} if $f$ is a homeomorphism onto its image. $f$ is called an {\bf open embedding} if it is an embedding with an open image; it is called a homeomorphism if $f$ is a homeomorphism of topological spaces and the underlying poset map is an isomorphism.

    \item The category $\uds{\bf Top}$ is the category of topological spaces stratified by regular posets whose morphisms are stratified embeddings.
    
    \item Let $G$ be a topological group. A $G$-action on an $A$-space $X$ is a continuous $G$-action on $X$ which preserves each stratum. In this case, the $\partial^\alpha A$-space $\partial^\alpha X$ has an induced $G$-action.
\end{enumerate}
\end{defn}

\subsubsection{Notations for corners}

In order to define stratified pseudomanifolds, we need to fix the notations for corners. Let $\epsilon>0$ be a positive real number. The model for collars is the space $[0, \epsilon)^k$ which is stratified by subsets of $\{1, \ldots, k\}$. More generally, for any finite set $S$, consider 
\beqn
[0, \epsilon)^S = \{(x_a)_{a\in S}\ |\ x_a \in [0, \epsilon)
\}
\eeqn
which is stratified by $A^{(S)}$. 

\begin{defn}\label{defn_A_manifold}
Let $A$ be a countable regular poset. 
\begin{enumerate}

\item An $A$-space $X$ is said to {\bf have corners} if for each $x \in X_\alpha$ and each maximal $\gamma \in A^{\max}$ with $\alpha \leq \gamma$, there exists a {\bf corner chart}, i.e., a stratified open embedding onto an open neighborhood of $x$ in $X$
\beqn
\xymatrix{  U_{x, \alpha} \times [0, \epsilon)^{\mb{Face}_\alpha^\gamma} \ar[r]^-\varphi \ar[d] & \partial^\gamma X \ar[d] \\
             \{ \alpha \} \times A^{(\mb{Face}_\alpha^\gamma)} \ar[r] &  \partial^\gamma A }
\eeqn
where $U_{x, \alpha}$ is an open neighborhood of $x$ in $X_\alpha$ (which is stratified by $\{\alpha\}$). Here the underlying poset map, whose image is $\{ \beta \in A\ |\ \alpha \leq \beta \leq \gamma\} \subset \partial^\gamma A$, is the inverse of the isomorphism \eqref{eqn_regular_poset}.

\item An {\bf $A$-stratified manifold} is an $A$-stratified space $X$ which has corners such that each $X_\alpha$ (the ``interior'' of $\partial^\alpha X$) is a topological manifold.\footnote{It follows from this condition that each $\partial^\alpha X$ is a topological manifold with boundary whose interior is $X_\alpha$.}  A {\bf smooth structure} on an $A$-stratified manifold $X$ consists of smooth structures on all $X_\alpha$ and a maximal atlas of corner charts which are smoothly compatible.
\end{enumerate}
\end{defn}

Then we can define the notion of smooth functions, tensors, vector bundles, Lie group actions, on smooth $A$-manifolds. We omit the details. Moreover, we obtain a regular stratification category $\uds{\bf Man}$ of topological manifolds and a regular stratification category $\uds{\bf SMan}$ of smooth manifolds. %

\subsubsection{Stratified pseudomanifolds}

\begin{defn}
An $n$-dimensional {\bf pseudomanifold} is a topological space $M$ with a partition $M = {\rm Int} M \sqcup \delta M$ satisfying
\begin{enumerate}
    \item ${\rm Int} M$ is open, dense, and itself an $n$-dimensional topological manifold.

    \item $\delta M$ admits a partition into topological manifolds of dimensions at most $n-2$.
\end{enumerate}
\end{defn}

\begin{defn}
Let $A$ be a countable regular poset. An {\bf $A$-stratified $n$-dimensional pseudomanifold} is an $A$-stratified space $X$ together with partitions $X_\alpha = {\rm Int} X_\alpha \sqcup \delta X_\alpha$ such that
\begin{enumerate}
    \item For each $\alpha \in A$, $X_\alpha = {\rm Int} X_\alpha \sqcup \delta X_\alpha$ is a pseudomanifold of dimension $n - \dep (\alpha)$. 

    \item For $x \in X_\alpha$ and each maximal $\gamma \in A^{\rm max}$ with $\alpha \leq \gamma$, there exists a corner chart
    \beqn
    U_{x, \alpha} \times [0, \epsilon)^{\mb{Face}_\alpha^\gamma } \to \partial^\gamma X
    \eeqn
    which sends ${\rm Int} U_{x, \alpha} \times [0,\epsilon)^{\mb{Face}_\alpha^\gamma}$ resp. $\delta U_{x, \alpha} \times [0,\epsilon)^{\mb{Face}_\alpha^\gamma}$ into ${\rm Int} \partial^\gamma X$ resp. $\delta( \partial^\gamma X)$.
\end{enumerate}
\end{defn}

Notice that a stratified zero-dimensional pseudomanifold is a discrete collection of points and a stratified 1-dimensional pseudomanifold is a disjoint union of 1-manifolds (with boundary).

\begin{defn}
An embedding of stratified pseudomanifolds from $X_1$ to $X_2$ is an embedding of stratified spaces $f: X_1 \to X_2$ such that $f({\rm Int} X_1) \subset {\rm Int}(X_2)$ and $f(\delta X_1) \subset \delta X_2$. The category of stratified pseudomanifolds, denoted by $\pman$, has objects being stratified pseudomanifolds and morphisms being embeddings. 
\end{defn}

One can check that $\pman$ is a regular stratification category (Definition \ref{defn_regular_stratification_category}).

\section{Kuranishi Spaces}

The Kuranishi theory for moduli spaces of pseudoholomorphic curves originated from the efforts of extending the definition of symplectic Gromov--Witten invariants and Hamiltonian Floer homology to full generality. It is based on the fundamental fact that the Cauchy--Riemann operator is elliptic. Fukaya--Ono \cite{Fukaya_Ono} firstly introduced the package of Kuranishi structures, alongside another finite-dimensional reduction approach of Li--Tian \cite{Li_Tian} (see also \cite{Siebert_virtual} \cite{Ruan_virtual}). While the method of Kuranishi structures was further developed by Fukaya--Oh--Ohta--Ono (see the comprehensive treatment in \cite{FOOO_Kuranishi}) together with closely related works by McDuff--Wehrheim \cite{MW_1, MW_2, MW_3}, there are alternate approaches with considerable novelties, including the algebraic version of Pardon \cite{Pardon_virtual} as well as the more recent ``global chart construction'' of Abouzaid--McLean--Smith \cite{AMS} (followed up by \cite{Hirschi_Swaminathan} and \cite{AMS2}). 

The AMS construction was adopted by the authors in \cite{Bai_Xu_2022} in the Gromov--Witten setting and \cite{Bai_Xu_Arnold} in the Hamiltonian Floer setting.

We use $\uds{\bf Kur}$, $\uds{\bf SKur}$, and $\uds{\bf SKur}^{\rm NC}$ to denote the category of Kuranishi charts in topological, smooth, or stably complex categories. When the discussion can be unified in the three cases, we simply use $\uds{\bf Kur}$. See below for the detailed definitions.

\subsection{Kuranishi spaces}

\begin{defn}\label{defn_K_chart}
Let $A$ be a countable regular poset.

\begin{enumerate}
\item An $A$-stratified {\bf (topological) Kuranishi space} (K-space for short) is a quadruple $K = (G, V, E, S)$ where $G$ is a compact Lie group, $V$ is an  $A$-manifold with a continuous $G$-action, $E \to V$ is a $G$-equivariant vector bundle, and $S: V \to E$ is a $G$-equivariant section, such that the $G$-action on $V$ has only finite stabilizers.%

\item A Kuranishi space $K = (G, V, E, S)$ is said to be {\bf smooth} if $V$ is a smooth  $A$-manifold, the $G$-action is smooth, and $E \to V$ is a smooth equivariant vector bundle (we do not impose any smoothness condition on $S$).

\end{enumerate}
\end{defn}

\subsubsection{Morphisms of Kuranishi spaces}

In order to define a regular stratification category (see Definition \ref{defn_regular_stratification_category}), we need to specify morphisms. We first define a general form of maps.

\begin{defn}
A {\bf strict map} of smooth/topological Kuranishi spaces from $K_1 = (G_1, V_1, E_1, S_1)$ to $K_2 = (G_2, V_2, E_2, S_2)$, denoted by $\iota_{21}: K_1 \to K_2$, consists of a Lie group homomorphism $\iota_{21}^G: G_1 \to G_2$ and a commutative diagram
\beqn
\xymatrix{   E_1 \ar[r]^{\iota_{21}^E} \ar[d]  &  E_2 \ar[d] \\
    V_1 \ar@/^1pc/[u]^{S_1} \ar[r]_{\iota_{21}^V} & V_2 \ar@/_1pc/[u]_{S_2} }
\eeqn
where $\iota_{21}^V$ is an equivariant smooth/continuous map and $\iota_{21}^E$ is an equivariant smooth/continuous bundle map.
\end{defn}

Notice that strict maps can be composed.

We list a few special kinds of strict maps.

\begin{defn}[Conjugation] \label{K_chart_conjugation} 
Let $K = (G, V, E, S)$ be a K-chart. The {\bf conjugation} of $K$ by $g \in G$ is the strict map $\iota_g: K \to K$ where $\iota_g^G: G \to G$ is the conjugation $h \mapsto g^{-1} h g$ and $\iota_g^V$ resp. $\iota_g^E$ is the action of the form $x \mapsto g^{-1} x$. Two strict maps $\iota_{21}, \iota_{21}': K_1 \to K_2$ are said to be {\bf conjugate} if there exists $g_2 \in G_2$ such that $\iota_{21}' = \iota_{g_2} \circ \iota_{21}$.
\end{defn}

\begin{defn}[Stabilization]
Let $K = (G, V, E, S)$ be an $A$-stratified K-chart. Let $\pi_N: N \to V$ be a $G$-equivariant continuous resp. smooth vector bundle and let $D \subset N$ be a $G$-invariant open (disk) subbundle.\footnote{When we talk about disk bundles, we always assume that it is a subset of a certain vector bundle, although it may not appear so in notations.} The {\bf stabilization} of $K$ by $D$, denoted by ${\rm Stab}_D (K)$, is the $A$-stratified K-chart 
\beqn
{\rm Stab}_D (K) = (G, D, \pi_N^* E \oplus \pi_N^* N, \pi_N^* S \oplus \tau_N)
\eeqn
where $\tau_N: N \to \pi_N^* N$ is the tautological section. By abuse of language, we also call the morphism
\beqn
\iota_D: K \to {\rm Stab}_D (K)
\eeqn
induced by the zero section $V \to D$ and the obvious bundle embedding $E \to \pi_N^* E \oplus \pi_N^* N$ a {\bf stabilization morphism}.
\end{defn}

\begin{defn}[Shrinking and open embedding]
\begin{enumerate}

\item Let $K = (G, V, E, S)$ be a K-chart. A {\bf shrinking} of $K$ is a K-chart $K' = (G, V', E', S')$ where $V' \subset V$ is a $G$-invariant open neighborhood of $S^{-1}(0)$, $E' = E|_{V'}$, and $S' = S|_{V'}$. 

\item A morphism $\iota_{21}: K_1 \to K_2$ is called an {\bf open embedding} if the underlying poset map is an isomorphism onto an Alexandrov open subset, the map $V_1 \to V_2$ is an isomorphism onto an open subset, and the bundle map is an isomorphism of vector bundles when restricted to the open subset identified with $V_1$. 
\end{enumerate}
\end{defn}

\begin{defn}[Group enlargement]
Let $K = (G, V, E, S)$ be a   K-chart. A {\bf group enlargement} of $K$ by a Lie group embedding $G \hookrightarrow G'$ is a K-chart
\beqn
G'\times_G K = (G', G'\times_G V, G'\times_G E, S')
\eeqn
where $S': G'\times_G V \to G'\times_G E$ is the section naturally induced from $S$. By abuse of language, the K-chart morphism
\beqn
K \to G'\times_G K
\eeqn
induced from the natural maps $V \to G'\times_G V$ and $E \to G'\times_G E$ is also called a group enlargement. 
\end{defn}

Finally, the morphism which will be used most frequently in this paper is a combination of the above more elementary kinds.

\begin{defn}[Morphisms of Kuranishi spaces]\label{chart_embedding}
A {\bf strict morphism} or {\bf strict embedding} of Kuranishi spaces from $K_1$ to $K_2$ is a strict map $\iota_{21} : K_1 \to K_2$ such that there exists a commutative diagram of strict maps
\beq\label{embedding_diagram}
\vcenter{ \xymatrix{  K_3 \ar[r]^{\iota_b}    &   K_4  \ar[d]^{\iota_c}   \\    K_1  \ar[u]^{\iota_a} \ar[r]_{\iota_{21}}    &  K_2      }
}
\eeq
where $\iota_a$ is a stabilization by a disk bundle, $\iota_b$ is a group enlargement, and $\iota_c$ is an open embedding. It implies in particular the group morphism $\iota_{21}^G$, the base map $\iota_{21}^V: V_1 \to V_2$, and the bundle map $\iota_{21}^E: E_1 \to E_2$, are all embeddings.
A conjugacy class of strict morphisms from $K_1$ to $K_2$ is called a {\bf morphism}. 
\end{defn}

\begin{rem}
The above definition basically says that a morphism $\iota_{21}: K_1 \to K_2$ is an embedding if (up to group enlargement) the domain map is an embedding and the restriction of $K_2$ to a neighborhood of the embedding image is isomorphic to a stabilization. One can also define a weaker ``infinitesimal version'' of K-chart embeddings. In the smooth category one can prove that the infinitesimal version is equivalent to the above definition. However, in the topological category it may not be so.
\end{rem}

One can check that compositions of embeddings are still embeddings, which ultimately follows from the fact that every $G'$-equivariant vector bundle over $G' / G$ for a group embedding $G \hookrightarrow G'$ is induced from a $G$-representation. Also notice that stabilization maps, open embeddings, and group enlargements are all special cases of embeddings.

\begin{defn}
The categories of regularly stratified topological and smooth Kuranishi spaces, denoted by $\uds{\bf Kur}$ and $\uds{\bf S Kur}$, are the categories whose objects are stratified Kuranishi charts in the topological and smooth categories and whose morphisms are embeddings. 
\end{defn}

\begin{prop}
$\uds{\bf Kur}$ and $\uds{\bf SKur}$ are regular stratification categories.
\end{prop}

\begin{proof}
Disjoint union and direct product make the category $\uds{\bf Kur}$ distributive monoidal, with monoidal unit being the singleton $K_0 = (\{0\}, {\rm pt}, 0, 0)$. On the other hand, one can obviously take pullbacks of any $A$-stratified Kuranishi space with respect to homogeneous poset maps $f: B \to A$. 
\end{proof}

\subsubsection{Free quotients}\label{subsec:free quotients}

There is another relation among Kuranishi spaces which will appear in the construction. This notion first appeared as {\it group enlargement} in \cite[Section 4]{AMS} and then showed up as {\it free quotient} in \cite[Lemma 3.5]{AMS2}.

\begin{defn}\label{defn_Kuranishi_free_quotient}
Let $K^\sim = (G^\sim, V^\sim, E^\sim, S^\sim)$ be a (topological or smooth) Kuranishi space where $G^\sim = G \times G'$ is the product of two compact Lie groups. Suppose the $G'$-action on $V^\sim$ is free. Then the {\bf free quotient} of $K^\sim$ by $G'$ is the Kuranishi space
\beqn
K:= K/G':= (G, V^\sim /G', E^\sim /G', S^\sim /G').
\eeqn
Notice that there is a natural strict morphism (which is not an embedding) from $K^\sim$ to $K$. We also call the morphism a free quotient by $G'$.
\end{defn}

\subsection{Rigidified embeddings}

In our geometric construction we can obtain objects and morphisms living in a more rigid category of Kuranishi charts. In short, we only allow the stabilization to be a product bundle coming from a $G$-representation.

\begin{defn}\label{defn_rigid_embedding}
Let $\iota_{21}: K_1 \to K_2$ be a strict embedding. 
\begin{enumerate}
    
\item A {\bf rigidification} (which may not exist) of $\iota_{21}$ consists of a $G_1$-representation $W_{21}$ (viewed as a $G_1$-vector bundle $\pi_{21}: \uds W_{21} \to V_1$), a $G_1$-invariant open neighborhood $D_{21} \subset W_{21}$ of the origin (viewed as a disk bundle over $V_1$), and an open embedding $\theta_{21}: G_2\times_{G_1} {\rm Stab}_{D_{21}}(K_1) \to K_2$ such that the diagram \eqref{embedding_diagram} is realized by 
\beqn
\vcenter{ \xymatrix{  {\rm Stab}_{D_{21}} (K_1) \ar[r]^-{\iota_b}  &  G_2 \times_{G_1} {\rm Stab}_{D_{21}}(K_1) \ar[d]^{\theta_{21}} \\
             K_1 \ar[u]^{ \iota_a} \ar[r]_{\iota_{21}}  &  K_2} }.
\eeqn
We denote the rigidification by $(D_{21} \subset W_{21}, \theta_{21})$.

\item Two rigidifications $(D_{21} \subset W_{21}, \theta_{21})$ and $(D_{21}' \subset W_{21}', \theta_{21}')$ are called {\bf equivalent} if there exists an isomorphism $\tau: W_{21} \cong W_{21}'$ of representations, which induces an isomorphism of K-charts
\beqn
\tau: {\rm Stab}_{D_{21}} K_1 \cong {\rm Stab}_{D_{21}'} K_1
\eeqn
such that after shrinking $D_{21}$ and $D_{21}'$ properly, as K-chart morphisms one has
\beqn
\theta_{21} = \theta_{21}' \circ \tau.
\eeqn

\item A {\bf rigidified embedding} from $K_1$ to $K_2$ is an embedding together with an equivalence class of rigidifications.
\end{enumerate}
\end{defn}

Note that not all embeddings can be rigidified. On the other hand, rigidified embeddings are well-behaved. We prove that the composition of rigidified embeddings is naturally rigidified.

\begin{lemma}
Let $\iota_{21}: K_1 \to K_2$ be an embedding with a rigidification $(D_{21} \subset W_{21}, \theta_{21} )$ and $\iota_{32}: K_2 \to K_3$ be an embedding with a rigidification $(D_{32}\subset W_{32}, \theta_{32} )$. Viewing $W_{32}$ as a $G_1$-representation via the group embedding $G_1 \to G_2$, define 
\beqn
W_{31}:= W_{21} \oplus W_{32},\ D_{31}:= D_{21} \oplus D_{32}.
\eeqn
Define $\theta_{31}: {\rm Stab}_{D_{31}}(K_1) \to K_3$ as the composition
\beqn
\xymatrix{ {\rm Stab}_{D_{31}} (K_1) \ar[rr]^-{\cong} & & {\rm Stab}_{D_{32}} ({\rm Stab}_{D_{21}} (K_1)) \ar[rr]^-{{\rm Stab}_{D_{32}} \theta_{21} } & & {\rm Stab}_{D_{32}} (K_2) \ar[rr]^-{\theta_{32}} & & K_3} .
\eeqn
Then $(D_{31} \subset W_{31}, \theta_{31})$ is a rigidification of $\iota_{32} \circ \iota_{21}$, whose equivalence class only depends on the equivalence classes of the rigidifications of $\iota_{32}$ and $\iota_{21}$.
\end{lemma}

\begin{proof}
This lemma holds essentially because a trivial bundle over the total space of a trivial bundle is canonically trivial. 
\end{proof}

\begin{rem}
If in Definition \ref{defn_rigid_embedding} we do not require the bundle $W_{21}$ to be trivial, then there is no well-defined composition for rigidified embeddings.
\end{rem}

\begin{rem}
One could criticize that rigidified embeddings are not anything natural. On the other hand, our geometric construction (the AMS construction) of Kuranishi flow categories allows us to obtain a topological Kuranishi flow category with rigidified embeddings as morphisms. This is very helpful in the smoothing procedure: once a stratum is smoothed, the total space of a trivial disk bundle is then canonically smoothed, giving a smooth tubular neighborhood in a higher stratum.
\end{rem}

\subsection{Stable complex structures and normal complex structures}
Next, we discuss stable and normal complex structures on Kuranishi spaces. Geometrically, such structures arise from index-theoretic considerations of (virtual) tangent bundles of moduli spaces of pseudoholomorphic curves. Such structures are also instrumental for incorporating the FOP perturbation scheme to define integral invariants.

\subsubsection{Stable morphisms of vector bundles}

We start with a few basic definitions. 

\begin{defn}\label{defn_stable_isomorphism}
Let $G$ be a Lie group, $V$ be a $G$-space, and $E \to V$ be an equivariant vector bundle. 

\begin{enumerate}

\item Let $F \to V$ be another $G$-equivariant vector bundle. A $G$-equivariant {\bf stable (iso)morphism}  from $E$ to $F$ is an equivalence class of triples $(R^-, R^+, \tau)$ where $R^\pm$ are a pair of real vector spaces (viewed as trivial representations of $G$) and 
\beqn
\tau: \uds R^- \oplus E \to \uds R^+ \oplus F
\eeqn
is a $G$-equivariant (iso)morphism of vector bundles. The equivalence relation is induced by $(R^-, R^+, \tau) \sim (R^- \oplus W, R^+ \oplus W, \tau \oplus {\rm Id}_{\uds W})$ where $W$ is another real vector space. Denote a stable (iso)morphism by 
\beqn
\tau: E \overset{s}{\to} F.
\eeqn

\item If $E, F$ are complex vector bundles, then a stable (iso)morphism $\tau: E \overset{s}{\to} F$ is called {\bf complex} if it is represented by a triple $(R^-, R^+, \tau)$ where $R^\pm$ are complex vector spaces and $\tau$ is a complex vector bundle map. 

\end{enumerate}
\end{defn}

One can verify that stable morphisms can be composed and stable isomorphisms admit inverses as stable isomorphisms. 

\subsubsection{Stable complex structures}

\begin{defn}\label{defn_stable_complex_structure}
Let $G$ be a Lie group, $V$ be a $G$-space, and $E \to V$ be a $G$-equivariant vector bundle.
\begin{enumerate}

\item A {\bf $G$-invariant stable complex structure} on $E$ consists of an equivalence class of pairs $(\tau, F)$ where $F \to V$ is a  $G$-equivariant complex vector bundle and $\tau: E \overset{s}{\to} F$ is a $G$-equivariant stable isomorphism. The equivalence relation is generated by the following relation: $(\tau_1, F_1)$ and $(\tau_2, F_2)$ are equivalent if there is a commutative diagram
    \beqn
    \xymatrix{ E \ar[r]^{\tau_1} \ar[d]_{{\rm Id}_E}  &  F_1 \ar[d] \\
     E \ar[r]_{\tau_2} & F_2}
    \eeqn
    of stable isomorphisms where the right vertical arrow is a complex stable isomorphism.

\item Suppose $E_1$ and $E_2$ are $G$-equivariant stably complex vector bundles over a $G$-space $V$. A $G$-equivariant stable (iso)morphism $f: E_1 \overset{s}{\to} E_2$ is said to be {\bf stably complex} if there exist representatives $(\tau_1, F_1)$ and $(\tau_2, F_2)$ of their stable complex structures and a commutative diagram of stable (iso)morphisms
\beqn
\xymatrix{  E_1 \ar[r]^{\tau_1} \ar[d]_{f}  & F_1 \ar[d] \\
 E_2 \ar[r]_{\tau_2} & F_2}
\eeqn
where the right vertical arrow is a $G$-equivariant complex stable (iso)morphism.

\item Let $V$ be a smooth $G$-manifold (with boundary and corners). A $G$-equivariant stable complex structure on $V$ is a $G$-equivariant stable complex structure on $TV$. Notice that for each stratum $\partial^\alpha V \subset V$, the stable complex structure on $TV$ induces a stable complex structure on $T(\partial^\alpha V)$ because the normal direction is trivial. 
\end{enumerate}
\end{defn}

Now we consider smooth Kuranishi space $K = (G, V, E, S)$ such that the $G$-action on $V$ has only finite stabilizers. In this setting, $TV$ contains a trivial subbundle isomorphic to $\uds{\mf g}$. On the other hand, in our construction, the obstruction bundle $E$ can always be made complex. This simplifies the definition of stable complex structures. 

\begin{defn}
A {\bf stably complex Kuranishi space} is a smooth Kuranishi space $K = (G, V, E, S)$ equipped with a $G$-invariant complex structure on $E$ and a $G$-invariant stable complex structure on $TV/\uds{\mf g}$. 

An {\bf open embedding} of stably complex Kuranishi spaces from $K_1$ to $K_2$ is a smooth open embedding $\iota_{21}: K_1 \to K_2$ such that 1) the vector bundle isomorphism 
\beqn
d\iota_{21}^V: TV_1/ \uds{\mf g}_1 \to TV_2/ \uds{\mf g}_2
\eeqn
is stably complex and 2) the vector bundle isomorphism $\iota_{21}^E: E_1 \to E_2$ is complex linear.
\end{defn}

We consider stabilizations of a stably complex Kuranishi space $K = (G, V, E, S)$ by a disk bundle $\pi_D: D \to V$ contained in a $G$-equivariant complex vector bundle $N \to V$, which is 
\beqn
{\rm Stab}_D(K) = (G, D, \pi_D^* E \oplus \pi_D^* N, \pi_D^* S \oplus \tau_N)
\eeqn
Then $\pi_D^* E \oplus \pi_D^* N$ carries a canonical complex structure. On the other hand, the bundle $TD/\uds{\mf g}$ is non-canonically isomorphic to $\pi_D^* TV/\uds{\mf g} \oplus \pi_D^* N$. An isomorphism can be induced from a $G$-invariant connection on the fiber bundle $D \to V$.

\begin{defn}
The {\bf stabilization} of a stably complex Kuranishi space $K = (G, V, E, S)$ by a disk bundle $D$ contained in a $G$-equivariant complex vector bundle $N \to V$ (equipped with a $G$-invariant complex linear connection) is the smooth stabilization ${\rm Stab}_D(K)$ with the induced stably complex structure on $TD/\uds{\mf g}$ and the complex structure on $\pi_D^* E \oplus \pi_D^* N$.
\end{defn}

\begin{defn}
A {\bf strict embedding} of stably complex Kuranishi spaces from $K_1$ to $K_2$ consists of a smooth strict embedding $\iota_{21}: K_1 \to K_2$ such that in the diagram \eqref{embedding_diagram}, the stabilization $\iota_a: K_1 \to K_3$ is provided by a disk bundle in a complex vector bundle and the open embedding $\iota_c: K_4 \to K_2$ is an open embedding of stably complex Kuranishi spaces.

An {\bf embedding} of stably complex Kuranishi spaces is a conjugacy class of strict embeddings.
\end{defn}

\begin{lemma}
(Strict) embeddings of stably complex Kuranishi spaces can be composed.
\end{lemma}

\begin{proof}
Essentially one needs to consider disk bundles over disk bundles to ensure the compatibility of stable complex structures. Given Ehresmann connections on $D_{21}\to V_1$ and $D_{32} \to V_2$, the horizontal distribution of $D_{21}$ can then be lifted to $D_{32}$, giving a Ehresmann connection on the iterated bundle.
\end{proof}

\begin{rem}
We could define the category 
\beqn
\uds{\bf SKur}^{\mb C}
\eeqn
whose objects are regularly stratified stably complex (in particular, smooth) Kuranishi spaces  and whose morphisms are embeddings. It is easy to see that this is a regular stratification category, and there is a natural forgetful functor
\beqn
\uds{\bf SKur}^{\mb C} \to \uds{\bf SKur}.
\eeqn
One can also consider the rigidified version. However, it is rather complicated to obtain stably complex Kuranishi spaces in our construction (although we could). Below we consider a more flexible one using normal complex structures, which suffices for our purpose.
\end{rem}

\subsubsection{Normally complex Kuranishi spaces}

Moduli spaces of pseudoholomorphic curves admit stable complex structure because of the nature of the Cauchy--Riemann operator. In Section \ref{section_NC_structure} we will follow the strategy of Abouzaid--Blumberg \cite{Abouzaid_Blumberg, Abouzaid_Blumberg_2024} to do the explicit construction. However, it is a bit cumbersome to construct stable complex structures compatibly. Instead, we will only construct the normal complex structure which is sufficient for our applications and where the compatibility is easier to achieve. 

Let $V$ be a $G$-space. Notice that a stable morphism of $G$-equivariant vector bundles $E, F \to V$ (Definition \ref{defn_stable_isomorphism}) induces the following well-defined linear maps. For each $x \in V$, let $G_x \subset G$ be the stabilizer. Then the fibers $E_x$ and $F_x$ split as 
\begin{align*}
&\ E_x = \mathring E_x \oplus \check E_x,\ &\ F_x = \mathring F_x \oplus \check F_x
\end{align*}
where $\mathring E_x$, $\mathring F_x$ are $G_x$-trivial summands and $\check E_x$, $\check F_x$ are $G_x$-nontrivial pieces. Then an equivalence class of stable isomorphisms induces well-defined $G_x$-equivariant maps
\beqn
\check \tau_x: \check E_x \to \check F_x.
\eeqn
(In particular, if $G$ is trivial, the maps $\check \tau_x$ are all trivial.)

\begin{defn}\label{defn_equivariant_NC_structure}
Let $G$ be a Lie group, $V$ be a $G$-space, and $E \to V$ be a $G$-equivariant vector bundle.
\begin{enumerate}
    \item Two stable complex structures on $E$, represented by $(\tau_1, F_1)$ and $(\tau_2, F_2)$, are {\bf normally equivalent}, if there exists a diagram of stable isomorphisms
    \beqn
    \xymatrix{ E \ar[r]^{\tau_1} \ar[d]_{{\rm Id}_E} &  F_1 \ar[d]\\
              E \ar[r]_{\tau_2}  & F_2  }
    \eeqn
    where the right vertical arrow is complex, such that for all $x \in V$, the following induced diagram of linear maps commutes.
    \beqn
    \xymatrix{ \check E_x \ar[r]^{\check \tau_{1, x}} \ar[d] & \check F_{1, x} \ar[d] \\
               \check E_x \ar[r]_{\check \tau_{2, x}} & \check F_{2, x} }
    \eeqn 

    \item A $G$-invariant {\bf normal complex structure} on $E$ is a normal equivalence class of stable complex structures. Notice that a normal complex structure on $E$ equips each $\check E_x$ with a $G_x$-invariant complex structure.

    \item Suppose $E_1, E_2 \to V$ are equipped with $G$-equivariant normal complex structures. A $G$-equivariant stable morphism $\tau: E_1 \overset{s}{\to} E_2$ is called {\bf normally complex} if for each $x \in V$, the induced linear map $\check E_{1, x} \to \check E_{2, x}$ is complex linear.
\end{enumerate}
\end{defn}

Now we define the notion of normal complex Kuranishi spaces. 

\begin{defn}\label{defn_NC_Kuranishi}
A {\bf normally complex Kuranishi space} (NC Kuranishi space for short) is a smooth Kuranishi space $K = (G, V, E, S)$ equipped with a $G$-invariant complex structure on $E$ and a $G$-invariant normal complex structure on $TV/{\mf g}$. An {\bf open embedding} of NC Kuranishi spaces from $K_1 = (G_1, V_1, E_1, S_1)$ to $K_2 = (G_2, V_2, E_2, S_2)$ is a smooth open embedding $\iota_{21}: K_1 \to K_2$ such that 1) the vector bundle isomorphism 
\beqn
d\iota_{21}: TV_1/{\mf g}_1 \to TV_2/{\mf g}_2
\eeqn
is normally complex and 2) the vector bundle isomorphism $\iota_{21}^E: E_1 \to E_2$ is complex-linear.
\end{defn}

One can again discuss stabilizations (by complex vector bundles), strict embeddings, conjugations, and embeddings of normally complex Kuranishi spaces. Then one obtains a regular stratification category called the category of normally complex Kuranishi spaces, denoted by 
\beqn
\uds{\bf SKur}^{\rm NC}
\eeqn
and the corresponding rigidified category $\uds{\bf SKur}^{\rm NC}_{\rm rig}$.

\section{Derived Orbifolds}

In this section we describe another regular stratification category called the category of derived orbifolds. We will consider as general as the smooth category, but not the topological category. 

\subsection{Stratified orbifolds}

All orbifolds, unless declared, are assumed to be smooth and effective. They are spaces which are locally modelled by $U / \Gamma$ where $\Gamma$ is a finite group and $U$ is an effective smooth $\Gamma$-manifold. Orbifolds with boundary and corners are locally modelled on $U/ \Gamma$ where $U$ is a $\Gamma$-manifold with boundary and corners. If $A$ is a regular poset (Definition \ref{defn_regular_poset}), then one can consider the notion of $A$-orbifold. Without giving the detailed definition, we remind the reader that a smooth effective $A$-orbifold contains charts of the form 
\beqn
U_\alpha \times [0, \epsilon)^{(\mb{Face}_\alpha^\gamma  )},\ {\rm where\ } \alpha \leq \gamma \in A^{\max},
\eeqn
where $U_\alpha$ is a $\Gamma$-manifold without boundary. Here we always assume that the isotropy group acts trivially on the normal factors. For example, if $X$ is a smooth $A$-manifold acted on by a compact Lie group $G$ with finite stabilizers, then the quotient $X/ G$ is naturally a smooth effective $A$-orbifold.

\begin{defn}[Derived orbifolds]
Fix a regular poset $A$. An $A$-stratified {\bf derived orbifold} is a triple ${\mc D} = ({\mc U}, {\mc E}, {\mc S})$ where ${\mc U}$ is an effective $A$-orbifold, ${\mc E} \to {\mc U}$ is a smooth orbifold vector bundle, and ${\mc S}: {\mc U} \to {\mc E}$ is a continuous section. A derived orbifold ${\mc D}$ is called {\bf compact} if the coarse space of ${\mc S}^{-1}(0)$ is a compact topological space.
\end{defn}

Similar to the case of Kuranishi spaces, we define embeddings of derived orbifolds as combinations of several more basic morphisms.

\begin{defn} Let ${\mc D}$ be an $A$-stratified derived orbifold. 

\begin{enumerate}

\item A {\bf morphism of derived orbifolds} from an $A_1$-stratified derived orbifold ${\mc D}_1 = ({\mc U}_1, {\mc E}_1, {\mc S}_1)$ to an $A_2$-stratified derived orbifold ${\mc D}_2 = ({\mc U}_2, {\mc E}_2, {\mc S}_2)$ consists of a morphism $i: A_1 \to A_2$ of $\regpos$, a smooth stratified orbifold map $\phi_{21}: {\mc U}_1 \to {\mc U}_2$, and a smooth orbibundle map $\wh\phi_{21}: {\mc E}_1 \to {\mc E}_2$ covering $\phi_{21}$ such that the following diagram commutes.
\beqn
\xymatrix{ {\mc E}_1 \ar[r]^{\wh\phi_{21}} \ar[d]  & {\mc E}_2 \ar[d] \\
           {\mc U}_1 \ar[r]_{\phi_{21}}  \ar@/^1pc/[u]^{{\mc S}_1} & {\mc U}_2 \ar@/_1pc/[u]_{{\mc S}_2} }
\eeqn
Such a morphism induces an orbispace map ${\mc S}_1^{-1}(0) \to {\mc S}_2^{-1}(0)$. 

\item An {\bf open embedding of derived orbifolds} from ${\mc D}_1$ to ${\mc D}_2$ is a morphism such that the poset map $i: A_1 \to A_2$ is an isomorphism onto an Alexandrov open subset, the orbifold map $\phi_{21}$ is an isomorphism onto an open subset, and the bundle map $\wh\phi_{21}$ is an isomorphism over the image of $\phi_{21}$.

\item Given a derived orbifold ${\mc D} = ({\mc U}, {\mc E}, {\mc S})$, for any open subset ${\mc U}' \subset {\mc U}$ with ${\mc S}^{-1}(0) \subset {\mc U}'$, the {\bf shrinking} of ${\mc D}$ to ${\mc U}'$, denoted by ${\mc D}|_{{\mc U}'}$, is the derived orbifold $({\mc U}', {\mc E}|_{{\mc U}'}, {\mc S}|_{{\mc U}'})$. It admits a natural open embedding into ${\mc D}$.

\item A disk bundle over an orbifold is an open neighborhood of the zero section of an orbifold vector bundle. A {\bf stabilization} of a derived orbifold ${\mc D} = ({\mc U}, {\mc E}, {\mc S})$ by a disk bundle ${\mc N} \to {\mc U}$ contained in a vector bundle $\pi_{\mc F}: {\mc F} \to {\mc U}$ is the derived orbifold 
\beqn
{\rm Stab}_{\mc N} {\mc D} = ({\mc N}, \pi_{\mc N}^* {\mc E} \oplus \pi_{\mc N}^* {\mc F}, \pi_{\mc N}^* {\mc S} \oplus \tau_{\mc F})
\eeqn
where $\tau_{\mc F}$ is the tautological section. We call the natural morphism
\beqn
{\mc D} \to {\rm Stab}_{\mc N}( {\mc D} )
\eeqn
the {\bf stabilization morphism} by ${\mc N}$.

\item An {\bf embedding of derived orbifolds} from ${\mc D}_1 = ({\mc U}_1, {\mc E}_1, {\mc S}_1)$ to ${\mc D}_2 = ({\mc U}_2, {\mc E}_2, {\mc S}_2)$ is a morphism ${\bm \phi}_{21}: {\mc D}_1 \to {\mc D}_2$ such that  there exists a disk bundle ${\mc N}_{21} \to {\mc U}_1$ and an open embedding $\theta_{21}: {\rm Stab}_{{\mc N}_{21}}( {\mc D}_1) \to {\mc D}_2$ extending $\phi_{21}$. 
\end{enumerate}
\end{defn}

\begin{defn}
The category of stratified derived orbifolds, denoted by $\uds{\bf dOrb}$, has objects being stratified derived orbifolds and morphisms being embeddings of derived orbifolds. 
\end{defn}

The following is readily checked. 

\begin{prop}
The category $\uds{\bf dOrb}$ is a regular stratification category (Definition \ref{defn_regular_stratification_category}). \qed
\end{prop}

One can see that there are functors
\beqn
\uds{\bf S Kur}  \to \uds{\bf dOrb} \to \uds{\bf Top}.
\eeqn
The first arrow sends $K = (G, V, E, S)$ to the quotient ${\mc D} = (V/G, E/G, S/G)$. Notice that a conjugacy class of Kuranishi space embeddings induces a well-defined embedding of derived orbifolds. The second arrow is induced by passing to the zero locus ${\mc S}^{-1}(0)$ of a derived orbifold ${\mc D} = ({\mc U}, {\mc E}, {\mc S})$.

\subsubsection{Rigidified embeddings}

We would like to introduce an analogue of rigidified embeddings for orbifolds. We need to relax the triviality requirement. A flat orbifold vector bundle is an orbifold vector bundle together with a flat connection. A key property of flat vector bundles is, if ${\mc E} \to {\mc U}$ is flat and ${\mc F} \to {\mc E}$ is flat, then ${\mc F}$ is canonically a flat bundle over ${\mc U}$. %

For example, if $G$ is a compact Lie group acting on a smooth manifold $V$ with only finite stabilizers and $W$ is a representation of $G$, then the vector bundle $(V\times W)/G$ is flat. %
Notice that the pullback of rigid vector bundles is still rigid.

\begin{defn}\label{defn_rigid_embedding_D}
Let $\phi_{21}: {\mc D}_1 \to {\mc D}_2   $ be an embedding of derived orbifolds. 
\begin{enumerate}

\item A {\bf rigidification} of $\phi_{21}$ consists of a disk bundle ${\mc N}_{21}$ contained in a flat vector bundle ${\mc F}_{21} \to {\mc U}_1$ and an open embedding
\beqn
\psi_{21}: {\rm Stab}_{{\mc N}_{21}}({\mc D}_1) \to {\mc D}_2
\eeqn
which extends $\phi_{21}$. 

\item Two rigidifications $({\mc N}_{21} \subset {\mc F}_{21}, \psi_{21})$, $({\mc N}_{21}'\subset {\mc F}_{21}', \psi_{21}')$ are equivalent if there is a flat bundle isomorphism $\tau: {\mc F}_{21} \cong {\mc F}_{21}'$ which identifies ${\mc N}_{21}$ with ${\mc N}_{21}'$ after an appropriate shrinking such that 
\beqn
\psi_{21}' = \psi_{21} \circ \tau.
\eeqn

\item A {\bf rigidified embedding} of derived orbifolds is an embedding together with an equivalence class of rigidifications. 
\end{enumerate}
\end{defn}

One can see that rigid derived orbifold embeddings can be composed. Therefore we can define a category of derived orbifolds with rigid embeddings as morphisms. This category is denoted by 
\beqn
\uds{\bf dOrb}_{\rm rig}.
\eeqn
Notice that there is a quotient functor
\beqn
\uds{\bf SKur}_{\rm rig} \to \uds{\bf dOrb}_{\rm rig}.
\eeqn

\subsection{Normal complex structures}

Normal complex (NC for short) structures are essential for constructing FOP transverse perturbations. We recall the definitions here. 

\subsubsection{Basic notions}

Suppose $G$ acts on a set $A$. A subgroup $H \subset G$ is called an $A$-essential subgroup, denoted by 
\beqn
H \subset_A G
\eeqn
if there exists $a \in A$ whose stabilizer is $H$. 

Given a real representation $V$ of a finite group $\Gamma$, the {\bf basic decomposition} of $V$ is the splitting
\beqn
V = V_\Gamma \oplus \check V_\Gamma
\eeqn
where $V_\Gamma$ is the direct sum of trivial subrepresentations and $\check V_\Gamma$ is the direct sum of nontrivial irreducible subrepresentations. An {\bf NC structure} on $V$ is a $\Gamma$-invariant complex structure on $\check V_\Gamma$. 

\subsubsection{NC structures}

\begin{defn}
Let $U$ be a $\Gamma$-manifold and $E \to U$ be a $\Gamma$-equivariant vector bundle. An {\bf NC structure} on $E$ consists of, for each $U$-essential subgroup $H \subset_U \Gamma$, a $H$-invariant complex structure $I^{\check E_H}$ on the subbundle
\beqn
\check E_H \subset E|_{U_H}
\eeqn
satisfying the following conditions.
\begin{enumerate}
    \item For each $g\in \Gamma$ which conjugates $H$ to $H':= g H g^{-1}$, the bundle isomorphism $\check E_H \to \check E_{g H g^{-1}}$ sends $I^{\check E_H}$ to $I^{\check E_{H'}}$.

    \item For each pair of $U$-essential subgroups $K \subset H \subset \Gamma$, one has the $K$-equivariant decomposition  
    \beqn
    \check E_H \cong \check E_K \oplus (\check E_H \cap E_K|_{U_H}).
    \eeqn
    We require that the restriction of $I^{\check E_H}$ to $\check E_K$ coincides with $I^{\check E_K}$.
\end{enumerate}

\end{defn}

\begin{defn}
Let ${\mc U}$ be an effective orbifold and ${\mc E} \to {\mc U}$ be an orbifold vector bundle. An {\bf NC structure} on ${\mc E}$ consists of, for each bundle chart $\hat C = (\Gamma, U, E)$ of ${\mc E}$, an NC structure ${\bm I}^E$ on the vector bundle $E$, such that chart embeddings respect those complex structures. An {\bf NC structure} on ${\mc U}$ is by definition an NC structure on $T{\mc U}$.

A {\bf normally complex derived orbifold} is a derived orbifold ${\mc D} = ({\mc U}, {\mc E}, {\mc S})$ together with NC structures on ${\mc U}$ and ${\mc E}$.
\end{defn}

\begin{defn}[Embeddings of NC derived orbifolds] \hfill
\begin{enumerate}

\item An {\bf open embedding of normally complex derived orbifolds} from ${\mc D}_1$ to ${\mc D}_2$ is an open embedding of derived orbifolds which respects normal complex structures. 

\item A stabilization of an NC derived orbifold is a stabilization by a disk bundle which is contained in a complex orbifold vector bundle.

\item An embedding of NC derived orbifolds $\iota_{21}: {\mc D}_1 \to {\mc D}_2$ is an embedding such that there exists a stabilization of ${\mc D}_1$ by a disk bundle ${\mc N}_{21} \to {\mc U}_1$ and an open embedding from ${\rm Stab}_{{\mc N}_{21}} ({\mc D}_1) \to {\mc D}_2$ which extends $\iota_{21}$. 

\item A {\bf rigidified embedding} of NC derived orbifolds from ${\mc D}_1$ to ${\mc D}_2$ is a rigidified embedding whose rigidification can be represented by a disk bundle in a rigid complex vector bundle. Two rigidifications are equivalent if they (after appropriate shrinking to smaller disk bundles) differ by an isomorphism of rigid complex vector bundles. 
\end{enumerate}
\end{defn}

Therefore we obtained the category $\uds{\bf dOrb}^{\rm NC}$ of normally complex derived orbifolds  with morphisms being embeddings and the category $\uds{\bf dOrb}^{{\rm NC}}_{\rm rig}$ with the same objects with morphisms being rigidified embeddings. One can see that they are both regular stratification categories.

\begin{lemma}
There are natural  functors of regular stratification categories
\begin{align*}
&\ \uds{\bf SKur}^{\rm NC} \to \uds{\bf dOrb}^{{\rm NC}},\ &\ \uds{\bf SKur}_{\rm rig }^{\rm NC} \to \uds{\bf dOrb}^{{\rm NC}}_{\rm rig}
\end{align*}
\end{lemma}

\begin{proof}
Let $K = (G, V, E, S)$ be a stably complex Kuranishi space. Then by definition, $E$ has a $G$-invariant complex structure and the vector bundle 
\beqn
TV/ \uds{\mf g}
\eeqn
has a $G$-invariant stable complex structure. This means there exists a real vector space $R$ and a complex structure on 
\beqn
TV/ \uds{\mf g} \oplus \uds R.
\eeqn
Consider the quotient chart
\beqn
K/G:= ({\mc U}, {\mc E}, {\mc S})
\eeqn
where ${\mc U} = V/ G$. Given any $x \in V$ with stabilizer $G_x \subset G$, one has the $G_x$-invariant decomposition 
\beqn
T_x V \cong T_x V_{G_x} \oplus N_x V_{G_x}.
\eeqn
The stable complex structure on $TV/\uds{\mf g}$ induces a complex structure on $N_x V_{G_x}$. The collection of these complex structures forms a normal complex structure on the quotient orbifold. One can easily check that the morphisms descend to morphisms in $\uds{\bf dOrb}^{{\rm NC}}$ and $\uds{\bf dOrb}^{{\rm NC}}_{\rm rig}$.
\end{proof}

\section{Collars}

Various inductive constructions require moduli spaces, which appear as, e.g., morphism spaces of flow categories, to have compatible collars. While one can possibly discuss collars in an abstract manner, we refrain from doing so. Instead, in this section we only consider the categories with geometric meanings. More precisely, in this section, the category $\uds{\bf C}$ belongs to the list $\uds{\bf Top}$, $\uds{\bf Man}$, $\uds{\bf SMan}$, $\pman$, $\uds{\bf Kur}$, $\uds{\bf dOrb}$ etc. Following an idea of Fukaya, we explain how to equip the geometric objects of interest with collars ``on the outside" which facilitates inductive constructions of various sorts.

\subsection{The subcategory of collared objects}

For the discussions in this section, let $\epsilon > 0$ be a positive real number. The model for collars is the space $[0, \epsilon)^k$ which is stratified by subsets of $\{1, \ldots, k\}$. More generally, for any finite set $S$, we can form the collar $[0, \epsilon)^S$ stratified by the set $A^{(S)}$ (Notation \ref{notation21}). These standard collars are objects of all the relevant regular stratification categories $\uds{\bf C}$ listed above.

\begin{defn}\label{defn_collared objects}
Let $\uds{\bf C}$ be either $\uds{\bf Top}$, $\uds{\bf Man}$, $\pman$, $\uds{\bf SMan}$, $\uds{\bf Kur}$, $\uds{\bf SKur}$, or $\uds{\bf dOrb}$. 
For any $\epsilon>0$, an {\bf $\epsilon$-collaring} of an $A$-stratified object $M \in {\rm Ob}\uds{\bf C}$ is a collection of open embeddings
\beqn
\partial^\alpha M \times [0, \epsilon)^{\mb{Face}_\alpha^\gamma} \to \partial^\gamma M,\ \alpha \leq \gamma
\eeqn
(called a collar embedding from $\alpha$ to $\gamma$) 
such that 1) its restriction to $\partial^\alpha M$ is the inclusion of $\partial^\alpha M$ into $M$ and 2) whenever $\alpha \leq \beta \leq \gamma$, the following diagram commutes.
\beqn
\xymatrix{   \partial^\alpha M \times [0, \epsilon)^{\mb{Face}_\alpha^\gamma}  \ar[r]   \ar[d]   &  \partial^\gamma M \\
            \partial^\alpha M \times [0, \epsilon)^{\mb{Face}_\alpha^\beta} \times [0, \epsilon)^{\mb{Face}_\beta^\gamma}    \ar[r]   & \partial^\beta M \times [0, \epsilon)^{ \mb{Face}_\beta^\gamma} \ar[u]  }
\eeqn

An $\epsilon$-collaring and an $\epsilon'$-collaring are {\bf germ-equivalent} if for a certain $\epsilon'' \leq \epsilon, \epsilon'$, for all $\alpha \leq \gamma$, the restrictions of the collar embeddings to $\partial^\alpha M \times [0, \epsilon'')^{\mb{Face}_\alpha^\gamma }$ are identical. A {\bf collared object} in $\uds{\bf C}$ is an object together with a germ-equivalence class of $\epsilon$-collarings. 
\end{defn}

\begin{defn}\label{defn:collared-morphism}
A morphism between two collared objects $K$ and $K'$ with underlying poset map $i: A \to A'$ is called {\bf collared} if there exists a sufficiently small $\epsilon>0$ such that for each $\alpha \leq \gamma$ in $A$, the morphism over the strata $\gamma \mapsto i(\gamma)$ in the collared neighborhoods of the strata $\alpha \mapsto i(\alpha)$ is given by
\beqn
\xymatrix{
\partial^{\alpha} K \times {\rm Collar}_\epsilon^{\mb{Face}_\alpha^{\gamma} } \ar[rr] \ar[d] &  &  \partial^{i(\alpha)} K' \times {\rm Collar}_\epsilon^{\mb{Face}_{i(\alpha)}^{i(\gamma)}}  \ar[d]\\
\partial^{\gamma} K \ar[rr]  & & \partial^{i(\gamma)} K'
}
\eeqn
where the map along the collar direction is the identity map after identifying $\mb{Face}_\alpha^{\gamma}\cong \mb{Face}_{i(\alpha)}^{i(\gamma)}$.
\end{defn}

In this way we obtain a subcategory 
\beqn
\outer \uds{\bf C}
\eeqn
of collared objects, whose morphisms are collared morphisms. Notice that the empty set and singletons are always collared. Moreover, the collar structures respect products and disjoint unions. Therefore, the category $\outer \uds{\bf C}$ is a regular stratification category (Definition \ref{defn_regular_stratification_category}).

\subsection{Outercollaring construction}\label{sub:outercollaring}

For manifolds with boundary, it is well-known that one can construct collars. However, it is not easy to find collars for more complicated objects, for example, smooth manifolds with corners. Here we would like to adopt the outercollaring construction from \cite[Section 17]{FOOO_Kuranishi} to avoid such discussions. Instead of finding collars in the given geometric object, we add collars outside. In general, for each category $\uds{\bf C}$ listed at the beginning of this section, we would like to construct monoidal functors called the {\bf $\epsilon$-outercollaring}
\beqn
\{\cdot\}: \uds{\bf C} \to \outer \uds{\bf C}.
\eeqn
Some care must be taken when we are in the smooth category.

\subsubsection{Outercollaring of stratified spaces}

We first discuss the outercollaring of objects in $\uds{\bf Top}$. Let $A$ be a regular poset and $X$ be an $A$-space. Upon choosing a nonnegative number $\epsilon \geq 0$, we would like to define a new $A$-space denoted by $\outer X$ called {\bf the outer collaring of $X$ of width $\epsilon$} such that the association $X \mapsto \outer X$ satisfies certain functorial properties (see Lemma \ref{outer_product}). 

We first assume that $A$ has a unique maximal element $\gamma$, i.e., $X$ has a single top stratum. Abbreviate $\mb{Face}_\alpha = \mb{Face}_\alpha^\gamma$. Define
\beqn
\outer X = \left( \bigsqcup_{\alpha \in A} \partial^\alpha X \times [-\epsilon, 0]^{\mb{Face}_\alpha} \right)/ \sim,
\eeqn
where  the equivalence relation $\sim$ is generated by the following relation: if $\alpha \leq \beta$ (which implies $\mb{Face}_\beta \subseteq \mb{Face}_\alpha$), we identify 
\beqn
\big( x, (t_i)_{i \in \mb{Face}_\beta} \big) \in \partial^\beta X \times [-\epsilon, 0]^{\mb{Face}_\beta}
\eeqn
with
\beqn
\big( y, (s_j)_{j\in \mb{Face}_\alpha} \big) \in \partial^\alpha X \times [-\epsilon, 0]^{\mb{Face}_\alpha}
\eeqn
if $x = y \in \partial^\alpha X$, $s_j = 0$ when $j \notin \mb{Face}_\beta$ and $t_j = s_j$ when $j \in \mb{Face}_\beta$. We call $\outer X$ the {\bf outer collaring} of $X$ of width $\epsilon$ (see Figure \ref{outercollar}).

\begin{center}
\begin{figure}[h]
\begin{tikzpicture}

\draw [white, fill = lightgray] (-8, 2) -- (-8, 0) -- (-6, 0) -- (-6, 2);

\draw [very thick] (-8, 2) -- (-8, 0) -- (-6, 0);
\node at (-7, 1) {\scriptsize $X$};

\node at (-8, 2.2) {\scriptsize $\partial^{\beta_2} X$};
\node at (-5.6, 0) {\scriptsize $\partial^{\beta_1} X$};
\node at (-8.2, -0.2) {\scriptsize $\partial^\alpha X$};

\draw [very thick, ->] (-5, 0.3) -- (-2.5, 0.3);

\node at (0, 2.2) {\scriptsize $\partial^{\beta_2} X$};
\node at (2.4, -0.05) {\scriptsize $\partial^{\beta_1} X$};

\draw [white, fill = lightgray] (0, 2) -- (0, 0) -- (2, 0) -- (2, 2);
\draw [very thick] (0, 2) -- (0, 0) -- (2, 0);

\node at (1, 1) {\scriptsize $X$};

\draw [very thick] (-2, 2) -- (-2, -2) -- (2, -2);

\draw [dotted] (-2, 0) -- (0, 0) -- (0, -2);

\node at (-1, -1) {\scriptsize $\partial^\alpha X \times [-\epsilon, 0]^2$};
\node at (1, -1) {\scriptsize $\partial^{\beta_1} X \times [-\epsilon, 0]$};
\node at (-1, 1) {\scriptsize $\partial^{\beta_2} X \times [-\epsilon, 0]$};
\end{tikzpicture}
\caption{The local picture of an outer collaring of a manifold with corners with strata $\alpha < \beta_1, \beta_2$. The gray area is in the top stratum.}\label{outercollar}
\end{figure}
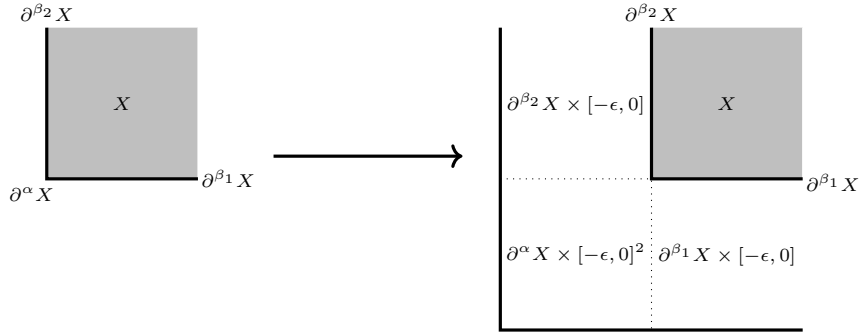
\end{center}

Now consider the general situation when the regular poset has arbitrarily many maximal elements. For each maximal $\gamma \in A^{\max}$, one has the outer collaring
\beqn
\outer (\partial^\gamma X).
\eeqn
Notice that if $\alpha \leq \gamma_1, \alpha \leq \gamma_2$, then there is a natural isomorphism 
\beqn
\outer \big( \partial^\alpha( \partial^{\gamma_1} X) \big) \cong \outer \big( \partial^\alpha ( \partial^{\gamma_2} X) \big).
\eeqn
Then define 
\beqn
\outer X:= \left( \bigsqcup_{\gamma \in A^{\max}} \outer (\partial^\gamma X) \right)/ \sim
\eeqn
where $\sim$ means we identify points in $\outer (\partial^{\gamma_1} X)$ and points in $\outer (\partial^{\gamma_2} X)$ with respect to the above natural isomorphisms. Figure \ref{outercollar2} describes an example of the outer collaring of a space $X$ with two top strata. 

\begin{center}
\begin{figure}[h]
\begin{tikzpicture}

\draw [white, fill = lightgray] (-10, 2) -- (-10, -2) -- (-8, -2) -- (-8, 2);

\draw [very thick] (-10, 2)-- (-10, -2);
\draw [very thick] (-10, -0) -- (-8, 0);
\node at (-9, 1) {\scriptsize $\partial^{\gamma_1} X$};
\node at (-9, -1) {\scriptsize $\partial^{\gamma_2} X$};

\node at (-10, 2.2) {\scriptsize $\partial^{\beta_1} X$};
\node at (-7.6, 0) {\scriptsize $\partial^{\beta_2} X$};
\node at (-10, -2.2) {\scriptsize $\partial^{\beta_3} X$};
\node at (-10.4, 0) {\scriptsize $\partial^\alpha X$};

\draw [very thick, ->] (-6.5, 0) -- (-3.5, 0);

\draw [white, fill = lightgray] (-2, 2.5) -- (-2, 0.5) -- (0, 0.5) -- (0, 2.5);
\draw [very thick] (-2, 2.5) -- (-2, 0.5) -- (0, 0.5);
\node at (-1, 1.5) {\scriptsize $\partial^{\gamma_1} X$};

\draw [white, fill = lightgray] (-2, -2.5) -- (-2, -0.5) -- (0, -0.5) -- (0, -2.5);
\draw [very thick] (-2, -2.5) -- (-2, -0.5) -- (0, -0.5);
\node at (-1, -1.5) {\scriptsize $\partial^{\gamma_2} X$};

\draw [very thick] (-2.5, 2.5) -- (-2.5, -2.5);

\draw [dotted] (-2.5, 0) -- (0, 0);
\draw [dotted] (-2.5, 0.5) -- (-2, 0.5) -- (-2, -0.5) -- (-2.5, -0.5);

\end{tikzpicture}
\caption{The local picture of an outer collaring of a composed manifold with corners with two top strata $\gamma_1, \gamma_2$ of dimension 2, three 1-dimensional strata $\beta_1, \beta_2, \beta_3$, and one zero-dimensional stratum $\alpha$.}\label{outercollar2}
\end{figure}
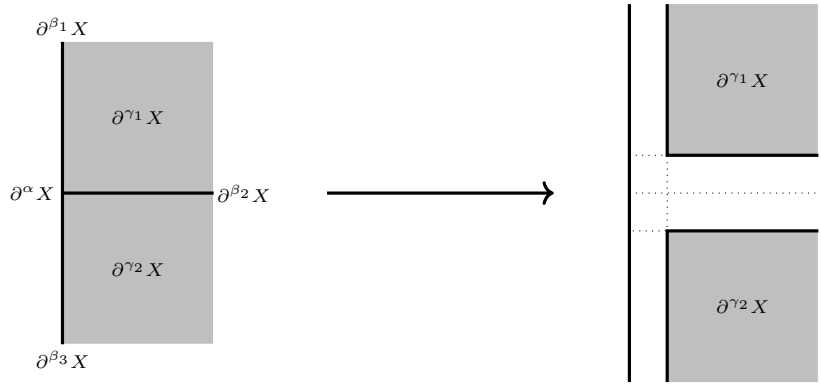
\end{center}

The following statement is straightforward from the definition.

\begin{lemma}\label{outer_product} The $\epsilon$-outer collaring is a regular stratification functor 
\beqn
\outer \{\cdot\}: \uds{\bf RegTop} \to \outer \uds{\bf RegTop}. \qed
\eeqn
\end{lemma}

The outer collar construction is compatible with group actions. If $G$ acts on $X$, then for each $g \in G$ and $(x, (t_i)_{i \in \mb{Face}_\alpha} )\in \partial^\alpha X \times [-\epsilon, 0]^{\mb{Face}_\alpha }$, define
\beqn
g \cdot (x, (t_i)_{i \in \bF_\alpha} ) = (gx, (t_i)_{i \in \bF_\alpha} ) \in \partial^\alpha X \times [-\epsilon, 0]^{\mb{Face}_\alpha } \subset \outer X.
\eeqn
It is easy to check that this gives a $G$-action on $\outer X$.

\subsection{Category of admissible smooth manifolds and orbifolds}\label{subsection_admissible_manifold}

The outercollaring of stratified topological manifolds is completely identical to the case of stratified spaces. However, the case with smooth manifolds is tricky.

We recall the notion of {\bf admissible manifolds with corners} introduced in \cite[Chapter 25]{FOOO_Kuranishi}. Let $V$ be a smooth manifold without boundary. For $t = (t_1, \ldots, t_k) \in [0, +\infty)^k$, denote 
\beqn
T_i:= e^{1/t_i}.
\eeqn 
Consider smooth functions 
\beqn
f: V \times [0, \epsilon)^k \to {\mb R}.
\eeqn
We consider the $C^m$-norm stratumwise using the above coordinates $T_i$ instead of $t_i$. Let 
\beqn
|f|_{C^m}
\eeqn
be the pointwise $C^m$-norm in the above sense. 

\begin{defn}(cf. \cite[Definition 25.3]{FOOO_Kuranishi}) We say $f: V \times [0, \epsilon)^k \to {\mb R}$ is {\bf admissible} if for each compact $K \subset V \times [0, \epsilon)^k$ and $m > 0$, there exist $\sigma(m, K)>0$ and $C(m, K)>0$ such that 
\beqn
\left| \frac{\partial f}{\partial T_i} \right|_{C^m} \leq C(m, K) e^{-\sigma(m, K) T_i},\ i = 1, \ldots, k.
\eeqn    
\end{defn}

Admissible functions are smooth (see \cite[Remark 25.5]{FOOO_Kuranishi}). Using the ring of admissible functions (see \cite[Lemma 25.13]{FOOO_Kuranishi}), we can consider corner charts on a smooth manifold with corners with admissible coordinate changes. In this way, one obtains a subcategory
\beqn
\uds{\bf SMan}_{\rm adm}
\eeqn
of $\uds{\bf SMan}$. Notice that the category $\outer \uds{\bf SMan}$ of collared smooth manifolds is automatically a subcategory of $\uds{\bf SMan}_{\rm adm}$.

Similarly, one can consider the admissible versions of the categories $\uds{\bf SKur}$ and $\uds{\bf dOrb}$. The key feature of the admissibility is that the outercollaring is well-defined (see \cite[Section 17.3]{FOOO_Kuranishi}). Indeed, for any $\epsilon>0$, there exist functors
\beqn
\outer \epsilon: \uds{\bf C}_{\rm adm} \to \outer \uds{\bf C}
\eeqn
where $\uds{\bf C}$ is either $\uds{\bf SMan}$, or $\uds{\bf SKur}$, or $\uds{\bf dOrb}$. For each admissible object $M$, let the result be called the {\bf $\epsilon$-outercollaring} of $M$. 

Finally, also notice that the standard collars are normally complex derived orbifolds. In fact manifolds are automatically normally complex orbifolds. Hence one has the collared versions
\begin{align*}
&\ \outer \uds{\bf dOrb}^{{\rm NC}},\ &\ \outer \uds{\bf dOrb}^{{\rm NC}}_{\rm rig}.
\end{align*}

\section{Flow Categories, Multimodules, and Homotopies}

\subsection{Flow categories}

The concept of flow categories was introduced by Cohen--Jones--Segal \cite{Cohen_Jones_Segal}. We need a variant of the original construction similar to \cite[Section 7]{pardon-VFC} and \cite[Section 7,8]{Abouzaid_axiomatic}. See also the formulations in \cite{Bai_Xu_Arnold}.

\begin{defn}[Flow categories]\label{def:flow-cat}\label{defn_flow_category}
Let $\uds{\bf C}$ be a regular stratification category (Definition \ref{defn_regular_stratification_category}). A $\uds{\bf C}$-enriched {\bf flow category} ${\mb F}$ consists of the following elements.

\begin{enumerate}

\item The set of objects ${\rm Ob}{\mb F}$, which is a countable poset (not assumed to be regular).

\item For any pair $p \leq q$ of objects, a morphism space $M_{pq}^{\mb F}$, which is an object of $\uds{\bf C}$ stratified by a normal poset $A_{pq}^{\mb F}$ (Definition \ref{defn_normal_poset}).

\item For any triple $p \leq r \leq q$ of objects, the composition map is defined as follows. It consists of an injective poset map homogeneous of degree $1$
\beqn
A_{pr}^{\mb F} \times A_{rq}^{\mb F} \to A_{pq}^{\mb F}
\eeqn
with image 
\beqn
A_{prq}^{\mb F} \subset A_{pq}^{\mb F},
\eeqn
together with a morphism in $\uds{\bf C}$
\beqn
\iota_{prq}^{\mb F}: M_{pr}^{\mb F} \times M_{rq}^{\mb F} \to M_{pq}^{\mb F} |_{A_{prq}^{\mb F}}\footnote{In later discussions, we will omit the notation $|_{A_{prq}^{\mb F}}$ for simplicity.}.
\eeqn

\end{enumerate}
Moreover, the following conditions must be satisfied.\footnote{The first four conditions are usually required for flow categories. The remaining conditions are imposed in order to obtain the correct chain level information.}
\begin{enumerate}

\item $M_{pq}^{\mb F} \neq \emptyset$ if and only if  $p \leq q$.

\item $M_{pp}^{\mb F}$ is a monoidal unit of $\uds{\bf C}$.

\item When $p = r$ resp. $r = q$, the composition $M_{pp}^{\mb F} \times M_{pq}^{\mb F} \to M_{pq}^{\mb F}$ resp. $M_{pq}^{\mb F} \times M_{qq}^{\mb F} \to M_{pq}^{\mb F}$ is the natural isomorphism from the monoidal structure.
 
\item {\bf (Associativity)} For any quadruple $p \leq r \leq s \leq q$ of objects, the following diagram commutes, where the arrows are defined using composition maps:
\beqn
\xymatrix{   &    M_{pr}^{\mb F} \times M_{rs}^{\mb F} \times M_{sq}^{\mb F}  \ar[ld] \ar[rd] & \\
 M_{ps}^{\mb F} \times M_{sq}^{\mb F}  \ar[rd]  &  &  M_{pr}^{\mb F} \times M_{rq}^{\mb F} \ar[ld] \\
& M_{pq}^{\mb F}  & }
\eeqn

\item If $p<r_1, r_2 < q$ and $A_{pr_1 q}^{\mb F} \cap A_{pr_2 q}^{\mb F} \neq \emptyset$, then either $r_1 \leq r_2$ or $r_2 \leq r_1$.

\item Maximal elements of each $A_{prq}^{\mb F}$ are true boundary strata of $A_{pq}^{\mb F}$ (Definition \ref{defn_normal_poset}). 
\end{enumerate}
\end{defn}

\begin{example}
A {\bf trivial flow category} enriched in $\uds{\bf C}$, denoted by ${\mb O}$, is the flow category whose object set has a single element $e$ with self-morphism space $M_{ee}^{{\mb O}}$ being a monoidal unit of $\uds{\bf C}$.
\end{example}

\subsection{Flow multimodules and homotopies}

Flow bimodules are abstract packages of moduli spaces considered in situations such as continuation maps, PSS and SSP maps, etc. There is a more general notion of flow multimodules allowing multiple flow categories as inputs, which model situations such as the pair-of-pants product.

\begin{defn}\label{defn_multimodule}
Let $\uds{\bf C}$ be a regular stratification category. Let ${\mb F}_1, \ldots, {\mb F}_m$ and ${\mb F}'$ be flow categories \footnote{One could consider more than one output; however, that will not be used in this paper.} enriched in $\uds{\bf C}$. A {\bf multimodule} over $({\mb F}_1, \ldots, {\mb F}_m; {\mb F}')$, denoted by ${\mb X}$, consists of the following objects.
\begin{enumerate}

\item For each tuple of objects $p_1\in {\rm Ob}{\mb F}_1, \ldots, p_m \in {\rm Ob}{\mb F}_m$, $p' \in {\rm Ob} {\mb F}'$, an object $M_{p_1 \cdots p_m; p'}^{\mb X}$ of $\uds{\bf C}$ with underlying normal poset $A_{p_1 \cdots p_m; p'}^{\mb X}$.

\item For any $i$ and a pair of objects $p_i \leq q_i$ in ${\mb F}_i$, an injective poset map homogeneous of degree $1$
\beqn
A^{{\mb F}_i}_{p_i q_i} \times A_{p_1 \cdots p_{i-1} q_i p_{i+1} \cdots p_m; p'}^{\mb X} \to A_{p_1 \cdots p_m; p'}^{\mb X}
\eeqn
with image
\beqn
A_{p_1 \cdots p_{i-1}(p_i q_i)p_{i+1} \cdots p_m; p'}^{\mb X} \subset A_{p_1 \cdots p_m; p'}^{\mb X},
\eeqn
together with a morphism in $\uds{\bf C}$
\beq\label{module_structural_map}
\iota_{p_1 \cdots p_{i-1}(p_i q_i) p_{i+1} \cdots p_m; p'}^{\mb X}: M_{p_i q_i}^{{\mb F}_i} \times M_{p_1 \cdots p_{i-1} q_i p_{i+1} \cdots p_m; p'}^{\mb X} \to M_{p_1 \cdots p_m; p'}^{\mb X}|_{A_{p_1 \cdots p_{i-1}(p_i q_i)p_{i+1} \cdots p_m; p'}^{\mb X}}.
\eeq

\item Similar to above, for a pair of objects $q' \leq  p'$ in ${\mb F}'$, an injective poset map homogeneous of degree $1$
\beqn
A_{p_1 \cdots p_m; q'}^{\mb X}\times A_{q'p'}^{{\mb F}'} \to A_{p_1 \cdots p_m; p'}^{\mb X}
\eeqn
with image
\beqn
A_{p_1 \cdots p_m; q'p'}^{\mb X} \subset A_{p_1 \cdots p_m; p'}^{\mb X},
\eeqn
together with a morphism in $\uds{\bf C}$
\beq\label{module_structural_map_2}
\iota_{p_1 \cdots p_m; q'p'}^{\mb X}: M_{p_1 \cdots p_m; q'}^{\mb X}\times M_{q'p'}^{{\mb F}'} \to M_{p_1 \cdots p_m; p'}^{\mb X} |_{A_{p_1 \cdots p_m; q'p'}^{\mb X}}.
\eeq
\end{enumerate}
These objects need to satisfy the following conditions.
\begin{enumerate}

\item When $p_i = q_i$ (recall $M_{pp}^{\mb F}$ is a monoidal unit), the morphism \eqref{module_structural_map} is the natural isomorphism
\beqn
M_{p_i p_i}^{{\mb F}_i} \times M_{p_1 \cdots p_m; p'}^{\mb X} \to M_{p_1 \cdots p_m; p'}^{\mb X}.
\eeqn
Similarly, when $q' = p'$, the morphism \eqref{module_structural_map_2} is the natural isomorphism
\beqn
M_{p_1 \cdots p_m; p'}^{\mb X} \times M_{p'p'}^{{\mb F}'} \to M_{p_1 \cdots p_m; p'}^{\mb X}.
\eeqn

\item The following associativity properties hold.
\begin{enumerate}
\item For $p_i < r_i < q_i$ in ${\mb F}_i$, the following diagram commutes.
\beq\label{module_associativity_1}
\vcenter{\xymatrix{  M_{p_i r_i}^{{\mb F}_i} \times M_{r_i q_i}^{{\mb F}_i} \times M_{p_1 \cdots p_{i-1} q_i p_{i+1} \cdots p_m; p'}^{\mb X} \ar[rr] \ar[d]  & & 
M_{p_i q_i}^{{\mb F}_i} \times M_{p_1 \cdots p_{i-1} q_i p_{i+1} \cdots p_m; p'}^{\mb X} \ar[d]\\
 M_{p_i r_i}^{{\mb F}_i} \times M_{p_1 \cdots p_{i-1} r_i p_{i+1} \cdots p_m; p'}^{\mb X} \ar[rr]  & & M_{p_1 \cdots p_m;p'}^{\mb X} }}
\eeq

\item For $i < j$, $p_i < q_i$ in ${\mb F}_i$ and $p_j < q_j$ in ${\mb F}_j$, the following diagram commutes.
\beq\label{module_associativity_2}
\vcenter{ \xymatrix{  M_{p_i q_i}^{{\mb F}_i} \times M_{p_j q_j}^{{\mb F}_j}  \times M_{p_1 \cdots p_{i-1} q_i p_{i+1} \cdots p_{j-1} q_j p_{j+1} \cdots p_m; p'}^{\mb X}  \ar[rr] \ar[d] & & 
M_{p_j q_j}^{{\mb F}_j} \times M_{p_1 \cdots p_{j-1} q_j p_{j+1} \cdots p_m;p'}^{{\mb X}} \ar[d] \\
 M_{p_i q_i}^{{\mb F}_i} \times M_{p_1 \cdots p_{i-1} q_i p_{i+1} \cdots p_m; p'}^{\mb X} \ar[rr]  & &  M_{p_1\cdots p_m;p'}^{\mb X}  }}
\eeq

\item For $q' < r' < p'$ in ${\mb F}'$, the following diagram commutes.
\beq\label{module_associativity_3}
\vcenter{ \xymatrix{  M_{p_1 \cdots p_m; q'}^{{\mb X}}\times M_{q'r'}^{{\mb F}'} \times M_{r'p'}^{{\mb F}'}  \ar[rr] \ar[d] & & M_{p_1 \cdots p_m; r'}^{{\mb X}} \times M_{r'p'}^{{\mb F}'} \ar[d]\\
M_{p_1 \cdots p_m; q'}^{{\mb X}}\times M_{q'p'}^{{\mb F}'} \ar[rr]  & & M_{p_1 \cdots p_m; p'}^{{\mb X}} }}
\eeq

\item For $p_i < q_i$ in ${\mb F}_i$ and $q' < p'$ in ${\mb F}'$, the following diagram commutes.
\beq
\vcenter{ \xymatrix{  M_{p_i q_i}^{{\mb F}_i} \times M_{p_1 \cdots p_{i-1} q_i p_{i+1} \cdots p_m; q'}^{{\mb X}} \times M_{q'p'}^{{\mb F}'}  \ar[rr] \ar[d]  & &    M_{p_1 \cdots p_m; q'}^{{\mb X}} \times M_{q'p'}^{{\mb F}'} \ar[d] \\
 M_{p_i q_i}^{{\mb F}_i} \times M_{p_1 \cdots p_{i-1} q_i p_{i+1} \cdots p_m; p'}^{{\mb X}}  \ar[rr] & &  M_{p_1 \cdots p_m; p'}^{{\mb X}}   }}
\eeq
\end{enumerate}

\item If $A_{p_1 \cdots p_{i-1} (p_i q_i) p_{i+1} \cdots p_m; p'}^{\mb X} \cap A_{p_1 \cdots p_{i-1}(p_i r_i) p_{i+1} \cdots p_m; p'}^{\mb X} \neq \emptyset$, then either $r_i \leq q_i$ or $q_i \leq r_i$. In the former case, one has 
\beqn
A_{p_1 \cdots p_{i-1} (p_i q_i) p_{i+1} \cdots p_m; p'}^{\mb X} \cap A_{p_1 \cdots p_{i-1}(p_i r_i) p_{i+1} \cdots p_m; p'}^{\mb X} = A_{p_1 \cdots p_{i-1}(p_i r_i q_i) p_{i+1} \cdots p_m; p'}^{{\mb X}}.
\eeqn
In the latter case, 
\beqn
A_{p_1 \cdots p_{i-1} (p_i q_i) p_{i+1} \cdots p_m; p'}^{\mb X} \cap A_{p_1 \cdots p_{i-1}(p_i r_i) p_{i+1} \cdots p_m; p'}^{\mb X} = A_{p_1 \cdots p_{i-1}(p_i q_i r_i) p_{i+1} \cdots p_m; p'}^{{\mb X}}.
\eeqn
\item If $A_{p_1\cdots p_m; q'p'}^{\mb X} \cap A_{p_1 \cdots p_m; r'p'}^{\mb X} \neq \emptyset$, then either $r' \leq q'$ or $q' \leq r'$. In the former case, one has 
\beqn
A_{p_1\cdots p_m; q'p'}^{\mb X} \cap A_{p_1 \cdots p_m; r'p'}^{\mb X} = A_{p_1\cdots p_m; r'q'p'}^{\mb X}.
\eeqn
In the latter case, 
\beqn
A_{p_1\cdots p_m; q'p'}^{\mb X} \cap A_{p_1 \cdots p_m; r'p'}^{\mb X} = A_{p_1 \cdots p_m; q'r'p'}^{\mb X}.
\eeqn

\item When $p_i < q_i$, the morphism \eqref{module_structural_map} has codimension 1 and maximal elements of $A_{p_1 \cdots p_{i-1} (p_i q_i) p_{i+1} \cdots p_m; p'}^{\mb X}$ are true boundary strata. Similarly, when $q' < p'$, the morphism \eqref{module_structural_map_2} has codimension 1 and maximal elements of $A_{p_1 \cdots p_m; q'p'}^{\mb X}$ are true boundary strata.

\end{enumerate}
\end{defn}

\begin{defn}
A multimodule over $({\mb F}; {\mb F}')$ is called a {\bf bimodule}, usually denoted by ${\mb B}$.
\end{defn}

\subsubsection{Products of multimodules}

In the above we only defined the situation with one output---this is the main scenario that will appear in this paper. On the other hand, there is no difficulty in extending the definition to multimodules with multiple output flow categories.

In general, let ${\mb F}_1, \ldots, {\mb F}_m; {\mb F}_1', \ldots, {\mb F}_n'$ be flow categories enriched over $\uds{\bf C}$. A multimodule over $({\mb F}_1, \ldots, {\mb F}_m; {\mb F}_1', \ldots, {\mb F}_n')$, denoted by ${\mb X}$, consists of objects
\beqn
M_{p_1 \cdots p_m; p_1' \cdots p_n'}^{\mb X} \in {\rm Ob} \uds{\bf C}
\eeqn
together with structural maps
\beqn
M_{p_i q_i}^{{\mb F}_i} \times M_{p_1 \cdots p_{i-1} q_i p_{i+1} \cdots p_m; p_1' \cdots p_n'}^{\mb X} \to M_{p_1 \cdots p_m; p_1' \cdots p_n'}^{\mb X}
\eeqn
and
\beqn
M_{p_1 \cdots p_m; p_1' \cdots p_{i-1}' q_i' p_{i+1}' \cdots p_n'}^{\mb X} \times M_{q_i' p_i'}^{{\mb F}_i'} \to M_{p_1 \cdots p_m; p_1' \cdots p_n'}^{\mb X}
\eeqn
which satisfy conditions similar to those in Definition \ref{defn_multimodule}. We omit the details.

It is straightforward to define products of multimodules. Suppose ${\mb X}$ is a multimodule over $({\mb F}_1, \ldots, {\mb F}_m; {\mb F}')$ and ${\mb Y}$ is a multimodule over $({\mb G}_1, \ldots, {\mb G}_n; {\mb G}')$. Then for objects $p_i \in {\rm Ob}{\mb F}_i$, $q_j \in {\rm Ob}{\mb G}_j$, $p' \in {\rm Ob}{\mb F}'$, $q' \in {\rm Ob}{\mb G}'$, the collection of objects
\beqn
M_{p_1 \cdots p_m q_1 \cdots q_n; p'q'}^{{\mb X} \times {\mb Y}}:=  M_{p_1 \cdots p_m; p'}^{\mb X} \times M_{q_1 \cdots q_n; q'}^{\mb Y}
\eeqn
define a multimodule ${\mb X} \times {\mb Y}$ over $({\mb F}_1, \cdots, {\mb F}_m, {\mb G}_1, \ldots, {\mb G}_n; {\mb F}', {\mb G}')$.

\subsubsection{Homotopies of multimodules}
For the following definition, denote by $[0,1]$ the homogeneous poset with one unique maximal element, thought as the total space of the unit interval, and two codimension $1$ element, which are the boundary points $\{0\}$ and $\{1\}$.

\begin{defn}\label{homotopy_modules}
Given two multimodules ${\mb X}_0$ and ${\mb X}_1$ over $({\mb F}_1, \ldots, {\mb F}_m; {\mb F}')$, a {\bf homotopy} of multimodules  from ${\mb X}_0$ to ${\mb X}_1$, denoted by ${\mb H}$, consists of the following objects.
\begin{enumerate}

\item For each tuple $(p_1, \ldots, p_m; p')$ of objects of involved flow categories, an object $M_{p_1 \cdots p_m; p'}^{{\mb H}}$ in $\uds{\bf C}$ stratified by a regular poset $A_{p_1 \cdots p_m; p'}^{\mb H} \times [0,1]$.

\item For a tuple $(p_1, \ldots, p_m; p')$, a morphism in $\uds{\bf C}$
\beqn
 \iota_{p_1\cdots p_m; p'}^{{\mb X}_0 \to {\mb H}} \sqcup \iota_{p_1 \cdots p_m; p'}^{{\mb X}_1\to {\mb H}}: M_{p_1 \cdots p_m; p'}^{{\mb X}_0} \sqcup M_{ p_1 \cdots p_m; p'}^{{\mb X}_1}  \to M_{p_1\cdots p_m; p'}^{\mb H}|_{A_{p_1 \cdots p_m; p'}^{\mb H} \times \{0\}} \sqcup M_{p_1\cdots p_m; p'}^{\mb H}|_{A_{p_1 \cdots p_m; p'}^{\mb H} \times \{1\}}.
\eeqn

\item For a tuple $(p_1, \ldots, p_m; p')$ and $p_i \leq q_i$ in ${\mb F}_i$, $q' \leq p'$, injective poset maps homogeneous of degree $1$
\beqn
A_{p_i q_i}^{{\mb F}_i} \times A_{p_1 \cdots p_{i-1} q_i p_{i+1} \cdots p_m; p'}^{\mb H}  \to A_{p_1 \cdots p_m; p'}^{{\mb H}}
\eeqn
and
\beqn
A_{p_1 \cdots p_m; q'}^{\mb H} \times A_{q'p'}^{{\mb F}'} \to A_{p_1 \cdots p_m; p'}^{\mb H},
\eeqn
together with morphisms in $\uds{\bf C}$ after pulling back the target to using the above poset maps (which we omit from the notations)
\beq\label{homotopy_composition_1}
\iota_{p_1 \cdots p_{i-1}(p_iq_i) p_{i+1}\cdots p_m; p'}^{\mb H}: M_{p_i q_i}^{{\mb F}_i} \times M_{p_1 \cdots p_{i-1} q_i p_{i+1} \cdots p_m; p'}^{\mb H}  \to M_{p_1 \cdots p_m; p'}^{{\mb H}}
\eeq
and 
\beq\label{homotopy_composition_2}
\iota_{p_1 \cdots p_m; q'p'}^{\mb H}: M_{p_1 \cdots p_m; q'}^{\mb H} \times M_{q'p'}^{{\mb F}'} \to M_{p_1 \cdots p_m; p'}^{\mb H}.
\eeq
\end{enumerate}
These objects need to satisfy the following conditions.
\begin{enumerate}

\item When $p_i = q_i$ resp. $q' = p'
$, the morphism \eqref{homotopy_composition_1} resp. \eqref{homotopy_composition_2} is the natural isomorphism from the monoidal structure of $\uds{\bf C}$.

\item {\bf (Associativity H1)} The following diagram commutes.
\beqn
\xymatrix{  M_{p_i q_i}^{{\mb F}_i} \times M_{p_1\cdots p_{i-1} q_i p_{i+1} \cdots p_m; q'}^{\mb H} \times M_{q'p'}^{{\mb F}'} \ar[r] \ar[d] &   M_{p_1 \cdots p_m; q'}^{\mb H} \times M_{q'p'}^{{\mb F}'} \ar[d]
\\
M_{p_i q_i}^{{\mb F}_i} \times M_{p_1 \cdots p_{i-1} q_i p_{i+1} \cdots p_m; p'}^{\mb H} \ar[r] & M_{p_1 \cdots p_m; p'}^{\mb H} }
\eeqn
Hence we have a well-defined morphism 
\beqn
M_{p_i q_i}^{{\mb F}_i} \times M_{p_1 \cdots p_{i-1} q_i p_{i+1} \cdots p_m; q'}^{\mb H} \times M_{q'p'}^{{\mb F}'} \to M_{p_1\cdots p_m; p'}^{\mb H}.
\eeqn

\item {\bf (Associativity H2)} For $\alpha = 0, 1$ and $i = 1, \ldots, m$, the following diagram commutes.
\begin{align*}
\xymatrix{  M_{p_i q_i}^{{\mb F}_i} \times M_{p_1 \cdots p_{i-1} q_i p_{i+1}\cdots p_m; q'}^{{\mb X}_\alpha} \times M_{q'p'}^{{\mb F}'}  \ar[r] \ar[d] & M_{p_1 \cdots p_m; q}^{{\mb X}_\alpha} \ar[d]  \\
  M_{p_i q_i}^{{\mb F}_i} \times M_{p_1 \cdots p_{i-1} q_i p_{i+1} \cdots p_m; q'}^{\mb H} \times M_{q'p'}^{{\mb F}'} \ar[r] &  M_{p_1 \cdots p_m; p'}^{{\mb H}} }
\end{align*}

\item {\bf (Associativity H3)}  For $p_i \leq r_i \leq q_i$ in ${\mb F}_i$ and $q'\leq r'\leq p'$ in ${\mb F}'$, the following diagram commutes.
\beq\label{homotopy_associativity_1}
\xymatrix{  M_{p_i r_i}^{{\mb F}_i} \times M_{r_i q_i}^{{\mb F}_i} \times M_{p_1 \cdots p_{i-1} q_i p_{i+1}\cdots p_m; q'}^{\mb H} \times M_{q'r'}^{{\mb F}'} \times M_{r'p'}^{{\mb F}'} \ar[r] \ar[d] &   M_{p_i r_i}^{{\mb F}_i} \times M_{p_1 \cdots p_{i-1} r_i p_{i+1} \cdots p_m; r'}^{\mb H} \times M_{r'p'}^{{\mb F}'} \ar[d] \\
M_{p_iq_i}^{{\mb F}_i} \times M_{p_1 \cdots p_{i-1} q_i p_{i+1}\cdots p_m; q'}^{{\mb H}} \times M_{q'p'}^{{\mb F}'} \ar[r] & M_{p_1\cdots p_m; p'}^{\mb H} }
\eeq

\item {\bf (Associativity H4)} For $i< j$, $p_i \leq q_i$ in ${\mb F}_i$ and $p_j \leq q_j$ in ${\mb F}_j$, the following diagram commutes.
\beqn
\xymatrix{  M_{p_i q_i}^{{\mb F}_i} \times M_{p_j q_j}^{{\mb F}_j} \times M_{p_1 \cdots p_{i-1} q_i p_{i+1} \cdots p_{j-1} q_j p_{j+1} \cdots p_m; p'}^{\mb H} \ar[r]  \ar[d]   &   M_{p_i r_i}^{{\mb F}_i} \times M_{p_1 \cdots p_{i-1} q_i p_{i+1} \cdots p_m; p'}^{\mb H} \ar[d]  \\
 M_{p_j q_j}^{{\mb F}_j} \times M_{p_1 \cdots p_{j-1} q_j p_{j+1} \cdots p_m; p'}^{\mb H} \ar[r]  & M_{p_1 \cdots p_m; p'}^{\mb H}} 
\eeqn
\end{enumerate}
\end{defn}

\subsection{Lifts of flow categories, multimodules, and homotopies}

The key issue of constructing Floer chain complex (and various chain maps, homotopies) is to regularize certain moduli spaces. Using the terminologies we just introduced, the associated flow categories (and multimodules, homotopies), whose morphism spaces are certain moduli spaces, are enriched in $\uds{\bf Top}$. The regularization requires a lift to a better category which allows manipulations in differential topological or algebraic topological sense. 

Let $\uds{\hat{\bf C}}$ and $\uds{\bf C}$ be two regular stratification categories (Definition \ref{defn_regular_stratification_category}). Let ${\bf F}: \uds{\hat{\bf C}} \to \uds{\bf C}$ be a covariant functor satisfying the natural compatibility conditions.
\begin{enumerate}
    \item ${\bf F}$ respects the distributive monoidal structure. \vspace{0.2cm}

    \item The diagram $    \xymatrix{ \underline{\hat{\bf C}} \ar[r]_{\bf F} \ar@<0pt> `u/10pt[r] `[rr] [rr] & \underline{\bf C} \ar[r]               & \underline{\bf RegPos}}$    commutes.

    \item ${\bf F}$ respects the pullback. Namely, if $\hat M$ is an $A$-stratified object of $\uds{\hat{\bf C}}$ and $B \subset A$, then one has the natural isomorphism
    \beqn
    {\bf F}(\hat M)|_B \cong {\bf F}(\hat M|_B).
    \eeqn
\end{enumerate}

\begin{defn}\label{defn:lifting}
Suppose ${\mb F}$ is a $\uds{\bf C}$-enriched flow category. A {\bf lift} of ${\mb F}$ with respect to ${\bf F}: \uds{\hat{\bf C}} \to \uds{\bf C}$ is a $\uds{\hat{\bf C}}$-enriched flow category $\hat{\mb F}$ satisfying the following conditions.
\begin{enumerate}
\item The set of objects are identical ${\rm Ob}\hat{\mb F} = {\rm Ob} {\mb F}$.

\item ${\bf F}(M_{pq}^{\hat {\mb F}}) = M_{pq}^{\mb F}$ and ${\bf F}_* \iota_{prq}^{\hat{\mb F}} = \iota_{prq}^{\mb F}$.
\end{enumerate}
\end{defn}

\begin{defn}\label{defn_bimodule_lifting}
Suppose ${\mb F}_1, \cdots, {\mb F}_m, {\mb F}'$ are $\uds{\bf C}$-enriched flow categories and ${\mb X}$ is a multimodule over $({\mb F}_1, \cdots, {\mb F}_m; {\mb F}')$. Suppose $\hat{\mb F}_1, \ldots, \hat{\mb F}_m, \hat{\mb F}'$ are lifts of ${\mb F}_1, \ldots, {\mb F}_m, {\mb F}'$ in $\uds{\hat{\bf C}}$ respectively. 
\begin{enumerate}

\item A {\bf lift} of ${\mb X}$ which {\bf extends} $\hat{\mb F}_1, \ldots, \hat{\mb F}_m; \hat{\mb F}'$ is a multimodule $\hat{\mb X}$ over $(\hat{\mb F}_1, \ldots, \hat{\mb F}_m; \hat{\mb F}')$ such that for any tuple of objects $p_1, \ldots, p_m; p'$ of the involved flow categories and $p_i \leq q_i$, $q' \leq p'$, one has 
\beqn
{\bf F}( M_{p_1\cdots p_m; p'}^{\hat{\mb X}}) = M_{p_1\cdots p_m; p'}^{\mb X}
\eeqn
and 
\begin{align*}
&\  {\bf F}_* \iota_{p_1 \cdots p_{i-1} (p_i q_i) p_{i+1} \cdots p_m; p'}^{\hat{\mb X}} = \iota_{p_1\cdots p_{i-1}(p_i q_i) p_{i+1} \cdots p_m; p'}^{\mb X},\ &\ {\bf F}_* \iota_{p_1 \cdots p_m; q'p'}^{\hat{\mb X}} = \iota_{p_1\cdots p_m; q'p'}^{\mb X}.
\end{align*}

\item Suppose $\hat{\mb X}_0$ and $\hat{\mb X}_1$ are two such lifts of multimodules ${\mb X}_0$ and ${\mb X}_1$ over $({\mb F}_1, \cdots, {\mb F}_m; {\mb F}')$ respectively. Let ${\mb H}$ be a homotopy from ${\mb X}_0$ to ${\mb X}_1$. A lift of ${\mb H}$ to $\uds{\hat{\bf C}}$ which extends $\hat{\mb F}$, $\hat{\mb F}'$, $\hat{\mb X}_0$, and $\hat{\mb X}_1$ is a homotopy $\hat{\mb H}$ from $\hat{\mb X}_0$ to $\hat{\mb X}_1$ such that for any tuple of objects $p_1, \cdots, p_m; p'$, $q', p' \in {\rm Ob}{\mb F}'$, one has
\beqn
{\bf F}(M_{p_1 \cdots p_m; p'}^{\hat{\mb H}}) = M_{p_1\cdots p_m; p'}^{\mb H},
\eeqn
\begin{align*}
&\ {\bf F}_*  \iota_{p_1 \cdots p_{i-1} (p_i q_i) p_{i+1} \cdots p_m; p'}^{\hat{\mb H}} = \iota_{p_1 \cdots p_{i-1}(p_i q_i) p_{i+1} \cdots p_m; p'}^{\mb H},\  &\ {\bf F}_*  \iota_{p_1 \cdots p_m; q'p'}^{\hat{\mb H}} = \iota_{p_1\cdots p_m; q'p'}^{\mb H},
\end{align*}
and 
\begin{align*}
&\ {\bf F}_* \iota_{p_1\cdots p_m; p'}^{\hat{\mb X}_0 \to \hat{\mb H}} = \iota_{p_1\cdots p_m; p'}^{{\mb X}_0 \to {\mb H}},\ &\ 
{\bf F}_* \iota_{p_1\cdots p_m; p'}^{\hat{\mb X}_1 \to \hat{\mb H}} = \iota_{p_1\cdots p_m; p'}^{{\mb X}_1 \to {\mb H}}.
\end{align*}
\end{enumerate}
\end{defn}

\subsection{Outercollaring of flow categories, multimodules, and homotopies}

Let $\uds{\bf C}$ be a regular stratification category which admits the outercollaring construction. Let ${\mb F}$ be a flow category enriched in $\uds{\bf C}$. Fix an $\epsilon>0$ (the {\bf outercollaring width}). Define the $\epsilon$-outercollaring of ${\mb F}$, which is a flow category $\outer   {\mb F}$ enriched in $\outer \uds {\bf C}$, as follows. The set of objects is the same as ${\mb F}$. For each pair $p, q \in {\rm Ob}{\mb F}$, the morphism space is
\beqn
M_{pq}^{\outer {}_\epsilon {\mb F}}:= \outer  M_{pq}^{\mb F}.
\eeqn
For $p\leq r \leq q$, the composition is the outercollaring
\beqn
\iota_{prq}^{\outer   {\mb F}}:= \outer {}_\epsilon \iota_{prq}^{\mb F}: \outer   ( M_{pr}^{\mb F} \times M_{rq}^{\mb F}) \cong \outer   M_{pr}^{\mb F} \times \outer  M_{rq}^{\mb F} \cong M_{pr}^{\outer   {\mb F}} \times M_{rq}^{\outer  {\mb F}} \to \outer  M_{pq}^{{\mb F}} \cong M_{pq}^{\outer   {\mb F}}.
\eeqn
It is straightforward to check that this set of data satisfies the definition of flow categories. 

The case of outercollaring of multimodules and homotopies are similar. We omit the details. We remark that one must choose the same outercollaring width of involved flow categories to obtain a corresponding outercollaring of multimodules and homotopies.

\subsection{Multimodule concatenations}\label{subsection_concatenation}

\subsubsection{Bimodule concatenations}

Before we discuss the more general multimodule concatenation, we first consider the bimodule case, which was discussed in \cite[Section 5]{Abouzaid_Blumberg_2024}. 

We first describe the intuitive picture. Suppose ${\mb F}, {\mb F}', {\mb F}''$ are flow categories enriched in $\uds{\bf Top}$, ${\mb B}$ is a bimodule over $({\mb F}; {\mb F}')$ and ${\mb B}'$ is a bimodule over $({\mb F}'; {\mb F}'')$. For example, ${\mb B}$ and ${\mb B}'$ are bimodules corresponding to continuation maps. We would like to define a bimodule ${\mb B}\circ {\mb B}'$ over $({\mb F}; {\mb F}'')$ which characterizes the features of compositions of continuation maps. Then given a pair of objects $p \in {\rm Ob}{\mb F}$ and $p'' \in {\rm Ob}{\mb F}''$, the space $M_{p;p''}^{{\mb B}\circ {\mb B}'}$ should be the quotient
\beq\label{concatenation_quotient}
\left( \bigsqcup_{p' \in {\rm Ob}{\mb F}'} M_{p;p'}^{{\mb B}} \times M_{p';p''}^{{\mb B}'} \right)/ \sim
\eeq
where the equivalence relation is generated by viewing the images of points of
\beqn
M_{p;p'}^{{\mb B}} \times M_{p'q'}^{{\mb F}'} \times M_{q';p''}^{{\mb B}'}
\eeqn
in either $M_{p;p'}^{{\mb B}} \times M_{p';p''}^{{\mb B}'}$ or $M_{p;q'}^{{\mb B}} \times M_{q';p''}^{{\mb B}'}$ the same point. 

However, this construction does not work for any category $\uds{\bf C}$. For example, when $\uds{\bf C} = \uds{\bf Kur}$, because of the symmetry groups, there is no obvious way to construct the above quotient in $\uds{\bf Kur}$. 

To rigorously define the notion of concatenations, one needs some preparations. Let $\mb{F}, \mb{F}', \mb{F}''$ be flow categories. Let $\mb{B}$ be a bimodule  over $({\mb F}; {\mb F}')$ and $\mb{B}'$ be a bimodule over $({\mb F}'; {\mb F}'')$. For each pair of objects $p\in {\rm Ob}{\mb F}$ and $p''\in {\rm Ob}{\mb F}''$, define the set
\beqn
A_{p;p''}^{{\mb B}' \circ {\mb B}}:= \left( \bigsqcup_{p' \in {\rm Ob} \mb{F}'}  A_{p;p'}^{\mb B} \times A_{p';p''}^{{\mb B}'} \right)/ \sim
\eeqn
where the equivalence relation $\sim$ is defined as follows. For any $p', q' \in {\rm Ob}{\mb F}'$ images of the two maps
\beqn
\xymatrix{ & A_{p;p'}^{\mb B} \times A_{p'q'}^{{\mb F}'} \times A_{q' ;p''}^{{\mb B}'}  \ar[ldd]_{{\rm Id}\times \iota_{p'q'; p''}^{{\mb B}'} } \ar[rdd]^{\iota_{p; p'q'}^{\mb B}\times {\rm Id} } \\ 
& & \\
   A_{p; p'}^{{\mb B}} \times A_{p'; p''}^{{\mb B}'}    & & A_{p; q'}^{{\mb B}} \times A_{q';p''}^{{\mb B}'} }
\eeqn
are identified. The equivalence relation $\sim$ is generated by this identification. 

\begin{lemma}\label{lemma_concatenation_equivalence}
Two pairs $(a', a'') \in A_{p;p'}^{{\mb B}} \times A_{p';p''}^{{\mb B}'}$ and $(b', b'')\in A_{p;q'}^{{\mb B}} \times A_{q';p''}^{{\mb B}'}$ are equivalent if and only if one of the following conditions holds: 1) $p' \leq q'$ and there exists $v \in A_{p'q'}^{{\mb F}'}$ such that 
\begin{align*}
&\ b' = \iota_{p;p'q'}^{\mb B}(a', v),\ &\ a'' = \iota_{p'q'; p''}^{{\mb B}'}(v, b'');
\end{align*}
2) $q' \leq p'$ and there exists $w \in A_{q'p'}^{{\mb F}'}$ such that 
\begin{align*}
&\ a' = \iota_{p; q'p'}^{\mb B}(b', w),\ &\ b'' = \iota_{q'p'; p''}^{{\mb B}'} (w, a'').
\end{align*}
\end{lemma}

\begin{proof}
    
We first check that the relation defined in the statement of this lemma, denoted by $\approx$, is an equivalence relation. The only non-obvious condition is the transitivity. Suppose we have three elements $(a', a'') \in A_{p; p'}^{\mb B} \times A_{p'; p''}^{{\mb B}'}$, $(b', b'') \in A_{p; q'}^{\mb B} \times A_{q'; p''}^{{\mb B}'}$, and $(c', c'') \in A_{p; r'}^{\mb B} \times A_{r'; p''}^{{\mb B}'}$. Suppose $(a', a'') \approx (b', b'')$ and $(b', b'') \approx (c', c'')$. Without loss of generality, we may assume one of the following holds.
\begin{enumerate}

\item $p' \leq q' \leq r'$. In this case, there exists $v \in A_{p' q'}^{{\mb F}'}$ and $w \in A_{q'r'}^{{\mb F}'}$ such that 
\begin{align*}
&\ b' = \iota_{p; p'q'}^{{\mb B}}(a', v),\ &\ a'' = \iota_{p'q'; p''}^{{\mb B}'}(v, b'')
\end{align*}
and 
\begin{align*}
&\ c' = \iota_{p; q'r'}^{{\mb B}}(b', w),\ &\ b'' = \iota_{q'r'; p''}^{{\mb B}'}(w, c'').
\end{align*}
Then by the associativity conditions for the bimodules, one has
\beqn
c' = \iota_{p; q'r'}^{\mb B}(b', w) = \iota_{p; q'r'}^{{\mb B}}(\iota_{p; p'q'}^{{\mb B}}(a', v), w) = \iota_{p; p'r'}^{{\mb B}} (a', \iota_{p'q'r'}^{{\mb F}'}(v, w))
\eeqn
and
\beqn
a'' = \iota_{p'q'; p''}^{{\mb B}'}(v, b'') = \iota_{p'q'; p''}^{{\mb B}'}(v, \iota_{q'r'; p''}^{{\mb B}'}(w, c'')) = \iota_{p'r'; p''}^{{\mb B}'}( \iota_{p'q'r'}^{{\mb F}'}(v, w), c'').
\eeqn
Hence $(a', a'') \approx (c',c'')$.

\item $p' \leq q'$ and $r' \leq q'$. In this case, there exists $v \in A_{p'q'}^{{\mb F}'}$ and $w\in A_{r'q'}^{{\mb F}'}$ such that 
\beqn
b = \iota_{p; p'q'}^{{\mb B}}(a, v),\ a' = \iota_{p'q'; p''}^{{\mb B}'}(v, b')
\eeqn
and
\beqn
b = \iota_{p; r'q'}^{{\mb B}}(c, w),\ c' = \iota_{r'q'; p''}^{{\mb B}'}(w, b').
\eeqn
Then we see 
\beqn
b \in \partial^{p; p'q'} A_{p;q'}^{{\mb B}} \cap \partial^{p; r'q'} A_{p; q'}^{{\mb B}} \neq \emptyset.
\eeqn
Then by the condition for bimodule/multimodules (see (3) and (4) of Definition \ref{defn_multimodule}), either $r' \leq q'$ or $q' \leq r'$. Without loss of generality, assume that $p' \leq r'$. Then one has
\beqn
b  \in A_{p;p'r'q'}^{\mb B}
\eeqn
As the poset maps underlying the structural maps are all injective (see Definition \ref{defn_flow_category} and Definition \ref{defn_multimodule}), it follows that there exists $u \in A_{p'r'}^{{\mb F}'}$ such that 
\beqn
c = \iota_{p; p'r'}^{{\mb B}}(a, u),\ v = \iota_{p'r'q'}^{{\mb F}'}(u, w). 
\eeqn
Then by using the associativity of the bimodules, one can conclude that
\beqn
a' = \iota_{p'q'; p''}^{{\mb B}'} (v, b') = \iota_{p'q'; p''}^{{\mb B}'} (\iota_{p'r'q'}^{{\mb F}'}(u, w), b') = \iota_{p'r'; p''}^{{\mb B}'}(u, \iota_{r'q'; p''}^{{\mb B}'}(w, b')) = \iota_{p'r'; p''}^{{\mb B}'}(u, c')
\eeqn
which implies $(a, a') \approx (c,c')$.
\end{enumerate}
Therefore, $\approx$ is an equivalence relation. 

Then one can see that the two relations $\sim$ and $\approx$ contain each other. Hence the lemma is proved. 
\end{proof}

The following corollary is an obvious consequence.

\begin{cor}\label{cor_module_concatenation}
The natural map
\beqn
A_{p; p'}^{\mb B} \times  A_{p'; p''}^{{\mb B}'} \to A_{p; p''}^{{\mb B}\circ {\mb B}'}
\eeqn
is injective.
\end{cor}

Denote the image of the above map by 
\beqn
A_{p; p'; p''}^{{\mb B} \circ {\mb B}'}.
\eeqn
More generally, assume $p_1' < \cdots < p_l'$ in ${\mb F}'$, denote
\beqn
A_{p; (p_1'\cdots p_l'); p''}^{{\mb B} \circ {\mb B}'} = \bigcap_{i =1}^l A_{p; p_i'; p''}^{{\mb B}\circ {\mb B}'}. 
\eeqn

The following corollary is also obvious.

\begin{cor}
If $A_{p; p'; p''}^{{\mb B} \circ {\mb B}'} \cap A_{p; q'; p''}^{{\mb B}\circ {\mb B}'} \neq \emptyset$, then either $p' \leq q'$ or $q' \leq p'$.
\end{cor}

We define a partial order between two elements $x, y\in A_{p; p''}^{{\mb B} \circ {\mb B}'}$ as follows. We denote $x \leq y$ if there exist $p' \in {\rm Ob}{\mb F}'$, $x', y' \in A_{p'p'}^{\mb B}$, $x'', y''\in A_{p';p''}^{{\mb B}'}$ such that $x' \leq y'$ and $x'' \leq y''$ and such that $x$ is represented by $(x', x'')$ and $y$ is represented by $(y', y'')$. 

\begin{lemma}\label{lemma_concatenation_partial_order}
This is a partial order. Moreover, the map $A_{p; p'}^{\mb B} \times A_{p'; p''}^{{\mb B}'} \to A_{p; p''}^{{\mb B} \circ {\mb B}'}$ is a poset map.
\end{lemma}

\begin{proof}
The only non-obvious condition is the transitivity. Suppose $x \leq y \leq z$, then there exists $p' \in {\rm Ob}{\mb F}'$ such that $x$ is represented by $(x', x'') \in A_{p; p'}^{{\mb B}} \times A_{p';p''}^{{\mb B}'}$, $y$ is represented by $(y', y'') \in A_{p; p'}^{\mb B}\times A_{p'; p''}^{{\mb B}'}$ such that $x' \leq y'$ and $x'' \leq y''$. There also exists $\tilde p' \in {\rm Ob}{\mb F}'$ such that $y$ is represented by $(\tilde y', \tilde y'') \in A_{p; \tilde p'}^{{\mb B}} \times A_{\tilde p'; p''}^{{\mb B}'}$ and $z$ is represented by $(\tilde z', \tilde z'') \in A_{p; \tilde p'}^{{\mb B}} \times A_{\tilde p'; p''}^{{\mb B}'}$ such that $\tilde y' \leq \tilde z'$ and $\tilde y'' \leq \tilde z''$. Since $(y', y'') \approx (\tilde y', \tilde y'')$, by Lemma \ref{lemma_concatenation_equivalence}, either $p' \leq \tilde p'$ or $\tilde p' \leq p'$. 

\begin{enumerate}

\item When $p' \leq \tilde p'$,  there exists $v \in A_{p'\tilde p'}^{{\mb F}'}$ such that 
\begin{align*}
&\ \tilde y' = \iota_{p; p'\tilde p'}^{{\mb B}}( y', v),\ &\  y'' = \iota_{p'\tilde p'; p''}^{{\mb B}'}(v, \tilde y'').
\end{align*}
Then define
\beqn
\tilde x' = \iota_{p; p'q'}^{{\mb B}}(x', v)
\eeqn
Then as $x' \leq y'$, one has $\tilde x' \leq \tilde y'$. Moreover, as $x'' \leq y'' = \iota_{p'q'; p''}^{{\mb B}'}(v, \tilde y'')$, it follows that $x'' = \iota_{p'q'; p''}^{{\mb B}'}(w, \tilde x'')$ for some $w$ and $\tilde x''$. Moreover, $\tilde x'' \leq \tilde y''$. It then follows that $(x', x'') \approx (\tilde x', \tilde x'')$ and $\tilde x' \leq \tilde y' \leq \tilde z'$ and $\tilde x'' \leq \tilde y'' \leq \tilde z''$. Hence $x \leq z$. 

\item When $\tilde p' \leq p'$, the argument is similar. There exists $w \in A_{\tilde p' p'}^{{\mb F}'}$ such that 
\begin{align*}
    &\ y' = \iota_{p; \tilde p' p'}^{\mb B}(\tilde y', w),\ &\ \tilde y'' = \iota_{\tilde p'p'; p''}^{{\mb B}'}(w, y'').
\end{align*}
Then define $\tilde x'' = \iota_{\tilde p'p'; p''}^{{\mb B}'}(w, x'')$. Then since $x'' \leq y''$, one has
\beqn
\tilde x'' = \iota_{\tilde p'p'; p''}^{{\mb B}'}(w, x'') \leq \iota_{\tilde p'p'; p''}^{{\mb B}'}(w, y'') = \tilde y'' \leq \tilde z''.
\eeqn
Moreover, as $x' \leq y' = \iota_{p; \tilde p'p'}^{\mb B}(\tilde y', w) \in A_{p; p'}^{\mb B}$, there exists $\tilde x' \in A_{p; \tilde p'}^{\mb B}$ and $u \in A_{\tilde p'p'}^{{\mb F}'}$ such that $x' =\iota_{p; \tilde p'p'}^{\mb B}(\tilde x', u)$ and $\tilde x' \leq \tilde y' \leq \tilde z'$. Hence $x \leq z$. 
\end{enumerate}
Lastly, it follows immediately from the definition of the partial order that the map $A_{p; p'}^{\mb B} \times A_{p'; p''}^{{\mb B}'} \to A_{p; p''}^{{\mb B} \circ {\mb B}'}$ is a poset map.
\end{proof}

\begin{lemma}
The poset $A_{p; p''}^{{\mb B}\circ {\mb B}'}$ is a normal poset.     
\end{lemma}

\begin{proof}
We first verify that this poset is homogeneous (see Definition \ref{defn_homogeneous_poset}). Let $x \in A_{p; p''}^{{\mb B} \circ {\mb B}'}$ be an arbitrary element. Notice that each $A_{p; p'; p''}^{{\mb B} \circ {\mb B}'} \cong A_{p; p'}^{{\mb B}}\times A_{p'; p''}^{{\mb B}'}$ is homogeneous. Hence if $x \in A_{p; p'; p''}^{{\mb B}\circ {\mb B}'}$, the depth of $x$ in this subset is well-defined, denoted by $d_{p'}(x)$. As maximal elements of $A_{p; p''}^{{\mb B} \circ {\mb B}'}$ are exactly maximal elements of $A_{p; p'; p''}^{{\mb B} \circ {\mb B}'}$ for all possible $p'$, we only need to show that $d_{p'}(x)$ is independent of $p'$. Suppose $x \in A_{p; p'; p''}^{{\mb B} \circ {\mb B}'} \cap A_{p; q'; p''}^{{\mb B} \circ {\mb B}'}$. Without loss of generality, we can assume $p' \leq q'$. Then $x \in A_{p; (p'q'); p''}^{{\mb B} \circ {\mb B}'}$. We need to show that $d_{p'}(x) = d_{q'}(x)$. 

{\it Claim.} Maximal elements of $A_{p; (p'q'); p''}^{{\mb B} \circ {\mb B}'}$ have depth $1$.

{\it Proof of the claim.} The depth is clearly at least one. If the depth is bigger than $1$, then there exists a non-maximal element $y$ such that $x< y$. Then there exists $r'$ different from $p', q'$ such that $x \in A_{p; r'; p''}^{{\mb B} \circ {\mb B}'}$. Then by the definition of the equivalence relation, the three elements $p', q', r'$ are totally ordered. In either case, $x$ cannot be a maximal element of $A_{p; (p'q'); p''}^{{\mb B} \circ {\mb B}'}$. \hfill {\it End of the proof of the claim.}

Therefore, both $d_{p'}(x)$ and $d_{q'}(x)$ are 1 higher than the depth of $x$ in the subset $A_{p; (p'q'); p''}^{{\mb B} \circ {\mb B}'}$. Hence $d_{p'}(x) = d_{q'}(x)$. Hence $A_{p; p''}^{{\mb B} \circ {\mb B}'}$ is a homogeneous poset.

Next we need to show that the poset is regular (see Definition \ref{defn_regular_poset}). Given any $\gamma \in A_{p; p''}^{{\mb B}\circ {\mb B}'}$, the three conditions for $\gamma$ in Definition \ref{defn_regular_poset} are about elements below $\gamma$. We may assume that $\gamma \in A_{p; p'; p''}^{{\mb B} \circ {\mb B}'}$, which is a regular poset. Hence all three conditions for $\gamma$ holds. Therefore, $A_{p; p''}^{{\mb B} \circ {\mb B}'}$ is regular.

Lastly, we verify that $A_{p; p''}^{{\mb B} \circ {\mb B}'}$ is normal (Definition \ref{defn_normal_poset}), i.e., there are only fake boundaries or true boundaries. Let $\gamma$ be an element of depth $1$.

{\it Claim. There exist at most two distinct $p' \in {\rm Ob}{\mb F}'$ such that $\gamma \in A_{p; p'; p''}^{{\mb B} \circ {\mb B}'}$.}

Suppose it is not the case, then there exists three distinct $p_1', p_2', p_3' \in {\rm Ob}{\mb F}'$ such that $\gamma$ is representable in $A_{p; p_i'}^{{\mb B}} \times A_{p_i'; p''}^{{\mb B}'}$ for $i = 1, 2, 3$. Then the three elements $p_1', p_2', p_3'$ are totally ordered; we may assume $p_1' < p_2' < p_3'$. Then it is easy to find an increasing path starting from $\gamma$ of length $2$. This contradicts the condition that $\gamma$ has depth $1$.

Now we discuss in two cases. First, suppose there exists a unique $p'$ such that $\gamma\in A_{p; p'; p''}^{{\mb B} \circ {\mb B}'}$, which is a normal poset. Hence $\gamma$ is adjacent to at most two maximal elements. Second, suppose $\gamma \in A_{p; p_i'; p''}^{{\mb B} \circ {\mb B}'}$ for $p_1' \neq p_2'$. Then we may assume $p_1' < p_2'$ and $\gamma \in A_{p; (p_1'p_2'); p''}^{{\mb B} \circ {\mb B}'}$. As $\gamma$ has depth 1, it is a maximal element of $A_{p; (p_1'p_2'); p''}^{{\mb B} \circ {\mb B}'}$. Hence by the properties of bimodule/multimodules (see (2) of Definition \ref{defn_multimodule}), it is a true boundary in both $A_{p; p_1'; p''}^{{\mb B} \circ {\mb B}'}$ and $A_{p; p_2'; p''}^{{\mb B}'}$. Hence $\gamma$ is adjacent to exactly two maximal elements in $A_{p; p''}^{{\mb B} \circ {\mb B}'}$. This finishes the proof that $A_{p; p''}^{{\mb B} \circ {\mb B}'}$ is normal.
\end{proof}

The above technical results allow us to do the following construction.

\begin{enumerate}

\item Let $\uds{\bf C}$ be either $\uds{\bf Top}$ or $\pman$. If we define 
\beqn
M_{p; p''}^{{\mb B} \circ {\mb B}'}:= \left( \bigsqcup_{p' \in {\rm Ob}{\mb F}'} M_{p; p'}^{{\mb B}} \times M_{p'; p''}^{{\mb B}'} \right)/ \sim
\eeqn
as topological spaces, then it is naturally stratified by $A_{p'; p''}^{{\mb B}'\circ {\mb B}}$. Moreover, if $\uds{\bf C} = \pman$, then $M_{p; p''}^{{\mb B}\circ {\mb B}'}$ is a stratified pseudomanifold. Moreover, one obtains a bimodule over $({\mb F}; {\mb F}'')$, denoted by ${\mb B} \circ {\mb B}'$. If the flow categories and bimodules are enriched in the collared version of $\uds{\bf C}$, then ${\mb B} \circ {\mb B}'$ is also collared.

\item If $\uds{\bf C} = \uds{\bf dOrb}$, let ${\mc U}_{p; p'}^{\mb B}$ and ${\mc U}_{p'; p''}^{{\mb B}'}$ be the underlying stratified orbifolds in $M_{p;p'}^{{\mb B}}$ and $M_{p'; p''}^{{\mb B}'}$ respectively. Then define
\beqn
{\mc U}_{p; p''}^{{\mb B}\circ {\mb B}'}:= \left( \bigsqcup_{p' \in {\rm Ob}{\mb F}'} {\mc U}_{p; p'}^{{\mb B}} \times {\mc U}_{p'; p''}^{{\mb B}'} \right) / \sim
\eeqn
in the similar way. Then it is a stratified orbifold. Moreover, the obstruction bundles ${\mc E}_{p;p'}^{\mb B} \to {\mc U}_{p;p'}^{\mb B}$ and ${\mc E}_{p'; p''}^{{\mb B}'} \to {\mc U}_{p'; p''}^{{\mb B}'}$ induce orbifold vector bundles ${\mc E}_{p; p''}^{{\mb B} \circ {\mb B}'} \to {\mc U}_{p; p''}^{{\mb B} \circ {\mb B}'}$, with induced Kuranishi sections 
\beqn
{\mc S}_{p;p''}^{{\mb B} \circ {\mb B}'}: {\mc U}_{p; p''}^{{\mb B} \circ {\mb B}'} \to {\mc E}_{p; p''}^{{\mb B}\circ {\mb B}'}.
\eeqn
Then one obtains objects $M_{p; p''}^{{\mb B}\circ {\mb B}'}$ in $\uds{\bf dOrb}$ together with structural maps making a bimodule ${\mb B}\circ {\mb B}'$ over $({\mb F}; {\mb F}'')$. Lastly, if the flow categories and bimodules are enriched in $\outer \uds{\bf dOrb}$, so is the concatenation ${\mb B} \circ {\mb B}'$.
\end{enumerate}

\subsubsection{General concatenations}

Now we consider concatenations in a general regular stratification category $\uds{\bf C}$. There is no canonical construction. We characterize the notion of concatenations in an axiomatic way. 

\begin{defn}\label{defn_general_concatenation}
Let ${\mb F}_1, \ldots, {\mb F}_m, {\mb G}_1, \ldots, {\mb G}_n, {\mb G}'$ be flow categories enriched in $\uds {\bf C}$. Let $i \in \{ 1, \ldots, n\}$. Let ${\mb X}$ resp. ${\mb Y}$ be a multimodule over $({\mb F}_1, \ldots, {\mb F}_m; {\mb G}_i)$ resp. over $({\mb G}_1, \ldots, {\mb G}_n; {\mb G}')$. Then a {\bf concatenation} of ${\mb X}$ and ${\mb Y}$ at ${\mb G}_i$, denoted by ${\mb X} \circ_i {\mb Y}$, consists of the following objects.
\begin{enumerate}
    \item A multimodule over $({\mb G}_1, \ldots, {\mb G}_{i-1}, {\mb F}_1, \ldots, {\mb F}_m, {\mb G}_{i+1}, \ldots, {\mb G}_n; {\mb G}')$.
    
    \item For each two tuples $(p_1, \ldots, p_m; r_i)$ and $(r_1, \ldots, r_n; r')$ of objects, a codimension zero morphism
    \beqn
    M_{p_1 \cdots p_m; r_i}^{\mb X} \times M_{r_1 \cdots r_n; r'}^{\mb Y} \to M_{r_1 \cdots r_{i-1} p_1 \cdots p_m r_{i+1} \cdots r_n; r'}^{{\mb X} \circ_i {\mb Y}}.
    \eeqn
    \end{enumerate}
    They are required to satisfy the following conditions.
\begin{enumerate}
    \item For each tuple $(r_1, \cdots, r_{i-1}, p_1, \cdots, p_m, r_{i+1}, \cdots, r_n; r')$, the object $M_{r_1 \cdots r_{i-1} p_1 \cdots p_m r_{i+1} \cdots r_n; r'}^{{\mb X} \circ_i {\mb Y}}$ is stratified by $A_{r_1 \cdots r_{i-1} p_1 \cdots p_m r_{i+1} \cdots r_n; r'}^{{\mb X} \circ_i {\mb Y}}$.

    \item The underlying poset maps of the structural maps of ${\mb X} \circ_i {\mb Y}$ is the same as the previously discussed special cases.

    \item The following diagram commutes.
    \beqn
    \xymatrix{ &    M_{p_1 \cdots p_m; r_i}^{{\mb X}} \times M_{r_i s_i}^{{\mb G}_i} \times M_{r_1 \cdots r_{i-1} s_i r_{i+1} \cdots r_n; r'}^{{\mb Y}} \ar[rd] \ar[ld]   & \\
     M_{p_1 \cdots p_m; s_i}^{{\mb X}} \times M_{r_1 \cdots r_{i-1} s_i r_{i+1}\cdots r_n; r'}^{\mb Y} \ar[rd] & &  M_{p_1 \cdots p_m; r_i}^{\mb X} \times M_{r_1 \cdots r_n; r'}^{\mb Y} \ar[ld] \\
     &  M_{r_1 \cdots r_{i-1} p_1 \cdots p_m r_{i+1} \cdots r_n; r'}^{{\mb X} \circ_i {\mb Y}}  &   }
    \eeqn

    \item The following diagram commutes.
    \beqn
    \xymatrix{
    M_{p_j q_j}^{{\mb F}_j} \times M_{p_1 \cdots p_{j-1} q_j p_{j+1} \cdots p_m; r_i}^{\mb X} \times M_{r_1 \cdots r_n; r'}^{{\mb Y}} \ar[r] \ar[d]   &    M_{p_j q_j}^{{\mb F}_j} \times M_{r_1 \cdots r_{i-1} p_1 \cdots p_{j-1} q_j p_{j+1} \cdots p_m r_{i+1} \cdots r_n; r'}^{{\mb X} \circ_i {\mb Y}} \ar[d] \\
    M_{p_1 \cdots p_m; r_i}^{{\mb X}} \times M_{r_1 \cdots r_n; r'}^{{\mb Y}} \ar[r] &   M_{r_1 \cdots r_{i-1} p_1 \cdots p_m r_{i+1} \cdots r_n; r'}^{{\mb X} \circ_i {\mb Y}} }
    \eeqn

    \item The following diagram commutes.
    \beqn
\xymatrix{
    M_{p_1 \cdots p_m; r_i}^{\mb X} \times M_{r_1 \cdots r_n; s'}^{\mb Y} \times M_{s'r'}^{{\mb G}'} \ar[r] \ar[d]   &  M_{r_1 \cdots r_{i-1} p_1 \cdots p_m r_{i+1} \cdots r_n; s'}^{{\mb X} \circ_i {\mb Y}} \times M_{s'r'}^{{\mb G}'} \ar[d]\\
    M_{p_1 \cdots p_m; r_i}^{{\mb X}} \times M_{r_1 \cdots r_n; r'}^{{\mb Y}} \ar[r] &   M_{r_1 \cdots r_{i-1} p_1 \cdots p_m r_{i+1} \cdots r_n; r'}^{{\mb X} \circ_i {\mb Y}} }
    \eeqn
    
\end{enumerate}
\end{defn}

We can see that the previously discussed canonical construction of concatenations satisfies the above definition. On the other hand, in many situations, the construction of the concatenation ${\mb X}\circ_i {\mb Y}$ is not canonical but involves choices one needs to make. 

\begin{rem}
    We would like to reiterate that in Definition \ref{defn_general_concatenation}, ${\mb X} \circ_i {\mb Y}$ is part of the data: we do not construct ${\mb X} \circ_i {\mb Y}$ using ${\mb X}$ and ${\mb Y}$. This is in contrast with the gluing construction in \cite[Section 5]{Abouzaid_Blumberg_2024}, where the concatenation is carried out for bimodules enriched in $\uds{\bf dOrb}$ by thickening up the boundary strata. In our applications, all ${\mb X} \circ_i {\mb Y}$ arise from geometric context, and the notations in Definition \ref{defn_general_concatenation} already ensure the desired chain-level statements.
\end{rem}

\subsection{Weak flow categories, bimodules, and homotopies}\label{weak_flow_category}

In certain cases it is convenient to use the notion of direct products of flow categories, which satisfy a weaker version of Definition \ref{defn_flow_category}. We make the following definitions in order to have a simple formulation of certain structures on multimodules involved in the equivariant pair-of-pants product. 

\begin{defn}\label{defn_weak_flow_category}
A {\bf weak flow category} enriched in $\uds{\bf C}$, denoted by ${\mb F}$, consists of a countable partially ordered set of objects ${\rm Ob}{\mb F}$ and a morphism space $M_{pq}^{\mb F}$ as an object of $\uds{\bf C}$ for each pair $p \leq q$, and composition maps $\iota_{prq}^{\mb F}: M_{pr}^{\mb F}\times M_{rq}^{\mb F} \to M_{pq}^{\mb F}$ (whose underlying poset map is not necessarily homogeneous of degree $1$) which satisfy (1)---(5) of Definition \ref{defn_flow_category}.
\end{defn}

\begin{example}
Let ${\mb F}_1, \ldots, {\mb F}_m$ be (weak) flow categories. Define ${\mb F}_1\times \cdots \times {\mb F}_m$ as a weak flow category with 
\beqn
{\rm Ob} {\mb F}_1\times \cdots \times {\mb F}_m = \prod_{j=1}^m {\rm Ob}{\mb F}_j,
\eeqn
morphism spaces
\beqn
M_{(p_1 \cdots p_m)(q_1 \cdots q_m)}^{{\mb F}_1\times \cdots \times {\mb F}_m} = \prod_{j=1}^m M_{p_j q_j}^{{\mb F}_j},
\eeqn
and product structural maps. 
\end{example}

We can also define the similar notion of weak multimodules. 

\begin{defn}\label{defn_weak_multimodule}
Let ${\mb F}_1, \ldots, {\mb F}_m; {\mb F}'$ be weak flow categories enriched in $\uds{\bf C}$. A {\bf weak multimodule} ${\mb X}$ over $({\mb F}_1, \ldots, {\mb F}_m; {\mb F}')$ consists of the same type of objects $M_{p_1 \cdots p_m; p'}^{\mb X}$ and structural maps (whose underlying poset map is not necessarily homogeneous of degree $1$) satisfying conditions (1)--(4) of Definition \ref{defn_multimodule}. 
\end{defn}

The following lemma is readily checked.

\begin{lemma}\label{lemma_multimodule_weak_bimodule}
Let ${\mb F}_1, \ldots, {\mb F}_m, {\mb F}'$ be flow categories enriched in $\uds{\bf C}$. Let ${\mb X}$ be a multimodule over $({\mb F}_1, \ldots, {\mb F}_m; {\mb F}')$. Then the objects $M_{p_1 \cdots p_m; p'}^{\mb X}$ and structural maps in ${\mb X}$ define a weak bimodule over $({\mb F}_1 \times \cdots \times {\mb F}_m; {\mb F}')$. 
\end{lemma}

The following notion is introduced merely to simplify the formulation of multimodule concatenations. It is not the same as the more geometric bimodule obtained via continuation maps.

\begin{defn}
Let ${\mb F}$ be a (weak) flow category. The {\bf naive identity bimodule} of ${\mb F}$, denoted by ${\mb I}^{\mb{FF}}$, is the weak bimodule over $({\mb F}; {\mb F})$ with
\beqn
M_{p;q}^{{\mb I}^{\mb{FF}}} = M_{pq}^{\mb F}
\eeqn
and structural maps induced from composition maps of ${\mb F}$.
\end{defn}

The case of homotopies is similar. We omit the definition and the corresponding lemma.

\subsubsection{Multimodule concatenations as weak bimodule concatenations}

Now we consider the situation of multimodule concatenations.  Let ${\mb X}$ be a multimodule over $({\mb F}_1, \ldots, {\mb F}_m; {\mb G}_i)$ and ${\mb Y}$ be a multimodule over $({\mb G}_1, \ldots, {\mb G}_n; {\mb G}')$. The concatenation ${\mb X} \circ_i {\mb Y}$ at ${\mb G}_i$ is defined previously as a multimodule over $({\mb G}_1, \cdots, {\mb G}_{i-1}, {\mb F}_1, \ldots, {\mb F}_m, {\mb G}_{i+1}, \ldots, {\mb G}_n; {\mb G}')$ (together with other data). As multimodules can be treated as weak bimodules (Lemma \ref{lemma_multimodule_weak_bimodule}), we can interpret the concatenation as a concatenation of a pair of weak multimodules.

First, notice that the notion of concatenations can be extended from multimodules to weak multimodules. Then we consider the product of weak bimodules
\beqn
{\mb I}^{{\mb G}_1{\mb G}_1} \times \cdots \times {\mb I}^{{\mb G}_{i-1} {\mb G}_{i-1}} \times {\mb X} \times {\mb I}^{\mb{G}_{i+1} \mb{G}_{i+1}} \times \cdots \times \mb{I}^{\mb{G}_n \mb{G}_n}
\eeqn
as a weak bimodule over $({\mb G}_1 \times \cdots \times {\mb G}_{i-1} \times {\mb F}_1 \times \cdots \times {\mb F}_m \times {\mb G}_{i+1} \times \cdots \times {\mb G}_n; {\mb G}_1 \times \cdots \times {\mb G}_n)$. Then the multimodule concatenation ${\mb X} \circ_i {\mb Y}$ is equivalent to the weak bimodule concatenation 
\beqn
({\mb I}^{{\mb G}_1{\mb G}_1} \times \cdots \times {\mb I}^{{\mb G}_{i-1} {\mb G}_{i-1}} \times {\mb X} \times {\mb I}^{\mb{G}_{i+1} \mb{G}_{i+1}} \times \cdots \times \mb{I}^{\mb{G}_n \mb{G}_n}) \circ {\mb Y}.
\eeqn
The purpose of introducing this formal notion is to simplify the definition of equivariance of multimodules when the symmetry group permutes the source flow categories.

\section{Algebraic Structures}\label{section6}\label{section_algebraic}

In this section we describe how to obtain algebraic structures such as chain complexes, chain maps, and chain homotopies from, respectively, flow categories, multimodules, and homotopies. When additional structures exist, these chain-level objects also carry corresponding structures such as equivariance and module structures over Novikov rings.

\subsection{Unobstructedness, local finiteness, grading, and orientations}

We start by considering the most basic properties of flow categories, multimodules, and homotopies that lead to chain-level structures. Notice that these structures will only produce chain complex etc. as abelian groups. To obtain objects as modules over various Novikov rings, one needs certain extra symmetries on the flow categories etc.

\subsubsection{Unobstructedness}

In geometric applications, the flow categories, multimodules, and homotopies always satisfy an additional property which does not have a general formulation. Here we specify this property.

\begin{defn}\label{defn_unobstructed}
Let ${\mb F}$ be a flow category enriched in $\uds{\bf Top}$. It is said to be {\bf unobstructed}\footnote{This terminology comes from Lagrangian Floer theory. If we formulate the Lagrangian Floer flow category, then in general the true boundary of moduli spaces contain configurations with disk bubbles, potentially causing obstruction to define a chain complex.} if for any triple $p < r < q$ of objects, the structural morphism
\beqn
\iota_{prq}^{\mb F}: M_{pr}^{\mb F} \times M_{rq}^{\mb F} \to M_{pq}^{\mb F}
\eeqn
is a homeomorphism onto $\partial^{prq} M_{pq}^{\mb F}$. Hence the union of these strata is the true boundary of $M_{pq}^{\mb F}$.

If $\uds{{\bf C}} \to \uds{\bf Top}$ is a functor of regular stratification categories and ${\mb F}$ is a flow category enriched in $\uds{\bf C}$, ${\mb F}$ is called {\bf unobstructed} if the induced flow category enriched in $\uds{\bf Top}$ is unobstructed.
\end{defn}

Similarly, we can define the notion of unobstructedness for multimodules and homotopies enriched in $\uds{\bf Top}$. We omit the details. In applications within this paper, all geometric flow categories and their lifts are unobstructed.

\subsubsection{Convergence}

We need certain properties guaranteeing the convergence of the chain-level objects in certain nonarchimedean completions. To incorporate the equivariant case, we allow a second functional which is bounded from one side to give a different completion.

\begin{defn}\label{defn_local_finite_flow}
A {\bf filtration} on a (weak) flow category ${\mb F}$ is a monotonic function
\beqn
{\mc A}: {\rm Ob}{\mb F} \to {\mb R}^k\ (k\geq 1)
\eeqn
where ${\mb R}^k$ has the product partial order. The filtration is called {\bf locally finite} if for each $p \in {\rm Ob}{\mb F}$ and any $c>0$,
    \beqn
\# \Big\{ q \in {\rm Ob}{\mb F}\ |\ p< q,\  {\mc A}_\alpha (p) - {\mc A}_\alpha (q) < c,\ \alpha = 1, \ldots, k \Big\} < \infty.
\eeqn
A (weak) flow category equipped with a locally finite filtration is called a {\bf locally finite} (weak) flow category.
\end{defn}

\begin{lemma}
If $({\mb F}_i, {\mc A}_i)$, $i = 1, \ldots, m$ are locally finite filtered flow categories, then the product $({\mb F}_1 \times \cdots \times {\mb F}_m, {\mc A}_1 \oplus \cdots \oplus {\mc A}_m)$ is a locally finite weak flow category.
\end{lemma}

\begin{proof}
Straightforward.
\end{proof}

\begin{rem}
In a typical Floer theoretic setup, the action ${\mc A}$ is ${\mb R}$-valued. In the equivariant case considered in Part 4 of this paper, the action will have two components. 
\end{rem}

Now we consider the local finiteness condition for multimodules. The condition we need is not only to guarantee a well-defined convergent chain map, but also that multimodule concatenations preserve such conditions. %

\begin{defn}\label{defn_local_finite_bimodule}
Let $({\mb F}, {\mc A})$ and $({\mb F}', {\mc A}')$ be $\uds{\bf C}$-enriched locally finite (weak) flow categories. Suppose ${\mc A}$ is ${\mb R}^k$-valued and ${\mc A}'$ is ${\mb R}^{k'}$-valued.  A (weak) bimodule ${\mb B}$ over $({\mb F}; {\mb F}')$ is called {\bf locally finite} if the following conditions are satisfied.

\begin{enumerate}

\item For any $p \in {\rm Ob}{\mb F}$, $c' \in {\mb R}$,
\beqn
 \# \Big\{ p' \in {\rm Ob}{\mb F}'\ |\  M_{p;p'}^{\mb B}  \neq \emptyset,\ \min_{1 \leq \alpha \leq k' } {\mc A}_\alpha' (p') \geq c' \Big\} < \infty.
\eeqn

\item For all $c \in {\mb R}$
\beqn
f^{\mb B}(c):= \sup \Big\{ \max_{1 \leq \alpha \leq k'} {\mc A}_\alpha'(p')\ |\ \max_{1 \leq \alpha \leq k} {\mc A}_\alpha(p) \leq c\ {\rm and}\ M_{p;p'}^{\mb B} \neq \emptyset \Big\} < +\infty,
\eeqn
and
\begin{multline*}
\lim_{\rho \to -\infty} h^{\mb B}(c, \rho)\\
:=\lim_{\rho \to -\infty} \sup \Big\{ \min_{1 \leq \alpha \leq k'} {\mc A}_\alpha'(p')\ |\ \max_{1 \leq \alpha \leq k}{\mc A}_\alpha(p) \leq c,\ \min_{1 \leq \alpha \leq k} {\mc A}_\alpha(p) \leq \rho,\ {\rm and}\ M_{p;p'}^{\mb B} \neq \emptyset \Big\}  = -\infty.
\end{multline*}

\end{enumerate}
\end{defn}

The case of multimodules can be obtained as a special case of weak bimodules.

\begin{defn}\label{defn_Novikov_multimodule}
Let $({\mb F}_i, {\mc A}_i)$, $i=1, \ldots, m$ and $({\mb F}', {\mc A}')$ be locally finite flow categories enriched in $\uds{\bf C}$. A {\bf locally finite multimodule} over $({\mb F}_1, \ldots, {\mb F}_m; {\mb F}')$ is a multimodule ${\mb X}$ which is a locally finite weak bimodule over $({\mb F}_1\times \cdots \times {\mb F}_m; {\mb F}')$. 
\end{defn}

We verify that concatenations of locally finite multimodules are locally finite.

\begin{lemma}\label{lemma_concatenation_Novikov}
Let $\uds{\bf C}$ be either $\uds{\bf Top}$ or $\uds{\bf dOrb}$. Let $({\mb F}_i, {\mc A}_i)$, $i = 1, \ldots, m$, $({\mb G}_j, {\mc B}_j)$, $j = 1, \ldots, n$, $({\mb G}', {\mc B}')$, be locally finite flow categories. Let ${\mb X}$ be a locally finite multimodule over $({\mb F}_1, \ldots, {\mb F}_m; {\mb G}_i)$ and ${\mb Y}$ be a locally finite multimodule over $({\mb G}_1, \ldots, {\mb G}_n; {\mb G}')$. Then the concatenation ${\mb X} \circ_i {\mb Y}$ is a locally finite multimodule. 
\end{lemma}

\begin{proof}
We verify conditions of Definition \ref{defn_local_finite_bimodule} one by one. Suppose (1) of Definition \ref{defn_local_finite_bimodule} does not hold for the concatenation ${\mb X} \circ_i {\mb Y}$. Then there exist 
\beqn
p_1 \in {\rm Ob}{\mb F}_1, \ldots, p_m \in {\rm Ob}{\mb F}_m, r_1 \in {\rm Ob}{\mb G}_1, \ldots, r_{i-1} \in {\rm Ob}{\mb G}_{i-1},\ r_{i+1} \in {\rm Ob}{\mb G}_{i+1}, \ldots, r_n \in {\rm Ob}{\mb G}_n,
\eeqn
$c' \in {\mb R}$ and infinitely many distinct $r_\nu'\in {\rm Ob}{\mb G}'$ such that 
\begin{align}\label{Novikov_concatenation_1}
&\ M_{r_1 \cdots r_{i-1} p_1 \cdots p_m r_{i+1} \cdots r_n; r'}^{{\mb X} \circ {\mb Y} } \neq \emptyset,\ &\ \min_{1 \leq \beta \leq l'} {\mc B}_\beta'(r_\nu') \geq c'
\end{align}
It follows that there exists a sequence $r_{i, \nu} \in {\rm Ob}{\mb G}_i$ such that 
\begin{align*}
&\ M_{p_1 \cdots p_m; r_{i,\nu}}^{\mb X} \neq \emptyset,\ &\ M_{r_1 \cdots r_{i-1} r_{i, \nu} r_{i+1} \cdots r_n; r_\nu'}^{\mb Y} \neq \emptyset.
\end{align*}
By choosing a subsequence, we may assume the terms of the sequence $r_{i, \nu}$ are either all identical, or all distinct. If all $r_{i, \nu}$ are identical, then it contradicts the local finiteness of ${\mb Y}$. Hence all $r_{i, \nu}$ are distinct. Then by the local finiteness of ${\mb X}$, for all $\alpha\in \{1, \ldots, l_i\}$, the sequence ${\mc B}_{i, \alpha}(r_{i, \nu})$ is bounded from above; moreover, by choosing a subsequence, for some $\alpha$ the sequence ${\mc B}_{i,\alpha}(r_{i,\nu}) \to -\infty$. Then by the local finiteness of ${\mb Y}$, one has
\beqn
\lim_{\nu \to \infty} \min_{1 \leq \beta \leq l'} {\mc B}_\beta'(r_\nu') = -\infty.
\eeqn
This contradicts \eqref{Novikov_concatenation_1}. Therefore, (1) holds for the concatenation ${\mb X} \circ_i {\mb Y}$. 

Now we verify (2) for the concatenation. 
Choose $c \in {\mb R}$. We need to show that 
\beqn
f^{{\mb X} \circ_i {\mb Y}}(c) < +\infty
\eeqn
where $f^{{\mb X} \circ_i {\mb Y}}$ is defined in Definition \ref{defn_local_finite_bimodule}. Consider any collection of objects $p_1, \ldots, p_m, r_1, \ldots, r_{i-1}, r_{i+1}, \ldots, r_n, r'$ such that
\beqn
\max_{1 \leq \alpha \leq k_i} {\mc A}_{i, \alpha}(p_i) \leq c,\ \max_{1 \leq \beta \leq l_j} {\mc B}_{j, \beta}(r_j) \leq  c
\eeqn
and 
\beqn
M_{r_1 \cdots r_{i-1} p_1 \cdots p_m r_{i+1} \cdots r_n; r'}^{{\mb X} \circ_i {\mb Y}} \neq \emptyset.
\eeqn
Then by the definition of concatenation, there exists $p' \in {\rm Ob}{\mb F}'$ such that 
\begin{align*}
&\ M_{p_1 \cdots p_m; r_i}^{{\mb X}} \neq \emptyset,\ &\ M_{r_1 \cdots r_n; r'}^{{\mb Y}} \neq \emptyset.
\end{align*}
Then 
\beqn
\max_{1 \leq \beta \leq l_i} {\mc B}_{i, \beta}(r_i) \leq f^{\mb X}(c).
\eeqn
Denote $c' = \max \{ c, f^{\mb X}(c)\}$. Then by the assumption, one has 
\beqn
\max_{1 \leq \beta \leq l'} {\mc B}_\beta'(r') \leq f^{\mb Y}(c').
\eeqn
This proves that 
\beqn
f^{{\mb X} \circ_i {\mb Y}}(c) \leq f^{{\mb Y}}(c') < +\infty.
\eeqn
Moreover, suppose 
\beqn
\rho  = \min \Big\{ \min_{1 \leq \alpha \leq k_j} {\mc A}_{j, \alpha}(p_j)\ |\ 1 \leq j \leq m \Big\}.
\eeqn
Then 
\beqn
\min_{1 \leq \beta \leq l_i} {\mc B}_{i, \beta} (r_i) \leq h^{{\mb X}}(c, \rho)
\eeqn
Denote
\beqn
\rho':= \min \Big\{ h^{{\mb X}}(c, \rho), \min_{1 \leq \beta \leq l_j} {\mc B}_{j, \beta}(r_j)\ |\ j \neq i \Big\}.
\eeqn
Then 
\beqn
\min_{1 \leq \beta \leq l'} {\mc B}_\beta'(r') \leq h^{\mb Y}(c', \rho')
\eeqn
which goes to $-\infty$ as $\rho' \to -\infty$. Therefore (2) of Definition \ref{defn_local_finite_bimodule} holds for the concatenation ${\mb X} \circ_i {\mb Y}$.
\end{proof}

Lastly we state the local finiteness of multimodule homotopies. 

\begin{defn}\label{defn_Novikov_homotopy}
Let $({\mb F}, {\mc A})$ and $({\mb F}', {\mc A}')$ be two locally finite weak flow categories. Suppose ${\mc A}$ is ${\mb R}^k$-valued and ${\mc A}'$ is ${\mb R}^{k'}$-valued. Suppose ${\mb B}_0$ and ${\mb B}_1$ are two locally finite weak bimodules over $({\mb F}; {\mb F}')$ (Definition \ref{defn_local_finite_bimodule}). Then a homotopy ${\mb H}$ from ${\mb B}_0$ to ${\mb B}_1$ is called a {\bf locally finite homotopy} if the following are satisfied.
\begin{enumerate}

\item For any $p \in {\rm Ob}{\mb F}$ and $c' \in {\mb R}$,
\beqn
\# \Big\{ p' \in {\rm Ob}{\mb F}'\ |\ M_{p; p'}^{\mb H} \neq \emptyset,\  \min_{1 \leq \alpha \leq k'} {\mc A}_\alpha'(p') \geq c' \Big\} < \infty.
\eeqn

\item For all $c \in {\mb R}$
\beqn
 \sup \Big\{ \max_{1 \leq \alpha \leq k'} {\mc A}_\alpha'(p')\ |\ \max_{1 \leq \alpha \leq k} {\mc A}_\alpha(p) \leq c\ {\rm and}\ M_{p;p'}^{\mb H} \neq \emptyset \Big\} < +\infty,
\eeqn
and
\beqn
\lim_{\rho \to -\infty} \sup \Big\{ \min_{1 \leq \alpha \leq k'} {\mc A}_\alpha'(p')\ |\ \max_{1 \leq \alpha \leq k}{\mc A}_\alpha(p) \leq c,\ \min_{1 \leq \alpha \leq k} {\mc A}_\alpha(p) \leq \rho,\ {\rm and}\ M_{p;p'}^{\mb H} \neq \emptyset \Big\}  = -\infty.
\eeqn

\end{enumerate}
\end{defn}

\begin{defn}
Let $({\mb F}_i, {\mc A}_i)$, $i=1, \ldots, m$ and $({\mb F}', {\mc A}')$ be locally finite flow categories enriched in $\uds{\bf C}$. Let ${\mb X}_0$ and ${\mb X}_1$ be locally finite multimodules over $({\mb F}_1, \cdots, {\mb F}_m; {\mb F}')$. A {\bf locally finite homotopy} from ${\mb X}_0$ to ${\mb X}_1$ is a locally finite weak homotopy between the corresponding weak bimodules.  
\end{defn}

\subsubsection{Grading}

Now we consider gradings on the flow categories, multimodules, and homotopies. Let $\uds{\bf C}$ be a regular stratification category which is either $\uds{\bf Kur}$, $\uds{\bf dOrb}$, or $\pman$. Define the notion of {\bf virtual dimension} by
\beqn
{\rm dim}^{\rm vir} C:= \left\{ \begin{array}{ll} {\rm dim} U - {\rm dim} G - {\rm rank} E,\ &\  C = (G, U, E, S) \in {\rm Ob} \uds{\bf Kur},\\
{\rm dim} {\mc U} - {\rm rank} {\mc E},\ &\ C = ({\mc U}, {\mc E}, {\mc S}) \in {\rm Ob} \uds{\bf dOrb},\\
{\rm dim} C,\ &\ C \in {\rm Ob} \pman
\end{array}\right. 
\eeqn
We allow different components of $C$ to have different virtual dimensions. 

Below we define gradings. In most of this paper, we will consider absolute, namely, ${\mb Z}$-gradings; in certain situations we will reduce to cyclic gradings.

\begin{defn}\label{defn_grading}
Let $N$ be a nonnegative integer. All flow categories and multimodules are assumed to be %
enriched in either $\uds{\bf Kur}$, $\uds{\bf dOrb}$, or $\pman$.
\begin{enumerate}

\item A (cohomological) {\bf ${\mb Z}$-grading} on a flow category ${\mb F}$ is a %
function $p \mapsto |p| \in {\mb Z}$ satisfying 
\beqn
{\rm dim}^{\rm vir} M_{pq}^{\mb F} = |q| - |p| - 1,\ \forall p, q \in {\rm Ob}{\mb F}.
\eeqn

\item Let ${\mb F}_1, \ldots, {\mb F}_m, {\mb F}'$ be flow categories equipped with ${\mb Z}$-gradings. A multimodule ${\mb X}$ over $({\mb F}_1, \ldots, {\mb F}_m; {\mb F}')$ is said to have cohomological degree $r$ if for all $p_i \in {\rm Ob} {\mb F}_i$ and $p' \in {\rm Ob}{\mb F}'$ 
\beqn
{\rm dim}^{\rm vir} M_{p_1\cdots p_m;p'}^{\mb X} = |p'| - |p_1| - \cdots - |p_m| - r.
\eeqn

\end{enumerate}
\end{defn}

\subsubsection{Orientations}\label{subsubsec:orientation}

Defining chain complexes in characteristics other than $2$ requires orientations. We do not include orientations as part of the structure on objects in the categories $\uds{\bf Kur}$, $\uds{\bf dOrb}$, or $\pman$ as it is cumbersome to discuss orientations on the corners. Instead, we regard it as an extra structure on the flow category, multimodule, or homotopy.

{\bf Convention.} Let $M$ be an oriented manifold with boundary. The induced boundary orientation is the one such that the 
\beqn
{\rm n}_{\rm out} \wedge {\rm or}_{\partial M} \in \det TM
\eeqn
coincides (up to positive scalar) with the orientation ${\rm or}_M$ of $M$.

\begin{defn}\label{defn711}
Let $C = ({\mc U}, {\mc E}, {\mc S})$ be an $A$-stratified derived orbifold with $A$ being a normal poset (Definition \ref{defn_normal_poset}).
\begin{enumerate}

\item For each element $\alpha \in A$, an orientation on $\partial^\alpha C$ is an orientation of the virtual vector bundle $T\partial^\alpha {\mc U} - {\mc E}|_{\partial^\alpha V}$, or equivalently an orientation on the determinant line 
\beqn
\det (\partial^\alpha C) := \det (T\partial^\alpha {\mc U}) \otimes \det ({\mc E}^*|_{T\partial^\alpha {\mc U}}).
\eeqn
Notice that for any face $\beta \in \mb{Face}_\alpha$, an orientation on $\partial^\alpha C$ induces an orientation on $\partial^\beta C$ by regarding the collar coordinate to point ``inward.''

\item An orientation on $C$ consists of orientations on all top strata $\partial^\gamma C$ such that for each fake boundary stratum $\beta$, the orientation on $\partial^\beta C$ induced from the (exactly) two adjacent top strata are opposite. 

\item Let $C' = ({\mc U}', {\mc E}', {\mc S}')$ be another $A'$-stratified derived orbifold. An embedding $\iota: C \to C'$ is said to be {\bf orientation-preserving} if the underlying poset map has degree zero (i.e. maps top strata to top strata) and if for each top stratum $\gamma \in A$ (which is mapped to $\gamma'\in A'$), the induced map of determinant lines
\beqn
\det (\partial^\gamma C) \to \det (\partial^{\gamma'} C')
\eeqn
is orientation-preserving.
\end{enumerate}
\end{defn}

On the other hand, an orientation on a stratified pseudomanifold $X$ is defined to consist of orientations on all top strata $\gamma \in A^{\rm max}$ satisfying the similar compatibility condition on fake boundaries.

\begin{defn}\label{defn_coherent_orientation}
Let ${\mb F}$ be a flow category enriched in $\uds{\bf dOrb}$ (resp. $\pman$). A {\bf coherent orientation} on ${\mb F}$ is a collection of orientations on the object $M_{pq}^{\mb F}$ such that the singletons $M_{pp}^{\mb F}$ are positively oriented and such that whenever $p< r < q$, the morphism
\beqn
\iota_{prq}^{\mb F}: M_{pr}^{\mb F} \times M_{rq}^{\mb F} \to \partial^{prq} M_{pq}^{\mb F}
\eeqn
is orientation-preserving up to factor
\beqn
(-1)^{|p|- |r|}.
\eeqn
Here by the definition of flow categories (Definition \ref{defn_flow_category}), the top strata of $\partial^{prq} M_{pq}^{\mb F}$ are true boundaries of $M_{pq}^{\mb F}$, hence have canonically induced orientations.

A flow category enriched in $\uds{\bf dOrb}$ or $\pman$ equipped with a coherent orientation is called {\bf oriented}.
\end{defn}

Similarly we have the notion of orientation for multimodules. Notice the signs appearing in the definition, which are put in order to obtain the correct chain maps.

\begin{defn}\label{defn_orientation_multimodule}
Let ${\mb F}_1, \ldots, {\mb F}_m, {\mb F}'$ be ${\mb Z}$-graded oriented flow categories enriched in $\uds{\bf dOrb}$ (resp. $\pman$). Let ${\mb X}$ be a degree zero multimodule over $({\mb F}_1, \ldots, {\mb F}_m; {\mb F}')$. Suppose all flow categories are equipped with coherent orientations. An {\bf coherent orientation} on ${\mb X}$ consists of orientations on the charts $M_{pq}^{{\mb X}}$ such that \begin{enumerate}

\item For each $i = 1, \ldots, m$, the morphism 
\beqn
M_{p_i q_i}^{{\mb F}_i} \times M_{p_1 \cdots p_{i-1} q_i p_{i+1}\cdots p_m; p'}^{{\mb X}} \to \partial M_{p_1 \cdots p_m; p'}^{{\mb X}}
\eeqn
is orientation-preserving up to the sign 
\beqn
(-1)^{ (\sum_{j<i} |p_j |) (|p_i| - |q_i|)}.
\eeqn
In particular, for $i=1$, the boundary orientation coincides with the product orientation.

\item The morphism 
\beqn
M_{p_1 \cdots p_m; q'}^{\mb X} \times M_{q' p'}^{{\mb F}'} \to \partial M_{p_1 \cdots p_m; p'}^{{\mb X}}
\eeqn
is orientation-preserving up to the sign
\beqn
(-1)^{{\rm dim} M_{p_1 \cdots p_m; q'}^{\mb X} + 1}.
\eeqn
 
\end{enumerate}

A multimodule over oriented flow categories equipped with coherent orientations is called an {\bf oriented multimodule}.
\end{defn}

We consider the orientation for concatenations. Let ${\mb X}$ be a multimodule over $({\mb F}_1, \ldots, {\mb F}_m; {\mb G}_i)$ and ${\mb Y}$ be a multimodule over $({\mb G}_1, \ldots, {\mb G}_n; {\mb G}')$. Recall that for each tuple $p_i \in {\rm Ob}{\mb F}_i$, $r_j \in {\rm Ob}{\mb G}_j$ and $r' \in {\rm Ob}{\mb G}'$, the top strata of $M_{r_1 \cdots r_{i-1} p_1 \cdots p_m  r_i \cdots r_n; r'}^{{\mb X} \circ_i {\mb Y}}$ is the disjoint union of the top strata of 
\beqn
\bigsqcup_{r_i} M_{p_1 \cdots p_m; r_i}^{\mb X} \times M_{r_1 \cdots r_n; r'}^{\mb Y}.
\eeqn
We define each of these products to be oriented by the product orientation times the sign
\beqn
(-1)^{\left( \sum_{j<i} |r_j|\right){\rm dim} M_{p_1 \cdots p_m; r_i}^{\mb X}} 
\eeqn
(if the concatenation happens at the $i+1$-th input). The moduli spaces of the concatenation necessarily have fake boundaries and one can check easily that the above defined orientation on the interior induces an orientation on the orientations in the sense of Definition \ref{defn711}. We then need to check that these orientations on individual moduli spaces define a coherent orientation of the multimodule.

\begin{lemma}\label{lemma_concatenation_orientation}
The orientation induced on ${\mb X}\circ_i {\mb Y}$ is coherent.
\end{lemma}

\begin{proof}
Fix a moduli space
\beqn
M_{r_1 \cdots r_{i-1} p_1 \cdots p_m r_{i+1} \cdots r_n; r'}^{{\mb X}\circ_i {\mb Y}}.
\eeqn
We first need to verify that the orientations on two open strata of this moduli space induce opposite orientations on the shared fake boundary stratum. Choose objects $r_i < s_i$ of ${\mb F}_i$. The two open strata are
\begin{align}
&\ M_{p_1 \cdots p_m; r_i}^{\mb X}\times M_{r_1 \cdots r_n; r'}^{\mb Y},\ &\ M_{p_1 \cdots p_m; s_i}^{\mb X}\times M_{r_1 \cdots r_{i-1} s_i r_{i+1} \cdots r_n; r'}^{\mb Y}.
\end{align}
The shared fake boundary stratum is 
\beqn
M_{p_1 \cdots p_m; r_i}^{\mb X} \times M_{r_i s_i}^{{\mb F}_i} \times M_{r_1 \cdots r_{i-1} s_i r_{i+1} \cdots r_m; r'}^{\mb Y}.
\eeqn
Use the product orientation as a reference. The orientation on this fake boundary induced from the first differs from the product orientation by 
\beqn
\left( \sum_{j<i}|r_j|\right) {\rm dim} M_{p_1 \cdots p_m; r_i}^{\mb X} + {\rm dim} M_{p_1 \cdots p_m; r_i}^{\mb X} + \left( \sum_{j<i} |r_j| \right) (|r_i| - |s_i|).
\eeqn
The second orientation on the fake boundary differs by a sign
\beqn
\left( \sum_{j<i} |r_j| \right)  {\rm dim} M_{p_1 \cdots p_m; s_i}^{\mb X} + {\rm dim} M_{p_1 \cdots p_m; r_i}^{\mb X} + 1.
\eeqn
One can see that they are indeed opposite.

Then one has to verify conditions listed in Definition \ref{defn_orientation_multimodule} for the concatenation. For breakings at some $j$-th input $s_j \in {\rm Ob} {\mb G}_j$ for $j<i$ and breakings at some $i+l-1$-th input $q_l \in {\rm Ob} {\mb F}_l$, the verification is straightforward. 

\end{proof}

\begin{defn}\label{defn_orientation_homotopy}
Let ${\mb F}_1, \ldots, {\mb F}_m, {\mb F}'$ be ${\mb Z}$-graded oriented flow categories enriched in $\uds{\bf dOrb}$ (resp. $\pman$), ${\mb X}_0, {\mb X}_1$ be oriented flow multimodules over $( {\mb F}_1, \ldots, {\mb F}_m; {\mb F}')$ of degree zero (see Definition \ref{defn_grading}). Let ${\mb H}$ be a ${\mb Z}$-graded homotopy from $ {\mb X}_0$ to $ {\mb X}_1$. A {\bf coherent orientation} on ${\mb H}$ consists of orientations on the charts $M_{pp'}^{{\mb H}}$ such that 
\begin{enumerate}

\item The embedding 
\beqn
(- M_{pp'}^{{\mb B}_0}) \sqcup M_{pp'}^{{\mb B}_1} \to \partial M_{pp'}^{{\mb H}}
\eeqn
is orientation-preserving.

\item The embedding
\beqn
M_{p_i q_i}^{{\mb F}_i} \times M_{p_1 \cdots p_{i-1} q_i p_{i+1} \cdots p_m; p'}^{{\mb H}} \to \partial M_{p_1 \cdots p_m; p'}^{{\mb H}}
\eeqn
is orientation-preserving up to the sign 
\beqn
(-1)^{ \left( \sum_{j<i} |p_j| \right) (|p_i|-|q_i|)}.
\eeqn

\item The embedding
\beqn
M_{p_1 \cdots p_m; q'}^{{\mb H}} \times M_{q'p'}^{{\mb F}'} \to \partial M_{p_1 \cdots p_m; p'}^{{\mb H}}
\eeqn
is orientation-preserving up to the sign
\beqn
(-1)^{{\rm dim} M_{p_1 \cdots p_m; q'}^{\mb H} + 1}.
\eeqn
\end{enumerate}
\end{defn}

\subsection{Algebraic structures without symmetry}\label{subsection_abstract_chain_complex}

We define chain complexes, chain maps, and chain homotopies associated to flow categories, bimodules, and homotopies enriched in $\pman$. We only impose the Novikov convergence condition but not the symmetries. Hence the resulting objects are only modules over ${\mb Z}$ but not the Novikov rings.

\subsubsection{The chain group}

We only state the construction for the ${\mb Z}$-graded case. The cyclically graded case can be produced in a similar fashion. 

To a ${\mb Z}$-graded locally finite weak flow category $({\mb F}, {\mc A})$ enriched in $\uds{\bf C}$ (see Definition \ref{defn_local_finite_flow}), one can associate a chain group
\beqn
C_* ( {\mb F}):= \left\{ \sum_{i=1}^\infty a_i p_i\ \left|\ \begin{array}{l} a_i \in {\mb Z},\ p_i \in {\rm Ob}{\mb F},\ \displaystyle \sup_i {\mc A}_\alpha (p_i) < +\infty \ \forall \alpha = 1, \ldots, m,\\
\forall c>0,\ \# \big\{ i\ |\ {\mc A}_\alpha (p_i) > -c\ \forall \alpha = 1, \ldots, k   \big\} < \infty \end{array}\right. \right\}.
\eeqn
It is not necessarily finitely generated over ${\mb Z}$. We call it the {\bf Floer chain group} associated to $({\mb F}, {\mc A})$. It carries a natural valuation
\beqn
{\mf v} = ({\mf v}_1, \ldots, {\mf v}_k): C_*({\mb F}) \setminus \{0\} \to {\mb R}^k,\ {\rm where}\ {\mf v}_\alpha \left(\sum_{i=1}^\infty a_i p_i \right) = \sup_i {\mc A}_\alpha(p_i).
\eeqn
For each $k \in {\mb Z}$, denote
\beqn
C_k ({\mb F}):= \left\{ \sum_{i=1}^\infty a_i p_i \in C_*({\mb F})\ |\ {\rm deg}(p_i) = k \right\}.\footnote{$C_*({\mb F})$ is not the algebraic direct sum of $C_k({\mb F})$, but rather a completion of the direct sum.
}
\eeqn

Suppose $({\mb F}', {\mc A}')$ and $({\mb F}'', {\mc A}'')$ are two weak Novikov flow categories. The Floer chain group $C_*({\mb F}' \times {\mb F}'')$ is a kind of completion of the tensor product chain groups of the factors. Indeed, on the algebraic (finite) tensor product
\beqn
C_*({\mb F}') \underset{\mb Z}{\otimes} C_*({\mb F}'')
\eeqn
there is the valuation ${\mf v}$  induced by ${\mf v}'$ and ${\mf v}''$. Then one can form the completion 
\beqn
C_*({\mb F}') \widehat{ \underset{\mb Z}{\otimes}} C_*({\mb F}'').
\eeqn

\subsubsection{The differential}

Now, we assume that the flow category ${\mb F}$ is enriched in $\pman$ and it is oriented. We define a linear map
\beqn
d: C_k( {\mb F}) \to C_{k-1}( {\mb F}),\ k \in {\mb Z}
\eeqn
in the following way. Note that the morphism space $ M_{pq}^{\mb F}$ is an oriented and compact pseudomanifold of dimension
\beqn
{\rm dim} M_{pq}^{\mb F} \equiv {\rm deg}(p) - {\rm deg} (q) - 1.
\eeqn
When the dimension is zero, $M_{pq}^{{\mb F}}$ consists of finitely many oriented points, which defines a count $n_{pq}^{\mb F} \in {\mb Z}$; in other cases, define $n_{pq}^{\mb F} = 0$. Define 
\beqn
d \left( \sum_{i=1}^\infty a_i p_i \right) = \sum_{i=1}^\infty \sum_{q\in {\rm Ob}{\mb F}} a_i n_{p_i q}^{\mb F} q.
\eeqn
The local finiteness property stated in Definition \ref{defn_local_finite_flow} implies that $d$ is well-defined.

\begin{prop}\label{prop_differential}
If ${\mb F}$ is unobstructed, then $d^2 = 0$. Therefore, $(C_*({\mb F}), d)$ gives rise to a ${\mb Z}$-graded chain complex of free abelian groups, called the {\bf complex} associated to ${\mb F}$.
\end{prop}

\begin{proof}
To show $d^2 = 0$, consider pairs $p, q$ with ${\rm deg} (p) - {\rm deg}(q) = 2$. By Definition \ref{defn_coherent_orientation}, the coefficients of $d^2(p)$ in $q$ is the oriented count of the union of products
\beqn
\bigcup_{p < r < q} M_{pr}^{\mb F} \times M_{rq}^{\mb F}.
\eeqn
Consider the 1-dimensional component of the pseudomanifold $M_{pq}^{\mb F}$. The boundary points only appear on the true boundary strata of $M_{pq}^{\mb F}$. By the definition of unobstructedness (Definition \ref{defn_unobstructed}), the true boundary strata are all the top strata of the products $M_{pr}\times M_{rq}$. By the coherence of the orientation (see Definition \ref{defn_coherent_orientation}), these counts cancel. Hence $d^2 = 0$. 
\end{proof}

If $C({\mb F}')$ and $C({\mb F}'')$ are complexes associated to flow categories ${\mb F}'$ and ${\mb F}''$, then the tensor product complex $C({\mb F}') \underset{\mb Z}{\otimes} C({\mb F}'')$ is defined. The following can be checked readily.

\begin{lemma}
The tensor product differential extends to the completion $\displaystyle 
C ({\mb F}') \widehat{\underset{\mb Z}{\otimes}} C ({\mb F}'')$.
\end{lemma}

\subsubsection{Chain maps and homotopies}

Now consider an oriented, unobstructed, ${\mb Z}$-graded Novikov multimodule ${\mb X}$ enriched in $\pman$ over $({\mb F}_1, \ldots, {\mb F}_m; {\mb F}')$. Let 
\beqn
m_{p_1 \cdots p_m; p'}^{\mb X} \in {\mb Z}
\eeqn
be the count of $M_{p_1 \cdots p_m; p'}^{{\mb X}}$ when the dimension is zero and $0$ otherwise. Define (formally) a linear map 
\beqn
\tilde \Phi^{\mb X}: C ({\mb F}_1) {\underset{{\mb Z}}{\otimes}} \cdots {\underset{{\mb Z}}{\otimes}} C ({\mb F}_m) \to C ({\mb F}')
\eeqn
by linearly extending (over ${\mb Z}$)
\beq\label{chain_map_sign}
\tilde \Phi^{\mb X} (p_1 \otimes \cdots \otimes p_m):= \sum_{p'\in {\rm Ob}{\mb F}'}  m_{p_1 \cdots p_m; p'}^{\mb X} p'.
\eeq
The Novikov property (Definition \ref{defn_local_finite_bimodule}) implies that $\tilde \Phi^{\mb X}$ is well-defined and extends to the completion.

\begin{prop}\label{prop_chain_map}
$\tilde \Phi^{\mb X}$ extends to the completion $\displaystyle C ({\mb F}_1) \widehat{\underset{{\mb Z}}{\otimes}} \cdots \widehat{\underset{{\mb Z}}{\otimes}} C ({\mb F}_m)$ and is a ${\mb Z}$-linear chain map of degree ${\rm index}({\mb X})$.
\end{prop}

\begin{proof}
The grading follows from the definition (see (2) of Definition \ref{defn_grading}). To prove that $\tilde \Phi^{\mb X}$ is a chain map, one needs to look at the orientation. By definition, the differential of the tensor product is determined by 
\beqn
d(p_1 \otimes \cdots \otimes p_m) = \sum_{i=1}^m (-1)^{ \sum_{j<i} |p_j|}  p_1 \otimes \cdots \otimes p_{i-1} \otimes d^{{\mb F}_i}(p_i) \otimes p_{i+1} \otimes \cdots \otimes p_m.
\eeqn
Look at one-dimensional components of the moduli space $M_{p_1 \cdots p_m; p'}^{\mb X}$ (see Definition \ref{defn_orientation_multimodule}), whose oriented boundary is 
\beqn
\bigsqcup_{i=1}^m \Big( (-1)^{(\sum_{j<i}|p_j|)(|p_i|-|q_i)|} M_{p_i q_i}^{{\mb F}_i} \times M_{p_1 \cdots p_{i-1} q_i p_{i+1} \cdots p_m; p'}^{{\mb X}}  \Big) \sqcup \Big((-1)^{{\rm dim} M_{p_1 \cdots p_m; q'}^{\mb X} + 1} M_{p_1 \cdots p_m; q'}^{{\mb X}} \times M_{q'p'}^{{\mb F}'}\Big).
\eeqn
Moreover, each component above is nonempty only if both factors are zero-dimensional. It follows that (also using \eqref{chain_map_sign})
\begin{multline*}
\sum_{q'} \langle \tilde\Phi^{\mb X}(p_1\otimes \cdots \otimes p_m), q'\rangle \langle d^{{\mb F}'}(q'), p' \rangle \\
= \sum_{i=1}^m (-1)^{\sum_{j<i}|p_j|} \sum_{q_i\in {\rm Ob}{\mb F}_i} \langle d^{{\mb F}_i}(p_i), q_i \rangle \langle \tilde\Phi^{{\mb X}}(p_1\otimes \cdots \otimes p_{i-1} \otimes q_i \otimes p_{i+1} \otimes \cdots \otimes p_m), p'\rangle
     \end{multline*}
     which means exactly that $\tilde \Phi^{{\mb X}}$ is a chain map.
\end{proof}

Now we consider the situation of multimodule concatenations (see Subsection \ref{subsection_concatenation}). Let 
\beqn
{\mb F}_1, \ldots, {\mb F}_m, {\mb G}_1, \ldots, {\mb G}_n, {\mb G}'
\eeqn
be flow categories enriched in $\pman$ which are oriented, graded, and Novikov. Let ${\mb X}$ be a multimodule over $({\mb F}_1, \ldots, {\mb F}_m; {\mb G}_i)$ which is oriented, graded, and Novikov. Let ${\mb Y}$ be a multimodule over $({\mb G}_1, \ldots, {\mb G}_n; {\mb G}')$ which is oriented, graded, and Novikov. Then one can consider the concatenation ${\mb X} \circ_i {\mb Y}$ which is a multimodule over $({\mb G}_1, \ldots, {\mb G}_{i-1}, {\mb F}_1, \ldots, {\mb F}_m, {\mb G}_{i+1}, \ldots, {\mb G}_n; {\mb G}')$. Lemma \ref{lemma_concatenation_Novikov} implies that the concatenation is Novikov. Lemma \ref{lemma_concatenation_orientation} says that the concatenation ${\mb X} \circ_i {\mb Y}$ is coherently oriented. Hence one obtains a chain map
\beqn
\tilde{\Phi}^{{\mb X} \circ_i {\mb Y}}: C({\mb G}_1) \widehat{\underset{\mb Z}{\otimes}} \cdots \widehat{\underset{\mb Z}{\otimes}} C({\mb G}_{i-1}) \widehat{\underset{\mb Z}{\otimes}} C ({\mb F}_1) \widehat{\underset{\mb Z}{\otimes}} \cdots \widehat{\underset{\mb Z}{\otimes}} C ({\mb F}_m) \widehat{\underset{\mb Z}{\otimes}} C ({\mb G}_{i+1}) \widehat{\underset{\mb Z}{\otimes}} \cdots \widehat{\underset{\mb Z}{\otimes}} C ({\mb G}_n) \to C ({\mb G}').
\eeqn

\begin{lemma}\label{lemma719}
There holds
\beqn
\tilde{\Phi}^{{\mb X} \circ_i {\mb Y}} (y_1, \cdots, y_{i-1}, x_1, \cdots x_m, y_{i+1}, \cdots, y_n) = \tilde{\Phi}^{\mb Y}(y_1, \cdots, y_{i-1}, \tilde{\Phi}^{\mb X}(x_1, \cdots, x_m), y_{i+1}, \cdots, y_n).
\eeqn
\end{lemma}

\begin{proof}
One only needs to verify this identity for $x_\alpha = p_\alpha$ and $y_\beta = r_\beta$ for generators of the involved chain complexes. Fix $r' \in {\rm Ob}{\mb G}'$. The $r'$-coefficient of the left is the signed count of the disjoint union
\beqn
M_{p_1 \cdots p_m; r_i}^{\mb X} \times M_{r_1 \cdots r_n; r'}^{\mb Y}.
\eeqn
By our orientation conventions (note that the moduli spaces are zero-dimensional), it is equal to 
\beqn
n_{p_1\cdots p_m; p'}^{\mb X} n_{p' r_1 \cdots r_n; r'}^{\mb Y}.
\eeqn
As only zero dimensional moduli spaces contribute, the sign above is positive; therefore, the induced count produces the right hand side.
\end{proof}

Lastly, suppose ${\mb X}_0, {\mb X}_1$ be two, unobstructed, oriented, ${\mb Z}$-graded Novikov multimodules over $({\mb F}_1, \ldots, {\mb F}_m; {\mb F}')$ enriched in $\pman$ of the same degree $n$, and ${\mb H}$ is an oriented, unobstructed, Novikov homotopy from ${\mb X}_0$ to ${\mb X}_1$. For each tuple $(p_1 \cdots p_m; p')$, the orientation on $M_{p_1 \cdots p_m;p'}^{\mb H}$ induces an integral count
\beqn
n_{p_1 \cdots p_m; p'}^{\mb H}\in {\mb Z}.
\eeqn
Define a linear map 
\beqn
\tilde \Psi^{\mb H}: C_*({\mb F}_1) \underset{{\mb Z}}{\otimes} \cdots \underset{{\mb Z}}{\otimes} C_*({\mb F}_m) \to C_*({\mb F}')
\eeqn
by ${\mb Z}$-linearly extending 
\beqn
p_1\otimes \cdots \otimes p_m \mapsto \sum_{p'} n_{p_1 \cdots p_m; p'}^{\mb H} p'.
\eeqn
Using the conditions stated in Definition \ref{defn_orientation_homotopy}, the following statement is immediate.

\begin{lemma}\label{lemma:homotopy}
$\tilde \Psi^{\mb H}$ extends to $C ({\mb F}_1) \widehat{\underset{{\mb Z}}{\otimes}} \cdots \widehat {\underset{{\mb Z}}{\otimes}} C({\mb F}_m)$ and is a chain homotopy from $\tilde \Phi^{{\mb X}_0}$ to $\tilde \Phi^{{\mb X}_1}$. 
\end{lemma}

\subsection{Symmetries on flow categories, multimodules, and homotopies}\label{subsection_Novikov_symmetry}

In the geometric situations considered in this paper, there are three types of symmetries on flow categories, multimodules, and homotopies. The first is the typical kind of Novikov symmetry, which is basically the symmetry related to changing cappings of contractible periodic orbits. The second and the third only appear in the equivariant case. The second type is the symmetry on the Borel model $S^\infty \cong E\zp$, given by the monoid ${\mb Z}_{\geq 0}$, which shifts the coordinates, rather than the group ${\mb Z}$. The third is symmetry by a finite group $\zp$.

\subsubsection{Twisted Novikov groups}

\begin{defn}[Novikov groups]\label{defn_Novikov_group} \hfill
\begin{enumerate}

\item A {\bf Novikov group} consists of a free abelian group $\mb\Pi$ together with morphisms
\begin{align*}
&\ \lambda = (\lambda_1, \ldots, \lambda_k): \mb\Pi \to {\mb R}^k,\ &\ \mu: \mb\Pi \to 2{\mb Z}
\end{align*}
A {\bf morphism} of Novikov groups from $(\mb\Pi, \lambda, \mu)$ to $(\mb\Pi', \lambda', \mu')$ consists of a group homomorphism $\rho: \mb\Pi \to \mb\Pi'$ and a linear map $\sigma: {\mb R}^k \to {\mb R}^{k'}$ such that the following diagrams commute.
\begin{align*}
&\ \xymatrix{  \mb\Pi \ar[r]^{\lambda} \ar[d]_\rho  & {\mb R}^k \ar[d]^\sigma \\
\mb\Pi' \ar[r]_{\lambda'}  & {\mb R}^{k'}}\ &\ \xymatrix{  \mb\Pi \ar[r]^{\mu} \ar[d]_\rho  &  {\mb Z} \ar[d]^{\rm Id} \\
\mb\Pi' \ar[r]_{\mu'}  &  {\mb Z} }
\end{align*}

\item The category of Novikov groups admits finite products, denoted by 
\beqn
(\mb \Pi, \lambda, \mu) \oplus (\mb\Pi', \lambda', \mu') = (\mb\Pi \oplus \mb\Pi', \lambda \oplus \lambda', \mu + \mu').
\eeqn

\item An {\bf twisted Novikov group}, denoted by $\mb\Pi \rtimes G$, consists of a Novikov group $\mb\Pi$ together with a finite group $G$ acting on $\mb\Pi$ by isomorphisms of Novikov groups. Notice that the semi-direct product $\mb\Pi \rtimes G$ is a not necessarily abelian.
\end{enumerate}
\end{defn}

\begin{example}
\begin{enumerate}

\item Let $(X, \omega)$ be a symplectic manifold. Then 
\beqn
\Pi = \frac{\pi_2(X)}{{\rm ker} \omega \cap {\rm ker} c_1}
\eeqn
together with the symplectic area $\lambda(a) = \omega(a)$ and twice the Chern number $\mu(a) = 2 c_1(a)$ is a Novikov group.

\item ${\mb Z}$ with the inclusion $\lambda: {\mb Z} \hookrightarrow {\mb R}$ and the multiplication by two $\mu: {\mb Z} \to 2{\mb Z}\hookrightarrow {\mb Z}$ is a Novikov group.

\item Any finite group is a twisted Novikov group.
\end{enumerate}
\end{example}

\subsubsection{Novikov rings}

Let $R$ be a commutative ring with unit, typically ${\mb Z}$ or $\fp$. Let $(\mb\Pi, \lambda, \mu)$ be a Novikov group. We define the Novikov ring $\Lambda_R^{\mb\Pi}$ as follows. Define
\beq\label{defn_Novikov_ring}
\Lambda_R^{\mb\Pi}:= \left\{  f: \mb\Pi \to R %
\ \left|\ \begin{array}{l} \forall \alpha = 1, \ldots, k,\  \displaystyle \inf \{ \lambda_\alpha(a)\ |\ f(a) \neq 0 \} > -\infty,\\
\forall c>0,\ \# \big\{ a \in \mb\Pi \ |\ f(a) \neq 0,\ \lambda_\alpha (a) < c,\ \alpha = 1, \ldots, k \big\} < \infty \end{array} 
  \right. \right\}.
\eeq

\begin{lemma}
$\Lambda_R^{\mb \Pi}$ is a graded unital commutative ring.
\end{lemma}

\begin{proof}
We can grade $\Lambda_R^{\mb \Pi}$ by dualizing the grading on ${\mb \Pi}$ induced by the morphism $\mu: {\mb \Pi} \to {\mb Z}$. The product structure is defined as follows. Choose $f', f'' \in \Lambda_R^{\mb\Pi}$. Define formally  $f'f'': \mb\Pi \to R$ by
\beq\label{Novikov_ring_multiplication}
f'f'' (a) = \sum_{a' + a'' = a} f'(a') f''(a'').
\eeq
which is commutative in $f'$ and $f''$. We claim that this is a finite sum. If not, then there is a sequence $(a_i', a_i'') \in \mb\Pi \times \mb\Pi$ such that $a_i' + a_i'' =a$ and $f'(a_i'), f''(a_i'') \neq 0$. Then 
\beqn
\lambda_\alpha(a_i') + \lambda_\alpha(a_i'') = \lambda_\alpha(a)
\eeqn
which is independent of $i$. Then
\beqn
\sup_i \lambda_\alpha(a_i'') = \lambda_\alpha(a) - \inf_i \lambda_\alpha (a_i') < \infty. 
\eeqn
The finiteness condition in the definition of $\Lambda_R^{\mb\Pi}$ shows that there are only finitely many distinct $a_i''$ in this sequence. As $a_i' = a - a_i''$, there are also only finitely many distinct $a_i'$ in this sequence. Hence the sum \eqref{Novikov_ring_multiplication} is indeed finite. Therefore, $f'f''$ is defined. It is also straightforward to check that $f'f''$ satisfies the definition of $\Lambda_R^{\mb\Pi}$. Using this definition, one can easily see that the unit is given by the map which assigns the identity element of ${\mb\Pi}$ to the multiplicative unit in $R$ and $0$ otherwise.
\end{proof}

We also define the subring
\beqn
\Lambda_{R,0}^{\mb\Pi}:= \left\{ f \in \Lambda_R^{\mb\Pi}\ |\ f(a) \neq 0 \Longrightarrow \lambda_\alpha (a_i) \geq 0,\ \alpha = 1, \ldots, k \right\}.
\eeqn

Now we consider the tensor product of Novikov rings. If $\mb\Pi_1$, $\mb\Pi_2$ are two Novikov groups, then the Novikov ring 
\beqn
\Lambda_R^{\mb\Pi_1 \oplus \mb\Pi_2}
\eeqn
is a certain completion of the tensor product $\Lambda_R^{\mb\Pi_1} \otimes_R \Lambda_R^{\mb\Pi_2}$. We denote the completed tensor product by $\widehat\otimes$. More generally, if $\mb\Pi_1, \ldots, \mb\Pi_m$ are Novikov groups, then there is a natural identification
\beqn
\Lambda^{\mb\Pi_1}_R \widehat{ \underset{R}{\otimes} } \cdots \widehat{ \underset{R}{\otimes} } \Lambda^{\mb\Pi_m}_R \cong \Lambda_R^{\mb\Pi_1 \oplus \cdots \oplus \mb\Pi_m}.
\eeqn

Lastly, suppose one has a morphism of Novikov groups $\rho: \mb\Pi \to \mb\Pi'$.  Then there is a natural ring map
\beqn
\Lambda_R^{\mb\Pi} \to \Lambda_R^{\mb\Pi'}
\eeqn
making the latter an algebra over $\Lambda_R^{\mb\Pi}$.

\subsubsection{Symmetries on flow categories}

We introduce the concepts related to symmetries on flow categories, multimodules, and homotopies which appear everywhere in the Floer-theoretic setting.

\begin{defn}\label{defn_equivariant_flow_category}
Let $\mb\Pi \rtimes G$ be a twisted Novikov group of rank $k$. Let ${\mb F}$ be a (weak) flow category enriched in $\uds{\bf C}$ equipped with a filtration ${\mc A}: {\rm Ob}{\mb F} \to {\mb R}^k$. A {\bf (free) ${\mb \Pi} \rtimes G$-action} on ${\mb F}$ consists of a (free) ${\mb\Pi} \rtimes G$-action on ${\rm Ob}{\mb F}$ together with $\uds{\bf C}$-isomorphisms
\beqn
T_a^{\mb F}: M_{pq}^{\mb F} \to M_{ap\ aq}^{\mb F},\ \forall a \in {\mb\Pi}\rtimes G
\eeqn
such that 
\begin{enumerate}

\item $T_0^{\mb F}$ is the identity morphism and for all $a, b \in \mb\Pi \rtimes G$, 
\beqn
T_{a+b}^{\mb F} = T_a^{\mb F} \circ T_b^{\mb F}.
\eeqn

\item For $a \in \mb\Pi \rtimes G$ and $p, q, r \in {\rm Ob}{\mb F}$, denote $p' = ap$, $r' = ar$, $q' = aq$. Then the following diagram commutes. 
\beqn
\xymatrix{  M_{pr}^{\mb F}\times M_{rq}^{\mb F} \ar[rr]^-{\iota_{prq}}  \ar[d]_{T_a^{\mb F} \times T_a^{\mb F}}  &  & M_{pq}^{\mb F}\ar[d]^{T_a^{\mb F}}\\
 M_{p' r'}^{\mb F} \times M_{r'q'}^{\mb F} \ar[rr]_-{\iota_{p'r'q'}}  &  &M_{p'q'}^{\mb F} }
\eeqn

\item For $a \in \mb\Pi \subset \mb\Pi \rtimes G$ and $p \in {\rm Ob}{\mb F}$ 
\begin{align*}
&\ {\mc A}(ap) = {\mc A}(p) - \lambda (a)\in {\mb R}^k,\ &\ {\rm deg}(ap) = {\rm deg}(p) - \mu(a).
\end{align*}
For $\theta \in G \subset \mb\Pi \rtimes G$, \begin{align*}
&\ {\mc A}(\theta p) = \theta({\mc A}(p)),\ &\ {\rm deg}(\theta p) = {\rm deg}(\theta p).
\end{align*}

\item Given $c>0$, the set of pairs
\beqn
\Big\{(p, q) \in {\rm Ob}{\mb F}\times {\rm Ob}{\mb F}\ |\ p<q, {\mc A}_\alpha(p) - {\mc A}_\alpha(q) < c\ \forall \alpha = 1, \ldots, k \Big\}
\eeqn
consists of finitely many $\mb\Pi \rtimes G$-orbits.

\end{enumerate}
\end{defn}

\begin{rem}
\begin{enumerate}

\item Notice that the last condition of Definition \ref{defn_equivariant_flow_category} implies that $({\mb F}, {\mc A})$ is a locally finite flow category. 

\item If $\mb\Pi_i$ are Novikov groups and $({\mb F}_i, {\mc A}_i)$ are $\mb\Pi_i$-equivariant Novikov flow categories, then the direct product
\beqn
({\mb F}_1 \times \cdots \times {\mb F}_m, {\mc A}_1 \times \cdots \times {\mc A}_m)
\eeqn
is a $\mb \Pi_1 \times \cdots \times \mb\Pi_m$-equivariant weak Novikov flow category. When all factors are equal to a $\mb\Pi$-equivariant Novikov flow category $({\mb F}, {\mc A})$, the power $({\mb F} \times \cdots \times {\mb F}, {\mc A} \times \cdots \times {\mc A})$ is $\mb\Pi^m \rtimes {\mb Z}_m$-equivariant where ${\mb Z}_m$ acts on $\mb\Pi^m$ by permuting the factors. 
\end{enumerate}
\end{rem}

\subsubsection{Symmetries on multimodules and homotopies}

\begin{defn}\label{defn_equivariant_bimodule}
Let $\mb\Pi \rtimes G$ and $\mb\Pi' \rtimes G$ be twisted Novikov groups. Let $({\mb F}, {\mc A})$ resp. $({\mb F}', {\mc A}')$ be a $\mb\Pi\rtimes G$-equivariant resp. $\mb\Pi' \rtimes G$-equivariant (weak) Novikov flow categories enriched in $\uds{\bf C}$. Let 
\beqn
\rho: \mb\Pi \to \mb\Pi'
\eeqn
be a $G$-equivariant morphism of Novikov groups, which induces a group morphism
\beqn
\rho: \mb\Pi \rtimes G \to \mb\Pi' \rtimes G.
\eeqn
A {\bf $\rho$-equivariant (weak) Novikov bimodule} over $({\mb F}; {\mb F}')$ is a (weak)  bimodule ${\mb B}$ together with $\uds{\bf C}$-isomorphisms
\beqn
T_a^{\mb B}: M_{p; p'}^{\mb B} \to M_{ap; \rho(a) p'}^{\mb B}\ \forall a \in \mb\Pi \rtimes G
\eeqn
satisfying 

\begin{enumerate}
    \item $T_{0}^{\mb B}$ is the identity and for $a, b\in \mb\Pi \rtimes G$, 
\beqn
T_{(a + b)}^{\mb B} = T_a^{\mb B} \circ T_b^{\mb B},
\eeqn

\item The following diagram commutes. 
\beqn
\xymatrix{ M_{pq}^{\mb F} \times  M_{q;p'}^{\mb B} \ar[rr] \ar[d]_-{T_a^{\mb F} \times T_a^{\mb B}}  & &     M_{p;p'}^{\mb B} \ar[d]^{T_a^{\mb B} } & &  M_{p;q'}^{\mb B}\times M_{q'p'}^{{\mb F}'} \ar[ll] \ar[d]^{T_a^{\mb B} \times T_{\rho(a)}^{{\mb F}'}} \\
 M_{ap; aq}^{{\mb F}}   \times M_{aq; \rho(a)p'}^{\mb B} \ar[rr] & &   M_{ap; \rho(a)p'}^{\mb B}  & & M_{ap; \rho(a) q'}^{{\mb B}} \times M_{\rho(a)q' \rho(a)p'}^{{\mb F}'}  \ar[ll]    }
\eeqn

\item For every constant $C>0$, the set 
\beqn
\left\{ (p,p') \in {\rm Ob}{\mb F} \times {\rm Ob}{\mb F}'\ \left| \  \begin{array}{c} M_{p; p'}^{\mb B} \neq \emptyset, \\
{\mc A}_1(p) + \cdots + {\mc A}_k(p) - {\mc A}_1'(p') - \cdots - {\mc A}_{k'}'(p') < C \end{array} \right. \right\}
\eeqn
consists of finitely many $\mb\Pi \rtimes G$-orbits.
\end{enumerate}
\end{defn}

The case of multimodules can be obtained from the case of weak bimodules. 

\begin{defn}\label{defn_equivariant_multimodule}
Let $\mb\Pi_i$, $\mb\Pi'$ be Novikov groups. Let ${\mb F}_i$ resp. ${\mb F}'$ be a $\mb\Pi_i$-equivariant resp. $\mb\Pi'$-equivariant flow categories enriched in $\uds{\bf C}$. Let 
\beqn
\rho: \mb\Pi_1 \times \cdots \times \mb\Pi_m \to \mb\Pi'
\eeqn
be a morphism of Novikov groups. A {\bf $\rho$-equivariant multimodule} over $({\mb F}_1, \ldots, {\mb F}_m; {\mb F}')$ is a multimodule ${\mb X}$, which can also be viewed as a weak bimodule over $({\mb F}_1\times \cdots \times {\mb F}_m; {\mb F}')$, together with $\uds{\bf C}$-isomorphisms
\beqn
T_a^{\mb X}: M_{p_1 \cdots p_m; p'}^{{\mb X}} \to M_{a_1 p_1 \cdots a_m p_m; \rho(a_1, \ldots, a_m) p'}^{\mb X}
\eeqn
for all $a = (a_1, \ldots, a_m) \in \mb\Pi_1 \times \cdots \times \mb\Pi_m$ satisfying the conditions for an equivariant weak Novikov bimodules given in Definition \ref{defn_equivariant_bimodule}. 
\end{defn}

\begin{rem}
The last condition of Definition \ref{defn_equivariant_bimodule} implies that ${\mb B}$ is locally finite (see Definition \ref{defn_local_finite_bimodule}).
\end{rem}

\begin{defn}
Let ${\mb F}, {\mb F}'$, and $\rho$ be as in Definition \ref{defn_equivariant_bimodule}. Let ${\mb B}_0$ and ${\mb B}_1$ be two $\rho$-equivariant bimodules over $({\mb F}; {\mb F}')$. 

\begin{enumerate}

\item A $\rho$-equivariant homotopy from ${\mb B}_0$ to ${\mb B}_1$ is a homotopy ${\mb H}$ from ${\mb B}_0$ to ${\mb B}_1$ together with $\uds{\bf C}$-isomorphisms
\beqn
T_a^{\mb H}: M_{p;p'}^{\mb H} \to M_{ap; \rho(a)p'}^{\mb H}\ \forall a\in \mb\Pi \rtimes G
\eeqn
satisfying 1) $T_{0}^{\mb H}$ is the identity, 2)
\beqn
T_{a+b}^{\mb H} = T_a^{\mb H} \circ T_b^{\mb H},
\eeqn
and 3) the following diagrams commute.
\beqn
\xymatrix{   {\mb M}_{p;p'}^{{\mb B}_0} \sqcup {\mb M}_{p;p'}^{{\mb B}_1} \ar[rr] \ar[d]_{T_a^{{\mb B}_0} \sqcup T_a^{{\mb B}_1} }   & & {\mb M}_{p;p'}^{\mb H}  \ar[d]^{T_a^{\mb H}} \\
   {\mb M}_{ap; \rho(a)p'}^{{\mb B}_0} \sqcup {\mb M}_{ap; \rho(a) p'}^{{\mb B}_1} \ar[rr] & & M_{ap; \rho(a)p'}^{\mb H}    }
\eeqn
\beqn
\xymatrix{ M_{pq}^{\mb F} \times  M_{q;p'}^{\mb H} \ar[rr] \ar[d]_-{T_a^{\mb F} \times T_a^{\mb H}}  & &     M_{p;p'}^{\mb H} \ar[d]^{T_a^{\mb H}}  & &  M_{p;q'}^{\mb H} \times M_{q'p'}^{{\mb F}'} \ar[ll] \ar[d]^{T_a^{\mb H} \times T_{\rho(a)}^{{\mb F}'}  } \\
 M_{ap\ aq}^{\mb F} \times M_{aq; \rho(a) p'}^{\mb H} \ar[rr] & &   M_{ap; \rho(a)p'}^{\mb H}  & & M_{ap; \rho(a)q'}^{\mb H}  \times M_{\rho(a)q'\ \rho(a)p'}^{{\mb F}'}   \ar[ll]    }
\eeqn
When $\mb\Pi = \mb\Pi'$ and $\rho$ is the identity, we simply say that ${\mb H}$ is $\mb\Pi \rtimes G$-equivariant.

\item Suppose in addition that $\mb\Pi \rtimes $ and $\mb\Pi' \rtimes$ are twisted Novikov groups and $\rho: \mb\Pi \to \mb\Pi'$ is an equivariant morphism of Novikov groups. Suppose ${\mb B}_0$ and ${\mb B}_1$ are $\rho$-equivariant Novikov bimodules. Then a $\rho$-equivariant Novikov homotopy is a $\rho$-equivariant homotopy ${\mb H}$ which satisfies, in addition, that for all constants $C>0$, the set 
\beqn
\left\{ (p,p') \in {\rm Ob}{\mb F} \times {\rm Ob}{\mb F}'\ \left| \  \begin{array}{c} M_{p; p'}^{\mb H} \neq \emptyset, \\
{\mc A}_1(p) + \cdots + {\mc A}_k(p) - {\mc A}_1'(p') - \cdots - {\mc A}_{k'}'(p') < C \end{array} \right. \right\}
\eeqn
consists of finitely many $\mb\Pi \rtimes G$-orbits.
\end{enumerate}
\end{defn}

Note that incorporating $G$-symmetry for flow multimodules and bimodule homotopies is straightforward. We do not spell out the details because they are not needed for this paper.

\subsection{Algebraic structures with symmetry}\label{subsection64}
In the previous two subsections, we provide separate treatments of convergence of structural maps of chain complexes and Novikov rings with geometric origin. As usual, in the concrete setting of Floer theory, the two structures interact with each other in the form that recapping of periodic orbits corresponds to multiplication by Novikov variables. Abstractly, equivariance under Novikov groups of flow categories, multimodules, and homotopies leads to linearity of differentials and chain maps under Novikov rings, and local finiteness of orbits under the Novikov group action is compatible with the convergence assumption. We record these classical observations here formally.

\subsubsection{Flow categories and bimodules}

Suppose $({\mb F}, {\mc A})$ is a ${\mb Z}$-graded, oriented, unobstructed, Novikov flow category enriched in $\pman$, then we have obtained previously the complex $C ({\mb F})$ as a free abelian group. 

\begin{prop}\label{prop_equivariant_differential}
Let $(\mb\Pi \rtimes G, \lambda, \mu)$ be a twisted Novikov group. If $({\mb F}, {\mc A})$ is $\mb\Pi \rtimes$-equivariant, then $C ({\mb F})$ has the structure of a ${\mb Z}$-graded $\Lambda_{\mb Z}^{\mb\Pi}$-module with a $G$-action and the differential is $G$-equivariant and $\Lambda_{\mb Z}^{\mb\Pi}$-linear.
\end{prop}

\begin{proof}
The $\mb\Pi$-equivariance condition on ${\mc A}$  implies that $C ({\mb F})$ has a natural $\Lambda_{\mb Z}^{\mb\Pi}$-module structure.  Moreover, the $\mb\Pi$-equivariance on the flow category implies that $d$ is a $\Lambda_{\mb Z}^{\mb\Pi}$-module map. In the presence of the finite group $G$, it follows from definition that the chain complex $C ({\mb F})$ is $G$-equivariant.
\end{proof}

The following proposition is straightforward.

\begin{prop}\label{prop_equivariant_chain_map}
Let $\mb\Pi \rtimes G$ and $\mb\Pi' \rtimes G$ be twisted Novikov groups with a $G$-equivariant morphism $\rho: \mb\Pi \to \mb\Pi'$. Let $({\mb F}, {\mc A})$ resp. $({\mb F}', {\mc A}')$ be a $\mb\Pi\rtimes G$-equivariant resp. $\mb\Pi' \rtimes G$-equivariant ${\mb Z}$-graded, oriented, unobstructed flow categories enriched in $\pman$. Let ${\mb B}$ be a $\rho$-equivariant graded, oriented, and unobstructed bimodule over $({\mb F}; {\mb F}')$. Then the chain map
\beqn
\tilde{\Phi}^{\mb B}: C ({\mb F}) \to C ({\mb F}')
\eeqn
extends to a $\Lambda^{{\mb \Pi}'}$-linear chain map
\beqn
\tilde \Phi^{\mb B}: C({\mb F}) \underset{\Lambda^{\mb\Pi}}{\otimes} \Lambda^{{\mb \Pi}'} \to C({\mb F}')
\eeqn
which is equivariant with respect to the $G$-actions.
\end{prop}

\subsubsection{Multimodules and homotopies}

Recall that by Proposition \ref{prop_chain_map}, given a ${\mb Z}$-graded, oriented, unobstructed multimodule, and locally finite ${\mb X}$ over $({\mb F}_1, \ldots, {\mb F}_m; {\mb F}')$ enriched in $\pman$, one has obtained a chain map 
\beqn
\tilde\Phi^{\mb X}: C ({\mb F}_1) \widehat{\underset{\mb Z}{\otimes}} \cdots \widehat{\underset{\mb Z}{\otimes}} C ({\mb F}_m) \to C({\mb F}')
\eeqn
To obtain a chain map as modules over the Novikov ring, one needs to make use of the equivariance condition.

\begin{prop}\label{prop_Novikov_linear_chain_map}
Suppose $({\mb F}_i, {\mc A}_i), i = 1, \ldots, m$ is equivariant with respect to the action of an (untwisted) Novikov group $(\mb\Pi_i, \lambda_i, \mu_i)$. Then the tensor product 
\beqn
C ({\mb F}_1) \widehat{\underset{\mb Z}{\otimes}} \cdots \widehat{\underset{\mb Z}{\otimes}} C ({\mb F}_m)
\eeqn
is a differential graded module over 
\beqn
\Lambda_{\mb Z}^{\mb\Pi_1} \widehat{\underset{\mb Z}{\otimes}} \cdots \widehat{\underset{\mb Z}{\otimes}} \Lambda_{\mb Z}^{\mb\Pi_m} \cong \Lambda_{\mb Z}^{\mb\Pi_1 \oplus \cdots \oplus \mb\Pi_m}.
\eeqn 
Suppose ${\mb X}$ is equivariant with respect to a Novikov group morphism $\rho: \mb\Pi_1 \oplus \cdots \oplus \mb\Pi_m \to \mb\Pi'$, then the chain map $\tilde \Phi^{\mb X}$ induces a $\Lambda_{\mb Z}^{\mb\Pi'}$-linear chain map
\beqn
\tilde \Phi^{\mb X}: \left( C_*({\mb F}_1) \widehat{\underset{\mb Z}{\otimes} } \cdots \widehat{ \underset{\mb Z}{\otimes}} C_*({\mb F}_m) \right) \underset{\Lambda_{\mb Z}^{\mb\Pi_1\oplus \cdots \oplus \mb\Pi_m}}{\otimes} \Lambda_{\mb Z}^{\mb\Pi'} \to C_*({\mb F}').
\eeqn
In addition, if ${\mb X}_0$ and ${\mb X}_1$ are two such multimodules and there exists a $\rho$-equivariant unobstructed and oriented homotopy from ${\mb X}_0$ to ${\mb X}_1$, then the corresponding chain maps $\Phi^{{\mb X}_0}$ and $\Phi^{{\mb X}_1}$ are homotopic as $\Lambda_{\mb Z}^{{\mb\Pi}'}$-linear chain maps.
\end{prop}

\begin{proof}
It follows from the same reasoning as in Lemma \ref{prop_equivariant_differential}.   
\end{proof}

\begin{rem}
There is a specific case corresponding to the equivariant pair-of-pants product which we would like to address in Section \ref{section_Steenrod} and Section \ref{section_Floer_Steenrod} in a more concrete setting rather than the current abstract setting.
\end{rem}

\subsection{Quantitative theories}

\subsubsection{Nonarchimedean linear algebra}

We recall some basic notions about nonarchimedean linear algebra systematically used in \cite{Usher_Zhang_2016}. Let $(\Pi, \lambda)$ be an ungraded Novikov group; the relevant discussion is insensitive to ${\mb Z}$-grading so we only consider a ${\mb Z}_2$-graded theory. Let ${\bf K}$ be a field. Then on the Novikov field $\Lambda_{\bf K}^\Pi$ there is the standard valuation 
\beqn
{\mf v}\left( \sum_{i=1}^\infty c_i T^{a_i} \right):= \inf \Big\{ \lambda(a_i)\ |\ c_i \neq 0 \Big\} \in {\mb R} \cup \{+\infty \}.
\eeqn
Then $\Lambda_{\bf K}^\Pi$ is a nonarchimedean field. Notice that in the setting of \cite{Usher_Zhang_2016} one requires $\Pi \subset {\mb R}$. We slightly relax this setting to better fit our discussion. To translate to the setting of \cite{Usher_Zhang_2016}, denote
\beq\label{eqn:r-value}
|\Pi|:= \lambda(\Pi) \subset {\mb R}.
\eeq
Then $\Lambda_{\bf K}^{|\Pi|}$ is an algebra over $\Lambda_{{\bf K}}^{\Pi}$.

A {\bf nonarchimedean normed vector space} of $\Lambda_{\bf K}^\Pi$ (cf. \cite[Definition 2.2]{Usher_Zhang_2016}) is a pair $(C, \ell)$ where $C$ is a $\Lambda_{\bf K}^\Pi$-vector space and $\ell: C \to \{-\infty\} \cup {\mb R}$ be a function satisfying
\begin{enumerate}
    \item $\ell (x) = -\infty \Longrightarrow x = 0$.

    \item $\ell (\lambda x) = \ell (x) - {\mf v}(\lambda)$.

    \item $\ell (x + y)  \leq \max \{ \ell(x),\ \ell(y)\}$.
\end{enumerate}
There is a notion called {\bf orthogonalizability} on nonarchimedean normed vector spaces (\cite[Definition 2.10]{Usher_Zhang_2016}) over $\Lambda_{{\bf K}}^{|\Pi|}$, which is satisfied by all concrete cases considered in Floer theory. Finally, Usher--Zhang defined the notion called {\bf Floer-type complex} which is crucial in their abstract discussion on barcodes. As the consideration does not necessarily involve grading, we slightly modify the original definition.

\begin{defn}(cf.\cite[Definition 4.1]{Usher_Zhang_2016})\label{defn_Floer_type_complex} Let $(\Pi, \lambda)$ be an ungraded Novikov group and ${\bf K}$ be a field.
\begin{enumerate}

\item A nonarchimedean normed vector space $(C, \ell_C)$ over $\Lambda_{\bf K}^\Pi$ is called {\bf orthogonalizable} if the tensor product $C \otimes \Lambda_{\bf K}^{|\Pi|}$ is orthogonalizable over $\Lambda_{\bf K}^{|\Pi|}$. 

\item A {\bf Floer-type complex} over $\Lambda_{\bf K}^\Pi$ consists of a ${\mb Z}_2$-graded orthogonalizable $\Lambda_{\bf K}^\Pi$-space $C$ together with an odd $\Lambda_{\bf K}^\Pi$-differential $d_C: C \to C$ satisfying
\beqn
\ell_C(d_C x) \leq \ell_C(x).
\eeqn
\end{enumerate}
\end{defn}

\begin{rem}
There is another notion called Floer--Novikov complex defined in \cite{Usher_2008}. It requires that the differential strictly decreases the nonarchimedean norm. On the other hand, that notion allows one to consider Novikov rings with coefficients in a general Noetherian ring and the exponents $\Gamma$ to be a general finite-rank free abelian group. In the Floer setup, the Floer complex can be made both a Floer--Novikov complex and a Floer-type complex with appropriate Novikov coefficients.
\end{rem}

Next we have two notions about chain maps and chain homotopies.

\begin{defn}(cf. \cite[Definition 4.5]{Usher_Zhang_2016})\label{defn_filtered_chain_map}
Let $C$ and $D$ be two Floer-type complexes over $\Lambda_{\bf K}^\Pi$. 
\begin{enumerate}

\item A linear map $\Phi: C \to D$ is said to {\bf preserve the filtration} if 
\beqn
\ell_D(\Phi(x)) \leq \ell_C(x)\ \forall x\in C.
\eeqn

\item $C$ and $D$ are said to be {\bf filtered isomorphic} if there is an isomorphism $\Phi: C \to D$ of chain complexes such that 
\beqn
\ell_D(\Phi(x)) = \ell_C(x)\ \forall x \in C.
\eeqn

\item Two chain maps $\Phi, \Psi: C \to D$ are said to be {\bf filtered chain homotopic} if 1) both $\Phi$ and $\Psi$ preserve the filtration and 2) there exists a homotopy $K: C \to D$ between $\Phi$ and $\Psi$ which also preserves the filtration. Two Floer-type complexes are {\bf filtered chain homotopy equivalent} if there exists a filtered chain homotopy $\Phi: C \to D$ and a filtered homotopy inverse $\Psi: D \to C$ such that $\Phi \circ \Psi$ is filtered chain homotopic to ${\rm Id}_D$ and $\Psi \circ \Phi$ is filtered chain homotopic to ${\rm Id}_C$.

\end{enumerate}
\end{defn}

One uses the following notion to compare Floer type complexes. 

\begin{defn}
Let $(C, \ell_C, d_C)$ and $(D, \ell_D, d_C)$ be two Floer-type complexes and $\delta \geq 0$. A {\bf $\delta$-quasiequivalence} between $C$ and $D$ is a quadruple $(\Phi, \Psi, K_C, K_D)$ where $\Phi: C \to D$, $\Psi: D \to C$ are chain maps, $K_C$, $K_D$ are chain homotopies as
\begin{align*}
    &\ \Psi \circ \Phi - {\rm Id}_C = d_C \circ K_C + K_C \circ d_C,\ &\ \Phi \circ \Psi - {\rm Id}_D = d_D \circ K_D + K_D \circ d_D
\end{align*}
satisfying for all $x \in C$ and $y \in D$, 
\begin{align*}
&\ \ell_D(\Phi(x)) \leq \ell_C (x) + \delta,\ &\ \ell_C(\Psi(y)) \leq \ell_D (y) + \delta,
\end{align*}
and
\begin{align*}
&\ \ell_C(K_C(x)) \leq \ell_C(x) + 2\delta,\ &\ \ell_D(K_D(y)) \leq \ell_D(y) + 2\delta.
\end{align*}
The {\bf quasiequivalence distance} between two Floer-type complexes is the infimum of $\delta$ for which there exists a $\delta$-quasiequivalence between them. 
\end{defn}

Below we show that the complex obtained from a Novikov flow category is indeed a Floer-type complex, once we use field coefficients.

\begin{prop}\label{prop738}
Let $\Pi$ be a finitely generated Novikov group and let $({\mb F}, {\mc A})$ be a $\Pi$-equivariant Novikov flow category enriched in $\pman$ satisfying in addition
\begin{enumerate}

\item The action $\lambda$ on $\Pi$ is ${\mb R}$-valued and ${\mc A}$ is ${\mb R}$-valued.

\item ${\rm Ob}{\mb F}/\Pi$ is finite.
\end{enumerate}
Then these data together with the differential on the chain complex $C( {\mb F})$ define a filtered Floer--Novikov complex in the sense of \cite[Definition 1.1]{Usher_2008}. 

Moreover, for any field ${\bf K}$, the complex 
\beqn
C({\mb F}) \otimes \Lambda_{\bf K}^{\Pi}
\eeqn
together with the norm $\ell:= {\mc A}$ is a Floer-type complex over $\Lambda_{\bf K}^{\Pi}$ in the sense of Definition \ref{defn_Floer_type_complex}.
\end{prop}

\begin{proof}
The first claim is straightforward where we take ${\rm Ob}(\mathring{\mb F})$ under the $\Pi$-action equipped with the action ${\mc A}$, and we use the incidence coefficients of the differential. For the second, see \cite[Section 12]{Usher_Zhang_2016}.
\end{proof}

\subsubsection{Filtration and spectral invariants}

For any Floer-type complex $(C, \ell_C, d_C)$, one can consider the filtration defined by 
\begin{align*}
&\ C^{\leq \tau} = \Big\{ x \in C\ |\ \ell_C (x) \leq \tau \Big\},\ &\ C^{<\tau}:= \Big\{ x \in C\ |\ \ell_C(x) < \tau \Big\}.
\end{align*}
As $d_C$ respects the norm, $C^{\leq \tau}$ is a subcomplex.  Then for any interval $I \subset {\mb R}$, one has a corresponding complex defined in a typical way. For example, when $I = [a, b]$, $C^I$ is the quotient complex
\beqn
C^I:= C^{\leq b}/ C^{< a}
\eeqn

For each $x \in H(C)$, define the {\bf spectral number} to be 
\beqn
c(x):= \inf \Big\{ \tau \in {\mb R}\ |\ x \in {\rm Im} \Big( H(C^{<\tau}) \to H(C) \Big) \Big\} \in {\mb R} \cup \{-\infty\}.
\eeqn
When the complex $C$ comes from a flow category ${\mb F}$ enriched in $\pman$, denote the spectral number by 
\beqn
c_{{\mb F}}(x)
\eeqn
where the coefficient ring is understood from the context.

Usher \cite{Usher_2008} provided an abstract argument showing the finiteness of spectral numbers and existence of spectral carriers.

\begin{thm}[cf. \cite{Usher_2008}]\label{thm_Usher}
Under the assumption of Proposition \ref{prop738}, for any Noetherian ring $R$, and $a \in H_*({\mb F}) \neq 0$, $\tau(a):= c_{\mb F}(a) > -\infty$. Moreover, $a$ is represented by a cycle contained in $C_*({\mb F})^{\leq \tau(a)}$.
\end{thm}

\begin{proof}
The theorem follows from Proposition \ref{prop738} and \cite[Theorem 1.3]{Usher_2008}.
\end{proof}

\subsubsection{Barcodes}

By the nonarchimedean Gram--Schmidt process established in \cite[Definition 6.3]{Usher_Zhang_2016}, given a Floer-type complex $(C, \ell_C, d_C) = (C, \ell, d)$, the map $d: C \to C$ of the ungraded complex factors through $C \to \ker(d) \subset C$, so we can find an orthogonal basis $y_1, \ldots, y_n$ of $C$ and an orthogonal basis of $x_1, \ldots, x_m$ of ${\rm ker} (d)$ such that one can write $d: C \to \ker(d)$ in the block-matrix form
\beqn
\left[ \begin{array}{cc} I_r & 0 \\
0 & 0 \end{array}\right]
\eeqn
where $r\leq n$ is the rank of $d$ and such that 
\beqn
\ell(y_1) - \ell(x_1) \geq \cdots \geq \ell(y_r) - \ell(x_r).
\eeqn
The {\bf verbose barcode} of $(C, d, \ell)$ is the multiset of elements $({\mb R}/ \Gamma) \times [0, +\infty]$ consisting of
\begin{enumerate}
    \item a pair $(\ell(x_i)\ {\rm mod}\ \Gamma, \ell(y_i) - \ell(x_i))$ for $i=1,\ldots,r$.

    \item a pair $(\ell(x_i)\ {\rm mod}\ \Gamma, \infty)$ for $i = r+1, \ldots, m$.
\end{enumerate}
Each element of this multiset is called a bar while the second variable is called the bar-length. The {\bf concise barcode} is the submultiset obtained by removing bars of length zero. It was proved (\cite[Theorem 7.1]{Usher_Zhang_2016}) that the verbose barcode of a Floer-type complex does not depend on the choice of the singular value decomposition, hence is well-defined.

\begin{thm}\label{thm_Usher_Zhang_barcodes} Let $C$ and $D$ be two Floer-type complexes over $\Lambda_{\bf K}^\Pi$. 
\begin{enumerate}

\item If $C$ and $D$ are filtered isomorphic, then they have identical verbose barcodes.

\item If $C$ and $D$ are filtered homotopy equivalent, then they have identical concise barcodes.
\end{enumerate}
\end{thm}

\begin{proof}
By \cite[Theorem A, B]{Usher_Zhang_2016}, this is true if the Novikov group $\Pi$ is contained in ${\mb R}$. In our setting, the conclusion holds after tensoring with $\Lambda_{\bf K}^{|\Pi|}$.
\end{proof}

\subsubsection{Continuity}

Now we consider the case of tensor product and discuss quantitative aspects of chain maps. Let $({\mb F}_1, {\mc A}_1), \ldots, ({\mb F}_m, {\mc A}_m)$ be locally finite flow categories enriched in $\pman$ so there are the associated chain complexes $C ({\mb F}_j)$. Recall that one can consider the completed tensor product (over ${\mb Z}$)
\beqn
C ({\mb F}_1) \widehat{\otimes} \cdots \widehat{\otimes} C ({\mb F}_m)
\eeqn
whose elements are formal sums
\beqn
\sum_{i=1}^\infty a_i \big( p_{i_1} \otimes \cdots \otimes p_{i_m}\big)
\eeqn
with $a_i \in R$ such that 
\beqn
\lim_{i \to \infty} \big( {\mc A}_1(p_{i_1}) + \cdots + {\mc A}_m( p_{i_m}) \big) = -\infty.
\eeqn
Then there is a similar filtration on the tensor product and for $\tau_1, \ldots, \tau_m \in {\mb R}$, an inclusion
\beqn
C^{\leq \tau_1}({\mb F}_1) \widehat{\otimes} \cdots \widehat{\otimes} C^{\leq \tau_m}({\mb F}_m) \to \big( C({\mb F}_1) \widehat{\otimes} \cdots \widehat{\otimes} C ({\mb F}_m )\big)^{\leq \tau_1 + \cdots + \tau_m}
\eeqn

\begin{defn}\label{defn_threshold}
Let ${\mb F}_1, \ldots, {\mb F}_m, {\mb F}'$ be Novikov flow categories with ${\mb R}$-valued actions ${\mc A}_1, \ldots, {\mc A}_m, {\mc A}'$. 

\begin{enumerate}

\item Let ${\mb X}$ be a Novikov multimodule over $({\mb F}_1, \ldots, {\mb F}_m; {\mb F}')$. The {\bf threshold} of ${\mb X}$ is the number
\beqn
\epsilon^{\mb X}:= \inf \Big\{ {\mc A}_1 (p_1) + \cdots + {\mc A}_m(p_m) - {\mc A}' (p')\ |\ M_{p_1 \cdots p_m; p'}^{\mb X}  \neq \emptyset \Big\} \in [-\infty, \infty].
\eeqn

\item Let ${\mb X}_0$ and ${\mb X}_1$ be two Novikov multimodules over $({\mb F}_1, \ldots, {\mb F}_m; {\mb F}')$. Let ${\mb H}$ be a Novikov homotopy from ${\mb X}_0$ to ${\mb X}_1$. The {\bf threshold} of ${\mb H}$ is the number
\beqn
\delta^{\mb H}:= \inf \Big\{ {\mc A}_1(p_1) + \cdots + {\mc A}_m (p_m) - {\mc A}'(p') \ |\ M_{p_1 \cdots p_m; p'}^{\mb H} \neq \emptyset \Big\} \in [-\infty, +\infty].
\eeqn
\end{enumerate}
\end{defn}

By the Novikov property of multimodules (see Definition \ref{defn_local_finite_flow}), $\epsilon^{\mb X} \neq -\infty$ and $\delta^{\mb H} \neq -\infty$.

\begin{prop}\label{prop_spectral_continuity}
In the situation of Definition \ref{defn_threshold}, the chain map
\beqn
\Phi^{\mb X}: C({\mb F}_1) \widehat{\otimes} \cdots \widehat{\otimes} C({\mb F}_m) \to C({\mb F}')
\eeqn
satisfies 
\beqn
\Phi^{\mb X} \left( \big( C ({\mb F}_1) \widehat{\otimes} \cdots \widehat{\otimes} C ({\mb F}_m) \big)^{\leq \tau} \right) \subset C_*^{ \leq \tau - \epsilon^{\mb X}} ({\mb F}').
\eeqn
The chain homotopy
\beqn
\Phi^{\mb H}: C({\mb F}_1) \widehat{\otimes}\cdots \widehat{\otimes} C({\mb F}_m) \to C({\mb F}')
\eeqn
satisfies
\beqn
\Phi^{\mb H} \left( \big( C ({\mb F}_1) \widehat{\otimes} \cdots \widehat{\otimes} C ({\mb F}_m) \big)^{\leq \tau} \right) \subset C_*^{ \leq \tau - \delta^{\mb H}} ({\mb F}').
\eeqn
\end{prop}

\begin{proof}
Straightforward check.
\end{proof}

\section{FOP Perturbations}\label{sec:FOP-perturb}

In this section we review the technique of FOP perturbations. The original idea was proposed by Fukaya--Ono in \cite{Fukaya_Ono_integer}, and further explored by Parker \cite{BParker_integer}. In \cite{Bai_Xu_2022} the authors rigorously established this technique and defined integer-valued Gromov--Witten invariants for general compact symplectic manifolds. In \cite{Bai_Xu_Arnold} the authors further proved the integral Arnold conjecture; see \cite{BPX} for another application.

\subsection{FOP perturbations}

\begin{defn}[Perturbations on derived orbifolds%
] \hfill
Given a derived orbifold ${\mc D} = ({\mc U}, {\mc E}, {\mc S})$, a {\bf perturbation} of ${\mc D}$ (or of ${\mc S}$) is a smooth section ${\mc S}': {\mc U} \to {\mc E}$ such that if ${\mc S}^{-1}(0)$ is compact, then $({\mc S}')^{-1}(0)$ is also compact.

\end{defn}

We recall the main ``black-box'' type theorem proved in \cite{Bai_Xu_2022}.

\begin{thm}\label{thm_FOP_property}\cite[Theorem 1.1]{Bai_Xu_2022}
Suppose ${\mc U}$ is an NC orbifold without boundary\footnote{When discussing FOP transverse sections for orbifolds with boundaries and corners, the behavior near boundary and corner strata are predetermined.} and ${\mc E} \to {\mc U}$ is an NC vector bundle. Let $\Gamma({\mc U}, {\mc E})$ be the space of smooth sections. Then there is a $C^0$-dense subset $\Gamma^{\rm FOP}({\mc U}, {\mc E}) \subset \Gamma({\mc U}, {\mc E})$ whose elements are called \emph{FOP transverse sections} satisfying the following properties.
\begin{enumerate}

\item {\bf (Classical transversality)} If ${\mc U}$ is a manifold, FOP transversality is equivalent to classical transversality.

\item {\bf (Locality Property)} The restrictions of FOP transverse sections to open subsets are still FOP transverse. 

\item {\bf (Extension Property)} For any pair of closed subsets ${\mc Y} \subset {\mc Y}'$ and open neighborhoods ${\mc V} \subset {\mc U}$ of ${\mc Y}$ and ${\mc V}' \subset {\mc U}$ of ${\mc Y}'$, if ${\mc S} \in \Gamma^{\rm FOP}({\mc V}, {\mc E}|_{{\mc V}})$, then there exists ${\mc S}' \in \Gamma^{\rm FOP}({\mc V}', {\mc E}|_{{\mc V}'})$ which agrees with ${\mc S}$ near ${\mc Y}$.

\item {\bf (Product Property)} Let ${\mc U}'$ be another normally complex  orbifold and ${\mc E}' \to {\mc U}'$ be a normally complex  vector bundle. Then the product map 
\beqn
\Gamma({\mc U}, {\mc E}) \times \Gamma({\mc U}', {\mc E}') \to \Gamma({\mc U} \times {\mc U}', {\mc E}\boxplus {\mc E}')
\eeqn
sends products of FOP transverse sections to FOP transverse sections.

\item {\bf (Stabilization Property)} If $\pi_{\mc F}: {\mc F} \to {\mc U}$ is rigid complex vector bundle\footnote{This condition was stated in \cite{Bai_Xu_2022} for more general normally complex vector bundles. When ${\mc F}$ is a rigid bundle, the total space ${\mc F}$ carries a canonical NC structure.}, then the stabilization map 
\beqn
\Gamma({\mc U}, {\mc E}) \to \Gamma({\mc F},  \pi_{\mc F}^* {\mc E} \oplus \pi_{\mc F}^* {\mc F}),\quad \quad {\mc S} \mapsto \pi_{\mc F}^* {\mc S} \oplus \tau_{\mc F},
\eeqn
where $\tau_{\mc F}: {\mc F} \to \pi_{\mc F}^* {\mc F}$ is the tautological section, sends FOP transverse sections to FOP transverse sections.

\item {\bf (Main stratum)} For each ${\mc S}\in \Gamma^{\rm FOP}({\mc U}, {\mc E})$, the closure 
\beqn
\ov{ {\mc S}^{-1}(0) \cap {\mc U}_{\rm free}}
\eeqn
is a pseudomanifold (whose dimension is ${\rm dim} {\mc U} - {\rm rank} {\mc E}$).
\end{enumerate}

\end{thm}

\begin{defn}
Let $\uds{\bf dOrb}^{\rm FOP}$ resp. $\uds{\bf dOrb}^{\rm FOP}_{\rm rig}$ be the category whose objects are 
\beqn
\mathring {\mc D}:=({\mc U}, {\mc E}, {\mc S}, {\mc S}')
\eeqn
where ${\mc D}:= ({\mc U}, {\mc E}, {\mc S})$ is an NC derived orbifold and ${\mc S}': {\mc U} \to {\mc E}$ is an FOP transverse perturbation. The morphism from $\mathring {\mc D}_1:= ({\mc U}_1, {\mc E}_1, {\mc S}_1, {\mc S}_1')$ to $\mathring {\mc D}_2:= ({\mc U}_2, {\mc E}_2, {\mc S}_2, {\mc S}_2')$ is an embedding of underlying stratified orbifolds covered by a bundle embedding ${\bm \iota}_{21} = (\iota_{21}, \wh\iota_{21}): ({\mc U}_1, {\mc E}_1) \to ({\mc U}_2, {\mc E}_2)$ such that it induces both a morphism of $\uds{\bf dOrb}^{\rm NC}$ resp. $\uds{\bf dOrb}^{\rm NC}_{\rm rig}$ from ${\mc D}_1:=({\mc U}_1, {\mc E}_1, {\mc S}_1)$ to ${\mc D}_2:=({\mc U}_2, {\mc E}_2, {\mc S}_2)$ and a morphism with ${\mc S}_i$ replaced by ${\mc S}_i'$.
\end{defn}

Notice that there are natural forgetful functors
\begin{align*}
 \vcenter{ \xymatrix{ & \uds{\bf dOrb}^{\rm FOP} \ar[ld] \ar[rd] & \\
\uds{\bf dOrb}^{{\rm NC}} & & {\pman} }}
\end{align*}
by forgetting the perturbation part ${\mc S}'$ resp. taking the closure of $({\mc S}')^{-1}(0) \cap {\mc U}_{\rm free}$.

By using Theorem \ref{thm_FOP_property} inductively, one can construct FOP transverse perturbations on stratified derived orbifolds which respect collars.

\begin{lemma}\label{lemma_FOP_induction}
Let $A$ be a regular poset and ${\mc D} = ({\mc U}, {\mc E}, {\mc S})$ be an $A$-stratified object of $\outer \uds{\bf dOrb}^{\rm NC}$. Let $A' \subset A$ be an Alexandrov closed subset and 
\beqn
\partial^{A'} {\mc D}
\eeqn
be the corresponding $A'$-stratified derived orbifold. Suppose ${\mc S}_{A'}: \partial^{A'} {\mc U} \to {\mc E}$ is an FOP transverse perturbation so that 
\beqn
(\partial^{A'} {\mc U}, {\mc E}, {\mc S}, {\mc S}_{A'})
\eeqn
is an $A'$-stratified object of $\outer \uds{\bf dOrb}^{\rm FOP}$. Then there exists an FOP transverse perturbation ${\mc S}': {\mc U} \to {\mc E}$ which extends ${\mc S}_{A'}$ such that $({\mc U}, {\mc E}, {\mc S}, {\mc S}')$ is an object of $\outer \uds{\bf dOrb}^{\rm FOP}$. 
\end{lemma}

\subsection{FOP perturbations on flow categories, multimodules, and homotopies}

Now we prove a major result of this paper stating the general procedure of constructing FOP perturbations on flow categories, multimodules, and homotopies.

\begin{thma}\label{thma_FOP} 
Flow categories, multimodules, and homotopies in this proposition are assumed to be unobstructed (Definition \ref{defn_unobstructed}) and locally compact (Definition \ref{defn_local_finite_flow}, \ref{defn_local_finite_bimodule}, \ref{defn_Novikov_multimodule}, and \ref{defn_Novikov_homotopy}). We also assume that all charts have compact zero locus.
\begin{enumerate}

\item {\bf (FOP perturbation for equivariant flow categories)} Let $\mb\Pi \rtimes G$ be a twisted Novikov group (Definition \ref{defn_Novikov_group}) and let $\tilde {\mb F}$ be a free $\mb\Pi \rtimes G$-equivariant flow category enriched in $\outer \uds{\bf dOrb}_{\rm rig}^{\rm NC}$. Then there exists a $\mb\Pi\rtimes G$-equivariant FOP perturbation $\mathring {\mb F}$ on $\tilde {\mb F}$.

\item {\bf (FOP perturbation for equivariant multimodules)} Let $\mb\Pi_1, \ldots, \mb\Pi_m; \mb\Pi'$ be Novikov groups and $G$ be a finite group acting on $\mb\Pi_1\times \cdots \times \mb\Pi_m$ and on $\mb\Pi'$. Let
\beqn
\rho: \mb\Pi_1\times \cdots \times \mb\Pi_m \to \mb\Pi'
\eeqn
be a $G$-equivariant morphism of Novikov groups. Let $\tilde {\mb F}_1, \ldots, \tilde {\mb F}_m$ be $\mb\Pi_i$-equivariant flow categories enriched in $\outer \uds{\bf dOrb}_{\rm rig}^{\rm NC}$ such that the weak flow category $\tilde {\mb F}_1 \times \cdots \times \tilde {\mb F}_m$ has a free $(\mb\Pi_1 \times \cdots \times {\mb \Pi}_m) \rtimes G$-action. Let $\tilde{\mb F}'$ be a free $\mb\Pi'\rtimes G$-equivariant flow category enriched in $\outer \uds{\bf dOrb}_{\rm rig}^{\rm NC}$. Let $\tilde {\mb X}$ be a $\rho$-equivariant multimodule over $(\tilde{\mb F}_1, \ldots, \tilde{\mb F}_m; \tilde {\mb F}')$. Let $\mathring {\mb F}_i$ resp. $\mathring {\mb F}'$ be a $\mb\Pi_i$-equivariant resp. $\mb\Pi'$-equivariant FOP perturbation on $\tilde {\mb F}_i$ resp. $\tilde {\mb F}'$  such that $\mathring {\mb F}_1 \times \cdots \times \mathring {\mb F}_m$ resp. $\mathring {\mb F}'$ is also $G$-equivariant. Then there exists a $\rho$-equivariant FOP perturbation $\mathring {\mb X}$ on $\tilde {\mb X}$ as an equivariant multimodule over $(\mathring {\mb F}_1, \ldots, \mathring {\mb F}_m; \mathring {\mb F}')$.

\item {\bf (FOP perturbation for equivariant multimodule homotopies)} Let $\mb\Pi_1, \ldots, \mb\Pi_m$, $\mb\Pi'$, $G$, $\rho$, $\tilde {\mb F}_1, \ldots, \tilde {\mb F}_m, \tilde {\mb F}'$, and $\mathring {\mb F}_1, \ldots, \mathring {\mb F}_m, \mathring {\mb F}'$ be as in the above situation. Let $\tilde {\mb X}_0$ and $\tilde {\mb X}_1$ be two $\rho$-equivariant multimodules over $(\tilde {\mb F}_1, \ldots, \tilde {\mb F}_m; \tilde {\mb F}')$ and $\tilde {\mb H}$ be a $\rho$-equivariant homotopy from $\tilde {\mb X}_0$ to $\tilde {\mb X}_1$. Let $\mathring {\mb X}_0$ and $\mathring {\mb X}_1$ be FOP perturbations on $\tilde {\mb X}_0$ and $\tilde {\mb X}_1$ respectively as $\rho$-equivariant multimodules over $(\mathring {\mb F}_1, \ldots, \mathring {\mb F}_m; \mathring {\mb F}')$. Then there exists a $\rho$-equivariant FOP perturbation $\mathring {\mb H}$ on $\tilde {\mb H}$ as a homotopy from $\mathring {\mb X}_0$ to $\mathring {\mb X}_1$.
\end{enumerate}
\end{thma}

\begin{proof}[Proof of Theorem \ref{thma_FOP}]

We construct the FOP perturbations inductively. In general, the Novikov conditions of Definition \ref{defn_local_finite_flow}, Definition \ref{defn_local_finite_bimodule}, and Definition \ref{defn_equivariant_flow_category} allow us to choose a beginning of the induction process. The freeness assumed in Theorem \ref{thma_FOP} allows us to achieve equivariant transversality.

---{\bf The case of flow categories}--- We say a pair of objects $p<q$ is minimal if $M_{pq}^{\tilde{\mb F}} \neq \emptyset$ and there does not exist a third object $r$ such that both $M_{pr}^{\tilde {\mb F}}$ and $M_{rq}^{\tilde {\mb F}}$ are nonempty. By (2) of Definition \ref{defn_equivariant_flow_category}, unless all morphism spaces are empty, minimal pairs exist. Moreover, the $\mb\Pi$-equivariance condition implies that minimal pairs form a union of $\mb\Pi$-orbits. 

Now we perturb $M_{pq}^{\tilde {\mb F}} = ({\mc U}_{pq}, {\mc E}_{pq}, {\mc S}_{pq})$ for minimal pairs $p<q$. We claim that ${\mc S}_{pq}^{-1}(0)$ is strictly contained in the union of top strata of ${\mc D}_{pq}$. If it is not the case, then there exists a codimension one stratum which intersects ${\mc S}_{pq}^{-1}(0)$. By the unobstructedness assumption, such a stratum comes from the product of two other morphism spaces. This contradicts the assumption that $p<q$ is a minimal pair. Hence our claim is correct. 

Then by Theorem \ref{thm_FOP_property}, there exists an FOP transverse perturbation ${\mc S}_{pq}'$. Then using the free $\mb\Pi$-action, we translate this perturbation to charts $M_{p'q'}^{\tilde{\mb F}}$ with all $(p', q')$ in the same $\mb\Pi$-orbit of $(p, q)$. We apply such perturbation constructions independently for all $\mb\Pi$-orbits of minimal pairs. 

Now we set up the induction. Recall that the Novikov flow category ${\mb F}$ carries an ${\mb R}^k$-valued action ${\mc A} = ({\mc A}_1, \ldots, {\mc A}_k)$. Define
\beqn
|{\mc A}| = {\mc A}_1 + \cdots + {\mc A}_m: {\rm Ob}{\mb F} \to {\mb R}.
\eeqn
Choose an increasing sequence $c_1 < c_2 < \cdots$ satisfying the following conditions.
\begin{enumerate}

    \item $c_l$ goes to infinity.

    \item For each pair $p<q$ with $c_{l-1} \leq |{\mc A}|(p) - |{\mc A}| (q) < c_l$, if $p<r<q$, then $|{\mc A}| (p) - |{\mc A}|(r) < c_{l-1}$ and $|{\mc A}| (r) - |{\mc A}| (q) < c_{l-1}$. 

    \item If $p<q$, then $|{\mc A}| (p)  - |{\mc A}| (q) \geq c_1$. 
\end{enumerate}

Now we state the following induction hypothesis: suppose for $c_l$ in this sequence, we have constructed FOP transverse perturbations
\beqn
{\mc S}_{pq}': {\mc U}_{pq} \to {\mc E}_{pq}
\eeqn
for all pairs $p<q$ with $| {\mc A} |(p) - |{\mc A}| (q) < c_l$ for all $\alpha = 1, \ldots, k$, which satisfy the following conditions.
\begin{enumerate}

\item The perturbations are $\mb\Pi$-equivariant.

\item Whenever $p<r<q$, the restriction of ${\mc S}_{pq}'$ to $\partial^{prq} {\mc U}_{pq}$ coincides with the stabilization of the product of ${\mc S}_{pr}'$ and ${\mc S}_{rq}'$ under the compositions of $\tilde {\mb F}$. 

\item If $p<q$ is minimal, then ${\mc S}_{pq}'$ is the one we have already constructed.
\end{enumerate}

Our choice of $c_1$ implies that the base case is established. Suppose the induction hypothesis holds for $c_{l-1}$. 

Consider a pair $p<q$ with $| {\mc A} | (p) - | {\mc A} | (q) \in [c_{l-1}, c_l)$. If the pair $(p, q)$ is minimal, then we take ${\mc S}_{pq}'$ to be the one we constructed before the induction setup. Suppose it is not minimal. By the induction hypothesis, for each $r$ with $p<r<q$, one has constructed FOP transverse perturbations
\begin{align*}
&\ {\mc S}_{pr}': {\mc U}_{pr} \to {\mc E}_{pr},\ &\ {\mc S}_{rq}': {\mc U}_{rq} \to {\mc E}_{rq}.
\end{align*}
As the structural map
\beqn
M_{pr}^{\mb F} \times M_{rq}^{\mb F} \to \partial^{prq} M_{pq}^{\mb F}
\eeqn
is a rigidified embedding of normally complex derived orbifolds, there is a canonical extension of the product ${\mc S}_{pr}'\times {\mc S}_{rq}'$ to $\partial^{prq} {\mc U}_{pq}$. The induction hypothesis also implies that for all such $r$, these canonical extensions agree on strata of the form $\partial^{pr_1 r_2 q} {\mc U}_{pq}$. Hence one obtains an object in $\outer \uds{\bf dOrb}_{\rm rig}^{\rm FOP}$ stratified by the true boundary of $A_{pq}^{\tilde {\mb F}}$. Then using Lemma \ref{lemma_FOP_induction}, there exists an FOP transverse perturbation ${\mc S}_\alpha: \partial^\alpha {\mc U}_{pq} \to {\mc E}_{pq}$ which is compatible with the existing ones. 

As the $\mb\Pi$-action is free, one can translate ${\mc S}_{pq}'$ to all pairs on the same $\mb\Pi$-orbit. Apply the same procedure for all $\mb\Pi$-orbits of pairs $p', q'$ with $|{\mc A}| (p') - | {\mc A}| (q') \in [c_{l-1}, c_l)$. This implies that the induction can be carried out. Hence the case of flow categories of Theorem \ref{thma_FOP} is established.

---{\bf The case of multimodules---} We follow a similar inductive process. Consider the set ${\rm Ob}\tilde {\mb F}_1\times \cdots \times {\rm Ob}\tilde{\mb F}_m \times {\rm Ob}\tilde{\mb F}'$, whose objects are denoted by $(p_1,\ldots, p_m; p')$. Then $\mb\Pi_1 \times \cdots \times \mb\Pi_m$ acts freely on this set via the morphism $\rho$; in the situation of item (3), the group $\Gamma$ also acts on this set and the $\Gamma$-action commutes with $\mb\Pi_1 \times \cdots \times \mb\Pi_m$. We say a tuple $(p_1, \ldots, p_m; p')$ is minimal if 1) there does not exist $i$ and $p_i < q_i$ in $\tilde{\mb F}_i$ such that $M_{p_1 \cdots p_{i-1} q_i p_{i+1} \cdots p_m; p'}^{\tilde {\mb X}} \neq \emptyset$ and 2) there does not exist $q' < p'$ in $\tilde{\mb F}'$ such that $M_{p_1 \cdots p_m; q'}^{\tilde {\mb X}} \neq \emptyset$. When $(p_1, \ldots, p_m; p')$ is minimal, similar to the flow category case, the unobstructedness assumption implies that one can construct an FOP transverse perturbation 
\beqn
{\mc S}_{p_1 \cdots p_m; p'}': {\mc U}_{p_1\cdots p_m; p'} \to {\mc E}_{p_1 \cdots p_m; p'}
\eeqn
which can be translated to orbits of $\mb\Pi_1\times \cdots \times \mb\Pi_m \times \Gamma$. 

Now we set up an induction scheme. For each tuple $(p_1, \ldots, p_m; p')$, define
\beqn
|p_1 \cdots p_m;p'| = |{\mc A}_1|(p_1) + \cdots + |{\mc A}_m|(p_m) - |{\mc A}'|(p')\in {\mb R}.
\eeqn

\begin{lemma}
There exists a sequences of real numbers
\beqn
c_1 < c_2 < \cdots
\eeqn
satisfying the following conditions.
\begin{enumerate}
    \item $c_l$ goes to infinity.

    \item If $|p_1 \cdots p_m; p'| \in [c_{l-1}, c_l)$ with $p_i < q_i$ in $\tilde{\mb F}_i$ such that $M_{p_1 \cdots p_{i-1} q_i p_{i+1} \cdots p_m; p'}^{\tilde{\mb X}} \neq \emptyset$, then $|p_1 \cdots p_{i-1} q_i p_{i+1} \cdots p_m; p'| < c_{l-1}$.

    \item If $|p_1 \cdots p_m; p'| \in [c_{l-1}, c_l)$ with $q' < p'$ in $\tilde{\mb F}'$ such that $M_{p_1 \cdots p_m; q'}^{\tilde{\mb X}} \neq \emptyset$, then $|p_1 \cdots p_m; q'| < c_{l-1}$.

    \item If $M_{p_1\cdots p_m; p'}^{\tilde{\mb X}} \neq \emptyset$, then $|p_1 \cdots p_m; p'| \geq c_1$. 
\end{enumerate}
\end{lemma}

\begin{proof}
It follows from the definition of equivariant Novikov flow categories (Definition \ref{defn_equivariant_flow_category}) that for each $i$, one has 
\beqn
\inf \Big\{ |{\mc A}_i| (p_i) - |{\mc A}_i|(q_i)\ |\ M_{p_iq_i}^{\tilde {\mb F}_i}\neq \emptyset \Big\} > 0.
\eeqn
Similarly
\beqn
\inf \Big\{ |{\mc A}'|(p') - |{\mc A}'|(q')  \ |\ M_{p'q'}^{\tilde {\mb F}'} \neq \emptyset \Big\} > 0.
\eeqn
Let $\epsilon$ be the minimum of these lower bounds. Then we can take $c_l = l \epsilon$. 
\end{proof}

Now we state the induction hypothesis: for a $c_l$ in this sequence, we have constructed FOP transverse perturbations
\beqn
{\mc S}_{p_1 \cdots p_m;p'}': {\mc U}_{p_1\cdots p_m; p'} \to {\mc E}_{p_1 \cdots p_m; p'}
\eeqn
for all tuples $(p_1, \ldots, p_m; p')$ satisfying $|p_1 \cdots p_m; p'| < c_l$ which satisfy the following conditions. 
\begin{enumerate}
\item The perturbations are $\mb\Pi_1 \times \cdots \times \mb\Pi_m$-equivariant. 

\item Whenever $p_i < q_i$ in $\tilde {\mb F}_i$, the morphism 
\beqn
\iota_{p_1 \cdots p_{i-1} (p_i q_i) p_{i+1} \cdots p_m; p'}^{\tilde {\mb X}}: M_{p_i q_i}^{\tilde {\mb F}_i} \times M_{p_1 \cdots p_{i-1} q_i p_{i+1} \cdots p_m; p'}^{\tilde {\mb X}} \to M_{p_1 \cdots p_m; p'}^{\tilde {\mb X}}
\eeqn
is also a morphism in $\outer \uds{\bf dOrb}_{\rm rig}^{\rm FOP}$. 

\item Whenever $q' < p'$ in $\tilde {\mb F}'$, the morphism
\beqn
\iota_{p_1 \cdots p_m; q'p'}^{\tilde {\mb X}}: M_{p_1 \cdots p_m; q'}^{\tilde {\mb X}} \times M_{q'p'}^{\tilde{\mb F}'} \to M_{p_1 \cdots p_m; p'}^{\tilde {\mb X}}
\eeqn
is also a morphism in $\outer \uds{\bf dOrb}_{\rm rig}^{\rm FOP}$.

\item When $(p_1, \ldots, p_m; p')$ is minimal, ${\mc S}_{p_1 \cdots p_m; p'}'$ is the one we have already constructed before the induction.
\end{enumerate}

Our choice of $c_1$ implies that the base case is established. Suppose the induction hypothesis holds for $c_{l-1}$. Consider a tuple $(p_1, \ldots, p_m;p')$ with $|p_1 \cdots p_m; p'| \in [c_{l-1}, c_l)$. If it is minimal, then the perturbation has been constructed. If not, then the induction hypothesis and the conditions on the sequence $c_l$ implies that an FOP transverse perturbation has been constructed compatibly on each lower stratum of ${\mc U}_{p_1\cdots p_m; p'}$. One can extend it using the collar structure to a neighborhood of the boundary. Then by (3) of Theorem \ref{thm_FOP_property}, one can extend to an FOP transverse perturbation ${\mc S}_{p_1\cdots p_m; p'}': {\mc U}_{p_1 \cdots p_m; p'}\to {\mc E}_{p_1\cdots p_m; p'}$.

Using the free $\mb\Pi_1\times \cdots \times \mb\Pi_m$-action, one can translate this perturbation to all tuples on the same orbit. Apply the same procedure to all orbits of tuples within this energy window. Therefore the induction can be carried out. Hence the multimodule case of Theorem \ref{thma_FOP} is established. 

The homotopy case of Theorem \ref{thma_FOP} is almost identical to the case of multimodules, hence is omitted.

\end{proof}

\subsection{FOP perturbation and multimodule concatenations}

Recall that although multimodule concatenation can be canonically done for multimodules enriched in $\uds{\bf dOrb}$, in practice we will have non-canonical concatenations (involving extra stabilizations).

\begin{prop}\label{prop_FOP_concatenation}
Let $\mathring {\mb F}_1, \ldots, \mathring {\mb F}_m, \mathring {\mb G}_1, \ldots, \mathring {\mb G}_n, \mathring {\mb G}'$ be Novikov flow categories enriched in $\outer \uds{\bf dOrb}_{\rm rig}^{\rm FOP}$. Let $\mathring {\mb X}$ be a Novikov multimodule over $(\mathring {\mb F}_1, \ldots, \mathring {\mb F}_m; \mathring {\mb G}_i)$ and $\mathring {\mb Y}$ be a Novikov multimodule over $(\mathring {\mb G}_1, \ldots, \mathring {\mb G}_n; \mathring {\mb G}')$. Let $\tilde {\mb X}$ resp. $\tilde {\mb Y}$ be their descents to $\outer \uds{\bf dOrb}_{\rm rig}^{\rm NC}$. Suppose $\tilde {\mb X} \circ_i \tilde {\mb Y}$ is a concatenation of $\tilde {\mb X}$ and $\tilde {\mb Y}$ at $\tilde {\mb G}_i$. Then this concatenation has a lift to $\outer \uds{\bf dOrb}_{\rm rig}^{\rm NC}$ as a concatenation $\mathring {\mb X} \circ_i \mathring {\mb Y}$. Moreover, the concatenation respects the functor
\beqn
\outer \uds{\bf dOrb}_{\rm rig}^{\rm FOP} \to \outer \pman.
\eeqn
\end{prop}

\begin{proof}
By Definition \ref{defn_general_concatenation}, the concatenation $\tilde {\mb X} \circ_i \tilde {\mb Y}$ contains codimension zero morphisms
\beq\label{FOP_concatenation_morphism}
M_{p_1 \cdots p_m; r_i}^{\tilde{\mb X}} \times M_{r_1 \cdots r_n; r'}^{\tilde {\mb Y}} \to M_{r_1 \cdots r_{i-1} p_1 \cdots p_m r_{i+1} \cdots r_n; r'}^{\tilde {\mb X} \circ_i \tilde {\mb Y}}
\eeq
where the three items above are objects of $\outer \uds{\bf dOrb}_{\rm rig}^{\rm NC}$ and the morphism is a collared rigidified embedding of NC derived orbifolds. Then by the {\bf (Product Property)} of Theorem \ref{thm_FOP_property}, the product of existing FOP perturbations on $M_{p_1 \cdots p_m; r_i}^{\tilde {\mb X}}$ and on $M_{r_1 \cdots r_n; r'}^{\tilde {\mb Y}}$ is an FOP perturbation on the source of the morphism \eqref{FOP_concatenation_morphism}. Then after certain stabilization by a rigid orbifold vector bundle, by the {\bf (Stabilization Property)} of Theorem \ref{thm_FOP_property}, the product has a canonical germ of FOP transverse extensions to the (corresponding open stratum) of the target of \eqref{FOP_concatenation_morphism}. The condition that the FOP perturbations are collared and respect the rigidifications guarantees that the extensions agree on adjacent strata. Therefore, one obtains a concatenation $\mathring {\mb X} \circ_i \mathring {\mb Y}$ lifting $\tilde{\mb X} \circ_i \tilde {\mb Y}$. It is easy to verify that the concatenation is compatible with the concatenation of multimodules enriched in $\pman$ and the forgetful functor $\outer\uds{\bf dOrb}_{\rm rig}^{\rm FOP} \to \outer \pman$.  
\end{proof}

\newpage

\part{GEOMETRIC FLOW CATEGORIES, BIMODULES, AND HOMOTOPIES}

\section*{Outline of Part 2}

In this part, we describe flow categories, bimodules, and homotopies that will appear in (most of) the settings of Hamiltonian Floer theory. We then state the main theorems about the construction of their lifts to the category of normally complex Kuranishi spaces, which induce lifts to normally complex derived orbifolds. Using the the counting scheme described in Section \ref{sec:FOP-perturb}, we prove the main theorems except for those in the $\zp$-equivariant setting. The technical constructions of the regularization of the moduli spaces are deferred to the next Part.

\section{Geometric Flow Categories}
\label{subsection_example}

\subsection{Morse flow category}

Given a closed oriented Riemannian manifold $(M, g)$ and a Morse function $f: M \to {\mb R}$, the Morse flow category, denoted by 
\beqn
\mb{M}(f, g)
\eeqn
and usually abbreviated as $\mb{M}$, has critical points of $f$ as its objects and the moduli spaces of unparametrized (broken) downward gradient trajectories as its morphism spaces. We always assume that $(f, g)$ is Morse--Smale and that near critical points $f$ is a quadratic function and $g$ is the Euclidean metric with respect to a certain coordinate system near each critical point. In this way, for any two critical points $p,q\in {\rm crit} f$, the moduli space $\ov{\mc M}{}_{pq}^{\mb M}$ is a smooth manifold with corners \cite{Wehrheim_Morse}.

The partial order of ${\rm Ob}\mb{M}$ is defined as
\beqn
p \leq q \Longleftrightarrow \ov{\mc M}{}_{pq}^{\mb M} \neq \emptyset
\eeqn
and the poset stratifying $\ov{\mc M}{}_{pq}^{\mb M}$ (when $p \neq q$)
\beqn
A_{pq}^{\mb M} = \{ p < r_1 < \cdots < r_k < q\}.
\eeqn
Then $\ov{\mc M}{}_{pq}^{\mb M}$ is an object of $\uds{\bf SMan}$.

We also define composition maps as the natural inclusion 
\beqn
\ov{\mc M}{}_{pr}^{\mb M} \times \ov{\mc M}{}_{rq}^{\mb M} \to \ov{\mc M}{}_{pq}^{\mb M}
\eeqn
which is a morphism of $\uds{\bf SMan}$. Therefore, we obtained a flow category enriched in $\uds{\bf SMan}$. In fact, the standard exponential decay estimate implies that the moduli spaces are indeed objects of $\uds{\bf SMan}_{\rm adm}$. Therefore, its outercollaring still lives in the same category.

The flow category $\mb{M}(f, g)$ admits a grading which assigns to each critical point $p$ the its Morse index. By assigning orientations on unstable manifolds and counting elements of zero-dimensional moduli spaces, one obtains a Morse--Smale--Witten complex $CM_*(f, g)$ over integers freely generated by critical points, whose homology is isomorphic to the singular homology of $M$.

\subsection{Floer flow category}

We now describe a central object considered in this paper. We will assign to each compact symplectic manifold together with a nondegenerate Hamiltonian and an almost complex structure a flow category, called the {\bf Floer flow category}, typically denoted by ${\mb F}$. More precisely, let $(X, \omega)$ be a compact symplectic manifold. Recall
\beqn
\Pi= \frac{ \pi_2(X)}{ {\rm ker} c_1 \cap {\rm ker} \omega}
\eeqn
which is an abelian group (monoid). We make it a Novikov monoid (Definition \ref{defn_Novikov_group}) by defining
\beqn
\lambda: \Pi \to {\mb R},\ \lambda(a) = \omega(a),
\eeqn
and
\beqn
\mu: \Pi \to {\mb Z},\ \mu(a) = 2 c_1(a).
\eeqn

Let $H: M \times S^1 \to X$ be a (time-dependent) $1$-periodic family of Hamiltonians which induces a 1-parameter family of Hamiltonian diffeomorphisms
\beqn
\phi_H^t \in {\rm Ham} (X, \omega).
\eeqn
We assume that $H$ is nondegenerate. Let $\tilde {\mc P}(H)$ be the set of pairs $(u, x)$ where $u: {\mb D}^2 \to X$ is a smooth map and $x = u|_{\partial \mb{D}^2}: S^1 \to X$ is a 1-periodic orbit of $H$. Define
\beqn
(u, x) \sim (u', x') \Longleftrightarrow x = x',\ \int_{{\mb D}^2} u^* \omega = \int_{{\mb D}^2} (u')^* \omega,\ \langle c_1(TX), u({\mb D}^2) \cup (-u' ({\mb D}^2)) \rangle = 0.
\eeqn
Define
\beqn
{\mc P}(H):= \tilde {\mc P}(H)/\sim
\eeqn
whose elements are called {\bf capped 1-periodic orbits}. Notice that there are well-defined action
\beqn
{\mc A}_H: {\mc P}(H) \to {\mb R}
\eeqn
and well-defined Conley--Zehnder index
\beqn
{\rm deg}: {\mc P}(H) \to {\mb Z}.
\eeqn
Moreover, $\Pi$ acts on ${\mc P}(H)$ freely with a finite quotient set, satisfying
\begin{align*}
&\ {\mc A}_H(ap) = {\mc A}_H(p) - \lambda(a),\ &\ {\rm deg}(ap) = {\rm deg} (p) - \mu(a).
\end{align*}

The Floer flow category ${\mb F}$ has object set
\beqn
{\rm Ob}{\mb F} = {\mc P}(H).
\eeqn
Choose an $\omega$-compatible almost complex structure $J$. For each pair $p, q \in {\mc P}(H)$, one has the moduli space $\ov{\mc M}{}_{pq}^{\mb F}$ of stable Floer trajectories. %
A partial order on ${\mc P}(H)$ is induced by:
\beqn
p\leq q \Longleftrightarrow \ov{\mc M}{}_{pq}^{\mb F} \neq \emptyset.
\eeqn
The moduli space $\ov{\mc M}{}_{pq}^{\mb F}$ is stratified by the regular poset $A_{pq}^{\mb F}$ which is (when $p \neq q$)
\beqn
A_{pq}^{\mb F}:= \{ p < r_1 < \cdots < r_k < q\}.
\eeqn
This provides a flow category which is denoted by $\mb{F}(X, \omega, H, J)$, often abbreviated by ${\mb F}$. At this moment this flow category is only enriched in $\uds{\bf Top}$.

We summarize the above setup in the following situation.

\begin{situationf}\label{situationf1}
$(X, \omega)$ is a compact symplectic manifold. ${\mb F}$ is the Floer flow category associated to a nondegenerate $1$-periodic Hamiltonian and a 1-periodic family of $\omega$-compatible almost complex structure $J$.
\end{situationf}

We summarize the usual properties of the Floer flow category.

\begin{prop}
In Situation \ref{situationf1}, the  Floer flow category ${\mb F}$, together with the symplectic action ${\mc A}_H$ and the Conley--Zehnder index, is a ${\mb Z}$-graded $\Pi$-equivariant Novikov flow category.
\end{prop}

\subsection{The pearly flow category}\label{subsection_pearly}

To compare Floer homology and Morse homology, we introduce a version of the pearly flow category. {\it A priori}, this flow category induces a chain complex which is not quasi-isomorphic to the Morse chain complex, as least not in the most naive way.

Let $(X, \omega, J)$ be a compact almost K\"ahler manifold. For each $a \in \Pi$, consider the moduli space 
\beqn
{\mc M}{}_{0,2}^{\mb R}(a)
\eeqn
which is a variant of the moduli space ${\mc M}{}_{0,2}(a)$ of $J$-holomorphic maps with two marked points in degree $d$. It is defined to be the moduli space of unparametrized Floer trajectories for the Hamiltonian $H \equiv 0$. By singularity removal, ${\mc M}{}_{0,2}^{\mb R}(a)$ can be alternatively described as the moduli space of $J$-holomorphic spheres of class $a$ with $2$ marked points such that the $2$ marked points are equipped with asymptotic markers with the same angular coordinate.

There is a natural compactification 
\beqn
\ov{\mc M}{}_{0,2}^{\mb R}(a)
\eeqn
which is naturally stratified by (ordered) decompositions $a = a_1 + \cdots + a_k$ with $\omega(a_i)>0$. The strata of the above moduli space $\ov{\mc M}{}_{0,2}^{\mb R}(a)$ can be labelled by the so-called {\it linear trees}. These are oriented trees whose vertices have valence at most 2, hence have a unique incoming vertex (the input) and a unique outgoing vertex (the output). The degrees $a_1, \ldots, a_k\in \Pi$ are viewed as extra decorations on the vertices. We then can allow the edges of the tree to acquire length, allowing the lengths to vary, and allowing the edges to break when the lengths approach infinity. Choose a Morse--Smale pair $(f, g)$ on $X$. We then add additional variables which are negative gradient segments of $(f, g)$ defined on the edges of varying length. For the input sphere resp. output sphere we also add a semi-infinite edge on which is supported a negative gradient ray. The natural matching conditions at nodes are required. By choosing an incoming resp. outgoing critical point $\uds x$ resp. $\uds y$ of $f$, one can consider the pearly moduli space and its compactification 
\beqn
\ov{\mc M}{}_{\uds x\uds y}^{\mb P} (a).
\eeqn

The underlying poset is not as simple as the Morse or Floer case. For $x = (\uds x, a)$ and $y = (\uds x, b)$ in ${\rm crit} f \times \Pi$, define $A_{xy}^{\mb P}$ to be the set of (possibly broken) linear trees   whose breakings are labelled by critical points of $f$, whose vertices are labelled by elements of $\Pi$ (with sum being $a-b$), and whose finite edges have either zero or positive lengths. This is a regular poset, whose maximal elements are those unbroken trees whose edges all have positive length. The depth function is given by the number of breakings plus the number of edges with length zero.

\begin{defn}
The {\bf pearly flow category} $\mb{P}$ associated to $X, J, f, g$ as above, is the flow category enriched in $\uds{\bf Top}$ which is defined as follows. 
\begin{enumerate}
\item The set of objects is ${\rm Ob}\mb{P} = {\rm crit} f \times \Pi$.

\item For any pair of objects $x = (\uds x, a)$ and $y = (\uds y, b)$, the morphism space is the moduli space
\beqn
\ov{\mc M}{}_{\uds x \uds y}^{\mb P} (a-b).
\eeqn

\item The composition map is given by the concatenation of pearly solutions. 
\end{enumerate}
\end{defn}

\begin{rem}
    The pearly flow category $\mb{P}$ should be distinguished from the pearly bimodule introduced in \cite[Section 4.4]{Bai_Xu_Arnold}, which is a flow bimodule between Morse flow categories. $\mb{P}$ here is used to set up the Morse--Bott homology of the constant Hamiltonian $H \equiv 0$.
\end{rem}

\subsection{Coupling of geometric flow categories}

Given symplectic manifolds $(X^j, \omega^j)$, $j = 1, \ldots, k$ and Floer flow categories $\mb{F}^j$ associated to pairs $(H^j, J^j)$, their coupling is defined to be the Floer flow category of the product manifold with the product Hamiltonians and almost complex structure. However, it is difficult to describe this construction on the abstract level for a collection of flow categories. The more abstract packaging would require building a symmetric monoidal structure on the (infinity) category of flow categories, cf. \cite{Abouzaid_Blumberg_2024} and its follow-up work. Here we describe certain ``coupling" that arises geometrically. 

{\bf Running Assumption.} Let $\mb{F}^j$, $j = 1, \ldots, k$ be a collection of flow categories. Each of them is either a Floer, Morse, or pearly flow category on a compact (symplectic) manifold.

We define their {\bf coupling}, denoted as
\beqn
{\mb F}:= {\mb F}^1 \ast \cdots \ast {\mb F}^k,
\eeqn
as follows. Its set of objects is
\beqn
{\rm Ob} {\mb F} = {\rm Ob}{\mb F}^1 \times \cdots \times {\rm Ob}{\mb F}^k
\eeqn
whose elements are written as ${\bf p} = (p^1, \ldots, p^k)$. The morphism space $M_{{\bf p}{\bf q}}$ is the moduli space of solutions ${\bf u} = (u^1, \ldots, u^k)$ where $u^j$ is a (stable) Floer/Morse/pearly trajectory connecting $p^j$ and $q^j$, depending on the geometric origin of $\mb{F}^j$. We take the quotient by the ``synchronized'' ${\mb R}$-action, i.e., for any $s \in {\mb R}$, we identify tuples of parametrized trajectories
\beqn
(u^1 (\cdot), \dots, u^k(\cdot)) \sim (u^1 (\cdot + s), \dots, u^k(\cdot + s)),
\eeqn
and compactify by allowing sphere bubbling and breaking of trajectories. Notice that this formulation means that breakings in each factor happen simultaneously. For example, when all ${\mb F}^j$ are Floer flow categories, this coupling is the same as the Floer flow category of the product manifold with respect to the product Hamiltonian and product almost complex structure. Notice that if ${\mb F}^j$ is equivariant with respect to an abelian group $\Pi^j$, then ${\mb F}$ is equivariant with respect to the natural action by the product $\Pi = \Pi^1 \times \cdots \times \Pi^k$.

In the discussion of equivariant Floer theory, one often needs to take self-coupling. For any $k$, denote
\beqn
{\mb F}^{\ast k}:= \underbrace{ {\mb F} \ast \cdots \ast {\mb F}}_{k},
\eeqn
which is equivariant with respect to the action by the product $\Pi \times \cdots \times \Pi$ of $k$ copies of $\Pi$. On the other hand, ${\mb Z}_k$ acts on this flow category by cyclically permuting the factors. This ${\mb Z}_k$-action commutes with the action by the diagonal $\Pi$. Hence this coupled flow category is equivariant with respect to an action by $\Pi \times {\mb Z}_k$.

\section{Geometric Multimodules}

Next, we describe bimodules and multimodules that lead to continuation maps, PSS and SSP maps, pair-of-pants products, and K\"unneth maps.

\subsection{Floer data and Floer domains}

\begin{defn}\label{defn_admissible_domain}
A {\bf punctured surface} is a compact oriented surface  (not necessarily connected) with finitely many points removed, usually denoted by $\Sigma$. We also include as part of data a decomposition of the punctures into {\bf negative ends} and {\bf positive ends}. The negative ends are usually denoted by $z_-^1, \ldots, z_-^{k_-}$, and the positive ends are usually denoted by $z_+^1, \ldots, z_+^{k_+}$. 

A {\bf cylindrical end} on a punctured surface $\Sigma$ around a negative resp. positive puncture $z_\pm^j$ is an orientation-preserving diffeomorphism between ${\mb R}^\pm \times S^1$ and a punctured neighborhood of $z_\pm^j$; two cylindrical ends are regarded as equivalent if they differ by a translation. An equivalence class is called a {\bf germ of cylindrical ends}.

A {\bf cylindrical surface} is a punctured  surface equipped with germs of cylindrical ends around all punctures. To save notations such an object is still denoted by $\Sigma$.
\end{defn}

\begin{defn}\label{defn_Floer_datum}
Let $(X, \omega)$ be a symplectic manifold. A {\bf Floer datum} over a cylindrical surface $\Sigma$, denoted by $\sigma = (j, J, \alpha)$, consists of a complex structure $j$ on $\Sigma$, a smooth family $J_z$ of $\omega$-compatible almost complex structures on $X$, and a smooth 1-form (called the Hamiltonian connection)
\beqn
\alpha \in \Gamma( X \times \Sigma, T^* \Sigma)
\eeqn
 such that 
\begin{enumerate}

\item For any cylindrical end $\varphi_\pm^j: {\mb R}^\pm \times S^1_t \to \Sigma$  at $z_\pm^j$, it is holomorphic with respect to $j$ near infinity where ${\mb R}^\pm \times S^1$ is equipped with the standard complex structure.

\item Near the infinity of the cylindrical end around $z_\pm^j$, we have $(\varphi_\pm^j)^* \alpha = H_\pm^j dt$ where $H_\pm^j$ is either a nondegenerate 1-periodic Hamiltonian or zero. In the former resp. latter case we say this end is of {\bf Floer type} resp. {\bf Morse type}.

\item Near the infinity of the cylindrical end around $z_\pm^j$, the almost complex structure $J_z$ coincides with a fixed almost complex structure $J_\pm^j$ on $X$.

\end{enumerate}

A cylindrical surface $\Sigma$ equipped with a Floer datum $\sigma$ is called a {\bf Floer domain}, which we still denote by $\Sigma$ for brevity.
\end{defn}

Given a Floer domain $(\Sigma, \sigma = (j, J, \alpha))$ and a smooth map $u: \Sigma \to X$, one can take the form
\beqn
\ov\partial_\sigma u:= (du + X_{\sigma}(u))^{0,1} \in \Omega^{0,1}(\Sigma, u^* TX)
\eeqn
where $X_{\sigma}$ is defined via 
\beqn
\omega(X_{\sigma}, -) = d_{X} \alpha,
\eeqn
where $d_{X}$ is the differentiation along the $X$ direction, and $(\cdot)^{0,1}$ is the projection onto the $(0, 1)$-part determined by the domain complex structure $j$ and the target almost complex structure $J$.

\subsection{Multimodules from smooth domains}

Now we can define a particular kind of multimodules coming from Floer domains. Let $\Sigma$ be a connected Floer domain with $m$ negative ends and $n$ positive ends such that all cylindrical ends are of Floer type (see Definition \ref{defn_Floer_datum}), with limiting Hamiltonians $H_1, \ldots, H_m$ on the $m$ negative ends and $H_1', \ldots, H_n'$ on the positive ends. Let $J_1,\ldots, J_m$ be the limiting almost complex structures on the negative ends and $J_1', \ldots, J_n'$ be the limiting almost complex structures on the positive ends. Let ${\mb F}_1, \ldots, {\mb F}_m$ resp. ${\mb F}_1', \ldots, {\mb F}_n'$ be Floer flow category associated to $(H_1, J_1), \ldots, (H_m, J_m)$ resp. $(H_1', J_1'), \ldots, (H_n', J_n')$.  Then the Floer domain $\Sigma$ defines a multimodule ${\mb X}^\Sigma$ over $({\mb F}_1, \ldots, {\mb F}_m; {\mb F}_1', \ldots, {\mb F}_n')$. More precisely, given any tuple ${\mf p} = (p_1, \ldots, p_m)$ of objects of ${\mb F}_1, \ldots, {\mb F}_m$ and ${\mf p}' = (p_1', \ldots, p_n')$ of objects of ${\mb F}_1', \ldots, {\mb F}_n'$, consider the space of solutions to 
\beqn
u: \Sigma \to X,\ \ov\partial_\sigma u = 0
\eeqn
which are asymptotic to $p_i$ near the $i$-th negative end and to $p_j'$ near the $j$-th positive end. We compactify this space by allowing sphere bubbles and Floer trajectory breakings, and denote the compactified space by 
\beqn
M_{{\mf p};{\mf p}'}^{{\mb X}^\Sigma}.
\eeqn
Then there are natural structural maps
\beqn
M_{p_i q_i}^{{\mb F}_i} \times M_{p_1\cdots p_{i-1} q_i p_{i+1} \cdots p_m; {\mf p}'}^{{\mb X}^\Sigma} \to M_{{\mf p};{\mf p}'}^{{\mb X}^\Sigma}
\eeqn
and
\beqn
M_{{\mf p}; p_1' \cdots p_{j-1}' q_j' p_{j+1}' \cdots p_n'}^{{\mb X}^\Sigma} \times M_{q_j' p_j'}^{{\mb F}_j'} \to M_{{\mf p};{\mf p}'}^{{\mb X}^\Sigma}.
\eeqn
It is straightforward to see that the collection of moduli spaces $M_{{\mf p}; {\mf p}'}^{{\mb X}^\Sigma}$ and the above maps satisfy the requirement for a multimodule over $({\mb F}_1, \ldots, {\mb F}_m; {\mb F}_1', \ldots, {\mb F}_n')$. Moreover, both ${\mb F}_i$ and ${\mb F}_j'$ have a free $\Pi$-action, hence $\prod {\mb F}_i$ resp. $\prod {\mb F}_j'$ has a free action by $\Pi^m$ resp. $\Pi^n$. Then it is easy to see that the multimodule ${\mb X}^\Sigma$ is equivariant with respect to the diagonal embeddings $\Pi \to \Pi^m$ and $\Pi \to \Pi^n$. We summarize the properties of this multimodule.

\begin{situationm}\label{situationm1}
$(X, \omega)$ is a compact symplectic manifold. ${\mb F}_1, \ldots, {\mb F}_m; {\mb F}'$ \footnote{Here $m$ is allowed to be $0$, the situation of the cigar bimodule described in Subsubsection \ref{subsection_cigar}.} are Floer flow categories associated to pairs $(H_1, J_1), \ldots, (H_m, J_m); (H', J')$. ${\mb X}$ is the multimodule over $({\mb F}_1, \ldots, {\mb F}_m; {\mb F})$ associated to a smooth Floer domain $\Sigma^{\mb X}$ with $m$ negative cylindrical ends and one positive cylindrical end. Moreover, the Floer datum on $\Sigma^{\mb X}$ coincides with $(H_j dt, J_j)$ on the $j$-th negative end and $(H' dt, J')$ on the positive end. 
\end{situationm}

\begin{prop}
In Situation \ref{situationm1}, the multimodule ${\mb X}$ is a $\Pi$-equivariant Novikov multimodule over $({\mb F}_1, \ldots, {\mb F}_m; {\mb F}')$.
\end{prop}

We give examples of a few concrete cases which appear naturally in applications.

\subsubsection{Bimodules for continuation maps} The first one is the bimodule coming from the continuation map construction. Given Floer flow categories $\mb{F}^\pm$ associated to a pair $(H_\pm, J_\pm)$ on the same symplectic manifold $(X, \omega)$, a homotopy of the Floer data
\beqn
(j_0, H_s dt, J_{s, t}), \quad \quad \quad \lim_{s \to \pm \infty} H_s = H_\pm,\ \lim_{s \to \pm} J_{s, t} = J_\pm
\eeqn
over the infinite cylinder, where $j_0$ is the standard complex structure on $\mb{R} \times S^1$, is a special case of Floer data over an arbitrary punctured sphere. Hence one obtains an $({\mb F}^-, {\mb F}^+)$-module, usually called the bimodule of continuation maps. 

There is a special one corresponding to the constant Floer datum on the infinite cylinder. 

\begin{defn}\label{defn_diagonal_bimodule}
The {\bf diagonal bimodule} of a Floer flow category $\mb{F}$ associated to $(H, J)$ is the bimodule for continuation maps associated to the $s$-independent Floer data $(j_0, Hdt, J)$ on the infinite cylinder. It is denoted by $\mb{\Delta}^{\mb{FF}}$.
\end{defn}

\subsubsection{Poincar\'e duality multimodule}

Consider a surface $\Sigma^{\rm PD}$ with two negative cylindrical ends and no positive ends. Let $\sigma^{\rm PD}$ be a Floer datum on $\Sigma^{\rm PD}$ which converges to $H_1 dt$ and $H_2 dt$ on the two negative ends. Let ${\mb F}_1$ and ${\mb F}_2$ be the Hamiltonian Floer flow categories associated to $H_1$ and $H_2$. On the other hand, consider the trivial flow category ${\mb O}$ whose objects are elements of $\Pi$ and the only morphisms are identities. Then we describe a multimodule ${\mb X}^{\rm PD}$ over $({\mb F}_1, {\mb F}_2; {\mb O})$. Given objects $p_1\in {\rm Ob}{\mb F}_1$, $p_2 \in {\rm Ob}{\mb F}_2$, and $a \in \Pi \cong {\rm Ob} {\mb O}$, consider the moduli spaces of maps $u: \Sigma^{\rm PD} \to X$ solving $\ov\partial_{\sigma^{\rm PD}} u= 0$ which converges to the (uncapped) orbits $\uds p_1$ and $\uds p_2$ at the two negative ends; moreover, the spherical class obtained by adding the cappings has symplectic area $a\in \Pi$. We compactify the space by allowing sphere bubbles as well as Floer breakings at the two negative ends, and denote the compactified space by 
\beqn
M_{p_1 p_2; a}^{{\mb X}^{\rm PD}}.
\eeqn
Then there are natural maps
\beqn
M_{p_1 q_1}^{{\mb F}_1} \times M_{q_1 p_2; a}^{{\mb X}^{\rm PD}} \to M_{p_1 p_2; a}^{{\mb X}^{\rm PD}}
\eeqn
and 
\beqn
M_{p_2 q_2}^{{\mb F}_2} \times M_{p_1 q_2; a}^{\mb{X}^{\rm PD}} \to M_{p_1 p_2; a}^{{\mb X}^{\rm PD}}
\eeqn
making it a multimodule over $({\mb F}_1, {\mb F}_2; {\mb O})$. We call it a {\bf Poincar\'e duality multimodule}.

\subsection{PSS type bimodules}\label{subsection_PSS}

In order to establish the relation between the Floer homology and the integral homology of the symplectic manifold, one typically considers the PSS (Piunikhin--Salamon--Schwarz \cite{PSS}) maps. In the current situation, we need to consider a few similar constructions. In terms of the formal terms, we describe bimodules among the three flow categories: Floer, Morse, and Pearly flow categories.  

\subsubsection{Deformed Morse-to-Morse bimodule}

Let $(f, g)$ be a Morse--Smale pair on the almost K\"ahler manifold $(X, \omega, J)$. Let $\uds{\mb M}$ be the Morse flow category associated to $(f, g)$ whose objects are critical points of $f$. Let ${\mb M}$ be the trivial product of $\uds{\mb M}$ with $\Pi$. Instead of the identity bimodule from ${\mb M}$ to itself, we consider a bimodule involving holomorphic spheres. Indeed, for each $a \in \Pi$, consider the space of $J$-holomorphic maps $u: S^2 \to X$ representing the class $a$ with two fixed marked points $z_-, z_+$. Denote its Gromov compactification by $\ov{\mc M}{}_{0,2}^{\rm parametrized} (X, J, a)$, consisting of genus zero stable $J$-holomorphic spheres with two marked points with energy $a$ with a rigid component. Let the two markings be $z_-$ and $z_+$. There are two evaluation maps
\beqn
\ev_\pm: \ov{\mc M}{}_{0,2}^{\rm parametrized} (X, J, a) \to X.
\eeqn
Then for each pair of ojects $x = (a_x, \uds x), y = (b_y, \uds y) \in {\rm Ob}{\mb M}$, define
\beqn
M_{x;y}^{\mb{MM}}:= \ev_-^{-1}( \ov{W_{\uds x}^{\rm us}}) \cap \ev_+^{-1}( \ov{W_{\uds y}^{\rm st}}) \cap \ov{\mc M}{}_{0,2}^{\rm parametrized} (X, J, a_x - b_y).
\eeqn
There are obvious structural maps coming from breaking of Morse trajectories:
\begin{align*}
    &\ M_{xz}^{\mb M} \times M_{z;y}^{\mb{MM}} \to M_{x;y}^{\mb{MM}},\ &\ M_{x;w}^{\mb{MM}}\times M_{wy}^{\mb{M}} \to M_{x;y}^{\mb{MM}}
\end{align*}
forming a bimodule over $({\mb M}, {\mb M})$. We call it the {\bf deformed Morse-to-Morse bimodule}, denoted by $\mb{\Delta}^{\mb{MM}}$.

\subsubsection{Floer-to-Morse and Morse-to-Floer bimodules}

We describe two bimodules relating the Morse flow category ${\mb M}$ and the Floer flow category ${\mb F}$, which will be denoted by ${\mb B}^{\mb{MF}}$ and $\mb{B}^{\mb{FM}}$. Fix a compact almost complex manifold $(X, \omega, J)$. Let $\uds{\mb M}$ be the Morse flow category associated to a Morse--Smale pair $(f, g)$ on $X$ whose objects are critical points of $f$. Let ${\mb M}$ be the trivial product of $\uds{\mb M}$ with $\Pi$. Let ${\mb F}$ be the Floer flow category associated to a nondegenerate Hamiltonian $H$. 

Let $\Sigma^{\mb{MF}}$ be ${\mb C}$ with $\infty$ viewed as a positive end with canonical cylindrical coordinates. Choose a Hamiltonian connection $\sigma^{\mb{MF}}$ which is equal to $H dt$ near $\infty$ and which is zero near $0$. Then for each object $p \in {\rm Ob}{\mb F}$, one can consider the equation 
\beqn
u: \Sigma^{\mb{MF}} \to X,\ \ov\partial_{\sigma^{\mb{MF}}} u = 0,
\eeqn
such that $u$ is asymptotic to the 1-periodic orbit underlying $p$ and the map $u$ represents the capped orbit $p$. Compactify this space by allowing sphere bubblings and Floer trajectory breakings on the positive ends. Denote the compactified moduli space by $\ov{\mc M}{}_p^{\mb{MF}}$. There is an evaluation map 
\beqn
\ev_0: \ov{\mc M}{}_p^{\mb{MF}} \to X
\eeqn
at the point $0 \in \Sigma^{\mb{MF}}$. Then for each pair of objects $x \in {\rm Ob}{\mb M}$, $p \in {\rm Ob}{\mb F}$, if $x = (a_x, \uds x)$ where $a_x \in \Pi$ and $\uds x \in {\rm crit} f$, then define
\beqn
M_{x; p}^{\mb{MF}}:= \ev_0^{-1}(\ov{W_{\uds x}^{\rm us}}) \cap \ov{\mc M}{}_{(-a_x)\cdot p}^{\mb{MF}}.
\eeqn
Here $W_{\uds x}^{\rm us}$ is the unstable manifold of $\uds x$ and $\ov{W_{\uds x}^{\rm us}}$ is its compactification. Then there are natural structural maps
\begin{align*}
&\ M_{xy}^{\mb{M}} \times M_{y;p}^{\mb{MF}} \to M_{x; p}^{\mb{MF}},\ &\ M_{x; q}^{\mb{MF}} \times M_{qp}^{\mb{F}} \to M_{x; p}^{\mb{MF}}
\end{align*}
defining a bimodule $\mb{B}^{\mb{MF}}$ over $({\mb M}; {\mb F})$.

Similarly, if we let $\Sigma^{\mb{FM}}$ be ${\mb C}$ with $\infty$ viewed as a negative end and if we choose a Hamiltonian connection $\sigma^{\mb{FM}}$ on this domain which is equal to $H dt$ on the cylindrical end, then we can define moduli spaces $M_{p;x}^{\mb{FM}}$ similarly, giving a bimodule ${\mb B}^{\mb{FM}}$ over $({\mb F};  {\mb M})$.

\subsubsection{The Pearly-to-Floer and Floer-to-Pearly bimodules}

Suppose we have fixed data for the Pearly flow category ${\mb P}$ and the Floer flow category on the almost K\"ahler manifold $(X, \omega, J)$. In particular, one has a Morse--Smale pair $(f, g)$. Then for the moduli space of solutions to 
\beqn
\ov\partial_{\sigma^{\mb{MF}}} u = 0,\ u: \Sigma^{\mb{MF}} \to X
\eeqn
we can take the ``pearly compactification'' in the same way as defining the pearly flow category. Indeed, we treat $0 \in \Sigma^{\mb{MF}}$ as a special marked point around which there is the standard cylindrical coordinates $(s, t)$. When compactifying the moduli space, sphere bubbles at $0$ are taken only modulo translation in $s$. Then one can insert gradient segments to obtain the corresponding pearly compactification. For each pair of objects $x \in {\rm Ob}{\mb P} = {\rm Ob}{\mb F}$, $p \in {\rm Ob}{\mb F}$, there is the moduli space $M_{x;p}^{\mb{PF}}$ together with structural maps
\begin{align*}
&\ M_{xy}^{\mb{P}}\times M_{y;p}^{\mb{PF}} \to M_{x;p}^{\mb{PF}},\ &\ M_{x; q}^{\mb{PF}} \times M_{qp}^{\mb{F}} \to M_{x;p}^{\mb{PF}}
\end{align*}
giving a bimodule ${\mb B}^{\mb{PF}}$ over $({\mb P}, {\mb F})$. Similarly, one also has the Floer-to-Pearly bimodule ${\mb B}^{\mb{FP}}$.

\subsubsection{Morse-to-Pearly and Pearly-to-Morse bimodules}

We describe a bimodule $\mb{B}^{\mb{PM}}$ over $({\mb P}, {\mb M})$ and a bimodule ${\mb{B}}^{\mb{MP}}$ over $({\mb M}, {\mb P})$. For simplicity, we assume ${\mb M}$ and ${\mb P}$ share the same Morse--Smale pair, though we can certainly take different ones. 

Given a pair of objects $x = (a_x, \uds x) \in {\rm Ob}{\mb P}$ and $y = (b_y, \uds y) \in {\rm Ob}{\mb M}$, consider the moduli space of $J$-holomorphic maps $u: {\mb R}\times S^1$ of symplectic area $a_x - b_y \in \Pi$ such that 
\beqn
u(-\infty) \in W_{\uds x}^{\rm us},\ u(+\infty) \in W_{\uds y}^{\rm st}.
\eeqn
When compactifying these moduli spaces, the bubbling at negative ends are taken modulo only the translation in ${\mb R}$ while bubbling at positive ends are taken modulo both translation and rotations; moreover, we also insert gradient segments at nodes connecting $-\infty$ and the principal component. In this way, one obtains a compact stratified space
\beqn
M_{x;y}^{\mb{PM}}
\eeqn
with structural maps
\begin{align*}
&\ M_{xw}^{\mb{P}} \times M_{w; y}^{\mb{PM}}\to M_{x;y}^{\mb{PM}},\ &\ M_{x; w}^{\mb{PM}}\times M_{wy}^{\mb{M}} \to M_{x;y}^{\mb{PM}}
\end{align*}
giving a bimodule ${\mb B}^{\mb{PM}}$ over $({\mb P}, {\mb M})$. Similarly one has a bimodule ${\mb B}^{\mb{MP}}$ over $({\mb M}, {\mb P})$ by taking ``half pearly'' compactifications at the positive end instead of the negative one.

\subsubsection{Pearly-to-Pearly bimodule}

\begin{figure}[h]
    \centering
    \includegraphics[width=0.9\linewidth]{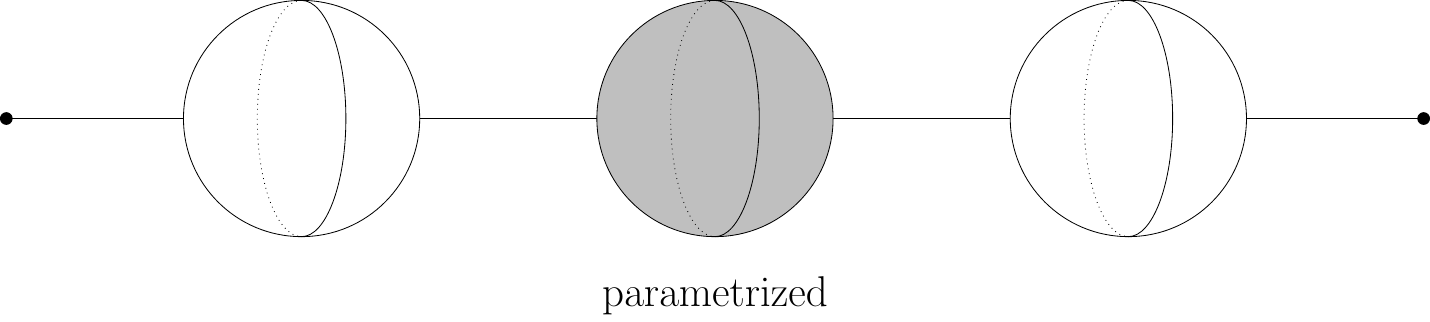}
    \caption{Configurations in $\mb{\Delta}^{\mb{PP}}$.}
    \label{Figure_PP}
\end{figure}

We describe a bimodule $\mb{\Delta}^{\mb{PP}}$ over $(\mb{P},\mb{P})$ where $\mb{P}$ is associated to a Morse--Smale pair $(f, g)$ on $(X, \omega, J)$. Given a pair of critical points $\uds x, \uds y \in {\rm crit} f$ and $a \in \Pi$, consider the space of $J$-holomorphic maps $u: {\mb R}\times S^1 \to X$ with symplectic area $a$ such that $u(-\infty) \in W_{\uds x}^{\rm us}$ and $u(+\infty) \in W_{\uds y}^{\rm st}$. We take the ``pearly compactification,'' namely, when sphere bubbles at $\pm\infty$, the bubbles are taken modulo only the ${\mb R}$-translation. Moreover, we insert (possibly broken) gradient segments at nodes connecting the distinguished marked points $\pm\infty$. Denote the compactified moduli space by  
\beqn
\ov{\mc M}{}_{\uds x;\uds y}^{\mb{PP}}(X, J, a).
\eeqn
A typical configuration of this moduli space is shown in Figure \ref{Figure_PP}.
Then for each pair of objects $x = (\uds x, a_x), y = (\uds y, a_y)$ of $\mb{P}$, define
\beqn
M_{x;y}^{\mb{PP}}:=  \ov{\mc M}_{\uds x\uds y}^{\mb{PP}}(X, J, a_x - a_y).
\eeqn
Then there are natural maps
\begin{align*}
&\ M_{xx'}^{\mb{P}} \times M_{x';y}^{\mb{PP}}\to M_{x;y}^{\mb{PP}},\ &\ M_{x;y'}^{\mb{PP}} \times M_{y' y}^{\mb{P}} \to M_{x;y}^{\mb{PP}}
\end{align*}
defining the structural maps of a bimodule over $(\mb{P},\mb{P})$. We call this bimoulde the {\bf Pearly-to-Pearly bimodule}, denoted by $\mb{\Delta}^{\mb{PP}}$. Notice that there is a ``classical part'' which corresponds to the diagonal bimodule of the Morse flow category.

\subsubsection{Cigar module}\label{subsection_cigar}

To construct a multiplicative identity of the Floer homology, we consider an analogue of the PSS moduli space where we do not have interior markings or Morse trajectories. Let $\Sigma^{\rm cigar} \cong {\mb C}$ be the open Riemann surface with a cylindrical end. Let $\mb{F}$ be the Floer flow category on $(X, \omega, J)$ associated to a nondegenerate Hamiltonian $H$. On the other hand, let ${\mb O}$ be the ``trivial'' flow category consisting of a single object $o \in {\rm Ob}{\mb O}$ with self-morphism space being a single point. 

Choose a Floer datum on $\Sigma^{\rm cigar}$ which coincides with $H$ near the infinity. Let ${\mb F}$ be the Floer flow category associated to $H$. The cigar module is a bimodule ${\mb B}^{\rm cigar}$ over $({\mb O}, {\mb F})$ whose morphism spaces consist of solutions to the Floer equation on $\Sigma$. More precisely, for each object $p \in {\rm Ob}{\mb F}$, consider the space of solutions to 
\beqn
u: \Sigma^{\rm cigar} \to X,\ \ov\partial_{\sigma^{\rm cigar}} u = 0 
\eeqn
such that $u$ converges to the 1-periodic orbit $\uds p$ underlying $p$ and such that $u$ represents the capped orbit $p$. Moreover, compactify this moduli space by allowing sphere bubbling and Floer trajectory breakings. Let the compactified space be $M_{o;p}^{\rm cigar}$. Then there are structural maps
\beqn
M_{o;q}^{\rm cigar} \times M_{qp}^{\mb F} \to M_{o;p}^{\rm cigar},\eeqn
making a bimodule over $({\mb O}; {\mb F})$. 

\subsection{Concatenations}\label{multimodule_concatenation}

We describe the situation corresponding to concatenating Floer data on different domains.

\begin{situationc}\label{situationc1}
$(X, \omega)$ is a compact symplectic manifold. ${\mb F}_1, \ldots, {\mb F}_m, {\mb G}_1, \ldots, {\mb G}_n, {\mb G}'$ are Floer flow categories associated to pairs $(H_1, J_1), \ldots, (H_m, J_m), (G_1, I_1), \ldots, (G_n, I_n), (G', I')$ respectively. ${\mb X}$ resp. ${\mb Y}$ is a multimodule over $({\mb F}_1, \ldots, {\mb F}_m; {\mb G}_i)$ resp. over $({\mb G}_1, \ldots, {\mb G}_n; {\mb G}')$ associated to smooth Floer domain $\Sigma^{\mb X}$ resp. $\Sigma^{\mb Y}$ as described in Situation \ref{situationm1}. 
\end{situationc}

\begin{prop}
In Situation \ref{situationc1}, the concatenation ${\mb X} \circ_i {\mb Y}$ is a $\Pi$-equivariant Novikov multimodule over $({\mb G}_1, \ldots, {\mb G}_{i-1}, {\mb F}_1, \ldots, {\mb F}_m, {\mb G}_{i+1}, \ldots, {\mb G}_n; {\mb G}')$. 
\end{prop}

\begin{proof}
Straightforward.
\end{proof}

\section{Homotopies of Multimodules}\label{section_homotopy}

In this section we descibe several geometrically defined bimodule homotopies. As we only use bimodule homotopies to show certain chain maps are canonically defined up to homotopy, we will not provide the most general treatment. Hence, we only give a minimal descriptions of bimodules which will be used in this paper. 

\subsection{Homotopies from varying Floer data}

Recall the notions of Floer data and Floer domains (Definition \ref{defn_Floer_datum}).

\begin{defn}
A smooth 1-parameter family of Floer domain consists of the following objects. 
\begin{enumerate}

\item A smooth genus zero surface $\Sigma$.

\item A smooth 1-parameter family of domain complex structures $j_\tau$, $\tau \in [0, 1]$, such that there exists an open and precompact region $\Sigma_0 \subset \Sigma$ and diffeomorphisms
\beqn
\phi_i^\pm: {\mb R}_\pm \times S^1 \to U_i^\pm \subset \Sigma \setminus \Sigma_0,\ i = 1, \ldots, l_\pm,
\eeqn
where $U_i^\pm$ are all open subsets of $\Sigma \setminus \Sigma_0$, such that $\phi_i^\pm$ is biholomorphic with respect to all $j_\tau|_{U_i^\pm}$. Each $U_i^\pm$ is called a positive/negative cylindrical end. Two cylindrical ends are germ-equivalent if they differ by a translation.

\item A smooth 1-parameter family of Hamiltonian connections $\sigma_\tau$ on $\Sigma$ such that on each $U_i^\pm$, using the cylindrical coordinates provided by $\phi_i^\pm$, we have $\sigma_\tau = H_i^\pm dt$ near $\infty$, where $H_i^\pm$ is a nondegenerate 1-periodic Hamiltonian on $X$.

\item A smooth 1-parameter family of domain-dependent, $\omega$-compatible almost complex structures $J_\tau$ on $X$ such that on each $U_i^\pm$, using the cylindrical coordinates, $J_\tau = J_i^\pm$ near $\infty$ where $J_i^\pm$ is a fixed almost complex structure on $X$.
\end{enumerate}
Denote such a 1-parameter family of Floer domains by $\Sigma_t^{\mb H}$, where ${\mb H}$ stands for homotopy.
\end{defn}

Given $\Sigma^{\mb H}$, which has the family of incoming Floer flow categories
\beqn
{\mb F}_1, \ldots, {\mb F}_m
\eeqn
and the collection of outgoing Floer flow categories 
\beqn
{\mb F}_1', \ldots, {\mb F}_n'.
\eeqn
For each pair of collections of objects 
\begin{align*}
&\ p_j \in {\rm Ob}{\mb F}_j,\ j = 1, \ldots, m,\ &\ p_i' \in {\rm Ob}{\mb F}_i',\ i=1, \ldots, n
\end{align*}
we consider the moduli space 
\beqn
M_{p_1 \cdots p_m; p_1' \cdots p_n'}^{\mb H}
\eeqn
which is the Gromov--Floer compactification of 
\beqn
\left\{ (\tau, u)\ \left|\ \begin{array}{l} \tau \in [0,1], u: \Sigma \to X,\   \ov\partial_{\sigma_\tau} u = 0,\\
\displaystyle \lim_{s \to -\infty} u|_{U_j} (s, t) =  p_j, j = 1, \ldots, m,\\
\displaystyle \lim_{s \to +\infty} u|_{U_i'} = p_i',\ i =1, \ldots, n  \end{array} \right. \right\}.
\eeqn
This is a compact topological space stratified by a regular poset $A_{p_1 \cdots p_m; p_1' \cdots p_n'}^{\mb H}$. 

Notice that the subset of pairs $(\tau, u)$ with $\tau = 0$ resp. $\tau =1$ defines a multimodule ${\mb X}_0$ resp. ${\mb X}_1$ over $({\mb F}_1, \cdots, {\mb F}_m; {\mb F}_1', \cdots, {\mb F}_n')$. Then there are natural structural maps
\beqn
M_{p_j q_j}^{{\mb F}_j} \times M_{p_1 \cdots p_{j-1} q_j p_{j+1} \cdots p_m; p_1' \cdots p_{i-1}' q_i' p_{i+1}' \cdots p_n'}^{\mb H} \times M_{q_i' p_i'}^{{\mb F}_i'} \to M_{p_1 \cdots p_m; p_1' \cdots p_n'}^{\mb H}
\eeqn
and 
\beqn
M_{p_1 \cdots p_m; p_1' \cdots p_n'}^{{\mb X}_0} \sqcup M_{p_1 \cdots p_m; p_1' \cdots p_n'}^{{\mb X}_1} \to M_{p_1\cdots p_m; p_1' \cdots p_n'}^{\mb H}.
\eeqn
One can then check that ${\mb H}$ together with these structural maps defines a homotopy from ${\mb X}_0$ to ${\mb X}_1$, enriched in $\uds{\bf Top}$.

We will mainly consider the case when $n=1$. We summarize the setting in the following situation.

\begin{situationh}\label{situationh1}
$(X, \omega)$ is a compact symplectic manifold. $\Sigma_t^{\mb H}, t \in [0, 1]$ is a smooth 1-parameter family of Floer domains with $m$ negative cylindrical ends and $1$ positive cylindrical end. ${\mb H}$ is a homotopy from the multimodule ${\mb X}_0$ over $({\mb F}_1, \ldots, {\mb F}_m; {\mb F}')$ defined by $\Sigma_0^{\mb H}$ to the multimodule ${\mb X}_1$ defined by $\Sigma_1^{\mb H}$. 
\end{situationh}

\subsection{Homotopies from gluing}\label{homotopy_gluing}

\begin{figure}[h]
    \centering
    \includegraphics[width=0.5\linewidth]{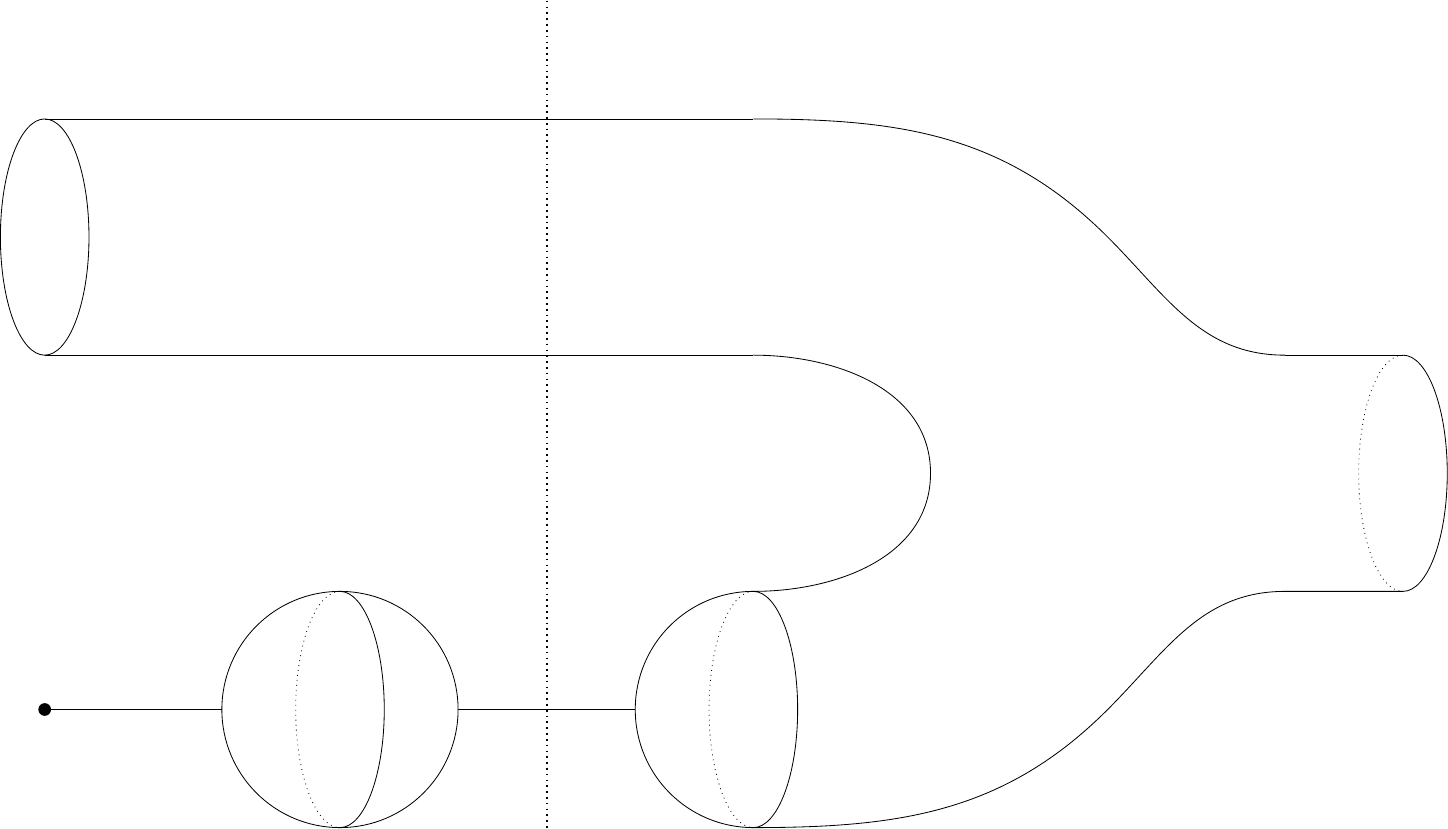}
    \caption{When breaking along the dotted line, the multimodule becomes the concatenation of two.}
    \label{fig:homotopy_gluing}
\end{figure}

Recall the notion of multimodule concatenations. Let $\Sigma^{\mb X}$ and $\Sigma^{\mb Y}$ be Floer domains defining Floer multimodules ${\mb X}$ and ${\mb Y}$ which can be composed to a multimodule ${\mb X} \circ_i {\mb Y}$. We assume that there is only one breaking on the surface $\Sigma^{\mb X} \circ_i \Sigma^{\mb Y}$, see Figure \ref{fig:homotopy_gluing}. Then depending on choosing the cylindrical ends on both sides of that breaking, there is a canonical way to construct a 1-parameter family of Floer domains $\Sigma^{{\mb X} \circ_S {\mb Y}}$ for $S \in [S_0, +\infty)$, where $S$ parametrizes the neck-length, such that when $S \to +\infty$, the Floer domain degenerates to the nodal $\Sigma^{\mb X} \circ_i \Sigma^{\mb Y}$. Denote this family by $\Sigma^{\mb H}$. Choose a smooth function $\tau \mapsto S$ which sends $\tau = 0$ to $S_0$ and $\tau = 1$ to $\infty$. Let $\Sigma^\tau$ be the corresponding Floer domain. For $\tau \in [0, 1)$, let ${\mb X} \circ_{S(\tau)} {\mb Y}$ be the multimodule coming from the Floer domain $\Sigma^\tau$.

Then for any collection of incoming objects ${\mf p}^-$ and outgoing objects ${\mf p}^+$, one has the moduli space $
M_{{\mf p}^-{\mf p}^+}^{\mb H}$, which is the natural compactification of 
\beqn
\Big\{ (\tau, u)\ |\ \tau \in [0, 1],\ u: \Sigma^\tau \to X,\ \ov\partial_{\sigma^\tau} u = 0,\ \lim_{s \to \pm \infty} u|_{U_\pm^i}(s, t) = p_i^\pm \Big\}.
\eeqn

As usual, we have the natural structural maps
\beqn
M_{{\mf p}^- {\mf q}^-}^{{\mf F}^-} \times M_{{\mf q}^- {\mf q}^+}^{\mb H} \times M_{{\mf q}^+ {\mf p}^+}^{{\mf F}^+} \to M_{{\mf p}^- {\mf p}^+}^{{\mb H}}
\eeqn
and 
\beqn
M_{{\mf p}^- {\mf p}^+}^{{\mb X}\circ_i {\mb Y}} \sqcup M_{{\mf p}^- {\mf p}^+}^{{\mb X} \circ_{S_0} {\mb Y}} \to M_{{\mf p}^- {\mf p}^+}^{\mb H}.
\eeqn
These structural maps provide a homotopy ${\mb H}$ from ${\mb X} \circ_{S_0} {\mb Y}$ to ${\mb X} \circ_i {\mb Y}$.

\begin{situationg}\label{situationg1}
$(X, \omega)$ is a compact symplectic manifold. ${\mb X}$ and ${\mb Y}$ are Floer multimodules as described in Situation \ref{situationc1}. $\Sigma_t^{\mb H}$, $t \in [0, 1]$ is a smooth family of Floer domains such that for $t >0$, $\Sigma_t^{\mb H}$ is smooth and $\Sigma_0$ is the concatenation $\Sigma^{\mb X} \# \Sigma^{\mb Y}$. ${\mb M}$ is the Floer multimodule associated to the Floer domain $\Sigma_1^{\mb H}$. ${\mb H}$ is the homotopy between ${\mb X} \circ_i {\mb Y}$ and ${\mb M}$.
\end{situationg}

\subsection{Homotopies associated to PSS type moduli spaces}\label{PSS_glue_homotopy}

On a symplectic manifold, consider the three geometric flow categories: a Floer flow category ${\mb F}$, a Morse flow category ${\mb M}$, and a pearly flow category ${\mb P}$. Consider the following diagram 
\beq\label{big_triangle}
\vcenter{ \xymatrix{    &  &      \mb{F}   \ar@(ur,ul)[]_{{\mb \Delta}^{\mb{FF}}}       \ar@<-.5ex>[lldd]_{{\mb B}^{\mb{FM}}} \ar@<.5ex>[rrdd]^{{\mb B}^{\mb{ FP}}}    &  &  \\
& & & & \\
\mb{M}  \ar@(l,d)[]_{{\mb \Delta}^{\mb{MM}}}       \ar@<-.5ex>[rruu]_{{\mb B}^{\mb{MF}}} \ar@<-.5ex>[rrrr]_{{\mb B}^{\mb{ MP}}}  & & & & \mb{P}    \ar@(d,r)[]_{{\mb \Delta}^{\mb{PP}}}        \ar@<.5ex>[lluu]^{{\mb B}^{\mb{PF}}} \ar@<-.5ex>[llll]_{{\mb B}^{\mb{PM}}}   }}
\eeq
where the arrows are bimodules. Notice that we can compose bimodules in this diagram following the arrows. We say a triangle in this diagram is {\bf commutative} if the composed bimodule from two edges is homotopic to the third bimodule. We will explain 1) why the concatenation ${\mb B}^{\mb{FP}}\circ {\mb B}^{\mb{PF}}$ resp. $\mb{B}^{\mb{FP}} \circ \mb{B}^{\mb{PF}}$ is homotopic to ${\mb \Delta}^{\mb{FF}}$ resp. $\mb{\Delta}^{\mb{PP}}$, whose construction actually can show all triangle whose concatenation happens at ${\mb F}$ or ${\mb P}$ commutes, 2) why the concatenation $\mb{B}^{\mb{MP}} \circ \mb{B}^{\mb{PM}}$ is homotopic to $\mb{\Delta}^{\mb{PP}}$, and 3) why the concatenation ${\mb B}^{\mb{FM}}\circ {\mb B}^{\mb{MF}}$ resp. $\mb{B}^{\mb{PM}} \circ \mb{B}^{\mb{MP}}$ is not (at least not obviously) homotopic to ${\mb \Delta}^{\mb{FF}}$ resp. $\mb{\Delta}^{\mb{MM}}$. 

\subsubsection{Floer-to-Pearly, Pearly-to-Floer bimodules, and their concatenations}

\begin{figure}[h]
    \centering
    \includegraphics[width=0.8\linewidth]{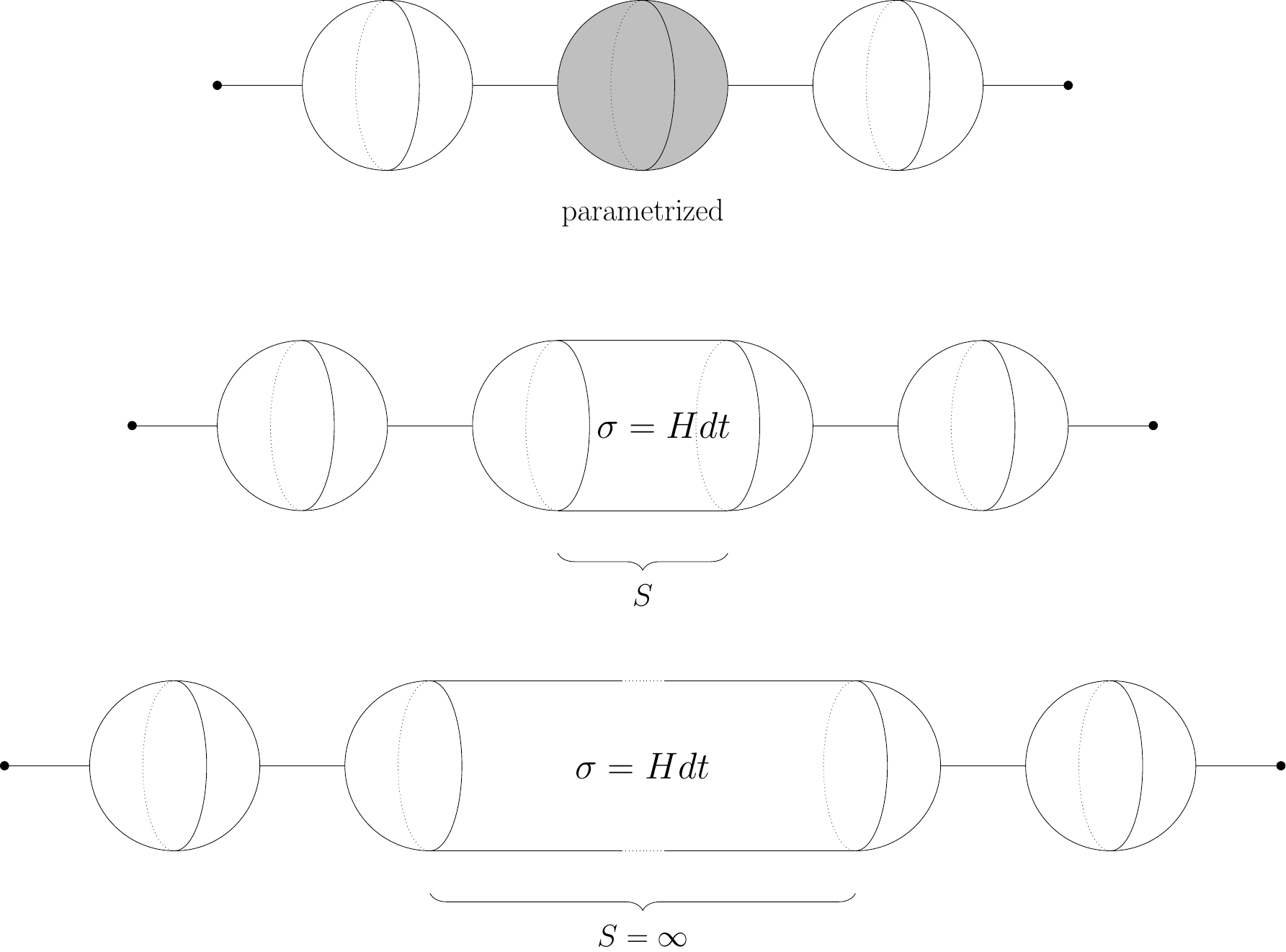}
    \caption{The homotopy from $\mb{\Delta}^{\mb{PP}}$ to $\mb{B}^{\mb{FP}} \circ \mb{B}^{\mb{PF}}$.}
    \label{figure_PFFP}
\end{figure}

Choose a 1-parameter family of Floer data on the infinite cylinder, parametrized by $\tau \in [0, 1)$, given by a Hamiltonian connection $\sigma_\tau$ supported in a region of length about $- 1/\log \tau$ when $\tau$ is close to $1$, and converges to zero when $\tau \to 0$, and such that when $\tau$ is close to $1$, over the long cylinder it is equal to $H dt$ where $H$ is a nondegenerate 1-periodic Hamiltonian. Given a pair $x = (\uds x, a_x), y = (\uds y, a_y)$ of objects of ${\mb P}$, consider the moduli space 
\beqn
\left\{ (\tau, u)\ \left|\ \begin{array}{c} \tau \in [0, 1),\ u: {\mb R}\times S^1\to X,\ \ov\partial_{\sigma_\tau} u = 0,\\ \displaystyle \lim_{s \to -\infty} u(s, t) \in W^{\rm us}_{\uds x},\ \lim_{s \to +\infty} u(s, t) \in W^{\rm st}_{\uds y},\ \omega(u) = a_x - a_y  \end{array} \right. \right\}.
\eeqn
We consider the pearly compactification at $\pm \infty$ and meanwhile allowing $\tau \to 1$, which means the map $u$ should break along Hamiltonian orbits as $- 1/\log \tau \to \infty$ (see Figure \ref{figure_PFFP}, where the neck length $- 1/\log \tau$ is denoted by $S$).

Then the restriction to the $\tau = 1$ slice gives moduli spaces of the bimodule concatenation $\mb{B}^{\mb{FP}} \circ \mb{B}^{\mb{PF}}$ (concatenation at Floer flow category). On the other hand, the restriction to the $\tau = 0$ slice gives moduli spaces of the bimodule $\mb{\Delta}^{\mb{PP}}$. Denote this homotopy by 
\beqn
{\mb H}^{\mb{PF} \tilde\circ \mb{FP}}
\eeqn

The concatenation at the Pearly flow category and its homotopy to $\mb{\Delta}^{\mb{FF}}$ is a bit more complicated. It is the concatenation of two homotopies.

\begin{figure}[h]
    \centering
\includegraphics[width=0.8\linewidth]{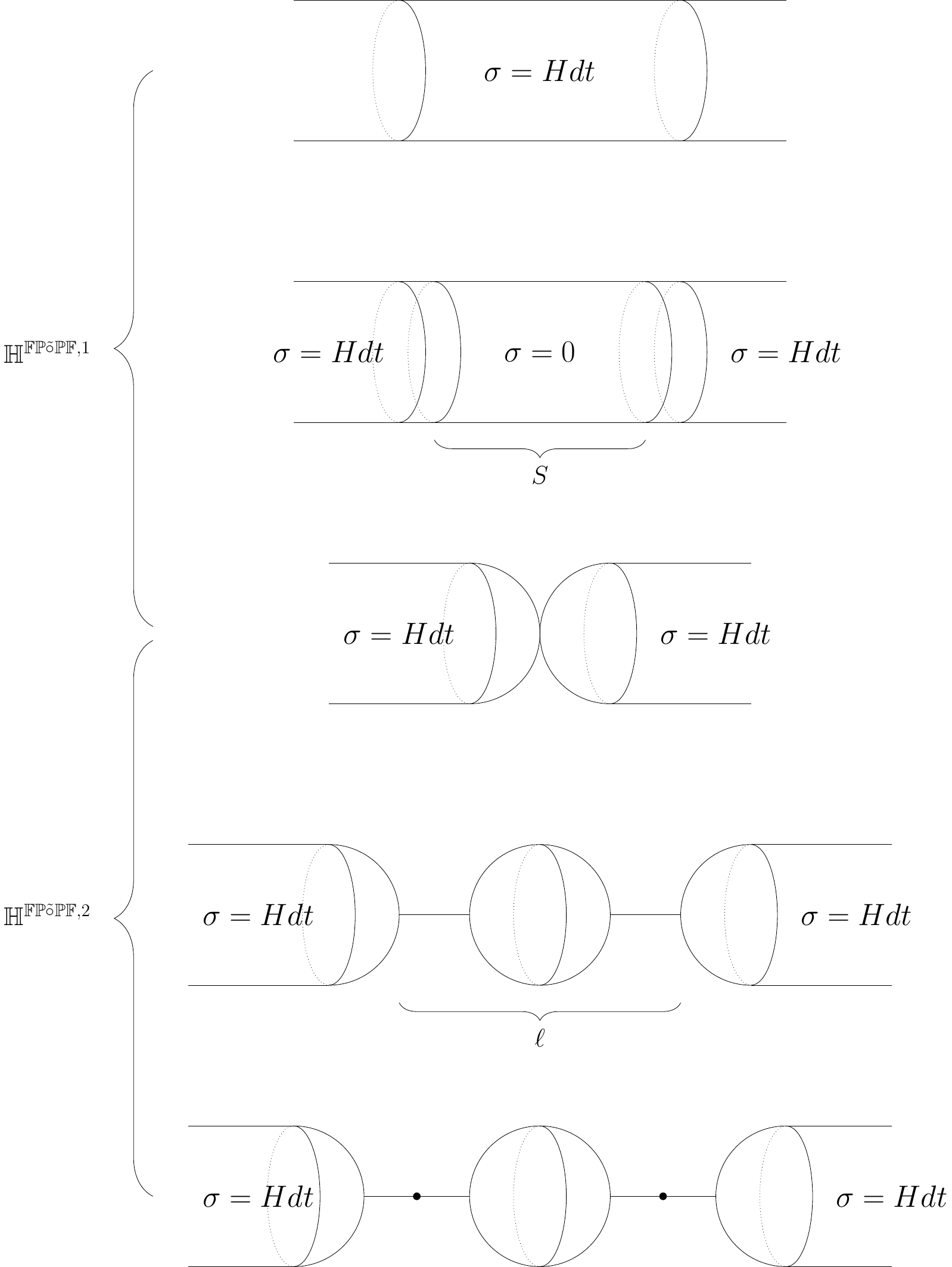}
    \caption{The homotopy between $\mb{\Delta}^{\mb{FF}}$ and $\mb{B}^{\mb{FP}}\circ \mb{B}^{\mb{PF}}$.}
    \label{figure_FPPF}
\end{figure}

First, consider the infinite cylinder equipped with a 1-parameter family of Floer data, parametrized by $\tau \in [0, 1)$ such that when $\tau = 0$, the Hamiltonian connection is the flat one $H dt$, which gives rise to the diagonal bimodule $\mb{\Delta}^{\mb{FF}}$, while as $\tau$ increases, there is a long cylinder of length growing to $\infty$ on which the Hamiltonian connection is trivial. Then for each pair of objects $p, p'\in {\rm Ob}{\mb F}$, consider the space
\beqn
\left\{ (\tau, u)\left| \begin{array}{cc} \tau \in [0, 1),\ u: {\mb R} \times S^1 \to X,\ \ov\partial_{\sigma_\tau} u = 0,\\
\displaystyle \lim_{s \to -\infty} u(s, t) = p,\ \lim_{s \to +\infty} u(s, t) = p' \end{array}  \right. \right\}
\eeqn
We compactify it naturally by adding sphere bubbles, Floer trajectory breakings at both negative and positive ends, and allowing cylindrical breakings when $\tau \to 1$. When the cylinder breaks off, we only quotient by the translation symmetry but not the rotation symmetry. This provides a compact moduli space
\beqn
M_{pp'}^{\mb{FP}\tilde\circ \mb{PF}, 1}
\eeqn
and gives a homotopy from $\mb{\Delta}^{\mb{FF}}$ to an intermediate bimodule, denoted by $\mb{B}^{\mb{F}\bullet \mb{F}}$, denoted by $\mb{H}^{\mb{FP}\tilde\circ \mb{PF}, 1}$ (see Figure \ref{figure_FPPF}).

Next, we introduce a homotopy from $\mb{B}^{\mb{F} \bullet \mb{F}}$ to $\mb{B}^{\mb{FP}\circ \mb{PF}}$.

First, given a nonnegative real number $\ell\geq 0$, we can consider the notion of {\bf telephone} (see Figure \ref{figure_FPPF}) of length $\ell$, which is a pearly-type configuration with two ends being solutions to a Floer equation while the surfaces in between are holomorphic cylinders connected by Morse flow segments of total length $\ell$. By compactifying the moduli space while allowing $\ell \to \infty$ (meaning the Morse trajectories will break one or more times), one obtains moduli spaces 
\beqn
M_{pp'}^{\mb{FP}\tilde\circ\mb{PF},2}
\eeqn
providing a homotopy from $\mb{B}^{\mb{F}\bullet\mb{F}}$ to $\mb{B}^{\mb{FP}}\circ \mb{B}^{\mb{PF}}$, denoted by $\mb{H}^{\mb{FP}\tilde\circ \mb{PF},2}$. Then define
\beqn
\mb{H} = \mb{H}^{\mb{FP}\tilde\circ\mb{PF},1} \cup \mb{H}^{\mb{FP}\tilde\circ\mb{PF},2}
\eeqn
which is a homotopy from $\mb{\Delta}^{\mb{FF}}$ to $\mb{B}^{\mb{FP}}\circ\mb{B}^{\mb{PF}}$. We call it the {\bf telephone homotopy}.

\begin{rem}
    The homotopy $\mb{H}^{\mb{FP}\tilde\circ\mb{PF},2}$ is used to relate two different constructions of the Morse--Bott theory for the constant Hamiltonian $H = 0$. The pearl flow category includes choosing a Morse function over the critical submanifolds and corresponds to the ``cascades model" (see, e.g., \cite{Bourgeois-thesis}) of Morse--Bott theory. In contrast, the natural output of the bimodule $\mb{H}^{\mb{FP}\tilde\circ\mb{PF},1}$ corresponds to the ``fiber product model" (see \cite{Austin-Braam} for an early incarnation). The point of our construction is to compare the ``Morse version" of Floer theory, coming from the flow category $\mb{F}$ that is associated with a nondegenerate Hamiltonian, with the Morse--Bott version, which uses pearly trajectories, using the fiber product model as an intermediate gadget.
\end{rem}

\subsubsection{Other concatenations along Pearly flow categories}

We can also construct a homotopy from $\mb{B}^{\mb{MF}}$ to the concatenation $\mb{B}^{\mb{MP}} \circ \mb{B}^{\mb{PF}}$, denoted by ${\mb H}^{\mb{MP} \tilde\circ \mb{PF}}$ and a homotopy from $\mb{B}^{\mb{FM}}$ to the concatenation $\mb{B}^{\mb{FP}} \circ \mb{B}^{\mb{PM}}$, denoted by $\mb{H}^{\mb{FP}\tilde\circ\mb{PM}}$. Since the two cases are dual to each other, we only explain the definition of the first homotopy. 

Consider the bimodule $\mb{B}^{\mb{MF}}$ given by the collection of moduli spaces
\beqn
M_{xp}^{\mb{MF}}
\eeqn
which is the natural compactification of the space of solutions to the Floer equation on the infinite cylinder with a Hamiltonian connection $\sigma$ which is zero near $-\infty$ and which is $H dt$ near $+\infty$. We deform $\sigma$ to a 1-parameter family $\sigma_\tau$, $\tau \in [0, 1]$ such that when $\tau \to 1$, the support of $\sigma_\tau$ is pushed to $+\infty$. Then as $\tau = 1$, one obtains a bimodule consisting of moduli spaces whose configurations are similar to the third picture of Figure \ref{figure_FPPF}. Then one can continue the homotopy by inserting gradient segments of lengths which vary from $0$ to $+\infty$. The very end of this homotopy is the concatenation $\mb{B}^{\mb{MP}} \circ \mb{B}^{\mb{PF}}$.

\subsubsection{Remarks on concatenation at the Morse flow category}\label{subsubsec:remark}

One can also form the concatenation at the Morse flow category, for example, $\mb{B}^{\mb{FM}} \circ\mb{B}^{\mb{MF}}$. In the category $\uds{\bf Top}$, it is homotopic to the bimodule $\mb{\Delta}^{\mb{FF}}$, firstly by gluing broken Morse trajectories and secondly by gluing holomorphic curves. However, it is hard to regularize such a homotopy. In particular, when gluing holomorphic curves at a node, unlike the Pearly case, the gluing has a complex gluing parameter. This does not fit very well with our framework. This is the reason why we need to go through the Pearly flow category in order to prove the invertibility of the PSS and SSP maps.

\section{Main Theorems on AMS Lifts}

In this section, we state the main results regarding regularizations of moduli spaces introduced in the previous section. The proof will be presented in the next part.

\subsection{Main statements for Floer flow categories}

We state the main theorem about global charts of Floer flow categories and other types of flow categories. The constructions depend on many choices. Most of the dependence of the choices can be absorbed via the global chart construction of bimodules; however, the dependence of the outercollaring width and the choice of the so-called integral action (Definition \ref{defn_integral_modification}) cannot be absorbed easily. Hence in the statement, we emphasize the dependence of these two choices.

\begin{thma}\label{thma_flow_chart}
In Situation \ref{situationf1}, given any outercollaring width $\epsilon>0$ and a choice of integral action on ${\mb F}$ (Definition \ref{defn_integral_modification}) on ${\mb F}$, the $\epsilon$-outercollaring of ${\mb F}$ admits a set of lifts $\tilde {\mb F}$ to $\outer \uds{\bf dOrb}^{\rm NC}_{\rm rig}$ (Definition \ref{defn_rigid_embedding}) as a $\Pi$-equivariant Novikov flow category. Each individual lift in this set is called an {\bf AMS lift} subject to the outercollaring width and the integral action. The AMS lifts admit coherent orientations in the sense of Definition \ref{defn_coherent_orientation}. Moreover, for any $\epsilon'>0$, the $\epsilon'$-outercollaring of $\tilde {\mb F}$ is an AMS lift of ${\mb F}$ subject to the outercollaring width $\epsilon + \epsilon'$ and the same integral action. 
\end{thma}

\begin{proof}
See Section \ref{section_AMS_proofs} except for the last statement, which is a straightforward consequence of the constructions.
\end{proof}

\begin{rem}
In fact what we mean by an AMS lift does not only mean the resulting flow category $\tilde {\mb F}$, but also all the choices made in the process of construction. When an follow-up construction depends on a given AMS lift, it means the construction {\it a priori} depends on the choices made for the AMS lift. 
\end{rem}

The (in)dependence on the choice of the integral action is characterized by the following theorem. Recall that for the Floer flow category ${\mb F}$, there is a bimodule $\mb{\Delta}^{\mb{FF}}$ over $({\mb F}; {\mb F})$ associated to the continuation map from $(J, H)$ to $(J, H)$ with the trivial Floer data.

\begin{prop}\label{thm_integral_action_comparison}
In Situation \ref{situationf1}, let $\tilde {\mb F}_1$ and $\tilde {\mb F}_2$ be two AMS lifts of ${\mb F}$ subject to the same outercollaring width and two different integral actions, provided by Theorem \ref{thma_flow_chart}. Then there exists a class of lifts of the corresponding outercollaring of $\mb{\Delta}^{\mb{FF}}$ to $\outer \uds{\bf dOrb}_{\rm rig}^{\rm NC}$ as a bimodule over $(\tilde {\mb F}_1; \tilde {\mb F}_2)$. Moreover, with respect to the coherent orientations on $\tilde {\mb F}_1, \tilde {\mb F}_2$, the lift is coherently oriented. %
\end{prop}

\begin{proof}
See Section \ref{section_comparison}.
\end{proof}

The (in)dependence on other choices for the AMS lift is included in the theorem on multimodules (Theorem \ref{thma_module_lift}).

\subsection{Multimodules and homotopies}

Next, we state the main theorems concerning regularizing flow multimodules and homotopies.

First, we formulate the theorem on lifts for multimodules. Due to technical reasons, we must consider the dependence on a specific choice, namely, that of integral actions (Definition \ref{defn_integral_modification} and Definition \ref{module_integral_action}) in the construction.

\begin{thma}\label{thma_module_lift}
Suppose ${\mb X}$ is a multimodule associated a smooth Floer domain $\Sigma^{\mb X}$ with $m \geq 0$ negative ends and $n \geq 0$ positive ends satisfying
\beqn
n \in \{0, 1\},\ n = 0 \Longrightarrow m = 2\ {\rm and}\ m = 0 \Longrightarrow n = 1.\footnote{This complicated range is to include the case for the Poincar\'e pairing and the unit of the pair-of-pants product.}
\eeqn
It is a multimodule over $({\mb F}_1, \ldots, {\mb F}_m; {\mb F}')$ where ${\mb F}_1, \ldots, {\mb F}_m$ are Floer flow categories and ${\mb F}'$ is a Floer flow category when $n=1$ and the trivial flow category when $n = 0$. When $m = 0$ and $n=1$, ${\mb X}$ is a cigar bimodule over $({\mb O}; {\mb F}')$.

Choose a common outercollaring width and compatible integral actions on the involved Floer flow categories (see Definition \ref{module_integral_action}). Let $\tilde {\mb F}_1, \ldots, \tilde {\mb F}_m, \tilde {\mb F}'$ be AMS lifts to $\outer \uds{\bf dOrb}^{\rm NC}_{\rm rig}$ of ${\mb F}_1, \ldots, {\mb F}_m, {\mb F}'$ subject to the fixed outercollaring width and these integral actions provided by Theorem \ref{thma_flow_chart} (when ${\mb F}' = {\mb O}$, take $\tilde {\mb F}' = {\mb O}$).

\begin{enumerate}

\item There exists a set of lifts of ${\mb X}$ to $\outer \uds{\bf dOrb}^{\rm NC}_{\rm rig}$, denoted by $\tilde {\mb X}$ (called an AMS lift subject to the chosen integral actions), which is a multimodule over $(\tilde {\mb F}_1, \ldots, \tilde {\mb F}_k; \tilde {\mb F}')$. Moreover, with respect to the coherent orientations on $\tilde {\mb F}_1, \ldots, \tilde {\mb F}_m, \tilde {\mb F}'$, the lift $\tilde {\mb X}$ is coherently oriented.

\item Any two AMS lifts of ${\mb X}$ (subject to the same set of integral actions) are homotopic as multimodules enriched in $\outer \uds{\bf dOrb}_{\rm rig}^{\rm NC}$ via a lift of the the trivial homotopy from the outercollaring of ${\mb X}$ to itself and the homotopy is coherently oriented.
\end{enumerate}
\end{thma}

\begin{proof}
See Section \ref{section_AMS_proofs}.
\end{proof}

Notice that Theorem \ref{thma_module_lift} depends on choosing an integral action. For comparing different choices, we consider a particular homotopy. Recall that $\mb{\Delta}^{{\mb F}_j\mb{F}_j}$, $\mb{\Delta}^{{\mb F}'{\mb F}'}$ are the diagonal bimodules. There is an obvious homotopy ${\mb H}$ from the concatenation (as concatenations of weak bimodules)
\beqn
\Big( \mb{\Delta}^{\mb{F}_1\mb{F}_1} \times \cdots \times \mb{\Delta}^{\mb{F}_m\mb{F}_m} \Big) \circ \mb{X} \circ \mb{\Delta}^{{\mb F}'\mb{F}'}
\eeqn
to ${\mb X}$ by simultaneously gluing at both positive and negative ends. When $n=0$ and $m=2$, we replace the above concatenation by
\beqn
\Big( \mb{\Delta}^{\mb{F}_1\mb{F}_1} \times \mb{\Delta}^{\mb{F}_2\mb{F}_2} \Big) \circ \mb{X}.
\eeqn

\begin{prop}\label{prop:module-double}
Under the assumption of Theorem \ref{thma_module_lift}, fix an outercollaring width. 

\begin{enumerate}

\item When $n = 1$, let
\beqn
{\mc A}_{\mb{F}_j}^{\Omega_{(s)}},\ {\mc A}_{{\mb F}'}^{\Omega_{(s)}},\ s = 0, 1
\eeqn
be two compatible sets of integral actions. Let $\tilde{\mb{X}}_{(s)}$ be an AMS lift subject to the fixed outercollaring width and one set of integral actions. Let $\tilde{\mb{\Delta}}_{(01)}^{\mb{F}_j\mb{F}_j}$, $\tilde{\mb{\Delta}}_{(10)}^{\mb{F}'\mb{F}'}$ be the lift provided by Proposition \ref{thm_integral_action_comparison}. Then there exists a lift of the outercollaring of ${\mb H}$ as a homotopy from 
\beqn
\Big( \tilde{\mb{\Delta}}_{(01)}^{\mb{F}_1\mb{F}_1} \times \cdots \times \tilde{\mb{\Delta}}_{(01)}^{\mb{F}_m\mb{F}_m} \Big) \circ \tilde{\mb X}_{(1)}  \circ \tilde{\mb{\Delta}}_{(10)}^{{\mb F}'\mb{F}'}
\eeqn
to $\tilde {\mb X}_{(0)}$.

\item When $n = 0$, $m = 2$, let 
\beqn
{\mc A}_{{\mb F}_1}^{\Omega_{(s)}},\ {\mc A}_{{\mb F}_2}^{\Omega_{(s)}},\ s = 0, 1
\eeqn
be two compatible pair of integral actions. Let $\tilde{\mb{X}}_{(s)}$ be an AMS lift subject to the fixed outercollaring width and one set of integral actions. Let $\tilde{\mb{\Delta}}_{(01)}^{\mb{F}_j\mb{F}_j}$, $j = 1, 2$ be the lift provided by Proposition \ref{thm_integral_action_comparison}. Then there exists a lift of the outercollaring of ${\mb H}$ as a homotopy from 
\beqn
\Big( \tilde{\mb{\Delta}}_{(01)}^{\mb{F}_1\mb{F}_1} \times \tilde{\mb{\Delta}}_{(01)}^{\mb{F}_2\mb{F}_2} \Big) \circ \tilde{\mb X}_{(1)}  
\eeqn
to $\tilde {\mb X}_{(0)}$.

\end{enumerate}
\end{prop}

\begin{proof}
See Section \ref{section_comparison}.
\end{proof}

Now we state the corresponding results for homotopies of multimodules from varying Floer data on smooth domains.

\begin{thma}\label{thma_homotopy_lift}
In Situation \ref{situationh1}, fix an outercollaring width and a compatible set of integral actions on involved Floer flow categories ${\mb F}_1, \ldots, {\mb F}_m, {\mb F}'$. Let $\tilde {\mb F}_1, \ldots, \tilde {\mb F}_m; \tilde {\mb F}'$ be AMS lifts of these flow categories subject to the fixed outercollaring width and the integral actions. Let $\tilde {\mb X}_0$ resp. $\tilde {\mb X}_1$ be AMS lifts of ${\mb X}_0$ resp. ${\mb X}_1$ subject to the same outercollaring width and the integral actions. Then there exists a lift of the outercollaring of ${\mb H}$ to $\outer \uds{\bf dOrb}_{\rm rig}^{\rm NC}$, denoted by $\tilde {\mb H}$, called an AMS lift subject to the outercollaring width and integral actions as a homotopy from $\tilde {\mb X}_0$ to $\tilde {\mb X}_1$.
\end{thma}

\begin{proof}
See Section \ref{section_AMS_proofs}.
\end{proof}

To establish algebraic relation such as the associativity of pair-of-pants product, one needs to consider the homotopy of multimodules corresponding to gluing of Floer domains. Recall the following setup. Let $\Sigma^{\mb X}$ be a Floer domain with $k$ negative ends and one positive end, defining a multimodule ${\mb X}$ over $({\mb F}_1, \ldots, {\mb F}_m; {\mb G}_i)$, and let $\Sigma^{\mb Y}$ be a Floer domain with $n$ negative ends and one positive ends defining a multimodule ${\mb Y}$ over $({\mb G}_1, \ldots, {\mb G}_n; {\mb G}')$. Here the source and target flow categories are all Floer flow categories associated to the limiting Hamiltonians on the two Floer domains. Recall that the concatenation ${\mb X} \circ_i {\mb Y}$ is a multimodule over $({\mb G}_1, \ldots, {\mb G}_{i-1}, {\mb F}_1, \ldots, {\mb F}_m, {\mb G}_{i+1}, \ldots, {\mb G}_n; {\mb G}')$.

\begin{thm}\label{thm_AMS_concatenation}
Fix a common outercollaring width. In Situation \ref{situationc1}, choose integral actions on the involved flow categories which are compatible with respect to ${\mb X}$ and ${\mb Y}$. Let $\tilde {\mb F}_1, \ldots, \tilde {\mb F}_m, \tilde {\mb G}_1, \ldots, \tilde {\mb G}_n, \tilde {\mb G}'$ be AMS lifts of ${\mb F}_1, \ldots, {\mb F}_m,  {\mb G}_1, \ldots, {\mb G}_n, {\mb G}'$ respectively subject to the fixed outercollaring width and the integral actions. Let $\tilde {\mb X}$ resp. $\tilde {\mb Y}$ be an AMS lift to $\outer \uds{\bf dOrb}_{\rm rig}^{\rm NC}$ subject to the fixed outercollaring width and the chosen integral actions provided by Theorem \ref{thma_module_lift}. Let $\tilde {\mb X}_0$ resp. $\tilde {\mb Y}_0$ be their descendants to $\outer \uds{\bf dOrb}$ which, by the discussion in Subsection \ref{multimodule_concatenation}, admit a natural concatenation $\tilde {\mb X}_0 \circ \tilde {\mb Y}_0$. Then there exists a concatenation $\tilde {\mb X} \circ_i \tilde {\mb Y}$ (see Definition \ref{defn_general_concatenation}) which is a lift of $\tilde {\mb X}_0 \circ_i \tilde {\mb Y}_0$ to $\outer \uds{\bf dOrb}_{\rm rig}^{\rm NC}$.\footnote{For the purpose of this paper, we do not need to prove that the concatenation $\tilde {\mb X} \circ_i \tilde {\mb Y}$ is well-defined up to homotopy.}
\end{thm}

\begin{proof}
See Section \ref{section_gluing}.
\end{proof}

On the other hand, let $\Sigma^{{\mb X} \circ_S {\mb Y}}$ be the family of Floer domains obtained from gluing along positive end of $\Sigma^{\mb X}$ matched with the negative end of $\Sigma^{\mb Y}$, each of which defines a multimodule ${\mb X}\circ_S {\mb Y}$. Recall (see Subsection \ref{homotopy_gluing}) that there is a homotopy ${\mb H}$ from ${\mb X} \circ_i {\mb Y}$ to ${\mb W}:= {\mb X} \circ_{S_0} {\mb Y}$ for $S_0 \gg 0$. 

\begin{thma}\label{thma_multimodule_gluing}
Fix a common outercollaring width. In Situation \ref{situationg1}, choose integral actions on involved flow categories which are compatible with respect to the homotopy ${\mb H}$. Suppose $\tilde{\mb F}_1, \ldots, \tilde {\mb F}_m,  \tilde {\mb G}_1, \ldots, \tilde {\mb G}_l, \tilde {\mb G}'$ are AMS lifts of the involved Floer flow categories subject to the fixed outercollaring width and the integral actions. Let $\tilde{\mb X}$ resp. $\tilde {\mb Y}$ be AMS lifts of the outercollaring of ${\mb X}$ resp. ${\mb Y}$ as multimodules over $(\tilde{\mb F}_1, \ldots, \tilde{\mb F}_m; \tilde{\mb G}_i)$ resp. over $( \tilde {\mb G}_1, \ldots, \tilde {\mb G}_l; \tilde {\mb G}')$ subject to the fixed outercollaring width and the chosen integral actions. Let $\tilde {\mb X} \circ_i \tilde {\mb Y}$ be a concatenation provided by Theorem \ref{thm_AMS_concatenation}. Let $\tilde {\mb M}$ be an AMS lift of the multimodule associated to the smooth domain $\Sigma^{{\mb X} \circ_S {\mb Y}}$ subject to the same outercollaring width and the chosen integral actions. Then there exists a class of lifts of the outercollaring of ${\mb H}$ to $\outer \uds{\bf dOrb}^{\rm NC}_{\rm rig}$ as a homotopy from $\tilde{\mb X} \circ_i \tilde {\mb Y}$ to $\tilde {\mb M}$. Each lift in this class is called an AMS lift of ${\mb H}$.
\end{thma}

\begin{proof}
See Section \ref{section_gluing}.
\end{proof}

\begin{rem}
    We do not discuss comparing non-homotopic choices in the constructions of lifts of homotopies. We only use homotopies to ensure well-defineness of operations induced by flow multimodules, so we stop our discussion on non-homotopic choices at multimodules.
\end{rem}

\subsection{Bimodules and homotopies associated to PSS type maps}

We state the main technical theorem towards the definition of the PSS and SSP maps (cf. Theorem \ref{thma_PSS}). 

\begin{thma}\label{thma_PSS_AMS}
Let ${\mb F}$ be Floer flow category on $(X, \omega)$ associated to a Hamiltonian $H$ and a time-independent almost complex structure $J$. Let ${\mb M}$ be a Morse flow category on $X$ associated to a Morse--Smale pair $(f, g)$ with the standard behavior near critical points. Let $\mb{B}^{\mb{MF}}$ resp. $\mb{B}^{\mb{FM}}$ be the Morse-to-Floer resp. Floer-to-Morse bimodule associated to certain Floer data on the infinite cylinder. Fix an outercollaring width. 

\begin{enumerate}

\item Choose an integral action ${\mc A}_{\mb F}^\Omega: {\rm Ob}{\mb F} \to {\mb Z}$ which is compatible with respect to $\mb{B}^{\mb{MF}}$ resp. to $\mb{B}^{\mb{FM}}$. Let $\tilde {\mb F}$ be an AMS lift of ${\mb F}$ to $\outer \uds{\bf dOrb}_{\rm rig}^{\rm NC}$ subject to the outercollaring width and the integral action. Then there exists a lift of the outercollaring of $\mb{B}^{\mb{MF}}$ resp. $\mb{B}^{\mb{FM}}$, denoted by $\tilde{\mb{B}}^{\mb{MF}}$ resp. $\tilde{\mb{B}}^{\mb{FM}}$ as a bimodule over $(\outer {\mb M}; \tilde {\mb F})$ resp. over $(\tilde {\mb F}; \outer {\mb M})$.

\item Let ${\mc A}_{\mb F}^{\Omega_{(s)}}$, $s = 0, 1$ be two integral actions which are compatible with respect to $\mb{B}^{\mb{MF}}$ resp. $\mb{B}^{\mb{FM}}$. Let $\tilde {\mb F}_{(s)}$ be AMS lifts of ${\mb F}$ subject to ${\mc A}_{\mb F}^{\Omega_{(s)}}$ and let $\tilde {\mb B}_{(s)}^{\mb{MF}}$ resp. $\tilde {\mb B}_{(s)}^{\mb{FM}}$ be AMS lift of $\mb{B}^{\mb{MF}}$ resp. $\mb{B}^{\mb{FM}}$ subject to ${\mc A}_{{\mb F}}^{\Omega_{(s)}}$. Let  $\tilde {\mb{\Delta}}_{(01)}^{\mb{FF}}$ be the AMS lift provided by Proposition \ref{thm_integral_action_comparison}. Then there exists a homotopy from $\tilde{\mb B}_{(1)}^{\mb{MF}}$ to the composition $\tilde{\mb B}_{(0)}^{\mb{MF}} \circ \tilde{\mb{\Delta}}_{(01)}^{\mb{FF}}$ resp. there exists a homotopy from $\tilde{\mb B}_{(0)}^{\mb{FM}}$ to the composition $\tilde{\mb{\Delta}}_{(01)}^{\mb{FF}} \circ \tilde{\mb B}_{(1)}^{\mb{FM}}$.
\end{enumerate}
\end{thma}

\begin{proof}
See Subsection \ref{subsection_Morse_Floer}.
\end{proof}

\begin{thm}\label{thm126}
Let ${\mb F}', {\mb M}', \mb{B}^{\mb{M}'\mb{F}'}$, $\mb{B}^{\mb{F}'\mb{M}'}$ be another set of Floer flow category, Morse flow category, Morse-to-Floer and Floer-to-Morse as in Theorem \ref{thma_PSS_AMS}. Let $\mb{B}^{\mb{FF}'}$ be a bimodule over $(\mb{F}; {\mb F}')$ associated to a Floer datum on the infinite cylinder. Let $\mb{B}^{\mb{MM}'}$ be a bimodule over $(\mb{M}; {\mb M}')$ associated to a one-parameter family of Morse data. Let ${\mb H}$ be a homotopy from $\mb{B}^{\mb{FF}'} \circ \mb{B}^{\mb{F}'\mb{M}'}$ to $\mb{B}^{\mb{FM}} \circ \mb{B}^{\mb{MM'}}$ by gluing the Floer/Morse data and varying in a one-parameter family. Then there exist AMS lifts $\tilde {\mb F}$ resp. $\tilde {\mb F}'$ of ${\mb F}$ resp. ${\mb F}'$ and AMS lifts $\tilde {\mb B}^{\mb{FM}}$ (as a bimodule over $(\tilde{\mb F}; \outer {\mb M})$ resp. $\tilde {\mb B}^{\mb{F}'{\mb M}'}$ (as a bimodule over $(\outer \mb{M}'; \tilde {\mb F}')$, an AMS lift $\tilde {\mb B}^{\mb{FF}'}$, and a lift of the outercollaring of ${\mb H}$ to $\outer \uds{\bf dOrb}_{\rm rig}^{\rm NC}$. 
\end{thm}

\begin{proof}
See Subsection \ref{subsection_Morse_Floer}.
\end{proof}

\subsubsection{PSS maps and pearly flow category}

Next, we consider the pearly flow category and its diagonal bimodule. 

\begin{thm}\label{thma_pearly_chart}
Let $(X, \omega)$ be a compact symplectic manifold. 
\begin{enumerate}

\item Let ${\mb P}$ be the pearly flow category associated to a Morse--Smale pair $(f, g)$ and a compatible almost complex structure $J$. Then upon choosing an integral symplectic form $\Omega$ which is tamed by $J$, any outercollaring of ${\mb P}$ admits a set of lifts $\tilde {\mb P}$ to $\outer \uds{\bf dOrb}_{\rm rig}^{\rm NC}$. Each individual lift is called an {\bf AMS lift} (subject to the outercollaring width and the integral symplectic form).

\item Let $\mb{\Delta}^{\mb{PP}}$ be the diagonal bimodule over $({\mb P}; {\mb P})$. Let $\tilde{\mb P}_0, \tilde {\mb P}_1$ be AMS lifts of ${\mb P}$ subject to the same outercollaring width and an integral symplectic form. Then there exists a lift of the outercollaring of $\mb{\Delta}^{\mb{PP}}$ to $\outer \uds{\bf dOrb}_{\rm rig}^{\rm NC}$ as a bimodule over $(\tilde {\mb P}_0; \tilde {\mb P}_1)$.
\end{enumerate}
\end{thm}

\begin{proof}
See Section \ref{AMS_pearly}.
\end{proof}

\subsubsection{The triangle of Morse, Floer and pearly flow categories}

Recall that in Subsection \ref{subsection_PSS}, besides the Morse-to-Floer and Floer-to-Morse bimodules, we also described a few other bimodules, namely the Pearly-to-Floer bimodule $\mb{B}^{\mb{PF}}$, the Floer-to-Pearly bimodule ${\mb B}^{\mb{FP}}$, the Pearly-to-Morse bimodule ${\mb B}^{\mb{PM}}$, and the Morse-to-Pearly bimodule ${\mb B}^{\mb{MP}}$. Moreover, one has the diagonal bimodules $\mb{\Delta}^{\mb{FF}}$ (which is a special case for the continuation map bimodule), $\mb{\Delta}^{\mb{PP}}$, and $\mb{\Delta}^{\mb{MM}}$ (a deformation of the diagonal bimodule of the classical Morse flow category). These bimodules are related by a triangle which commutes up to homotopies of bimodules. 

\begin{thma}\label{thma_PSS_SSP}
Fix a common outercollaring width and an integral action ${\mc A}_{\mb F}^\Omega$ compatible with all involved bimodules and homotopy. There exist the following objects.
\begin{enumerate}

\item An AMS lift $\tilde {\mb F}$ of ${\mb F}$ (subject to the outercollaring width and the integral action).

\item An AMS lift $\tilde {\mb P}$ of $\mb{P}$ (subject to the outercollaring width and the integral symplectic form).

\item Lifts of the outercollaring of ${\mb B}^{\mb{MF}}$, $ \mb{B}^{\mb{FM}}$, $ \mb{B}^{\mb{FP}}$, $ \mb{B}^{\mb{FP}}$, $ \mb{B}^{\mb{MP}}$, $ \mb{B}^{\mb{PM}}$, $ \mb{\Delta}^{\mb{FF}}$, $ \mb{\Delta}^{\mb{PP}}$, $ \mb{\Delta}^{\mb{MM}}$, denoted by 
\beqn
\tilde{\mb B}^{\mb{MF}},\ \tilde{\mb{B}}^{\mb{FM}}, \ \tilde{\mb{B}}^{\mb{FP}},\ \tilde{\mb{B}}^{\mb{FP}},\ \tilde {\mb{B}}^{\mb{MP}},\ \tilde{\mb{B}}^{\mb{PM}}, \ \tilde{\mb{\Delta}}^{\mb{FF}},\ \tilde{\mb{\Delta}}^{\mb{PP}},\ \tilde{\mb{\Delta}}^{\mb{MM}}
\eeqn
which are bimodules over two of the three flow categories $\tilde {\mb F}, \tilde {\mb P}, \outer {\mb M}$.

\end{enumerate}
Moreover, the following diagram commutes in the following sense.
\beq\label{triangle_lift}
\xymatrix{      &  & &  \outer {\mb M} \ar@/^0.7pc/[lllddd]^{\tilde{\mb B}^{\mb{MF}}} \ar@(ur,ul)[]_{\tilde{\mb \Delta}^{\mb{MM}}}      \ar@/^0.7pc/[rrrddd]^{\tilde{\mb B}^{\mb{MP}}}      & & & \\
& & & & & &  \\
& & & & & & \\
              \tilde{\mb F}    \ar@(l,d)[]_{\tilde{\mb \Delta}^{\mb{FF}}}  \ar@/^0.7pc/[rrrrrr]^{\tilde{\mb B}^{\mb{FP}} }   \ar@/^0.7pc/[rrruuu]^{\tilde{\mb B}^{\mb{FM}} }  & &  & &  & & 
 \tilde{\mb P}     \ar@(d,r)[]_{\tilde{\mb \Delta}^{\mb{PP}} }  \ar@/^0.7pc/[llllll]^{\tilde{\mb B}^{\mb{PF}} } \ar@/^0.7pc/[llluuu]^{\tilde{\mb B}^{\mb{PM}} }  }
\eeq
All triangles are commutative up to homotopy of bimodules enriched in $\outer \uds{\bf dOrb}_{\rm rig}^{\rm NC}$, except for triangles whose vertex of concatenation is $\outer {\mb M}$.
\end{thma}

\begin{proof}
    See Section \ref{subsection_PSS_triangle}.
\end{proof}

\begin{rem}
As discussed in Section \ref{subsubsec:remark}, the regularization of homotopies related to gluing Morse trajectories does not fit well with the current framework.
\end{rem}

\section{Proof of the Main Theorems}

In this section we prove (most of) the main theorems using previously stated theorems on Kuranishi lifts of various flow categories, bimodules, homotopies and the theorem on FOP perturbations (Theorem \ref{thma_FOP}). The proofs are mostly formal consequences of the existence of the regularizations and coherent FOP perturbations thereon, following the general theme of Floer theory which turns bordifications of moduli spaces into algebraic structures.

\subsection{Proof of Theorem \ref{thma_Floer_complex}}

Suppose we are in Situation \ref{situationf1}. Theorem \ref{thma_flow_chart} provides, depending on a choice, an AMS lift $\tilde {\mb F}$ enriched in $\outer \uds{\bf dOrb}^{{\rm NC}}_{\rm rig}$ together with a coherent orientation. As the Floer flow category is an unobstructed (Definition \ref{defn_unobstructed}) and $\Pi$-equivariant Novikov flow category  (Definition \ref{defn_local_finite_flow} and Definition \ref{defn_equivariant_flow_category}), Theorem \ref{thma_FOP} provides a $\Pi$-equivariant lift, i.e., an FOP perturbation $\mathring {\mb F}$ on $\tilde {\mb F}$, which is a flow category enriched in $\outer \uds{\bf dOrb}^{\rm FOP}_{\rm rig}$. Then apply the forgetful functor (i.e., taking the closure of the isotropy-free locus of the Kuranishi section)
\beqn
\outer \uds{\bf dOrb}^{\rm FOP}_{\rm rig} \to \pman,
\eeqn
one obtains a flow category enriched in $\pman$. Using the coherent orientation and the  ${\mb Z}$-grading by the Conley--Zehnder index, by Proposition \ref{prop_differential} and Proposition \ref{prop_equivariant_differential}, one obtains a canonically associated chain complex
\beqn
CF_*(H, J; \Xi):= C_*( \mathring {\mb F} )
\eeqn
generated by objects of $\mb{F}$, which has the structure of a graded $\Lambda_{\mb Z}^\Pi$-module. Here $\Xi$ denotes all choices made in constructing the AMS lift and the FOP perturbation; in particular, $\Xi$ includes the information of the outercollaring width and the integral action. We call this process the {\bf basic construction for Floer complex}.

We prove that the Floer chain complex has a well-defined chain homotopy type. Let $\Xi_{(0)}$, $\Xi_{(1)}$ be two different sets of choices, including outercollaring widths, integral actions, and choices made for the two AMS lifts $\tilde {\mb F}_{(0)}$, $\tilde {\mb F}_{(1)}$ and the FOP perturbations $\mathring {\mb F}_{(0)}$, $\mathring {\mb F}_{(1)}$. 

\begin{lemma}\label{lemma131}
There exists a chain homotopy equivalence 
\beqn
\Psi_{(01)}: CF_*(H, J, \Xi_{(0)}) \to CF_*(H, J, \Xi_{(1)})
\eeqn
whose chain homotopy type is well-defined.
\end{lemma}

\begin{proof}
By the part of Theorem \ref{thma_flow_chart} about outercollaring width, one can reduce the comparison to the case that the outercollaring widths are equal. Then by Proposition \ref{thm_integral_action_comparison}, there exists a lift of the outercollaring of the diagonal bimodule $\mb{\Delta}^{\mb{FF}}$ to $\outer \uds{\bf dOrb}_{\rm rig}^{\rm NC}$ as a bimodule over $(\tilde{\mb F}_{(0)}; \tilde {\mb F}_{(1)})$. Then by Theorem \ref{thma_FOP}, the FOP perturbations $\mathring {\mb F}_{(0)}$ and $\mathring {\mb F}_{(1)}$ can be extended to an FOP perturbation $\mathring{\mb{\Delta}}^{\mb{FF}}$ on $\tilde{\mb{\Delta}}^{\mb{FF}}$. Then by Proposition \ref{prop_equivariant_chain_map}, one obtains the chain map
\beqn
\Psi_{(01)}: CF_*(H, J, \Xi_{(0)}) \to CF_*(H, J, \Xi_{(1)}).
\eeqn
Using the part of Proposition \ref{thm_integral_action_comparison} stating that two lifts of ${\mb \Delta}^{\mb{FF}}$ are homotopic and Theorem \ref{thma_FOP}, one can prove that the homotopy type of $\Psi_{(01)}$ is well-defined.

It remains to prove that $\Psi_{(01)}$ is a homotopy equivalence. Indeed, for each object $p \in {\rm Ob}{\mb F}$, the bimodule moduli space $M_{p;p}^{\mb{FF}}$ consists of a single point corresponding to the constant trajectory at $p$. This moduli space has no lower strata, hence coincides with its outercollaring. Moreover, it is transversely cut out, so the FOP perturbation can be chosen to be the same as the original Kuranishi section. We can then choose the FOP perturbation $\mathring{\mb{\Delta}}^{\mb{FF}}$ inductively such that we do not perturb for $M_{p;p}^{\mb{FF}}$. For such a choice of FOP perturbation, the count $n_{p;p}^{\mb{FF}}$ needed to define the corresponding coefficient of the chain map is $1$. With respect to the natural identification of the chain groups $CF_*(H, J, \Xi_{(0)}) \cong CF_*(H, J, \Xi_{(1)})$, $\Psi_{(01)}$ is a higher order perturbation of the identity in terms of $T$-valuation. Therefore, $\Psi_{(01)}$ is an invertible chain map, hence a chain homotopy equivalence. 
\end{proof}

This finishes the proof of Theorem \ref{thma_Floer_complex}.

\subsection{General construction of chain maps and homotopics}\label{subsection132}

\subsubsection{Chain maps associated to smooth Floer domains}

Suppose we are in Situation \ref{situationm1}. Choose an outercollaring width and a compatible set of integral actions
\begin{align*}
&\ {\mc A}_{{\mb F}_j}^{\Omega}: {\rm Ob}{\mb F}_j \to {\mb Z},\ &\ {\mc A}_{{\mb F}'}^{\Omega}: {\rm Ob}{\mb F}' \to {\mb Z}.
\end{align*}
Let $CF_*(H_j, J_j, \Xi_j)$ resp. $CF_*(H', J', \Xi')$ be Floer chain complexes associated to the Floer flow categories ${\mb F}_j$ resp. ${\mb F}'$, which depend on the choices $\Xi_j$ resp. $\Xi'$, including the information of the chosen outercollaring width, the integral actions, the construction of the AMS lifts $\tilde {\mb F}_j$ reps. $\tilde{\mb F}'$, and FOP perturbations $\mathring {\mb F}_j$ resp. $\mathring {\mb F}'$. Then by Theorem \ref{thma_module_lift}, there exists an AMS lift $\tilde{\mb X}$ of ${\mb X}$ subject to the outercollaring width and the integral actions as a multimodule over $(\tilde {\mb F}_1, \ldots, \tilde {\mb F}_m; \tilde {\mb F}')$. By Theorem \ref{thma_FOP}, the FOP perturbations $\mathring {\mb F}_j$ and $\mathring {\mb F}'$ can be extended to an FOP perturbation $\mathring {\mb X}$ on $\tilde {\mb X}$. Then by Proposition \ref{prop_equivariant_chain_map}, one obtains a chain map
\beqn
\Phi^{\mb{X}}(\Theta): \bigotimes_{1 \leq j \leq m} CF_*(H_j, J_j, \Xi_j) \to CF_*(H', J', \Xi')
\eeqn
where $\Theta$ contains the information of all choices made in order to define such a chain map.

In practice, the multimodule ${\mb X}$ could be either a bimodule needed for establishing the continuation map, the multimodule needed for defining the pair-of-pants product, or the multimodule for defining the Poincar\'e duality etc. We need to prove that it is well-defined up to homotopy. Again, one can reduce the comparison to the case when the outercollaring widths are the same.

\begin{lemma}\label{lemma132}
Fix the outercollaring width and the compatible set of integral actions on the involved flow categories. Let $CF_*(H_j, J_j, \Xi_j)$ resp. $CF_*(H', J', \Xi')$ be Floer chain complexes associated to these flow categories. Let 
\beqn
\Phi^{\mb{X}}(\Theta_0),\ \Phi^{{\mb X}}(\Theta_1): \bigotimes_{1\leq j \leq m} CF_*(H_j, J_j, \Xi_j) \to CF_*(H', J', \Xi')
\eeqn
be two chain maps defined via two differet sets of choices made for two AMS lifts $\tilde {\mb X}_0$ resp. $\tilde {\mb X}_1$ and two FOP perturbations $\mathring {\mb X}_0$ resp. $\mathring {\mb X}_1$. Then $\Phi^{{\mb X}}(\Theta_0)$ and $\Phi^{{\mb X}}(\Theta_1)$ are chain homotopic.
\end{lemma}

\begin{proof}
Consider the trivial homotopy $\mb{H}$ from ${\mb X}$ to itself. The set of integral actions remain compatible with respect to ${\mb H}$. Therefore, by part (2) of Theorem \ref{thma_module_lift}, there exists a homotopy $\tilde {\mb H}$ from $\tilde {\mb X}_0$ to $\tilde {\mb X}_1$ as a lift of the outercollaring of ${\mb H}$. By Theorem \ref{thma_FOP}, the FOP perturbations $\mathring {\mb F}_j$, $\mathring {\mb F}'$, $\mathring {\mb X}_0$, $\mathring {\mb X}_1$ can be extended to an FOP perturbation $\mathring {\mb H}$ on $\tilde {\mb H}$. Then by Proposition \ref{prop_Novikov_linear_chain_map} one obtains a chain homotopy between $\Phi^{{\mb X}}(\Theta_0)$ and $\Phi^{{\mb X}}(\Theta_1)$.
\end{proof}

Next we compare different choices of integral actions.

\begin{lemma}\label{lemma133}
Fix an outercollaring width. Let 
\beqn
{\mc A}_{{\mb F}_j}^{\Omega_{(s)}},\ {\mc A}_{{\mb F}'}^{\Omega_{(s)}},\ s = 0, 1
\eeqn
be two compatible sets of integral actions. Let
\beqn
CF_*(H_j, J_j, \Xi_{j, (s)}),\ CF_*(H', J', \Xi_{(s)}')
\eeqn
be Floer complexes whose constructions are subject to one of the sets of compatible integral actions, including the construction of the AMS lifts $\tilde {\mb F}_{j, (s)}$ and $\tilde {\mb F}_{(s)}'$ and FOP perturbations $\mathring {\mb F}_{j, (s)}$, $\mathring {\mb F}_{(s)}'$. Let 
\beqn
\Phi^{{\mb X}}(\Theta_{(s)}): \bigotimes_{1 \leq j \leq m} CF_*(H_j, J_j, \Xi_{j, (s)}) \to CF_*(H', J', \Xi_{(s)}')
\eeqn
be the chain map defined via two different AMS constructions subject to the two sets of integral actions, including the AMS lifts $\tilde {\mb X}_{(s)}$ and the FOP perturbation $\mathring {\mb X}_{(s)}$. 

Moreover, let 
\beqn
\Psi_{j, (01)}: CF_*(H_j, J_j, \Xi_{j, (0)}) \to CF_*(H_j, J_j, \Xi_{j, (1)})
\eeqn
and 
\beqn
\Psi_{(10)}': CF_*(H', J', \Xi_{(1)}') \to CF_*(H', J', \Xi_{(0)}')
\eeqn
be the chain maps induced from the continuation maps provided by Proposition \ref{thm_integral_action_comparison}. Then there exists a chain homotopy between $\Phi^{\mb X}(\Theta_{(0)})$ and the composition
\beqn
\Psi_{(10)}' \circ \Phi^{\mb X}(\Theta_{(1)}) \circ \Big( \Psi_{1, (01)} \otimes \cdots \otimes \Psi_{m, (01)} \Big).
\eeqn
\end{lemma}

\begin{proof}
This is essentially a consequence of Proposition \ref{prop:module-double}. Indeed, the assumption of this lemma can be plugged into the assumption of Proposition \ref{prop:module-double}. Then one obtains a homotopy $\tilde {\mb H}$ from $\tilde {\mb X}_{(0)}$ to the concatenation
\beqn
\Big( \tilde{\mb{\Delta}}_{(01)}^{\mb{F}_1 \mb{F}_1} \times \cdots \times \tilde{\mb{\Delta}}_{(01)}^{\mb{F}_m \mb{F}_m} \Big) \circ \tilde {\mb X}_{(1)} \circ \tilde{\mb{\Delta}}_{(10)}^{\mb{F}'\mb{F}'}.
\eeqn
Then by Theorem \ref{thma_FOP}, the existing FOP perturbations can be extended to an FOP perturbation $\mathring {\mb H}$ on $\tilde {\mb H}$. The conclusion then follows. 
\end{proof}

\subsubsection{Homotopy}

We prove the following basic chain-level consequence of homotopies of multimodules from varying Floer data.

\begin{lemma}\label{lemma134}
In Situation \ref{situationh1}, fix an outercollaring width and a compatible set of integral actions on involved Floer flow categories. Let 
\begin{align*}
&\ CF_*(H_j, J_j, \Xi_j),\ &\ CF_*(H', J', \Xi')
\end{align*}
be Floer chain complexes provided by Theorem \ref{thma_Floer_complex} subject to the outercollaring width and the integral actions. Let 
\beqn
\Phi^{{\mb X}_0}(\Theta_0),\ \Phi^{{\mb X}_1}(\Theta_1): \bigotimes_{1\leq j \leq m} CF_*(H_j, J_j, \Xi_j) \to CF_*(H', J', \Xi')
\eeqn
be chain maps provided by the basic construction of chain maps. Then these two chain maps are homotopic.
\end{lemma}

\begin{proof}
The construction of the Floer chain complexes involves choosing AMS lifts $\tilde {\mb F}_1, \ldots, \tilde {\mb F}_m; \tilde {\mb F}'$ and FOP perturbations $\mathring {\mb F}_1, \ldots, \mathring {\mb F}_m; \mathring {\mb F}'$. Then the chain maps are constructed by choosing AMS lifts $\tilde {\mb X}_0$ and $\tilde {\mb X}_1$ 
as multimodules over $(\tilde {\mb F}_1, \ldots, \tilde {\mb F}_m; \tilde {\mb F}')$, and FOP perturbations $\mathring {\mb X}_0$, $\mathring {\mb X}_1$ as multimodules over $(\mathring {\mb F}_1, \ldots, \mathring {\mb F}_m; \mathring {\mb F}')$. Then by Theorem \ref{thma_homotopy_lift}, there exists an AMS lift $\tilde {\mb H}$ as a homotopy from $\tilde {\mb X}_0$ to $\tilde {\mb X}_1$. Then by Theorem \ref{thma_FOP}, the existing FOP perturbations can be extended to an FOP perturbation $\mathring {\mb H}$ on $\tilde {\mb H}$. Then Proposition \ref{prop_Novikov_linear_chain_map} implies that the chain maps are homotopic via a concrete homotopy from $\mathring {\mb H}$.
\end{proof}

\subsubsection{Concatenation and gluing}

Now we consider the situation of concatenation. 

\begin{lemma}\label{lemma135}
In Situation \ref{situationg1}, fix an outercollaring width and a 1-parameter family of lateral lines in the domains. Suppose
\begin{align*}
&\ {\mc A}_{{\mb F}_1}^{\Omega},\ldots, {\mc A}_{{\mb F}_m}^\Omega,\ &\ {\mc A}_{{\mb G}_1}^\Omega,\ldots, {\mc A}_{{\mb G}_n}^\Omega, {\mc A}_{{\mb G}'}^\Omega
\end{align*}
are integral actions on the involved flow categories which are compatible with respect to ${\mb X}$, ${\mb Y}$, ${\mb M}$, and the homotopy ${\mb H}$. Subject to these choices, let
\beqn
CF(H_1, J_1, \Xi_1), \ldots, CF(H_m, J_m, \Xi_m)
\eeqn
and 
\beqn
CF(G_1, I_1, \Lambda_1), \ldots, CF(G_n, I_n, \Lambda_n), CF(G', I', \Lambda')
\eeqn
be Floer chain complexes provided by Theorem \ref{thma_Floer_complex}. Let 
\beqn
\Phi^{{\mb X}}: \bigotimes_{1\leq j \leq m} CF_*(H_j, J_j, \Xi_j) \to CF_*(G_i, I_i, \Lambda_i),
\eeqn
\beqn
\Phi^{\mb Y}: \bigotimes_{1 \leq l \leq n} CF(G_l, I_l, \Lambda_l) \to CF(G', I', \Lambda'),
\eeqn
and
\beqn
\Phi^{{\mb M}}: \bigotimes_{l<i} CF(G_l, I_l, \Lambda_l) \otimes \bigotimes_{1 \leq j \leq m} CF(H_j, J_j, \Xi_j) \otimes \bigotimes_{l>i} CF(G_l, I_l, \Lambda_l) \to CF(G', I', \Lambda')
\eeqn
be chain maps provided by the general construction of Floer chain maps. Then the last chain map is homotopic to the composition
\beqn
(y_1, \cdots, y_{i-1}, x_1, \ldots, x_m, y_{i+1}, \ldots, y_n) \mapsto \Phi^{{\mb Y}} (y_1, \ldots, y_{i-1}, \Phi^{\mb X}(x_1,\ldots, x_m), y_{i+1}, \ldots, y_n).
\eeqn
\end{lemma}

\begin{proof}
Subject to the integral actions and outercollaring widths, the construction of the Floer chain complexes depends on the AMS lifts $\tilde {\mb F}_1, \ldots, \tilde {\mb F}_m$ and $\tilde {\mb G}_1, \ldots, \tilde {\mb G}_n, \tilde {\mb G}'$ and FOP perturbations $\mathring {\mb F}_1, \ldots, \mathring {\mb F}_m$ and $\mathring {\mb G}_1, \ldots, \mathring {\mb G}_n, \mathring {\mb G}'$. The chain maps $\Phi^{\mb X}$, $\Phi^{\mb Y}$, and $\Phi^{\mb M}$ also depend on corresponding AMS lifts and FOP perturbations. One can also choose concatenations $\tilde {\mb X} \circ_i \tilde {\mb Y}$ and $\mathring {\mb X} \circ_i \mathring {\mb Y}$ (provided by Theorem \ref{thm_AMS_concatenation} and Proposition \ref{prop_FOP_concatenation}). 
Then by Theorem \ref{thma_multimodule_gluing}, there exists a lift $\tilde {\mb H}$ of the outercollaring of ${\mb H}$ as a homotopy from $\tilde {\mb X} \circ_i \tilde {\mb Y}$ to $\tilde {\mb M}$. By Theorem \ref{thma_FOP}, one can find an FOP perturbation $\mathring {\mb H}$ on $\tilde {\mb H}$ as a homotopy from $\tilde {\mb X} \circ \tilde {\mb Y}$ to $\tilde {\mb M}$. Then one obtains the desired chain homotopy.
\end{proof}

\subsection{Proof of Theorem \ref{thma_continuation_map}}\label{subsection_proof_thma}

Assume we are in Situation \ref{situationf1}. Above, we have defined the Floer chain complex 
\beqn
CF_*(H, J, \Xi)
\eeqn
which appears to depend on the outercollaring width, the integral action, choices made for constructing the AMS lift $\tilde {\mb F}$, and the FOP perturbation $\mathring {\mb F}$. We have already proved (see Lemma \ref{lemma131}) that the chain homotopy type only depends on $(H, J)$, but not on the system of choices $\Xi$. It remains to construct the continuation map.

\subsubsection{Definition of the continuation map}

Let ${\mb F}'$ be another Floer flow category associated to another pair $(H', J')$. Choose a 1-parameter family of Floer data interpolating between $(H, J)$ and $(H', J')$, resulting in a bimodule ${\mb B}$ from ${\mb F}$ to ${\mb F}'$. Then we are in a special case of Situation \ref{situationm1}. Then, using the above construction, upon choosing an outercollaring width and a pair of integral actions on ${\mb F}$ and ${\mb F}'$ compatible with respect to ${\mb B}$, after AMS constructions of $\tilde {\mb F}$, $\tilde {\mb F}'$,$\tilde {\mb B}$, and choosing FOP perturbations $\mathring {\mb F}$, $\mathring {\mb F}'$, and $\mathring {\mb B}$, one obtains a chain map
\beqn
\Phi^{\mb B}(\Theta): CF_*(H, J, \Xi) \to CF_*(H', J', \Xi')
\eeqn
whose homotopy class only depends on the bimodule ${\mb B}$ and the choices of the integral actions (see Lemma \ref{lemma132}). 

\subsubsection{Independence of the integral actions}

We first would like to reduce the dependence to ${\mb B}$ but not the integral actions. Indeed, if ${\mc A}_{{\mb F}}^{\Omega_{(s)}}$, ${\mc A}_{{\mb F}'}^{\Omega_{(s)}}$, $s = 0, 1$ are two compatible pairs of integral actions, then by Lemma \ref{lemma133}, the diagram
\beqn
\xymatrix{ CF_*(H, J, \Xi_{(0)}) \ar[rr]^-{\Phi^{{\mb B}}(\Theta_{(0)}) } \ar[d]_{\Psi_{(01)}^{\mb F}}     &    &CF_*(H', J', \Xi_{(0)}')  \\
          CF_*(H, J, \Xi_{(1)}) \ar[rr]_-{\Phi^{{\mb B}}(\Theta_{(1)}) }  &   &   CF_*(H', J', \Xi_{(1)}')  \ar[u]_{\Psi_{(10)}^{{\mb F}'}}  }
\eeqn
commutes up to homotopy. Therefore, the homotopy class only depends on the bimodule ${\mb B}$. 

\subsubsection{Independence of the interpolation of Floer data}

Next we would like to prove the homotopy independence of the bimodule ${\mb B}$. Suppose ${\mb B}_0$ and ${\mb B}_1$ are two such bimodules over $({\mb F}; {\mb F}')$. Since ${\mb B}_0$ and ${\mb B}_1$ are both defined via varying Floer data, one can construct a homotopy ${\mb H}$ from ${\mb B}_0$ to ${\mb B}_1$. Then by Proposition \ref{prop163}, one can find a pair of integral actions ${\mc A}_{{\mb F}}^\Omega$, ${\mc A}_{{\mb F}'}^{\Omega}$ compatible with respect to ${\mb H}$. This leads to the setup of Lemma \ref{lemma134}. Let 
\begin{align*}
&\ CF_*(H, J, \Xi),\ &\ CF_*(H', J', \Xi')
\end{align*}
be Floer chain complexes constructed subject to these integral actions and let 
\beqn
\Phi^{{\mb B}_0}(\Theta_0),\ \Phi^{{\mb B}_1}(\Theta_1): CF_*(H, J, \Xi) \to CF_*(H', J', \Xi')
\eeqn
be continuation maps associated to these two bimodules constructed as above. Then by Lemma \ref{lemma134}, the two continuation maps are homotopic. This proves that the continuation map does not depend on the bimodule either. 

\subsubsection{Composition of continuation maps}

We prove that the composition of continuation maps is homotopic to another continuation map.

\begin{prop}\label{prop134}
Suppose ${\mb F}_0, {\mb F}_1, {\mb F}_2$ are Floer flow categories associated to three pairs $(H_i, J_i)$ for $i = 0, 1, 2$. Let ${\mb B}_{01}$ resp. ${\mb B}_{12}$ be a bimodule over $({\mb F}_0; {\mb F}_1)$ resp. over $({\mb F}_1; {\mb F}_2)$ associated to Floer data $\sigma_{01}$ resp. $\sigma_{12}$ on the infinite cylinder. Let ${\mb H}$ be a homotopy associated to gluing $\sigma_{01}$ and $\sigma_{12}$, as a homotopy from the composition ${\mb B}_{01}\circ {\mb B}_{12}$ to a bimodule ${\mb B}_{02}$ over $({\mb F}_0; {\mb F}_2)$ associated to a Floer datum $\sigma_{02}$. Then there exist continuation maps
\beqn
\Phi^{{\mb B}_{01}}: CF_*(H_0, J_0, \Xi_0) \to CF_*(H_1, J_1, \Xi_1)
\eeqn
\beqn
\Phi^{{\mb B}_{12}}: CF_*(H_1, J_1, \Xi_1) \to CF_*(H_2, J_2, \Xi_2)
\eeqn
and 
\beqn
\Phi^{\mb{B}_{02}}: CF_*(H_0, J_0, \Xi_0) \to CF_*(H_2, J_2, \Xi_2)
\eeqn
such that the last one is homotopic to the composition $\Phi^{{\mb B}_{12}} \circ \Phi^{\mb{B}_{01}}$.
\end{prop}

\begin{proof}
By Proposition \ref{prop163}, we can choose integral actions on the three Floer categories which are compatible with respect to ${\mb B}_{01}$, ${\mb B}_{12}$, ${\mb B}_{02}$, and the homotopy ${\mb H}$. Then we are in the situation of Lemma \ref{lemma135}, which implies the desired homotopy.
\end{proof}

\subsubsection{Continuation maps are homotopy equivalences}

\begin{cor}
Any chain-level continuation map is a homotopy equivalence. 
\end{cor}

\begin{proof}
Using the notations in Proposition \ref{prop134}, take ${\mb F}_2 = {\mb F}_0$ and denote ${\mb B}_{12}$ by ${\mb B}_{10}$. Let ${\mb B}_{02}$ be the diagonal bimodule $\mb{\Delta}^{\mb{F}_0\mb{F}_0}$ associated to the trivial Floer datum on the cylinder. Choose a homotopy ${\mb H}$ from $\mb{B}_{01} \circ \mb{B}_{10}$ to $\mb{\Delta}^{\mb{F}_0\mb{F}_0}$ by gluing Floer data and varying in a one-parameter family, Proposition \ref{prop134} provides a chain homotopy between the composition $\Phi_{01} \circ \Phi_{10}$ and the chain map associated to an AMS lift of $\mb{\Delta}^{\mb{F}_0\mb{F}_0}$. Notice that the latter has the form 
\beqn
\Phi^{{\mb \Delta}^{\mb{FF}}} = {\rm Id} + {\rm higher\ order\ in\ } T
\eeqn
hence is an invertible chain map. Therefore, by a similar argument switching the roles of $\Phi_{01}$ and $\Phi_{10}$, the continuation maps $\Phi_{01}$ and $\Phi_{10}$ are homotopy inverses to each other. 
\end{proof}

Lastly, we prove 

\begin{cor}
The chain map from $CF_*(H, J, \Xi) \to CF_*(H, J, \Xi)$ induced from an AMS lift of the diagonal bimodule $\mb{\Delta}^{\mb{FF}}$ is  chain homotopic to the identity.
\end{cor}

\begin{proof}
The same argument of the proof of Lemma \ref{lemma131} also applies here, resulting in a chain map associated to an AMS lift of the bimodule $\mb{\Delta}^{\mb{FF}}$ of the form
\beqn
\Phi = {\rm Id} + A.
\eeqn
Notice that ${\rm Id}$ is also a chain map in this case. Now since $\mb{\Delta}^{\mb{FF}} \circ\mb{\Delta}^{\mb{FF}}$ is homotopic to $\mb{\Delta}^{\mb{FF}}$, Proposition \ref{prop134} implies that 
\beqn
({\rm Id} + A)^2 - ({\rm Id} + A) = A ({\rm Id} + A) = d\rho + \rho d
\eeqn
for some linear map $\rho$. As ${\rm Id} + A$ is an invertible chain map, one has 
\beqn
A = (d \rho + \rho d)({\rm Id} + A)^{-1} = d( \rho ({\rm Id} + A)^{-1}) + ( \rho ({\rm Id} + A)^{-1}) d
\eeqn
meaning that ${\rm Id} +A$ is chain homotopic to ${\rm Id}$.
\end{proof}

\subsection{Proof of Theorem \ref{thma_product}}\label{subsection_proof_product}

\subsubsection{Definition of the pair-of-pants products}

The setting of Theorem \ref{thma_product} falls into Situation \ref{situationm1}. Let $\Sigma^{\mb X}$ be a smooth pair-of-pants with two negative ends and one positive end equipped with a Floer datum which are asymptotic to $H_1dt$ and $H_2dt$ on the negative ends and to $H'dt$ on the positive end. To simplify the argument, use a fixed almost complex structure $J$; more general cases can be treated with the help of continuation maps. Hence we omit the reference to almost complex structures. Let ${\mb F}_1, {\mb F}_2, {\mb F}'$ be the Floer flow categories associated to $H_1$, $H_2$, and $H'$ respectively.

We fix an outercollaring width and will not mention it anymore in this subsection. Let ${\mb X}$ be the multimodule associated to $\Sigma^{\mb X}$. Then by Proposition \ref{prop163}, one can choose a triple of integral actions ${\mc A}_{{\mb F}_1}^\Omega, {\mc A}_{{\mb F}_2}^\Omega$, ${\mc A}_{{\mb F}'}^\Omega$. Then one obtains Floer chain complexes 
\beqn
CF_*(H_1, \Xi_1), CF_*(H_2, \Xi_2), CF_*(H', \Xi')
\eeqn
which depend on the AMS lifts $\tilde {\mb F}_1, \tilde {\mb F}_2, \tilde {\mb F}'$ subject to the integral actions and their FOP perturbations $\mathring {\mb F}_1, \mathring {\mb F}_2, \mathring {\mb F}'$. Upon choosing a pair of lateral lines in $\Sigma^{\mb X}$, the basic construction provides a chain map
\beqn
\Phi^{\mb X}(\Theta): CF_*(H_1, \Xi_1) \otimes CF_*(H_2, \Xi_2) \to CF_*(H', \Xi')
\eeqn
which depends on an AMS lift $\tilde {\mb X}$ and an FOP perturbation $\mathring {\mb X}$. 

\subsubsection{Well-defined homotopy class}

Now we prove Part (1) of Theorem \ref{thma_product}. Similar to the case of continuation maps, by Lemma \ref{lemma132}, the homotopy class of $\Phi^{{\mb X}}(\Theta)$ does not depend on choices made for constructing the AMS lift $\tilde {\mb X}$ and FOP perturbation $\mathring {\mb X}$, so {\it a priori} only depends on the integral actions and the choices of lateral lines. These dependence is removed (up to homotopy) by Lemma \ref{lemma133}. Hence the homotopy class of $\Phi^{{\mb X}}(\Theta)$ {\it a priori} only depends on the Floer data defining ${\mb X}$. Lastly, if ${\mb X}_0$ and ${\mb X}_1$ are multimodules defined by two different Floer data on the pair-of-pants (including different cylindrical ends), then one can build a homotopy ${\mb H}$ from ${\mb X}_0$ to ${\mb X}_1$ by varying the Floer data. Then Lemma \ref{lemma134}  implies that the corresponding chain maps $\Phi^{{\mb X}_0}(\Theta_0)$ and $\Phi^{{\mb X}_1}(\Theta_1)$ are also homotopic.

\subsubsection{Compatibility with the continuation maps}

We only treat the case with composing continuation maps at the positive end, and the other cases follow from the same reasoning. Let ${\mb F}''$ be a fourth Floer flow category associated to $(H'', J'')$. We assume that the almost complex structure is fixed and drop it from the notation; the general case when $J$ is domain-dependent can be treated using continuation maps. Let $\mb{B}^{\mb{F}'\mb{F}''}$ be a bimodule over $({\mb F}'; {\mb F}'')$ provided by continuation map Floer data on the cylinder. Consider the concatenation ${\mb X} \circ \mb{B}^{\mb{F}'\mb{F}''}$. Then by gluing the Floer domains, one obtains a homotopy ${\mb H}$ from ${\mb X} \circ \mb{B}^{\mb{F}'\mb{F}''}$ to a multimodule ${\mb Y}$ over $({\mb F}_1, {\mb F}_2; {\mb F}'')$ defined on the pair-of-pants with another Floer datum. 

Then by Proposition \ref{prop163}, one can choose integral actions on involved Floer flow categories which are compatible with respect to $\mb{X}$, $\mb{B}^{\mb{F}'\mb{F}''}$, ${\mb Y}$, and ${\mb H}$. Then we are in the situation of Lemma \ref{lemma135}, which provides specific chain complexes 
\beqn
CF_*(H_1, \Xi_1), CF_*(H_2, \Xi_2), CF_*(H', \Xi'), CF_*(H'', \Xi'')
\eeqn
chain maps
\beqn
\begin{split}
\Phi^{\mb X}: CF_*(H_1, \Xi_1) \otimes CF_*(H_2, \Xi_2) \to &\ CF_*(H', \Xi'),\\
\Phi^{\mb{F}'\mb{F}''}: CF_*(H', \Xi') \to &\  CF_*(H'', \Xi''),\\
\Phi^{\mb Y}: CF_*(H_1, \Xi_1) \otimes CF_*(H_2, \Xi_2) \to &\ CF_*(H'', \Xi'').
\end{split}
\eeqn
The Lemma \ref{lemma135} shows that $\Phi^{\mb Y}$ is chain homotopic to the composition $\Phi^{\mb{F}'\mb{F}''} \circ \Phi^{\mb X}$. This is Theorem \ref{thma_product} (2).

\subsubsection{Associativity}

Associativity follows from Lemma \ref{lemma134} and \ref{lemma135}. In the notations of chain complexes, we drop the dependence on $J$. 

\begin{figure}[h]
    \centering
    \includegraphics[width=0.9\linewidth]{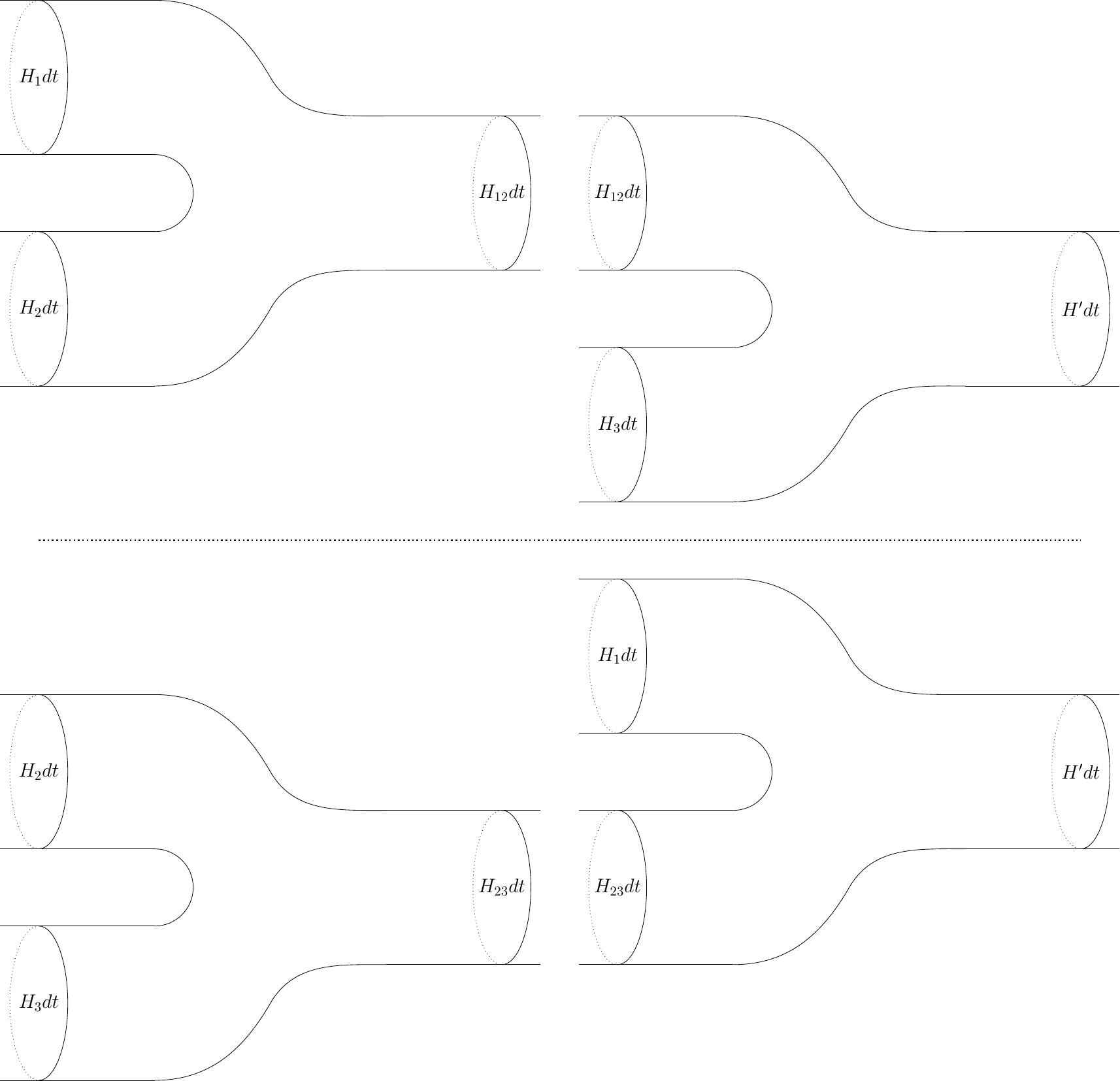}
    \caption{The homotopy between two concatenations.}
    \label{Figure_associativity}
\end{figure}

Choose Hamiltonians $H_1, H_2, H_3, H'$, $H_{12}$, $H_{23}$ independently. One can build a homotopy of two concatenated multimodules as shown in Figure \ref{Figure_associativity}. One can choose integral actions on all involved Floer flow categories which are compatible with respect to this homotopy. Then apply Lemma \ref{lemma135} twice, one obtains the following diagram in which each triangle commutes up to homotopy.
\beqn
\xymatrix{ CF_*(H_1, \Xi_1) \otimes CF_*(H_2, \Xi_2) \otimes CF_*(H_3, \Xi_3)   \ar[rd] \ar[d]   \ar[r] &  CF_*(H_1, \Xi_1) \otimes CF_*(H_{23}, \Xi_{23}) \ar[d]   \\
CF_*(H_{12}, \Xi_{12}) \otimes CF_*(H_3, \Xi_3) \ar[r] &   CF_*(H', \Xi')  }
\eeqn
Here the diagonal arrow is provided by a multimodule with three source flow categories. As the homotopy types of continuation maps are well-defined, the associativity is proved. This is Theorem \ref{thma_product} (3). 

\subsubsection{Identity element}

There is a natural multiplicative identity of the pair-of-pants product coming from the cigar bimodule (see Subsection \ref{subsection_cigar}). Let ${\mb F}$ be the Floer flow category associated to a nondegenerate Hamiltonian $H$ on the almost K\"ahler manifold $(X, \omega, J)$. On the other hand, let ${\mb O}$ be the trivial flow category consisting of only one object. Then the cigar bimodule is a bimodule ${\mb B}^{\mb{OF}}$ associated to a Floer data on ${\mb C}$ which is asymptotic to $H dt$ near $\infty$. Then this is a special case of Situation \ref{situationm1}. 

By Proposition \ref{prop163}, one can choose an integral action ${\mc A}_{\mb F}^\Omega$ on ${\mb F}$ which is compatible with respect to $\mb{B}^{\mb{OF}}$. Choose an outercollaring width. Let 
\beqn
CF_*(H, J, \Xi)
\eeqn
be the associated Floer chain complex depending an AMS lift $\tilde{\mb F}$ subject to the outercollaring width, the integral action, and an FOP perturbation $\mathring {\mb F}$. The basic construction for multimodules then provides a chain map
\beqn
\Phi^{\mb{OF}}(\Theta): {\mb Z} \cong C_*({\mb O})  \to CF_*(H, J, \Xi)
\eeqn
which depends on an AMS lift $\tilde {\mb B}^{\mb{OF}}$ and an FOP perturbation $\mathring {\mb B}^{\mb{OF}}$. 

Following the same argument as before, one can prove that the homotopy class of the chain map $\Phi^{\mb{OF}}(\Theta)$ is well-defined, independent of all choices as well as the Floer data responsible for the bimodule ${\mb B}^{\mb{OF}}$. Moreover, $\Phi^{\mb{OF}}(\Theta)$ is compatible with continuation maps up to homotopy. As the differential of $C_*(\mb{O})$ is trivial, the image of $1 \in {\mb Z}$ is closed and its homology class is well-defined. Therefore, one obtains a well-defined element 
\beqn
{\bf 1} \in HF_*(X, \omega; \Lambda_{\mb Z}^\Pi).
\eeqn

We prove that ${\bf 1}$ is the multiplicative identity. This follows from a particular concatenation of multimodules and Lemma \ref{lemma135}. Consider the concatenation of the bimodule $\mb{B}^{\mb{OF}}$ with a multimodule ${\mb X}$ over $({\mb F}, {\mb F}'; {\mb F}'')$ for pair-of-pants product. One can see that by gluing the Floer data, the concatenation is homotopic to a bimodule for continuation map $\mb{B}^{\mb{F}'\mb{F}''}$, hence a special case of Situation \ref{situationg1}. Therefore, by following similar argument as previous situations where gluing is concerned, one can prove that ${\bf 1}$ is a multiplicative identity. We omit the details. This is Theorem \ref{thma_product} (4) and we conclude the proof.

\subsection{Proof of Theorem \ref{thma_Poincare}}\label{subsection_PD}

We establish the Floer-theoretic Poincar\'e duality. Consider the punctured Riemann sphere $\Sigma^{\rm PD}$ which has two negative ends and no positive end. A Floer datum whose restrictions to the two negative ends are $(H_idt, J_i)$, $i = 1, 2$, leads to the multimodule ${\mb X}^{PD}$ over $({\mb F}_1, {\mb F}_2; \tilde {\mb O})$  where $\tilde {\mb O}$ is the $\Pi$-covering of the trivial flow category ${\mb O}$. In particular, ${\rm Ob}\tilde {\mb O} \cong \Pi$. 

Now we define hte Poincar\'e pairing on the chain level. By Proposition \ref{prop163}, one can choose integral actions on ${\mb F}_1$, ${\mb F}_2$ compatible with respect to ${\mb X}^{\rm PD}$. Fix an outercollaring width. By Theorem \ref{thma_Floer_complex}, one obtains chain complexes 
\beqn
CF_*(H_1, J_1, \Xi_1),\ CF_*(H_2, J_2, \Xi_2).
\eeqn
Then the general construction of chain maps provided in Subsection \ref{subsection132} provides 
\beqn
\Phi^{\rm PD}: CF_*(H_1, J_1, \Xi_1) \otimes CF_*(H_2, J_2, \Xi_2) \to C_*(\tilde {\mb O}) \cong \Lambda^\Pi.
\eeqn

In the same way as proving independence of choices, one can see that $\Phi^{\rm PD}$ has a well-defined homotopy type which does not depend on choices made for the AMS construction and FOP perturbations. It is also compatible with continuation maps, and compatible with the comparison between different choices of integral actions (here we do not need lateral lines). Hence the Poincar\'e pairing is well-defined. 

We need to prove the compatibility between the Poincar\'e pairing and the pair-of-pants product. In fact this is similar to the proof of the associativity of the pair-of-pants product, using the concatenation of a multimodule ${\mb X}$ (for the pair-of-pants products) and ${\mb X}^{\rm PD}$ and a gluing construction, utilizing Lemma \ref{lemma134} and Lemma \ref{lemma135}. We omit the details.

Lastly, we prove the nondegeneracy of the Poincar\'e pairing with field coefficients. We only need to prove it on the homology level. We make use of the PSS map $\Phi^{\mb{MF}}$ regarded as a map on homology. One can see that the composition
\beqn
\Phi^{\rm PD} \circ (\Phi^{\mb{MF}}\otimes \Phi^{\mb{MF}})
\eeqn
is a bilinear pairing on $H_*(M; \Lambda_{\bf K}^\Pi)$. The pairing has a leading order term which can be identified with the classical Poincar\'e pairing on the integral homology of $X$. Hence the induced pairing is still nondegenerate on the Morse homology (tensoring with $\Lambda_{\bf K}^\Pi$). As the PSS map is an isomorphism on homology (Theorem \ref{thma_PSS_isomorphism}), it implies that the Poincar\'e pairing on the Floer homology is also nondegenerate. This finishes the proof of Theorem \ref{thma_Poincare}.

\subsection{Proof of Theorem \ref{thma_PSS}}

The definition of the PSS and the SSP maps are simply special cases of general chain maps induced from bimodules. Fix a Floer flow category ${\mb F}$ on $(X, \omega)$ associated to a pair $(H, J)$. To simplify the argument we assume $J$ is time-independent; the general case can be treated with continuation maps. The same data provide the diagonal bimodule $\mb{\Delta}^{\mb{FF}}$. Choose a Morse--Smale pair $(f, g)$ on $M$. Then one has the Morse flow category ${\mb M}$ associated to $(f, g)$. Choosing an interpolation between $H$ and $0$, one obtains a Floer-to-Morse bimodule $\mb{B}^{\mb{FM}}$; a backward interpolation provides a Morse-to-Floer bimodule $\mb{B}^{\mb{MF}}$.

Fix a common outercollaring width. Choose an integral action ${\mc A}_{\mb F}^\Omega$ on ${\mb F}$ which are compatible with all bimodules and homotopies involved in the discussion. Then Theorem \ref{thma_Floer_complex} provides a Floer complex $CF_*(H, J, \Xi)$ depending on choices, which involve an AMS lift $\tilde {\mb F}$ subject to the outercollaring width and the integral action, and an FOP perturbation $\mathring {\mb F}$. Then by Theorem \ref{thma_PSS_AMS}, one can find AMS lifts $\tilde {\mb B}^{\mb{MF}}$ and $\tilde{\mb B}^{\mb{FM}}$ (which are enriched in $\outer\uds{\bf dOrb}^{\rm NC}_{\rm rig}$) of $\mb{B}^{\mb{MF}}$ and $\mb{B}^{\mb{FM}}$. Notice that $\outer {\mb M}$ can be viewed as enriched in $\outer \uds{\bf dOrb}_{\rm rig}^{\rm FOP}$. Then the FOP perturbation $\mathring {\mb F}$ can be extended, by Theorem \ref{thma_FOP}, to an FOP perturbation $\mathring {\mb B}^{\mb{MF}}$ and an FOP perturbation $\mathring {\mb B}^{\mb{FM}}$. Then one obtains chain maps
\begin{align*}
&\ \Phi^{\mb{MF}}: CM_*(f, g)  \to CF_*(H, J,  \Xi),\ &\ \Phi^{\mb{FM}}:  CF_*(H, J, \Xi) \to CM_*(f, g).
\end{align*}
We call them the {\bf chain-level PSS map} and the {\bf chain-level SSP map}. 

Using similar argument as the case of continuation maps, based on part of Theorem \ref{thma_PSS_AMS} and Theorem \ref{thm126}, one can show that the homotopy classes of the PSS and SSP maps are independent of the choices made for constructing the AMS lifts or the FOP lifts. Moreover, using homotopies realized by gluing Floer cylinders, one can show that the PSS and SSP maps are compatible with continuation maps. Namely, the composition of a chain-level PSS map resp. a chain-level SSP map with a chain-level continuation map (in the right way) is homotopic to a chain-level PSS resp. a chain-level SSP map. This finishes the proof of Theorem \ref{thma_PSS}.

\subsection{Proof of Theorem \ref{thma_PSS_isomorphism}}

To prove Theorem \ref{thma_PSS_isomorphism}, what remains to show is that the PSS and SSP maps are chain homotopy equivalences. We remind the reader again that the PSS and SSP maps are not known to be homotopy inverses to each other. To make the comparison, we need to go through the pearly flow category $\mb{P}$. 

Keep the data fixed in the proof of Theorem \ref{thma_PSS}. Then the Morse--Smale pair $(f, g)$ and the almost complex structure $J$ also determines the pearly flow category $\mb{P}$ (see Subsection \ref{subsection_pearly}). Meanwhile, one has the diagonal bimodules 
\beqn
\mb{\Delta}^{\mb{FF}}, \mb{\Delta}^{\mb{MM}},\ \mb{\Delta}^{\mb{PP}}
\eeqn
and the Morse-to-pearly bimodule $\mb{B}^{\mb{MP}}$, pearly-to-Morse bimodule $\mb{B}^{\mb{PM}}$. Using the same Floer data on the infinite cylinder defining the bimodules $\mb{B}^{\mb{MF}}$, $\mb{B}^{\mb{FM}}$ also provides the Floer-to-pearly bimodule $\mb{B}^{\mb{FP}}$ and the pearly-to-Floer bimodule $\mb{B}^{\mb{PF}}$. Moreover, one can build homotopies realizing the commutativity-up-to-homotopy of each triangle of the diagram \eqref{big_triangle}. 

Fix an outercollaring width. Choose an integral action on ${\mb F}$ which is compatible with respect to all involved bimodules and homotopies. Then upon making choices, one obtains the Floer chain complex $CF_*(H, J, \Xi)$ by Theorem \ref{thma_Floer_complex} and the PSS and SSP maps as in the proof of Theorem \ref{thma_PSS}. 

Moreover, by Theorem \ref{thma_pearly_chart}, upon making choices, one can construct an AMS lift of the outercollaring of ${\mb P}$ to $\outer \uds{\bf dOrb}_{\rm rig}^{\rm NC}$, denoted by $\tilde {\mb P}$. By choosing an FOP perturbation $\mathring {\mb P}$ (granted by Theorem \ref{thma_FOP}), one obtains a chain complex
\beqn
CP_*(f, g, J, \Theta)
\eeqn
where $\Theta$ represents all the choices, including the outercollaring width, the integral symplectic form contained in ${\mc A}_{\mb F}^\Omega$, the choices made for constructing $\tilde {\mb P}$, and the FOP perturbation. Then Theorem \ref{thma_pearly_chart} also provides a lift $\tilde{\mb \Delta}^{\mb{PP}}$ of the outercollaring of $\mb{\Delta}^{\mb{PP}}$ as a bimodule over $(\tilde {\mb P}; \tilde {\mb P})$. Then Theorem \ref{thma_FOP} allows us to extend the FOP perturbations $\mathring {\mb P}$ to an FOP perturbation $\mathring{\mb{\Delta}}^{\mb{PP}}$. Then one obtains a chain map
\beqn
\Phi^{\mb{PP}}: CP_*(f, g, J, \Theta) \to CP_*(f, g, J, \Theta).
\eeqn

Based on the above constructions and choices, Theorem \ref{thma_PSS_SSP} guarantees the existence of chain maps
\begin{align*}
&\ \Phi^{\mb{FP}}: CF_*(H, J, \Xi) \to CP_*(f, g, J, \Theta),\ &\ \Phi^{\mb{PF}}: CP_*(f, g, J, \Theta) \to CF_*(H, J, \Xi).
\end{align*}
and
\begin{align*}
&\ \Phi^{\mb{PM}}: CP_*( f, g, J, \Theta) \to CM_*(f, g) \otimes \Lambda^\Pi,\ &\ \Phi^{\mb{MP}}: CM_*(f, g)\otimes \Lambda^\Pi \to CP_*(f, g, J, \Theta).
\end{align*}

\begin{lemma}\label{lemma139}
$\Phi^{\mb{FM}}$ is chain homotopic to $\Phi^{\mb{PM}} \circ \Phi^{\mb{FP}}$ and $\Phi^{\mb{MF}}$ is chain homotopic to $ \Phi^{\mb{PF}} \circ \Phi^{\mb{MP}}$.
\end{lemma}

\begin{proof}
One has fixed a specific homotopy $\mb{H}$ from the concatenation $\mb{B}^{\mb{FP}} \circ \mb{B}^{\mb{PM}}$ to $\mb{B}^{\mb{FM}}$. By Theorem \ref{thm_AMS_concatenation}, there exists a concatenation $\tilde{\mb{B}}^{\mb{FP}}\circ \tilde{\mb{B}}^{\mb{PM}}$ enriched in $\outer \uds{\bf dOrb}_{\rm rig}^{\rm NC}$. Then by Theorem \ref{thm_AMS_gluing_1}, there exists a lift $\tilde {\mb H}$ of the outercollaring of ${\mb H}$ to $\outer \uds{\bf dOrb}_{\rm rig}^{\rm NC}$ as a homotopy from $\tilde{\mb{B}}^{\mb{FP}}\circ \tilde{\mb{B}}^{\mb{PM}}$ to $\tilde {\mb B}^{\mb{FM}}$. Moreover, Proposition \ref{prop_FOP_concatenation} allows us to have a concatenation $\mathring{\mb{B}}^{\mb{FP}}\circ \mathring{\mb{B}}^{\mb{PM}}$. Then Theorem \ref{thma_FOP} allows us to extend this concatenation and the FOP perturbation $\mathring {\mb{B}}^{\mb{FM}}$ to an FOP perturbation $\mathring{\mb H}$ on $\tilde {\mb H}$. Then one obtains a chain homotopy between $\Phi^{\mb{FM}}$ to $\Phi^{\mb{PM}}\circ \Phi^{\mb{FP}}$. The other chain homotopy can be established in the same way. 
\end{proof}

\begin{lemma}
$\Phi^{\mb{FP}}$ and $\Phi^{\mb{PF}}$ are invertible after passing to homology.
\end{lemma}

\begin{proof}
Similar to the above lemma, by Theorem \ref{thma_PSS_SSP}, using FOP perturbations, one obtains a chain homotopy between $\Phi^{\mb{PF}} \circ \Phi^{\mb{FP}}$ and the chain-level continuation map associated to the diagonal bimodule $\mb{\Delta}^{\mb{FF}}$, which is invertible after passing to homology. Via the reversed gluing along the Floer flow category, one obtains a chain homotopy between $\Phi^{\mb{FP}} \circ \Phi^{\mb{PF}}$ and the chain map $\Phi^{\mb{PP}}$ associated to the diagonal bimodule $\mb{\Delta}^{\mb{PP}}$, which is also invertible. Hence $\Phi^{\mb{PF}}$ and $\Phi^{\mb{FP}}$ are both invertible after passing to homology.
\end{proof}

\begin{lemma}
$\Phi^{\mb{PM}}$ and $\Phi^{\mb{MP}}$ are invertible chain maps.
\end{lemma}

\begin{proof}
Notice that the objects of $\mb{P}$ and $\mb{M}$ are identical. Moreover, both $\Phi^{\mb{PM}}$ and $\Phi^{\mb{MP}}$ have leading order term being the identity map. Hence both these two chain maps are invertible. 
\end{proof}

Now by the previous three lemmas, $\Phi^{\mb{FM}}$ and $\Phi^{\mb{MF}}$ are chain homotopic to invertible chain maps. Hence they are chain homotopy equivalences. This finishes the proof of Theorem \ref{thma_PSS_isomorphism}.

\subsection{Proof of Theorem \ref{thma_semipositive}}\label{subsection_proof_semipositive}

By Hofer--Salamon \cite{Hofer_Salamon} and Ono \cite{Ono_1995} (see also \cite{Floer_Hofer_Salamon} and \cite{McDuff_Salamon_2004}), one can choose a time-independent $J$ and suitable nondegenerate $H$ satisfying
\begin{enumerate}
    \item All simple nonconstant $J$-holomorphic sphere $u: S^2 \to X$ is transverse.

    \item When ${\rm deg} (p) - {\rm deg}(q) \leq 1$, $M_{pq}^{{\mb F}}$ is transverse and consists of stable Floer trajectories without sphere bubbles.

    \item When ${\rm deg} (p) - {\rm deg}(q) = 0$, $M_{pq}^{{\mb F}}$ consists of maps defined over smooth cylinders.
\end{enumerate}

To proceed, one needs to carefully make choices when constructing the AMS lift $\tilde {\mb F}$, especially the smooth structure. Choose an outercollaring width and an integral action on ${\mb F}$. Let $\hat{\mb F}$ be a relatively NC AMS lift of ${\mb F}$ to $\outer \uds{\bf S^{\rm rel}Kur}_{\rm rig}^{\rm NC}$. Let $K_{pq} = (G_{pq}, V_{pq}/B_{pq}, E_{pq}, S_{pq})$ be the corresponding relative Kuranishi space. 

Notice that for pair $p<q$, there is an open subset 
\beqn
\mathring M_{pq}^{\mb F} \subset M_{pq}^{\mb F}
\eeqn
corresponding to smooth Floer trajectories. By the above assumption on $H$ and $J$, we see that when ${\rm deg}(p) - {\rm deg}(q) \leq 0$, $M_{pq}^{\mb F} = \mathring M_{pq}^{\mb F}$. There is also the corresponding open subset
\beqn
\mathring V_{pq}/ \mathring B_{pq}.
\eeqn
We know that $\mathring V_{pq}$ has a canonical smooth structure such that the map $\mathring V_{pq} \to \mathring B_{pq}$ is a smooth submersion. Moreover, $\mathring S_{pq}:=S_{pq}|_{\mathring V_{pq}}$ is smooth. 

Recall that the stable smoothing of $\hat{\mb F}$ should be given by a $G_{pq}$-smoothing on stabilizations $\hat V_{pq}:= V_{pq} \times W_{pq}$ where $W_{pq}$ is a linear representation of $G_{pq}$. The open set $\mathring V_{pq}\times W_{pq}$ has a canonical smooth structure. 

\begin{lemma}\label{lemma1312}
There exists a stable smoothing $\hhat {\mb F}$ of $\hat{\mb F}$ satisfying the following conditions. If $p<q$ and ${\rm deg}(p) \leq {\rm deg}(q)$, then $\hat V_{pq}$ is equivariantly diffeomorphic to $\mathring V_{pq} \times W_{pq}$ near $S_{pq}^{-1}(0)$. In particular, the stabilization $\hat S_{pq}$ is smooth near its zero locus and transverse in the classical sense. 
\end{lemma}

\begin{proof}
This is because the smoothing can be done locally and extend (cf. \cite[Theorem B.3]{Bai_Xu_Arnold}).
\end{proof}

Now consider the AMS lift $\tilde{\mb F}$ of the outercollaring of ${\mb F}$. In order to define the Floer chain complex using FOP transverse perturbations, one needs to construct the perturbations inductively. 

\begin{lemma}\label{lemma_semipositive_FOP}
There exists an FOP perturbation $\mathring {\mb F}$ of $\tilde{\mb F}$ satisfying the following conditions. Let ${\mc S}_{pq}': {\mc U}_{pq} \to {\mc E}_{pq}$ be the FOP transverse perturbation. If $p<q$ and ${\rm deg}(p) - {\rm deg}(q) \leq 0$, then $({\mc S}_{pq}')^{-1}(0) = ({\mc S}_{pq})^{-1}(0)$ and ${\mc S}_{pq}'$ agrees with ${\mc S}_{pq}$ near the zero locus.
\end{lemma}

\begin{proof}
For the proof, except for the conditions listed in the statement, we further ask that for any $p<q$ such that ${\rm deg}(p) < {\rm deg}(q)$, the FOP transverse perturbation satisfies $\| {\mc S}'_{pq} - {\mc S}_{pq} \|_{C^0} < \frac12 \| {\mc S}_{pq}\|_{C^0}$. We consider the inductive procedure provided in the proof of Theorem \ref{thma_FOP}. Consider a pair $p<q$. Let ${\mc D}_{pq} = ({\mc U}_{pq}, {\mc E}_{pq}, {\mc S}_{pq})$ be the normally complex derived orbifold obtained from the AMS construction.

For the base case of the induction argument, we assume that $(p, q)$ is  minimal, i.e., $M_{pq}^{\mb F}$ does not contain broken trajectories. Then we can choose an FOP transverse perturbation ${\mc S}_{pq}'$ independent of other pairs. In this case, if ${\rm deg}(p) - {\rm deg}(q) <0$, note that ${\mc S}_{pq}$ is nowhere vanishing. Using the $C^0$-density of FOP sections (cf. Theorem \ref{thm_FOP_property}), via a partition-of-unity argument, we can choose a smooth FOP transverse perturbation ${\mc S}'_{pq}$ such that $\| {\mc S}'_{pq} - {\mc S}_{pq} \|_{C^0} < \frac12 \| {\mc S}_{pq}\|_{C^0}$. Then, in particular, ${\mc S}_{pq}$ is nowhere vanishing. If ${\rm deg}(p) = {\rm deg}(q)$, then ${\mc S}_{pq}$ vanishes at finitely many points which correspond to smooth solutions to the Floer equation. In this case, our stable smoothing for such Kuranishi spaces carry the original smooth structure (see Lemma \ref{lemma1312}), the original section ${\mc S}_{pq}$ is smooth near the zero locus and transverse in the classical sense. Moreover, ${\mc U}_{pq}$ is a manifold near its zero locus. Hence {\bf (Classical Transversality)} condition of Theorem \ref{thm_FOP_property}, ${\mc S}_{pq}$ is also FOP transverse near the zero loucs. Hence one can choose ${\mc S}_{pq}'$ to agree with ${\mc S}_{pq}$ near its zero locus without creating more zeroes. If ${\rm deg}(p) - {\rm deg}(q) \geq 1$, we choose an arbitrary FOP transverse perturbation ${\mc S}_{pq}'$.  

Inductively, suppose we have chosen FOP transverse perturbations ${\mc S}_{rs}'$ satisfying our requirement for all pairs $r, s$ such that ${\mc A}_H(r) - {\mc A}_H(s) < {\mc A}_H(p) - {\mc A}_H(q)$. To continue the induction, we analyze the induced perturbation near the boundary strata of ${\mc U}_{pq}$. We first consider the case when ${\rm deg}(p) - {\rm deg}(q) <0$. For any triple $p<r<q$, one of the factors of $M_{pr}^{\mb F} \times M_{rq}^{\mb F} $ has negative virtual dimension. Without loss of generality, we can assume that ${\rm deg}(p) - {\rm deg}(r) <0$. By the inductive hypothesis, we know that $\| {\mc S}_{pr} - {\mc S}_{pr}' \|_{C^0} < \frac12 \| {\mc S}_{pr} \|_{C^0}$. By the compatibility requirement, the putative perturbation ${\mc S}_{pq}'$ is obtained from ${\mc S}_{pr}' \times {\mc S}_{rq}'$ via a stabilization, which implies that the condition $\| {\mc S}'_{pq} - {\mc S}_{pq} \|_{C^0} < \frac12 \| {\mc S}_{pq}\|_{C^0}$ is already met near the boundary stratum $\partial^{prq}{\mc U}_{pq}$. Because this induced section is nowhere vanishing, it is FOP transverse. Then we can apply the $\textbf{(Extension Property)}$ of Theorem \ref{thm_FOP_property} together with the $C^0$-density of FOP transverse perturbations to construct ${\mc S}'_{pq}$ which argees with the induced section near $\partial^{prq}{\mc U}_{pq}$ with $\| {\mc S}'_{pq} - {\mc S}_{pq} \|_{C^0} < \frac12 \| {\mc S}_{pq}\|_{C^0}$. In particular, we see that ${\mc S}'_{pq}$ is nowhere vanishing. 

As for the case ${\rm deg}(p) = {\rm deg}(q)$, if we have $p<r<q$, we can similarly analyze the induced section near $\partial^{prq}{\mc U}_{pq}$ to see that it is nowhere vanishing near $\partial^{prq}{\mc U}_{pq}$. As for the interior part of ${\mc U}_{pq}$, similar to the above discussion, we know that the original section ${\mc S}_{pq}$ vanishes at finitely many points, near which the Kuranishi spaces carry the original smooth structure with respect to which ${\mc S}_{pq}$ is smooth, thereby being an FOP transverse section. So we can choose the perturbation ${\mc S}'_{pq}$ to agree with ${\mc S}_{pq}$ near ${\mc S}^{-1}_{pq}(0)$ and to satisfy the inductive compatibility near the boundary strata such that $({\mc S})^{-1}_{pq}(0) = ({\mc S}')^{-1}_{pq}(0)$. 

Finally, if ${\rm deg}(p) - {\rm deg}(q) \geq 1$, we choose ${\mc S}'_{pq}$ such that it satisfies the compatibility condition near the boundary strata as in the proof of Theorem \ref{thma_FOP}. This finishes the inductive step and we conclude the proof.
\end{proof}

\begin{proof}[Proof of Theorem \ref{thma_semipositive}]
By Lemma \ref{lemma_semipositive_FOP}, there exists an FOP transverse perturbation $\mathring {\mb F}$ such that the associated counts (with respect to the natural orientations) of zeros coincides with the counts of zero-dimensional moduli spaces using the classical method. Then the chain complexes defined using the FOP perturbation agrees with the classical one. 
\end{proof}

\subsection{Proof of Theorem \ref{thma_spectral_invariants}}\label{subsection_proof_spectral}

We verify each item of Theorem \ref{thma_spectral_invariants}. In fact the verification is similar to the case of Hamiltonian Floer spectral invariants defined via the Floer chain complex over ${\mb Q}$ (see for example \cite{FOOO_spectral}). The argument here is a rather straightforward adaptation to the current setting.%

\begin{enumerate}
    \item The {\bf (Spectrality)} follows from Usher's abstract result (Theorem \ref{thm_Usher}). Notice that the argument requires that $CF_* (H, J, \Xi)$ is a Floer--Novikov complex, which is proved in Proposition \ref{prop738}. 

    \item The {\bf (Shifting property)} follows from the fact that the multiplication by $a \in \Lambda_R^\Pi$ is an isomorphism of chain complexes from $CF_*(H, J, \Xi)^{\leq \tau}$ to $CF_*(H, J, \Xi)^{\leq \tau - \mf{v}(a)}$.

    \item The {\bf (Lipschitz Continuity)} is a straightforward consequence of Theorem \ref{thma_Floer_filtered} (2), where we recall that the Hofer distance between Hamiltonians $H_1$ and $H_2$ is
    $$\int_0^1\max_{x \in M}(H_1(x,t) - H_2(x,t))dt - \int_0^1 \min_{x \in M}(H_1(x,t) - H_2(x,t))dt,$$
    which is an upper bound for the threshold (cf. Definition \ref{defn_threshold}) of the continuation maps.

    \item The {\bf (Isotopy invariance)} following from the continuation map construction. Let $H_1, H_2$ induce the same Hamiltonian isotopy $\tilde \phi \in \wham(X, \omega)$ and $J_1$, $J_2$ be two $\omega$-compatible almost complex structure. Then there exists a Floer data on the infinite cylinder connecting $H_1 dt$ and $H_2 dt$ and a family of compatible almost complex structures connecting $J_1$ and $J_2$. Let $\mb{F}_1, \mb{F}_2$ be the associated Floer flow categories. Then the Floer data described here defines a bimodule ${\mb B}_{12}$ from ${\mb F}_1$ to ${\mb F}_2$. Theorem \ref{thma_module_lift} together with the FOP perturbation provides a chain map 
    \beqn
    \Phi^{\mb{B}_{01}}: CF_*(H_1, J_1, \Xi_1) \to CF_*(H_2, J_2, \Xi_2).
    \eeqn
    The Hamiltonian connection can be chosen to be flat precisely because the Hamiltonian isotopies induced from $H_1$ and $H_2$ are isotopic. Hence the threshold of the chain map $\Phi^{\mb{B}_{01}}$ (Definition \ref{defn_threshold}) is zero. The same is true for a continuation map in the reversed direction. Hence we obtain the equality between the spectral numbers. 

    \item For the {\bf (Symplectic invariance)} property, fix a symplectomorphism $\rho: X \to X$. Let $J_\rho$ and $H_\rho$ be the pullback almost complex structure and Hamiltonian. Let ${\mb F}$ resp. ${\mb F}_\rho$ be the Floer flow category associated to $(H, J)$ resp. $(H_\rho, J_\rho)$. Then $\rho$ induces an isomorphism of Floer flow categories $\mb{F} \to \mb{F}_\rho$. The the AMS+FOP construction for ${\mb F}$ can be pulled back to $\mb{F}_\rho$, providing an isomorphic chain complex $CF_*(H_\rho, J_\rho, \Xi_\rho)$. This isomorphism also preserves the symplectic action and the filtration on the chain complexes. 

    \item The {\bf (Triangle inequality)} follows from the chain-level pair-of-pants product construction provided by Theorem \ref{thma_product}. Let $H_1, H_2$ be nondegenerate Hamiltonians generating $\tilde \phi_1, \tilde \phi_2\in \wham(X, \omega)$ such that $\tilde\phi_1 \circ \tilde\phi_2$ is also nondegenerate. Then there is a composition $H_1 \natural H_2$ generating the composition $\tilde\phi_1 \circ \tilde \phi_2$. Moreover, on the pair-of-pants $\Sigma^{\rm pants}$, there exists a flat Hamiltonian connection $\sigma$ which converges to $H_1 dt$, $H_2 dt$ on the two negative ends and converges to $H_3dt:= H_1 \natural H_2 dt$ on the positive ends. Let ${\mb X}$ be the corresponding multimodule, which is an equivariant Novikov multimodule. Let ${\mb F}_i$, $i = 1, 2, 3$ be the Floer flow categories associated to $H_i$. Then it is well-known that the curvature of $\sigma$ gives the threshold of the Novikov multimodule (Definition \ref{defn_threshold}), and the curvature can be taken to be zero. Therefore, by Proposition \ref{prop_spectral_continuity}, one has 
    \beqn
    c_{{\mb F}_3}( {\mf x}_1 \star {\mf x}_2) \leq c_{{\mb F}_1} ({\mf x}_1) + c_{\mb{F}_2} ({\mf x}_2).
    \eeqn

\item The {\bf (Poincar\'e duality)} property follows from the chain-level construction of the Poincar\'e duality and the spectral invariant continuatity result (Proposition \ref{prop_spectral_continuity}). The argument also follows that of \cite[Lemma 5.1]{Ostrover_2006}.

Fix an almost complex structure $J$ and omit it from the argument from now on. Let $H$ be a nondegenerate Hamiltonian on $X$ and ${\mb F}$ be the associated Floer flow category. Let $\ov{H}$ be the reversal of $H$, which has an induced Floer flow category $\ov{\mb F}$. Then on the domain $\Sigma^{\mb{PD}}$ which has two negative ends and no positive end, there exists a flat Hamiltonian connection $\sigma^{\mb{PD}}$ which is $H dt$ on one end and $\ov{H} dt$ on the other. Consider this specific Poincar\'e duality multimodule ${\mb X}^{\mb{PD}}$ (over $({\mb F}, \ov{{\mb F}}; {\mb O})$). Notice that the involved moduli spaces are identical to moduli spaces for the diagonal bimodule $\mb{\Delta}^{\mb{F}{\mb F}}$. 

Notice that ${\mb O}$ contains a special object, or rather a subcategory $\uds{\mb O}$ corresponding to the trivial element of $\Pi$. Let ${\mb X}_0^{\rm PD}$ be the restriction to $\uds {\mb O}$, which is a multimodule over $({\mb F}, \ov{\mb F}; \uds {\mb O})$. Let
\beqn
\Phi_0^{\rm PD}: CF(H) \otimes CF(\ov{H}) \to {\bf K}
\eeqn
be the chain-level pairing. Notice that $\Phi_0^{\rm PD}$ only comes from trivial Floer cylinders, hence is obviously a nondegenerate pairing (in both variables). 

Now we start to prove Poincar\'e duality property of spectral invariants. Fix $x \in HF(X)$. 
\begin{enumerate}

    \item We first prove that 
    \beqn
    c_H(x) \geq - \inf \Big\{ c_{\ov{H}}(y)\ |\ \Phi_0^{\rm PD}(x,y) \neq 0 \Big\}
    \eeqn
    Indeed, the multimodule ${\mb X}_0^{\mb{PD}}$ has a threshold $0$ (see Definition \ref{defn_threshold}). So for any $y \in HF(X)$, if $\Phi_0^{\rm PD}(x,y) \neq 0$, then there exist cycles ${\mf x}$ and ${\mf y}$ representing them such that the chain-level pairing between ${\mf x}$ and ${\mf y}$ is nonzero. Then by Proposition \ref{prop_spectral_continuity}, one has
    \beqn
    {\mc A}_H({\mf x}) + {\mc A}_{\ov{H}}({\mf y}) \geq 0.
    \eeqn
    As ${\mc A}_H({\mf x})$ and ${\mc A}_{\ov{H}}({\mf y})$ can approximate $c_H(x)$ and $c_{\ov{H}}(y)$ respectively, one has
\beqn
\Phi_0^{\rm PD}(x,y) \neq 0 \Longrightarrow c_H(x) + c_{\ov{H}} (y) \geq 0.
\eeqn
Then  one has 
\beqn
c_H(x) \geq \sup \Big\{ - c_{\ov{H}}(y)\ |\ \Phi_0^{\rm PD}(x,y) = 0 \Big\}.
\eeqn

\item We prove that for all $\epsilon>0$,
\beqn
c_H(x) < \epsilon - \inf \Big\{ c_{\ov{H}}(y)\ |\ \Phi_0^{\rm PD}(x,y) \neq 0 \Big\}.
\eeqn
Choose $\epsilon$ such that $\alpha:= \epsilon - c_H(x) \notin {\rm Spec}(\ov{H})$. By energy consideration, the restriction of the chain-level pairing $\Phi_0^{\rm PD}$ to 
\beqn
CF(H)^{<\beta} \otimes CF(\ov{H})^{<\gamma} \to {\bf K}
\eeqn
vanishes when $\beta + \gamma < 0$. Hence one obtains a well-defined nondegenerate pairing
\beqn
\Phi_{0,\alpha}^{\rm PD}: HF(H)^{>-\alpha} \otimes HF(\ov{H})^{<\alpha} \to {\bf K}.
\eeqn 
As $-\alpha = c_H(x) - \epsilon < c_H(x)$, 
\beqn
x \notin {\rm Im} \Big( HF(H)^{<-\alpha} \to HF(X) \Big).
\eeqn
Therefore, the image of $x$ under
\beqn
HF(X) \to HF(H)^{>-\alpha}
\eeqn
denoted by $x_\alpha$, is nonzero. By the nondegeneracy of the pairing $\Phi_{0, \alpha}^{\rm PD}$, there exists $y_\alpha \in HF(\ov{H})^{<\alpha}$ such that $\Phi_{0,\alpha}^{\rm PD}(x_\alpha, y_\alpha) \neq 0$. Choose a cycle ${\mf x} \in CF(H)^{\leq c_H(x)}$ representing $x$ and a cycle ${\mf y} \in CF(\ov{H})^{<\alpha}$ representing $y_\alpha$. Then the chain level pairing $\Phi_0^{\rm PD}({\mf x}, {\mf y})$ is nonzero. As ${\mf x}$ is closed, ${\mf y}$ is not a boundary. Hence the map $HF(\ov{H})^{<\alpha} \to HF(X)$ does not annihilate $y_\alpha$. Denote the image by $y$. Then $\Phi_0^{\rm PD}(x,y) \neq 0$. Then one has 
\beqn
c_{\ov{H}}(y) \leq \alpha = \epsilon - c_H(x)\Longrightarrow c_H(x) \leq \epsilon - c_{\ov{H}}(y) \leq \epsilon - \inf \Big\{ c_{\ov{H}}(y)\ |\ \Phi_0^{\rm PD}(x,y) \neq 0 \Big\}.
\eeqn    
\end{enumerate}

Combining (a) and (b), we conclude that Poincar\'e duality holds.

\end{enumerate}

\newpage

\part{THE AMS CONSTRUCTION}

\section*{Outline of Part 3}

We go through the AMS construction of lifts of various geometric flow categories, multimodules, and homotopies to the category of (collared rigidified) smooth normally complex Kuranishi spaces.

\section{The Flow Categories, Multimodules, and Homotopies of Domains}\label{section_AMS_domains}

The Abouzaid--McLean--Smith (AMS) construction of global charts for genus zero Gromov--Witten theory \cite{AMS} uses certain auxiliary moduli spaces of rational curves in projective spaces.  To extend such construction to the case of Floer moduli spaces, we consider an analogy, which is an auxiliary flow category. The final goal of this section is to describe a flow category enriched in an intermediate category. The discussion in this section and subsequent ones has certain overlap with \cite[Section 5]{Bai_Xu_Arnold}, but we include these expositions to make the paper self-contained and various technical aspects more streamlined.

\subsection{Monotone flow categories, multimodules, and homotopies}

\subsubsection{Monotone flow category}\label{subsubsec:mono-flow-cat}

We introduce the technical notion of monotone flow categories to simplify the discussions. %

We introduce the following normal poset used frequently in this paper. For each $d \geq 0$, define
\beqn
A_d = \left\{ \begin{array}{cc} \Big\{ (d_0, \ldots, d_l)\ |\ d_i > 0,\ d_0 + \cdots + d_l = d \Big\},\ &\ d \neq 0,\\
\{ (0) \},\ &\ d = 0.
\end{array} \right.
\eeqn
The depth function is the number $l$. There are also homogeneous inclusions of degree $1$
\beq\label{flow_poset_inclusion}
A_d \times A_{d'} \to A_{d + d'}
\eeq
defined by the natural concatenation.

\begin{defn}\label{defn_monotone_flow}
Let $\uds {\bf C}$ be a regular stratification category. A {\bf monotone flow category} enriched in $\uds{\bf C}$, denoted by ${\mb f}$, consists of %
for each $d \geq 0$ an object $
M_d^{\mb f}$ in $\uds{\bf C}$  together with morphisms
\beqn
\iota_{dd'}^{\mb f}: M_d^{\mb f} \times M_{d'}^{\mb f} \to M_{d + d'}^{\mb f}
\eeqn
satisfying the following conditions.

\begin{enumerate}
    \item For $d  = 0$, $M_d^{\mb f}$ is a singleton in $\uds{\bf C}$.

\item The associativity: for $d, d', d'' \geq 0$, the following diagram commutes.
\beqn
\xymatrix{ M_d^{\mb f} \times M_{d'}^{\mb f} \times M_{d''}^{\mb f} \ar[r] \ar[d] & M_{d + d'}^{\mb f} \times M_{d''}^{\mb f} \ar[d] \\
 M_d^{\mb f} \times M_{d' + d''}^{\mb f} \ar[r] & M_{ d + d' + d''}^{\mb f} }
 \eeqn
    
\end{enumerate}
\end{defn}

\begin{rem}
Notice that given a monotone flow category $\mb{f}$, one obtains an ordinary flow category, denoted by ${\mb F}$, whose object set is ${\mb Z}$. The partial order is the opposite of the natural one on ${\mb Z}$. The space of morphisms from $a$ to $b$ (with $a \geq b$) is 
\beqn
M_{ab}^{\mb F} = M_{a-b}^{\mb f}. 
\eeqn
Moreover, the flow category ${\mb F}$ has a free action by ${\mb Z}$.
\end{rem}

\begin{rem}
One can generalize the notion of monotone flow categories by replacing ${\mb Z}_{\geq 0}$ with certain other abelian monoids. One option is ${\mb Z}_{\geq 0}^k$, where the $k=2$ case will be used to prove Proposition \ref{thm_integral_action_comparison}. The corresponding flow category has object set being ${\mb Z}^2$. The other is the effective cone
\beqn
\Pi^{\rm eff} = \{ a \in \Pi \ |\ \omega(a) > 0 \} \cup \{0\}
\eeqn
of a symplectic manifold. The corresponding flow category has object set being $\Pi$.
\end{rem}

\subsubsection{Monotone flow multimodules and homotopies}

Now we discuss the case of multimodules.

\begin{defn}\label{defn_monotone_multimodule}
Let $\mb{f}_1, \ldots, \mb{f}_m$ resp. ${\mb f}'$ be monotone flow categories %
enriched in a regular stratification category $\uds{\bf C}$. A {\bf monotone flow multimodule} over $(\mb{f}_1, \ldots, \mb{f}_m; \mb{f}')$, denoted by ${\mb x}$, consists of a  collection of objects $M_d^{\mb x}$ of $\uds{\bf C}$ indexed by $d \geq 0$ (the underlying poset may not be $A_d$),  %
and for $d_j, d, d' \in {\mb Z}_{\geq 0}$ morphisms in $\uds{\bf C}$
\beqn
\iota_{d_j; d}^{\mb x}: M_{d_j}^{{\mb f}_j} \times M_d^{\mb x} \to M_{d_j + d}^{\mb x}\ \ \ \ {\rm resp.}\ \ \ \ \iota_{d; d'}^{\mb x}: M_{d}^{{\mb x}}\times M_{d'}^{{\mb f}'} \to M_{d + d'}^{\mb x}.
\eeqn
They need to satisfy the following conditions.

\begin{enumerate}

\item When $d_j = 0$ resp. $d' = 0$, the morphism $\iota_{d_j; d}^{\mb x}$ resp. $\iota_{d; d'}^{\mb x}$ are the canonical isomorphism.

\item When $d_j > 0$ resp. $d'>  0$, the morphism $\iota_{d_j; d}^{\mb x}$ resp. $\iota_{d; d'}^{\mb x}$ has underlying poset map being homogeneous of degree $1$. 

\item For each $d_j, d, d'\geq 0$, the following diagram commutes.
\beqn
\xymatrix{   M_{d_j}^{{\mb f}_j} \times M_d^{\mb x} \times M_{d'}^{{\mb f}' } \ar[rr] \ar[d] &  &M_{d_j}^{{\mb f}_j} \times M_{d + d'}^{{\mb x}} 
 \ar[d]^{\iota_{d_j; d + d'}^{\mb x}} \\
M_{d_j + d}^{\mb x} \times M_{d'}^{{\mb f}' } \ar[rr]_{\iota_{d_j + d; d'}^{\mb x}}   &  & 
  M_{d_j + d + d'}^{\mb x} }
\eeqn
Hence we may combine the morphisms $\iota_{d_j; d}^{\mb x}$ and $\iota_{d; d'}^{\mb x}$ as one
\beqn
M_{d_j}^{{\mb f}_j} \times M_d^{\mb \xi}\times M_{d'}^{{\mb f}' } \to M_{d_j + d + d'}^{\mb x}.
\eeqn

\item For $d_j, e_j, d, d', e'\geq 0$, the following diagram commutes. 
    \beqn
    \xymatrix{  M_{e_j}^{{\mb f}_j} \times M_{d_j}^{{\mb f}_j} \times M_d^{{\mb x}} \times M_{d'}^{{\mb f}' } \times M_{e'}^{{\mb f}' } \ar[r] \ar[d] & M_{e_j + d_j}^{{\mb f}_j} \times M_d^{{\mb x}}  \times M_{d' + e'}^{{\mb f}' } \ar[d] \\
    M_{e_j}^{{\mb f}_j} \times M_{d_j + d + d'}^{\mb x} \times M_{e'}^{{\mb f}' }  \ar[r] & M_{ e_j + d_j + d + d' + e'}^{\mb x}
     }
    \eeqn

    \item For $i \neq j$, $d_i, d_j \geq 0$ and $d \geq 0$, the following diagram commutes.
    \beqn
    \xymatrix{ M_{d_j}^{{\mb f}_j} \times M_{d_i}^{{\mb f}_i}  \times M_{d}^{\mb x}\ar[rr] \ar[d]  & &   M_{d_j}^{{\mb f}_j} \times M_{d_i + d}^{{\mb x}}\ar[d] \\
    M_{d_i}^{{\mb f}_i} \times M_{d_j + d}^{{\mb x}} \ar[rr]  & & M_{d_i + d_j + d}^{\mb x} }
    \eeqn
    \end{enumerate}
\end{defn}

Notice that the monotone multimodule ${\mb x}$ induces an ordinary multimodule ${\mb X}$ over $({\mb F}_1, \ldots, {\mb F}_m; {\mb F}')$ with 
\beqn
M_{a_1 \cdots a_m; a'}^{\mb X}: =  \left\{ \begin{array}{cc} M_{a_1 + \cdots + a_m - a'}^{\mb x},\ &\ a_1 + \cdots + a_m \geq a',\\
\emptyset,\ &\ {\rm otherwise} \end{array}\right.
\eeqn
and structural maps induced from structural maps of ${\mb x}$. Moreover, it is equivariant with respect to the natural group morphism 
\beqn
\underbrace{{\mb Z} \times \cdots \times {\mb Z}}_{m} \to {\mb Z}.
\eeqn

The case of homotopies is similar and stated as follows.

\begin{defn}\label{defn_monotone_homotopy}
 Let ${\mb f}_1, \ldots, {\mb f}_m, {\mb f}'$ be as in Definition \ref{defn_monotone_multimodule}. Let ${\mb x}_0$, ${\mb x}_1$ be two monotone multimodules over $({\mb f}_1, \ldots, {\mb f}_m; {\mb f}')$ %
A {\bf monotone homotopy} from ${\mb x}_0$ to ${\mb x}_1$, denoted by ${\mb h}$, consists of, for each $d \geq 0$, an object $M_d^{\mb h}$ of $\uds{\bf C}$, a morphism
\beqn
\iota_d^{{\mb x}_0 \to {\mb h}} \sqcup \iota_d^{{\mb x}_1 \to {\mb h}}: M_d^{{\mb x}_0} \sqcup M_d^{{\mb x}_1} \to M_d^{\mb h}
\eeqn
of codimension 1, and for each $d_j\geq 0$ resp. $d' \geq 0$, a morphism
\beqn
\iota_{d_j; d}^{\mb h}: M_{d_j}^{{\mb f}_j} \times M_d^{\mb h} \to X_{d_j + d}^{\mb h}\ \ \ \ {\rm resp.}\ \ \ \ \iota_{d; d'}^{\mb h}: M_d^{\mb h} \times M_{d'}^{\mb{f}'} \to  M_{d + d'}^{\mb h}
\eeqn
which satisfy the following conditions.
\begin{enumerate}
    \item When $d_j = 0$ resp. $d' = 0$, the morphisms $\iota_{d_j; d}^{{\mb h}}$ resp. $\iota_{d; d'}^{\mb h}$ are the canonical isomorphisms.

    \item When $d_j > 0$ resp. $d'>0$, the morphisms $\iota_{d_j; d}^{{\mb h}}$ resp. $\iota_{d; d'}^{\mb h}$ have codimension 1.

    \item For $d_j, d, d'\geq 0$, the following diagram commutes.
    \beqn
    \xymatrix{ M_{d_j}^{{\mb f}_j} \times M_d^{\mb h} \times M_{d'}^{{\mb f}'} \ar[r] \ar[d]   &   M_{d_j + d}^{\mb h}\times M_{d'}^{{\mb f}'} \ar[d] \\
     M_{d_j}^{{\mb f}_j} \times M_{d + d'}^{{\mb h}} \ar[r]     &   M_{d_j + d + d'}^{{\mb h}} }
    \eeqn
    Hence there is a well-defined composition
    \beqn
     M_{d_j}^{{\mb f}_j} \times M_d^{\mb h} \times M_{d'}^{{\mb f}'} \to M_{d_j + d + d'}^{\mb{h}}.
     \eeqn

    \item For $d_j, d, d'\geq 0$, $\alpha = 0, 1$, the following diagram commutes.
    \beqn
    \xymatrix{
    M_{d_j}^{{\mb f}_j} \times M_{d}^{{\mb x}_\alpha} \times M_{d'}^{{\mb f}'} \ar[r] \ar[d]  &    M_{d_j}^{{\mb f}_j} \times M_d^{{\mb h}} \times M_{d'}^{{\mb f}'} \ar[d] \\
    M_{d_j + d + d'}^{{\mb x}_\alpha} \ar[r] & M_{d_j + d + d'}^{\mb h} }
    \eeqn

    \item For $d_j, e_j, d, d', e' \geq 0$, the following diagram commutes.
    \beqn
    \xymatrix{
    M_{e_j}^{{\mb f}_j} \times M_{d_j}^{{\mb f}_j} \times M_d^{{\mb h}} \times M_{d'}^{{\mb f}'} \times M_{e'}^{{\mb f}'} \ar[r] \ar[d]    &    M_{e_j}^{{\mb f}_j} \times M_{d_j + d + d'}^{{\mb h}} \times M_{e'}^{{\mb f}'} \ar[d]\\
    M_{e_j + d_j}^{{\mb f}_j} \times M_d^{{\mb h}} \times M_{d' + d'}^{{\mb f}'} \ar[r]  & M_{e_j + d_j + d + d' + d'}^{\mb h} }
    \eeqn

    \item For $d_i, d_j \geq 0$, $i \neq j$ and $d \geq 0$, the following diagram commutes.
    \beqn
    \xymatrix{  M_{d_i}^{{\mb f}_i} \times M_{d_j}^{{\mb f}_j} \times M_d^{\mb h}  \ar[r] \ar[d]  &  M_{d_i}^{{\mb f}_i} \times M_{d_j + d}^{\mb h} \ar[d]\\
     M_{d_j}^{{\mb f}_j} \times M_{d_i + d}^{\mb h} \ar[r] & M_{d_i + d_j + d}^{\mb h}
    }
    \eeqn
\end{enumerate}
\end{defn}

Again, a monotone homotopy induces an ordinary homotopy ${\mb H}$ from ${\mb X}_0$ to ${\mb X}_1$ with  
\beqn
M_{a_1 \cdots a_m; a'}^{\mb H}:=  \left\{ \begin{array}{cl} M_{a_1 + \cdots + a_m - a'}^{\mb h},&\ a_1 + \cdots + a_m \geq a',\\
\emptyset,\ &\ {\rm otherwise}.
\end{array}\right.
\eeqn

\subsection{The category of curves} 

The following definition may not be the most canonical one. However, as we will consider only very concrete objects, the {\it ad hoc} description suffices our purpose. We should remind the reader that the following definition goes beyond the algebro-geometric situation as the base $B$ is typically only a stratified smooth manifold.

\begin{defn}\label{defn_curve_1}
An {\bf equivariant families of curves} is a triple ${\mc C} = ({\mc G}, B, C)$ where ${\mc G}$ is a (not necessarily reductive) complex Lie group containing a maximal compact subgroup $G\subset {\mc G}$, $B$ and $C$ are stratified smooth ${\mc G}$-manifolds, with a ${\mc G}$-equivariant map $C \to B$, and the structure of a prestable genus zero complex curve on each fiber $C_\phi \subset C$. This triple needs to satisfy the following conditions.
\begin{enumerate}
    \item For any $g \in {\mc G}$ and $\phi \in B$, the map $g: C_\phi \to C_{g\phi}$ is an isomorphism of complex curves.

    \item Let $\mathring C \subset C$ be the complement of fiberwise nodes. Then $\mathring C$ is a stratified ${\mc G}$-manifold and the map $\mathring C \to B$ is a submersion.
\end{enumerate}
\end{defn}

For each object ${\mc C} = ({\mc G}, B, C)$, we introduce
\begin{align*}
&\ P_{\mc C} = {\mc G}/G,\ &\ Q_{\mc C} = {\rm Lie}({\mc G})/ {\rm Lie} G.
\end{align*}

We can also add additional structures to an equivariant family of curves, such as a collection of marked points, which are precisely ${\mc G}$-equivariant sections of $C \to B$ avoiding nodes and each other. Another structure we will use all the time is the cylindrical ends around markings and nodes, which are choices of germs of cylindrical coordinates. 

We define morphisms of the category of equivariant families of curves. 

\begin{defn}\label{defn_curve_2} Let ${\mc C} = ({\mc G}, B, C)$ be an equivariant family of curves.
\begin{enumerate}

\item A {\bf group enlargement} of ${\mc C} = ({\mc G}, B, C)$ is an equivariant family of curves ${\mc C}' = ({\mc G}', B', C')$ where ${\mc G} \to {\mc G}'$ is a complex Lie group embedding which sends the compact group $G$ into $G'$ with $B' = {\mc G}' \times_{\mc G} B$ and $C' = {\mc G}' \times_{\mc G} C$. We denote the group enlargement as such by 
\beqn
{\mc G}'\times_{\mc G} {\mc C}.
\eeqn

\item A {\bf strict morphism} from ${\mc C}_1 = ({\mc G}_1, B_1, C_1)$ to ${\mc C}_2 = ({\mc G}_2, B_2, C_2)$, denoted by $\zeta_{21}: {\mc C}_1 \to {\mc C}_2$, consists of a complex Lie group embedding $\zeta_{21}^{{\mc G}}: {\mc G}_1 \to {\mc G}_2$ which restricts to a group map $\zeta_{21}^G: G_1 \to G_2$ and an equivariant commutative diagram
\beqn
\xymatrix{ C_1 \ar[r]^{\zeta_{21}^C}  \ar[d] & C_2 \ar[d] \\
 B_1 \ar[r]_{\zeta_{21}^B} & B_2 }
\eeqn
satisfying 1) the map $\zeta_{21}^B: B_1 \to B_2$ is a stratified smooth map; 2) the map $\zeta_{21}^C: C_1 \to C_2$ is continuous such that the restricted map $\mathring C_1 \to \mathring C_2$ is smooth and such that for each $\phi_1 \in B_1$ sent to $\phi_2 \in B_2$, the induced map $\mathring C_{\phi_1} \to \mathring C_{\phi_2}$ is an isomorphism (of complex curves together with extra structures such as germs of cylindrical ends and marked points). Notice that a strict morphism ${\mc C}_1 \to {\mc C}_2$ induces a strict morphism from ${\mc G}_2 \times_{{\mc G}_1} {\mc C}_1$ into ${\mc C}_2$. 

\item A strict morphism is called a {\bf strict embedding} if the induced map ${\mc G}_2\times_{{\mc G}_1} {\mc C}_1 \to {\mc C}_2$ is an isomorphism.

\item A {\bf unitary conjugation} of $({\mc G}, B, C)$ is a strict morphism induced from an element $g \in G$ which lifts the conjugation action of $g$ on ${\mc G}$.

\item An {\bf embedding} of equivariant families of curves is a unitary conjugacy class of strict embeddings, where we declare that two strict embeddings ${\mc C}_1 \to {\mc C}_2$ are equivalent if one can be taken to another via a unitary conjugation of $g \in G_2$.

\item The category of equivariant families of curves, denoted by $\uds{\bf Curve}$, is the category whose objects are equivariant families of curves and whose morphisms are embeddings.

\end{enumerate}
\end{defn}

Notice that the category $\uds{\bf Curve}$ is monoidal: given two object ${\mc C}_i = ({\mc G}_i, B_i, C_i)$, the product ${\mc C}_1 \times {\mc C}_2$ is $({\mc G}_1 \times {\mc G}_2, B_1 \times B_2, C_1 \sqcup C_2)$ where the fiber of $C_1 \sqcup C_2$ over $(\phi_1, \phi_2) \in B_1\times B_2$ is the disjoint union $C_{\phi_1}\sqcup C_{\phi_2}$. The initial object, which is $\emptyset$-stratified (see Definition \ref{defn_regular_stratification_category}), is the object with ${\mc G}$ being the trivial group, $B$ being the singleton, and $C$ being empty. 

Notice that the category $\uds{\bf Curve}$ also allows outercollaring.

\subsubsection{Lateral lines}

We introduce the notation of lateral lines, which are used for gluing domains of solutions to Floer equations in various flow-category constructions.

\begin{defn}\label{defn:lateral-line}
\begin{enumerate}

\item Let $\Sigma$ be a smooth genus zero curve with cylindrical ends. A {\bf lateral line} in $\Sigma$ is a real curve ${\mb R} \cong \ell \subset \Sigma$ satisfying: near $\pm \infty$, $\ell$ is contained in a cylindrical end of $\Sigma$ and in cylindrical coordinates, is parametrized by $(s, 0)$. In particular, for domains of Floer trajectories, this amounts to an ${\mb R}$-orbit of the translation action on the cylinder.

\item Let ${\mc C} = ({\mc G}, B, C)$ be an equivariant family of curves equipped with marked points and cylindrical ends. A family of lateral lines is a collection of lateral lines $\ell_\phi \subset C_\phi$ such that $\ell_{g\phi} = g(\ell_\phi)$ for all $\phi \in B$ and $g \in {\mc G}$.
\end{enumerate}
\end{defn}

By abuse of notations, we use $\uds{\bf Curve}$ to denote the category of equivariant family of curves equipped with cylindrical ends and lateral lines. It is easy to see that lateral lines can be concatenated (as we always require the angular coordinate to be $0 \in S^1$), which also guarantees their existence by induction.

\subsubsection{Rigidifications}

To fit into the construction of Kuranishi charts and FOP perturbations, we need to introduce the more refined category $\uds{\bf Curve}_{\rm rig}$. Comparing to other ``rigidified'' morphisms in categories such as $\uds{\bf dOrb}$ or $\uds{\bf Kur}$, the difference between the domain and target of an embedding in $\uds{\bf Curve}$ is caused typically by the enlargement of the complex Lie group. So the rigidification becomes more specific.

\begin{defn}\label{defn_curve_rigidification}
Let $\zeta_{21}: ({\mc G}_1, B_1, C_1) \to ({\mc G}_2, B_2, C_2)$ be a strict embedding in $\uds{\bf Curve}$. A {\bf rigidification} of $\zeta_{21}$ consists of an orthogonal representation $Q_{21}$ of $G_1$ and
a germ of $\zeta_{21}^G$-equivariant maps
\beqn
\theta_{21}: B_1 \times Q_{21}^\epsilon \to B_2
\eeqn
(where $G_1$ acts diagonally on the domain) whose restriction to $B_1 \times \{0\}$ coincides with $\zeta_{21}^B$, such that the map
\beqn
\begin{split}
G_2 \times_{G_1} ( B_1 \times Q_{21}^\epsilon) \to &\ B_2,\\
[g_2, (b_1, \eta_{21})] \mapsto &\ g_2 \theta_{21} (b_1, \eta_{21})
\end{split}
\eeqn
is a diffeomorphism onto an open neighborhood of the image of $\zeta_{21}^B$. A {\bf strict rigidified embedding} is a strict embedding together with a rigidification.
\end{defn}

Notice that the group embedding $\iota_{21}^{\mc G}: {\mc G}_1 \to {\mc G}_2$ induces a linear injection
\beqn
Q_{{\mc C}_1} \to Q_{{\mc C}_2}
\eeqn
and $Q_{21}$ is necessarily isomorphic to its cokernel.

Each group element $g_2\in G_2$ induces a strict rigidified embedding from ${\mc C}_2$ to itself.

\begin{lemma}
Strict rigidified embeddings can be composed.
\end{lemma}

\begin{proof}
Suppose ${\mc C}_i = ({\mc G}_i, B_i, C_i)$, $i = 1, 2, 3$ are three equivariant families. Let $\zeta_{21}: {\mc C}_1 \to {\mc C}_2$ and $\zeta_{32}: {\mc C}_2 \to {\mc C}_3$ be strict embeddings with rigidifications 
\begin{align*}
&\ \theta_{21}: B_1\times Q_{21}^\epsilon \to {\mc G}_2,\ &\  \theta_{32}: B_2 \times Q_{32}^\epsilon \to {\mc G}_3.
\end{align*}
The composed rigidification consists of the representation $Q_{21} \oplus Q_{32}$, the neighborhood $Q_{21}^\epsilon \oplus Q_{32}^\epsilon$, and the map 
\beqn
\begin{split}
\theta_{31}: B_1 \times Q_{21}^\epsilon \times Q_{32}^\epsilon \to &\ {\mc G}_3,\\
(b_1, \eta_{21}, \eta_{32}) \mapsto &\ \theta_{32} \Big( \theta_{21}(b_1, \eta_{21}), \eta_{32} \Big).
\end{split}
\eeqn
The property that the induced map is a diffeomorphism onto a neighborhood of $\zeta_{31}^B(B_1)$ is an obvious consequence of the corresponding properties for the two rigidified embeddings.
\end{proof}

\begin{defn}
$\uds{\bf Curve}_{\rm rig}$ has morphisms being conjugacy classes of strict rigidified embeddings.
\end{defn}

Notice that we can naturally define the collared version of $\uds{\bf Curve}_{\rm rig}$, denoted by 
\beqn
\outer \uds{\bf Curve}_{\rm rig}.
\eeqn

\subsubsection{Stable complex structure}

We would like to define a category $\uds{\bf Curve}^{\mb C}$ of equivariant families of curves whose bases have equivariant stable complex structures. 

\begin{defn}
A {\bf stable complex structure} on an equivariant family of curves ${\mc C} = ({\mc G}, B, C)$ is a $G$-equivariant stable complex structure of the tangent bundle $TB$ (see Definition \ref{defn_stable_complex_structure}). Recall that it consists of a ${\mc G}$-equivariant complex vector bundle $I \to B$ and a ${\mc G}$-equivariant stable isomorphism
    \beqn
    \tau: TB \overset{s}{\to} I.
    \eeqn
\end{defn}

To define morphisms of stably complex objects, note the following fact. Let $\zeta_{21}: ({\mc G}_1, B_1, C_1) \to ({\mc G}_2, B_2, C_2)$ be a strict embedding of equivariant families of curves. Notice that there are ${\mc G}_2$-equivariant vector bundles
    \begin{align*}
    &\ {\mc G}_2\times_{{\mc G}_1} TB_1,\ &\ \pi^* T({\mc G}_2/{\mc G}_1)
    \end{align*}
    where $\pi: {\mc G}_2\times_{{\mc G}_1} B_1 \to {\mc G}_2/{\mc G}_1$ temporarily denotes the natural projection. 
    Via $\zeta_{21}$, one can see that there is a ${\mc G}_2$-equivariant bundle isomorphism
    \beq\label{eqn142}
    TB_2 \cong \pi^* T({\mc G}_2/{\mc G}_1) \oplus ( {\mc G}_2 \times_{{\mc G}_1} TB_1)
    \eeq
    where the first summand is complex. Hence a stable complex structure on ${\mc C}_1$ naturally induces a stable complex structure on ${\mc C}_2$ via the strict embedding $\zeta_{21}$. One can see that the induced one is invariant under unitary conjugations of strict embeddings. Hence one can make the following definition.

\begin{defn}
A strict embedding of equivariant families of curves with stably complex structures from $({\mc G}_1, B_1, C_1, I_1, \tau_1)$ to $({\mc G}_2, B_2, C_2, I_2, \tau_2)$ is a strict embedding $\zeta_{21}$ from $({\mc G}_1, B_1, C_1)$ to $({\mc G}_2, B_2, C_2)$ such that the bundle isomorphism \eqref{eqn142} respects the ${\mc G}_2$-equivariant stable complex structures. The category $\uds{\bf Curve}^{\mb C}$ consists of objects being equivariant families of curves with stable complex structures whose morphisms are unitary conjugacy classes of strict embeddings. 
\end{defn}

One can also define rigidifications for strict embeddings in $\uds{\bf Curve}^{\mb C}$. Indeed, suppose
\beqn
({\mc G}_1, B_1, C_1, I_1, \tau_1),\ ({\mc G}_2, B_2, C_2, I_2, \tau_2) \in {\rm Ob} \uds{\bf Curve}^{\mb C}
\eeqn
and 
\beqn
\zeta_{21}: 
({\mc G}_1, B_1, C_1) \to ({\mc G}_2, B_2, C_2)
\eeqn
is a strict embedding. Suppose
\beqn
\theta_{21}: B_1 \times Q_{21}^\epsilon \to {\mc G}_2
\eeqn
is a rigidification of $\zeta_{21}$ as a strict embedding in $\uds{\bf Curve}$. Consider the image
\beqn
B_2^\epsilon:= {\rm Im}( G_2 \times_{G_1}( B_1 \times Q_{21}^\epsilon) \to B_2)
\eeqn
which is a $G_2$-invariant open subset of $B_2$. There is a projection
\beqn
B_2^\epsilon \to G_2\times_{G_1}B_1.
\eeqn
We say that the rigidification $\theta_{21}$ respects the stable complex structure if over $B_2^\epsilon$, the stable complex structure on $B_2$ is pulled back from its restriction to $G_2\times_{G_1} B_1$ via the above projection. This provides a category $\uds{\bf Curve}_{\rm rig}^{\mb C}$. One also has the corresponding collared version $\outer \uds{\bf Curve}_{\rm rig}^{\mb C}$.

\subsection{Construction of the domain flow categories}

The purpose of the rest of this section is to construct monotone flow categories (Definition \ref{defn_monotone_flow}) %
enriched in the category $\outer \uds{\bf Curve}_{\rm rig}^{\mb C}$, as well as monotone flow bimodules and homotopies.

\subsubsection{Some Lie groups and homogeneous spaces}

We will consider variants of $U(d)$ and $GL(d)$. Inside $PGL(d+1)$, there are two complex subgroups
\beqn
{\mc G}_d^-:= PGL^-(d+1) = \left\{ \left[ \begin{array}{cc} 1 & * \\ 0 & * \end{array}\right] \in PGL(d+1) \right\}
\eeqn
and
\beqn
{\mc G}_d^+:= PGL^+(d+1) = \left\{ \left[ \begin{array}{cc} * & 0 \\ * & 1  \end{array}\right] \in PGL(d+1)  \right\}.
\eeqn
In this paper we will use exclusively the positive version\footnote{This preference is due to the fact that the domains of multimodules always have a unique output. However, for Floer flow categories, using the negative version is also possible.}, hence we simply denote ${\mc G}_d = {\mc G}_d^+$. The unitary group $G_d:= U(d)$ embeds into ${\mc G}_d$ by 
\beqn
U(d) \ni g \mapsto \left[ \begin{array}{cc} g & 0 \\ 0 & 1 \end{array} \right] \in {\mc G}_d.
\eeqn
One needs to consider certain specific group embeddings for different $d$'s. The embedding
\beqn
G_{d_1} \times G_{d_2} = U(d_1) \times U(d_2) \to G_{d_1 + d_2} = U(d_1 + d_2)
\eeqn
is obvious. For the corresponding complex Lie group ${\mc G}_d$, there are corresponding embeddings. For example, for $d_1 = d_2 = 2$, the embedding ${\mc G}_2 \times {\mc G}_2 \to {\mc G}_4$ reads 
\beqn
\left( \left[ \begin{array}{ccc} a_{11} & a_{12} & 0 \\ a_{21} & a_{22} & 0 \\ a_{31} & a_{32} & 1 \end{array}\right], \left[ \begin{array}{ccc} b_{11} & b_{12} & 0 \\ b_{21} & b_{22} & 0 \\  b_{31} & b_{32} & 1 \end{array}\right] \right) \mapsto \left[ \begin{array}{ccccc} a_{11} & a_{12} & 0 & 0 & 0 \\ a_{21} & a_{22} & 0 & 0 & 0 \\
0 & 0 & b_{11} & b_{12} & 0 \\ 0 & 0 & b_{21} & b_{22} & 0 \\ a_{31} & a_{32} & b_{31} & b_{32} & 1 \end{array}\right].
\eeqn

We set up the following notations for later purposes. Denote
\beqn
P_d:= {\mc G}_d / G_d.
\eeqn
The tangent space at the identity coset is the quotient of Lie algebras
\beqn
Q_d:= {\mf g}_d^{\mb C}/ {\mf g}_d.
\eeqn
There are canonical $G_d$-invariant metrics on $Q_d$ making the linear embeddings
\beqn
Q_{d_1} \oplus Q_{d_2} \to Q_{d_1 + d_2}
\eeqn
isometric. Denote 
\beq\label{eqn143}
Q_{d_1, d_2}:= Q_{d_1 + d_2}/ ( Q_{d_1}\oplus Q_{d_2} ).
\eeq

\subsubsection{Moduli spaces of cylinders in projective spaces}

For $d \geq 1$, we will describe an object ${\mc C}_d = ({\mc G}_d, B_d, C_d)$ in $\uds{\bf Curve}$. Consider the moduli space of holomorphic maps 
\beqn
u: {\mb R}\times S^1 \to \mb{P}^d%
\eeqn
of degree $d$, modulo the translation in ${\mb R}$-direction (but not the $S^1$-rotation). Let 
\beqn
\ov{\mc M}{}_{0,2}^{\mb R}(\mb{P}^d, d)
\eeqn
be the natural compactification of this moduli space (cf. \cite[Definition 5.14]{Bai_Xu_Arnold} for an explanation of this notation). Then the group $PU(d+1)$ acts on $\ov{\mc M}{}_{0,2}^{\mb R}(\mb{P}^d, d)$. Let 
\beqn
\ov{\mc M}{}_{0,2}^{{\mb R}, +}({\mb P}^d, d) \subset 
\ov{\mc M}{}_{0,2}^{\mb R} ({\mb P}^d, d)
\eeqn
be the subset satisfying the constraint:
\beq\label{positive_constraint}
u (z_+) = [0, \ldots, 0, 1]
\eeq
which is invariant under $G_d = U(d)$, where $z_+ = \{\infty\} \in \mathbb{P}^1 \supseteq {\mb R}\times S^1$. Let 
\beqn
B_{d} \subset \ov{\mc M}{}_{0,2}^{{\mb R},+}(\mb{P}^{d}, d)
\eeqn
be the open subset of curves $u$ whose image is not contained in any hyperplane of $\mb{P}^d$. It is easy to calculate the dimension of $B_d$, which is 
\beqn
{\rm dim}_{\mb R} B_d = 2d(d+1)- 1 = {\rm dim}_{\mb R} {\mc G}_d -1.
\eeqn
Further, let $C_d \to B_d$ be the universal curve; notice that fibers also carry canonical germs of cylindrical ends.

The moduli space $B_{d}$ is stratified by the poset $A_d$ (cf. the notations introduced in Section \ref{subsubsec:mono-flow-cat}): for any $\tau  = (d_0 \cdots d_r)\in A_d$, the stratum $\partial^\tau B_{d}$ is the subspace of curves in $\mb{P}^d$ which breaks into $r+1$ stable cylinders such that the $i$-th one has degree $d_i$. Therefore, one obtains an object
\beqn
{\mc C}_{d} = ({\mc G}_d, B_d, C_d)
\eeqn
in the category $\uds{\bf Curve}$. 

\begin{rem}
There is a refinement of the above stratification, namely, the stratification by combinatorial types, and each stratum is preserved by the ${\mc G}_d$-action. 
\end{rem}

\subsubsection{Orientations}

The smooth part of the moduli space $\ov{\mc M}_{0,2}({\mb P}^d, d)$ is a complex manifold, which is canonically oriented. The interior of $B_d$ has a natural $S^1$-action by reparametrizing the cylinders so that the quotient ${\rm Int} B_d/S^1$ is identified with an open subset of $\ov{\mc M}_{0,2}({\mb P}^d, d)$. Hence we can orient $B_d$ by regarding the infinitesimal $S^1$-direction as positively oriented.

\subsubsection{Concatenation of cylinders}

We would like to define the structural maps, i.e., morphisms 
\beqn
\zeta_{d, d'}: {\mc C}_{d} \times {\mc C}_{d'} \to {\mc C}_{d + d'},\ d, d' > 0
\eeqn
in the category $\uds{\bf Curve}$ such that they satisfy associativity. In this way, one obtains a monotone flow category enriched in $\uds{\bf Curve}$. 

We describe the definition as follows. Given two points $\phi \in B_d$, $\phi' \in B_{d'}$ represented by maps $u: \Sigma \to {\mb P}^d$, $u': \Sigma' \to {\mb P}^{d'}$, let $\Sigma \# \Sigma'$ be the prestable cylinder obtained by identifying the output of $\Sigma$ with the input $\Sigma'$. The map $u$ can be represented by sections
\beqn
h_0, \ldots, h_d \in H^0(\Sigma, u^* {\mc O}_{{\mb P}^{d}}(1)).
\eeqn
We abbreviate $L_u:= u^* {\mc O}_{{\mb P}^d}(1)$ and extend it by trivial bundle to $\Sigma'$. Then $h_0, \ldots, h_d$ are sections over $\Sigma \# \Sigma'$. Similarly, the map $u'$ can be represented by sections $h_0', \ldots, h_{d'}'$ of the pullback bundle $L_{u'} \to \Sigma'$ which can be extended to the whole of $\Sigma \# \Sigma'$. 

For each $\phi' \in B_{d'}$, we choose the lateral line, a (broken) real curve
\beq\label{lateral_line}
\ell_{\phi'} \subset C_{\phi'} 
\eeq
which is parametrized by $(s, 0)$ in the cylindrical coordinates in all cylindrical components. Along $\ell_{\phi'}$, the parallel transport with respect to the Chern connection on $L_{u'}$ sends the (nonzero) vector $h_{d'}(z_+) \in L_{u'}|_{z_+}$ to a nonzero vector $\tau_{u'} \in L_{u'}|_{z_-}$. Then consider $1 + d + d'$ sections of $L_u \otimes L_{u'} \to \Sigma \# \Sigma'$ given by 
\beqn
h_0 \otimes \tau_{u'}, \ldots, h_{d-1} \otimes \tau_{u'}, h_d \otimes h_0', \ldots, h_d \otimes h_{d'}.
\eeqn
It is easy to check that it is base-point-free, hence defines a map 
\beqn
(\phi \# \phi'): \Sigma \# \Sigma' \to \mb{P}^{d + d'}
\eeqn
which represents a point in $B_{d + d'}$. Moreover, it is straightforward to see that the equivalence class of $u\# u'$ is independent of the choice of the representatives $u$ and $u'$. Therefore, one obtains a well-defined map 
\beqn
\zeta_{d,d'}^{B}: B_d \times B_{d'} \to B_{d + d'}.
\eeqn

\begin{lemma}
The map $\zeta_{d, d'}^B: B_d  \times B_{d'}  \to B_{d + d'} $ is smooth and equivariant with respect to the group embedding 
\beqn
\zeta_{d, d'}^{\mc G}: {\mc G}_{d}  \times {\mc G}_{d'}  \to {\mc G}_{d + d'} ,
\eeqn
and satisfying the associativity, i.e., the following diagram commutes.
\beqn
\xymatrix{  B_{d_1}  \times B_{d_2}  \times B_{d_3}  \ar[r] \ar[d]   & B_{d_1 + d_2} \times B_{d_3} \ar[d] \\
B_{d_1}  \times B_{d_2 + d_3}  \ar[r]   &  B_{d_1 + d_2 + d_3}   }
\eeqn
Furthermore, the maps $\zeta_{d, d'}^B$ can be lifted canonically to a strict morphism
\beqn
 \zeta_{d, d'}: {\mc C}_{d}  \times {\mc C}_{d'}  \to {\mc C}_{d + d'} 
\eeqn
in the category $\uds{\bf Curve}$  such that the following diagram commutes.
\beqn
\xymatrix{    {\mc C}_{d_1}  \times {\mc C}_{d_2}  \times {\mc C}_{d_3}  \ar[r]\ar[d]    & {\mc C}_{d_1 + d_2}  \times {\mc C}_{d_3}  \ar[d]\\
{\mc C}_{d_1}  \times {\mc C}_{d_2 + d_3}  \ar[r] & {\mc C}_{d_1 + d_2 + d_3}  }
\eeqn
\end{lemma}

\begin{proof}
All elements in the definition depend smoothly on the variables $\phi$ and $\phi'$. From the explicit definition of $u\#u'$ given above, the equivariance follows immediately. The associativity can also be checked directly. See \cite[\S 5.2.4]{Bai_Xu_Arnold} for more details.
\end{proof}

We need to prove that the morphisms are strict embeddings into the corresponding stratum.

\begin{prop}
The morphism $\zeta_{d, d'}$ is a strict embedding into the stratum $\partial^{d, d'} {\mc C}_{d + d'} $.
\end{prop}

\begin{proof}
What we need to show is that the induced strict morphism 
\beqn
{\mc G}_{d + d'}  \times_{{\mc G}_d  \times {\mc G}_{d'} } \left( {\mc C}_{d}  \times {\mc C}_{d'}   \right) \to \partial^{d, d'} {\mc C}_{d + d'} 
\eeqn
is an isomorphism. Indeed, given any curve in $\partial^{d, d'} B_{d + d'} $, it is represented by a holomorphic map $u: \Sigma \# \Sigma' \to \mb{P}^{d + d'}$ such that the restriction of $u$ to $\Sigma$ resp. $\Sigma'$ has degree $d$ resp. $d'$. After applying a transformation in ${\mc G}_{d + d'} $, one can make $u(\Sigma)$ be contained in the special ${\mb P}^d \subset {\mb P}^{d + d'}$ with last $d'$ coordinates vanish, and make $u(\Sigma')$ be contained in a ${\mb P}^{d'} \subset {\mb P}^{d + d'}$ which is spanned by the last $d'$ coordinates as well as a line contained in the special ${\mb P}^d$. Such curves are exactly in the image of $\zeta_{d, d'}^B$. 
\end{proof}

\subsubsection{Rigidification}

Lastly we construct the rigidifications (see Definition \ref{defn_curve_rigidification}).

\begin{prop}\label{prop_curve_rigidification}
The monotone flow category $\outer \mb{Dom}$ by taking the outer-collaring of the monotone flow category which has objects ${\mc C}_d$ and morphisms $\zeta_{d, d'}$ admits a lift to $\outer \uds{\bf Curve}_{\rm rig}$. 
\end{prop}

\begin{proof}
We first explain why a single embedding $\zeta_{d_1, d_2}: {\mc C}_{d_1} \times {\mc C}_{d_2} \to \partial^{d_1, d_2} {\mc C}_{d_1 + d_2}$ can be rigidified (before outer-collaring). (See Definition \ref{defn_curve_rigidification}) the representation of $G_1 = G_{d_1}\times G_{d_2}$ is $Q_{d_1, d_2}$ given by \eqref{eqn143}. Roughly, $Q_{d_1, d_2}$ is the normal direction to $G_{d_1 + d_2}( {\mc G}_{d_1} \times {\mc G}_{d_2})$ inside ${\mc G}_{d_1 + d_2}$. To be more specific, one can choose a $G_d$-invariant germ of maps
\beqn
\theta_{d_1, d_2}: Q_{d_1, d_2} \to {\mc G}_{d_1 + d_2}
\eeqn
sending $0$ to the identity (such as the restriction of an exponential map) which is orthogonal to $G_{d_1+d_2}( {\mc G}_{d_1} \times {\mc G}_{d_2})$. 

By the definition of embeddings in $\uds{\bf Curve}$, there is an isomorphism
\beqn
{\mc G}_{d_1+d_2} \times_{{\mc G}_{d_1} \times {\mc G}_{d_2}} ({\mc C}_{d_1}\times {\mc C}_{d_2} ) \to \partial^{d_1, d_2} {\mc C}_{d_1 + d_2}.
\eeqn
Then the map 
\beqn
\begin{split}
Q_{d_1, d_2} \times (B_{d_1} \times B_{d_2}) &\ \to \partial^{d_1, d_2} B_{d_1 + d_2}\\
(v, x_1, x_2) & \ \mapsto \theta_{d_1, d_2}(v) (\zeta_{d_1, d_2}^B (x_1, x_2))
\end{split}
\eeqn
is a local diffeomorphism after equivariantize by $G_{d_1 + d_2}$.

We need to choose rigidifications for all such embeddings so that they satisfy the requirement for flow categories. For this one can allow $\theta_{d_1, d_2}$ to depend on $(x_1, x_2) \in B_{d_1} \times B_{d_2}$ so $\theta_{d_1, d_2}$ is a germ of equivariant bundle maps
\beqn
Q_{d_1, d_2} \times B_{d_1} \times B_{d_2} \to {\mc G}_{d_1 + d_2} \times B_{d_1} \times B_{d_2}
\eeqn
over $B_{d_1} \times B_{d_2}$. All such rigidifications of the embedding $\zeta_{d_1, d_2}$ are equivariantly homotopic. The collaring region allows us to connect different choices.
\end{proof}

We choose such a lift of $\outer \mb{Dom}$ to $\outer \uds{\bf Curve}_{\rm rig}$ and still denote it by $\outer \mb{Dom}$. The effect of such a lift to later constructions needs to be addressed in the following way. Roughly, the resulting Kuranishi spaces of the AMS construction have symmetry groups being compact Lie groups $G$ while the domains we constructed here have symmetry groups being corresponding complex Lie groups ${\mc G}$. This redundance will be removed by using a notion called {\it group reduction} (see Definition \ref{defn_group_reduction}), which later becomes part of the Kuranishi section. 

Consider a strict embedding 
\beqn
\zeta_{21}: {\mc C}_1 = ({\mc G}_1, B_1, C_1) \to {\mc C}_2 = ({\mc G}_2, B_2, C_1)
\eeqn
of equivariant families of curves. Denote by
\beqn
\zeta_{21}^P: P_1 = {\mc G}_1/ G_1 \to P_2 = {\mc G}_2/ G_2
\eeqn
the induced embedding of homogeneous spaces and 
\beqn
\zeta_{21}^Q: Q_1 = {\mf g}_1^{\mb C} / {\mf g}_1 \to Q_2 = {\mf g}_2^{\mb C}/ {\mf g}_2
\eeqn
the linearization of $\zeta_{21}^P$. 

\begin{lemma}
Let $\outer \mb{Dom}$ be the monotone flow category enriched in $\outer \uds{\bf Curve}_{\rm rig}$ we constructed above. There exist a collection of germs of $G_d\times G_d$-equivariant  maps
\beqn
\lambda_d: \outer B_d \times Q_d \to \outer B_d\times P_d
\eeqn
as bundle maps over $\outer B_d$ satisfying the following conditions.
\begin{enumerate}

\item For each $\phi\in \outer B_d$, the restriction of $\lambda_d$ to $\{\phi\}\times Q_d$ is a local diffeomorphism onto an open neighborhood of $[1] \in P_d$.

    \item $\lambda_d$ is collared, i.e., $\lambda_d$ is determined over the collar neighborhood by its restriction to the boundary strata via the projections.

    \item For each decomposition $d = d_1 + d_2$, the following diagram commutes.
    \beqn
    \xymatrix{ \big( \outer B_{d_1} \times Q_{d_1} \big) \times \big(\outer B_{d_2} \times Q_{d_2} \big) \ar[rr]^-{\lambda_{d_1}\times \lambda_{d_2}} \ar[d]_{\zeta_{d_1,d_2}^B \times \zeta_{d_1, d_2}^Q} & & \big( \outer B_{d_1} \times P_{d_1}\big)\times \big( \outer B_{d_2} \times P_{d_2}\big) \ar[d]^{\zeta_{d_1, d_2}^B \times \zeta_{d_1, d_2}^P} \\
            \partial^{d_1, d_2} \outer B_d \times Q_d \ar[rr]_-{\lambda_{d_1 + d_2}}   &  & \partial^{d_1, d_2} \outer B_d \times P_d }
    \eeqn

\item For each decomposition $d = d_1 + d_2$, $x_1 \in \outer B_{d_1}$, $x_2 \in \outer B_{d_2}$, $q_1 \in Q_{d_1}$, $q_2 \in Q_{d_2}$, and $v \in Q_{d_1, d_2}$ with $|v|$ sufficiently small, one has 
\beq\label{eqn146}
\lambda_d ( \zeta_{d_1, d_2}^B (x_1, x_2), (q_1, q_2, v) ) = (1, \theta_{d_1, d_2}(x_1, x_2, v)) \cdot \Big( \zeta_{d_1, d_2}^B(x_1, x_2), \zeta_{d_1, d_2}^P ( \lambda_{d_1}(x_1, q_1), \lambda_{d_2}(x_2, q_2)) \Big)
\eeq
where $\cdot$ denots the ${\mc G}_d \times {\mc G}_d$-action on $B_d \times P_d$.

\end{enumerate}
\end{lemma}

\begin{proof}
The proof, which is similar to that of Proposition \ref{prop_curve_rigidification}, is based on induction and homotopies over the collar region. For $d = 1$, choose an arbitrary $\lambda_1$ which realizes the local diffeomorphism $Q_1$ with $P_1$. The formula \eqref{eqn146} allows us to define $\lambda_d$ on boundary. The condition that $\lambda_d$ is collared implies that the boundary values of $\lambda_d$ determine its values near the boundary. One can then extend to the interior of $\outer B_d$. 
\end{proof}

\subsubsection{Stable complex structure}

As the first step towards constructing the stable complex structure for the AMS construction, we have the following ``natural'' stable complex structure on the monotone flow category of curves. 

\begin{prop}\label{prop_domain_stable_complex}
The monotone flow category $\outer {\mb{Dom}}$ (enriched in $\outer \uds{\bf Curve}_{\rm rig}$) admits a natural lift to $\outer \uds{\bf Curve}_{\rm rig}^{\mb{C}}$.
\end{prop}

\begin{proof}
Notice that the smooth part $\ov{\mc M}{}_{0,2}^* (\mb{CP}^d, d) \subset \ov{\mc M}{}_{0,2}(\mb{CP}^d, d)$ has a $PGL(d+1)$-invariant complex structure. The constraint at the negative or the positive marking does not affect this complex structure. Notice that the interior of $B_d$ has a projection map $B_d \to \ov{\mc M}{}_{0,2}(\mb{CP}^d, d)$ whose fibers are the $S^1$-rotations (which has no fixed point); moreover, the ${\mc G}_d$-action sends $S^1$-orbits to $S^1$-orbits. Hence we can define $I_{B_d} \to B_d$ as the pullback of the tangent bundle of $\ov{\mc M}{}_{0,2}^* (\mb{CP}^d, d)$ so that there is a ${\mc G}_d$-equivariant vector bundle isomorphism
\beqn
TB_d \oplus \uds{\mb R} \cong I_{B_d}
\eeqn
where the trivial bundle $\uds {\mb R}$ is the subbundle coming from the $S^1$-orbits. 

Now we consider the structural maps. Given $d_1, d_2\geq 1$, recall one has the map
\beqn
\zeta_{d_1, d_2}^B: B_{d_1} \times B_{d_2} \to \partial^{d_1, d_2} B_{d_1 + d_2}.
\eeqn
The differential reads 
\beqn
TB_{d_1} \oplus TB_{d_2} \to T \partial^{d_1, d_2} B_{d_1 + d_2}.
\eeqn
On the other hand, one can see from the construction of $\zeta_{d_1, d_2}^B$ that it induces a complex analytic embedding
\beqn
\ov{\mc M}{}_{0,2}^*(\mb{CP}^{d_1}, d_1) \times \ov{\mc M}{}_{0,2}^*(\mb{CP}^{d_2}, d_2) \to \ov{\mc M}{}_{0,2}^* (\mb{CP}^{d_1 + d_2}, d_1 + d_2).
\eeqn
Therefore, the map $\zeta_{d_1, d_2}^B$ respects the corresponding stable complex structures.
\end{proof}

\subsubsection{Orientation and concatenation}

Recall that each $B_d$ is oriented in a canonical way. Hence for $d = d_1 + d_2$, there is a product orientation on $B_{d_1} \times B_{d_2}$. As ${\mc G}_d$ is a complex Lie group, ${\mc G}_d\times_{{\mc G}_{d_1}\times {\mc G}_{d_2}} (B_{d_1}\times B_{d_2})$ has the induced product orientation. On the other hand, it is a codimension one boundary of $B_d$, hence has a corresponding boundary orientation.

\begin{lemma}\label{lemma_domain_orientation}
The product orientation and the boundary orientation coincide.
\end{lemma}

\begin{proof}
We only consider the interior of the moduli spaces $B_{d_1}$ and $B_{d_2}$. Recall that $B_{d_1}/S^1$ and $B_{d_2}/S^1$ are open subsets of complex manifolds. There are also orientation-preserving gluing maps
\beqn
B_{d_1}/S^1 \times \Delta \times B_{d_2}/S^1 \to B_d/S^1
\eeqn
where $\Delta$ is a small disk of gluing parameters. By explicitly writing down the gluing map in terms of representative and the effect of reparametrizations, one can see that the induced map
\beqn
B_{d_1}\times [0, \epsilon) \times B_{d_2} \to B_d
\eeqn
is orientation preserving. Because all $B_d$ are odd-dimensional, our convention on induced boundary orientation leads to the conclusion.
\end{proof}

\subsection{Multimodules and homotopies for smooth domains}\label{subsection_domain_multimodule}

In this subsection, we describe the multimodules of domain curves which will be used for the AMS construction for moduli spaces for pair-of-pants products. 

We would like to specify a class of constructions as follows. Fix a common outercollaring width. Choose (independently) lifts of the outercollared domain flow category $\outer \mb{Dom}_1, \ldots, \outer \mb{Dom}_m, \outer \mb{Dom}'$ to $\outer \uds{\bf Curve}_{\rm rig}^{\mb C}$.\footnote{When we treat the cigar bimodule and the Poincar\'e duality multimodule, one does not need $\outer \uds{\mb{Dom}}'$. However the method can be easily modified to suit those cases.}

\begin{enumerate}
    \item For any smooth genus zero curve $\Sigma$ with $m$ negative cylindrical ends and one positive cylindrical end, we construct a class of monotone multimodule $\outer \mb{Dom}^\Sigma$ over $(\outer \mb{Dom}_1, \ldots, \outer \mb{Dom}_m; \outer \mb{Dom}')$ enriched in $\outer \uds{\bf Curve}_{\rm rig}^{\mb C}$.

    \item Let $\Sigma_t$, $t \in [0, 1]$ be a smooth family of smooth genus zero curves with $m$ negative cylindrical ends and one positive cylindrical end. Let $\outer \mb{Dom}^{\Sigma_0}$, $\outer \mb{Dom}^{\Sigma_1}$ be two monotone multimodules as above. We construct a class of homotopies $\outer \mb{Dom}^{\Sigma_I}$ from $\outer \mb{Dom}^{\Sigma_0}$ to $\outer \mb{Dom}^{\Sigma_1}$ enriched in $\outer \uds{\bf Curve}_{\rm rig}^{\mb C}$.    
    \end{enumerate}

We first describe the underlying posets of the monotone multimodule. For each $d$, define
\beqn
A_d^\Sigma:= \left\{ (r_1, \ldots, r_m; d_0; r')\ \left| \begin{array}{c} r_j \in A_{d_j},\ j = 1, \ldots, m,\ r' \in A_{d'},\ d_j, d'  > 0,\ d_0 \geq 0,\\
d_1 + \cdots + d_m + d_0 + d' = d 
\end{array} \right.\right\}.
\eeqn
Then there are the natural poset maps
\beqn
\prod_{j=1}^m A_{d_j} \times A_d^\Sigma \times A_{d'} \to A_{d_1 + \cdots + d_m + d + d'}^\Sigma.
\eeqn
These are the poset maps underlying the structural maps of the monotone multimodule we will define. 

Now we describe, for each $d \geq 0$, an object ${\mc C}_d^\Sigma = ({\mc G}_d^\Sigma, B_d^\Sigma, C_d^\Sigma) \in \uds{\bf Curve}$. We set ${\mc G}_d^\Sigma = PGL^+(d+1)$ and $G_d^\Sigma = U^+(d)$, the same as before. For the space $B_d^\Sigma$, consider the moduli space of holomorphic maps $u: \Sigma \to \mb{P}^d$ of degree $d$. Let 
\beqn
{\mc M}^\Sigma({\mb P}^d, d)
\eeqn
be the subspace satisfying the positive constraint
\beqn
u  (z_+) = [0, \ldots, 0, 1].
\eeqn
Let $\ov{\mc M}{}^\Sigma({\mb P}^d, d)$ be the natural compactification (where cylindrical bubbles are taken modulo only translation but not rotation). Let 
\beqn
B^\Sigma_d \subset \ov{\mc M}{}^\Sigma(\mb{P}^d, d)
\eeqn
be the subset of maps whose images are not entirely contained in any hyperplane. Notice that $B_d^\Sigma$ is naturally stratified by the poset $A_d^\Sigma$. Moreover, ${\mc G}_d^\Sigma$ acts on $B_d^\Sigma$. Let $C_d^\Sigma \to B_d^\Sigma$ be the universal curve, which is a ${\mc G}_d^\Sigma$-equivariant family. Define
\beqn
{\mc C}_d^\Sigma = ({\mc G}_d^\Sigma, B_d^\Sigma, C_d^\Sigma)
\eeqn
which is an object in $\uds{\bf Curve}$.

\subsubsection{Structural maps}

To describe the structural maps, we need to choose the corresponding lateral lines on $\Sigma$.\footnote{Notice that there are homotopically non-equivalent choices of lateral lines. The comparison between two such choices of lateral lines on $\Sigma$ will be given in Subsection \ref{subsection_comparison_lateral_line}.} We also need to restrict to the case that there is only one output (in the positive version). The morphisms we need to describe are
\beqn
\prod_{j=1}^m {\mc C}_{d_j} \times {\mc C}_d^\Sigma  \to {\mc C}_{d_1 + \cdots + d_m + d}^\Sigma.
\eeqn
and 
\beqn
{\mc C}_d^\Sigma \times {\mc C}_{d'} \to {\mc C}_{d +d'}^\Sigma.
\eeqn
There is not much difficulty in defining the second map, which is similar to the previous case. On the other hand, there is a subtlety in defining the first structural map. In fact, the underlying group map is not well-defined; it depends on the ordering of the cylindrical end. However, different choices are related by conjugations. 

Assume for simplicity $l = 2$, which is the case of the pair-of-pants product.  For any $d_1, d_2$, we define
\beqn
{\mc G}_{d_1} \times {\mc G}_{d_2} \times {\mc G}_d^\Sigma  \to   {\mc G}_{d_1 + d_2 + d}^\Sigma
\eeqn
to be the same morphism described before. More explicity, when $d^1 = d^2 = d = 2$, the map is
\begin{multline*}
\left( \left[ \begin{array}{ccc} a_{11} & a_{12} & 0 \\ a_{21} & a_{22} & 0 \\
a_{31} & a_{32} & 1  \end{array} \right], \left[ \begin{array}{ccc} b_{11} & b_{12} & 0 \\  b_{21} & b_{22} & 0 \\
b_{31} & b_{32} & 1  \end{array} \right], \left[ \begin{array}{ccc} c_{11} & c_{12} & 0 \\ c_{21} & c_{22} & 0 \\ c_{31} & c_{32} & 1  \end{array} \right] \right) \\
\mapsto \left[ \begin{array}{ccccccc} a_{11} & a_{12} & 0 & 0 & 0 & 0 & 0 \\
a_{21} & a_{22} & 0 & 0 & 0 & 0 & 0 \\
0 & 0 & b_{11} & b_{12} & 0 & 0 & 0 \\
0 & 0 & b_{21} & b_{22} & 0 & 0 & 0 \\
0 & 0 & 0 & 0 & c_{11} & c_{12} & 0 \\
0 & 0 & 0 & 0 & c_{21} & c_{22} & 0 \\
a_{31} & a_{32} & b_{31} & b_{32} & c_{31} & c_{32} & 1 
\end{array} \right].
\end{multline*}

Similarly, we also define the concatenation map. Consider the concatenation at the two negative ends of $\Sigma$. Let $\phi^j \in B_{d^j}$, $j = 1, 2$, be represented by maps $u^j: \Sigma^j \to \mb{P}^{d^j}$. Then they pullback line bundles $L_{u^1}, L_{u^2} \to \Sigma$ with frames
\beqn
(h_0^1, \ldots, h_{d^1}^1),\ (h_0^2, \ldots, h_{d^2}^2).
\eeqn
The positive constraint is that
\begin{align*}
&\ h_0^1(z_+) = \cdots = h_{d^1-1}^1(z_+) = 0,\ &\  h_0^2(z_+) = \cdots = h_{d^2-1}^2(z_+) = 0.
\end{align*}
On the other hand, suppose the point $\phi \in B_d^\Sigma$ is represented by a map $u: \Sigma \to {\mb P}^d$, which pulls back a bundle $L_u \to \Sigma$ having a frame
\beqn
h_0, \ldots, h_d
\eeqn
satisfying 
\beqn
h_0(z_+) = \cdots = h_{d-1} (z_+) = 0.
\eeqn
The parallel transport along the two lateral lines induces nonzero vectors $\tau^1 \in L_u(z_-^1)$, $\tau^2 \in L_u(z_-^2)$. Then consider the prestable curve 
\beqn
\tilde \Sigma:= (\Sigma^1 \sqcup \Sigma^2) \# \Sigma.
\eeqn
We extend the bundles $L_{u^1} \to \Sigma^1$, $L_{u^2} \to \Sigma^2$, and $L_u \to \Sigma$ to other components by trivial bundle extension. Then all the sections extend by constants. Then define $L = L_{u^1} \otimes L_{u^2} \otimes L_u$ and consider the collection of sections (written in three groups)
\beqn
\Big( h_0^1 \otimes h_{d^2}^2 \otimes \tau^1, \ldots, h_{d^1-1}^1 \otimes h_{d^2}^2 \otimes \tau^1; h_{d^1}^1 \otimes h_0^2 \otimes \tau^2, \ldots, h_{d^1}^1 \otimes h_{d^2-1}^2 \otimes \tau^2; h_{d^1}^1 \otimes h_{d^2}^2 \otimes h_0, \ldots, h_{d^1}^1 \otimes h_{d^2}^2 \otimes h_d \Big).
\eeqn
Then it is easy to see that this collection forms a basis of $H^0(L)$. In particular, it is base-point-free. Hence it defines a holomorphic map
\beqn
\tilde u: \Sigma \to \mb{P}^{d^1 + d^2 + d}.
\eeqn
One can also see that this is independent of the choices of the representatives. Hence one obtains a well-defined map
\beqn
\zeta_{d^1, d^2; d}^B: B_{d^1} \times B_{d^2} \times B_d^\Sigma \to B_{d^1 + d^2 + d}^\Sigma.
\eeqn
It is easy to see that the map is equivariant with respect to the group embedding.

The associativity is more subtle than the previous case. Choose $d^1_a, d^1_b, d_a^2, d_b^2, d$. Denote $\tilde d = d_b^1 + d_b^2 + d$. Let's compare the following two elements
\beqn
\zeta_{d^1_a + d_b^1, d^2_a + d_b^2; d}^B \left( \zeta_{d_a^1, d_b^1}^B(\phi^1_a, \phi^1_b), \zeta_{d_a^2, d_b^2}^B (\phi^2_a, \phi^2_b);  \phi \right)\ {\rm vs.}\ \zeta_{d_a^1, d_a^2; \tilde d}^B \left( \phi_a^1, \phi_a^2; \zeta_{\phi_b^1, \phi_b^2; d}^B (\phi_b^1, \phi_b^2; \phi) \right).
\eeqn
They indeed differ by a conjugation. In fact, the diagram
\beqn
\xymatrix{  U(d_a^1) \times U(d_b^1) \times U(d_a^2)\times U(d_b^2) \times U(d)  \ar[r] \ar[d] &  U(d_a^1 + d_b^1) \times U(d_a^2 + d_b^2) \times U(d)\ar[d]\\
U(d_a^1) \times U(d_a^2) \times U(\tilde d) \ar[r] & U(d_a^1 + d_a^2 + d_b^1 + d_b^2 + d )}
\eeqn
commutes up to an element in $U(d_a^1 + d_a^2 + d_b^1 + d_b^2 + d)$ which is the matrix permuting corresponding coordinates. Therefore, the above diagram commutes up to a conjugation. 

One can extend the above construction to the general case when $\Sigma$ is smooth with multiple inputs.

\begin{prop}
The collections of equivariant family of curves ${\mc C}_{d}^\Sigma$ together with the structural maps
\beqn
{\mc C}_{d_1} \times \cdots \times {\mc C}_{d_m} \times {\mc C}_d^\Sigma \times {\mc C}_{d'} \to {\mc C}_{d_1+ \cdots + d_m + d + d'}^\Sigma
\eeqn
define a monotone multimodule over $(\mb{Dom}, \ldots, \mb{Dom}; \mb{Dom})$.
\end{prop}

\subsubsection{Rigidification and stable complex structure}

We can construct a rigidification of any outercollaring of the multimodule $\mb{Dom}^{\Sigma}$ in the same way as constructing the rigidification of the domain flow category (see Proposition \ref{prop_curve_rigidification}). 

\begin{prop}\label{prop_domain_multimodule}
Given lifts of $\outer \mb{Dom}_1, \ldots, \outer \mb{Dom}_m, \outer \mb{Dom}'$ to $\outer \uds{\bf Curve}_{\rm rig}^{\mb C}$, there exists a lift of $\outer \mb{Dom}^\Sigma$ to $\outer \uds{\bf Curve}_{\rm rig}^{\mb C}$ as a monotone multimodule over the given lifts of the monotone flow categories. Moreover, for all $d$ there exist ${\mc G}_d$-invariant orientations on $B_d^\Sigma$ satisfying
\begin{enumerate}
    \item For each $i \in \{1, \ldots, m\}$ and $d_i > 0$, the map
    \beqn
    {\mc G}_{d_i + d} \times_{{\mc G}_{d_i} \times {\mc G}_d} (B_{d_i} \times B_d^\Sigma) \to \partial^{d_i, d} B_{d + d_i}^\Sigma
    \eeqn
    is orientation preserving.

    \item For each $d' > 0$, the map
    \beqn
    {\mc G}_{d + d'} \times_{{\mc G}_d \times {\mc G}_{d'}} (B_d^\Sigma \times B_{d'}) \to \partial^{d, d'} B_{d + d'}^\Sigma
    \eeqn
    is orientation reversing. 
\end{enumerate}
\end{prop}

\begin{proof}
The construction of the rigidification can be carried out as done in the proof of Proposition \ref{prop_curve_rigidification}. The construction of stable complex structure can be carried out similarly as in the proof of Proposition \ref{prop_domain_stable_complex}. As for the orientation, the argument is similar to Lemma \ref{lemma_domain_orientation}. However, notice that $B_d$ has an odd dimension and $B_d^\Sigma$ has an even dimension, which make the difference in the two cases as one has to interchange the $[0, \epsilon)$ factor with either $B_{d_i}$ or $B_d^\Sigma$. 
\end{proof}

\subsubsection{Homotopies}

Flow homotopies associated to a 1-parameter family of variations of domain Floer data is indeed quite trivial. Indeed, let $\Sigma_t$, $t \in [0, 1]$ be a smooth family of Floer domains with $m$ negative ends and $1$ positive end. Then there are $m+1$ copies of the monotone flow categories $\mb{Dom}$ corresponding to the $m$ inputs and 1 output. Suppose $\Sigma_t$ is equipped with a smooth family of lateral lines. Then $\Sigma_0$ and $\Sigma_1$ correspond to monotone multimodules $\mb{Dom}^{\Sigma_0}$ and $\mb{Dom}^{\Sigma_1}$ over $(\mb{Dom}, \ldots, \mb{Dom}; \mb{Dom})$. Moreover, by allowing $t$ to vary, one naturally obtains a monotone homotopy $\mb{Dom}^{\mb H}$ from $\mb{Dom}^{\Sigma_0}$ to $\mb{Dom}^{\Sigma_1}$. Upon choosing any outercollaring width, these items become corresponding collared versions. The following proposition demonstrates that the choices of rigidifications and stable complex structures can also be carried inductively. 

\begin{prop}
Let $\outer\mb{Dom}_1, \ldots, \outer \mb{Dom}_m; \outer\mb{Dom}'$ be monotone flow categories enriched in $\outer \uds{\bf Curve}_{\rm rig}^{\mb C}$ obtained from the monotone flow category $\mb{Dom}$ by outercollaring and possibly different rigidifications. Let 
\begin{align*}
&\ \outer \mb{Dom}^{\Sigma_0},\ &\ \outer \mb{Dom}^{\Sigma_1}
\end{align*}
be monotone multimodules over $(\outer \mb{Dom}_1, \ldots, \outer \mb{Dom}_m; \outer \mb{Dom}')$ enriched in $\outer \uds{\bf Curve}_{\rm rig}^{\mb C}$ constructed above. Then the outercollaring $\outer \mb{Dom}^{\mb H}$ admits a lift to $\outer \uds{\bf Curve}_{\rm rig}^{\mb C}$ as a monotone homotopy from $\outer \mb{Dom}^{\Sigma_0}$ to $\outer \mb{Dom}^{\Sigma_1}$.
\end{prop}

\begin{proof}
Routine verification of the inductive constructions. 
\end{proof}

\begin{figure}[h]
    \centering
    \includegraphics[scale = 0.6]{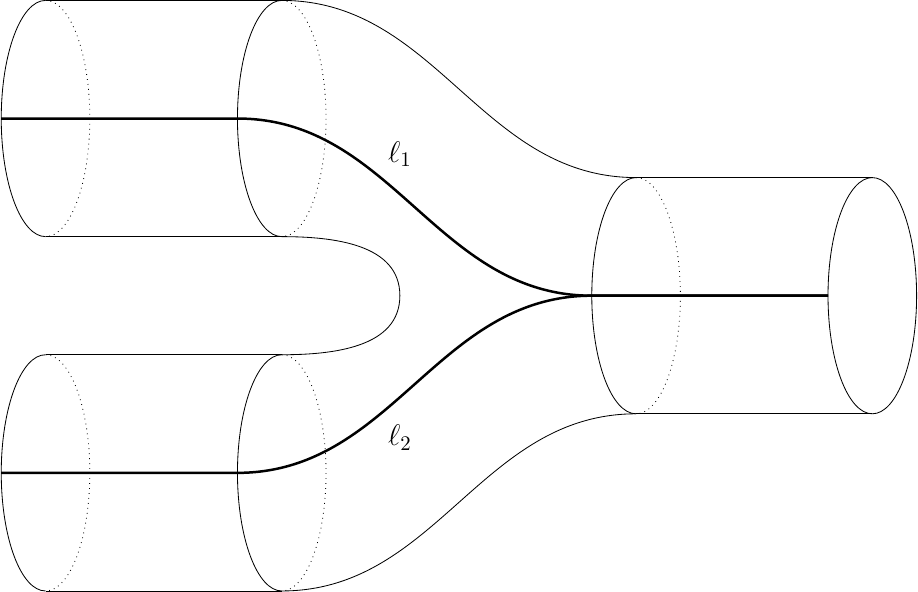}
    \caption{Lateral lines in a pair-of-pants.}
    \label{fig:lateral_lines}
\end{figure}

\subsubsection{The case of cigar bimodule}

Let $\Sigma^{\rm cigar}$ be the complex plane ${\mb C}$ equipped with a cylindrical end. It is used to construct the cigar bimodule and eventually lead to the multiplicative identity of the pair-of-pants product. The construction of the corresponding monotone bimodule, denoted by $\outer \mb{Dom}^{\rm cigar}$, is slightly different from the general case. 

We shall regard $\outer \mb{Dom}^{\rm cigar}$ not as a bimodule, but only a {\bf right module} over $\outer \mb{Dom}$. This means it consists of objects
\beqn
{\mc C}_d^{\rm cigar} \in \outer \uds{\bf Curve}_{\rm rig}^{\mb C}
\eeqn
together with structural maps
\beqn
{\mc C}_d^{\rm cigar}\times {\mc C}_{d'} \to {\mc C}_{d + d'}^{\rm cigar}
\eeqn
which satisfy the obvious associativity relation. The construction of the object ${\mc C}_d^{\rm cigar}$ is the same as before. Indeed, consider the moduli space of degree $d$ maps
\beqn
\phi: \Sigma^{\rm cigar} \to \mb{CP}^d
\eeqn
satisfying the constraint
\beqn
\phi(\infty) = [0, \ldots, 0, 1]\in \mb{P}^d
\eeqn
and such that the image of $\phi$ is not contained in any hyperplane. After compactification in the standard way, one obtains a moduli space $B_d^{\rm cigar}$ equipped with a ${\mc G}_d$-action and a universal curve $C_d^{\rm cigar} \to B_d^{\rm cigar}$, hence provides an object
\beqn
{\mc C}_d^{\rm cigar} = ({\mc G}_d, B_d^{\rm cigar}, C_d^{\rm cigar})
\eeqn
stratified by the poset
\beqn
A_d^{\rm cigar} = \{ (d_0, d_1, \ldots, d_l)\ |\ d_0 + \cdots + d_l = d,\ d_i \geq 1 \}.
\eeqn

The difference from the previous case is the definition of the structural maps. In the current situation, we do not need a lateral line in $\Sigma^{\rm cigar}$. Other constructions are the same. This provides a monotone right module $\mb{Dom}^{\rm cigar}$ over $\mb{Dom}$. Then after outercollaring and choosing a lift of $\outer \mb{Dom}$ to $\outer \uds{\bf Curve}_{\rm rig}^{\mb C}$, one can find a rigidification of the outercollaring $\outer \mb{Dom}^{\rm cigar}$.

\subsubsection{The case of Poincar\'e duality multimodule}

Now let $\Sigma^{\rm PD}$ be the smooth Floer domain with two negative ends and no positive ends, which is responsible for defining the Poincar\'e duality for Floer theory. One can similar define a monotone multimodule $\outer \mb{Dom}^{\rm PD}$ enriched in $\outer \uds{\bf Curve}_{\rm rig}^{\mb C}$ as above, whose construction slightly differ from previous cases. 

In this case, the structure we build is better called a double left module, meaning a collection of objects
\beqn
{\mc C}_d^{\rm PD} \in \outer \uds{\bf Curve}_{\rm rig}^{\mb C}
\eeqn
together with two structural maps
\begin{align*}
&\ {\mc C}_{d_1} \times {\mc C}_{d}^{\rm PD} \to {\mc C}_{d_1 + d}^{\rm PD},\ {\mc C}_{d_2} \times {\mc C}_d^{\rm PD} \to {\mc C}_{d_2 + d}^{\rm PD}
\end{align*}
corresponding to breakings at the two negative ends, which satisfy similar associativity conditions as a general monotone multimodule.

We first consider the moduli space of degree $d$ stable maps
\beqn
\phi: \Sigma^{\rm PD} \to \mb{CP}^d
\eeqn
whose images are not contained in any hyperplane. As the cylindrical ends are on the negative side, in the current convention we do not impose any evaluation constraint. After the standard compactification, one obtains such a moduli space
\beqn
B_d^{\rm PD}
\eeqn
which has a ${\mc G}_d$-action and a universal curve $C_d^{\rm PD} \to B_d^{\rm PD}$, providing an object
\beqn
{\mc C}_d^{\rm PD} = ({\mc G}_d, B_d^{\rm PD}, C_d^{\rm PD}).
\eeqn

The structural maps are also defined similarly, but without using lateral lines. We omit the details. Hence one obtains a monotone double left module $\mb{Dom}^{\rm PD}$ over $(\mb{Dom}, \mb{Dom})$. Then after outercollaring, upon choosing two independent lifts $\outer \mb{Dom}_1, \outer \mb{Dom}_2$ to $\outer \uds{\bf Curve}_{\rm rig}^{\mb{C}}$, one can construct a lift $\outer \mb{Dom}^{\rm PD}$ to $\outer \uds{\bf Curve}_{\rm rig}^{\mb{C}}$. We omit the details.

\section{Relatively Smooth AMS lifts}\label{section_Floer_lift}

In this section we carry out perhaps the most complex part of the AMS construction. We demonstrate how to lift an outercollaring of a Floer flow category ${\mb F}$ to $\outer \uds{\bf S^{\rm rel}Kur}_{\rm rig}$, the category of (collared and rigidified) relatively smooth Kuranishi spaces. In fact, we will construct a lift to $\outer \uds{\bf Kur}_{\rm rig}$, the category of collared and rigidified topological Kuranishi spaces. The lift to relatively smooth version is automatic. The collared-and-rigidified feature is always preserved in the construction as it is a basic prerequisite for inductively construct perturbations.

\subsection{Integral modification of topological energy}

We generalize the method of Rezchikov \cite{Rezchikov_Arnold} of modifying the topological energy. First, there is a way to modify the symplectic form.

\begin{defn}\label{defn_integral_modification}
Let ${\mb F}$ be the Floer flow category associated to a nondegenerate Hamiltonian $H$ on $(X, \omega, J)$. An {\bf integral action} on ${\mb F}$ is a  function
\beqn
{\mc A}_{\mb F}^\Omega: {\rm Ob}{\mb F} \to {\mb Z}
\eeqn
satisfying the following conditions.
\begin{enumerate}

\item There exists an integral symplectic form $\Omega \in \Omega^2(X)$ tamed by $J$ such that ${\rm ker}\Omega = {\rm ker} \omega \subset \pi_2(X)$ and for each $A \in \Pi$ and $p \in {\rm Ob}{\mb F}$
\beqn
{\mc A}_{\mb F}^\Omega(A\cdot p) = {\mc A}^\Omega_{\mb F} (p) - \Omega(A).
\eeqn

\item There exists $C>0$ such that whenever $p \neq q$ and $M_{pq}^{\mb F} \neq \emptyset$, one has
\beqn
0 <  C \left(  {\mc A}_H(p) - {\mc A}_H(q) \right) \leq {\mc A}^\Omega_{\mb F}(p) - {\mc A}^\Omega_{\mb F}(q).
\eeqn
\end{enumerate}
\end{defn}

We need similar notions for multimodules and homotopies.

\begin{defn}\label{module_integral_action}
Consider a smooth genus zero Riemann surface with cylindrical ends $\Sigma$ together with a family of Floer data $\sigma_w$ parametrized continuously by $w$ in a compact topological space such that for all $w$, the Floer data is asymptotic to $(H_i dt, J_i)$ on the $i$-th negative end ($i = 1, \ldots, m$) and asymptotic to $(H_j' dt, J_j')$ on the $j$-th positive end ($j=1, \ldots, n$). A compatible collection of integral actions consists of integral actions 
\begin{align*}
&\ {\mc A}_{{\mb F}_i}^\Omega: {\rm Ob} {\mb F}_i \to {\mb Z}\  (i=1, \ldots, m),\ &\ {\mc A}_{{\mb F}_j'}^\Omega: {\rm Ob}{\mb F}_j' \to {\mb Z}
\end{align*}
on the Floer flow categories ${\mb F}_i$ resp. ${\mb F}_j'$ associated to $(H_i, J_i)$ resp. $(H_j', J_j')$ satisfying the following conditions.
\begin{enumerate}

\item The integral actions induce the same integral symplectic form $\Omega$.

\item For each parameter $w$ and each point $z \in \Sigma$, $\Omega$ is tamed by the almost complex structure $J_{w, z}$.

\item For objects $p_i \in {\rm Ob}{\mb F}_i$ and $p_j' \in {\rm Ob}{\mb F}_j'$, let $M_{p_1\cdots p_m; p_1'\cdots p_n'}$ denote temporarily the moduli space of stable solutions of Floer equation for all parameter $s$ asymptotic to the labeling objects. Then
\beqn
M_{p_1 \cdots p_m; p_1'\cdots p_n'} \neq \emptyset \Longrightarrow \sum_{i=1}^m {\mc A}_{{\mb F}_i}^\Omega(p_i) > \sum_{j=1}^n {\mc A}_{{\mb F}_j'}^\Omega(p_j')
\eeqn
\end{enumerate}
\end{defn}

One can then specialize to multimodules and homotopies. 

\begin{defn}
\begin{enumerate}

\item In Situation \ref{situationm1}, integral actions ${\mc A}_{{\mb F}_1}^\Omega, \ldots, {\mc A}_{{\mb F}_m}^\Omega; {\mc A}_{{\mb F}'}^\Omega$ on the involved flow categories are called {\bf compatible} with respect to ${\mb X}$ if it is compatible with respect to the fixed Floer data $\sigma^{\mb X}$ on $\Sigma^{\mb X}$.

\item In Situation \ref{situationh1}, integral actions ${\mc A}_{{\mb F}_1}^\Omega, \ldots, {\mc A}_{{\mb F}_m}^\Omega, {\mc A}_{{\mb F}'}^\Omega$ are called {\bf compatible} with respect to ${\mb H}$ if it is compatible with respect to the 1-parameter family of Floer data. 

\end{enumerate}
\end{defn}

\begin{prop}\label{prop163} In the setting of Definition \ref{module_integral_action}, there exist compatible integral actions on involved Floer flow categories.
\end{prop}

\begin{proof}
It was proved in \cite[Lemma 16]{Rezchikov_Arnold} that given any collection $H_1, \ldots, H_m$ of nondegenerate 1-periodic Hamiltonians on $(X, \omega)$, there exists an open subset $U\subset X$ containing the images of all 1-periodic orbits of all $H_i$ and for any $\epsilon>0$, a rational symplectic form $\Omega \in \Omega^2(M)$ satisfying the following conditions.
\begin{enumerate}
    \item $\Omega|_U = \omega|_U$.

    \item $\| \Omega - \omega\|_{C^0(X)} \leq \epsilon$.

    \item ${\rm ker} \Omega = {\rm ker} \omega \subset \pi_2(X)$.
\end{enumerate} 

Now consider a Floer domain $\Sigma$ with the Hamiltonian connection $\sigma$.  Let $\Sigma_0 \subset \Sigma$ be a compact region such that $\Sigma \setminus \Sigma_0$ is a union of half cylinders. Let $(s, t)$ be the cylindrical coordinates and let $H_i dt$ be the restriction of $\sigma$ to the $i$-th cylindrical end $U_i$, where $H_i$ is nondegenerate. 

Let $U \subset X$ be the open subset provided by Rezchikov's lemma. Let $\Omega$ be a symplectic form which agrees with $\omega$ within $U$ and which is $C^0$-close to $\omega$ so $\Omega$ tames $J_{w, z}$ for all almost complex structures appearing.  Then there exists $\epsilon_0>0$ satisfying 
\beq\label{eqn161}
\forall i\ \ \forall  \rho: S^1 \to X,\ \int_{S^1} |\rho'(t) - X_{H_{i, t}}(\rho(t))|^2 dt \leq \epsilon_0 \Longrightarrow {\rm Im}(\rho) \subset U.
\eeq
For each capped 1-periodic orbit $p$ of $H_i$ where the capping is a map $u: {\mb D}^2 \to X$, define
\beqn
{\mc A}_{H_i}^\Omega(\gamma, u) = -\int_{{\mb D}^2} u^* \Omega - \int_{S^1} H_i(\gamma(t)) dt.
\eeqn

Let $|\cdot|_\Omega$ be the norm on the tangent bundle $TX$ determined by $\Omega$ and $J_{w, z}$. For any map $u: \Sigma \to X$, define
\beqn
\nabla^{\sigma, \Omega} (u) = \nabla u + X_\sigma^\Omega(u)
\eeqn
where $X_\sigma^\Omega$ is the 1-form on $\Sigma$ taking values in Hamiltonian vector fields associated to $\Omega$. In any local coordinate $z = s + {\bf i} t$ one can write $\sigma = \sigma_s ds + \sigma_t dt$ and
\beqn
\nabla^{\sigma, \Omega} u = \nabla_s^{\sigma, \Omega} u ds + \nabla_t^{\sigma, \Omega} u dt.
\eeqn
The difference of the covariant derivatives between $\Omega$ and $\omega$ is
\beqn
\nabla^{\sigma, \Omega} u - \nabla^{\sigma, \omega} u = (X_{\sigma_s}^\Omega(u) - X_{\sigma_s}^\omega(u)) ds + (X_{\sigma_t}^\Omega(u) - X_{\sigma_t}^\omega(u)) dt =: X_{\sigma_s}^{\Omega- \omega}(u) ds + X_{\sigma_t}^{\Omega - \omega} (u) dt
\eeqn
which is linear in the difference $\Omega - \omega$. 

Define the energy density of the map $u$ as
\beqn
\frac{1}{2} | \nabla^{\sigma, \omega} (u)|_\Omega^2 {\rm dvol}_\Sigma.
\eeqn
Then one has the calculation (see \cite[Lemma 2.2.1]{McDuff_Salamon_2004})
\beqn
\frac{1}{2} | \nabla^{\sigma, \omega} (u)|_\Omega^2 {\rm dvol}_\Sigma = | (\nabla^{\sigma, \omega} (u))^{0,1}|^2 {\rm dvol}_\Sigma + \frac{1}{2} \Big( \Omega(\nabla_s^{\sigma, \omega} u, \nabla_t^{\sigma, \omega} u) + \Omega( J \nabla_s^{\sigma, \omega} u, J \nabla_t^{\sigma, \omega}u)  \Big) ds dt.
\eeqn
We observe that
\beq\label{eqn:trivial}
\begin{split}
&\ |\Omega( J \nabla_s^{\sigma, \omega} u, J \nabla_t^{\sigma, \omega}u) - \Omega(\nabla_s^{\sigma, \omega} u, \nabla_t^{\sigma, \omega} u)| \\
= &\ |\Omega( J \nabla_s^{\sigma, \omega} u, J \nabla_t^{\sigma, \omega}u) - \omega( J \nabla_s^{\sigma, \omega} u, J \nabla_t^{\sigma, \omega}u) + \omega( J \nabla_s^{\sigma, \omega} u, J \nabla_t^{\sigma, \omega}u) - \Omega(\nabla_s^{\sigma, \omega} u, \nabla_t^{\sigma, \omega} u)| \\
= &\ |\Omega( J \nabla_s^{\sigma, \omega} u, J \nabla_t^{\sigma, \omega}u) - \omega( J \nabla_s^{\sigma, \omega} u, J \nabla_t^{\sigma, \omega}u) + \omega( \nabla_s^{\sigma, \omega} u, \nabla_t^{\sigma, \omega}u) - \Omega(\nabla_s^{\sigma, \omega} u, \nabla_t^{\sigma, \omega} u)| \\
\leq &\ \| \omega - \Omega \|_{C^0(X)} \cdot | \nabla^{\sigma, \omega} (u)|_\Omega^2
\end{split}
\eeq

Now suppose $u$ solves $(\nabla^{\sigma, \omega} u)^{0,1} = 0$. Then one has 
\beqn
\begin{split}
&\ \Omega(\nabla_s^{\sigma, \omega} u, \nabla_t^{\sigma, \omega} u)\\
= &\  \Omega \Big( \nabla_s^{\sigma, \Omega} u + X_{\sigma_s}^{\omega-\Omega}(u), \nabla_t^{\sigma, \Omega} u + X_{\sigma_t}^{\omega-\Omega}(u) \Big)\\
 = &\ \Omega( \nabla_s^{\sigma, \Omega} u , \nabla_t^{\sigma, \Omega} u) + \Omega( \nabla_s^{\sigma, \Omega} u, X_{\sigma_t}^{\omega - \Omega} (u)) + \Omega( X_{\sigma_s}^{\omega - \Omega}(u), \nabla_t^{\sigma, \Omega} u ) + \Omega( X_{\sigma_s}^{\omega - \Omega}(u), X_{\sigma_t}^{\omega - \Omega}(u))\\
 = &\ \Omega( \nabla_s^{\sigma, \Omega}u, \nabla_t^{\sigma, \Omega}u) + \Omega( \nabla_s^{\sigma, \omega}u, X_{\sigma_t}^{\omega-\Omega}(u)) + \Omega( X_{\sigma_s}^{\omega-\Omega}(u), \nabla_t^{\sigma, \omega}(u)) + \Omega( X_{\sigma_s}^{\omega-\Omega}(u), X_{\sigma_t}^{\Omega- \omega}(u)) 
 \end{split}
\eeqn
Notice that if $u$ converges to capped 1-periodic orbits $p_i$ at the $i$-th cylindrical end, then 
\beqn
\int_\Sigma \Omega(\nabla_s^{\sigma, \Omega} u, \nabla_t^{\sigma, \Omega} u) ds dt =  \sum_{i=1}^m \pm {\mc A}_{H_i}^\Omega(p_i) + \int_\Sigma F_\sigma(u)
\eeqn
where the first term only depends on the homotopy class of $u$ and $F_\sigma$ does not depend on the symplectic form. Then one has
\begin{multline*}
\sum_{i=1}^m \pm {\mc A}_{H_i}^\Omega(p_i) = \int_\Sigma \Omega(\nabla_s^{\sigma, \Omega} u, \nabla_t^{\sigma, \Omega} u) ds dt-  \int_\Sigma F_\sigma(u) \\
- \int_\Sigma \Big(  \Omega( \nabla_s^{\sigma, \omega} u, X_{\sigma_t}^{\omega-\Omega}(u) + \Omega( X_{\sigma_s}^{\omega-\Omega} (u), \nabla_t^{\sigma, \omega}u) - \Omega( X_{\sigma_s}^{\omega-\Omega}(u), X_{\sigma_t}^{\omega-\Omega}(u)) \Big) dsdt.
\end{multline*}
We would like to estimate the last term. First, on the compact region $\Sigma_0$, one has 
\begin{multline*}
\left| \int_{\Sigma_0} \Big(  \Omega( \nabla_s^{\sigma, \omega} u, X_{\sigma_t}^{\omega-\Omega}(u) + \Omega( X_{\sigma_s}^{\omega-\Omega} (u), \nabla_t^{\sigma, \omega}u) - \Omega( X_{\sigma_s}^{\omega-\Omega}(u), X_{\sigma_t}^{\omega-\Omega}(u)) \Big) dsdt \right| \\
\leq C (\epsilon) \| \omega - \Omega\|_{C^0(X)}^2 + \epsilon \int_{\Sigma_0} | \nabla^{\sigma, \omega} (u) |_\Omega^2 {\rm dvol}_\Sigma
\end{multline*}
for any $\epsilon>0$ and a suitable $C(\epsilon)>0$. On the other hand, on each half cylinder $U_i\cong I_i \times S^1$ where $I_i = {\mb R}_\pm$, contained in $\Sigma \setminus \Sigma_0$, we may take $(s, t)$ to be the cylindrical coordinate. Then we know that $\sigma_s = 0$. Hence we are left with 
\beq\label{eqn153}
\int_{U_i} \Omega(\nabla_s^{\sigma, \omega}u, X_{H_i}^{\omega - \Omega}(u)) ds dt.
\eeq
Recall the pre-chosen constant $\epsilon_0$ (see \eqref{eqn161}). Consider 
\beqn
I_i(\epsilon_0):= \Big\{ s \in I_i \ |\ \int_{S^1} |\nabla_s^{\sigma, \omega}(u)(s, t)|_\omega^2 dt \geq \epsilon_0 \Big\}.
\eeqn
Then $|I_i(\epsilon)| \leq E^\omega(u) / \epsilon_0$. Since $\Omega = \omega$ near 1-periodic orbits, by our assumption, the integrand of \eqref{eqn153} vanishes outside $I_i(\epsilon_0)\times S^1$. Then by Cauchy--Schwarz inequality, one has 
\begin{multline*}
\left| \int_{U_i} \Omega( \nabla_s^{\sigma, \omega}u, X_{\sigma_t}^{\omega-\Omega}(u)) ds dt \right| = \left| \int_{I_i(\epsilon_0) \times S^1} \Omega(\nabla_s^{\sigma, \omega} u, X_{\sigma_t}^{\omega-\Omega}(u)) ds dt \right|  \\
\leq C \|\omega-\Omega\|_{C^0(X)} \sqrt{ |I_i(\epsilon_0)|} \sqrt{ \int_{I_i(\epsilon_0)\times S^1} |\nabla_s^{\sigma, \omega} u |_\omega^2 ds dt } \\
\leq C \| \omega- \Omega\|_{C^0(X)} \frac{ E^\omega(u)}{\sqrt{\epsilon_0}}. 
\end{multline*}
On the other hand, when $\Omega$ is close to $\omega$, the norms $|\cdot|_\omega$ and $|\cdot|_\Omega$ are compatible. Therefore, when $\Omega$ is $C^0$-close to $\omega$, one has the estimate
\beqn
\sum_{i=1}^m \pm {\mc A}_{H_i}^\Omega(p_i) \geq  \frac{1}{2} E^\Omega(u) - C
\eeqn
for some constant $C>0$ by combining the above computations with \eqref{eqn:trivial} to rearrange the inequalities, provided that $\epsilon$ is sufficiently small.

Further, to make the symplectic action rational, we define
\beqn
{\mc A}_{{\mb F}_i}^\Omega(p_i):= {\mc A}_{H_i}^\Omega(p_i) + \delta_{p_i}
\eeqn
for a small real number $\delta_{p_i}$ making the above rational. Then one still has 
\beqn
\sum_{i=1}^m \pm {\mc A}_{\mb{F}_i}^\Omega(p_i) \geq \frac{1}{2} E^\omega (u) - C
\eeqn
after changing $C$ a bit. Then we argue in two cases.

\begin{enumerate}

\item When $\Sigma$ is an infinite cylinder, $\sigma = H dt$ and $J$ only depends on $t \in S^1$, one can take $C = 0$ and obtain instead
\beqn
{\mc A}_{{\mb F}}^\Omega(p) - {\mc A}_{{\mb F}}^\Omega(q) \geq \frac{1}{2} E^\omega (u) = \frac{1}{2} ({\mc A}_H^\omega(p) - {\mc A}_H^\omega(q)). 
\eeqn
This provides an integral action on the associated Floer flow category.

\item In other cases, one can add rational constants to ${\mc A}_{H_i}^\Omega$ to make the sum of $\pm {\mc A}_{H_i}^\Omega$ always positive.
\end{enumerate}
Finally, we choose a common large positive integral multiple of the ${\mc A}^{\Omega}_{\mb F}$ constructed above so that the action becomes integral valued and the symplectic form $\Omega$ becomes integral.
\end{proof}

\subsubsection{Pullback the domain flow category, multimodule, and homotopy}\label{subsubsec:pullback-curve}

Now choose an integral action  ${\mc A}_{\mb F}^\Omega:  {\rm Ob}{\mb F} \to {\mb Z}$ with associated symplectic form $\Omega$. For each pair of objects $p, q \in {\rm Ob}{\mb F}$, abbreviate $d = d_{pq} = {\mc A}_{\mb F}^\Omega(p) - {\mc A}_{\mb F}^\Omega (q)$. Inside the poset $A_{d}$ (of positive partitions of $d$), there is an Alexandrov open subset 
\beqn
A_{pq}^\Omega \subset A_d
\eeqn
whose elements are decompositions 
\beqn
d = d_0 + \cdots + d_l
\eeqn
such that there exist objects $r_1, \ldots, r_l\in {\rm Ob} {\mb F}$ such that 
\beqn
d_0 = {\mc A}_{\mb F}^\Omega (p) - {\mc A}_{\mb F}^\Omega ({r_1}),\ d_1 = {\mc A}_{\mb F}^\Omega (r_1) - {\mc A}_{\mb F}^\Omega (r_2),\ \ldots, d_l =  {\mc A}_{\mb F}^\Omega (r_l) -{\mc A}_{\mb F}^\Omega (q).
\eeqn

Now we use these data to ``pull back'' the monotone flow categories $\outer \mb{Dom}$ constructed before. Define a flow category enriched in $\outer \uds{\bf Curve}_{\rm rig}^{\mb C}$ whose objects are the same as ${\mb F}$, and whose morphism from $p$ to $q$ is the restriction of $\outer {\mc C}_{d}$ to the corresponding open subset
\beqn
\outer B_{pq}^\Omega \subset \outer B_d.
\eeqn
Denote this object by 
\beqn
{\mc C}_{pq}^\Omega = ({\mc G}_{pq}^\Omega, B_{pq}^\Omega, C_{pq}^\Omega).
\eeqn
When $\Omega$ is understood from the context, we often drop it from the notation. Then the restriction of the structural maps of the flow category $\outer \mb{Dom}$ gives a flow category
\beqn
\outer \mb{Dom}^{\Omega, {\mb F}}
\eeqn
enriched in $\outer \uds{\bf Curve}_{\rm rig}^{\mb C}$ whose objects are the same as objects of ${\mb F}$ and whose morphisms are the families $\outer {\mc C}_{pq}^\Omega$.

Now we summarize the construction as an update of Situation \ref{situationf1}.

\begin{situationf}\label{situationf2}
In Siutation \ref{situationf1}, one fixes an outercollaring width. $\outer \mb{Dom}$ is the monotone flow category enriched in $\outer \uds{\bf Curve}_{\rm rig}^{\mb C}$. Furthermore, let $\outer \mb{Dom}^{\Omega, {\mb F}}$ be the flow category enriched in $\outer \uds{\bf Curve}_{\rm rig}^{\mb C}$ subject to a chosen integral action on ${\mb F}$.
\end{situationf}

Similarly, let ${\mb X}$ be a multimodule associated to a Floer domain $\Sigma$. Choose an outercollaring width and a compatible set of integral actions on involved flow categories 
\begin{align*}
&\ {\mc A}_{{\mb F}_j}^\Omega: {\rm Ob}{\mb F}_j \to {\mb Z},\ &\ {\mc A}_{{\mb F}'}^\Omega: {\rm Ob}{\mb F}' \to {\mb Z}.
\end{align*}
Then one can choose an integer $n^{\mb X}$ such that for all tuples of objects $(p_1, \ldots, p_m; p')$, one has
\beq\label{eqn:constant-shift}
d:= d_{p_1 \cdots p_m; p'}:= \sum_{i=1}^m {\mc A}^\Omega_{{\mb F}_i}(p_i) - {\mc A}^\Omega_{{\mb F}'}(p') + n^{\mb X} > 0.
\eeq
Let 
\beqn
A_{p_1 \cdots p_m; p'}^{\Omega, {\mb X}} \subset 
A_d
\eeqn
be the Alexandrov open subset of partitions of $d$ corresponding to possible degenerations of $\Sigma$ and let 
\beqn
{\mc C}_{p_1\cdots p_m; p'}^{\mb X}
\eeqn
be the restriction to the open subset corresponding to $A_{p_1 \cdots p_m; p'}^{\Omega, \mb{X}}$. Then this provides a multimodule over $(\outer \mb{Dom}^{\Omega, {\mb F}_1}, \ldots, \outer \mb{Dom}^{\Omega, {\mb F}_m}; \outer \mb{Dom}^{\Omega, {\mb F}'} )$ enriched in $\outer \uds{\bf Curve}_{\rm rig}^{\mb C}$, denoted by 
\beqn
\outer \mb{Dom}^{\Omega, {\mb X}}.
\eeqn

\begin{rem}
    The constant $n^{\mb X}$ should be thought as shifting the integral action ${\mc A}^{\Omega}_{\mb F'}$ by $-n^{\mb X}$ such that for all $(p_1, \ldots, p_m; p')$, we have
\beqn
d:= d_{p_1 \cdots p_m; p'}:= \sum_{i=1}^m {\mc A}^\Omega_{{\mb F}_i}(p_i) - ({\mc A}^\Omega_{{\mb F}'}(p') - n^{\mb X}) > 0.
\eeqn
The positive integer $d$ will be the degree of an auxiliary line bundle over ${\mb P}^1$.
\end{rem}

\begin{situationm}\label{situationm2}
In Situation \ref{situationm1}, one fixes an outercollaring width and a compatible list of integral actions on the involved Floer categories. $\outer \mb{Dom}^{\Omega, {\mb F}_1}, \ldots, \outer \mb{Dom}^{\Omega, {\mb F}_m}; \outer \mb{Dom}^{\Omega, {\mb F}'}$ are flow categories enriched in $\outer \uds{\bf Curve}_{\rm rig}^{\mb C}$ subject to the chosen integral actions as described in Situation \ref{situationf2}. $\outer \mb{Dom}^{\Omega, {\mb X}}$ is the multimodule over $(\outer \mb{Dom}^{\Omega, {\mb F}_1}, \ldots, \outer \mb{Dom}^{\Omega, {\mb F}_m}; \outer \mb{Dom}^{\Omega, {\mb F}'})$.
\end{situationm}

The case of Situation \ref{situationh1} is similar and we only describe the construction below without details.

\begin{situationh}\label{situationh2}
In Situation \ref{situationh1}, one fixes an outercollaring width and a compatible list of integral actions on involved Floer flow categories. 
\beqn
\outer \mb{Dom}^{\Omega, {\mb F}_1},\ldots, \outer \mb{Dom}^{\Omega, {\mb F}_m}; \outer \mb{Dom}^{\Omega, {\mb F}'}
\eeqn
are flow categories enriched in $\outer \uds{\bf Curve}_{\rm rig}^{\mb C}$ subject to the chosen integral actions as described in Situation \ref{situationf2};
\beqn
\outer \mb{Dom}^{\Omega, {\mb X}_0}, \outer \mb{Dom}^{\Omega, {\mb X}_1}
\eeqn
are multimodules over them associated to $\Sigma_0$ and $\Sigma_1$ respectively as described in Situation \ref{situationm2}. $\outer \mb{Dom}^{\Omega, {\mb H}}$ is a homotopy from $\outer \mb{Dom}^{\Omega, {\mb X}_0}$ to $\outer \mb{Dom}^{\Omega, {\mb X}_1}$ subject to the fixed outercollaring width and the integral actions.
\end{situationh}

\subsection{Framed maps}\label{subsection152}

In this subsection, we describe a suitable package describing flow categories, multimodules, and homotopies of infinite-dimensional spaces of maps into the target symplectic manifold.

\subsubsection{Category of families of framed maps}

The following definitions are rather artificial. However, they help us organize various constructions and package required properties. We use the following notations a lot: given any group $G$, let $\hat G = G \times G$. We often denote the first factor by $G_L$ and the second factor by $G_R$.

\begin{defn}\label{defn156}
\begin{enumerate}

\item The category $\uds{\bf Map}^{\rm fr}$ has objects
\beqn
\hat {\mc C}:= ({\mc G}, B, C, Y)
\eeqn
where $({\mc G}, B, C)$ is an object of $\uds{\bf Curve}$, $Y$ is a ${\mc G}\times {\mc G}$-space together with a map $\pi: Y  \to B$ which is equivariant with respect to the projection ${\mc G}_L \times {\mc G}_R \to {\mc G}$ onto the first factor and the action by ${\mc G}_R$ on $Y$ is free.  

\item Let $\hat {\mc C}_i = ({\mc G}_i, B_i, C_i, Y_i)$, $i = 1, 2$ be two objects of $\uds{\bf Map}^{\rm fr}$. A strict embedding from $\hat{\mc C}_1$ to $\hat {\mc C}_2$, denoted by $\hat \zeta_{21}$, consists of a strict embedding
\beqn
\zeta_{21}: {\mc C}_1 \to {\mc C}_2
\eeqn
in $\uds{\bf Curve}$ and a map $\zeta_{21}^Y: Y_1 \to Y_2$ equivariant with respect to $\hat {\mc G}_1 \to \hat{\mc G}_2$ such that 
\begin{enumerate}
    \item The induced map
    \beqn
    \hat {\mc G}_2 \times_{\hat {\mc G}_1} Y_1 \to Y_2
    \eeqn
    is an homeomorphism onto the corresponding union of strata. 
    
    \item The following diagram commutes.
\beqn
\xymatrix{ B_1  \ar[d]_{\zeta_{21}^B}   &  Y_1  \ar[d]^{\zeta_{21}^Y} \ar[l]    \\
B_2 & Y_2 \ar[l]    }
\eeqn
\end{enumerate}

\item A {\bf rigidification} of a strict embedding $\hat\zeta_{21}$ is a rigidification of the underlying strict embedding $\zeta_{21}$ in $\uds{\bf Curve}$.

\item The notion of unitary conjugations and embeddings are defined similar to many previous cases. We omit the details. Let $\uds{\bf Map}^{\rm fr}$ be the category of framed maps with morphisms being embeddings and $\uds{\bf Map}^{\rm fr}_{\rm rig}$ be the rigidified version.
\end{enumerate}
\end{defn}

\subsubsection{Group reductions} 

By definition, Kuranishi spaces have symmetry group being a compact Lie group $G$, while the corresponding AMS domain has the symmetry group being a complex Lie group ${\mc G}$. The redundancy is reduced by choosing certain group reductions. 

\begin{defn}[The category of families of reduced framed maps]\label{defn_group_reduction} \hfill
\begin{enumerate}

    \item Let $\hat {\mc C} = ({\mc G}, B, C, Y)$ be an object of $\uds{\bf Map}^{\rm fr}$. A {\bf group reduction} of $\hat {\mc C}$ is a map
\beqn
\Lambda:  Y \to {\mc G}/ G
\eeqn
which is ${\mc G}_L$-invariant and ${\mc G}_R$-equivariant. 

\item The category of {\bf equivariant families of reduced framed maps}, denoted by $\uds{\bf map}^{\rm fr}$ consists of objects $({\mc G}, B, C, Y, \Lambda)$ where $({\mc G}, B, C, Y)$ is an object of $\uds{\bf Map}^{\rm fr}$ and $\Lambda$ is a group reduction. A {\bf strict embedding} $\hat\zeta_{21}: \hat {\mc C}_1 \to \hat{\mc C}_2$ in $\uds{\bf map}^{\rm fr}$ consists of a strict embedding of $\uds {\bf Map}^{\rm fr}$, which induces a map
\beqn
\zeta_{21}^P: {\mc G}_1/G_1 \to {\mc G}_2/G_2
\eeqn
such that following diagram commutes.
\beqn
\xymatrix{  Y_1 \ar[r]^-{\Lambda_1} \ar[d]_{\zeta_{21}^Y}  &   ({\mc G}_1/G_1) \ar[d]^{\zeta_{21}^P} \\
            Y_2 \ar[r]_-{\Lambda_2}     & ({\mc G}_2/G_2)  }
\eeqn
The notions of unitary conjugations and embeddings are defined in a way similar to previously considered situations. Morphisms in $\uds{\bf map}^{\rm fr}$ are embeddings, i.e., strict embeddings modulo unitary conjugations.

\item A rigidification of a strict embedding $\hat\zeta_{21}$ as considered above is a rigidification of the underlying strict embedding in $\uds{\bf Curve}$, and we denote the resulting category by $\uds{\bf map}^{\rm fr}_{\rm rig}$.

\end{enumerate}
\end{defn}

As usual, we indicate the outercollaring using the notation $\outer$.

\subsubsection{Modification of energy density}\label{subsubsec:energy-modify}

Given an integral action ${\mc A}_{\mb{F}}^\Omega$ on ${\mb F}$, there are also the induced modification of the energy density. Given any smooth map $u: {\mb R}\times S^1 \to X$ which converges to a pair of objects $p$ and $q$ with ${\mc A}_H(p) > {\mc A}_H(q)$, define 
\beq\label{integral_energy_density}
\Omega_u:= \Omega( \partial_s u, \partial_t u - X^{\Omega}_{H_t}(u)) ds \wedge dt \in \Omega^2({\mb R}\times S^1).
\eeq
Then its integral is the integer ${\mc A}_{\mb F}^\Omega(p) - {\mc A}_{\mb F}^\Omega(q)$. We may regard $\Omega_u$ as a modified energy density, which is everywhere nonnegative when $u$ is a Floer cylinder. Similarly, given a smooth map $u: S^2 \to X$ representing a sphere class $A \in \pi_2(X)$, define
\beqn
\Omega_u:= u^* \Omega \in \Omega^2(S^2).
\eeqn
Its integral is $\Omega(A)$; when $u$ is a $J$-holomorphic sphere, $\Omega_u$ is everywhere nonnegative.

In general, given smooth map $u: \Sigma \to X$ over a Floer domain with associated flow multimodudle ${\mb X}$ over $({\mb F}_1, \dots, {\mb F}_m, {\mb F}')$, recall that we have chosen an integer $n^{\mb X} > 0$. In this setting, we consider the $2$-form
\beqn
\Omega_u := u^* \Omega + F^{\Omega}_{\alpha}(u) + n^{\mb X} d{\rm vol}_{\Sigma},
\eeqn
where $\alpha$ is the Hamiltonian connection from the Floer datum, and $F^{\Omega}_{\alpha}(u)$ denotes its curvature with respect to the symplectic form $\Omega$. If $u$ solves the Floer equation, $\Omega_u$ is everywhere nonnegative, and $\int_{\Sigma} \Omega_u$ recovers $d_{p_1 \cdots p_m;p'}$ from \eqref{eqn:constant-shift}. The same construction applies to maps from multimodule homotopies.

\subsubsection{Floer admissible maps
and framed maps}

Now we provide the concrete flow categories, multimodules, and homotopies enriched in $\uds{\bf Map}^{\rm fr}$. 

We first consider the flow category case. Let ${\mb F}$ be a Floer flow category associated to $(X, \omega, J)$ and the Hamiltonian $H$. Let ${\mc A}_{\mb F}^\Omega: {\rm Ob}{\mb F} \to {\mb Z}$ be an integral action with associated integral symplectic form $\Omega$. Recall that one has obtained a flow category $\outer \mb{Dom}^{\Omega, {\mb F}}$ with the same set of objects as ${\mb F}$ enriched in $\outer \uds{\bf Curve}_{\rm rig}$. For each $p, q \in {\rm Ob}{\mb F}$, one has the object
\beqn
{\mc C}_{pq}^{\mb F} = ({\mc G}_{pq}^{\mb F}, B_{pq}^{\mb F}, C_{pq}^{\mb F}) \in {\rm Ob} \uds{\bf Curve}_{\rm rig}.
\eeqn
We define a $\hat{\mc G}_{pq}^{\mb F}$-space $Y_{pq}^{\mb F}$ as follows. For each $\phi \in B_{pq}^{\mb F}$, consider continuous maps $u: C_\phi \to X$ which are smooth away from nodes and which converge exponentially fast towards 1-periodic orbits of $H$. We require that $u$ converges to the orbits $p$ resp. $q$ at the negative resp. positive end of $C_\phi$, and that the cappings satisfy $p \# u = q$. In particular, the homotopy type of the map $u$ is constrained by the underlying curve; for example, if $C_\phi$ has a sphere component of degree $d'$, then the restriction of $u$ to that component has the same degree $d'$ measured by the modified integral symplectic form $\Omega$. For $\phi$ and $u$ as above, there is a Hermitian holomorphic line bundle
\beqn
L_u  \to C_\phi
\eeqn
determined by the 2-form $\Omega_u \in \Omega^2(C_\phi)$ (see \eqref{integral_energy_density}) which decays exponentially along cylindrical end. The degree of $L_u$ is $d_{pq}$. Then consider frames 
\beqn
F = (f_0, \cdots, f_{d_{pq}})
\eeqn
of the space $H^0(L_u)$ modulo ${\mb C}^*$ such that 1) at the negative end $z_-$ one has 
\beqn
[f_0(z_-), \ldots, f_{d_{pq}}(z_-)] = [1, 0, \ldots, 0] \in \mb{CP}^{d_{pq}}
\eeqn
and 2) the Hermitian matrix
\beqn
\left( \int_{C_\phi} \langle f_\mu, f_\nu \rangle \Omega_u \right)_{\mu\nu}
\eeqn
is positive definite, where $\langle \cdot , \cdot \rangle$ is the fiberwise Hermitian metric on the bundle $L_u$. Notice that the set of  such frames $F$ is a ${\mc G}_{pq}$-torsor. Then let $Y_{pq}^{\mb F}$ be the set of all tuples $(\phi, u, F)$ described as above. 

The set $Y_{pq}^{\mb F}$ has a natural topology and admits a natural forgetful map $\pi_0: Y_{pq}^{\mb F} \to B_{pq}^{\mb F}$. There is another map
\beqn
\pi_1: Y_{pq}^{\mb F} \to B_{pq}^{\mb F},\ (\phi, u, F) \mapsto \phi_F 
\eeqn
where $\phi_F: C_\phi \to \mb{CP}^{{\bf d}_{pq}}$ is the Kodaira embedding induced by the frame. Moreover, $Y_{pq}^{\mb F}$ has an action by $\hat {\mc G}_{pq} \cong {\mc G}_{pq} \times {\mc G}_{pq}$, where the first factor acts by reparametrizing the domain and the second factor acts by linearly transforming the frame $F$. One can see that the tuple
\beqn
({\mc G}_{pq}^{\mb F}, B_{pq}^{\mb F}, C_{pq}^{\mb F}, Y_{pq}^{\mb F})
\eeqn
together with the maps $\pi_0, \pi_1$ is an object of $\uds{\bf Map}^{\rm fr}$. 

The cases of multimodules and homotopies are similar, using the modified energy from Section \ref{subsubsec:energy-modify}. One obtains on the level of spaces
\begin{align*}
&\ Y_{p_1 \cdots p_m; p'}^{{\mb X}},\ &\ 
Y_{p_1 \cdots p_m; p'}^{{\mb H}}
\end{align*}
as well as the corresponding objects in $\uds{\bf Map}^{\rm fr}$.

\subsubsection{The structural maps}

We would like to make the spaces of framed maps into flow categories, multimodules, and homotopies enriched in $\uds{\bf Map}^{\rm fr}$. We only describe the flow category case. %
Given $(\phi, u, F) \in Y_{pq}^{{\mb F}}$ and $(\phi', u', F') \in Y_{qr}^{{\mb F}}$, one has defined the composition 
\beqn
\phi \# \phi':= \zeta_{prq}^B(\phi, \phi') \in B_{d_{pq} + d_{qr}} = B_{pr}
\eeqn
using the structural maps of the monotone flow category $\mb{Dom}^k$ with the corresponding curve $C_{\phi \# \phi'}$. The maps $u$ and $u'$ obviously define a Floer admissible map
\beqn
u \# u': C_{\phi \# \phi'} \to X.
\eeqn
Then for the two frames $F = (f_0, \ldots, f_d)$ and $F' = (f_0', \ldots, f_{d'}')$, we will define $F\# F'$ consisting of holomorphic sections of $L_{u\# u', H} \to C_{\phi \# \phi'}$. We can extend $L_{u, H}$ resp. $L_{u', H}$ to the other component of $C_{\phi \# \phi'}$ by trivial bundle and denote the extension still by the same notation. Then one has the natural isomorphism
\beqn
L_{u \# u', H} = L_{u, H} \otimes L_{u', H}.
\eeqn
Let $\tau_u \in L_u|_{z_+}$ be the parallel transport of $f_0(z_-)$ along the lateral line with respect to the connection on $L_u \to C_\phi$. Then in the same way as we define the concatenating $\phi \# \phi'$, consider the frame
\beqn
F \# F' = \Big( f_0 \otimes f_0', \ldots, f_d \otimes f_0', \tau_u \otimes f_1', \ldots, \tau_u \otimes f_{d'}' \Big).
\eeqn
Hence we have obtained a map 
\beqn
\hat\zeta_{prq}: Y_{pr}^{\mb F} \times Y_{rq}^{\mb F} \to \partial^{prq} Y_{pq}^{\mb F}.
\eeqn
We can also prove the associativity of the composition maps as the case of the monotone flow category $\mb{Dom}$. Hence we obtained a flow category. 

\begin{prop}\label{prop:reference-diagram}
\begin{enumerate}

\item The spaces $Y_{pq}^{\mb F}$ together with the composition maps form a flow category $\mb{Map}_{\mb F}^{\mb{fr}}$ enriched in $\outer \uds{\bf Map}^{\rm fr}_{\rm rig}$ which lifts the flow category $\mb{Dom}^{\Omega, {\mb F}}$. 

\item The spaces $Y_{p_1\cdots p_m; p'}^{{\mb X}, {\rm fr}}$ together with the structural maps form a multimodule $\mb{Map}_{\mb X}^{\mb{fr}}$ enriched in $\outer \uds{\bf Map}^{\rm fr}_{\rm rig}$ over $(\mb{Map}_{{\mb F}_1}^{\mb{fr}}, \cdots, \mb{Map}_{{\mb F}_m}^{\mb{fr}}; \mb{Map}_{{\mb F}'}^{\mb{fr}})$.

\item The spaces $Y_{p_1 \cdots p_m; p'}^{{\mb H}, {\rm fr}}$ together with the structural maps form a homotopy $\mb{Map}_{\mb H}^{\mb{fr}}$ enriched in $\outer \uds{\bf Map}^{\rm fr}_{\rm rig}$ from $\mb{Map}_{{\mb X}_0}^{\mb{fr}}$ to $\mb{Map}_{\mb{X}_1}^{\mb{fr}}$.
\end{enumerate}
\end{prop}

\begin{proof}
We just need to point out that the difference between the image of ${\rm Map}_{pr}^{\rm fr}\times {\rm Map}_{rq}^{\rm fr}$ and the stratum $\partial^{prq} {\rm Map}_{pq}^{\rm fr}$ lies in the difference of the domains. The concatenations of frame maps into $X$ do not affect the rigidification.
\end{proof}

\subsubsection{Group reductions}

We have the following proposition showing the existence of compatible systems of group reductions.

\begin{prop}\label{prop_group_reduction}
Let ${\mb F}$, ${\mb F}_1, \ldots, {\mb F}_m, {\mb F}', {\mb X}, {\mb X}_0, {\mb X}_1, {\mb H}$ be as above. Let $\outer \mb{Map}_{\mb F}^{\mb{fr}}, \cdots, \outer \mb{Map}_{{\mb H}}^{\mb{fr}}$ be the corresponding flow categories, multimodules, and homotopy enriched in enriched in $\outer \uds{\bf Map}^{\rm fr}_{\rm rig}$.
\begin{enumerate}

\item There exists a lift of $\outer \mb{Map}_{\mb F}^{\mb{fr}}$ to $\outer \uds{\bf map}^{\rm fr}_{\rm rig}$, denoted by $\outer \mb{map}{}_{\mb{F}}^{\mb{fr}}$.

\item Given lifts $\outer {\mb{map}}{}_{{\mb F}_1}^{\mb{fr}}, \cdots, \outer  {\mb{map}}{}_{\mb{F}_m}^{\mb{fr}}, \outer {\mb{map}}{}_{{\mb F}'}^{\mb{fr}}$ of $\outer \mb{Map}_{{\mb F}_1}^{\mb{fr}}, \cdots, \outer \mb{Map}_{{\mb F}_m}^{\mb{fr}}, \outer \mb{Map}_{\mb{F}'}^{\mb{fr}}$ respectively, there exists a lift of $\outer {\mb{Map}}_{\mb X}^{\mb{fr}}$ to $\outer \uds{\bf map}^{\rm fr}_{\rm rig}$ as a multimodule over $(\outer {\mb{map}}{}_{{\mb F}_1}^{\mb{fr}}, \cdots, \outer  {\mb{map}}{}_{\mb{F}_m}^{\mb{fr}}; \outer  {\mb{map}}{}_{{\mb F}'}^{\mb{fr}})$, which is denoted by $\outer {\mb{map}}{}_{{\mb X}}^{\mb{fr}}$.

\item Given lifts $\outer {\mb{map}}{}_{{\mb F}_1}^{\mb{fr}}, \cdots, \outer   {\mb{map}}{}_{\mb{F}_m}^{\mb{fr}}, \outer {\mb{map}}{}_{{\mb F}'}^{\mb{fr}}$ of $\outer \mb{Map}_{{\mb F}_1}^{\mb{fr}}, \cdots, \outer \mb{Map}_{{\mb F}_m}^{\mb{fr}}, \outer \mb{Map}_{\mb{F}'}^{\mb{fr}}$ to $\outer \uds{\bf map}^{\rm fr}_{\rm rig}$ and lifts $\outer { \mb{map}}{}_{{\mb X}_0}^{\mb{fr}}$, $\outer {\mb{map}}{}_{{\mb X}_1}^{\mb{fr}}$ of $\outer \mb{Map}_{\mb{X}_0}^{\mb{fr}}, \outer \mb{Map}{}_{{\mb X}_1}^{\mb{fr}}$ respectively, as multimodules over the tuple $(\outer {\mb{map}}{}_{{\mb F}_1}^{\mb{fr}}, \cdots, \outer {\mb{map}}{}_{\mb{F}_m}^{\mb{fr}}; \outer {\mb{map}}{}_{{\mb F}'}^{\mb{fr}})$, there exists a lift of $\outer \mb{Map}{}_{{\mb H}}^{\mb{fr}}$ to $\outer \uds{\bf map}^{\rm fr}_{\rm rig}$ as a homotopy from $\outer{ \mb{map}}{}_{{\mb X}_0}^{\mb{fr}}$ to $\outer {\mb{map}}{}_{{\mb X}_1}^{\mb{fr}}$.
\end{enumerate}
\end{prop}

\begin{proof}
We prove by induction. The construction of the base situation was provided in \cite{AMS}, whose construction is recalled as follows. Given $p, q$ and a framed map $(\phi, u, F) \in {\rm Map}_{pq}^{\mb{F},{\rm fr}}$, consider the $(d_{pq}+1)\times (d_{pq} +1)$ positive-definite Hermitian matrix
\beq\label{frame_pairing}
(H_F)_{ab} = \int_{C_\phi} \langle f_a, f_b \rangle \Omega_u.
\eeq
As frames are well-defined up to ${\mb C}^*$-rescaling, this Hermitian matrix is well-defined up to positive rescaling, hence defines a map 
\beqn
\Lambda_{pq}^{\mb F}: Y_{pq}^{{\mb F}} \to {\mc G}_{pq}/ G_{pq} = P_{pq}.
\eeqn
This map is obviously equivariant with respect to the ${\mc G}_R$-action; moreover, as ${\mc G}_L$ acts by reparametrizing the map and pulling back the framings, $\Lambda_{pq}^{\mb F}$ is ${\mc G}_L$-invariant. Therefore, one obtains a group reduction on $Y_{pq}^{\mb F}$. 

If $Y_{pq}^{\mb F}$ has no lower stratum, then its outercollaring is itself and hence one obtains the base case of the induction. Suppose inductively, one has constructed compatible collared group reductions
\beqn
\Lambda_{rs}^{\mb F}: \outer Y_{rs}^{\mb F} \to P_{rs}
\eeqn
for all pairs $r, s$ satisfying ${\mc A}(r) - {\mc A}(s) < {\mc A}(p) - {\mc A}(q)$. Notice that for each triple $p< r < q$, one has a strict embedding
\beqn
\outer Y_{pr}^{\mb F} \times \outer Y_{rq}^{\mb F} \to \partial^{prq} \outer Y_{pq}^{\mb F}
\eeqn
which is equivariant with respect to the complex group embedding $\hat {\mc G}_{pr}\times \hat {\mc G}_{rq} \to \hat {\mc G}_{pq}$. Hence one obtains an extension
\beqn
\vcenter{ \xymatrix{  \outer Y_{pr}^{\mb F} \times \outer Y_{rq}^{\mb F} \ar[r]^-{\hat \zeta_{prq}} \ar[d]_{\Lambda_{pr}^{\mb F} \times \Lambda_{rq}^{\mb F} } & \partial^{prq} \outer Y_{pq}^{\mb F} \ar[d]^{\Lambda_{pq}^{\mb F}}  \\
 P_{pr} \times P_{rq}  \ar[r]_-{\zeta_{prq}^P} &  P_{pq} }}.
\eeqn
The compatibility assumption and the $\hat {\mc G}_{pq}$-equivariance imply that one obtains a group reduction on the boundary of $\outer Y_{pq}^{\mb F}$ which is collared. We also have a group reduction in $Y_{pq}^{\mb F}$ (viewed as a subset away from boundary) given by \eqref{frame_pairing}. Then using the product structure on the collar region (where the actions by the complex Lie groups fix the collar coordinates) and choosing a homotopy, one can interpolate over the collar region to get a continuous group reduction on the whole of $\outer Y_{pq}^{\mb F}$. 

The constructions for the cases of multimodules and homotopies are completely the same and we omit the details.
\end{proof}

\subsubsection{Linearizing the group reduction}

We would like to regard the group reduction as a section of a vector bundle rather than a nonlinear fiber bundle. This requires linearizing the homogeneous space ${\mc G}/G$, which has a contractible set of choices. 

Recall that we have constructed a collection of germs of maps
\beqn
\lambda_d: \outer B_d \times Q_d \to \outer B_d \times P_d
\eeqn
By abuse of notation, let $Q_d^\epsilon \subset Q_d$ be a $G_d$-invariant disk bundle over which $\lambda_d$ is defined and $P_d^\epsilon \subset P_d$ the image of $Q_d^\epsilon$. Denote
\beqn
\lambda_{pq}: \outer B_{pq}^{\mb F} \times Q_{pq} \to \outer B_{pq}^{\mb F} \times P_{pq}
\eeqn
the pullback of $\lambda_{d_{pq}}$. Then define
\beqn
\outer Y_{pq}^{{\mb F}, \epsilon}:= (\Lambda_{pq}^{\mb F})^{-1} (P_{pq}^\epsilon)
\eeqn
and define the composition
\beqn
\lambda_{pq}:= \left( \xymatrix{ \outer Y_{pq}^{\mb{F}, \epsilon}  \ar[r]^-{\Lambda_{pq}^{\mb F}}   &     P_{pq}^\epsilon \ar[r]^{\lambda_{pq}^{-1}} & Q_{pq}^\epsilon \ar[r]& Q_{pq} } \right).
\eeqn
Viewing this as a section of a certain trivial vector bundle, this map will become part of the Kuranishi sections. We call such a section $\lambda_{pq}$ a {\bf linearized group reduction}.

\subsection{The topological Kuranishi lifts}\label{subsection163}

\subsubsection{Thickening data}

Given an equivariant family ${\mc C} = ({\mc G}, B, C)$ and an almost complex manifold $(X, J)$. Recall $\mathring C \subset C$ is the complement of nodes and markings. Consider the vector bundle
\beqn
\Lambda^{0,1}_{\mathring C/ B} \otimes TX \to \mathring C \times X.
\eeqn 

\begin{defn}\label{defn_thickening_datum}
A {\bf thickening datum} on an equivariant family of curves ${\mc C} = ({\mc G}, B, C)$ consists of a pair $(W, \nu)$ where $W$ is a finite-dimensional complex unitary representation\footnote{We restrict to complex representation in order to simplify the discussion of normal complex structures.} of the maximal compact subgroup $G \subset {\mc G}$ and 
\beqn
\nu: W \to \Gamma_c( \mathring C \times X, \Lambda^{0,1}_{\mathring C/B} \otimes TX)
\eeqn
is a $G$-equivariant linear map. Here $\Gamma_c$ stands for smooth sections with compact support.
\end{defn}

Notice that given a thickening datum $(W, \nu)$ on ${\mc C}$, for each $\phi \in B$ and a smooth map $u: C_\phi \to X$, there is an induced linear map
\beqn
\nu_{\phi, u}: W \to \Omega^{0,1}_c( \mathring C_\phi, u^* TX)
\eeqn
given by restricting $\nu (e)$ to the graph of $u$ in $\mathring C_\phi \times X$. These maps will be used to thicken up the moduli space of Floer equation.

\begin{defn}
The {\bf category of curves-with-thickening data}, denoted by $\uds{\bf Thick}$ is defined as follows. The objects are tuples
\beqn
({\mc G}, B, C, W, \nu )
\eeqn
where $({\mc G}, B, C)$ is an equivariant family of curves and $(W, \nu )$ is a thickening datum on ${\mc C}$.

A {\bf strict embedding of perturbations} consists of a strict embedding 
\beqn
\zeta_{21} = (\zeta_{21}^{\mc G}, \zeta_{21}^B, \zeta_{21}^C): ({\mc G}_1, B_1, C_1) \to ({\mc G}_2, B_2, C_2)
\eeqn
(see Definition \ref{defn_curve_2}) and a  $\zeta_{21}^G$-equivariant isometric linear embedding $\zeta_{21}^W: W_1 \to W_2$ such that 1) the following diagram commutes;
\beqn
\xymatrix{    W_1 \ar[r]^-{\nu_1} \ar[d]_{\zeta_{21}^W}  &  \Gamma_c( \mathring C_1 \times X, \Lambda_{\mathring C_1/B_1}^{0,1} \otimes TX)  \ar[d]\\
  W_2 \ar[r]_-{\nu_2}  &   \Gamma_c( \mathring C_2 \times X, \Lambda_{\mathring C_2/B_2}^{0,1} \otimes TX) }
\eeqn
(here the right vertical arrow is induced from the morphism $\zeta_{21}$) and 2) the restriction of $\nu_2$ to the orthogonal complement $W_1^\bot \subset W_2$ vanishes. 

The conjugation of $({\mc G}, B, C, W, \nu)$ is a strict embedding induced by $g \in G$. An {\bf embedding} of curvs-with-thickening-datum is a conjugacy class of strict embeddings.
\end{defn}

The category $\uds{\bf Thick}$ is a regular stratification category with a forgetful functor $\uds{\bf Thick} \to \uds{\bf Curv}$. One can easily formulate the category $\outer \uds{\bf Thick}$ of collared curves-with-thickening-datum, together with forgetful functor 
\beqn
\outer\uds{\bf Thick} \to \outer\uds{\bf Curve}.
\eeqn

We need to consider the rigidified version. Recall that for a strict embedding 
\beqn
\zeta_{21}: ({\mc G}_1, B_1, C_1) \to ({\mc G}_2, B_2, C_2)
\eeqn
a rigidification is a germ of $G_1$-equivariant map $\theta_{21}: B_1 \times Q_{21} \to {\mc G}_2$ satisfying an appropriate condition. 

\begin{defn}
Let $(\zeta_{21}, \zeta_{21}^W)$ be a strict embedding in $\uds{\bf Thick}$ from $({\mc G}_1, B_1, C_1, W_1, \nu_1)$ to $({\mc G}_2, B_2, C_2, W_2, \nu_2)$. A {\bf rigidification} of $(\zeta_{21}, \zeta_{21}^W)$ consists of a rigidification $\theta_{21}: B_1 \times Q_{21} \to {\mc G}_2$ of $\zeta_{21}$ satisfying the following condition. Fix $b_1 \in B_1$. The rigidification induces identifications
\beqn
C_{\zeta_{21}^B(b_1)} \cong C_{\theta_{21}(b_1, q)( \zeta_{21}^B (b_1))}
\eeqn
for all $q \in Q_{21}$ sufficient close to the origin. Then we can regard $\nu_2$ as a family of linear maps
\beqn
\nu_2(b_1, q): W_2 \to \Gamma_c( \mathring C_{\zeta_{21}(b_1)} \times X, \Lambda^{0,1} \otimes TX).
\eeqn
We require that these family only depends on $b_1$. In particular, for $q$ close to the origin, $\nu_2$ vanishes on the orthogonal complement of $W_1$.
\end{defn}

\begin{defn}
A {\bf thickening datum} on a flow category/multimodule/homotopy enriched in $\uds{\bf Curve}$ is a lift in $\uds{\bf Thick}$. If the flow category/multimodule/homotopy is enriched in $\outer \uds{\bf Curve}_{\rm rig}$, then a perturbation is called {\bf collared and rigidified} if it is a lift in $\outer\uds{\bf Thick}_{\rm rig}$.
\end{defn}

In the following steps, one needs to choose perturbations satisfying additional properties. The construction of such perturbations will be provided in due course.

Now we can describe the AMS lift of the Floer flow categories, bimodules, and homotopies to the category of topological Kuranishi spaces.

\subsubsection{The thickenings}

We describe the thickened moduli spaces as candidates of the Kuranishi regularizations of the moduli spaces. We first consider the flow category case. Return to the domain flow category $\outer \mb{Dom}^{\Omega, {\mb F}}$ enriched in $\outer\uds{\bf Curve}_{\rm rig}$ whose morphism spaces are equivariant families of curves ${\mc C}_{pq}^{\mb F}$. A thickening datum of $\outer \mb{Dom}^{\Omega, {\mb F}}$ is denoted by $\mb{W}^{\mb F}$ which consists of thickening data
\beqn
\nu_{pq}^{\mb F}: W_{pq}^{\mb F} \to \Gamma_c ( \mathring C_{pq}^{\mb F} \times X, \Lambda_{\mathring C_{pq}^{\mb F}/ B_{pq}^{\mb F}}^{0,1}\otimes TX).
\eeqn

\begin{defn}\label{defn_thickened_moduli}
Given any collared and rigidified thickening datum $\mb{W}^{\mb F}$ on the flow category $\outer \mb{Dom}^{\Omega, {\mb F}}$ and the flow category $\outer {\mb{map}}{}_{{\mb F}}^{\mb{fr}}$ of reduced framed maps, for each pair $p, q$ of objects in the Floer flow category ${\mb F}$, define the following items.

\begin{enumerate}

\item The {\bf thickened moduli space} is
\beqn
V_{pq}:= V_{pq}^{\mb F}:= \left\{ (\phi, u, F, e)\ \left|\ \begin{array}{c} (\phi, u, F) \in Y_{pq}^{{\mb F}, \epsilon},\ e \in W_{pq}^{\mb F},\\
\ov\partial_H u + \nu_{pq}^{\mb F} (e) = 0,\ \phi = \phi_F \in \outer B_{pq}^{\mb F}. \end{array} \right. \right\}.
\eeqn

\item The $G_{pq}$-action on $V_{pq}$ is the product of the diagonal $G_{pq}$-action on $Y_{pq}^{\mb{F}, \epsilon}$ space of framed maps and the action on $W_{pq}$.

\item The {\bf obstruction bundle} is the trivial bundle
\beqn
E_{pq}:= E_{pq}^{\mb F}:= \uds Q_{pq}^{\mb F} \oplus 
\uds W_{pq}^{\mb F} \to V_{pq}.
\eeqn

\item The {\bf Kuranishi section} is (given by) the equivariant map
\beqn
S_{pq}^{\mb F} (\phi, u,  F, e) = (\lambda_{pq}^{\mb F}(\phi, u, F), e).
\eeqn

\end{enumerate}
\end{defn}

\begin{thm}\label{thm_footprint}
The natural ``footprint'' map
\beqn
 (S_{pq}^{\mb F})^{-1}(0)/ G_{pq}^{\mb F} \to \outer \ov{\mc M}_{pq}^{\mb F}
\eeqn
which sends the orbit of $(\phi, u, F, 0)$ to the equivalence class of the stable Floer trajectory $u: C_\phi \to X$ is an isomorphism of orbispaces.
\end{thm}

\begin{proof}
We drop ${\mb F}$ from the notations for simplification. We also assume $k=1$. We prove that the footprint map is surjective. Given any point of $\outer \ov{\mc M}_{pq}$ represented by a stable Floer trajectory $u: \Sigma \to X$ (assume that this point is not in the outer collar; the case that it is in the outer collar can be deduced similarly). Then $u$ determines a unitary line bundle $L_u  \to \Sigma$ whose Chern connection has curvature $\Omega_u$. Let $F$ be a frame of $L_u$ satisfying the negative (or positive) constraint, which determines a holomorphic map $\phi_F: \Sigma \to {\mb P}^{d_{pq}}$ and hence can be identified with a point $\phi_F \in B_{pq}$. The domain $\Sigma$ can be then identified with $C_{\phi_F}$. Let $u_F: \mathring C_{\phi_F} \to X$ be the pullback map. As two different frames $F$ and $F'$ differ by an element of the group ${\mc G}_{pq, R}$, one obtains a ${\mc G}_{pq}$-orbit (where we regard ${\mc G}_{pq}$ as the diagonal of ${\mc G}_{pq, L} \times {\mc G}_{pq, R}$) of framed maps 
\beqn
( \phi_F, u_F, F) \in Y_{pq}^{\mb F}.
\eeqn
By the definition of a structural group reduction scheme, there is a unique $G_{pq}$-orbit (for the diagonal action) of framed maps such that 
\beqn
\Lambda_{pq}^{\mb F}(\phi_F, u_F, F) = [1] \in P_{pq}
\eeqn
Then $(\phi_F, u_F, F) \in Y_{pq}^{\mb{F}, \epsilon}$. Hence $\lambda_{pq}^{\mb F}(\phi_F, u_F, F) = 0$. This $G_{pq}$-orbit is then sent by the footprint map to the given point of $\ov{\mc M}_{pq}$. Therefore the footprint map is surjective.

We prove that the footprint map is injective. Suppose $(\phi, u, F, 0)$ and $(\phi', u', F', 0)$ are sent to the same point of $\outer \ov{\mc M}_{pq}$. By definition, there is a domain isomorphism $C_\phi \cong C_{\phi'}$, which must be given by an element $g \in {\mc G}_{pq}$. As $\phi = \phi_F$, $\phi' = \phi_{F'}$, it follows that 
\beqn
g (\phi, u, F, 0) = (\phi', u', F',0).
\eeqn
As in the same ${\mc G}_{pq}$-orbit there can only be one normalized $G_{pq}$-oribt, it follows that $g \in G_{pq}$. This implies that the footprint map is injective.
\end{proof}

Now we turn to the case of multimodules ${\mb X}:= \mb{X}^\Sigma$ associated to a connected Floer domain $\Sigma$ over $({\mb F}_1, \ldots, {\mb F}_m; {\mb F}')$. Previously one has constructed a multimodule $\outer{\mb{Dom}}^{\mb X}$ enriched in $\outer\uds{\bf Curve}_{\rm rig}$, together with the group reductions
\beqn
\lambda_{p_1 \cdots p_m; p'}^{\mb X}: \outer Y_{p_1 \cdots p_m; p'}^{{\mb X}} \to {\mc G}_{p_1\cdots p_m; p'}/ G_{p_1 \cdots p_m; p'}
\eeqn
which extends the group reductions on the involved flow categories of framed maps.

\begin{defn}\label{defn_module_thickening}
Given a thickening datum $\mb{W}^\Sigma$ on the multimodule $\outer \mb{Dom}^{\Omega, {\mb X}}$ (as a multimodule over the flow categories equipped with thickening data, enriched in $\outer \uds{\bf Thick}_{\rm rig}$), for objects $p_1, \ldots, p_m; p'$, define the following objects.

\begin{enumerate}
    \item The {\bf thickened moduli space} is 
    \beqn
    V_{p_1 \cdots p_m; p'}:= V_{p_1 \cdots p_m; p'}^{\mb X}:= \left\{ (\phi, u, F, e)\ \left| \begin{array}{l}  (\phi, u, F) \in \outer Y_{p_1 \cdots p_m; p'}^{{\mb X}, \epsilon},\ e \in W_{p_1 \cdots p_m; p'}^{\mb X},\\
    \ov\partial_\sigma u + \nu_{p_1 \cdots p_m; p'}^{\mb X} (e) = 0,\ \phi = \phi_F \in \outer B_{p_1 \cdots p_m; p'}^{\mb{X}}. \end{array} \right.  \right\}
\eeqn

\item The $G_{p_1 \cdots p_m; p'}$-action on $V_{p_1 \cdots p_m; p'}$ is the restriction of the diagonal action on $\outer Y_{p_1 \cdots p_m; p'}^{\mb{X}, \epsilon} \times W_{p_1 \cdots p_m; p'}^{\mb X}$. 

\item The {\bf obstruction bundle} is the trivial bundle
\beqn
E_{p_1 \cdots p_m; p'}:=E_{p_1 \cdots p_m; p'}^{\mb X}:=  \uds Q_{p_1 \cdots p_m; p'} \oplus \uds W_{p_1 \cdots p_m; p'} \to V_{p_1 \cdots p_m; p'}.
\eeqn

\item The {\bf Kuranishi section} is (equivalent to) the map
\beqn
S_{p_1\cdots p_m; p'}^{\mb X} (\phi, u, F, e) \mapsto (\lambda_{p_1 \cdots p_m; p'}^{\mb X} (\phi, u, F), e)
\eeqn
where $\lambda_{p_1 \cdots p_m; p'}$ is the linearized group reduction. 
\end{enumerate}
\end{defn}

The case of homotopies is the same. We omit the detailed discussions. We will use ${\mb H}$ to indicate the relationship with a multimodule homotopy ${\mb H}$ coming from the Floer theory background.

\subsubsection{Transversality}

Next we deal with transversality. For this we need to set up a certain Fredholm problem. For each pair of objects $p, q$, define
\beqn
V_{pq}':= \Big\{ (\phi, u, e) |\ \phi \in \outer B^{\mb F}_{pq}, e \in W_{pq}^{\mb F}, u: C_\phi \to X \text{ satisfies }\ov\partial_H u + \nu_{pq}^{\mb F} (e) = 0 \Big\}.
\eeqn

\begin{lemma}
There is a canonical $G_{pq}$-equivariant isomorphism $V_{pq}' \cong V_{pq}$.
\end{lemma}

\begin{proof}
We show the natural map $V_{pq} \to V_{pq}'$ is bijective. Indeed, given any $(\phi, u, e) \in V_{pq}'$, the set of all frames $F$ on it are related by an element of ${\mc G}_R$ which linearly transforms the frame. Each frame $F$ induces a map $\phi_F: C_\phi \to \mb{P}^{\bf d}$. The condition that $\phi_F = \phi$ specifies a unique frame. Hence the map $V_{pq} \to V_{pq}'$ is surjective. As it is obviously injective and $G_{pq}$-equivariant, the lemma follows.
\end{proof}

For a fixed $\phi$, let $V_{pq}'(\phi) \subset V_{pq}'$ be the preimage of $\phi$. Then $V_{pq}'(\phi)$ is the zero locus of a smooth Fredholm section defined on a certain Banach manifold of maps from $C_\phi$ to $X$; the concrete choice of the Banach manifold depends on certain Sobolev parameters. 

\begin{defn}
We say that $(\phi, u, e) \in V_{pq}'(\phi)$ is {\bf transverse} if the linearization of the corresponding Fredholm section is surjective (which is independent of the Sobolev parameters). We say that the perturbation $\mb{W}$ is {\bf transverse} if all points of $V_{pq}'(\phi)$ for all $p,q$ are transverse. 
\end{defn}

We can choose the thicikening data appropriately so that transversality holds.

\begin{prop}\label{prop_perturbation} We omit the dependence on $\Omega$ from the notations.
\begin{enumerate}

\item Given any integral modification of the topological energy, there exists a transverse thickening datum $\mb{W}^{\mb F}$ on the flow category $\outer \mb{Dom}^{{\mb F}}$, i.e., a transverse lift to $\outer \uds{\bf Thick}_{\rm rig}$.

\item Consider the multimodule $\outer \mb{Dom}^{\mb{X}}$ over $(\outer \mb{Dom}^{{\mb F}_1}, \cdots \outer \mb{Dom}^{{\mb F}_m}; \outer \mb{Dom}^{{\mb F}'})$. Suppose $\mb{W}^{{\mb F}_1}$, $\cdots$, $\mb{W}^{{\mb F}_m}; \mb{W}^{{\mb F}'}$ are transverse thickening data on these flow categories. Then there exists a transverse thickening datum ${\mb W}^{\mb X}$ which is a multimodule over $(\mb{W}^{{\mb F}_1}, \cdots \mb{W}^{{\mb F}_m}; \mb{W}^{{\mb F}'})$.

\item Suppose $\mb{W}^{{\mb X}_0}$ resp. $\mb{W}^{{\mb X}_1}$ are two transverse thickening data on ${\mb X}_0$ resp. $\mb{X}_1$ as multimodules over $(\mb{W}^{{\mb F}_1}, \cdots, \mb{W}^{{\mb F}_m}; \mb{W}^{{\mb F}'})$. Then there exists a lift of $\outer \mb{Dom}^{{\mb H}}$ to $\outer \uds{\bf Thick}_{\rm rig}$, denoted by $\mb{W}^{{\mb H}}$, as a homotopy from $\mb{W}^{{\mb X}_0}$ to $\mb{W}^{{\mb X}_1}$.
\end{enumerate}
\end{prop}

\begin{proof}
See Subsubsection \ref{subsubsection_perturbation}.
\end{proof}

Transversality implies the regularity of the thickened moduli spaces. 

\begin{prop}
Suppose we choose transverse thickening data as claimed in Proposition \ref{prop_perturbation}.
\begin{enumerate}
    \item For each $p, q \in {\rm Ob}{\mb F}$, the thickened moduli space $V_{pq}^{\mb F}$ is an $A_{pq}^{\mb F}$-stratified topological $G_{pq}^{\mb F}$-manifold near $(S_{pq}^{\mb F})^{-1}(0)$. 

    \item Each thickened moduli space $V_{p_1 \cdots p_m; p'}^{\mb X}$ is an $A_{p_1 \cdots p_m; p'}^{\mb X}$-stratified topological $G_{p_1 \cdots p_m; p'}^{\mb X}$-manifold near $(S_{p_1 \cdots p_m; p'}^{\mb X})^{-1}(0)$.

    \item Each thickened moduli space $V_{p_1 \cdots p_m; p'}^{\mb H}$ is an $A_{p_1 \cdots p_m; p'}^{\mb H}$-stratified topological $G_{p_1 \cdots p_m; p'}^{\mb H}$-manifold near $(S_{p_1 \cdots p_m; p'}^{\mb H})^{-1}(0)$.
    \end{enumerate}
\end{prop}

\begin{proof}
We only discuss the flow category case. By definition, transversality implies that each $V_{pq}'(\phi) \cong V_{pq}$ is a topological manifold near where $e = 0$ for all $\phi \in B_{pq}$. The gluing argument implies that $V_{pq}'$ is a stratified manifold. Since the thickening datum is collared, $V_{pq}'$ is also collared. Finally the canonical identification $V_{pq}' \cong V_{pq}$ shows $V_{pq}$ is a collared manifold.
\end{proof}

\subsubsection{The topological AMS lifts}\label{subsubsec:top-AMS}

Now assuming transversality, by abuse of notations, the ``regular locus'' of $V_{pq}$ is still denoted by $V_{pq}$ (similar treatment for thickenings for multimodules and homotopies). Then one obtains topological Kuranishi spaces
\beqn
K_{pq} = (G_{pq}, V_{pq}, E_{pq}, S_{pq}),\ p, q \in {\rm Ob}{\mb F}.
\eeqn
Together with the structural maps, one obtains a lift of $\outer {\mb F}$ to the category $\outer \uds{\bf Kur}$, denoted by $\hat {\mb F}$. 

\begin{lemma}\label{lemma_Kuranishi_rigidification}
The topological AMS lift $\hat{\mb F}$ can be canonically lifted to $\outer \uds{\bf Kur}_{\rm rig}$.
\end{lemma}

\begin{proof}
We only need to construct rigidifications of the embeddings
\beqn
\iota_{prq} : K_{pr} \times K_{rq} \to \partial^{prq} K_{pq}.
\eeqn
By definition, one needs to find an orthogonal representation $N_{prq}$ of $G_{pr}\times G_{rq}$ and a germ of open embedding
\beqn
G_{pq}\times_{G_{pr}\times G_{rq}} ( K_{pr} \times K_{rq} \times N_{prq}) \to \partial^{prq} K_{pq}
\eeqn
which extends $\iota_{prq}$. Let
\beqn
W_{prq} \subset W_{pq}
\eeqn
be the orthogonal complement of $W_{pr} \oplus W_{rq} \subset W_{pq}$. Moreover, the underlying strict embedding in $\outer \uds{\bf Curve}$ is also rigidified via a representation $Q_{prq} \subset Q_{pq}$. Set
\beqn
N_{prq}:= W_{prq} \oplus Q_{prq}.
\eeqn
Then one can check that there is a natural strict open embedding
\beqn
G_{pq} \times_{G_{prq}} \left( {\rm Stab}_{N_{prq}^\epsilon} ( K_{pr} \times K_{rq}) \right) \to \partial^{prq} K_{pq}
\eeqn
(for a sufficiently small ball $N_{prq}^\epsilon \subset N_{prq}$) which extends $\iota_{prq}^{\mb F}$.
\end{proof}

\subsubsection{Multimodules and homotopies}

Similarly, for the flow multimodule ${\mb X}$ over $({\mb F}_1, \ldots, {\mb F}_m; {\mb F}')$ induced from a Floer domain $\Sigma$, choose integral modifications 
\begin{align*}
&\ {\mc A}_{{\mb F}_i}^\Omega: {\rm Ob}{\mb F}_i \to {\mb Z}^k,\ &\ {\mc A}_{{\mb F}'}^\Omega: {\rm Ob}{\mb F}' \to {\mb Z}^k
\end{align*}
(and the constant $n_{\mb X}\in {\mb Z}$, see Definition \ref{module_integral_action}). Make choices for 
\begin{align*}
&\ \outer {\mb{map}}{}_{{\mb F}_i}^{\mb{fr}},\ &\ \outer {\mb{map}}{}_{{\mb{F}'}}^{\mb{fr}},
\end{align*}
namely group reductions, by Lemma \ref{prop_group_reduction}, one can find a multimodule
\beqn
\outer {\mb{map}}{}_{{\mb X}}^{\mb{fr}}
\eeqn
over $(\outer {\mb{map}}{}_{\mb{F}_1}^{\mb{fr}}, \cdots \outer {\mb{map}}{}_{{\mb F}_m}^{\mb{fr}}; \outer {\mb{map}}{}_{{\mb F}'}^{\mb{fr}})$ enriched in the category $\outer \uds{\bf map}{}_{\rm rig}^{\rm fr}$. 

In parallel, choose collared and rigidified  thickening data $\mb{W}^{{\mb F}_i}$ resp. $\mb{W}^{\mb{F}'}$ over $\outer \mb{Dom}^{\Omega, {\mb F}_i}$ resp. $\outer \mb{Dom}^{{\bm \Omega}, {\mb F}'}$, by Proposition \ref{prop_perturbation}, one can find a thickening datum ${\mb W}^{{\mb X}}$ as a multimodule over $(\mb{W}^{{\mb F}_1}, \ldots, \mb{W}^{{\mb F}_m}; \mb{W}^{\mb{F}'})$. Then the data $\outer {\mb{map}}{}_{\mb X}^{\mb{fr}}$ and $\mb{W}^{\mb X}$ induce the collection of topological Kuranishi spaces
\beqn
K_{p_1 \cdots p_m; p'}^{\mb X} = (G_{p_1 \cdots p_m; p'}^{\mb X}, V_{p_1 \cdots p_m; p'}^{\mb X}, E_{p_1 \cdots p_m; p'}^{\mb X}, S_{p_1 \cdots p_m; p'}^{\mb X}).
\eeqn
This provides a lift of $\outer {\mb X}$ to $\outer \uds{\bf Kur}$, denoted by $\wh{\mb X}$ and called a {\bf topological AMS lift} of the multimodule $\outer {\mb X}$. It is a multimodule over $(\wh {\mb F}_1, \ldots, \wh {\mb F}_m; \wh {\mb F}')$. Similar to Lemma \ref{lemma_Kuranishi_rigidification}, this lift can be further lifted (in a natural way) to $\outer \uds{\bf Kur}{}_{\rm rig}$. 

The case of homotopies is similar and omited.

\subsubsection{Proof of Proposition \ref{prop_perturbation}}\label{subsubsection_perturbation}

The existence of a transverse thickening datum for a single moduli space was (essentially) given in \cite[Section 4.1]{AMS2}. To construct a transverse thickening datum on a flow category, one needs a more systematic treatment. 

All the constructions towards the proof is {\it ad hoc}.

\begin{defn}
Let $G$ be a Lie group.
\begin{enumerate}

\item A {\bf $\mu$-representation} of a Lie group $G$ is a direct system of complex representations $(W_\mu)_{\mu \geq 1}$ of $G$ where $W_\mu \to W_{\mu+1}$ is an embedding. A $\mu$-representation is said to be {\bf exhaustive} if any finite-dimensional representation $R$ of $G$ can be embedded into $W_\mu$ for $\mu$ sufficiently large. 

\item A single representation $W$ is regarded as a $\mu$-representation with $W_\mu \equiv W$.

\item A $\mu$-representation is called {\bf finite-dimensional} if all $W_\mu$ are finite-dimensional; is called {\bf unitary} if all $W_\mu$ are equipped with $G$-invariant Hermitian inner products and the embeddings $W_\mu \to W_{\mu+1}$ are isometric.

\item A {\bf strict embedding} of $\mu$-representations from $(W_\mu)$ to $(W_\mu')$ is a sequence of $G$-equivariant linear isometric embeddings $W_\mu \to W_{\mu}'$ such that 
\beqn
\xymatrix{  W_\mu \ar[r] \ar[d] & W_\mu' \ar[d] \\
 W_{\mu+1} \ar[r] & W_{\mu+1}' }
 \eeqn
 \end{enumerate}
\end{defn}

Recall the definition of the category $\uds{\bf Curve}$. One can define another category $\uds{\bf Curve}_\mu$ whose objects are $({\mc G}, B, C, (W_\mu))$ where $({\mc G}, B, C)$ is an object of $\uds{\bf Curve}$ and $(W_\mu)$ is a finite-dimensional, exhaustive, and unitary $\mu$-representation of $G \subset {\mc G}$. One can similarly define strict embeddings; the morphisms are conjugacy classes of strict embeddings. Then there is an obvious forgetful functor
\beqn
\uds{\bf Curve}_\mu \to \uds{\bf Curve}.
\eeqn

\begin{prop}\label{system_representation}
\hfill
\begin{enumerate}

\item There exist a lift of the monotone flow category $\mb{Dom}$ to $\uds{\bf Curve}_\mu$, denoted by $\mb{Dom}_\mu$.

\item For each smooth Floer domain $\Sigma$ with $m$ negative cylindrical ends and one positive cylindrical end\footnote{The case with zero positive cylindrical ends, e.g., the case for Poincar\'e pairing, can be proved as well.}, a lift of $\mb{Dom}^\Sigma$ to $\uds{\bf Curve}_\mu$, denoted by $\mb{Dom}^\Sigma_\mu$, which is a monotone  multimodule over $(\mb{Dom}_\mu, \ldots, \mb{Dom}_\mu; \mb{Dom}_\mu)$.

\item For any 1-parameter family of Floer domains $\Sigma^I$ and the associated homotopy $\mb{Dom}^I$ from $\mb{Dom}^{\Sigma_0}$ to $\mb{Dom}^{\Sigma_1}$, a lift to $\uds{\bf Curve}_\mu$, denoted by $\mb{Dom}_\mu^I$, which is a monotone homotopy from $\mb{Dom}_\mu^{\Sigma_0}$ to $\mb{Dom}_\mu^{\Sigma_1}$.
\end{enumerate}
\end{prop}

\begin{proof}
Consider the standard unitary representation $W_d = {\mb C}^d$ of $G_d \cong U(d)$. The induced representation of the exterior algebra $\Lambda W_d$ is also unitary. Then define the $\mu$-representation with 
\beqn
W_{d, \mu}:= \left( \bigoplus_{i=1}^\mu (\Lambda W_d)^{\otimes i} \right)^{\oplus 2^\mu}
\eeqn
whose Hermitian structure is the direct sum of the Hermitian structures of summands. The inclusion of $W_{d, \mu}$ into $W_{d, \mu+1}$ is the composition
\begin{multline*}
W_{d, \mu} = \left( \bigoplus_{i=1}^\mu (\Lambda W_d)^{\otimes i} \right)^{\oplus 2^\mu} \to 
\left( \bigoplus_{i=1}^{\mu+1} (\Lambda W_d)^{\otimes i} \right)^{\oplus 2^\mu} \oplus \{0\} \\
\subset 
\left( \bigoplus_{i=1}^{\mu+1} (\Lambda W_d)^{\otimes i} \right)^{\oplus 2^\mu} \oplus
\left( \bigoplus_{i=1}^{\mu+1} (\Lambda W_d)^{\otimes i} \right)^{\oplus 2^\mu} = W_{d, \mu+1}
\end{multline*}
which is isometric. As any finite-dimensional representation of $U(d)$ is contained in a sufficiently higher tensor power $(\Lambda W_d)^{\otimes i}$, one obtains a $\mu$-representation.

Then for any decomposition $d = d_0 + d_1$, one has the canonical isometry
\beqn
W_{d_0} \oplus W_{d_1} \cong W_d
\eeqn
which is equivariant with respect to $U(d_0) \times U(d_1) \to U(d)$. This map induces an equivariant isometric embeddings of the exterior algebras
\beqn
\Lambda W_{d_0} \oplus \Lambda W_{d_1} \to \Lambda W_d
\eeqn
and hence isometric embeddings
\beqn
W_{d_0, \mu} \oplus W_{d_1, \mu} \to W_{d, \mu}.
\eeqn
There is an obvious associativity property of such embeddings, which induces well-defined maps for each partition $d = d_0 + \cdots + d_l$
\beqn
W_{d_0, \mu} \oplus \cdots \oplus W_{d_l, \mu} \to W_{d, \mu}.
\eeqn
Then for each $d$, the quadruple $({\mc G}_d, B_d, C_d, (W_{d, \mu}))$ is an object of $\uds{\bf Curve}_\mu$ and they form a monotone flow category enriched in this category. The cases of multimodules and homotopies are the same; notice that in these cases, there are no well-defined strict embeddings but up to conjugacy, the embeddings are well-defined.
\end{proof}

We explain what we will do next. We explain for the flow category case. We will choose sufficiently large $\mu_{pq}$ for each pair $p, q$ inductively, such that $W_{pq, \mu_{pq}}$ can cover all possible cokernels of the linearized operator of the Floer equation. 

First we recall the following definition and lemma from \cite{AMS2} with slight modification.

\begin{defn}
Let $G$ be a compact Lie group acting smoothly on a manifold $B$ and $\pi: V \to B$ be a smooth complex $G$-vector bundle. A {\bf finite-dimensional approximation scheme} is a finite-dimensional $\mu$-representation $(V_\mu)$ and a morphism of $\mu$-representations $\lambda_\mu: V_\mu \to C_c^\infty(B, V)$  such that the union of $\lambda_\mu(V_\mu)$ over $\mu$ is dense with respect to the $C^\infty_{\rm loc}$-topology.
\end{defn}

\begin{lemma}\cite[Lemma 4.2]{AMS2}
There exists a finite-dimensional approximation scheme for given $G, B, V$. 
\end{lemma}

\begin{cor}
When $G = U(d)$, the finite-dimensional approximation scheme can be taken such that $V_\mu = W_{d, \mu}$ where $(W_{d, \mu})$ is the $\mu$-representation specified in the proof of Proposition \ref{system_representation}.
\end{cor}

\begin{proof}
This is because any $V_\mu$ is contained in some $W_{\mu'}$ for $\mu'$ sufficiently large.
\end{proof}

Now we can prove Proposition \ref{prop_perturbation}. 

\begin{proof}[Proof of Proposition \ref{prop_perturbation}]
We can choose a collection of collared thickening data
\beqn
\nu_{pq, \mu}: W_{pq, \mu} \to C_c^\infty( \outer \mathring C_{pq} \times X, \Lambda^{0,1}_{\outer \mathring C_{pq}/ \outer B_{pq}} \otimes TX)
\eeqn
such that they are compatible with the structural maps of the flow category. Then using the compactness of each $M_{pq}$, there exists $\mu_{pq}$ such that the image of $\nu_{pq, \mu_{pq}}$, when restricted to the domain of each solution, covers the cokernel of the linearization. Then we can choose $\mu_{pq}$ inductively. Then define $W_{pq}:= W_{pq, \mu_{pq}}$ and $\nu_{pq}:= \nu_{pq, \mu_{pq}}$. This provides a transverse thickening datum on the flow category. The inductive construction on multimodules and homotopies are exactly the same. 
\end{proof}

\subsection{Relatively smooth AMS lifts}\label{subsection_relative_smooth}

The topological AMS lifts actually admit refinement in the so-called relatively smooth category (see \cite{Hirschi_Swaminathan}). The discussion in this category is important for eventually obtaining smooth Kuranishi spaces. We introduce the category of Kuranishi spaces in this version. 

\subsubsection{Relative spaces and relative smooth structures}

\begin{defn}[Relative spaces]\label{defn_relative_space} \hfill   
\begin{enumerate}
\item A (stratified) {\bf relative space} consists of a pair of $A$-stratified topological spaces $Y$, $S$ and a stratified continuous map $\pi: Y \to S$. When $\pi$ is understood from the context, we denote a relative space by $Y/S$.
 
\item A {\bf relative map} from $Y/S$ to $Z/T$ is a commutative diagram
\beqn
\xymatrix{  Y \ar[r] \ar[d] & Z \ar[d] \\
            \pi(Y)  \ar[r]  & \pi(Z) }.
            \eeqn
When $S = T$, without further declaration, we will assume that the underlying map from $\pi(Y)$ to $\pi(Z)$ is the restriction of the identity map of $S$. 

\item A {\bf relative homeomorphism} is a relative map $\varphi: Y/S \to Z/T$ consisting of a homeomorphism from $Y$ to $Z$ and a homeomorphism from $S$ to $T$.

\item If $G$ is a topological group, a $G$-action on relative space is a continuous action via relative homeomorphisms. In particular, the map $\pi: Y \to S$ is equivariant.

\item Let $D$ be another topological space. The {\bf horizontal stabilization} of a relative space $Y/S$ by $D$ is the relative space 
\beqn
{\rm Stab}_D^{\rm hor}(Y/S) = (Y\times D)/(S \times D).
\eeqn
The {\bf vertical stabilization} of $Y/S$ is the relative space
\beqn
{\rm Stab}_D^{\rm vert}(Y/S) = (Y \times D)/S.
\eeqn
Notice that both horizontal and vertical stabilizations are functors on the categories of relative spaces. Moreover, given $D_1, D_2$, there is a natural transformation
\beqn
{\rm Stab}_{D_1}^{\rm hor} \circ {\rm Stab}_{D_2}^{\rm vert} \to {\rm Stab}_{D_2}^{\rm vert} \circ {\rm Stab}_{D_1}^{\rm hor}.
\eeqn

\end{enumerate}
\end{defn}

\begin{rem}
    A fiber bundle is a canonically a relative space. In particular, we denote by $\uds{\mb R}^k/S$ the relative space associated to the trivial ${\mb R}^k$ bundle over $S$.
\end{rem}

\begin{defn}
Let $U/S \subset \uds{\mb R}^m/S$ and $V/T \subset \uds{\mb R}^n/T$ be open subsets. A continuous map $\varphi: U/S \to V/T$ is called rel-$C^{\infty, k}$ if the map $F: U \to {\mb R}^n$, determined by 
\beqn
\varphi(v, x) = (F(v, x), \uds\varphi(x))
\eeqn
has all partial derivatives in $v \in {\mb R}^m$ and those up order $k$ varies continuously in $x \in S$ uniformly on compact subsets of ${\mb R}^m$. A rel-$C^{\infty, k}$-map is called a rel-$C^{\infty, k}$-diffeomorphism if it has a rel-$C^{\infty, k}$-inverse. $\varphi$ is called {\bf relatively smooth}, or rel-$C^\infty$, if it is rel-$C^{\infty, k}$ for all $k$.
\end{defn}

\begin{defn}
A {\bf relative smooth structure} on $Y/S$ consists of an atlas $\varphi_\alpha: U_\alpha/S \to \uds{\mb R}^n/S$ which are rel-$C^{\infty}$-compatible. More precisely 1) $Y = \cup U_\alpha$ and 2) for any $\alpha, \beta$, the map 
\beqn
\varphi_\alpha \circ \varphi_\beta^{-1}: \varphi_\beta( (U_\alpha \cap U_\beta)/S) \to \varphi_\alpha( (U_\alpha \cap U_\beta)/S)
\eeqn
is a rel-$C^{\infty}$-diffeomorphism.
\end{defn}

Notice that if $Y/S$ has a relative smooth structure, then each fiber of $Y$ has a well-defined smooth structure. Then there is a well-defined rel-$C^\infty$ vector bundle
\beqn
T^{\rm vert} Y/S \to Y,
\eeqn
called the {\bf vertical tangent bundle}.

\subsubsection{The category of relative smooth Kuranishi spaces}

\begin{defn}[Category of relatively smooth Kuranishi spaces] \hfill
\begin{enumerate}

\item A {\bf relative smooth Kuranishi space} is a quadruple $K/B = (G, V/B, E, S)$ where $(G, V, E, S)$ is a topological Kuranishi space, $B$ is a smooth $G$-manifold, together with a $G$-invariant rel-$C^\infty$ structure on $V/B$.

\item A {\bf stabilization} of a relative smooth Kuranishi space $K/B$ by a $G$-equivariant disk bundle
\beqn
D = D^{\rm hor} \oplus D^{\rm vert} \subset N^{\rm hor} \oplus N^{\rm vert}
\eeqn
is the relative Kuranishi space
\beqn
{\rm Stab}_D(K/B) = (G, \tilde V/ \tilde B, \tilde E, \tilde S)
\eeqn
where $\tilde V / \tilde B$ is the stabilization of the relative space $V/B$ by $D$ (which carries a canonically induced relative smooth structure), $\tilde E$ and $\tilde S$ are the natural stabilization. 

\item A {\bf strict embedding} of relative smooth Kuranishi spaces from $K_1/B_1 = (G_1, V_1/B_1, E_1, S_1)$ to $K_2/B_2 = (G_2, V_2/B_2, E_2, S_2)$ consists of a strict embedding 
\beqn
\iota_{21}: K_1 \to K_2
\eeqn
(which includes a group embedding $\iota_{21}^G: G_1 \to G_2$) and an $\iota_{21}^G$-equivariant smooth embedding $B_1 \to B_2$.

\item An {\bf embedding} of relative Kuranishi spaces from $K_1/B_1$ to $K_2/B_2$ is a conjugacy class of strict embeddings. 

\item Let $\uds{\bf S^{\rm rel} Kur}$ be the category of relatively smooth Kuranishi spaces whose morphisms are embeddings. One can also consider the rigidified and collared versions of this category.
\end{enumerate}
\end{defn}

Notice that there is an obvious functor
\beqn
\uds{\bf S^{\rm rel} Kur} \to \uds{\bf Kur}
\eeqn
by sending
\beqn
(G, V/B, E, S) \mapsto (G, V, E, S).
\eeqn

\begin{prop}[Relative smooth structures on AMS lifts]\hfill
\begin{enumerate}

\item In Situation \ref{situationf1}, any topological AMS lift of the Floer flow category is canonically lifted to $\outer \uds {\bf S^{\rm rel} Kur}_{\rm rig}$.

\item In Situation \ref{situationm1}, let $\hat{\mb X}$ be the AMS lift of the outercollaring of ${\mb X}$ to $\outer \uds{\bf Kur}_{\rm rig}$ as a multimodule over the topological AMS lifts $(\hat{\mb F}_1, \ldots, \hat{\mb F}_m; \hat{\mb F}')$ enriched in $\outer \uds{\bf Kur}_{\rm rig}$. Then $\hat {\mb X}$ is canonically lifted to $\outer \uds{\bf S^{\rm rel} Kur}_{\rm rig}$.

\item In Situation \ref{situationh1}, let $\hat {\mb H}$ be a topological AMS lift of the outercollaring of ${\mb H}$ to $\outer \uds{\bf Kur}_{\rm rig}$. Then $\hat {\mb H}$ is canonically enriched in $\outer \uds{\bf S^{\rm rel} Kur}_{\rm rig}$.

\end{enumerate}
\end{prop}

\begin{proof}
We only discuss the case about the flow category. The other cases can be proved in the same manner. For each pair $p<q$ of objects in ${\mb F}$, the fiberwise smooth structure of $\outer V_{pq} \to \outer B_{pq}$ follows from the implicit function theorem. The fact that these fiberwise structures vary in the rel-$C^\infty$ sense can be proved using the method of \cite{Swaminathan_relative} (which treats the case of Gromov--Witten theory). Alternatively, it follows from the proof of \cite[Proposition 6.2]{Bai_Xu_Arnold}.
\end{proof}

We summarize the construction so far in the following situations.

\begin{situationf}\label{situationf3}
In Situation \ref{situationf2}, the domain flow category $\outer \mb{Dom}^{\Omega, {\mb F}}$ induces an associated flow category $\outer \mb{Map}_{\mb F}^{\mb{fr}}$. Choose a transverse thickening datum ${\mb W}^{\mb F}$ and a group reduction, which gives a relatively smooth AMS lift $\hat{\mb F}$ enriched in $\outer \uds{\bf S^{\rm rel} Kur}_{\rm rig}$.
\end{situationf}

\begin{situationm}\label{situationm3}
In Situation \ref{situationm2}, $\hat{\mb X}$ is a relatively smooth AMS lift subject to an outercollaring width and a compatible list of integral actions on the involved flow categories. $\hat {\mb X}$ is a multimodule over $(\hat {\mb F}_1, \ldots, \hat {\mb F}_m; \hat {\mb F}')$ where each member of this tuple is a relative smooth AMS lift of an involved Floer flow category subject to the outercollaring width and the chosen integral action. 
\end{situationm}

\begin{situationh}\label{situationh3}
In Situation \ref{situationh2}, $\hat {\mb H}$ is a relatively smooth AMS lift of the homotopy ${\mb H}$ enriched in $\outer \uds{\bf S^{\rm rel} Kur}_{\rm rig}$ subject to an outercollaring width and a compatible list of integral actions on the involved Floer flow categories based on the homotopy $\outer \mb{Dom}^{\Omega, {\mb H}}$ enriched in $\outer \uds{\bf Curve}_{\rm rig}^{\mb C}$. The two boundaries of $\hat {\mb H}$ are relatively smooth AMS lifts of multimodules $\mb{X}_0$ resp. $\mb{X}_1$. 
\end{situationh}

\section{Vertical Normal Complex Structures and orientations}\label{section_NC_structure}

\subsection{Main statements about normal complex structures}

Recall that one has defined the category of relatively smooth Kuranishi spaces $\uds{\bf S^{\rm rel}Kur}$. Given any object $K/B = (G, V/B, E, S)$, one has a well-defined vertical tangent bundle
\beqn
T^{\rm vt} V \to V
\eeqn
which is a $G$-equivariant relatively smooth vector bundle. Recall the definition of equivariant normal complex structures (Definition \ref{defn_equivariant_NC_structure}). Given a relatively smooth Kuranishi space $K/B = (G, V/B, E, S)$, a {\bf relative NC structure} on $K/B$ is by definition a $G$-equivariant NC structure on $T^{\rm vt} V \to V$. One can then define the category of relative NC Kuranishi spaces, denoted by 
\beqn
\uds{\bf S^{\rm rel} Kur}^{\rm NC}.
\eeqn
We can also consider the rigidified and collared versions. We omit the detailed definitions.

Now we state our main theorems about the existence of fiberwise stable complex structures of Floer flow categories, multimodules, and homotopies. 

\begin{thm}\label{thm171}
In Situation \ref{situationf3}, the relatively smooth AMS lift $\hat {\mb F}$ admits a set of lifts, denoted temporarily by $\hhat {\mb F}$, to $\outer \uds{\bf S^{\rm rel} Kur}_{\rm rig}^{\rm NC}$. We call such a lift a {\bf relatively NC AMS lift}. 
\end{thm}

\begin{proof}
See Subsection \ref{subsection173}.
\end{proof}
\begin{rem}
In the following theorem, what we call a relative NC AMS lifts of $\hat {\mb F}$ actually involves more data than the lifts to $\outer \uds{\bf S^{\rm rel} Kur}_{\rm rig}^{\rm NC}$; they also include certain choices made for the construction. The construction of the relative NC AMS lifts for multimodules actually depends on these choices, not only on the resulting lifts of the flow categories.
\end{rem}

\begin{thm}\label{thm173}
In Situation \ref{situationm3}, suppose $\hhat {\mb F}_1, \ldots, \hhat {\mb F}_m; \hhat {\mb F}'$ are relatively NC AMS lifts of them provided by Theorem \ref{thm171}. Then $\hat {\mb X}$ admits a set of lifts, denoted temporarily by $\hhat {\mb X}$, to $\outer \uds{\bf S^{\rm rel} Kur}_{\rm rig}^{\rm NC}$ as a multimodule over $(\hhat {\mb F}_1, \ldots, \hhat {\mb F}_m; \hhat {\mb F}')$. Each of such $\hhat{\mb X}$ is called a {\bf relatively NC AMS lifts} (of ${\mb X}$ or $\hat {\mb X}$).
\end{thm}

\begin{proof}
See Subsection \ref{subsection174}.
\end{proof}

\begin{thm}
In Situation \ref{situationh3}, suppose $\hhat {\mb F}_1, \ldots, \hhat {\mb F}_m, \hhat {\mb F}'$ are relatively NC AMS lifts of the involved Floer flow categories, $\hhat {\mb X}_0$ and $\hhat {\mb X}_1$ are relatively NC AMS lifts of the two multimodules subject to the fixed outercollaring width and the fixed list of compatible integral actions. Then the relatively smooth AMS lift $\hat{\mb H}$ admits a lift, denoted by $\hhat {\mb H}$ temporarily, to $\outer \uds{\bf S^{\rm rel} Kur}_{\rm rig}^{\rm NC}$ as a homotopy from $\hhat {\mb X}_0$ to $\hhat {\mb X}_1$. We call $\hhat {\mb H}$ a relatively NC AMS lift (of ${\mb H}$ or $\hat {\mb H}$).
\end{thm}

\begin{proof}
A tedious reformulation of the proofs of Theorem \ref{thm171} and Theorem \ref{thm173}. There are no essential new ingredients. 
\end{proof}

\subsection{Abouzaid--Blumberg's construction for a single moduli space}

We first review the construction of \cite[Section 11]{Abouzaid_Blumberg} in the case of a single moduli space from the Floer flow category. Let ${\mb F}$ be the Floer flow category associated to the Hamiltonian $H$. Consider a pair of objects $p<q$ such that there is no intermediate object $r$ with $p<r<q$. Let $K_{pq} = (G_{pq}, V_{pq}, E_{pq}, S_{pq})$ be an AMS global chart of $M_{pq}^{\mb F}$, viewed as an object of $\uds{\bf S^{\rm rel} Kur}$. We can ignore $E_{pq}$ and $S_{pq}$ for the discussion but should remember the reletive structure $V_{pq}/B_{pq}$. We will also drop the outercollaring for a moment. 

\subsubsection{ A complex-linear Cauchy--Riemann operator and an interpolation}

We first specify a family of complex linear Cauchy--Riemann operator for points in $V_{pq}$. Fix an almost complex structure $J_0$ and a $J_0$-linear connection $\nabla^{TX}$ such that its holonomy along each 1-periodic orbit of $H_t$ is trivial. Choose $\tau > 2$. 

We first look at an individual smooth map $u: {\mb R} \times S^1 \to X$ which converges exponentially to 1-periodic orbits of $H$. Define the space
\beqn
\wt W^{1,\tau}(u^* TX) \subset W^{1,\tau}_{\rm loc}(u^* TX)
\eeqn
to be the space of sections which differ from $\nabla^{TX}$-parallel sections along the limiting periodic orbits by $W^{1,\tau}$-small terms near infinity. Then the $J_0$-linear connection $\nabla^{TX}$ induces a complex-linear Cauchy--Riemann operator
\beqn
D_u^{\mb C}: \wt W^{1,\tau}(u^*TX) \to L^\tau (\Lambda^{0,1} \otimes u^* TX).
\eeqn
As $u$ converges exponentially to periodic orbits, it is easy to see that $D_u^{\mb C}$ is Fredholm.

The operator $D_u^{\mb C}$ and the linearization of the Floer equation can be connected by a particular 1-parameter family of Fredholm operators modulo certain finite-dimensional corrections. Choose a cut-off function 
\beqn
\chi: {\rm Dom}(u) \cong {\mb R} \times S^1 \to [0, 1]
\eeqn
which only depends on the ${\mb R}$-variable such that $\chi(s) = 0$ for $s \ll 0$ and $\chi(s) = 1$ for $s \gg 0$. For each object $p$, let $u_p: {\mb R} \times S^1 \to X$ be the corresponding constant solution to the Floer equation (which does not depend on the capping). Then there is an associated linearized operator
\beqn
D_{u_p}(\xi) = \partial_s \xi + J_t \left( \nabla_t \xi - \nabla_\xi X_{H_t}(p) \right).
\eeqn

Let 
\beqn
\wt W^{1,\tau}({\mb R}\times S^1, u_p^* TX) \subset W^{1,\tau}_{\rm loc}({\mb R}\times S^1, u_p^* TX)
\eeqn
be the subspace of $W^{1,\tau}_{\rm loc}$-sections which satisfies the same asymptotic condition near $-\infty$ as $\wt W^{1,\tau}(u^* TX)$ and which is $W^{1,\tau}$-small near $+\infty$. Then consider the interpolation
\beq\label{operator_dp}
\begin{split}
D_p: \wt W^{1,\tau} ({\mb R}\times S^1, u_p^* TX) \to &\  L^\tau ({\mb R}\times S^1, \Lambda^{0,1}\otimes u_p^* TX)\\
\xi\mapsto &\  ( 1- \chi) (\nabla^{TX} \xi)_{J_0}^{0,1} + \chi D_{u_p} ( \xi).
\end{split}
\eeq
Then $D_p$ is a real-linear Fredholm operator. Notice that from the definition, $D_p$ only depends on the underlying 1-periodic orbit but not the capping.

We now choose, for each object $p$ (actually, as discussed above, only depending on the underlying orbit), a finite-dimensional complex vector space $R_p^+$, a linear map
\beq\label{eqn_nup}
\nu_p: R_p^+ \to C_0^\infty( {\mb R}\times S^1, \Lambda^{0,1} \otimes u_p^* TX)
\eeq
such that the map
\beqn
\nu_p \oplus D_p: R_p^+ \oplus \wt W^{1,\tau} ({\mb R}\times S^1, u^* TX) \to L^\tau ({\mb R}\times S^1, \Lambda^{0,1}\otimes u_p^* TX)
\eeqn
is surjective. Denote
\beqn
R_p^-:= {\rm Ker}( \theta_p \oplus D_p) \subset R_p^+ \oplus \wt W^{1,\tau}({\mb R}\times S^1, u^* TX).
\eeqn
Notice that the pair $(R_p^-, R_p^+)$ does not depend on the Sobolev exponent $\tau$. 

Now we would like to construct interpolations between $D_u$ and $D_u^{\mb C}$, up to factors coming from $D_p$ and $D_q$, in a 1-parameter family. The basic idea has been described in \cite[Section 11]{Abouzaid_Blumberg} and \cite[Appendix B]{Abouzaid_Blumberg_2024}. However, we need to be more careful when generalizing to the case of multimodules. 

More precisely, we would like to construct an interpolation  between 
\beqn
D_p \oplus D_u: \wt W^{1,\tau }(u_p^* TX) \oplus W^{1,\tau }(u^* TX) \to L^\tau (\Lambda^{0,1} \otimes u_p^* TX) \oplus L^\tau ( \Lambda^{0,1} \otimes u^* TX)
\eeqn
and 
\beqn
D_u^{\mb C} \wt \oplus D_q: \wt W^{1, \tau}(u^* TX) \wt \oplus \wt W^{1,\tau }(u_q^* TX) \to L^\tau (\Lambda^{0,1}\otimes u_q^* TX).
\eeqn
Here $\wt \oplus$ means the subset of the direct sum satisfying the matching condition. 

To build the interpolation, consider the family of cut-off functions
\beqn
\chi_\rho: {\mb R} \to [0, 1],\ \chi_\rho (s) = \chi( s- \rho).
\eeqn
Then for all $t \in {\mb R}$, define
\beq\label{family_CR_operator}
\begin{split}
D_{u,\rho}: \wt W^{1,\tau} (u^* TX) \to &\ L^\tau (\Lambda^{0,1}\otimes u^* TX),\\
\xi \mapsto &\ (1 - \chi_\rho) (\nabla^{TX} \xi)_{J_0}^{0,1} + \chi_\rho D_u (\xi).
\end{split}
\eeq
Notice that even if $u$ has bubble components, these operators are still defined. Then intuitively, we have
\begin{align*}
&\ \lim_{\rho \to -\infty} D_{u, \rho} \approx D_p \oplus D_u,\ &\ \lim_{\rho \to +\infty} D_{u, \rho } \approx D_u^{\mb C} \wt \oplus D_q.
\end{align*}
This is the underlying idea of interpolating the linearized operator $D_u$ and a complex linear operator $D_u^{\mb C}$, up to the two operators $D_p$ and $D_q$ which are independent of $u$.

To better formulate this construction over the whole Kuranishi chart $V_{pq}$, we need to introduce an auxiliary space. Denote
\beqn
\wt B_{pq}:= B_{pq} \times [-\infty, +\infty],
\eeqn
viewed as a product of stratified spaces where $[-\infty, +\infty]$ has one open stratum and two boundary strata. We define an auxiliary universal curve $\wt C_{pq} \to \wt B_{pq}$ whose fiber $\wt C_{\phi, \rho}$ at $(\phi, \rho)$ with $\phi \in B_{pq}$ and $\rho \in (-\infty, +\infty)$ is still $C_\phi$, while for $\rho = -\infty$ resp. $\rho = +\infty$, $\wt C_{\phi, \rho}$ has the additional component of the domain of the map $u_p$ resp. the domain of $u_q$. One can then define the additional space $\wt V_{pq}$ which is identical to $V_{pq} \times [-\infty, +\infty]$. However, for each point $(\wt\phi, \wt u, e, F) \in \wt V_{pq}$, we regard $u$ as a map $\wt u: \wt C_{\phi, \rho} \to X$ where if $\rho = \pm\infty$, $\wt u$ maps the extra cylindrical component to the constant loop at the negative or the positive end. Then $\wt V_{pq} \to \wt B_{pq}$ is still relatively smooth. The operator $D_{u,  \rho}$ can be viewed as a Cauchy--Riemann operator $D_{\wt u}$ acting on the bundle $\wt u^* TX$ (and extended to suitable Sobolev completions). We remark that, this setup facilitates the gluing construction: when $\rho \in \pm\infty$, if the operator $D_{u, \rho}$ is surjective, then each vector in the kernel can be glued to an element in the kernel at a nearby point.  

With the above viewpoint understood, using basic gluing construction, one obtains the following technical lemma. Let 
\beqn
\mathring {\mc E}_{pq}^0 \to \wt V_{pq}
\eeqn
be the infinite-dimensional vector bundle whose fiber over $\wt x = (\wt \phi, \wt u, e, F)$ is the space
\beqn
\mathring {\mc E}_{\wt x}^0 = \Omega_c^{0,1}( \wt C_{\phi, \rho}, \wt u^* TX).
\eeqn
Notice that when $\rho = \pm \infty$, one has the splitting
\beq\label{eqn172}
\mathring {\mc E}_{\wt x}^0 = \Omega_c^{0,1}( {\mb R}\times S^1, u_p^* TX) \oplus \Omega_c^{0,1}( C_\phi, u^* TX)
\eeq
or
\beq\label{eqn173}
\mathring {\mc E}_{\wt x}^0 = \Omega_c^{0,1}( C_\phi, u^* TX) \oplus \Omega_c^{0,1}({\mb R}\times S^1, u_q^* TX).
\eeq

\begin{lemma}\label{lemma175}
There exist a unitary representation $W_{pq}$ of $G_{pq}$ and a continuous $G_{pq}$-equivariant bundle map over $\wt V_{pq} \times X$
\beqn
\wt\nu_{pq}: \uds R_p^+ \oplus \uds W_{pq} \oplus \uds R_q^+ 
\to \mathring {\mc E}_{pq}^0
\eeqn
satisfying the following conditions. Let $\wt x = (\wt \phi, \wt u, e, F) \in \wt V_{pq}$ where $\wt \phi = (\phi, \rho ) \in \wt B_{pq}$.
\begin{enumerate}

    \item When $\rho = -\infty$, with respect to the splitting \eqref{eqn172}, $\wt\nu_{pq}$ has the block form
    \beqn
    \wt \nu_{pq} = \left[ \begin{array}{ccc} \nu_p &  0 & 0 \\
                     0 & 0 & 0 \end{array}\right]
    \eeqn
    where $\nu_p$ is specified in \eqref{eqn_nup}.
    
    \item When $\rho = +\infty$, with respect the splitting \eqref{eqn173}, $\wt\nu_{pq}$ has the block form
    \beqn
    \wt\nu_{pq} = \left[ \begin{array}{ccc} 0 &  \nu_{\wt x} & 0 \\
     0 & 0 &  \nu_q \end{array}\right]
     \eeqn
     where
    \beqn
    \nu_{\wt x}: W_{pq} \to \Omega_c^{0,1}(C_\phi, u^* TX)
    \eeqn
    is a complex-linear linear map.
    
    \item For each $\wt x \in \wt V_{pq}$, the image of $\wt\nu_{pq}$ at $\wt x$ is transverse to the image of the operator $D_{\wt u}$. 
\end{enumerate}
\end{lemma}

\begin{proof}
Basic Fredholm theory and a portion of gluing argument.
\end{proof}

\subsubsection{An extended index bundle}

Now assume $W_{pq}$ and $\wt\nu_{pq}$ is given as in Lemma \ref{lemma175}, consider the family of linear maps parametrized by $\wt x = (\wt \phi, \wt u, e, F) \in \wt V_{pq}$:
\beqn
\nu_{\wt x} \oplus D_{\wt u}: \big( R_p^+ \oplus W_{pq} \oplus R_q^+ \big) \oplus \wt W^{1,\tau}(\wt u^* TX) \to L^ \tau (\Lambda^{0,1}\otimes \wt u^* TX).
\eeqn

\begin{prop}\label{prop1711}
The union of kernels of $\wt \nu_{\wt x} \oplus D_{\wt u}$ form a continuous \footnote{In fact it is a relative smooth vector bundle.} vector bundle 
\beqn
\wt I_{pq}  \to \wt V_{pq}.
\eeqn
\end{prop}

\begin{proof}
The proof is based on the typical gluing construction. However, one needs to maintain the linear structure in the gluing procedure. The pregluing is constructed using parallel transport on the target manifold $X$ and cut-paste on the domain, which are linear in the input. The use of the implicit function theorem can also respect the linear feature.
\end{proof}

Hence the vector bundles on the two boundary pieces, $\partial^\pm \wt V_{pq}:= V_{pq}\times \{\pm\infty\} \cong V_{pq}$ are (non-canonically) isomorphic. The family of kernels selects a homotopy class of isomorphisms. By the specific characterization of $\wt\nu_{pq}$ given in Lemma \ref{lemma175}, one can see that 
\beqn
\wt I_{pq} |_{\partial^- \wt V_{pq}} \cong \uds R_p^- \oplus T^{\rm vt} V_{pq} \oplus \uds W_{pq} \oplus \uds R_q^+;
\eeqn
and
\beqn
\wt I_{pq} |_{\partial^+ \wt V_{pq}} \cong \uds R_p^+ \oplus I_{pq}^{\mb C} \oplus \uds R_q^-.
\eeqn
Here $I_{pq}^{\mb C}$ is formed by the kernels of the complex-linear  operator 
\beqn
\nu_{\wt x} \oplus D_u^{\mb C}: W_{pq} \oplus \wt W^{1,\tau}(C_\phi, u^* TX) \to L^\tau( C_\phi, u^* TX)
\eeqn
hence is a complex vector bundle.

We can then choose a continuous bundle isomorphism
\beqn
\uds R_p^- \oplus T^{\rm vt} V_{pq} \oplus \uds W_{pq} \oplus \uds R_q^+ \cong \uds R_p^+ \oplus I_{pq}^{\mb C} \oplus \uds R_q^-.
\eeqn
Notice that $G_{pq}$ acts trivially on the vector spaces $R_p^\pm$ and $R_q^\pm$. As $W_{pq}$ is complex, one can see that this is a stable isomorphism from $T^{\rm vt} V_{pq}$ to $I_{pq}^{\mb C}$ and hence induces an NC structure on $T^{\rm vt} V_{pq}$.

\subsection{Proof of Theorem \ref{thm171}}\label{subsection173}

\subsubsection{An extended Kuranishi flow category}

We extend the previous construction to all moduli spaces involved in a single Floer flow category ${\mb F}$. To construct the interpolation in a compatible way for all moduli spaces, we utilize the lateral lines (see Definition \ref{lateral_line}). Recall that for each pair of objects $p\leq q$, in the AMS construction there is a $G_{pq}$-manifold $B_{pq}$ parametrizing certain curves in $\mb{CP}^{d_{pq}}$. Consider the moduli space
\beqn
\ov{\mc M}{}_{0, 2, +}^{\mb R}(\mb{CP}^{d_{pq}}, d_{pq})
\eeqn
of 1-marked stable cylinders with the interior marked points lying on the lateral lines of domains. There is the forgetful map with ${\mb R}$ as generic fibers
\beqn
\ov{\mc M}{}_{0,2,+}^{\mb R}(\mb{CP}^{d_{pq}}, d_{pq}) \to \ov{\mc M}{}_{0,2}^{\mb R}(\mb{CP}^{d_{pq}}, d_{pq}).
\eeqn
Let $\wt B_{pq}$ be the preimage of $B_{pq}$; let $\wt \phi$ denote a typical element of $\wt B_{pq}$ with the underlying point $\phi \in B_{pq}$. We use the convention that this interior marked point does not go to sphere bubbles but only stay in cylindrical components; namely, it can coincide with a node. Let 
\beqn
\wt C_{pq} \to \wt B_{pq}
\eeqn
be the pullback of $C_{pq} \to B_{pq}$,  which still gives an equivariant family of curves
\beqn
\wt {\mc C}_{pq}:= ({\mc G}_{pq}, \wt B_{pq}, \wt C_{pq}).
\eeqn
Let $C_{\wt \phi} \subset \wt C_{pq}$ be the fiber over $\wt \phi \in \wt B_{pq}$. Then there is a natural map 
\beqn
C_{\wt \phi} \to C_\phi
\eeqn
which collapses at most one constant cylindrical component. Notice that for each $\wt \phi \in \wt B_{pq}$, there is a special cylindrical component
\beqn
C_{\wt \phi}^\star \subset C_{\wt \phi}
\eeqn
which contains the marked point.

Notice that there two disjoint closed subsets 
\begin{align*}
&\ \partial^- \wt B_{pq},\ &\ \partial^+ \wt B_{pq}
\end{align*}
corresponding to configurations where the marked points is at the very left resp. very right of the curve. Notice that one has canonical identifications
\beqn
\partial^- \wt B_{pq} \cong B_{pq} \cong \partial^+ \wt B_{pq}.
\eeqn

Notice that $\wt B_{pq}$ and $\wt C_{pq}$ can also be organized as morphism spaces of a flow category. However, as they are objects in a special category which does not have an obvious symmetric monoidal structure, therefore, we will describe the flow category structure in a more {\it ad hoc} way, without specifying an underlying category of stratified objects. Define a ``product'' $\wt B_{pr} \boxtimes \wt B_{rq}$ via the diagram
\beqn
\vcenter{ \xymatrix{ \wt B_{pr} \boxtimes \wt B_{rq} \ar[r] \ar[d]  &   \wt B_{pq} \ar[d] \\
          B_{pr}\times B_{rq} \ar[r]_-{\zeta_{prq}^B}   &  B_{pq} } },
\eeqn
which we require to be a pullback. The ``product" universal curves also be defined as a pullback, which admits a domain map
\beqn
\wt C_{pr} \boxtimes \wt C_{rq} \to \wt B_{pr} \boxtimes \wt B_{rq}.
\eeqn

Now suppose we are in Situation \ref{situationf3} where $\hat {\mb F}$ is a relativley smooth AMS lift of the Floer flow category ${\mb F}$ enriched in $\outer \uds{\bf S^{\rm rel} Kur}_{\rm rig}$. For each pair of $p \leq q$, let 
\beqn
\wt V_{pq} \to \wt B_{pq}
\eeqn
be the pullback of $V_{pq} \to B_{pq}$. Then one can define similarly $\wt V_{pr} \boxtimes \wt V_{rq}$ and corresponding maps
\beqn
\wt V_{pr} \boxtimes \wt V_{rq} \to \partial^{prq} \wt V_{pq}.
\eeqn
Denote
\beqn
\wt V_{pr_1 \cdots r_l q} :=  \wt V_{pr_1} \boxtimes \cdots \boxtimes \wt V_{r_l q},
\eeqn
which is well-defined due to the associativity assumption of compositions of the flow category. Notice that there is also a corresponding stratum 
\beqn
\partial^{pr_1 \cdots r_l q} \wt V_{pq} \subset \wt V_{pq}
\eeqn
which contains $\wt V_{pr_1 \cdots r_l q}$.

\subsubsection{The construction of the normal complex structure}\hfill

---{\bf The family of Cauchy--Riemann operators}---Then we would like to build a family of Fredholm operators indexed by $\wt x = (\wt \phi, \wt u, e, F) \in \outer \wt V_{pq}$. Fix $\tau > 2$. The construction appears to depend on $\tau$ but the resulting normal complex structure does not. Notice that for each $\wt x \in \wt V_{pq}$, there are a pair of Banach spaces
\begin{align*}
&\ {\mc E}_{\wt x}^1:= \wt W^{1,\tau}(C_{\wt \phi}, \wt u^* TX),\ &\ {\mc E}_{\wt x}^0:= L^{\tau}(C_{\wt \phi}, \Lambda^{0,1} \otimes \wt u^* TX).
\end{align*}
We regard them as ``vector bundles'' over $\wt V_{pq}$ only in a set-theoretic sense, denoted by 
\begin{align*}
&\ {\mc E}_{pq}^1\to \wt V_{pq},\ &\ {\mc E}_{pq}^0 \to \wt V_{pq}.
\end{align*}
However, we keep in mind that the induced index-theoretic constructions over the finite-dimensional moduli spaces is nevertheless well-defined as topological objects, and we can perform gluing constructions over the nodal locus to make sense of them.

We give an {\it ad hoc} construction of a family of Cauchy--Riemann operators between these spaces. First we construct a family of cut-off functions.

\begin{lemma}\label{lemma_cut_off_function}
There exist a collection of $G_{pq}$-invariant cut-off functions
\beqn
\chi_{pq}: \outer \wt C_{pq} \to [0, 1]
\eeqn
for all pairs of objects $p \leq q$ of ${\mb F}$ satisfying the following conditions.
\begin{enumerate}

    \item $\chi_{pq} \equiv 1$ on components to the right of $C_{\wt \phi}^\star$ and $\wt \chi_{pq} \equiv 0$ on components to the left of $C_{\wt \phi}^\star$.

    \item On $C_{\wt \phi}^\star$, $\chi_{pq} = 0$ for $s \ll 0$ and $\wt \chi_{pq} = 1$ for $s \gg 0$ where $s$ is a cylindrical coordinate.

    \item $\chi_{pq}$ is collared.

    \item $\chi_{pq}$ respects the flow category structure. Namely, whenever $p< r < q$ and $\wt \phi \in \outer \wt C_{pr} \boxtimes \outer \wt C_{rq}$, if the marked point is on the $pr$-section then on this part $\chi_{pq}$ coincides with $\chi_{pr}$; if the marked point is on the $rq$-section then on this section $\chi_{pq}$ coincides with $\chi_{rq}$.

    \item $\chi_{pq}$ respect the rigidification. $\chi_{pq}$ does not vary with the normal direction near $\outer \wt C_{pr}\boxtimes \outer \wt C_{rq}$ inside $\partial^{prq} \wt C_{pq}$. 

\end{enumerate}
\end{lemma}

\begin{proof}
An inductive construction similar to the proof of Theorem \ref{thma_FOP}, which is standard at this point.
\end{proof}

Now for each $\wt x = (\wt \phi, \wt u, e, F) \in \wt V_{pq}$, consider the Cauchy--Riemann operator
\beqn
D_{\wt x}: {\mc E}_{\wt x}^1 \to {\mc E}_{\wt x}^0
\eeqn
given by the same formula as \eqref{family_CR_operator}. We would like to consider the index bundles. Let
\beqn
\mathring {\mc E}_{\wt x}^0 \subset {\mc E}_{\wt x}^0
\eeqn
be the subspace of compactly supported smooth sections and let 
\beqn
\mathring {\mc E}_{pq}^0 \to \wt V_{pq}
\eeqn
denote the corresponding ``subbundle.''

\begin{defn}\label{defn_bundle_interpolator}
An {\bf index bundle interpolator} over the Kuranishi flow category $\hat {\mb F}$ consists of the following items.
\begin{enumerate}

\item For each object $p \in {\rm Ob}{\mb F}$, a finite-dimensional complex vector space $R_p^+$ and a complex-linear map
\beqn
\nu_p: R_p^+ \to \Omega_c^{0,1} ({\mb R}\times S^1, u_p^* TX)
\eeqn
which is transverse to the image of $D_p$. We require that $R_p$ and $\nu_p$ only depends on the $\Pi$-orbit, i.e., the underlying uncapped orbit $\uds p$. Denote
\beqn
R_p^- = {\rm ker} (\nu_p \oplus D_p) \subset R_p^+ \oplus  \wt W^{1,\tau}({\mb R}\times S^1, u_p^* TX),
\eeqn
which is only a real vector space.

\item For each $p \leq q$, a finite-dimensional Hermitian representations $W_{pq}$ of $G_{pq}$ and a $G_{pq}$-equivariant complex-linear bundle maps over $\wt V_{pq}$
\beqn
\wt\nu_{pq}: \uds R_p^+ \oplus \uds W_{pq} \oplus \uds R_q^+ \to \mathring {\mc E}_{pq}^0.
\eeqn

\item For each triple $p < r < q$, a complex-linear map
\beqn
\wt\iota_{prq}: W_{pr} \oplus R_r^+ \oplus W_{rq} \to W_{pq}
\eeqn
which is equivariant with respect to the group homomorphism $G_{pr}\times G_{rq} \to G_{pq}$.
\end{enumerate}
They are required to satisfy the following conditions.
\begin{enumerate}
    \item {\bf (The original index bundle)} For each $\wt x  = (\wt \phi, \wt u, e, F) \in \partial^- \wt V_{pq}$, with respect to the natural decomposition
    \beqn
    \mathring {\mc E}_{\wt x}^0 \cong \Omega^{0,1}_c({\mb R}\times S^1, u_p^* TX) \oplus \Omega^{0,1}_c(C_\phi, u^* TX)
    \eeqn
    we have
    \beqn
    \wt\nu_{pq} = \left[ \begin{array}{ccc}  \nu_p & 0 & 0  \\  0 &  0  & 0 \end{array}\right]
    \eeqn

    \item {\bf (The complex index bundle)} For each $\wt x  = (\wt \phi, \wt u, e, F) \in \partial^+ \wt B_{pq}$, with respect to the natural decomposition
    \beqn
    \mathring {\mc E}_{\wt x}^0 \cong 
    \Omega_c^{0,1}( C_\phi, u^* TX) \oplus \Omega_c^{0,1}({\mb R}\times S^1, u_q^* TX)
    \eeqn
    we have
    \beqn
    \wt\nu_{pq} = \left[ \begin{array}{ccc} 0  &  \nu_{\wt x} &  0 \\
          0 & 0 & \nu_q \end{array}\right].
    \eeqn
    Here $\nu_{\wt x}: W_{pq} \to \Omega_c^{0,1}(C_\phi, u^* TX)$ is a complex-linear map.

    \item {\bf (Flow category structure I)} The above conditions imply that $\wt\nu_{pr}$ and $\wt\nu_{rq}$ determine a bundle map 
    \beqn
    \wt\nu_{prq}: (\uds R_p^+ \oplus \uds W_{pr} \oplus \uds R_r^+ \oplus \uds W_{rq} \oplus \uds R_q^+ ) |_{\wt V_{prq}} \to  \mathring {\mc E}_{pr}^0 \boxplus \mathring {\mc E}_{rq}^0.
    \eeqn
    We require that the following diagram commutes.
    \beqn
    \xymatrix{ \big( \uds R_p^+ \oplus \uds W_{pr} \oplus \uds R_r^+ \oplus \uds W_{rq} \oplus \uds R_q^+ \big)|_{\wt V_{prq}} \ar[rr]^-{\wt\nu_{prq}} \ar[dd]_-{{\rm Id}_{R_p^+} \oplus \wt\iota_{prq} \oplus {\rm Id}_{R_q^+}}   &  &  \mathring {\mc E}_{pr}^0 \boxplus \mathring {\mc E}_{rq}^0 \ar[dd] 
    \\ & &  \\ \uds R_p^+ \oplus \uds W_{pq} \oplus \uds R_q^+|_{ \partial^{prq} \wt V_{pq}}  \ar[rr]_-{\wt\nu_{pq}} & &  \mathring {\mc E}_{pq}^0|_{\partial^{prq} \wt V_{pq} }}
    \eeqn

    \item {\bf (Flow category structure II)} Whenever $p<r < s < q$, the following diagram commutes. 
    \beqn
    \xymatrix{   W_{pr} \oplus R_r^+ \oplus W_{rs} \oplus R_s^+ \oplus W_{sq} \ar[rrr]^-{\wt\iota_{prs} \oplus {\rm Id}_{R_s^+} \oplus {\rm Id}_{W_{sq}}} \ar[d]_{{\rm Id}_{W_{pr} \oplus {\rm Id}_{R_r^+} \oplus \wt\iota_{rsq}}}   & &   & W_{ps} \oplus R_s^+ \oplus W_{sq} \ar[d]^{\wt\iota_{psq}} \\
           W_{pr} \oplus R_r^+ \oplus W_{rq} \ar[rrr]_-{\wt\iota_{prq}} &  & & W_{pq} }
    \eeqn

    \item {\bf (Transversality)} For each $\wt x \in \wt V_{pq}$, the image $\wt \nu_{pq}$ contained in $\mathring {\mc E}_{\wt x}^0 \subset {\mc E}_{\wt x}^0$ is transverse to the image of $D_{\wt x}$ (which is a condition true for all $\tau>2$). 

    \item {\bf (Collar)} The maps $\wt\nu_{pq}$ satisfy the natural collaring condition. Namely, near each stratum $\partial^\alpha \wt V_{pq}$, $\wt\nu_{pq}$ is independent of the collar coordinates.

    \item {\bf (Rigidification)} $\wt \nu_{pq}$ does not vary in the normal direction to $\wt V_{prq}$  specified by the rigidfication.

    \item {\bf (Novikov equivariance)} These objects are equivariant with respect to the natural $\Pi$-action.
    
\end{enumerate}
\end{defn}

\begin{lemma}
There exists an index bundle interpolator on $\hat{\mb F}$.
\end{lemma}

\begin{proof}
This is an inductive construction. The spaces $R_p^+$ and the maps $\nu_p$ can be chosen independently for all objects. The construction of Abouzaid--Blumberg recalled previously provides the required space $W_{pq}$ and the bundle map $\wt\nu_{pq}$ for $p<q$ which does not have intermediate objects. Suppose for a pair of objects $p< q$, one has constructed such structures for all $r<s$ with $d_r -d_s < d_p - d_q$, then for $W_{pq}$, one can choose initially $W_{pq}$ to be the direct sum of all $W_{pr} \oplus R_p^+ \oplus W_{rq}$ for intermediate objects $r$. The value of $\wt\nu_{pq}$ on the boundary is then determined by first three conditions of Definition \ref{defn_bundle_interpolator} as well as the rigidification condition. The collar condition allows us to extend slightly into the interior of $\wt V_{pq}$. The transversality achieved on the boundary strata is carried over to such an extension. To achieve global transversality, one may need to enlarge $W_{pq}$ by taking direct sum with another complex vector space.  
\end{proof}

Once the index bundle interpolator is chosen, one automatically obtains a system of vector bundles over the spaces $\wt V_{pq}$. Indeed, for each $\wt x\in \wt V_{pq}$, denote
\beqn
\wt I_{pq}|_{\wt x}:= {\rm ker} ( \wt \nu_{pq}|_{\wt x} \oplus D_{\wt x}) \subset {\mc E}_{pq}^1|_{\wt x}.
\eeqn
They form a topological vector bundle $\wt I_{pq} \to \wt V_{pq}$. We consider its restriction to $\partial^\pm \wt V_{pq}$. By the first two conditions of Definition \ref{defn_bundle_interpolator}, one has $G_{pq}$-equivariant identifications
\beqn
\wt I_{pq}|_{\partial^- \wt V_{pq}} \cong \uds R_p^- \oplus T^{\rm vt} V_{pq} \oplus \uds W_{pq} \oplus \uds R_q^+
\eeqn
and 
\beqn
\wt I_{pq}|_{\partial^+ \wt V_{pq}} \cong \uds R_p^+ \oplus I_{pq}^{\mb C} \oplus \uds R_q^-.
\eeqn
Here $I_{pq}^{\mb C}$ has fiber over $\wt x$ being
\beqn
{\rm ker} ( \nu_{\wt x} \oplus D_u^{\mb C})
\eeqn
hence a $G_{pq}$-equivariant complex vector bundle. Notice that $\partial^- \wt V_{pq} \cong V_{pq} \cong \partial^+ \wt V_{pq}$ canonically and the inclusions of $V_{pq}$ into $\wt V_{pq}$ from both sides are $G_{pq}$-equivariant homotopy equivalences. Hence the two restriction bundles are equivariantly isomorphic.

To construct a normal complex structure on the flow category, one must choose such bundle isomorphisms in a way compatible with the flow category structure. Notice that the product $\wt V_{prq} \cong \wt V_{pr}\boxtimes \wt V_{rq}$, it is the union
\beqn
\wt V_{pr}\times \partial^- \wt V_{rq} \underset{\partial^+ \wt V_{pr}\times \partial^- \wt V_{rq}}{\cup} \partial^+ \wt V_{pr} \times \wt V_{rq}.
\eeqn
Hence $\wt I_{pr}$ and $\wt I_{rq}$ determine a vector bundle
\beqn
\wt I_{prq} \to \wt V_{prq}
\eeqn
which is $\wt I_{pr} \oplus \uds R_r^- \oplus T^{\rm vt} V_{rq} \oplus \uds W_{rq} \oplus \uds R_q^+$ on the first piece and which is $\uds R_p^+ \oplus I_{pq}^{\mb C} \oplus \wt I_{rq}$ on the second piece.

\begin{lemma}\label{lemma1713}
There exists a system of $G_{pq}$-equivariant real vector bundle isomorphisms
\beqn
\wt I_{pq}|_{\partial^- \wt V_{pq}} \cong \wt I_{pq}|_{\partial^+ \wt V_{pq}}
\eeqn
satisfying the following conditions.

\begin{enumerate}

\item Whenever $p<r < q$, the following diagram commutes.
\beqn
\xymatrix{  \wt I_{prq}|_{\partial^- \wt V_{pr}\times \partial^- \wt V_{rq}}      \ar[r] \ar[d]  &  \wt I_{prq}|_{\partial^+ \wt V_{pr} \times \partial^- \wt V_{rq}}  \ar[r] & \wt I_{prq}|_{\partial^+ \wt V_{pr} \times \partial^+ \wt V_{rq}} \ar[d]\\
\wt I_{pq}|_{\partial^- \wt V_{prq}} \ar[rr] & & \wt I_{pq}|_{\partial^+ \wt V_{prq}}} 
\eeqn

\item The isomorphisms are collared and respect the rigidification.

\end{enumerate}
\end{lemma}

\begin{proof}
Induction and homotopy triviality of the interpolating interval $[0,1]$.
\end{proof}

Now as $I_{pq}^{\mb C}$ is complex and $R_p^\pm$, $R_q^\pm$ are trivial representations of $G_{pq}$, the isomorphism induces an NC structure on $T^{\rm vt} V_{pq}$. The compatibility condition given in Lemma \ref{lemma1713} implies that the structural maps of the flow category respect the vertical NC structures. Therefore, one obtains a lift of $\hat {\mb F}$ to $\outer \uds {\bf S^{\rm rel} Kur}_{\rm rig}^{\rm NC}$. This finishes the proof of Theorem \ref{thm171}.

\subsection{Proof of Theorem \ref{thm173}}\label{subsection174}

Let $\Sigma^{\mb X}$ be the smooth Floer domain in Situation \ref{situationm3}. One still needs to relate the linearization of the Floer equation with a complex-linear operator on $\Sigma^{\mb X}$ via a certain interpolation. However, in order to be compatible with chosen interpolations for the involved flow categories, the interpolation here needs to depend on generally $m$ parameters instead of just one, where $m$ is the number of negative cylindrical ends, or more generally, the number of lateral lines in $\Sigma^{\mb X}$.

Let $\ell_1, \ldots, \ell_m \subset \Sigma^{\mb X}$ be the lateral lines. Choose diffeomorphisms
\beqn
\rho_j: \ell_j \cong {\mb R}
\eeqn
such that on cylindrical ends, $\rho_j$ coincides with the radial coordinates up to a shift. Let 
\beqn
\wt B_{p_1 \cdots p_m; p'}^{\mb X}
\eeqn
be the version of the moduil $B_{p_1 \cdots p_m; p'}^{\mb X}$ whose elements are $(\phi, z_1, \ldots, z_m)$ where $z_j$ is contained in the $j$-th lateral line $\ell_j \subset C_\phi$ and
\beqn
{\rm either}\ \rho_j(z_j) \leq 0\ \forall j\ \ {\rm or}\ \rho_1(z_1) = \cdots = \rho_m (z_m) \geq 0.
\eeqn
Compactify the spaces to allow the marked points to escape from cylindrical ends: notice that they can either escape from negative ends independently, or escape from the positive end simultaneously. Then the foregetful map
\beqn
\wt B_{p_1 \cdots p_m; p'}^{\mb X} \to B_{p_1 \cdots p_m; p'}
\eeqn
has generic fibers being homeomorphic to a $m$-dimensional complex. Notice that the functions $\rho_j$ extend to a continuous and $G_{p_1\cdots p_m; p'}^{\mb X}$-invariant map
\beqn
\rho: \wt B_{p_1 \cdots p_m; p'}^{\mb X} \to [-\infty, +\infty]^m
\eeqn
such that either all $\rho_j$ are negative or all $\rho_j$ are nonnegative and identical.

Moreover, one can pullback the universal curve to 
\beqn
\wt C_{p_1 \cdots p_m; p'}^{\mb X} \to \wt B_{p_1 \cdots p_m; p'}^{\mb X}.
\eeqn

Similar to before, one has ``structural maps''
\begin{align*}
&\ \wt B_{p_i q_i}^{{\mb F}_i} \boxtimes \wt B_{p_1 \cdots p_{i-1} q_i p_{i+1}\cdots p_m; p'}^{\mb X} \to \wt B_{p_1 \cdots p_m; p'}^{\mb X},\ &\ \wt B_{p_1 \cdots p_m; q'}^{\mb X}\boxtimes \wt B_{q' p'}^{{\mb F}'} \to \wt B_{p_1 \cdots p_m; p'}^{\mb X}
\end{align*}
which satisfy the usual compatibility condition.

We choose compatible cut-off functions, whose existence follow again from an induction argument.
\begin{lemma}
Suppose $\chi_{p_j q_j}^{\mb{F}_j}$ and $\chi_{p' q'}^{{\mb F}'}$ have been chosen independently as given by Lemma \ref{lemma_cut_off_function}. Then there exists a family of $G_{p_1 \cdots p_m; p'}^{\mb X}$-invariant continuous functions
\beqn
\chi_{p_1 \cdots p_m; p'}^{{\mb X}}: \wt C_{p_1 \cdots p_m; p'}^{{\mb X}} \to [0, 1]
\eeqn
satisfying the following conditions. We characterize the property by considering its restriction to the fiber at each point $\wt\phi \in \wt B_{p_1 \cdots p_m; p'}^{\mb X}$. Denote the restriction by $\chi_{\wt\phi}: C_{\wt\phi} \to [0, 1]$. 

\begin{enumerate}
    \item When $\rho_1, \ldots, \rho_m \leq 0$, the support of $\chi_{\wt\phi}$ is the disjoint union of $m$ connected cylindrical regions contained in the $m$ negative ends of $C_{\wt\phi}$.

    \item When $\rho_1= \ldots = \rho_m \gg 0$, the support of $\chi_{\wt\phi}$ is contained in the positive cylindrical end.

    \item $\chi_{p_1 \cdots p_m; p'}^{\mb X}$ respects the structural maps of the multimodule about the concatenation from the negative side. More precisely, when $p_j < q_j$ and $\wt\phi \in \wt B_{p_1 \cdots p_{j-1}(p_j q_j) p_{j+1} \cdots p_m; p'}^{\mb X}$, $\chi_{\wt\phi}$ is determined by $\chi_{p_j q_j}^{{\mb F}_j}$ and $\chi_{p_1 \cdots p_{j-1} q_j p_{j+1} \cdots p_m; p'}^{\mb X}$ in the canonical way.

    \item $\chi_{p_1 \cdots p_m; p'}^{\mb X}$ respects the structural map of the multimodule about concatenation from the positive side. More precisely, if $\wt\phi \in \wt B_{p_1 \cdots p_m; q'p'}^{\mb X}$, then $\chi_{\wt \phi}$ is canonically determined by $\chi_{p_1 \cdots p_m; q'}^{\mb X}$ and $\chi_{q'p'}^{{\mb F}'}$.
 
    \item $\chi_{p_1\cdots p_m; p'}^{\mb X}$ is collared. Namely, near each stratum, it is constant in collar coordinates.

    \item $\chi_{p_1 \cdots p_m; p'}^{\mb X}$ respects the rigidification. 

    \item The collection of functions are equivariant with respect to the action by $\Pi^m$. Namely, if $a_1, \ldots, a_m \in \Pi$ and $q_1 = a_1 p_1, \ldots, q_m = a_m p_m$, $q' = (a_1 + \cdots + a_m) p'$, with respect to the identification $\wt C_{p_1 \cdots p_m; p'}^{\mb X} \cong \wt C_{q_1 \cdots q_m; q'}^{\mb X}$, one has 
    \beqn
    \chi_{p_1 \cdots p_m; p'}^{\mb X} = \chi_{q_1 \cdots q_m; q'}^{\mb X}.
    \eeqn
    \end{enumerate}
\end{lemma}

Now similar to the flow category case, for each $\wt x = (\wt \phi, \wt u, e, F) \in \wt V_{p_1\cdots p_m; p'}$, there are a pair of Banach spaces
\begin{align*}
&\ {\mc E}_{\wt x}^{1},\ &\ {\mc E}_{\wt x}^{0}
\end{align*}
which form equivariant ``vector bundles''
\begin{align*}
&\ {\mc E}_{p_1 \cdots p_m; p'}^1,\ &\ {\mc E}_{p_1 \cdots p_m; p'}^0.
\end{align*}
The cut-off functions also induce a family of Cauchy--Riemann operators
\beqn
D_{\wt x}: {\mc E}_{\wt x}^1 \to {\mc E}_{\wt x}^0.
\eeqn

Before we state the definition of index bundle interpolators for the multimodule case, we specify certain boundary strata of $\wt V_{p_1 \cdots p_m; p'}^{\mb X}$. For $j = 1, \ldots, m$, there is a corresponding codimension 1 stratum 
\beqn
\partial^j \wt V_{p_1 \cdots p_m; p'}
\eeqn
corresponding to configurations where the $j$-th marked point goes to the very left of the underlying curve. There is also a positive end
\beqn
\partial^+  \wt V_{p_1 \cdots p_m; p'} \cong V_{p_1 \cdots p_m; p'}
\eeqn
corresponding to configurations where all the $m$ marked points go to the very right of the underlying curve. Notice that we also have a canonical identification
\beqn
\partial^- \wt V_{p_1 \cdots p_m; p'}:= \bigcap_{j=1}^m \partial^j \wt V_{p_1 \cdots p_m; p'} \cong V_{p_1 \cdots p_m; p'}.
\eeqn

We consider the compactly supported version of ${\mc E}_{p_1 \cdots p_m; p'}^0$, which is denoted by 
\beqn
\mathring {\mc E}_{p_1 \cdots p_m; p'}^0.
\eeqn

We give the analogue of Definition \ref{defn_bundle_interpolator}.

\begin{defn}\label{defn_bundle_interpolator_module}
Suppose we have chosen index bundle interpolators over $\hat {\mb F}_1, \ldots, \hat {\mb F}_m; \hat {\mb F}'$ independently. An {\bf index bundle interpolator} over $\hat {\mb X}$ (compatible with the existing index bundle interpolators) consists of the following items. 
\begin{enumerate}

\item A collection of finite-dimensional unitary representations $W_{p_1 \cdots p_m; p'}^{\mb X}$ of $G_{p_1 \cdots p_m; p'}^{\mb X}$.

\item A collection of equivariant complex-linear bundle maps
\beqn
\wt \nu_{p_1 \cdots p_m; p'}^{\mb X}: \bigoplus_{j=1}^m \uds R_{p_j}^+ \oplus \uds W_{p_1 \cdots p_m; p'}^{\mb X} \oplus \uds R_{p'}^+ \to \mathring {\mc E}_{p_1 \cdots p_m; p'}^0
\eeqn

\item A collection of equivariant complex linear maps
\beqn
\wt\iota_{p_1 \cdots p_{i-1} (p_i q_i) p_{i+1} \cdots p_m; p'}: W_{p_i q_i}^{\mb{F}_i} \oplus R_{q_i}^+ \oplus W_{p_1 \cdots p_{i-1} q_i p_{i+1} \cdots p_m; p'} \to W_{p_1 \cdots p_m; p'}.
\eeqn

\item A collection of equivariant complex linear maps
\beqn
\wt \iota_{p_1 \cdots p_m; q'p'}: W_{p_1 \cdots p_m; q'}^{\mb X} \oplus R_{q'}^+ \oplus W_{q'p'}^{{\mb F}'} \to W_{p_1 \cdots p_m; p'}^{\mb X}.
\eeqn
\end{enumerate}
They are required to satisfy the following conditions.
\begin{enumerate}

     \item For each $j$ and $\wt x = (\wt \phi, \wt u, e, F) \in \partial^j \wt V_{p_1 \cdots p_m; p'}$, $\wt\nu_{\wt x}$ sends $R_{p_j}^+$ via $\nu_{p_j}$ into the component over the domain of $u_{p_j}$ while sends the complement of $R_{p_j}^+$ into other components.

    \item {\bf (The original index bundle)} For each $\wt x = (\wt \phi, \wt u, e, F) \in \partial^- \wt V_{p_1 \cdots p_m; p'}$, with respect to the decomposition
    \beqn
    \mathring {\mc E}_{\wt x}^0 \cong \bigoplus_{j=1}^m \Omega_c^{0,1}({\mb R}\times S^1, u_{p_j}^* TX) \oplus \Omega_c^{0,1}( C_\phi, u^* TX)
    \eeqn
    we have
    \beqn
    \wt\nu_{p_1 \cdots p_m; p'}^{\mb X} = \left[  \begin{array}{ccc} \bigoplus_{j=1}^m \nu_{p_j} & 0 & 0 \\
    0 & 0 & 0 \end{array}\right].
    \eeqn

    \item {\bf (The complex index bundle)} For each $\wt x = (\wt \phi, \wt u, e, F) \in \partial^+ \wt B_{p_1 \cdots p_m; p'}^{\mb X}$, with respect to the decomposition 
    \beqn
    \wt\nu_{p_1 \cdots p_m; p'}^{\mb X} = \left[ \begin{array}{ccc} 0 & \nu_{\wt x} & 0 \\   0 & 0 & \nu_{p'}\end{array} \right]
\eeqn
where 
\beqn
\nu_{\wt x}: W_{p_1 \cdots p_m; p'}^{\mb X} \to \Omega_c^{0,1}(C_\phi, u^* TX)
\eeqn
is a complex-linear map.

    \item {\bf (Flow multimodule structure I)} For each $j$ and $p_j < q_j$, the above conditions induces the following bundle maps
\beqn
\uds R_{p_j}^+ \oplus \uds W_{p_j q_j} \oplus \uds R_{q_j}^+ \oplus \uds W_{p_1 \cdots p_{j-1}q_j p_{j+1} \cdots p_m; p'}^{\mb X} \oplus \uds R_{p'}^+  \to \wt V_{p_1 \cdots p_{j-1}(p_j q_j) p_{j+1} \cdots p_m; p'}.
\eeqn
Then the following diagram commutes.
\beqn
\xymatrix{ \displaystyle \bigoplus_{i=1}^m \uds R_{p_i}^+ \oplus \uds W_{p_j q_j} \oplus \uds R_{q_j}^+ \oplus \uds W_{p_1 \cdots p_{j-1}q_j p_{j+1} \cdots p_m; p'}^{\mb X} \oplus \uds R_{p'}^+  \ar[rr] \ar[d] & &  \mathring {\mc E}_{p_1 \cdots p_m; p'}^0 |_{\wt V_{p_1 \cdots p_{j-1}(p_j q_j) p_{j+1} \cdots p_m; p'}} \ar[d] \\
      \displaystyle  \bigoplus_{i=1}^m \uds R_{p_i}^+ \oplus \uds W_{p_1 \cdots p_m; p'}^{\mb X} \oplus \uds R_{p'}^+ \ar[rr]   & &   \mathring {\mc E}_{p_1 \cdots p_m; p'}^0}
\eeqn

\item {\bf (Flow multimodule structure II)} For each $i$ and $p_i < r_i < q_i$, the following diagram of linear maps commutes.
\beqn
\xymatrix{  W_{p_i r_i}^{{\mb F}_i} \oplus R_{r_i}^+ \oplus W_{r_i q_i}^{{\mb F}_i} \oplus R_{q_i}^+ \oplus W_{p_1 \cdots p_{i-1} q_i p_{i+1} \cdots p_m; p'}^{\mb X} \ar[rrr] \ar[d] & & &  W_{p_i r_i}^{{\mb F}_i} \oplus R_{r_i}^+ \oplus W_{p_1 \cdots p_{i-1} r_i p_{i+1} \cdots p_m; p'}^{\mb X} \ar[d] \\
   W_{p_i q_i}^{{\mb F}_i} \oplus R_{q_i}^+ \oplus W_{p_1 \cdots p_{i-1} q_i p_{i+1} \cdots p_m; p'}^{{\mb X}} \ar[rrr] &  & &     W_{p_1 \cdots p_m; p'}^{{\mb X}}}.
\eeqn

\item {\bf (Flow multimodule structure III)} For each $i < j$, $p_i< q_i$, $p_j < q_j$, the following diagram of linear maps commutes.
\beqn
\xymatrix{  W_{p_i q_i}^{{\mb F}_i} \oplus R_{q_i}^+ \oplus W_{p_j q_j}^{{\mb F}_j} \oplus R_{q_j}^+ \oplus W_{p_1 \cdots p_{i-1} q_i p_{i+1} \cdots p_{j-1} q_j p_{j+1} \cdots p_m; p'}^{\mb X}  \ar[d] \ar@/^2pc/[rdd] & \\
   W_{p_i q_i}^{{\mb F}_i} \oplus R_{q_i}^+ \oplus W_{p_1 \cdots p_{i-1} q_i p_{i+1} \cdots p_m; p'}^{\mb X} \ar@/_2pc/[rdd]  & \\
   &  W_{p_j q_j}^{{\mb F}_j} \oplus R_{q_j}^+ \oplus W_{p_1 \cdots p_{j-1} q_j p_{j+1} \cdots p_m; p'}^{\mb X} \ar[d] \\
    & W_{p_1 \cdots p_m; p'}^{{\mb X}} & } 
\eeqn

\item {\bf (Flow multimodule structure IV)} For $q' < p'$, the following diagram commutes.
\beqn
\xymatrix{ \displaystyle   \bigoplus_{i=1}^m \uds R_{p_i}^+ \oplus \uds W_{p_1 \cdots p_m; q'}^{\mb X} \oplus \uds R_{q'}^+ \oplus \uds W_{q'p'}^{{\mb F}'} \oplus \uds R_{p'}^+ \ar[rr]  \ar[d]   &  &  \mathring {\mc E}_{p_1 \cdots p_m; p'}^0 \ar[d]  \\
        \displaystyle     \bigoplus_{i=1}^m \uds R_{p_i}^+ \oplus \uds W_{p_1 \cdots p_m; p'}^{\mb X} \oplus \uds R_{p'}^+ \ar[rr]   & &  \mathring {\mc E}_{p_1 \cdots p_m; p'}^0    }
\eeqn

\item {\bf (Flow multimodule structure V)} For $q' < r' < p'$, the following diagram commutes.
\beqn
\xymatrix{ W_{p_1 \cdots p_m; q'}^{\mb X} \oplus R_{q'}^+ \oplus W_{q'r'}^{{\mb F}'} \oplus R_{r'}^+ \oplus W_{r'p'}^{{\mb F}'} \ar[rrr] \ar[d] & & &  W_{p_1 \cdots p_m; r'}^{{\mb X}} \oplus R_{r'}^+ \oplus W_{r'p'}^{{\mb F}'} \ar[d]\\
     W_{p_1 \cdots p_m; q'}^{{\mb X}} \oplus R_{q'}^+ \oplus W_{q'p'}^{{\mb F}'} \ar[rrr] & & &  W_{p_1 \cdots p_m; p'}^{{\mb X}}}
\eeqn

\item {\bf (Flow multimodule structure VI)} For $p_j< q_j$ and $q' < p'$, the following diagram commutes.
\beqn
\xymatrix{ W_{p_j q_j}^{{\mb F}_j} \oplus R_{q_j}^+ \oplus W_{p_1 \cdots p_{j-1} q_j p_{j+1} \cdots p_m; q'}^{{\mb X}} \oplus R_{q'}^+ \oplus W_{q'p'}^{{\mb F}'} \ar[rrr] \ar[d] & & &   W_{p_1 \cdots p_m; q'}^{{\mb X}} \oplus R_{q'}^+ \oplus W_{q'p'}^{{\mb F}'}\ar[d]\\
W_{p_j q_j}^{{\mb F}_j} \oplus R_{q_j}^+ \oplus W_{p_1 \cdots p_{j-1} q_j p_{j+1} \cdots p_m; p'}^{{\mb X}} \ar[rrr] & & & W_{p_1 \cdots p_m; p'}^{{\mb X}}}
\eeqn

    \end{enumerate}
\end{defn}

\begin{lemma}\label{lemma1718}
Given index bundle interpolators for $\hat {\mb F}_1, \ldots, \hat {\mb F}_m; \hat {\mb F}'$, there exists an index bundle interpolator for $\hat {\mb X}$.
\end{lemma}

\begin{proof}
Induction. 
\end{proof}

Then we obtain a family of continuous real $G_{p_1 \cdots p_m; p'}^{\mb X}$-equivariant vector bundles
\beqn
\wt I_{p_1 \cdots p_m; p'} \to \wt V_{p_1 \cdots p_m; p'}
\eeqn
whose fiber over $\wt x$ is the kernel of the operator 
\beqn
\wt \nu_{p_1 \cdots p_m; p'}^{\mb X} \oplus D_{\wt x}.
\eeqn
We examine its restrictions to specific lower strata. By one of the conditions of Definition \ref{defn_bundle_interpolator_module}, one has 
\beqn
\wt I_{p_1 \cdots p_m; p'}|_{\partial^- \wt V_{p_1 \cdots p_m; p'}} \cong \bigoplus_{j=1}^m \uds R_{p_j}^- \oplus \uds W_{p_1 \cdots p_m; p'}^{\mb X} \oplus \uds R_{p'}^+ \oplus T^{\rm vt} V_{p_1 \cdots p_m; p'}
\eeqn
which is a stabilization of the original vertical tangent bundle. On the other hand,
\beqn
\wt I_{p_1 \cdots p_m; p'}|_{\partial^+ \wt V_{p_1 \cdots p_m; p'}} \cong \bigoplus_{j=1}^m \uds R_{p_j}^+ \oplus I_{p_1 \cdots p_m; p'}^{\mb C} \oplus \uds R_{p'}^-
\eeqn
where $I_{p_1 \cdots p_m; p'}^{\mb C}$ is fiberwise the kernel of the complex-linear operator
\beqn
\nu_{\wt x} \oplus D_u^{\mb C}.
\eeqn
Since both the negative and positive inclusions
\beqn
V_{p_1 \cdots p_m; p'} \cong \partial^\pm \wt V_{p_1 \cdots p_m; p'} \to \wt V_{p_1 \cdots p_m; p'}
\eeqn
are equivariant homotopy equivalences, the above two restriction bundles are isomorphic. One needs to choose a compatible system of isomorphisms in order to obtain a relative normal complex structure on the Kuranishi multimodule.

We consider certain codimension 1 strata of the thickenings. 
\begin{enumerate}
    \item Suppose $p_j < q_j$ in ${\mb F}_j$. Then there is a corresponding relative space $\wt V_{p_1 \cdots p_{j-1}(p_j q_j) p_{j+1} \cdots p_m; p'} \to \wt B_{p_j q_j}\boxtimes \wt B_{p_1 \cdots p_{j-1}q_j p_{j+1} \cdots p_m; p'}$. Notice that there are further lower strata of the form
    \beqn
    \partial^- \wt V_{p_1 \cdots p_{j-1} (p_j q_j) p_{j+1} \cdots p_m; p'} \cong \partial^- \wt V_{p_j q_j} \times \partial^- \wt V_{p_1 \cdots p_{j-1} q_j p_{j+1} \cdots p_m; p'}
    \eeqn
    corresponding to where the $j$-th marked point is on the very left,
    \beqn
    \partial^+ \wt V_{p_1 \cdots p_{j-1} (p_j q_j) p_{j+1} \cdots p_m; p'} \cong \partial^+ \wt V_{p_j q_j} \times \partial^+ \wt V_{p_1 \cdots p_{j-1} q_j p_{j+1} \cdots p_m; p'}
    \eeqn
    where all the marked points are on the very right. Notice that here we have the usual direct product but not $\boxtimes$. There is also an intermediate strata
    \beqn
    \partial^0 \wt V_{p_1 \cdots p_{j-1} (p_j q_j) p_{j+1} \cdots p_m; p'} \cong \partial^+ \wt V_{p_j q_j} \times \partial^- \wt V_{p_1 \cdots p_{j-1} q_j p_{j+1} \cdots p_m; p'}
    \eeqn
    corresponding to configurations where the $j$-th marked point is at the breaking at $q_j$ and all other marked points are on the very left. Notice that all these strata are isomorphic with $V_{p_1 \cdots p_{j-1}(p_j q_j) p_{j+1} \cdots p_m; p'}$. Notice that, the bundle $\wt I_{p_1 \cdots p_m; p'}$, when restricted to $\wt V_{p_1 \cdots p_{j-1} (p_j q_j) p_{j+1} \cdots p_m; p'}$, contains a subbudnle $\wt I_{p_1 \cdots p_{j-1}(p_j q_j) p_{j+1} \cdots p_m; p'}$ coming from the index bundle interpolators, whose restrictions to the above three strata are canonically identified with 
    \begin{multline*}
    \wt I_{p_1 \cdots p_{j-1} (p_j q_j) p_{j+1} \cdots p_m; p'}|_{\partial^- \wt V_{p_1 \cdots p_{j-1} (p_j q_j) p_{j+1} \cdots p_m; p'}} \\
    \cong \wt I_{p_j q_j}|_{\partial^- \wt V_{p_j q_j}} \oplus \wt I_{p_1 \cdots p_{j-1} q_j p_{j+1} \cdots p_m; p'}|_{\partial^- \wt V_{p_1 \cdots p_{j-1} q_j p_{j+1} \cdots p_m; p'}},
    \end{multline*}
\begin{multline*}
    \wt I_{p_1 \cdots p_{j-1} (p_j q_j) p_{j+1} \cdots p_m; p'}|_{\partial^0 \wt V_{p_1 \cdots p_{j-1} (p_j q_j) p_{j+1} \cdots p_m; p'}} \\
    \cong \wt I_{p_j q_j}|_{\partial^+ \wt V_{p_j q_j}} \oplus \wt I_{p_1 \cdots p_{j-1} q_j p_{j+1} \cdots p_m; p'}|_{\partial^- \wt V_{p_1 \cdots p_{j-1} q_j p_{j+1} \cdots p_m; p'}},
    \end{multline*}
    and
    \begin{multline*}
    \wt I_{p_1 \cdots p_{j-1} (p_j q_j) p_{j+1} \cdots p_m; p'}|_{\partial^+ \wt V_{p_1 \cdots p_{j-1} (p_j q_j) p_{j+1} \cdots p_m; p'}} \\
    \cong \wt I_{p_j q_j}|_{\partial^+ \wt V_{p_j q_j}} \oplus \wt I_{p_1 \cdots p_{j-1} q_j p_{j+1} \cdots p_m; p'}|_{\partial^+ \wt V_{p_1 \cdots p_{j-1} q_j p_{j+1} \cdots p_m; p'}}.
    \end{multline*}

    \item Suppose $q' < p'$ in ${\mb F}'$. Then there is a corresponding relative space $\wt V_{p_1 \cdots p_m; q'p'} \to \wt B_{p_1 \cdots p_m; q'} \boxtimes \wt B_{q'p'}$. The three further lower strata are 
    \beqn
    \begin{split}
        \partial^- \wt V_{p_1 \cdots p_m; q'p'} \cong &\ \partial^- \wt V_{p_1 \cdots p_m; q'} \times \partial^- \wt V_{q'p'},\\
        \partial^0 \wt V_{p_1 \cdots p_m; q'p'} \cong &\ \partial^+ \wt V_{p_1 \cdots p_m; q'} \times \partial^- \wt V_{q'p'},\\
        \partial^+ \wt V_{p_1 \cdots p_m; q'p'} \cong &\ \partial^+ \wt V_{p_1 \cdots p_m; q'} \times \partial^+ \wt V_{q'p'}.
    \end{split}
    \eeqn
    The extended index bundle $\wt I_{p_1 \cdots p_m; p'}$, when restricted to then, has the subbundle $\wt I_{p_1 \cdots p_m; q'p'}$ and is identified respectively as 
    \beqn
    \begin{split}
       \wt I_{p_1 \cdots p_m; q'p'}|_{\partial^- \wt V_{p_1 \cdots p_m; q'p'}} \cong &\  \wt I_{p_1 \cdots p_m; q'}|_{\partial^- \wt V_{p_1 \cdots p_m; q'}} \oplus \wt I_{q'p'}|_{\partial^- \wt V_{q'p'}},\\
       \wt I_{p_1 \cdots p_m; q'p'}|_{\partial^0 \wt V_{p_1 \cdots p_m; q'p'}} \cong &\ \wt I_{p_1 \cdots p_m; q'}|_{\partial^+ \wt V_{p_1 \cdots p_m; q'}} \oplus \wt I_{q'p'}|_{\partial^- \wt V_{q'p'}},\\
        \wt I_{p_1 \cdots p_m; q'p'}|_{\partial^+ \wt V_{p_1 \cdots p_m; q'p'}} \cong &\ \wt I_{p_1 \cdots p_m; q'}|_{\partial^+ \wt V_{p_1 \cdots p_m; q'}} \oplus \wt I_{q'p'}|_{\partial^+ \wt V_{q'p'}}.
    \end{split}
    \eeqn
\end{enumerate}

\begin{lemma}\label{lemma1714}
Given bundle isomorphisms 
\begin{align*}
    &\ \wt I_{p_j q_j}|_{\partial^- \wt V_{p_j q_j}} \cong \wt I_{p_j q_j} |_{\partial^+ \wt V_{p_j q_j}},\ &\ \wt I_{p'q'}|_{\partial^- \wt V_{p' q'}} \cong \wt I_{p'q'}|_{\partial^+ \wt V_{p'q'}}
\end{align*}
satisfying conditions listed in Lemma \ref{lemma1713}, there exists a system of $G_{p_1 \cdots p_m; p'}$-equivariant real vector bundle isomorphisms
\beqn
\wt I_{p_1 \cdots p_m; p'}|_{\partial^- \wt V_{p_1 \cdots p_m; p'}} \cong \wt I_{p_1 \cdots p_m; p'}|_{\partial^+ \wt V_{p_1 \cdots p_m; p'}}
\eeqn
satisfying the following conditions.
\begin{enumerate}
    \item When $p_j < q_j$, the following diagram commutes

    \item When $q' < p'$, the following diagram commutes.
    \beqn
    \xymatrix{ \wt I_{p_1 \cdots p_m; q'p'}|_{\partial^- \wt V_{p_1 \cdots p_m; q'} \times \partial^- \wt V_{q'p'}} \ar[d] \ar[rr] &  &   \wt I_{p_1 \cdots p_m; q'p'}|_{\partial^+ \wt V_{p_1 \cdots p_m; q'} \times \partial^+ \wt V_{q'p'}} \ar[d]\\
    \wt I_{p_1 \cdots p_m; p'}|_{\partial^- \wt V_{p_1 \cdots p_m; q'p'}} \ar[rr] & & \wt I_{p_1 \cdots p_m; p'}|_{\partial^+ \wt V_{p_1 \cdots p_m; q'p'}}
    }
    \eeqn

    \item The isomorphisms are collared.

    \item The isomorphisms respect the rigidifications.
\end{enumerate}
\end{lemma}

\begin{proof}
Induction.
\end{proof}

Now using the identifications of the extended index bundles over the negative and positive strata, one obtains the equivariant bundle isomorphisms
\beqn
\bigoplus_{j=1}^m \uds R_{p_j}^- \oplus \uds W_{p_1 \cdots p_m; p'}^{\mb X} \oplus \uds R_{p'}^+ \oplus T^{\rm vt} V_{p_1 \cdots p_m; p'} \cong \bigoplus_{j=1}^m \uds R_{p_j}^+ \oplus I_{p_1\cdots p_m; p'}^{\mb C} \oplus \uds R_{p'}^-
\eeqn
In particular, this gives an equivariant normal complex structure (see Definition \ref{defn_equivariant_NC_structure}) of the vertical tangent bundle $T^{\rm vt} V_{p_1 \cdots p_m; p'}$. Hence the chart $K_{p_1 \cdots p_m; p'}$ becomes an object of $\outer \uds{\bf S^{\rm rel} Kur}_{\rm rig}^{\rm NC}$. The conditions of Lemma \ref{lemma1714} imply that the structural maps are also morphisms in this category. Hence one obtains a lift of the Kuranishi flow category $\hat {\mb F}$ to $\outer \uds{\bf S^{\rm rel} Kur}_{\rm rig}^{\rm NC}$. This finishes the proof of Theorem \ref{thm173}.

\subsection{Orientations}\label{subsection165}

We explain how to obtain coherent orientations on the AMS construction. We give orientations on the tangent bundles of the domain moduli spaces as well as the vertical tangent bundles. After the stable smoothing to be done in the next section, they induce coherent orientations on the AMS lifts.

We first look at the operator $D_p$ given by \eqref{operator_dp}. Its Fredholm index is congruent (mod $2$) to the cohomological degree of the periodic orbit $p$, denoted by $|p|$. Let $\det D_p$ be its real determinant line. 

Consider the flow category case. For each pair of objects $p, q$, one can use the 1-parameter family of operators \eqref{family_CR_operator} to obtain an interpolation; in particular an isomorphism between determinant lines
\beqn
\det D_p \otimes \det D_u \cong \det (D_u^{\mb C} \oplus D_q) \cong \det D_q. 
\eeqn
The last isomorphism comes from the fact that $D_u^{\mb C}$ is complex. We may regard the orientation of the moduli space as the family of isomorphismss on the determinant lines
\beqn
\det D_u: \det D_p \to \det D_q.
\eeqn
In terms of the AMS lift, this means one obtains a family of isomorphisms (well-defined up to homotopy) 
\beqn
\sigma_{pq}: \det D_p \to \det D_q
\eeqn
over the thickened moduli $V_{pq}$. The construction implies
\beqn
\sigma_{rq} \circ \sigma_{pr} = \sigma_{pq}.
\eeqn

The case with multimodules and homotopies are essentially the same. For a Floer multimodule ${\mb X}$, using the interpolation constructed before, one obtains a family of isomorphisms
\beqn
\sigma_{p_1 \cdots p_m; p'}^{\mb X}: \bigotimes_{i=1}^m \det D_{p_i} \to \det D_{p'}.
\eeqn
By considering indices, one has that
\beqn
\sigma_{p_1 \cdots p_m; p'}^{\mb X} = (-1)^{( \sum_{j<i}|p_j|) (|p_i|- |q_i|)} \sigma_{p_1 \cdots p_{i-1} q_i p_{i+1} \cdots p_m; p'}^{\mb X} \circ \sigma_{p_i q_i}^{{\mb F}_i} 
\eeqn
and
\beqn
\sigma_{p_1 \cdots p_m; p'}^{\mb X} = \sigma_{q'p'}^{{\mb F}'}\circ \sigma_{p_1 \cdots p_m;q'}^{\mb X}.
\eeqn
The case of homotopy is similar and we omit it.

\section{Stable Smoothings}

We proceed with the last step of constructing the AMS lifts. The original argument of such abstract smoothing of moduli spaces was contained in \cite{AMS}. The inductive argument in the setting of Floer theories was first systematically discussed in \cite{Bai_Xu_Arnold}.

\subsection{Main results about stable smoothings}

In the construction we will frequently stabilize thickened moduli spaces obtained in the AMS construction. To simplify the discussion, we always stabilize by trivial complex vector bundles. We make the following formal definition.

\begin{defn}
The category $\uds{\bf Kur}^+$ has objects being $(K, R)$ where $K = (G, V, E, S)$ is a Kuranishi space and $R$ is a finite-dimensional unitary complex representation of $G$. A strict embedding from $(K_1, R_1)$ to $(K_2, R_2)$ consists of a strict embedding of Kuranishi spaces $\iota_{21}: K_1 \to K_2$ and an $\iota_{21}^G$-equivariant isometric linear embedding $\iota_{21}^R: R_1 \to R_2$. Morphisms of $\uds{\bf Kur}^+$ are unitary conjugacy classes of strict embeddings. 
\end{defn}

Similarly, one can consider the collared and rigidified versions. Then there are two natural functors from $\uds{\bf Kur}^+$ to $\uds{\bf Kur}$: the first is to forget the representation $R$; the other is to take the stabilization by $R$. What we would consider is the version
\beqn
\outer\uds{\bf S^{\rm rel}Kur}_{\rm rig}^+
\eeqn
for collared, rigidified, relatively smooth Kuranishi spaces.

\begin{defn}\label{defn_stabilizer}
A {\bf stabilizer} of a flow category/multimodule/homotopy enriched in $\uds{\bf Kur}$ is a lift to $\uds{\bf Kur}^+$.
\end{defn}

\begin{rem}
In the AMS construction, a stabilizer can be viewed as adding trivial pieces to the thickening data. 
\end{rem}

Notice that there are two natural functors from $\outer \uds{\bf S^{\rm rel} Kur}_{\rm rig}^+$ to $\outer \uds{\bf S^{\rm rel} Kur}_{\rm rig}$. The first is to forget the stabilizer. The second is to stabilize each Kuranishi space by the corresponding trivial vector bundle via the vertical stabilization (see Definition \ref{defn_relative_space}). We often freely switch between the stabilizer, which is a system of representations, and the corresponding vertical stabilization, viewed as total space of the corresponding product equivariant vector bundles.

\begin{defn}
A {\bf smoothing} of a flow category/multimodule/homotopy enriched in $\outer \uds{\bf Kur}_{\rm rig}$ is a lift to $\outer \uds{\bf SKur}_{\rm rig}$. A {\bf stable smoothing} of a flow category/multimodule/homotopy enriched in $\outer \uds{\bf Kur}_{\rm rig}$ consists of a lift to $\outer \uds{\bf Kur}_{\rm rig}^+$ and a smoothing of the corresponding stabilization.
\end{defn}

\begin{thm}(Stable Smoothing Theorem)
\begin{enumerate}

\item In Situation \ref{situationf1}, suppose $\hat{\mb F}$ is a relatively smooth AMS lift of ${\mb F}$ subject to a certain outercollaring width. Then there exists a stable smoothing of $\hat {\mb F}$. In other words, there exists a vertical stabilization $\hat{\mb F}^+$ and a  smoothing of $\hat {\mb F}^+$. We still denote the flow category enriched in $\outer \uds{\bf SKur}_{\rm rig}$ by $\hat {\mb F}^+$ and call it an {\bf AMS stable smoothing} of $\hat {\mb F}$.

\item In Situation \ref{situationm1}, suppose $\hat{\mb F}_1, \cdots, \hat {\mb F}_m; \hat {\mb F}'$ are relatively smooth AMS lifts and $\hat {\mb X}$ is a relatively smooth AMS lift as a multimodule over $(\hat{\mb F}_1, \ldots, \hat{\mb F}_m; \hat {\mb F}')$.  Suppose $\hat {\mb F}_1^+, \ldots, \hat {\mb F}_m^+, \hat {\mb F}'{}^+$ are AMS stable smoothings of $\hat {\mb F}_1, \ldots, \hat {\mb F}_m, \hat {\mb F}'$ respectively. Then there exists a class of stable smoothings of $\hat {\mb X}$ as multimodules over $(\hat {\mb F}_1^+, \ldots, \hat {\mb F}_m^+; \hat {\mb F}'{}^+)$. These stable smoothings are called {\bf AMS stable smoothings} of $\hat {\mb X}$. 

\item In Situation \ref{situationh3}, let $\hat {\mb F}_1^+, \cdots, \hat {\mb F}_m^+, \hat {\mb F}'{}^{+}$ be AMS stable smoothings of the relatively smooth AMS lifts $\hat {\mb F}_1, \ldots, \hat {\mb F}_m, \hat {\mb F}'$. Let $\hat {\mb X}_0^+$ resp. $\hat {\mb X}_1^+$ be AMS stable smoothings of $\hat {\mb X}_0$ resp. $\hat {\mb X}_1$ as multimodules over $(\hat {\mb F}_1^+, \ldots, \hat {\mb F}_m^+; \hat {\mb F}'{}^+)$. Then the relatively smooth AMS lift $\hat {\mb H}$ has a class of stable smoothings as a homotopy from $\hat {\mb X}_0^+$ to $\hat {\mb X}_1^+$. Such stable smoothings are called {\bf AMS stable smoothings} of $\hat {\mb H}$.
\end{enumerate}
\end{thm}

\subsection{Preliminaries}

\subsubsection{Microbundles}

Recall that a {\bf microbundle} of rank $k$ over a topological space   $M$, denoted by ${\ms F}$, is a topological space ${\ms F}$ together with two structural maps $p: {\ms F} \to M$, $i: M\to {\ms F}$ such that $p\circ i = {\rm Id}_M$ and such that near the image of $i$, locally ${\ms F}$ is an ${\mb R}^k$-bundle. Any real vector bundle $F \to M$ together with the bundle projection and the zero section induces a microbundle $F_\mu$. If $M$ is a topological manifold, then its tangent microbundle ${\ms T}_\mu M$ has the total space $M \times M$ with $p$ the projection onto the first factor and $i$ the diagonal embedding of $M$.

A {\bf morphism} of microbundles over $M$ from ${\ms F}$ to ${\ms F}'$ is a germ of continuous maps from a neighborhood of $i(M) \subset {\ms F}$ to ${\ms F}'$ which commutes with the structural maps. A morphism is called an {\bf isomorphism} if it can be represented by a homeomorphism. An {\bf isotopy} of microbundle isomorphisms from $\phi_0: {\ms F} \to {\ms F}'$ to $\phi_1: {\ms F} \to {\ms F}'$ is a microbundle isomorphism from $\tilde {\ms F} \to \tilde {\ms F}'$ over $M \times [0, 1]$ whose restriction to $M\times \{i\}$ coincides with $\phi_i$ (as germs) for $i = 0, 1$.

In particular, a {\bf vector bundle reduction} of a microbundle ${\ms F} \to M$ consists of a vector bundle $F \to M$ and a microbundle isomorphism $F_\mu \to {\ms F}$. When a topological manifold $M$ is equipped with a smooth structure, denoted by $M_\alpha$, then there is a vector bundle reduction $(TM_\alpha)_\mu \to {\ms T}_\mu M$: a concrete choice could be the exponential map associated to a Riemannian metric. The isotopy class of such a vector bundle reduction is canonical. Hence it is a necessary condition for a topological manifold to admit a smooth structure that its tangent microbundle admits a vector bundle reduction. 

All the above notions admit equivariant extensions with respect to actions by compact Lie groups. In particular, when a topological $G$-manifold $M$ is equipped with a smooth structure $M_\alpha$ such that the $G$-action is smooth, then there is a well-defined isotopy class of $G$-equivariant vector bundle reduction $(TM_\alpha)_\mu \to {\ms T}_\mu M$.

\begin{rem}
We need the notion of smooth microbundles over smooth manifolds. These are microbundles whose total spaces and structural maps are all in the $C^\infty$ category. One can also define smooth microbundle maps and smooth vector bundle reductions. Equivariant versions can be formulated as well. For example, when $M$ is already a smooth manifold, the tangent microbundle is smooth and any vector bundle reduction of ${\ms T}_\mu M$ using a smooth Riemannian metric is smooth.
\end{rem}

\subsubsection{Stable $G$-smoothings}

Now we recall notions about equivariant smoothings. Let $M$ be a topological manifold with an action by a compact Lie group $G$. A {\bf $G$-smoothing} of $M$ is a smooth structure $\alpha$ on $M$ such that the $G$-action is smooth. Two $G$-smoothings $\alpha_0$ and $\alpha_1$ are called {\bf equivalent} if the identity map $M_{\alpha_0} \to M_{\alpha_1}$ is a diffeomorphism; they are called {\bf isotopic} if the identity map is homotopic through $G$-equivariant homeomorphisms to a diffeomorphism. More generally, a {\bf stable $G$-smoothing} of $M$ consists of an orthogonal $G$-representation $R$ and a $G$-smoothing on $M \times R$. Notice that given any stable $G$-smoothing $(M\times R)_\alpha$ and any other orthogonal $G$-representation $R'$, the product smooth structure $(M\times R)_\alpha \times R'$ gives another stable $G$-smoothing. Two stable $G$-smoothings $(M\times R_1)_{\alpha_1}$ and $(M \times R_2)_{\alpha_2}$ are {\bf stably isotopic} if there exist orthogonal $G$-representations $R_1', R_2'$ such that $R_1 \oplus R_1' \cong R_2 \oplus R_2'$ as orthogonal $G$-representations and such that the two $G$-smoothings $(M \times R_1)_{\alpha_1} \times R_1'$ and $(M\times R_2)_{\alpha_2} \times R_2'$ on $M \times R_1 \times R_1' = M\times R_2 \times R_2'$ are isotopic.

Notice that given a $G$-smoothing $M_\alpha$ on $M$, there is an induced isotopy class of $G$-equivariant vector bundle reductions $(TM_\alpha)_\mu \to {\ms T}_\mu M$. In view of stable smoothings, one also need to define the notion of stable isotopy of stable vector bundle reductions.

\begin{defn}\cite[Page 300]{Lashof_1979}
Let ${\ms F} \to M$ be a $G$-microbundle. A {\bf stable $G$-vector bundle reduction} of ${\ms F}$ consists of a $G$-vector bundle $F$ and an orthogonal $G$-representation $R$ together with a $G$-microbundle isomorphism
\beqn
\phi: F_\mu \to {\ms F} \oplus \uds R_\mu.
\eeqn
Two stable $G$-vector bundle reductions $\phi_i: (F_i)_\mu \to {\ms F} \oplus (\uds R_i)_\mu$, $i = 0, 1$, are called {\bf stably isotopic} if there exist orthogonal $G$-representations $R_0', R_1'$ such that $R_0\oplus R_0' \cong R_1 \oplus R_1'$ and the stabilizations $\phi_0 \oplus {\rm Id}_{\uds R_0'}$ and $\phi_1 \oplus {\rm Id}_{\uds R_1'}$ are isotopic vector bundle reductions of ${\ms F} \oplus (\uds R_0)_\mu \oplus (\uds R_0')_\mu \cong {\ms F} \oplus (\uds R_1)_\mu \oplus (\uds R_1')_\mu$.
\end{defn}

\begin{rem}
In the above discussions, we can replace ``orthogonal representations'' by ``unitary representations'' to suit our purpose. Any orthogonal representation can be stabilized to a unitary representation.
\end{rem}

Now given any stable $G$-smoothing $(M\times R)_\alpha$ of $M$, there is an isotopy class of $G$-vector bundle reductions $T(M\times R)_\alpha \to {\ms T}_\mu( M \times R)$. The restriction of this reduction to $ M\times \{0\}$ gives a stable isotopy class of $G$-vector bundle reductions of ${\ms T}_\mu M$. Lashof's theorem says that this association has an inverse. 

\begin{thm}\cite{Lashof_1979}\label{thm_Lashof} Let $M$ be a topological $G$-manifold with finitely many orbit types. Then the above map which sends each stable isotopy class of stable $G$-smoothings to a stable isotopy class of $G$-vector bundle reductions of ${\ms T}_\mu M$ is invertible. In particular, if ${\ms T}_\mu M$ admits a $G$-vector bundle reduction, then $M$ admits a stable $G$-smoothing.
\end{thm}

\subsubsection{Micro-connections} 

Let $\pi: V \to B$ be a $G$-equivariant relatively smooth manifold such that $B$ is a smooth $G$-manifold. Then there is a well-defined $G$-microbundle over $V$:
\beqn
\pi^* {\ms T}_\mu B \oplus {\ms T}^{\rm vt}_\mu V.
\eeqn
To facilitate the stable smoothing of $V$, we need to identify the above microbundle with the tangent microbundle of $V$. Notice that there are natural microbundle maps
\beqn
{\ms T}^{\rm vt}_\mu V \to {\ms T}_\mu V \to \pi^* {\ms T}_\mu B.
\eeqn
Hence such an identification can be regarded as a microbundle version of connections.

\begin{defn}
Let $V/B$ be a $G$-equivariant topological submersion. A {\bf micro-connection} is a $G$-equivariant microbundle map $\nabla_\mu: {\ms T}_\mu V \to {\ms T}_\mu^{\rm vt} V$ satisfying
\begin{enumerate}
    \item The restriction of $\nabla^\mu$ to ${\ms T}_\mu^{\rm vt} V$ is the identity map.

    \item $\nabla_\mu$ is represented by a germ of the map $\phi: \Delta^+(V) \to V$ such that for any $(b, b') \in (\pi\times \pi)(\Delta^+(V)) \subset B\times B$ and $p \in \pi^{-1}(b) \subset V$, the restriction of $\phi$ to $\Delta^+(V) \cap (\{p\} \times \pi^{-1}(b'))$ is a homeomorphism onto an open neighborhood of $(p, p)$ in $\{p\}\times \pi^{-1}(b))$.
\end{enumerate}
\end{defn}

Notice that a micro-connection induces a microbundle isomorphism ${\ms T}_\mu M \cong \pi^* {\ms T}_\mu B \oplus {\ms T}_\mu^{\rm vt} M$.

\begin{lemma}\label{splitting_lemma}
Let $C \subset V$ be a $G$-invariant closed  subset and $U \subset V$ be a $G$-invariant open neighborhood of $C$. Let $D \subset V$ be another $G$-invariant closed set. Suppose $\nabla_\mu: {\ms T}_\mu U \to {\ms T}_\mu^{\rm vt} U$ is a $G$-invariant micro-connection. Then there exists a $G$-invariant open neighborhood $W$ of $C \cup D$ and a $G$-invariant micro-connection $\nabla_\mu': {\ms T}_\mu W \to {\ms T}_\mu^{\rm vt} W$ which coincides with $\nabla_\mu$ near $C$.
\end{lemma}

\begin{proof}
This is a restatement of the extension lemma \cite[Lemma 4.24]{AMS} and a relative version of \cite[Proposition 4.25]{AMS}.
\end{proof}

\subsection{Pre-smoothed relatively smooth Kuranishi spaces}

To describe the structures one needs to perform the stable smoothing, we enhance the category of relatively smooth Kuranishi spaces to include various specific add-ons.

\begin{defn}\label{defn:pre-smoothing}
Let $K = (G, V/B, E, S)$ be an object of $\outer \uds{\bf S^{\rm rel} Kur}_{\rm rig}$. 
\begin{enumerate}

\item A {\bf pre-smoothing} on $K$ consists of 
\begin{enumerate}

\item A $G$-invariant Riemannian metric $g^{TB}$ on $B$ and a $G$-invariant metric $h^{T^{\rm vt} V}$ on $T^{\rm vt} V$.

\item A micro-connection $\nabla_\mu$ on $V/B$.
\end{enumerate}

\item A pre-smoothing is called {\bf collared} if for each stratum $V_\alpha/B_\alpha$, the following holds.
\begin{enumerate}
    \item Near $B_\alpha$, $g^{TB}$ is the product of $g^{TB_\alpha}$ and the standard metric on the collar.

    \item Near $V_\alpha$, $h^{T^{\rm vt} V}$ is constant in collar coordinates.

    \item Near $V_\alpha$, $\nabla_\mu$ factors as 
    \beqn
    \xymatrix{ {\ms T}_\mu V \ar[r] & {\ms T}_\mu V_\alpha \ar[r] & {\ms T}_\mu^{\rm vt} V_\alpha}
    \eeqn
    where the first arrow is induced from the collar structure and the second arrow is the restriction of the micro-connection on the stratum $V_\alpha/B_\alpha$.
\end{enumerate}

\item A pre-smoothed Kuranishi space is a relatively smooth Kuranishi space equipped with a pre-smoothing.

\item A stabilization of a pre-smoothed Kuranishi space $K$ by a trivial $G$-equivariant disk bundle $D = D^{\rm hor} \oplus D^{\rm vt}$ is the stabilized Kuranishi space ${\rm Stab}_D K$ together with the naturally induced pre-smoothing (where $D^{\rm hor}$ and $D^{\rm vt}$ are equipped with linear metrics).

\item A {\bf rigidified strict embedding} of pre-smoothed Kuranishi spaces from $K_1$ to $K_2$ consists of a rigidified strict embedding $\iota_{21}: K_1 \to K_2$ (with normal disk bundle $D_{21} = D_{21}^{\rm hor} \oplus D_{21}^{\rm vt}$ equipped with the direct sum of linear metrics) such that the induced germ of strict open embedding preserves the pre-smoothings.

\item A {\bf rigidified embedding} is a unitary conjugacy class of rigidified strict embeddings. Let $\outer \uds{\bf S^{\rm pre} Kur}_{\rm rig}$ be the regular stratification category whose objects are collared, pre-smoothed Kuranishi spaces and whose morphisms are rigidified embeddings. There is an obvious forgetful functor 
\beqn
\outer \uds{\bf S^{\rm pre} Kur}_{\rm rig} \to \outer \uds{\bf S^{\rm rel} Kur}_{\rm rig}.
\eeqn

\item A {\bf pre-smoothing} of a flow category/multimodule/homotopy enriched in $\outer \uds{\bf S^{\rm rel} Kur}_{\rm rig}$ is a lift to $\outer \uds{\bf S^{\rm pre} Kur}_{\rm rig}$.
\end{enumerate}
\end{defn}

\begin{lemma}\label{lemma_presmoothing}\hfill
\begin{enumerate}
    \item In Situation \ref{situationf3}, the relatively smooth AMS lift $\hat {\mb F}$ admits a pre-smoothing.

    \item In Situation \ref{situationm3}, suppose the AMS lifts $\hat {\mb F}_1 ,\ldots, \hat {\mb F}_m; \hat {\mb F}'$ of the involved Floer flow categories are equipped with pre-smoothings, then $\hat {\mb X}$ admits a pre-smoothing as a multimodule enriched in $\outer \uds{\bf S^{\rm pre} Kur}_{\rm rig}$.

    \item In Situation \ref{situationh3}, suppose the AMS lifts $\hat {\mb F}_1, \ldots, \hat {\mb F}_m, \hat {\mb F}', \hat {\mb X}_0, \hat {\mb X}_1$ are equipped with pre-smoothings so $\hat {\mb X}_0$ and $\hat {\mb X}_1$ are multimodules over the tuple of pre-smoothed flow categories, then $\hat {\mb H}$ admits a pre-smoothing as a homotopy enriched in $\outer \uds{\bf S^{\rm pre} Kur}_{\rm rig}$.
\end{enumerate}
\end{lemma}

\begin{proof}
The collared and rigidified structure of the relatively smooth AMS lifts allows us to choose invariant horizontal and vertical Riemannian metrics inductively. On the other hand, Lemma \ref{splitting_lemma} allows us to inductively construct micro-connections. 
\end{proof}

Once one obtains a pre-smoothed Kuranishi space $K = (G, V/B, E, S)$, one obtains a $G$-vector bundle reduction of the tangent microbundle. More precisely, the exponential map of $g^{TB}$ induces a reduction 
\beqn
(TB)_\mu \to {\ms T}_\mu B.
\eeqn
The vertical exponential map of $h^{T^{\rm vt}V}$ induces a reduction
\beqn
(T^{\rm vt} V)_\mu \to {\ms T}_\mu^{\rm vt} V.
\eeqn
Together with the micro-connection, one obtains
\beqn
(\pi^* TB \oplus T^{\rm vt} V)_\mu \to {\ms T}_\mu V.
\eeqn
Notice that such $G$-vector bundle reductions are compatible with morphisms of $\outer \uds{\bf S^{\rm pre} Kur}_{\rm rig}$. 

\subsubsection{Post-smoothed Kuranishi spaces}

Definition \ref{defn:pre-smoothing} ensures the existence of coherent vector bundle lifts. The next definition is designed to make the stabizers explicit in order to carry out the inductive smoothing construction.

\begin{defn}
Let $K = (G, V/B, E, S)$ be an object of $\outer \uds{\bf S^{\rm pre} Kur}_{\rm rig}$. A {\bf post-smoothing} on $K$ consists of a $G$-smoothing of $V$, a collared smooth $G$-invariant vector bundle reduction
\beqn
\exp: (TV)_\mu \to {\ms T}_\mu V,
\eeqn
and a collared $G$-equivariant vector bundle isomorphism
\beqn
TV \cong \pi^* TB \oplus T^{\rm vt} V
\eeqn
such that $\exp$ and the $G$-vector bundle reduction induced from the pre-smoothing are in the same stable $G$-isotopy class on the interior of each stratum.
\end{defn}

Now we state the results about constructing stable smoothings. Since we need to constantly stabilize charts, we would like to use the convenient system of unitary representations of $G_d \cong U(d)$ specified by Proposition \ref{system_representation}. 

\begin{lemma}\label{lemma_postsmoothing}\hfill 
\begin{enumerate}

\item In Situation \ref{situationf3}, let $\hat {\mb F}$ be a pre-smoothed AMS lift of ${\mb F}$ provided by Lemma \ref{lemma_presmoothing}. Then there exists a stabilization $\hat {\mb F}^+$ and a post-smoothing of $\hat {\mb F}^+$. Each such flow category is called a post-smoothed AMS lift of ${\mb F}$.

\item In Situation \ref{situationm3}, let $\hat {\mb F}_1, \ldots, \hat {\mb F}_m; \hat {\mb F}'$ be pre-smoothed AMS lifts of the involved Floer flow categories. Let $\hat {\mb X}$ be a pre-smoothed AMS lift of the multimodule ${\mb X}$ provided by Lemma \ref{lemma_presmoothing}. Suppose $\hat {\mb F}_1^+, \ldots, \hat {\mb F}_m^+; \hat {\mb F}'{}^+$ are post-smoothings of $\hat {\mb F}_1, \ldots, \hat {\mb F}_m; \hat {\mb F}'$. Then there exists a stabilization $\hat {\mb X}^+$ of $\hat {\mb X}$ and a post-smoothing of it as a multimodule over $(\hat {\mb F}_1^+, \ldots, \hat {\mb F}_m^+; \hat {\mb F}'{}^+)$.

\item In Situation \ref{situationh3}, let $\hat {\mb F}_1^+, \ldots, \hat {\mb F}_m^+; \hat{\mb F}'{}^+$ be post-smoothed AMS lifts of the involved flow categories and $\hat {\mb X}_0^+$, $\hat {\mb X}_1^+$ be post-smoothed AMS lifts of the two multimodules ${\mb X}_0$ and ${\mb X}_1$. Let $\hat {\mb H}$ be a pre-smoothed AMS lift of the homotopy provided by Lemma \ref{lemma_presmoothing}. Then there exists a stabilization $\hat {\mb H}^+$ of $\hat {\mb H}$ and a post-smoothing of it as a homotopy from $\hat {\mb X}_0^+$ to $\hat {\mb X}_1^+$.
\end{enumerate}
\end{lemma}

\begin{proof}
Again, it suffices to demonstrate the argument in the case of flow categories. Recall the system of unitary representations $W_{d, \mu}$. For each pair $p<q$, denote $W_{pq, \mu}:= W_{d_{pq}, \mu}$ where $d_{pq}$ is the integer specified via the choice of integral action. Notice that once we choose a system of positive integers $\mu_{pq}$ satisfying 
\beq\label{mu_inequality}
\mu_{pq} > \max \{ \mu_{pr}, \mu_{rq}\ |\ p < r < q \}
\eeq
then one obtains a stabilization $\hat {\mb F}^+$ with each chart $K_{pq}$ vertically stabilized by $W_{pq, \mu_{pq}}$. The main claim is that one can choose such a stabilizer together with a post-smoothing of $\hat {\mb F}^+$ satisfying the condition on stable isotopy classes of $G$-vector bundle reductions. 

We construct $\mu_{pq}$ and the post-smoothing inductively. For the base case when $V_{pq}$ has no boundary strata, the existence of $W_{pq, \mu_{pq}}$ and the post-smoothing is a consequence of Lashof's theorem (Theorem \ref{thm_Lashof}), since every finite-dimensional representation of $G_{pq}$ is contained in a certain $W_{pq, \mu}$ for $\mu$ sufficiently large. Choose such $\mu_{pq}$ and denote $V_{pq}^+ = V_{pq} \times W_{pq, \mu_{pq}}$. Then we can choose a smooth $G_{pq}$-invariant Riemannian metric on $V_{pq}^+$ and obtain the exponential map $\exp_{pq}$ as a smooth vector bundle reduction of the tangent microbundle. 

Now we state the induction hypothesis. Given $p<q$, suppose we have chosen $\mu_{rs}$ for all $r, s$ satisfying $d_{rs} < d_{pq}$ satisfying \eqref{mu_inequality} and post-smoothings such that the induced embeddings of the smoothed Kuranishi spaces are morphisms of $\outer \uds{\bf SKur}_{\rm rig}$. \footnote{This claim is a purpose of using rigidified embeddings: the trivialized normal bundle induces a canonical smoothing in the tubular neighborhood of the embedding. If we do not require the normal bundle to be trivialized, then one needs to make the painful choices of smooth bundle structures.}

Now we can construct the stable smoothing on $V_{pq}$. Choose $\mu_{pq}'$ sufficiently large such that the already stabilized boundary strata can all be included into $V_{pq} \times W_{pq, \mu_{pq}'}$. On the other hand, the interior ${\rm Int} V_{pq}$ can be stably smoothed by Theorem \ref{thm_Lashof}, meaning that one can choose $\mu_{pq}''$ sufficiently large and a smoothing of ${\rm Int} V_{pq} \times W_{pq, \mu_{pq}''}$. By choosing $\mu_{pq}$ greater than both $\mu_{pq}'$ and $\mu_{pq}''$, one can assume $\mu_{pq}' = \mu_{pq}''$. Notice that the induction hypothesis implies a canonically induced smoothing near the boundary. On overlaps, these two smoothings induce the same stable $G$-isotopy class of microbundle reductions. Hence by the relative smoothing result proved in \cite[Appendix B]{Bai_Xu_Arnold}, one can interpolate the two smoothings in the collar region, possibly after a further enlargement of $\mu_{pq}$, so that near the boundary the smoothing is still the one canonically induced from other already-stabilized charts. Therefore, one obtains a stabilization by $W_{pq, \mu_{pq}}$, a smoothing of this stabilization, so that all the chart embeddings are morphisms of $\outer \uds{\bf SKur}_{\rm rig}$. Further, we choose the smooth vector bundle reduction on the stabilized chart $V_{pq}  \times W_{pq, \mu_{pq}}$. 
\end{proof}

\subsection{Normally complex AMS lifts}

Finally we are at the last step of proving Theorem \ref{thma_flow_chart} and counterparts for Floer multimodules and homotopies.

\begin{thm}\label{thm_stable_complex}\hfill
\begin{enumerate}

\item In Situation \ref{situationf3}, let $\hat {\mb F}$ be a post-smoothed AMS lift of ${\mb F}$. Then $\hat {\mb F}$, viewed as a lift in $\outer \uds{\bf SKur}_{\rm rig}$, admits a class of lifts to $\outer \uds{\bf SKur}_{\rm rig}^{\rm NC}$. Each of such lifts, denoted by $\hhat {\mb F}$ temporarily, is called an {\bf NC AMS lift} of (either ${\mb F}$ or $\hat {\mb F}$).

\item In Situation \ref{situationm1}, let $\hat {\mb F}_1, \ldots, \hat {\mb F}_m, \hat {\mb F}'$ be post-smoothed AMS lifts of involved Floer flow categories and let $\hat {\mb X}$ be a post-smoothed AMS lift of the multimodule ${\mb X}$. Suppose $\hhat {\mb F}_1, \ldots, \hhat {\mb F}_m; \hhat {\mb F}'$ are NC AMS lifts of $\hat {\mb F}_1, \ldots, \hat {\mb F}_m; \hat {\mb F}'$. Then there exists a class of lifts of $\hat {\mb X}$ to $\outer \uds{\bf SKur}_{\rm rig}^{\rm NC}$ as multimodules over $(\hhat {\mb F}_1, \ldots, \hhat {\mb F}_m; \hhat {\mb F}')$. Each of such lifts, denoted by $\hhat {\mb X}$ temporarily, is called an {\bf NC AMS lift} of (either ${\mb X}$ or $\hat {\mb X}$).

\item In Situation \ref{situationh1}, let $\hat {\mb F}_1, \ldots, \hat {\mb F}_m, \hat {\mb F}, \hat {\mb X}_0, \hat {\mb X}_1, \hat {\mb H}$ be post-smoothed AMS lifts of the involved flow categories, multimodules, and homotopy. Suppose $\hhat {\mb F}_1, \ldots, \hhat {\mb F}_m, \hhat {\mb F}', \hhat {\mb X}_0, \hhat {\mb X}_1$ be NC AMS lifts. Then there exists a lift of $\hat {\mb H}$ to $\outer \uds{\bf SKur}_{\rm rig}^{\rm NC}$ as a homotopy from $\hhat {\mb X}_0$ to $\hhat {\mb X}_1$.
\end{enumerate}
\end{thm}

\begin{proof}
Recall that the post-smoothing of a relatively smooth Kuranishi space $K = (G, V/B, E, S)$ still remembers the vector bundle $\pi^* TB \oplus T^{\rm vt} V$ although the relative structure is destroyed. As the post-smoothing also includes the identification between this vector bundle and the smooth tangent bundle, we may assume $TV$ has the canonical splitting. 

Now notice that $T^{\rm vt} V$ is equipped with the $G$-equivariant NC structure while $TB$ is equipped with the $G$-equivariant stable complex structure. This does not yet provide an NC structure on $TV$ since $TV$ contains the subbundle of the infinitesimal $G$-action. To obtain the NC structure, one needs to make the obstruction bundle complex first. 

Recall that from the construction of the AMS lifts, the obstruction bundle is always trivial and decomposes as 
\beqn
E = \uds W \oplus \uds Q
\eeqn
where $W$ is a complex representation of $G$ and $Q$ is a real representation of $G$. However, we observe that in all the Kuranishi spaces we obtained in the construction, the direct sum 
\beqn 
{\mf g} \oplus Q
\eeqn
is naturally isomorphic to the Lie algebra of the complex Lie group ${\mc G}$. Therefore, we can stabilize all charts by the Lie algebra. Denote by 
\beqn
\tilde K:= {\rm Stab}_{\mf g} (K) = (G, \tilde V, \tilde E, \tilde S).
\eeqn
In this way, the obstruction bundles of all involved Kuranishi spaces are equivariant complex vector bundles. Let
\beqn
{\mf g}_{\tilde V} \subset T \tilde V
\eeqn
be the subbundle spanned by infinitesimal $G$-actions. Then by Definition \ref{defn_NC_Kuranishi}, it remains to construct $G$-equivariant NC structures on the bundles $T \tilde V/ {\mf g}_{\tilde V}$. Indeed, for each chart $\tilde K = (G, \tilde V, \tilde E, \tilde S)$ after the stabilization, let $\pi: \tilde V \to V$ be the projection. Then one has
\beqn
T\tilde V = \pi^* TV \oplus \uds {\mf g}
\eeqn
where the previous constructions provide $TV$ with NC structures. Therefore, it surfaces to define compatible equivariant bundle isomorphisms
\beqn
\pi^* TV \cong (\pi^* TV \oplus \uds {\mf g})/ {\mf g}_{\tilde V}.
\eeqn

One chooses such isomorphisms in the following way. Choose a compatible system of $G$-invariant metrics on $V$ inductively. Then one obtains an equivariant splitting
\beqn
TV \cong {\mf g}_V^\bot \oplus {\mf g}_V.
\eeqn
Notice that ${\mf g}_V$ is isomorphic with the trivial bundle $\uds {\mf g}$. Then define
\beqn
\pi^* TV \cong \pi^* ({\mf g}_V^\bot) \oplus \pi^* {\mf g}_V \to \pi^*(TV \oplus \uds {\mf g})
\eeqn
by sending $\pi^* {\mf g}_V$ into $\pi^* TV$ while sending $\pi^* {\mf g}_V$ onto $\uds{\mf g}$. It is then easy to see that the composition
\beqn
\pi^* TV \to \pi^* (TV \oplus \uds {\mf g}) \to \pi^* (TV \oplus \uds {\mf g}) / {\mf g}_{\tilde V}
\eeqn
is an isomorphism. This equips the Kuranishi spaces $G$-equivariant NC structures. 
\end{proof}

\subsection{Coherent orientations}

For all involved flow categories and all objects $p \in {\mb Ob}{\mb F}$, choose a trivialization of the determinant line $\det D_p$. In last section, the isomorphism $\sigma_{pq}^{\mb F}$ hence determines an orientation on the vertical tangent bundle $T^{\rm vt} V_{pq}^{\mb F}$. Similarly, there are orientations on the vertical tangent bundle $T^{\rm vt} V_{p_1\cdots p_m; p'}^{\mb X}$ for a multimodule. 

After stable smoothing, one has equivariant bundle isomorphisms
\beqn
TV_{pq} \cong \pi_{pq}^* TB_{pq} \oplus T^{\rm vt} V_{pq}.
\eeqn
Both summands have been given certain orientations. We then orient $V_{pq}$ using the above decomposition.

One can then verify the coherence requirement of Definition \ref{defn_coherent_orientation}. Indeed, for $p<r<q$, the isomorphism
\beqn
\Big( \pi_{pr}^* TB_{pr}^{\mb F} \oplus T^{\rm vt} V_{pr}^{\mb F} \Big) \oplus \Big( \pi_{rq}^* T B_{rq}^{\mb F} \oplus T^{\rm vt} V_{rq}^{\mb F} \Big) \to \pi^* TB_{pq}^{\mb F} \oplus T^{\rm vt} V_{pq}^{\mb F}
\eeqn
is orientation preserving up to a sign
\beqn
(-1)^{{\rm rank} T^{\rm vt} V_{rq}} = (-1)^{|p| - |r|}
\eeqn
because $B_{rq}^{\mb F}$ is odd-dimensional.

Suppose the AMS lifts of ${\mb F}_1, \ldots, {\mb F}_m; {\mb F}'$ are oriented as above. W orient the multimodule ${\mb X}$ as follows. For each tuple $(p_1, \ldots, p_m; p')$, orient the tangent bundle of the thickening $V_{p_1 \cdots p_m; p'}^{\mb X}$ via the direct sum decomposition
\beqn
TV_{p_1 \cdots p_m;p'}^{\mb X} \cong \pi_{p_1 \cdots p_m; p'}^* TB_{p_1\cdots p_m; p'}^{\mb F} \oplus T^{\rm vt} V_{p_1 \cdots p_m; p'}^{\mb X}.
\eeqn
Notice that the domain moduli is even-dimensional. Then one can see that with respect to breakings, the isomorphism
\beqn
TV_{p_i q_i}^{{\mb F}_i} \oplus TV_{p_1 \cdots p_{i-1} q_i p_{i+1} \cdots p_m; p'}^{\mb X} \to TV_{p_1 \cdots p_m; p'}^{\mb X}
\eeqn
has a sign $(-1)^{(\sum_{j<i} |p_i|)(|p_i|- |q_i|)}$ while the isomorphism
\beqn
T V_{p_1 \cdots p_m; q'}^{\mb X} \oplus TV_{q'p'}^{{\mb F}'} \cong TV_{p_1 \cdots p_m; p'}^{\mb X}
\eeqn
has a sign $(-1)^{{\rm rank} T^{\rm vt} V_{p_1\cdots p_m; q'}^{\mb X}}$.

\section{Proofs of the Main Theorems on AMS Lifts}\label{section_AMS_proofs}

In this section we summarize the previous constructions as formal proofs of the main theorems about the AMS lifts for Floer flow categories, multimodules and homotopies to the category of (rigidified and collared) normally complex derived orbifolds $\outer\uds{\bf dOrb}_{\rm rig}^{\rm NC}$.

\subsection{Proof of Theorem \ref{thma_flow_chart}}

Assume we are in Situation \ref{situationf1}. The construction is illustrated and summarized in Figure \ref{figure_roadmap}. Each step of the roadmap is completed by a concrete theorem or proposition in previous sections. Upon taking the quotient orbifold, we obtain the proof of Theorem \ref{thma_flow_chart} except the orientation part. 

\begin{center}
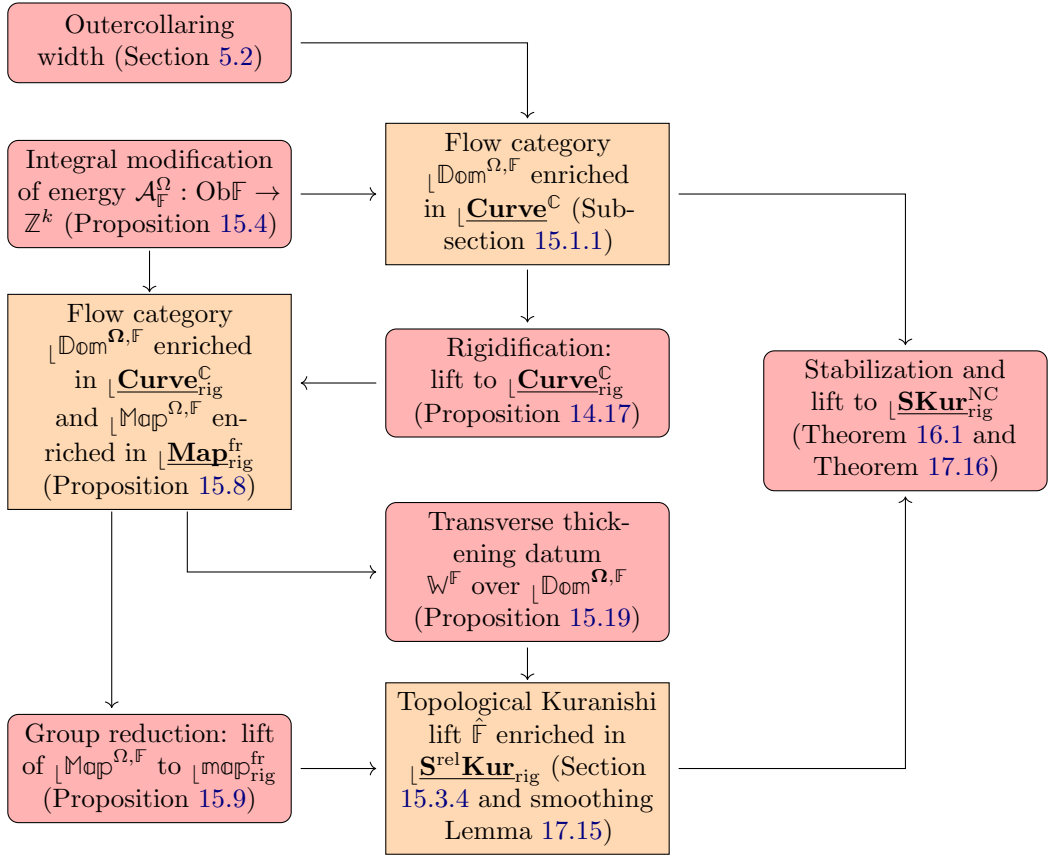
\begin{figure}[h]

\begin{tikzpicture}

\node at (0,0) [choice] {Outercollaring width (Section \ref{sub:outercollaring})};

\node at (0,-2) [choice] {Integral modification of energy ${\mc A}_{\mb F}^\Omega: {\rm Ob}{\mb F} \to {\mb Z}^k$ (Proposition  \ref{prop163})};

\draw [->] (2, 0) -- (5, 0) -- (5, -1);

\draw [->] (2, -2) -- (3, -2);

\node at (5, -2) [result] {Flow category $\outer \mb{Dom}^{\Omega, {\mb F}}$ enriched in $\outer \uds{\bf Curve}^{\mb C}$ (Subsection \ref{subsubsec:pullback-curve})};

\draw [->] (5, -3) -- (5, -3.7);

\node at (5, -4.5) [choice] {Rigidification: lift to $\outer \uds{\bf Curve}_{\rm rig}^{\mb C}$ (Proposition \ref{prop_curve_rigidification})};

\draw [->] (0, -2.75) -- (0, -3.3);

\draw [->] (3, -4.5) -- (2, -4.5);

\node at (0, -4.75) [result] {Flow category $\outer \mb{Dom}^{\bm{\Omega}, {\mb F}}$ enriched in $\outer \uds{\bf Curve}_{\rm rig}^{\mb C}$ and $\outer \mb{Map}^{\Omega, \mb{F}}$ enriched in $\outer \uds{\bf Map}^{\rm fr}_{\rm rig}$ (Proposition \ref{prop:reference-diagram})};

\draw [->] (-0.5, -6.2) -- (-0.5, -8.8);

\draw [->] (0.5, -6.2) -- (0.5, -7) -- (3, -7);

\node at (5, -7) [choice] {Transverse thickening datum $\mb{W}^{\mb F}$ over $\outer \mb{Dom}^{\bm{\Omega}, {\mb F}}$ (Proposition \ref{prop_perturbation})};

\node at (0, -9.6) [choice] {Group reduction: lift of $\outer \mb{Map}^{\Omega, {\mb F}}$ to $\outer \mb{map}_{\rm rig}^{\rm fr}$ (Proposition \ref{prop_group_reduction})};

\draw [->] (2, -9.6) -- (3, -9.6);

\draw [->] (5, -8) -- (5, -8.43);

\node at (5, -9.6) [result] {Topological Kuranishi lift $\hat{\mb F}$ enriched in $\outer \uds{\bf S^{\rm rel}Kur}_{\rm rig}$ (Section \ref{subsubsec:top-AMS} and smoothing Lemma \ref{lemma_postsmoothing})};

\draw [->] (7, -2) -- (10, -2) -- (10, -4);
\draw [->] (7, -9.6) -- (10, -9.6)-- (10, -6);

\node at (10, -5) [choice] {Stabilization and lift to $\outer \uds{\bf SKur}_{\rm rig}^{\rm NC}$ (Theorem \ref{thm171} and Theorem \ref{thm_stable_complex})};

\end{tikzpicture}
\caption{Roadmap of the AMS construction  for Floer flow categories. Here the pink, rounded rectangles involve choices to be made and the orange rectangles are items canonically induced from previous steps.}\label{figure_roadmap}
\end{figure}
\end{center}

Now we explain how to obtain a coherent orientation for the AMS lift. Given a pair of objects $p, q \in {\rm Ob}{\mb F}$, one needs to specify an orientation on the virtual vector bundle
\beqn
T(V_{pq}/G_{pq}) - E_{pq}/G_{pq}.
\eeqn
Recall that $E_{pq} = W_{pq} \oplus Q_{pq}$ where $W_{pq}$ comes from the thickening data and $Q_{pq}$ corresponds to ${\mf g}_{pq}^{\mb C}/ {\mf g}_{pq}$. As the thickening space $W_{pq}$ is complex, it can be neglected for orientation considerations. Moreover, if we stabilize the virtual vector bundle by ${\mf g}_{pq}$, then as $Q_{pq} \oplus {\mf g}_{pq} \cong {\mf g}^{\mb C}$ is complex and hence canonically oriented, it is equivalent to orient the vector bundle
\beqn
T(V_{pq}/G_{pq}) \oplus \uds{\mf g}_{pq}/G_{pq}
\eeqn
where $\uds{\mf g}_{pq} \to V_{pq}$ is the trivial bundle with fibres being ${\mf g}_{pq}$ and $\uds{\mf g}_{pq}/G_{pq}$ is the induced orbifold vector bundle. 

On the other hand, the stable smoothing operation preserves the ($C^0$) isomorphism
\beqn
TV_{pq} \cong \pi_{pq}^* TB_{pq} \oplus T^{\rm vert} V_{pq}
\eeqn
where both summands have been given $G_{pq}$-invariant orientations. Moreover, there is a summand ${\mf g}_{pq} \subset TV_{pq}$ corresponding to the infinitesimal $G_{pq}$-actions. Then the direct sum orientation descends to an orientation on the bundle
\beqn
T(V_{pq}/G_{pq}) \oplus \uds{\mf g}_{pq}/G_{pq}
\eeqn
which is what we need.

It remains to verify that the individual orientations on the derived orbifolds are coherent. Suppose $p<r<q$ in ${\mb F}$. Recall that the domain moduli spaces $B_{pq}$ are oriented in a way which respects breakings (see Lemma \ref{lemma_domain_orientation}). On the other hand, the orientations on the vertical tangent bundles $T^{\rm vt} V_{pq}$ satisfy 
\beqn
\det T^{\rm vt} V_{pq}|_{\partial^{prq} V_{pq}} \cong \det T^{\rm vt} V_{pr} \otimes \det T^{\rm vt} V_{rq}
\eeqn
as oriented bundles. As $B_{pq}$ are always odd-dimensional, there holds
\begin{multline*}
\det \Big( \pi_{pq}^* TB_{pq} \oplus T^{\rm vt} V_{pq}\Big)|_{\partial^{prq} V_{pq}} \cong \det \Big( \pi_{pr}^* TB_{pr} \oplus \pi_{rq}^* TB_{rq} \oplus Q_{prq} \Big) \otimes \det \Big( T^{\rm vt} V_{pr} \oplus T^{\rm vt} V_{rq} \oplus W_{prq} \Big) \\
\cong (-1)^{{\rm rank} T^{\rm vt} V_{pr}} \det \Big( \pi_{pr}^* TB_{pr} \oplus T^{\rm vt} V_{pr} \Big)  \otimes \det \Big( \pi_{rq}^* TB_{rq}\oplus T^{\rm vt} V_{rq}  \Big).
\end{multline*}
Notice that here $Q_{prq}$ and $W_{prq}$ are both complex. The rank of the vertical bundle here is congruent mod 2 to the difference $|p| - |r|$ of cohomological gradings. This means the orientations on the derived orbifolds specified above provide a coherent orientation on $\tilde {\mb F}$.

\subsection{Proof of Theorem \ref{thma_module_lift}}

We also summarize the construction of AMS lifts of multimodules associated to a smooth Floer domain.

\begin{prop}\label{prop181}
In Situation \ref{situationm1}, fix an outercollaring width and choose a list of integral actions on the involved flow categories which are compatible with respect to ${\mb X}$. Let $\hat {\mb F}_1, \ldots, \hat{\mb F}_m; \hat {\mb F}'$ be AMS lifts to $\outer \uds{\bf SKur}_{\rm rig}^{\rm NC}$ (subject to the fixed outercollaring width and the integral actions). Then there exists a lift of the outercollaring of ${\mb X}$ to $\outer \uds{\bf SKur}_{\rm rig}^{\rm NC}$ as a multimodule over $(\hat{\mb F}_1, \ldots, \hat{\mb F}_m; \hat{\mb F}')$.
\end{prop}

\begin{proof}
To save notations, we write the proof in the case of bimodules and write ${\mb X} = {\mb B}$. The case of multimodules is the same. Abbreviate ${\mb F}_1$ by ${\mb F}$. The constructions of $\hat {\mb F}$ and $\hat{\mb F}'$ start with the monotone flow categories $\outer \mb{Dom}$ and $\outer \mb{Dom}'$. Then by Proposition \ref{prop_domain_multimodule}, there exists a monotone bimodule $\outer \mb{Dom}^{\Sigma^{\mb X}}$ over $(\outer \mb{Dom}; \outer \mb{Dom}')$. Then using the compatible integral actions, one obtains the corresponding bimodule $\outer \mb{Dom}^{\Omega, {\mb X}}$ over $(\outer \mb{Dom}^{\Omega, {\mb F}}; \outer \mb{Dom}^{\Omega, {\mb F}'})$. The choices of transverse thickening data for $\hat {\mb F}$ and $\hat {\mb F}'$ can be extended, by Proposition \ref{prop_perturbation}, to a transverse thickening datum ${\mb W}^{\mb X}$. The choices of the group reductions for $\hat {\mb F}$ and $\hat {\mb F}'$ can be extended, by Proposition \ref{prop_group_reduction}, to a group reduction on $\outer \mb{Map}_{{\mb X}}^{\mb{fr}}$. Then by Definition \ref{defn_module_thickening}, one obtains the collections of thickenings $V_{p;p'}^{\mb B}$ which provides (in a canonical way) a lift of the outercollaring of ${\mb B}$ to $\outer \uds{\bf S^{\rm rel} Kur}_{\rm rig}$ as a bimodule over $(\hat {\mb F}; \hat{\mb F}')$. The remaining structures, including the relative NC structure, the stable smoothing, and the NC structure can all be constructed using the inductive procedure. We omit the details. 
\end{proof}

The two exceptional cases, i.e., the multimodules for Poincar\'e pairing and the cigar bimodule can be established similarly. We only give the statement without detailed proofs.

\begin{prop}\label{prop182}
Let ${\mb F}_1, {\mb F}_2$ be Floer flow categories associated to $(H_1, J_1)$ and $(H_2, J_2)$ respectively. Let ${\mb X}^{\rm PD}$ be the Poincar\'e multimodule over $({\mb F}_1, {\mb F}_2; {\mb O})$ associated to a Floer data interpolating on the infinite cylinder. Fix an outercollaring width and choose integral actions ${\mc A}_{{\mb F}_1}^\Omega$, ${\mc A}_{{\mb F}_2}^\Omega$ which are compatible with respect to this multimodule. Let $\hat {\mb F}_1$ resp. $\hat {\mb F}_2$ be AMS lifts of ${\mb F}_1$ resp. ${\mb F}_2$ to $\outer \uds{\bf SKur}_{\rm rig}^{\rm NC}$ subject to the outercollaring width and the chosen integral actions. Then there exists a lift of the outercollaring of ${\mb X}^{\rm PD}$ to $\outer \uds{\bf SKur}_{\rm rig}^{\rm NC}$ as a multimodule over $(\hat {\mb F}_1, \hat {\mb F}_2; {\mb O})$, which we call an AMS lift subject to the outercollaring width and the chosen integral actions.
\end{prop}

\begin{prop}\label{prop183}
Let ${\mb F}$ be a Floer flow category associated to $(H,J)$. Let ${\mb B}^{\rm cigar}$ be the cigar bimodule over $({\mb O}; {\mb F})$ associated to a Floer data on the complex plane. Let $\hat {\mb F}$ be an AMS lift of ${\mb F}$ subject to an outercollaring width and an integral action ${\mc A}_{\mb F}^\Omega$ which is compatible with respect to the bimodule ${\mb B}^{\rm cigar}$. Then there exists a lift of the outercollaring of ${\mb B}^{\rm cigar}$ to $\outer \uds{\bf SKur}_{\rm rig}^{\rm NC}$ which is a bimodule over $({\mb O}; \hat {\mb F})$, which we call an AMS lift of ${\mb B}^{\rm cigar}$ subject to the outercollaring width and the chosen integral action. 
\end{prop}

Propositions \ref{prop181}--\ref{prop183} lead to the first item of Theorem \ref{thma_module_lift} (without the orientation claim). To establish the second item, consider the trivial homotopy ${\mb H}$ from ${\mb X}$ to itself. Notice that the chosen integral actions are still compatible with respect to ${\mb H}$. Then for any two AMS lifts $\hat {\mb X}_0$ and $\hat{\mb X}_1$ subject to the same outercollaring width and the same set of integral actions, one can construct a lift of the outercollaring of ${\mb H}$ to $\outer \uds{\bf SKur}_{\rm rig}^{\rm NC}$ by using the corresponding homotopy version of each step of the AMS construction. We omit the details. After passing to derived orbifolds, one obtains the second item of Theorem \ref{thma_module_lift} (without the orientation claim).

To obtain coherent orientations, one can follow the same argument as in the proof of Theorem \ref{thma_flow_chart}. The signs   come  from the properties of the orientations on the domain moduli spaces (see Proposition \ref{prop_domain_multimodule}) as well as the orientations on vertical tangent bundles (see Subsection \ref{subsection165}). 

\subsection{Proof of Theorem \ref{thma_homotopy_lift}}

The proof is simply a different organization of the previous constructions for flow categories and multimodules. Notice that it is important to choose integral actions which are compatible throughout the consideration. Then each stage of the AMS construction, starting from the domain flow category, rigidifications, choosing thickening data and group reductions, constructing relative NC structures, stable smoothing, and smooth NC structures, can all be inductively accomplished. The statement of Theorem \ref{thma_homotopy_lift} then follows without any new difficulties than flow category and multimodule cases.

\section{Concatenation and Gluing---Proof of Theorem \ref{thma_multimodule_gluing}}\label{section_gluing}

In this section we prove Theorem \ref{thm_AMS_concatenation} and Theorem \ref{thma_multimodule_gluing}. Note that, although the statement of Theorem \ref{thm_AMS_concatenation} is on the level of derived orbifolds, the constructions of the rigidification and NC structure require us to go through the AMS construction on the level of Kuranishi spaces. 

\subsection{Concatenation and gluing of domains}

\subsubsection{Concatenations of monotone multimodules}\label{subsection_domain_concatenation}

We make the analgoue of Definition \ref{defn_general_concatenation} for the case of monotone multimodules. Let ${\mb f}_1, \ldots, {\mb f}_m, {\mb g}_1, \ldots, {\mb g}_n, {\mb g}'$ be monotone flow categories enriched in $\uds{\bf C}$. Let ${\mb x}$ resp. ${\mb y}$ be monotone multimodules over $({\mb f}_1, \ldots, {\mb f}_m; {\mb g}_i)$ resp. over $({\mb g}_1, \ldots, {\mb g}_n; {\mb g}')$. For each $d\geq 0$, let 
\begin{align*}
    &\ A_d^{\mb x}\ &\ A_d^{\mb y}
\end{align*}
be the underlying regular posets of the objects $M_d^{\mb x}$ and $M_d^{\mb y}$ respectively. Define 
\beq\label{monotone_concatenation_poset}
A_d^{{\mb x} \circ_i {\mb y}}:= \left( \bigsqcup_{ d' + d'' = d} A_{d'}^{\mb x} \times A_{d''}^{\mb y} \right)/ \sim 
\eeq
where the equivalence relation $\sim$ is similar to the case discussed in Subsection \ref{subsection_concatenation}. More precisely, the equivalence relation is generated by 
\beqn
\xymatrix{  & A_{d_1}^{\mb x} \times A_d \times A_{d_2}^{\mb y} \ar[ld] \ar[rd]  & \\
A_{d_1 + d}^{\mb x} \times A_{d_2}^{\mb y} & & A_{d_1}^{\mb x} \times A_{d + d_2}^{\mb y}}.
\eeqn
Then for each pair $d, e \geq 0$, there is a codimension zero subset
\beqn
A_d^{{\mb x}} \times A_e^{{\mb y}} \cong A_{d; e}^{{\mb x} \circ_i {\mb y}} \subset A_{d+e}^{{\mb x}\circ_i {\mb y}}.
\eeqn

\begin{defn}\label{defn_monotone_concatenation}
A {\bf concatenation} of ${\mb x}$ and ${\mb y}$ at ${\mb g}_i$  consists of the following objects.
\begin{enumerate}
    \item A monotone multimodule denoted by ${\mb x} \circ_i {\mb y}$ over $({\mb g}_1, \ldots, {\mb g}_{i-1}, {\mb f}_1, \ldots, {\mb f}_m, {\mb g}_{i+1}, \ldots, {\mb g}_n; {\mb g}')$ consisting of objects
    \beqn
    M_d^{{\mb x} \circ_i {\mb y}} \in {\rm Ob} \uds{\bf C}
    \eeqn
    together with structural maps. 
    
    \item For $d, e \geq 0$, codimension zero morphisms
    \beqn
    M_{d}^{\mb x} \times M_{e}^{\mb y} \to M_{d+e}^{{\mb x} \circ_i {\mb y}}
    \eeqn
    whose underlying poset map the inclusion $A_d^{{\mb x}} \times A_e^{{\mb y}} \to A_{d+e}^{{\mb x} \circ_i {\mb y}}$.
\end{enumerate}
They are required to satisfy the following conditions. 
\begin{enumerate}
\item The objects $M_d^{{\mb x} \circ_i {\mb y}}$ is stratified by the poset $A_d^{{\mb x} \circ_i {\mb y}}$ defined by \eqref{monotone_concatenation_poset}.

\item For any $d, d_i, e \geq 0$, the following diagram commutes.
\beqn
\xymatrix{   &    M_d^{\mb x} \times M_{d_i}^{{\mb g}_i} \times M_{e}^{\mb y} \ar[ld] \ar[rd]   & \\
M_{d + d_i}^{\mb x} \times M_e^{\mb y} \ar[rd]   & & M_d^{\mb x} \times M_{d_i + e}^{\mb y} \ar[ld] \\
& M_{ d+ d_i + e}^{{\mb x} \circ_i {\mb y}} &  }
\eeqn

\item The following diagram commutes.
\beqn
\xymatrix{ M_{d_j}^{{\mb f}_j} \times M_d^{{\mb x}} \times M_e^{\mb y} \times M_{e'}^{{\mb g}'} \ar[r] \ar[d]   &   M_{d_j + d}^{\mb x} \times M_{e + e'}^{\mb y} \ar[d] \\
  M_{d_j}^{{\mb f}_j} \times M_{d+e}^{{\mb x} \circ_i {\mb y}} \times M_{e'}^{{\mb g}'} \ar[r]  & M_{d_j + d + e + e'}^{{\mb x} \circ_i {\mb y}}}
  \eeqn
\end{enumerate}
\end{defn}

\subsubsection{Concatenations of domains}

Now we consider the concatenation of two monotone multimodules associated to smooth curves.

\begin{prop}\label{prop_domain_concatenation}
In Situation \ref{situationc1}, let 
\beqn
\outer \mb{Dom}^{{\mb F}_1},  \ldots, \outer \mb{Dom}^{{\mb F}_m},\ \outer \mb{Dom}^{{\mb G}_1}, \ldots, \outer \mb{Dom}^{{\mb G}_n},\ \outer \mb{Dom}^{{\mb G}'}
\eeqn
be the monotone flow categories enriched in $\outer \uds{\bf Curve}_{\rm rig}^{\mb C}$ subject to the same outercollaring width. Let $\outer \mb{Dom}^{\Sigma^{\mb X}}$ resp. $\outer \mb{Dom}^{\Sigma^{{\mb Y}}}$ be monotone multimodules over $(\outer \mb{Dom}^{{\mb F}_1}, \ldots, \outer \mb{Dom}^{{\mb F}_m}; \outer \mb{Dom}^{{\mb G}_i})$ resp. over $(\outer \mb{Dom}^{{\mb G}_1}, \ldots, \outer \mb{Dom}^{{\mb G}_n}; \outer \mb{Dom}^{{\mb G}'})$ constructed in Subsection \ref{subsection_domain_multimodule}. Then there exists a concatenation $\outer \mb{Dom}^{\Sigma^{\mb X}} \circ_i \outer \mb{Dom}^{\Sigma^{{\mb Y}}}$ (see Definition \ref{defn_monotone_concatenation}). 
\end{prop}

\begin{proof}
The underlying monotone multimodules enriched in $\uds{\bf Curve}$ have a canonical concatenation. Indeed, let $\Sigma^{{\mb X} \circ_i {\mb Y}}$ be the singular curve with two components $\Sigma^{{\mb X}}$ and $\Sigma^{{\mb Y}}$ glued along the cylindrical end of the concatenation. Notice that the lateral lines in $\Sigma^{\mb X}$ and $\Sigma^{\mb Y}$ glue together to lateral lines in the singular domain. Then the same construction of $\mb{Dom}^{\Sigma}$ can be applied to $\Sigma^{{\mb X}\circ_i {\mb Y}}$, providing a monotone multimodule $\mb{Dom}^{\Sigma^{\mb X}} \circ_i \mb{Dom}^{\Sigma^{\mb Y}}$. It is straightforward to check that the outercollaring of $\mb{Dom}^{\Sigma^{\mb X}} \circ_i \mb{Dom}^{\Sigma^{\mb Y}}$ is also a concatenation of the corresponding outercollarings of $\mb{Dom}^{\Sigma^{\mb X}}$ and $\mb{Dom}^{\Sigma^{\mb Y}}$. 

The remaining is to construct the rigidifications and the stable complex structures. For simplicity, we assume that ${\mb X}$ is a bimodule $\mb{B}$ over $({\mb F}; {\mb F}')$ and ${\mb Y}$ is a bimodule ${\mb B}'$ over $({\mb F}'; {\mb F}'')$; the general case is essentially the same with only additional notational complexities. 

We first consider the stable complex structure on the concatenation. Recall the case of smooth domains is obtained in Proposition \ref{prop_domain_stable_complex} where the stable complex structure is naturally induced. The same fact holds for the concatenated domain. Therefore, the natural concatenation $\mb{Dom}^{{\mb B}} \circ \mb{Dom}^{{\mb B}'}$ is in fact in $\uds{\bf Curve}^{\mb C}$. 

Now we consider rigidifications. We first rigidify the embeddings
\beqn
{\mc C}_{d}^{\mb{B}} \times {\mc C}_{d'}^{\mb{B}'} \to {\mc C}_{d + d'}^{{\mb B} \circ {\mb B}'}.
\eeqn
Notice that this is a codimension zero embedding and there is a target top stratum 
\beqn
{\mc C}_{d; d'}^{{\mb B} \circ {\mb B}'} \subset {\mc C}_{d + d'}^{{\mb B} \circ {\mb B}'}
\eeqn
such that 
\beqn
{\mc G}_{d + d'} \times_{{\mc G}_d \times {\mc G}_{d'}} ({\mc C}_{d}^{\mb{B}} \times {\mc C}_{d'}^{\mb{B}'}) \cong {\mc C}_{d; d'}^{\mb{B}\circ \mb{B}'}.
\eeqn 
Recall that a rigidification of the above codimension zero embedding (see Definition \ref{defn_curve_rigidification}) is a map 
\beqn
\theta_{d; d'}: B_{d}^{\mb{B}} \times B_{d'}^{\mb{B}'} \times Q_{d; d'}^\epsilon \to B_{d; d'}^{\mb{B} \circ \mb{B}'}.
\eeqn
These maps have to incorporate the existing one defining rigidifications of $\outer \mb{Dom}^{\mb{B} }$ and $\outer \mb{Dom}^{\mb{B}'}$. One can construct these maps inductively, from the deepest strata up. 

First, if $B_{d; d'}^{\mb{B} \circ \mb{B}' }$ has no lower strata, then there is no constraint coming from rigidifications of $\outer \mb{Dom}^{\mb B}$ and $\outer\mb{Dom}^{{\mb B}'}$ and one can choose the rigidification arbitrarily. 

In general one needs to construct $\theta_{d; d'}$ inductively. We only explain how to construct it when $B_{d; d' }^{\mb{B} \circ \mb{B}' }$ has only a codimension one stratum given by a decomposition $d = d_1 + d_2$. This codimension one stratum is shared with another open stratum $B_{d_1; d_2 + d'}^{\mb{B} \circ \mb{B}' }$. The rigidifications for these two open strata (with different symmetry groups) must be compatible along the shared fake boundary. We first choose rigidification of the codimension 1 embedding
\beqn
B_{d_1}^{\mb B}\times B_{d_2}^{\mb B}\times B_{d'}^{{\mb B}'} \to \partial^{d_1; d_2; d'} B_{d + d'}^{{\mb B} \circ {\mb B}'}.
\eeqn
Notice that the embedding can be factorized in two ways
\beqn
\xymatrix{  &   \partial^{d_1 d_2} B_d^{{\mb B}} \times B_{d'}^{{\mb B}'}  \ar[rd] &  \\
B_{d_1}^{\mb B}\times B_{d_2}^{\mb B}\times B_{d'}^{{\mb B}'}  \ar[ru] \ar[rd] &   &  \partial^{d_1; d_2; d'} B_{d + d'}^{{\mb B} \circ {\mb B}'}  \\
&   B_{d_1}^{{\mb B}} \times \partial^{d_2 d'} B_{d_2 + d'}^{{\mb B}'}  \ar[ru] & 
}
\eeqn
where the two arrows on the left are already rigidified. Consider
\beqn
Q_{d_1; d_2; d'} = Q_{d_1 + d_2 + d'}/ Q_{d_1} \oplus Q_{d_2} \oplus Q_{d'}.
\eeqn
The rigidification on this lower stratum should be a germ of maps
\beqn
\theta_{d_1; d_2; d'}: \Big( B_{d_1}^{\mb B}\times B_{d_2}^{\mb B}\times B_{d'}^{{\mb B}'} \Big) \times Q_{d_1; d_2; d'}  \to \partial^{d_1; d_2; d'} B_{d + d'}^{{\mb B} \circ {\mb B}'} 
\eeqn
Notice that there is a natural inclusion
\beqn
Q_{d_1 d_2} \oplus Q_{d_2 d'} \subset Q_{d_1; d_2; d'}
\eeqn
where the value of the rigidification is already determined. One then only has to extend it to the complement. Choosing such an extension equivariantly, one hence obtains the rigidification $\theta_{d_1; d_2; d'}$ on this lower stratum. 

Notice that this lower stratum is shared by two top strata. One uses the collar regions in both directions to extend to a rigidification for both separately. 
\end{proof}

\subsubsection{Gluing of domains}

Now suppose we are in Situation \ref{situationh1}. Suppose on the family of domains $\Sigma_t$ ($t \in [0, 1]$) one chooses a smooth family of lateral lines. Here $t = 0$ corresponds to the broken domain. Then without making further choices, one has the monotone flow categories
\beqn
\mb{Dom}^{\mb{F}_1}, \ldots, \mb{Dom}^{\mb{F}_m}, \mb{Dom}^{{\mb G}_1}, \ldots, \mb{Dom}^{\mb{G}_n}, \mb{Dom}^{\mb{G}'}
\eeqn
monotone multimodules
\begin{align*}
&\ \mb{Dom}^{\Sigma^{\mb X}},\ &\ \mb{Dom}^{\Sigma^{\mb Y}}.
\end{align*}
As demonstrated above, the concatenation $\mb{Dom}^{\Sigma^{\mb X}} \circ_i \mb{Dom}^{\Sigma^{\mb Y}}$ within the category $\uds{\bf Curve}$ is canonically defined. Moreover, by varying $t$, one also obtains a monotone homotopy $\mb{Dom}^{\mb H}$ from the concatenation to $\mb{Dom}^{\mb{M}}$ which is the multimodule associated to $\Sigma_1$. 

What we do below is the construction of the further choices such as rigidifications and stable complex structure. 

\begin{prop}\label{prop_domain_gluing}
Fix an outercollaring width. Let $\outer \mb{Dom}^{\mb{F}_1}, \ldots, \outer \mb{Dom}^{\mb{F}_m}$, $\outer \mb{Dom}^{\mb{G}_1}, \ldots, \outer \mb{Dom}^{\mb{G}_n}$, $\outer \mb{Dom}^{\mb{G}'}$, $\outer \mb{Dom}^{\Sigma^{\mb X}}$, $\outer \mb{Dom}^{\Sigma^{\mb Y}}$, $\outer \mb{Dom}^{\Sigma^{\mb X}} \circ_i \outer \mb{Dom}^{\Sigma^{\mb Y}}$ be as in Proposition \ref{prop_domain_concatenation}. Let $\outer \mb{Dom}^{\mb{M}}$ be a monotone multimodule associated to the smooth domain $\Sigma_1$ enriched in $\outer \uds{\bf Curve}_{\rm rig}^{\mb C}$. Then for the outercollaring $\outer \mb{Dom}^{\mb H}$, there exists a lift to $\outer \uds{\bf Curve}_{\rm rig}^{\mb C}$ as a monotone homotopy from $\outer \mb{Dom}^{\Sigma^{\mb X}} \circ_i \outer \mb{Dom}^{\Sigma^{\mb Y}}$ to $\outer \mb{Dom}^{\mb{M}}$.
\end{prop}

\begin{proof}
Again, it only remains to construct the rigidifications and stable complex structures. The construction of rigidification is still based on the same induction principle. Then stable complex structure can be chosen after choosing a rigidification.
\end{proof}

\subsection{Relatively smooth AMS lifts for concatenations and gluing}

In this subsection we treat the topological AMS construction for concatenations of Floer multimodules associated to smooth Floer domains.

\subsubsection{Concatenation}

Assume we are in Situation \ref{situationc1}. We have constructed relatively smooth AMS lifts $\hat{\mb X}$ and $\hat{\mb Y}$ to $\outer \uds{\bf S^{\rm rel} Kur}_{\rm rig}$ respectively. The purpose of this subsection is to demonstrate how to extend the constructions of $\hat{\mb X}$ and $\hat{\mb Y}$ to a lift of the concatenation ${\mb X} \circ_i {\mb Y}$ to $\outer \uds{\bf S^{\rm rel} Kur}_{\rm rig}$.

\begin{thm}\label{thm_AMS_concatenation_1}
In Situation \ref{situationc1}, let $\hat {\mb X}$ and $\hat {\mb Y}$ be relatively smooth AMS lifts of the outercollarings of ${\mb X}$ and ${\mb Y}$ respectively, subject to the same outercollaring width and integral actions which induce the same integral symplectic form on $X$. Then there exists a set of concatenations $\hat {\mb X} \circ_i \hat {\mb Y}$ enriched in $\outer \uds{\bf S^{\rm rel} Kur}_{\rm rig}$.
\end{thm}

\begin{proof}
As in Subsection \ref{subsection_domain_concatenation}, we only describe the construction in the case of concatenations of bimodules. The case of multimodules differs only in notation. Denote ${\mb X} = {\mb B}_{01}$ and ${\mb Y} = {\mb B}_{12}$ where ${\mb B}_{01}$ is a bimodule over $({\mb F}_0; {\mb F}_1)$ and $\mb{B}_{12}$ is a bimodule over $({\mb F}_1; \mb{F}_2)$, where ${\mb F}_i$ is the Floer flow category associated to Hamiltonian $H_i$. Let $\Sigma$ be the broken domain with two cylinders equipped with Floer datum $\sigma$ from the concatenation of the two Floer data for $\mb{B}_{01}$ and $\mb{B}_{12}$. 

First consider the concatenation of domains. The AMS constructions of $\hat {\mb F}_0, \hat {\mb F}_1, \hat {\mb F}_2$ and $\hat {\mb B}_{01}$, $\hat {\mb B}_{12}$ start with choosing monotone flow categories $\outer \mb{Dom}_0, \outer \mb{Dom}_1, \outer \mb{Dom}_2$ enriched in $\outer \uds{\bf Curve}_{\rm rig}^{\mb C}$ as considered in Section \ref{section_AMS_domains}, as well as monotone bimodules $\outer \mb{Dom}_{01}, \outer \mb{Dom}_{12}$ respectively. Then Proposition \ref{prop_domain_concatenation} allows us to choose a concatenation $\outer \mb{Dom}_{01} \circ \outer \mb{Dom}_{12}$ as a monotone bimodule over $(\outer \mb{Dom}_0; \outer \mb{Dom}_2)$. 

Now using the chosen integral actions, one can obtain pullbacks 
\begin{align*}
&\ \outer \mb{Dom}^{{\mb F}_i}, i = 0, 1, 2;\ &\ \outer \mb{Dom}^{{\mb B}_{01}},\ \outer \mb{Dom}^{{\mb B}_{12}}
\end{align*}
as ordinary flow categories and bimodules. The corresponding bimodules of framed maps are also obtained
\begin{align*}
&\ \outer \mb{Map}_{\mb{B}_{01}}^{\rm fr},\ &\ \outer \mb{Map}_{\mb{B}_{12}}^{\rm fr}
\end{align*}
with a concatenation 
\beqn
\outer \mb{Map}_{\mb{B}_{01}}^{\rm fr} \circ \outer \mb{Map}_{\mb{B}_{12}}^{\rm fr}
\eeqn
canonically determined by the concatenation $\outer\mb{Dom}^{{\mb B}_{01}}\circ \outer \mb{Dom}^{\mb{B}_{12}}$.

We need to show that the two further systems of choices, i.e., group reductions and thickening data, can be extended to the concatenation. 

We first explain why the group reduction can be naturally extended. Indeed, given a triple of objects $p_0, p_1, p_2$ of the three involved flow categories, one has the morphism
\beq\label{eqn165}
\outer {\rm Map}_{p_0; p_1}^{{\rm fr}} \times \outer {\rm Map}_{p_1; p_2}^{{\rm fr}} \to \outer {\rm Map}_{p_0; p_2}^{{\rm fr}}
\eeq
in the category of framed maps. One can see that if we choose any strict embedding representing the conjugacy class, then the strict embedding induces an open embedding onto the closure of a top stratum
\beqn
\hat {\mc G}_{p_0; p_2} \times_{\hat {\mc G}_{p_0; p_1}\times \hat {\mc G}_{p_1; p_2}} \Big( \outer {\rm Map}_{p_0; p_1}^{{\rm fr}} \times \outer {\rm Map}_{p_1; p_2}^{{\rm fr}} \Big) \to \outer {\rm Map}_{p_0; p_2}^{{\rm fr}}
\eeqn
(recall $\hat {\mc G} = {\mc G}_L \times {\mc G}_R$ where ${\mc G}_L$ acts by reparametrizing the maps and frames while ${\mc G}_R$ acts by linearly rearranging the frames). As the group reduction is invariant under the first group action and equivariant under the second action, the open embedding naturally induces an extension of the group reduction. One can check that the extension is well-defined on fake boundaries of $\outer {\rm Map}_{p_0; p_2}^{\rm fr}$ (which are adjacent to two such top strata). 

The extension of the thickening datum is also canonical in a neighborhood of the morphism \eqref{eqn165}, because the embedding is rigidified. Using the rigidification of this embedding, we can extend the thickening datum in a way which is constant in each normal fiber specified by the rigidification. As one can always shrink the resulting charts, up to shrinking, the extension of the thickening datum is canonically extended. 

Then one can write down the thickenings in the same way as for a multimodule with a smooth domain. We can verify easily that one obtains a topological Kuranishi lift $\hat {\mb B}$ which is also a concatenation $\hat {\mb B}_{01} \circ \hat {\mb B}_{12}$ within the category $\outer \uds{\bf Kur}_{\rm rig}$. The relative smooth structure is canonically given and one can also easily see that the concatenation is indeed in $\outer \uds{\bf S^{\rm rel} Kur}_{\rm rig}$.
\end{proof}

\subsubsection{Gluing}

\begin{thm}\label{thm_AMS_gluing_1}
In Situation \ref{situationg1}, fix an outercollaring width and integral actions on the involved flow categories compatible with respect to ${\mb H}$. Let $\hat{\mb X}$ and $\hat{\mb Y}$ be relatively smooth AMS lifts of ${\mb X}$ and ${\mb Y}$ subject to the outercollaring width and the integral actions. Let $\hat{\mb X} \circ \hat{\mb Y}$ be a concatenation provided by Theorem \ref{thm_AMS_concatenation_1}. Let $\hat{\mb M}$ be an AMS lift of the multimodule ${\mb M}$. Then there exists a set of lifts of the outercollaring of the homotopy ${\mb H}$ to $\outer \uds{\bf S^{\rm rel} Kur}_{\rm rig}$ as a homotopy from $\hat{\mb X} \circ \hat {\mb Y}$ to $\hat {\mb M}$.    
\end{thm}

\begin{proof}
At this stage, it follows from the same strategy as the proof of Theorem \ref{thma_homotopy_lift}. The choices of rigidifications of domains are provided by Proposition \ref{prop_domain_gluing}. The choices of group reductions and thickening datum can also be interpolated over the homotopy. The relative smooth structures are canonically obtained.  
\end{proof}

\subsection{Normal complex structure and stable smoothing}

Finally, we sketch the discussion on normal complex structures and stable smoothing of the objects arising in concatenation and gluing, which will finish the proof of Theorem \ref{thm_AMS_concatenation} and Theorem \ref{thma_multimodule_gluing}.

\begin{proof}[Proof of Theorem \ref{thm_AMS_concatenation}]
We continue after Theorem \ref{thm_AMS_concatenation_1}. Suppose the relatively smooth AMS lifts $\hat {\mb F}_1, \ldots, \hat {\mb F}_m$,  $\hat{\mb G}_1, \ldots, \hat {\mb G}_n$, $\hat {\mb G}'$, and the multimodules $\hat {\mb X}$, $\hat {\mb Y}$ are lifted to $\outer \uds{\bf S^{\rm rel} Kur}_{\rm rig}^{\rm NC}$ using the construction of Section \ref{section_NC_structure}. We would like to see if the composition $\hat {\mb X} \circ_i \hat {\mb Y}$ provided by Theorem \ref{thm_AMS_concatenation_1} can also be lifted to $\outer \uds{\bf S^{\rm rel} Kur}_{\rm rig}^{\rm NC}$ as a composition. In fact, the vertical tangent bundles contained in the composition $\hat {\mb X} \circ_i \hat {\mb Y}$ are canonically determined by the corresponding vertical tangent bundles. Hence the vertical NC structures are canonically induced. Therefore, we see that $\hat {\mb X} \circ_i \hat {\mb Y}$ is lifted to $\outer \uds{\bf S^{\rm rel} Kur}_{\rm rig}^{\rm NC}$. 

The stable smoothing of the concatenation is also automatically induced. We briefly explain for the bimodule case. Let ${\mb F}_0, {\mb F}_1, {\mb F}_2$ and $\mb{B}_{01}$, $\mb{B}_{12}$ be as in the proof of Theorem \ref{thm_AMS_concatenation_1}. For each triple of objects $p_0, p_1, p_2$, one has the embedding
\beqn
K_{p_0; p_1}^{\mb{B}_{01}} \times K_{p_1; p_2}^{\mb{B}_{12}} \to K_{p_0; p_2}^{\mb{B}_{01} \circ \mb{B}_{12}}.
\eeqn
The difference essentially lies in the difference between the complex Lie groups ${\mc G}_{p_0; p_1}\times {\mc G}_{p_1; p_2}$ and ${\mc G}_{p_0; p_2}$. In terms of Kuranishi charts, the unitary part of the difference is realized by changing the symmetry groups while the non-unitary part of the difference is realized by having a further stabilization. Therefore, if both $K_{p_0; p_1}^{\mb{B}_{01}}$ and $K_{p_1; p_2}^{\mb{B}_{12}}$ is equipped with a (stable) smoothing, then the corresponding top stratum of $K_{p_0; p_2}^{\mb{B}_{01} \circ \mb{B}_{12}}$ has an induced stable smoothing (possibly after doing further stabilizations). 

Lastly, we combine the stable complex structures of the domain moduli space discussed in the proof of Proposition \ref{prop_domain_concatenation}, the NC structures on the vertical tangent bundles to obtain an NC structure on the quotient derived orbifolds.
\end{proof}

\begin{proof}[Proof of Theorem \ref{thma_multimodule_gluing}]
We continue after Theorem \ref{thm_AMS_gluing_1}. Suppose the relatively smooth AMS lifts $\hat {\mb F}_1, \ldots, \hat {\mb F}_m$, $\hat {\mb G}_1, \ldots, \hat {\mb G}_n$, $\hat {\mb G}'$, $\hat {\mb X}$, $\hat {\mb Y}$, $\hat {\mb M}$ are lifted to $\outer \uds{\bf S^{\rm rel} Kur}_{\rm rig}^{\rm NC}$. Suppose we also have a concatenation $\hat {\mb X} \circ_i \hat {\mb Y}$ enriched in $\outer \uds{\bf S^{\rm rel} Kur}_{\rm rig}^{\rm NC}$. The homotopy $\hat{\mb H}$ is also lifted to $\outer \uds{\bf S^{\rm rel} Kur}_{\rm rig}$. We first can lift $\hat {\mb H}$ to $\outer \uds{\bf S^{\rm rel} Kur}_{\rm rig}^{\rm NC}$. This can be done by using the same argument as using a family of Cauchy--Riemann operators interpolating between the linearization of the Floer equation and a complex-linear one. We omit the details, which would be a tedious repetition of the flow category and multimodule case already considered. 

Next, suppose $\hat {\mb F}_1, \ldots, \hat {\mb F}_m$, $\hat {\mb G}_1, \ldots, \hat {\mb G}_n$, $\hat {\mb G}'$, $\hat {\mb X}$, $\hat {\mb Y}$, $\hat {\mb M}$ are also lifted to $\outer \uds{\bf SKur}_{\rm rig}^{\rm NC}$ after further stabilization and smoothing. The stable smoothing and NC structures of $\hat {\mb H}$ can be constructed using the same inductive argument as we always did. The result is obtained by passing to the orbifold quotient.
\end{proof}

\section{Comparing Non-homotopic Choices}\label{section_comparison}

The comparison between different choices which are not homotopic requires more setup. Choices in this category include
\begin{enumerate}
    \item The choice of integral actions (see Definition \ref{defn_integral_modification}) in construction the AMS lifts.

    \item Different choices of lateral lines in Floer domains.

\end{enumerate}

We remove the dependence on these two choices uniformly via a doubling construction (cf. \cite[Section 6.10]{AMS}). Here, different integral actions give rise to different moduli spaces of curves of projective spaces; as for different choices of lateral lines under the same choice of integral actions, they induce different composition laws. Despite the different effects, we can handle them in the same way.

\subsection{Curves in products of projective spaces}

We can generalize the notion of monotone flow categories/multimodules/homotopies (see Definition \ref{defn_monotone_flow}, Definition \ref{defn_monotone_multimodule}, Definition \ref{defn_monotone_homotopy}) in the following way. Fix $k \geq 2$ and consider the monoid ${\mb Z}_{\geq 0}^k$ instead of ${\mb Z}_{\geq 0}$. Then for each ${\bf d} = (d^1, \ldots, d^k) \in {\mb Z}_{\geq 0}^k$, there is a natural poset
\beqn
A_{\bf d} = \Big\{ ({\bf d}_0, \ldots, {\bf d}_l)\ |\ {\bf d}_i \in {\mb Z}_{\geq 0}^k \setminus \{0\},\ {\bf d} = {\bf d}_0 + \cdots + {\bf d}_l \Big\}.
\eeqn
A {\bf monotone flow category of rank $k$} enriched in $\uds {\bf C}$, denoted by ${\mb f}$, consists of a collection of $A_{\bf d}$-stratified objects
\beqn
M_{\bf d}^{\mb f},\ \forall {\bf d} \in {\mb Z}_{\geq 0}^k
\eeqn
and a collection of morphisms
\beqn
M_{{\bf d}_1}^{\mb f} \times M_{{\bf d}_2}^{\mb f} \to \partial^{{\bf d}_1, {\bf d}_2} M_{{\bf d}_1 + {\bf d}_2}^{\mb f}
\eeqn
which satisfy the natural associativity condition, similar to that imposed in Definition \ref{defn_monotone_flow}. One can see any monotone flow category ${\mb f}$ of rank $k$ induces an ordinary flow category ${\mb F}$ with ${\rm Ob}{\mb F} \cong {\mb Z}^k$. Similarly, one can define monotone multimodules and homotopies of higher ranks as we did in Definition \ref{defn_monotone_multimodule} and Definition \ref{defn_monotone_homotopy}: we leave it to the interested reader to spell out the details.

Now we would like to build rank $k$ monotone flow category $\outer \mb{Dom}^k$ enriched in $\outer \uds{\bf Curve}_{\rm rig}^{\mb C}$. Indeed, for ${\bf d} = (d^1, \ldots, d^k)$, define
\beqn
{\mb P}^{\bf d} = {\mb P}^{d^1} \times \cdots \times {\mb P}^{d^k}.
\eeqn
Each element is an equivalence class of maps with $k$ components. Indeed, one can consider the moduli space $\ov{\mc M}{}_{0, 2}^{\mb R}(\mb{P}^{\bf d}, {\bf d})$ of cylinders in $\mb{P}^{\bf d}$ of multidegree ${\bf d}$ and the open subset 
\beqn
B_{\bf d}  \subset \ov{\mc M}{}_{0, 2}^{\mb R}(\mb{P}^{\bf d}, {\bf d})
\eeqn
of cylinders satisfying 1) the evaluations of the $i$-th component at  $+\infty$ is  $[0, \cdots, 0, 1] \in \mb{P}^{d^i}$ and 2) the image of the $i$-th component is not contained in any hyperplanes of $\mb{P}^{d^i}$. Then $B_{\bf d} $ is an $A_{\bf d}$-stratified manifold with a smooth action by the complex Lie group
\beqn
{\mc G}_{\bf d} := {\mc G}_{d^1}  \times \cdots \times {\mc G}_{d^k} 
\eeqn
which contains a maximal compact subgroup
\beqn
G_{\bf d} := G_{d^1}  \times \cdots \times G_{d^k}  \cong U(d^1) \times \cdots \times U(d^k).
\eeqn
Similar to the rank 1 case, one obtains a collection of objects
\beqn
{\mc C}_{\bf d}    = ({\mc G}_{\bf d} , B_{\bf d} , C_{\bf d} ) \in {\rm Ob} \uds{\bf Curve}
\eeqn
where $C_{\bf d}  \to B_{\bf d} $ is the universal curve. 

One can define similar structural maps to obtain a monotone flow category, denoted by 
\beqn
\mb{Dom}^k.
\eeqn
Similarly, for any smooth Floer domain $\Sigma^{\mb X}$, one obtains a monotone multimodule $\mb{Dom}^{k, \Sigma^{\mb X}}$ over $(\mb{Dom}^k, \ldots, \mb{Dom}^k; \mb{Dom}^k)$.

A slight difference from the previous rank $1$ case is that we choose not to construct rigidifications and stable complex structures at this stage (although we could do so). We would like to do the construction later after pulling back so that certain compatibility holds directly from our construction. At this point, we only take outercollarings and the monotone flow categories/multimodules are enriched in $\outer \uds{\bf Curve}$.

\subsection{Double-framed thickenings}

We may carry out the roadmap of the AMS construction for the higher-rank approach which facilitates the comparison of different integral actions. We only need to consider the $k=2$ case and call such a construction a ``double-framed'' AMS construction.

We briefly go through the construction for Floer flow categories. Assume one is in Situation \ref{situationf1}. Fix an outercollaring width. Let 
\beqn
{\mc A}_{\mb{F}}^{\Omega_1},\ {\mc A}_{\mb{F}}^{\Omega_2}: {\rm Ob}{\mb F} \to {\mb Z}
\eeqn
be two integral actions on ${\mb F}$. Denote their direct sum by
\beqn
{\mc A}_{\mb F}^{\Omega_1 \oplus \Omega_2}: {\rm Ob}{\mb F} \to {\mb Z}^2.
\eeqn
Then for each pair $p<q$ of objects of ${\mb F}$, define
\beqn
{\bf d}_{pq} = (d_{pq}^1, d_{pq}^2) = {\mc A}_{\mb F}^{\Omega_1 \oplus \Omega_2}(p) - {\mc A}_{\mb F}^{\Omega_1 \oplus \Omega_2} (q)\in {\mb Z}_{\geq 0}^2.
\eeqn
Then for each pair $p<q$, let
\beqn
{\mc C}_{pq}  = ({\mc G}_{pq}, B_{pq}, C_{pq}), 
\eeqn
be the restriction of ${\mc C}_{{\bf d}_{pq}}$ to the strata whose associated degrees can be realized by integralized actions of broken Floer trajectories of $H$. This provides a pullback 
\beqn
\outer \mb{Dom}^{\Omega_1 \oplus \Omega_2, {\mb F}}
\eeqn
which is a flow category enriched in $\outer \uds{\bf Curve}$ whose objects are those of ${\mb F}$.

One can write down the associated flow category of framed maps in the same way as the rank $1$ case. More precisely, given objects $p<q$, the space of framed maps consists of tuples $(\phi, u, F)$ where $\phi \in B_{pq}$, $u: C_\phi \to X$ is a smooth map satisfying the exponential decay condition, and $F = (F_1, F_2)$ consists of two frames where $F_i$, called an $\Omega_i$-frame along $u$, is a frame of the holomorphic line bundle $L_{u, \Omega_i} \to C_\phi$ satisfying the positive constraint.

The following fact is important for the comparison. For each $p<q$, write
\beqn
{\mc G}_{pq} = {\mc G}_{pq}^1 \times {\mc G}_{pq}^2 = {\mc G}_{d_{pq}^1} \times {\mc G}_{d_{pq}^2}.
\eeqn

\begin{lemma}\label{lemma_modification_comparison}
$B_{pq}$ is a ${\mc G}_{pq}^1$-equivariant principal ${\mc G}_{pq}^2$-bundle over $B_{pq}^1$.
\end{lemma}

\begin{proof}
Abbreviate
\beqn
d^i = {\mc A}_{\mb{F}}^{\Omega_i}(p) - {\mc A}_{\mb F}^{\Omega_i}(q),\ i = 1, 2.
\eeqn
Each element of $B_{pq}$ is represented by a holomorphic map $(u^1, u^2): C_\phi \to \mb{P}^{d^1} \times \mb{P}^{d^2}$. Then there is an obvious free ${\mc G}_{d^2}$-action on $B_{pq}$ by linearly transforming the second component. Moreover, fix $\phi^1 \in B_{pq}^1$, all holomorphic maps into ${\mb P}^{d^2}$ which have the prescribed degrees on each irreducible component are in the same free ${\mc G}_{pq}^2$-orbit. They come from projective embeddings of the irreducible components into ${\mb P}^{d^2}$ induced by a basis of global sections of the line bundle with prescribed degree, with matching condition on the nodes. This amounts to looking at all framings with positive evaluation constraint, thereby producing ${\mc G}_{pq}^2$.
\end{proof}

We denote the natural projection
\beqn
\pi_i: \outer B_{pq} \to \outer B_{pq}^i,\ i = 1, 2.
\eeqn
We will use these projections to pull back all choices made for the rank 1 construction to the rank 2 case and take direct sums.

\begin{cor}
The lifts of $\outer \mb{Dom}^{\Omega_i, {\mb F}}$ to $\outer \uds{\bf Curve}_{\rm rig}^{\mb C}$ canonically induce a lift of $\outer \mb{Dom}^{\Omega_1 \oplus \Omega_2, {\mb F}}$. 
\end{cor}

\begin{proof}
By Definition \ref{defn_curve_rigidification}, the rigidification for $\outer \mb{Dom}^{\Omega_i, {\mb F}}$ is given by maps
\beqn
\theta_{prq}^i: \outer B_{pr}^i \times \outer B_{rq}^i \times Q_{prq}^i \to {\mc G}_{pq}^i
\eeqn
satisfying certain conditions. One can simply pullback $\theta_{prq}^i$ to $\outer B_{pr}\times \outer B_{rq}$. More precisely, define
\beqn
\begin{split}
\theta_{prq}: \outer B_{pr} \times \outer B_{rq} \to &\ Q_{prq} \\
(\phi_{pr}, \phi_{rq}) \mapsto &\ \Big( \theta_{prq}^1( \pi_1(\phi_{pr}), \pi_1(\phi_{rq})), \theta_{prq}^2( \pi_2(\phi_{pr}), \pi_2(\phi_{rq})) \Big).
\end{split}
\eeqn
These maps define a rigidification of $\outer \mb{Dom}^{\Omega_1 \oplus \Omega_2, {\mb F}}$.

On the other hand, the stable complex structure can also be obtained by using such pullbacks, whose compatibility can be checked easily.
\end{proof}

Now we assume that $\outer \mb{Dom}^{\Omega_1 \oplus \Omega_2, {\mb F}}$ is enriched in $\outer \uds{\bf Curve}_{\rm rig}^{\mb C}$. Of course, such a lift does not necessarily come from pullbacks, but for our purpose, it suffices to consider this special situation.

Similarly, for a flow multimodule ${\mb X}$ over $({\mb F}_1, \ldots, {\mb F}_m; {\mb F})$, if each of the flow categories are equipped with integral actions for $i = 1, 2$,
\beqn
{\mc A}_{{\mb F}_1}^{\Omega_i}, \ldots, {\mc A}_{{\mb F}_m}^{\Omega_i}, {\mc A}_{{\mb F}'}^{\Omega_i}.
\eeqn 
we obtain a multimodule $\outer \mb{Dom}^{\Omega_1 \oplus \Omega_2, {\mb X}}$ over $(\outer \mb{Dom}^{\Omega_1 \oplus \Omega_2, {\mb F}_1}, \ldots, \outer \mb{Dom}^{\Omega_1 \oplus \Omega_2, {\mb F}_m}; \outer \mb{Dom}^{\Omega_1 \oplus \Omega_2, {\mb F}'})$ enriched in $\outer \uds{\bf Curve}_{\rm rig}^{\mb C}$.

\subsection{Flow category situation}\label{subsec:flow-double}
We spell out the comparison between different integral actions on the same flow category in detail.

\subsubsection{Framed maps and group reductions}
We can write down the flow category of framed maps
\beqn
\outer \mb{Map}^{\mb{fr}}_{\mb F}
\eeqn
which assigns to each pair $p<q$ a space $\outer {\rm Map}_{pq}^{\rm fr}$ consisting of triples $(\phi, u, F)$ where $\phi \in \outer B_{pq}$, $u: C_\phi \to X$ is a component-wise smooth map which extends continuously over nodes and which converges exponentially fast to periodic orbits, and $F = (F_1, F_2)$ is a double-frame where $F_i$ is a frame of the holomorphic line bundle $L_{u, \Omega_i} \to C_\phi$ satisfying the pre-chosen positive or negative constraint.

Notice that there are still natural projection maps
\beq\label{projection}
\begin{split}
\pi_i: \outer {\rm Map}_{pq}^{\rm fr} \to &\ \outer {\rm Map}_{pq, i}^{\rm fr}\\
(\phi, u, F) \mapsto & (\pi_i(\phi), u, F_i).
\end{split}
\eeq

Now denote
\beqn
\hat {\mc G}_{pq} = {\mc G}_{L, pq} \times {\mc G}_{R, pq} \cong {\mc G}_{pq} \times {\mc G}_{pq}.
\eeqn
It acts on $\outer {\rm Map}_{pq}^{\rm fr}$ as before, where ${\mc G}_{L, pq}$ reparametrizes the maps and the frames, while ${\mc G}_{R, pq}$ linearly transforms the frames. Recall that a group reduction (Definition \ref{defn_group_reduction}) consists of a compatible system of $\hat {\mc G}_{pq}$-equivariant maps
\beqn
\Lambda_{pq}: \outer {\rm Map}_{pq}^{\rm fr} \to {\mc G}_{pq}/ G_{pq}.
\eeqn
If we have chosen independently group reductions for $\outer \mb{Map}_{{\mb F}_i}^{\mb{fr}}$, then using the projection map \eqref{projection}, we can pull them back to obtain a group reduction on the double-framed construction $\outer \mb{Map}_{\mb{F}}^{\mb{fr}}$.

\subsubsection{Thickening data}

The last step of constructing the topological AMS lift is to choose transverse thickening data. By slightly generalizing Definition \ref{defn_thickening_datum}, a thickening datum on the flow category $\outer \mb{Dom}^{\Omega_1\oplus \Omega_2, {\mb F}}$ consists of a collection of unitary $G_{pq}$-representations and $G_{pq}$-equivariant linear maps 
\beqn
\nu_{pq}: W_{pq} \to \Gamma_c( \mathring C_{pq} \times X, \Lambda_{\mathring C_{pq}/ B_{pq}}^{0,1}\otimes TX).
\eeqn
satisfying the natural compatibility condition with the structural maps. With or without transversality condition, one can write down the thickening (cf. Definition \ref{defn_thickened_moduli}) 
\beqn
V_{pq} = \Big\{ (\phi, u, e, F)\ |\ (\phi, u, F) \in \outer {\rm Map}_{pq}^{\rm fr}, e \in W_{pq},\ \phi_F = \phi \in \outer B_{pq},\ \ov\partial_H u + \nu_{pq}(e) = 0 \Big\}.
\eeqn
The obstruction bundle $E_{pq}$ is still the direct sum of the trivial bundles $W_{pq} \oplus Q_{pq}$ and the Kuranishi section $S_{pq}: V_{pq} \to E_{pq}$ is the direct sum of the (linearized) group reduction $\lambda_{pq}$ and the forgetful map $(\phi, u, e, F) \mapsto e$. There is the natural footprint map
\beqn
S_{pq}^{-1}(0)/ G_{pq} \to \outer M_{pq}^{\mb F}
\eeqn
sending the $G_{pq}$-orbit of $(\phi, u, e, F)$ to the equivalence class of $u: C_\phi \to X$.

Then suppose one has chosen independently thickening datum
\beqn
\nu_{pq}^i: W_{pq}^i \to \Gamma_c( \mathring C_{pq}^i \times X, \Lambda_{\mathring C_{pq}^i / B_{pq, i}}^{0, 1} \otimes TX),\ i = 1, 2,
\eeqn
then one can pullback each of them
\beqn
(\pi_i)^* \nu_{pq}^i: W_{pq}^i \to \Gamma_c( \mathring C_{pq} \times X, \Lambda_{\mathring C_{pq}/ B_{pq}}^{0,1}\otimes TX).
\eeqn

\begin{lemma}\label{lemma193}
If $\nu_{pq}^i$ is transverse, so is $\pi_i^* \nu_{pq}^i$.
\end{lemma}

\begin{proof}
This follows from the observation that if $W_{pq}^i$ is mapped surjectively onto the cokernel of the linearized Cauchy--Riemann operator, so is the pullback.
\end{proof}

\subsubsection{Finalizing proof}

Now we would like to compare the two lifts $\hat {\mb F}$ using a rank 1 construction and $\hat {\mb F}^\sim$ using a rank 2 construction where the thickening datum is pulled back from the choices made for $\hat {\mb F}$. 

Recall the notion of free quotient of Kuranishi spaces (Definition \ref{defn_Kuranishi_free_quotient}).

\begin{defn}
Let $\hat{\mb F}^\sim$ be a flow category enriched in $\uds{\bf Kur}$ with morphism space 
\beqn
K_{pq}^\sim = (G_{pq}^\sim, V_{pq}^\sim, E_{pq}^\sim, S_{pq}^\sim).
\eeqn
A {\bf free quotient} of $\hat{\mb F}^\sim$ consists of for each pair of objects $p, q \in {\rm Ob}{\mb F}$ a decomposition 
\beqn
G_{pq}^\sim = G_{pq} \times G_{pq}'
\eeqn
of compact Lie groups satisfying 1) the conjugacy class of group embeddings $G_{pr}^\sim \times G_{rq}^\sim \to G_{pq}^\sim$ respect the decomposition and 2) $G_{pq}'$ acts freely on $V_{pq}^\sim$. Notice that one obtains a new flow category $\hat{\mb F}$ with morphism spaces $K_{pq}$ being the free quotient of $K_{pq}^\sim$ by $G_{pq}'$ (see Definition \ref{defn_Kuranishi_free_quotient}). By abuse of notation, we also call the resulting flow category consisting of these free quotients a free quotient of $\hat{\mb F}^\sim$.
\end{defn}

Notice that in the above definition, if $\hat{\mb F}^\sim$ is enriched in $\uds{\bf SKur}$, then the resulting flow category $\tilde {\mb F}^\sim$ enriched in $\uds{\bf dOrb}$ obtained by taking the quotient by the full group is isomorphic to the one obtained from the free quotient obtained by taking the residual group action.

\begin{prop}\label{prop195}
Assume we are in Situation \ref{situationf1}. Fix an outercollaring width. Let ${\mc A}_{\mb F}^{\Omega_i}: {\rm Ob}{\mb F} \to {\mb Z}$, $i = 1, 2$ be two integral actions. Let $\hat {\mb F}$ be an AMS lift of ${\mb F}$ (of rank $1$, subject to an outercollaring width and ${\mc A}_{\mb F}^{\Omega_i}$). Then there exists a rank 2 AMS lift $\hat {\mb F}^\sim$ (subject to the same outercollaring width and the rank 2 integral action) and another rank $2$ lift of the outercollaring ${\mb F}$, denoted by $\hat {\mb F}$ such that 1) $\hat {\mb F}_1$ is isomorphic to a free quotient of $\hat {\mb F}$ and 2) $\hat {\mb F}^\sim$ is isomorphic to a stabilization of $\hat {\mb F}$. The same applies to an AMS lift $\hat{\mb F}_2$ of ${\mb F}$ subject to an outercollaring width and ${\mc A}_{\mb F}^{\Omega_i}$.
\end{prop}

\begin{proof}
We provide a construction that works equally well for both $\hat{\mb F}_1$ and $\hat{\mb F}_2$: their relations with our $\hat{\mb F}^\sim$ are completely symmetric. We will focus on explaining the situation for $\hat{\mb F}_1$.

Recall the roadmap of constructing $\hat {\mb F}_1$. First, consider the flow category $\outer \mb{Dom}^{\Omega_1, {\mb F}}$ enriched in $\outer \uds{\bf Curve}_{\rm rig}^{\mb C}$. By Lemma \ref{lemma_modification_comparison}, there is another flow category $\mb{Dom}^{\Omega_1 \oplus \Omega_2, {\mb F}}$ whose corresponding moduli spaces are principal bundles over those in $\outer \mb{Dom}^{\Omega_1, {\mb F}}$. Then the flow category $\hat{\mb F}$ is naturally induced from it in the following way. The differences between $\hat{\mb F}$ and $\hat{\mb F}_1$ lie in: 1) the compact Lie group associated to objects $p < q$ is $G^1_{pq} \times G^2_{pq}$; 2) the space $B^1_{pq}$ is replaced by a principal $G^2_{pq}$-bundle over $B^1_{pq}$, which is a subspace of $B_{pq}$ discussed in Lemma \ref{lemma193} where we replace each fiber of the principal ${\mc G}^2_{pq}$-bundle $B_{pq} \to B^1_{pq}$ by the maximal compact subgroup, i.e., looking at the unitary framings on the second factor; 3) the obstruction bundle and Kuranishi section are pulled back to the principal $G^2_{pq}$-bundle constructed in 2). By construction, we see that $\hat{\mb F}_1$ is isomorphic to the free quotient of $\hat{\mb F}$.

On the other hand, we define $\hat {\mb F}^\sim$ taking full advantage of the double framing. This time, along each map $u: C_\phi \to X$ there are two sets of frames, one for the line bundle induced from $\Omega_1$ and the other for the line bundle induced from $\Omega_2$. The group reduction chosen for $\hat {\mb F}_1$ can be extended to this ``double-framed'' version. This time, we consider the total space $B_{pq}$ for each $p < q$, acted on by $G_{pq} = G^1_{pq} \times G^2_{pq}$. As for the obstruction bundle, we first pull back the thickening data ${\mb W}^{\mb F}_1$ and ${\mb W}^{\mb F}_2$, which is still transverse by Lemma \ref{lemma193}, then we take the direct sum with $\uds Q_{pq}^1$ and $\uds Q_{pq}^2$ so that the obstruction bundle is 
\beqn
\uds W_{pq}^{{\mb F}_1} \oplus \uds W_{pq}^{{\mb F}_2} \oplus \uds Q_{pq}^1 \oplus \uds Q_{pq}^2
\eeqn
where the $\uds Q_{pq}^2$-component of the Kuranishi section is induced from the group reduction for the $\Omega_2$-frames. Therefore, in the category $\outer \uds{\bf S^{\rm rel} Kur}_{\rm rig}$, the statement of Proposition \ref{prop195} holds, where we construct the rigidification by pullback.

To complete the construction of $\hat {\mb F}^\sim$, one needs to extend the statement to the category $\outer \uds{\bf SKur}_{\rm rig}^{\rm NC}$. Notice that the pullback is induced from maps between domain moduli. Hence the construction of the relative NC structure for $\hat {\mb F}_1$ can be canonically lifted to $\hat {\mb F}^\sim$. To extend the smoothing, notice that for each $p, q$, the ambient space associated with $p < q$ in $\hat{\mb F}$ is a topological $G_{pq}^2$-principal bundle over the ambient space associated with $p < q$ in $\hat{\mb F}_1$. One can always use the group actions to construct smooth bundle structures inductively, hence $\hat {\mb F}$ is smoothed. Then using the stabilization structure between $\hat {\mb F}$ and $\hat {\mb F}^\sim$, one can induce a canonical smooth structure on $\hat {\mb F}^\sim$. Hence the statement is also true within the category $\outer \uds{\bf SKur}_{\rm rig}$. Lastly, the NC structure construction is quite elementary and can also be respected. 
\end{proof}

\begin{proof}[Proof of Proposition \ref{thm_integral_action_comparison}]
    Using the notations introduced in Proposition \ref{prop195}, recall that the flow categories $\tilde{\mb F}_1, \tilde{\mb F}_2$ enriched in $\outer \uds{\bf dOrb}_{\rm rig}^{\rm NC}$ are obtained as quotients of the flow categories $\hat{\mb F}_1, \hat{\mb F}_2$ enriched in $\outer \uds{\bf SKur}_{\rm rig}^{\rm NC}$. We can use the intermediate flow category $\hat{\mb F}^{\sim}$ to provide the desired construction. Indeed, recall that the Floer flow category requires passing to the quotient under the ${\mb R}$-translation action if the energy is strictly greater than $0$. To obtain the lift of the diagonal bimodule $\mb{\Delta}^{\mb{FF}}$, we take the normally complex Kuranishi spaces before taking the quotient. By passing to the orbifold quotient, we finish the proof.
\end{proof}

In fact we can also prove that given a pair of integral actions ${\mc A}_{\mb F}^{\Omega_i}$, $i = 1, 2$, the double-framed construction together with FOP perturbation produces a chain complex $C({\mb F})$ which is well-defined up to homotopy given by a continuation map also defined via a double-framed construction. The proof of Proposition \ref{thm_integral_action_comparison} can be viewed as a special case of this fact.

\subsection{Multimodules}\label{subsection_comparison_lateral_line}
Compared with flow categories, for multimodules and homotopies, the lateral lines specified in Definition \ref{defn:lateral-line} may not be homotopically unique. Note that under the double-frame construction, we need to use the lateral lines from two different versions simultaneously to define concatenation maps of the framings. They exactly facilitate the construction of the doubly-thickened moduli spaces.

\begin{proof}[Proof of Proposition \ref{prop:module-double}]
We only present the proof in the first situation; the second situation (the Poincar\'e duality multimodule) only differs in notations. Fix an outercollaring width. For any involved flow category ${\mb F}$, the double-framed construction produces AMS lifts $\hat {\mb F}^{\rm double}$ (to $\outer \uds{\bf SKur}^{\rm NC}_{\rm rig}$), which depends on choices made in the construction. Two particular choices as described above are the ``pullbacks'' from each single-framed versions, denoted by 
\begin{align*}
&\ \hat {\mb F}_{(0)}^{{\rm double}},\ &\ \hat {\mb F}_{(1)}^{{\rm double}}.
\end{align*}

Now turn to the multimodule ${\mb X}$. Starting from $\outer \mb{Dom}^{\Omega_1 \oplus \Omega_2, {\mb X}}$. For any collection of double-framed AMS lifts
\beqn
\hat {\mb F}_1^{\rm double},\ \ldots, \hat {\mb F}_m^{\rm double},\ \hat {\mb F}^{'{\rm double}}
\eeqn
then one can construct a double-framed AMS lift 
\beqn
\hat {\mb X}^{\rm double}
\eeqn
as a multimodule over $(\hat {\mb F}_1^{\rm double},\ldots, \hat{\mb F}_m^{\rm double}; \hat {\mb F}^{'{\rm double}})$. Moreover, the construction is compatible with continuation maps connecting different double-framed AMS lifts of each individual flow category. 

Notice that one has two particular situations of $\hat {\mb X}^{\rm double}$ coming from pulling back the single-framed construction $\hat {\mb X}^{\rm double}_{(s)}$ for $s = 0, 1$, which are, respectively, multimodules over $(\hat {\mb F}_{1, (s)}^{{\rm double}},\ldots, \hat {\mb F}_{m, (s)}^{\rm double}; \hat {\mb F}^{'{\rm double}}_{(s)})$. Now let $\hat{\mb\Delta}^{{\mb F}_j\mb{F}_j}_{(01)}$ and $\hat{\mb\Delta}^{{\mb F}'{\mb F}'}_{(10)}$ be the double-framed AMS lift of the diagonal bimodules. Then the double-framed concatenation/gluing produces a homotopy from
\beqn
\left( \hat {\mb\Delta}^{{\mb F}_1{\mb F}_1}_{(01)} \times \cdots \times \hat {\mb \Delta}^{{\mb F}_m {\mb F}_m}_{(01)} \right) \circ \hat {\mb X}_{(1)}^{\rm double} \circ \hat {\mb \Delta}^{{\mb F}' {\mb F}'}_{(10)}
\eeqn
to $\hat {\mb X}_{(1)}^{\rm double}$. This requires extending the proof of Theorem \ref{thm_AMS_concatenation} provided in the last section to the double-framed setting, which is nothing more complicated. Passing to the quotient via the functor $\outer \uds{\bf SKur}^{\rm NC}_{\rm rig}$ to $\outer \uds{\bf dOrb}^{\rm NC}_{\rm rig}$, one obtains a homotopy exactly expected in Proposition \ref{prop:module-double}.
\end{proof}

\section{AMS Construction for Other Flow categories, Multimodules, and Homotopies}\label{section_PSS_AMS}

In this section we discuss the AMS construction for flow categories etc. which are not from smooth Floer domains.

\subsection{Pearly flow category}\label{AMS_pearly}

Recall that the Pearly flow category $\mb{P}$ is associated to an almost K\"ahler manifold $(X, \omega, J)$ and a Morse--Smale pair $(f, g)$ on $X$. The object set is 
\begin{align*}
&\ {\rm Ob} \mb{P} = {\rm crit} f \times \Pi\ &\  {\rm where}\ \Pi = \frac{\pi_2(X)}{{\rm ker} c_1 \cap {\rm ker} \omega}.
\end{align*}
An object is written as $x = (\uds x, a_x)$ where $\uds x \in {\rm crit} f$ and $a_x \in \Pi$. We will construct AMS lifts of outercollarings of $\mb{P}$ as declared in Theorem \ref{thma_pearly_chart}. The main difference from the case of the Floer flow category is that the boundary strata of moduli spaces are given by fiber products instead of direct products.

\subsubsection{Fibered topological spaces, Kuranishi spaces, and derived orbifolds}

\begin{defn}\label{defn191}
Let $Y$ be a smooth manifold (without boundary and corners).
\begin{enumerate}

\item Let $\uds{\bf Top}(Y)$ be the regular stratification category whose objects are pairs $(N, f_-, f_+)$ where $N\in {\rm Ob} \uds{\bf Top}$ and $f_\pm: N \to Y$ is a pair of continuous maps. A morphism in $\uds{\bf Top}(Y)$ from $(N_1, f_{1, \pm})$ to $(N_2, f_{2, \pm})$ is a morphism of stratified spaces $\rho: N_1 \to N_2$ such that the following diagram commutes.
\beqn
\xymatrix{ N_1 \ar[r]^{f_{1, \pm}} \ar[d]_\rho & Y \ar[d]^{{\rm Id}_Y}\\
           N_2 \ar[r]_{f_{2, \pm}} & Y }
\eeqn
The  (non-symmetric) product of $(N_1, f_{1, \pm})$ and $(N_2, f_{2, \pm})$ is the fiber product
\beqn
\left( N_1 \underset{Y}{\times} N_2,\ \pi_1^* f_{1, -}, \ \pi_2^* f_{2, +} \right):= \left( N_1\ {}_{f_{1, +}}\underset{Y}{\times}{}_{f_{2, -}} N_2,\ \pi_1^* f_{1, -},\ \pi_2^* f_{2, +} \right).
\eeqn
The monoidal unit is $(Y, {\rm Id}_Y, {\rm Id}_Y)$. 

\item Let $\uds{\bf Kur}(Y)$ be the regular stratification category whose objects are $(K, f_-, f_+)$ where $K = (G, V, E, S)$ is a Kuranishi space and $f_\pm: V \to Y$ is a pair of $G$-invariant topological submersions. Let $\uds{\bf SKur}(Y)$ be the corresponding category of smooth Kuranishi spaces together with smooth submersions to $Y$. Let $\uds{\bf SKur}^{\rm NC}(Y)$ be the subcategory whose objects are smooth normally complex Kuranishi spaces with $G$-invariant submersions. Let $\uds{\bf S^{\rm rel} Kur}(Y)$ be the category of relatively smooth Kuranishi spaces $K/B = (G, V/B, E, S)$ together with $G$-invariant maps $f_\pm: V \to Y$ such that they are fiberwise submersions. Let $\uds{\bf S^{\rm rel}Kur}^{\rm NC}(Y)$ be the subcategory where the relative Kuranishi spaces are equipped with relative NC structures. 

\item Let $\uds{\bf dOrb}(Y)$ be the category whose objects are pairs $({\mc D}, f_-, f_+)$ where ${\mc D} = ({\mc U}, {\mc E}, {\mc S})$ is a derived orbifold and $f_\pm: {\mc U} \to Y$ are smooth submersions. Let $\uds{\bf dOrb}^{\rm NC}(Y)$ be the category whose objects are $({\mc D}, f_-, f_+)$ where ${\mc D}$ is an NC derived orbifold and $f_\pm: {\mc U} \to Y$ are smooth submersions.

\end{enumerate}
\end{defn}

\begin{rem}
Unlike the presumed stable complex versions, the monoidal structure given by fiber products in the NC versions ($\uds{\bf SKur}^{\rm NC}(Y)$, $\uds{\bf S^{\rm rel}Kur}^{\rm NC}(Y)$, $\uds{\bf dOrb}^{\rm NC}(Y)$) does not require the almost complex structure on $Y$ because $Y$ is assumed to be a manifold.
\end{rem}

Notice that one can also talk about collared and rigidified versions of $\uds{\bf Kur}(Y), \uds{\bf SKur}(Y)$, and $\uds{\bf dOrb}(Y)$.

One needs to generalize the discussion of stable smoothings of topological Kuranishi spaces to the category $\uds{\bf Kur}(Y)$.

\begin{defn}
A {\bf stable smoothing} of an object $(K, f_\pm) \in \uds{\bf Kur}(Y)$ where $K = (G, V, E, S)$ consists of a unitary representation $W$ of $G$, a $G$-smoothing on $V \times W$, a $G$-invariant extension $\tilde f_\pm: V \times W \to Y$ of $f_\pm: V \cong V \times \{0\} \to Y$ such that $\tilde f_\pm$ are smooth submersions.
\end{defn}

\begin{lemma}[Stable smoothing with maps] (cf. \cite[Lemma 4.5]{AMS}) \label{lemma_smoothing_map}
Let $V$ be a smooth $G$-manifold and $f_0: V \to Y$ be a continuous $G$-invariant map. Then there exists a unitary representation $W$ of $G$, a $G$-smoothing on $V \times W$ which is in the same stable isotopy class of the given smooth structure on $V$, and a smooth submersion $\tilde f: V \times W \to Y$ whose restriction to $V \times \{0\}$ coincides with $f_0$.
\end{lemma}

\begin{proof}
The proof is essentially that of \cite[Lemma 4.5]{AMS}, which we reproduce using our language. Choose a Riemannian metric on $Y$. %
Then one can choose a smooth $G$-invariant approximation $f: V \to Y$ and a smooth $G$-invariant function $\epsilon: V \to {\mb R}_+$ such that for all $v \in V$, there is a unique shortest geodesic of length shorter than $\epsilon(v)$ connecting $f_0 (v)$ and $f (v)$, and such that the exponential map $\exp_{f(v)} \xi$ is well-defined for $\xi$ in the $\epsilon(v)$-ball of $T_{f(v)} Y$. Then there is a smooth equivariant bundle map $\rho: f^* TY \to f^* TY$ compressing $f^* TY$ diffeomorphically into the $\epsilon$-disk bundle. Therefore the map 
\beqn
\tilde f: f^* TY \to Y,\ \tilde f(v, \xi) = \exp_{f(v)} ( \rho(\xi)).
\eeqn
is a smooth submersion. Choose another smooth (real) $G$-vector bundle $E \to V$ such that 
\beqn
E \oplus f^* TY \cong V \times W
\eeqn
for some unitary representation $W$ of $G$. $\tilde f$ clearly extends to a smooth $G$-equivariant map from $V\times W$ to $Y$, which is still a smooth submersion.

We would like to define a self-homeomorphism of $V \times W$ to pullback the smooth structure on $V \times W$. Indeed, we can write\beqn
f_0 (v) = \exp_{f (v)} \xi_0 (v),\ \xi_0 (v) \in T_{f(v)} Y.
\eeqn
Define a nonlinear bundle map
\beqn
F: f^* TY \to f^* TY,\ (v, \xi) \mapsto (v, \xi + \xi_0(v)).
\eeqn
Extend it to a smooth map $V \times W \to V \times W$ which is the identity on the summand $E$. Then we can see the following diagram commutes.
\beqn
\xymatrix{    V \ar[r]^{f_0} \ar[d]_{v \mapsto (v, 0)} & Y \\
            V\times W \ar[r]_F  & V \times W  \ar[u]_{\tilde f} }
\eeqn
Therefore, $f_0$ is extended to a smooth submersion from $V \times W$ to $Y$. Lastly, the smoothing on $V \times W$ is in the same stable isotopy class because the choices made in the proof, especially the smooth approximation $f$ of $f_0$, can all be made isotopic.
\end{proof}

\subsubsection{The monotone flow category of holomorphic cylinders}

Let $(X, \omega, J)$ be an almost K\"ahler manifold. Denote
\beqn
\Pi^{\rm eff}:= \Big\{ a \in \Pi\ |\ \omega(a) > 0 \Big\} \cup \{0\}
\eeqn
which is an abelian monoid whose Grothendieck construction is $\Pi$. 

We would like to generalize the notion of monotone flow categories. Recall that a monotone flow category enriched in $\uds {\bf C}$ contains objects $M_d^{\mb f}$ indexed by $d \in {\mb Z}_{\geq 0}$. In the current case, one can define a generalization consisting of objects $M_a^{\mb f}$ indexed by $a\in \Pi^{\rm eff}$, together with structural maps from the fiber product of $M_a^{\mb f}$ and  $M_b^{\mb f}$ over a space to $M_{a+b}^{\mb f}$ satisfying a version of associativity condition. We describe such a monotone flow category $\mb{cyl}(X)$ enriched in $\uds{\bf C} = \uds{\bf Top}(X)$ by defining
\beqn
M_a^{\mb{cyl}}(X):= \left( \bigsqcup_{\omega(A) = a} \ov{\mc M}{}_{0,2}^{\mb R}(X, J, A),\ \ev_-, \ \ev_+ \right) \in {\rm Ob} \uds{\bf Top}(X),
\eeqn
where $\ev_{\pm}$ are the evaluation maps at the marked points. Then for $a, b \in \Pi^{\rm eff}$, there is a natural map
\beqn
M_a^{\mb{cyl} }(X) \underset{X}{\times} M_b^{\mb{cyl} }(X) \to M_{a + b}^{\mb{cyl} }(X).
\eeqn
These structures define the monotone flow category $\mb{cyl}(X)$. By taking the Grothendieck construction, it induces an ordinary flow category
\beqn
\mb{Cyl}(X)
\eeqn
called the {\bf cylinder flow category} enriched in $\uds{\bf Top}(X)$ whose object set is $\Pi$.

To carry out the AMS construction, we first choose an integral symplectic form $\Omega \in \Omega^2(X)$ satisfying 1) $J$ is tamed by $\Omega$ and 2) ${\rm Ker} \Omega = {\rm ker} \omega\subset \pi_2(X)$. Then for each $a \in \Pi^{\rm eff}$, there is an associated integer $d = d_a$. As $J$ is tamed by $\Omega$, $d_a \geq 0$.

\begin{prop}\label{prop215}
There exists a lift $\wh{\mb{cyl}}  (X)$ of the outercollaring $\outer \mb{cyl}  (X)$ to $\outer \uds{\bf SKur}^{\rm NC}_{\rm rig}(X)$ as a monotone flow category supported over $\Pi^{\rm eff}$.
\end{prop}

\begin{proof}
We first describe a monotone flow category of curves $\outer \mb{Dom}^\Omega$ enriched in $\outer \uds{\bf Curve}_{\rm rig}^{\mb C}$, which is the restriction of $\outer \mb{Dom}$ (specified in Section \ref{section_AMS_domains}) to open sets induced from $\Omega$. For each $d$, let
\beqn
A_d^\Omega \subset A_d
\eeqn
to be the Alexandrov open subset of partitions $d = d_0 + \cdots + d_l$ such that for each $d_l$, there exists $a_l \in \Pi^{\rm eff}$ with $d_l = \Omega(a_l)$. Recall that the monotone flow category $\outer \mb{Dom}$ gives objects
\beqn
{\mc C}_d = ({\mc G}_d, B_d, C_d) \in {\rm Ob} \uds{\bf Curve}.
\eeqn
For $a \in \Pi^{\rm eff}$, let $B_a^\Omega \subset B_{\Omega(a)}$ be the union of strata prescribed by $A_d^\Omega$, which is open and ${\mc G}_d$-invariant. Define
\beqn
{\mc C}_a^\Omega \in {\rm Ob} \outer \uds{\bf Curve}_{\rm rig}^{\mb C}
\eeqn
to be the restriction of ${\mc C}_{\Omega(a)}$ to $B_a^\Omega$. Then one obtains a monotone flow category $\outer \mb{Dom}^{\Omega, \mb{cyl}}$. 

Then one has a corresponding monotone flow category of framed maps consisting of objects
\beqn
{\rm Map}_a^{\rm fr} (X) \to \outer B_a^{\Omega}
\eeqn
whose fiber of $\phi \in \outer B_a^\Omega$ is the space of pairs $(u, F)$ where $u: \mathring C_\phi \to X$ is a smooth map on each irreducible component which extends continuously to $C_\phi$ and which converges exponentially at cylindrical nodes and markings and where $F$ is a frame of the holomorphic line bundle $L_{u, \Omega} \to C_\phi$.

One can also choose a collection of transverse thickening data
\beqn
\nu_a: W_a \to \Gamma_c( \outer \mathring C_a^\Omega \times X, \Lambda_{\outer \mathring C_a^\Omega/ \outer B_a^\Omega}^{0,1}\otimes TX)
\eeqn
which are collared and respect the structural maps. Here $W_a$ is a collection of unitary representations of $G_{\Omega(a)}$. On the other hand, one can choose a group reduction scheme, given by a collection of continuous maps
\beqn
\lambda_a: {\rm Map}_a^{\rm fr}(X) \to Q_a\cong {\mc G}_{\Omega(a)}/ G_{\Omega(a)}.
\eeqn
Then one obtains a corresponding thickening
\beqn
K_a = (G_a, V_a, E_a, S_a) \in {\rm Ob} \outer \uds{\bf Kur}
\eeqn
where 
\beqn
V_a = \Big\{ (\phi, u, F, e)\ |\ (\phi, u, F) \in {\rm Map}_a^{\rm fr}(X),\ e \in W_a,\ \phi = \phi_F \in \outer B_a^\Omega,\ \ov\partial_J u + \nu_a (e) = 0    \Big\},
\eeqn
with the natural $G_a$-action, $E_a$ is the trivial bundle with fiber $W_a \oplus Q_a$, and $S_a$ is the section given by 
\beqn
S_a (\phi, u, F, e) = (e, \lambda_a (\phi, u, F)).
\eeqn
Similar to Theorem \ref{thm_footprint}, the natural map
\beqn
S_a^{-1}(0)/G_a \to \outer M_a^{\mb{cyl} } (X)
\eeqn
is an isomorphism of orbispaces. By an abuse of notation, let $V_a$ denote the regular locus for the operator $\ov\partial_J u + \nu_a (e)= 0$, i.e., the subspace such that the linearization of this equation is surjective. Then, it contains the zero locus by the transversality assumption. Then one obtains a lift of the monotone flow category $\outer \mb{cyl} (X)$ to the category $\outer \uds{\bf Kur}_{\rm rig}(X)$. In the same way as argued in Subsection \ref{subsection_relative_smooth}, the lift is in fact enriched in $\outer \uds{\bf S^{\rm rel}Kur}_{\rm rig}(X)$. We denote this lift by 
\beqn
\widehat{\mb{cyl}} (X).
\eeqn

The construct of normal complex structure is slightly easier than the case of Floer flow categories. Indeed, because we do not have Hamiltonian perturbations, the (vertical) linearized operators are the Cauchy--Riemann operators
\beqn
D_{u, e}= \left[ D_u\ov\partial_J + D_u \nu_a(u, e) + \nu_a(u) \right]: \Omega^0 (C_\phi, u^* TX) \oplus W_a \to \Omega^{0,1}(C_\phi, u^* TX)
\eeqn
where $W_a$ is a complex representation. Recall that the canonical linearization is calculated by the Levi--Civita connection on $X$ induced from the symplectic form and the almost complex structure, which is not necessarily complex-linear. However, one can choose a 1-parameter family of linear operators connecting $D_{u, e}$ to a complex linear one, without using the complicated interpolation argument as in the Floer case. This means for each individual moduli space, one has a canonical homotopy class of isomorphisms between the vertical tangent bundle and a complex vector bundle. One can choose concrete isomorphisms inductively. This provides a lift to $\outer \uds{\bf S^{\rm rel}Kur}_{\rm rig}^{\rm NC}(X)$.

The smoothing procedure is slightly different from the previous case as we also need to make the evaluation maps smooth. However, by Lemma \ref{lemma_smoothing_map}, one can perform stable smoothing on $(K_a, \ev_\pm)$ so that not only the Kuranishi spaces become smooth, but also the evaluation maps become smooth submersions. Then following the same inductive procedure, after stabilization, one obtains a lift to the category $\outer\uds{\bf SKur}_{\rm rig} (X)$. Together with the stable complex structures on the domain flow categories, this is a lift to $\outer \uds{\bf SKur}_{\rm rig}^{\rm NC}(X)$.
\end{proof}

\begin{proof}[Proof of Theorem \ref{thma_pearly_chart}]
We first describe the lift of the pearly flow category ${\mb P}$. Choose a lift $\widehat{\mb{cyl}}  (X)$ provided by Proposition \ref{prop215}. By taking quotient, one obtains a monotone flow category enriched in $\outer\uds{\bf dOrb}^{\rm NC}_{\rm rig}(X)$, which induces an ordinary flow category, denoted by $\widetilde{\mb{Cyl}} (X)$, whose objects are elements of $\Pi$. 

Then one can obtain a lift of the corresponding outercollaring of ${\mb P}$ by ``inserting flow lines.'' We do this inductively. In the lowest energy case, one only needs to add flow rays at the negative and positive ends to obtain Kuranishi charts of moduli spaces labelled by a pair of critical points. In higher energy case, there are configurations where in the interior there are gradient segments of positive length. An inductive scheme can be designed in a rather elementary way to produce inductively all charts for involved moduli spaces. Hence one obtains a lift of the outercollaring of ${\mb P}$ to $\outer \uds{\bf dOrb}_{\rm rig}^{\rm NC}$ (inserting flow lines does not affect the NC structure).

Next, we discuss the regularization of the bimodule $\mb{\Delta}^{\mb{PP}}$. In fact, one can take different Morse--Smale pairs and different almost complex structures to define $\mb{P}$ and consider the bimodule as an interpolation between them. However, we choose the same one to simplify the notations.

The construction is very similar to the regularization of the Pearly flow category. First, there is a monotone bimodule $\mb{\Delta}^{\mb{cyl}  (X)}$ given by holomorphic cylinders which have a parametrized principal component. This is a monotone bimodule over $(\mb{cyl} (X); \mb{cyl}  (X))$. 

Then after outercollaring, following the same procedure as the proof of Proposition \ref{prop215}, one obtains an AMS lift  to $\outer \uds{\bf SKur}_{\rm rig}^{\rm NC}(X)$, denoted by $\wh{\mb\Delta}^{\mb{cyl} (X)}$. After taking quotient, one obtains a bimodule $\widetilde{\mb\Delta}^{\mb{cyl}^{\mb P} (X)}$ over $(\widetilde{\mb{cyl}}  (X); \widetilde{\mb{cyl}}  (X))$ enriched in $\outer \uds{\bf dOrb}^{\rm NC}_{\rm rig}(X)$. Then using the fact that evaluation maps are smooth submersions, by inserting Morse flow lines, one obtains a lift of $\outer \mb{\Delta}^{\mb{PP}}$ as a bimodule over $(\tilde{\mb P}; \tilde{\mb P})$. This proves the second item of Theorem \ref{thma_pearly_chart}.
\end{proof}

\subsection{Deformed Morse-to-Morse bimodule}

The regularization of the Morse-to-Morse bimodule was discussed in \cite{Bai_Xu_Arnold} using the original perturbation from \cite{AMS}. 

\begin{prop}
Upon choosing an integral symplectic form $\Omega$, there exists a lift of $ \Delta^{\mb{MM}}$ of the deformed Morse-to-Morse bimodule to $\uds{\bf dOrb}^{\rm NC}_{\rm rig}$ as a bimodule over $(\mb{M}; \mb{M})$ (without outercollaring).
\end{prop}

\begin{proof}
For each $d$, consider the moduli space $\ov{\mc M}_{0,2}(\mb{CP}^d, d)$ of 2-pointed stable holomorphic spheres in $\mb{CP}^d$. Choose the negative convention; the other choice works as well. Let $B_d = B_d^-$ be the subset of curves whose images are not contained in any hyperplane and whose value at the negative end is $[1, 0, \ldots, 0]$. Then consider the curve ${\mc C}_d = ({\mc G}_d, B_d, C_d)$. Notice that this is a family of curves without boundary or corners. 

Now choose an integral symplectic form $\Omega$. For each $a \in \Pi^{\rm eff}$, denote $d = d_a$. Consider the space of framed maps
\beqn
{\rm Map}_d^{\rm fr} = \Big\{ (\phi, u, F)\ |\ \phi \in B_d, u: C_\phi \to X,\ F = (f_0, \ldots, f_d)\ {\rm is\ a\ basis\ of\ } H^0(L_{u, \Omega}) \Big\}.
\eeqn
Choose a group reduction scheme $\lambda_d: {\rm Map}_d^{\rm fr} \to Q_d$ and a transverse perturbation 
\beqn
\iota: W \to \Gamma_c( \mathring C_d \times X, \Lambda^{0,1}_{\mathring C_d / B_d} \otimes TX).
\eeqn
This induces a topological Kuranishi space 
\beqn
K_d = (G_d, V_d, E_d, S_d)
\eeqn
with ambient space
\beqn
V_d = \Big\{ (\phi, u, F, e)\ |\ (\phi, u, F) \in {\rm Map}_d,\ e \in W_a,\ \phi = \phi_F \in B_d,\ \ov\partial_J u + \iota (e) = 0    \Big\},
\eeqn
obstruction bundle $E_d = W \oplus Q_d$, and Kuranishi section 
\beqn
S_d (\phi, u, F, e) = (e, \lambda_d (\phi, u, F)).
\eeqn
It also comes with a {\rm ref}-$C^{\infty, 1}$-structure on $V_d/B_d$ as well as continuous evaluation maps $\ev_\pm: V_d \to X$. Then using Lemma \ref{lemma_smoothing_map}, one can find a stable smoothing of $V_d$ making the evaluation maps submersive.

Last, using the submersiveness of the evaluation maps, one obtains a lift of the bimodule $\mb{\Delta}^{\mb{MM}}$ by taking fiber products with stable/unstable manifolds of the Morse flow.
\end{proof}

\begin{rem}
In \cite{Bai_Xu_Arnold} one considered the Morse-to-Floer and Floer-to-Morse bimodules and their concatenation along the Floer flow category. However, their concatenation along the Morse flow category is difficult to regularize (the codimension two gluing is hard to fit into the framework of flow bimodules). In this paper, we do not consider these two bimodules and their concatenations. The Morse-to-Morse bimodule will arise from the concatenations of the Morse-to-Pearly and the Pearly-to-Morse bimodules, whose regularization will be addressed later.
\end{rem}

\subsection{Pearly-to-Morse and Morse-to-Pearly bimodules}

In this subsection we regularize the two bimodules shown in the title. To simplify the notations, we will only use a single integral symplectic form $\Omega \in \Omega^2(X)$ which is a multiple of a rational approximation of $\omega$ which also tames $J$. The case with different choices of $\Omega$ can be discussed in a similar way as the case of pair-of-pants multimodules.

\subsubsection{The domains}

The model for the Morse flow category can be viewed as a monotone flow category over the trivial monoid $\{0\}$. Let such a trivial monotone flow category be $\mb{o}$, which has only one object $M_0^{\mb{o}}$ chosen to be a singleton in the category $\uds{\bf Curve}$. We would like to define a monotone flow bimodule over ${\mb Z}_{\geq 0}$. For each $d$, we consider the family of curves
\beqn
{\mc C}_d^{\mb{PM}} = ({\mc G}_d, B_d^{\mb{PM}}, C_d^{\mb{PM}})
\eeqn
where the space of curves $B_d^{\mb{PM}}$ in $\mb{CP}^d$ is different from previous cases. Indeed, consider the moduli space of parametrized holomorphic cylinders
\beqn
u: {\mb R} \times S^1 \to \mb{CP}^d
\eeqn
of degree $d$. When compactifying this space, the bubbling near the negative and positive ends is treated differently. Near the negative ends, the broken pieces are considered modulo only the real translation; near the positive end, the broken pieces are considered modulo ${\mb C}^*$, i.e., translations and rotations. The domain of a point in the compactification still has two distinguished points: $z_-$ and $z_+$. Then let $B_d^{\mb{PM}}$ be the subspace in the compactification satisfying 1) the image is not contained in any hyperplane and 2) the map satisfies the positive:
\beqn
u(z_+) = [0, 0, \ldots, 0, 1] \in {\mb P}^d.
\eeqn

The structural map of this bimodule only contains the negative part as we view it as a bimodule with one side the trivial flow category: namely, given $d, d'$, there is a morphism
\beqn
{\mc C}_d \times {\mc C}_{d'}^{\mb{PM}} \to {\mc C}_{d+d'}^{\mb{PM}}.
\eeqn
The construction is similar to before and hence skipped. One can check that the collection of equivariant families of curves $\mc{C}_d^{\mb{PM}}$ defines a monotone bimodule.

\subsubsection{Cylinders in $X$}

We would like to construct monotone bimodules enriched in  $\uds{\bf dOrb}_{\rm rig}^{\rm NC}(X)$. For the Morse side, there is still the trivial monotone flow category with a single object whose only morphism space is the singleton, i.e., the triple
\beqn
(X, {\rm Id}, {\rm Id})
\eeqn
where $X$ is regarded as a normally complex derived orbifold. On the other hand, one has constructed the lift $\wt{\mb{cyl}}  (X)$ of $\mb{cyl} (X)$ to $\outer \uds{\bf dOrb}^{\rm NC}_{\rm rig}(X)$. By considering parametrized cylinders in $X$, one can easily formulate a monotone bimodule
\beqn
\mb{cyl}^{\mb{PM}}(X)
\eeqn
over $(\mb{cyl}(X); \mb{o})$ consisting of moduli spaces of parametrized cylinders in $X$ with two marked points where the compactification at negative and positive marked points is taken differently.

\begin{prop}\label{prop218}
Given an AMS lift $\widetilde{\mb{cyl}} (X)$ of the outercollaring $\outer \mb{cyl} (X)$, there exists a lift of the outercollaring $\outer \mb{cyl}{}^{\mb{PM}}(X)$ to $\outer \uds{\bf dOrb}^{\rm NC}_{\rm rig}(X)$, denoted by $\widetilde{\mb{cyl}}{}^{\mb{PM}}(X)$, as a monotone bimodule over $(\widetilde{\mb{cyl}} (X); \mb{o})$.
\end{prop}

\begin{proof}[Sketch of Proof] It is completely parallel to the proof of Proposition \ref{prop215}.
\end{proof}

\begin{cor} Given an AMS lift $\tilde {\mb P}$ of the pearly flow category ${\mb P}$ subject to a fixed outercollaring width and an integral symplectic form, there exists a lift $\tilde{\mb B}^{\mb{PM}}$ resp. $\tilde{\mb B}^{\mb{MP}}$ of the corresponding collaring of $\mb{B}^{\mb{PM}}$ resp. $\mb{B}^{\mb{MP}}$ to $\outer \uds{\bf dOrb}_{\rm rig}^{\rm NC}$ as a bimodule over $(\tilde{\mb P}; \outer \mb{M})$ resp. over $(\outer {\mb M}; \tilde {\mb P})$. Moreover, any two such lifts are homotopic. 
\end{cor}

\begin{proof}
We work over the domain bimodule $\wt{\mb{cyl}}{}^{\mb{PM}}(X)$ (and in parallel $\wt{\mb{cyl}}^{}{\mb{MP}}(X)$). Extend thickening data used in the construction of $\tilde{\mb P}$ to thickening data over $\wt{\mb{cyl}}^{}{\mb{PM}}(X)$. The construction of NC structures and stable smoothing can all be reproduced in the same way as before. This provides the AMS lifts of outercollarings of $\mb{B}^{\mb{PM}}$ and $\mb{B}^{\mb{MP}}$. 
\end{proof}

\subsection{Pearly-to-Floer and Floer-to-Pearly bimodules}

The AMS constructions for the bimodules $\mb{B}^{\mb{FP}}$ and $\mb{B}^{\mb{PF}}$ are slight variants of that for $\mb{\Delta}^{\mb{PP}}$. We describe in detail the construction for $\mb{B}^{\mb{PF}}$. 

\subsubsection{A new category of fibered spaces}

In order to deal with the mixture of fiber products (in the Pearly part) and direct products (in the Floer part), we introduce the space $X^+ = X \sqcup \{{\rm pt}\}$ which is a smooth manifold of mixed dimensions. Then similar to Definition \ref{defn191}, one can consider the category $\uds{\bf Top}(X^+)$ whose objects are spaces with continuous maps to $X^+$, as well as the categories $\uds{\bf Kur}(X^+)$ and $\uds{\bf dOrb}(X^+)$. Then there are natural functors
\beqn
\vcenter{ \xymatrix{     \uds{\bf Top}(X) \ar[rd] & & \uds{\bf Top} \ar[ld]\\
 & \uds{\bf Top}( X^+ ) & } }.
\eeqn
Similar functors exist for Kuranishi spaces or derived orbifolds. 

Now consider the flow category $\mb{Cyl} (X)$ (obtained from the monotone flow category $\mb{cyl} (X)$), which can be viewed as enriched in $\uds{\bf Top}(X^+)$. Similarly, the Floer flow category ${\mb F}$ associated to a Hamiltonian $H$ and an almost complex structure $J$ is also enriched in $\uds{\bf Top}(X^+)$. Then the Floer data used for defining the bimodules $\mb{B}^{\mb{PF}}$ resp. $\mb{B}^{\mb{FP}}$ also induces a bimodule
\beqn
\mb{Cyl}^{\mb{PF}}(X)\ {\rm resp.}\ \mb{Cyl}{}^{\mb{FP}}(X)
\eeqn
over $(\mb{Cyl} (X); \mb{F})$ resp. over $(\mb{F}; \mb{Cyl} (X))$.

\begin{prop}\label{prop:21-10}
Given an AMS lift $\widetilde{\mb{Cyl}} (X)$ of $\outer \mb{Cyl} (X)$ to $\outer \uds{\bf dOrb}_{\rm rig}^{\rm NC}(X)$ and an AMS lift $\tilde{\mb F}$ of $\outer \mb{F}$ to $\outer \uds{\bf dOrb}_{\rm rig}^{\rm NC}$, there exists a lift (called an AMS lift) of the outercollaring $\outer \mb{Cyl}{}^{\mb{PF}}(X)$ resp. $\outer \mb{Cyl}{}^{\mb{FP}}(X)$ to $\outer \uds{\bf dOrb}{}_{\rm rig}^{\rm NC}(X^+)$, denoted by $\widetilde{\mb{Cyl}}{}^{\mb{PF}}(X)$ resp. $\widetilde{\mb{Cyl}}{}^{\mb{FP}}(X)$.
\end{prop}

\begin{proof}[Sketch of Proof] There is no essential difference from the previous cases. 
\end{proof}

\begin{cor}
Given an AMS lift $\tilde{\mb P}$ of the outercollaring of the pearly flow category $\mb{P}$ and an AMS lift $\tilde{\mb F}$ of the Floer flow category, there exists a lift of $\outer \mb{B}^{\mb{PF}}$ resp. $\outer \mb{B}^{\mb{FP}}$ to $\outer \uds{\bf dOrb}_{\rm rig}^{\rm NC}$. 
\end{cor}

\begin{proof}
It follows from the same arguments as we have seen many times, building on Proposition \ref{prop:21-10}. The only difference is that we should insert Morse flow lines on the pearly side to complete the construction.
\end{proof}

\subsection{Morse-to-Floer and Floer-to-Morse bimodules}\label{subsection_Morse_Floer}

Now we consider the AMS lift of the bimodules $\mb{B}^{\mb{MF}}$ and $\mb{B}^{\mb{FM}}$, which are responsible for the PSS and SSP maps claimed in Theorem \ref{thma_PSS}. The claim of such AMS lifts is stated as Theorem \ref{thma_PSS_AMS}. As the PSS map on the homology level is supposed to be canonical, we also need to compare different choices. 

\begin{proof}[Proof of Theorem \ref{thma_PSS_AMS}]We sketch the construction for $\mb{B}^{\mb{MF}}$. The other case has no difference other than a simple reversal. See \cite[Section 7]{Bai_Xu_Arnold} for a more explicit discussion in a slightly different framework.

By removing the Morse flow ray, the Morse-to-Floer bimodule ${\mb B}^{\mb{MF}}$ becomes the cigar bimodule $\mb{B}^{\rm cigar}$, which is enriched in ${\rm Top}(X^+)$. Recall that it is a bimodule over $({\mb O}; {\mb F})$ where ${\mb O}$ is the $\Pi$-covering of the trivial flow category. Then the AMS construction is similar to that for the cigar bimodule $\mb{B}^{\rm cigar}$. When construct an AMS lift $\tilde {\mb B}^{\rm cigar}$ to $\outer \uds{\bf dOrb}_{\rm rig}^{\rm NC}$ (which is covered by Theorem \ref{thma_module_lift}), especially for smoothing, one can make the evaluation maps at the incoming marked point smooth submersions onto $X$. Hence one obtains a lift, still denoted by $\tilde {\mb B}^{\rm cigar}$, enriched in $\outer \uds{\bf dOrb}_{\rm rig}^{\rm NC}(X^+)$. Then by adding gradient rays, one obtains a lift $\tilde {\mb B}^{\mb{MF}}$ to $\uds{\bf dOrb}_{\rm rig}^{\rm NC}$ as a bimodule over $(\outer \mb{M}; \tilde {\mb F})$. 

For the second item of Theorem \ref{thma_PSS_AMS}, one also needs to compare different choices of integral actions. This can be done using the double-framed construction introduced in Section \ref{section_comparison} in a similar way as proving Proposition \ref{prop:module-double}. We omit the details.
\end{proof}

\begin{proof}[Proof of Theorem \ref{thm126}]
For the concatenation $\mb{B}^{\mb{FF}'} \circ \mb{B}^{\mb{F'}\mb{M}'}$ at the Floer flow category and its homotopy to $\mb{B}^{\mb{F}\mb{M}'}$, the corresponding AMS construction is similar to what is required for establishing the naturality of the continuation map. For the concatenation $\mb{B}^{\mb{FM}} \circ \mb{B}^{\mb{MM}'}$ at the Morse flow category $\mb{M}$, the construction is even simpler. We omit the details. 
\end{proof}

\subsection{Proof of Theorem \ref{thma_PSS_SSP}}\label{subsection_PSS_triangle}

We have constructed the lifts of flow categories and bimodules listed in the statement of Theorem \ref{thma_PSS_SSP}. We choose a common outercollaring width and a compatible set of integral actions. To finish the proof of Theorem \ref{thma_PSS_SSP}, one needs to establish corresponding lifts of homotopies involving certain concatenations and gluings. Notice that we only glue at the Floer flow category or the pearly flow category, but not the Morse flow category. We need to construct such homotopies for many triangles of the diagram \eqref{triangle_lift}. However, the constructions for each triangle have no substantial differences. Hence we only give proofs for two typical cases.

\subsubsection{Concatenation at the Floer flow category}

We prove the following statement.

\begin{lemma}
Fix an outercollaring width and an integral action ${\mc A}_{\mb F}^\Omega$ which is compatible with all involved Floer data. Let $\hat {\mb F}$ be an AMS lift of ${\mb F}$ subject to the outercollaring width and ${\mc A}_{\mb F}^\Omega$. Let $\hat {\mb P}$ be an AMS lift of the pearly flow category subject to the same outercollaring width and the integral symplectic form $\Omega$. Let $\hat {\mb B}^{\mb{PF}}$ be an AMS lift of ${\mb B}^{\mb{PF}}$ and $\hat {\mb B}^{\mb{FP}}$ be an AMS lift of $\mb{B}^{\mb{FP}}$. Let $\hat{\mb{\Delta}}^{\mb{PP}}$ be an AMS lift of the diagonal bimodule. Then there exists a concatenation $\hat {\mb B}^{\mb{PF}} \circ \hat {\mb B}^{\mb{FP}}$ and a homotopy from this concatenation to $\hat{\mb\Delta}^{\mb{PP}}$.
\end{lemma}

\begin{proof}
The construction of a concatenation $\hat{\mb B}^{\mb{PF}} \circ \hat{\mb B}^{\mb{FP}}$ is no different from other situations involving concatenating at a Floer flow category. The construction of the lift of the homotopy to $\outer \uds{\bf SKur}_{\rm rig}^{\rm NC}$ also only differs in notation. 
\end{proof}

\subsubsection{Concatenation at the pearly flow category}

The concatenation at the pearly flow category is slightly more complicated. We consider the reversed concatenation ${\mb B}^{\mb{FP}} \circ \mb{B}^{\mb{PF}}$. 

Notice that the homotopy between the concatenation $\mb{B}^{\mb{FP}}\circ \mb{B}^{\mb{PF}}$ and $\mb{\Delta}^{\mb{FF}}$ is the union of two parts: the first part is when there are gradient flow segments of a positive total length (where the length is a parameter for the homotopy); the second part is when there are no gradient segments. See Figure \ref{figure_FPPF} for an illustration. The boundary between the two parts can be viewed as another concatenation
\beqn
\mb{Cyl}^{\mb{FP}}(X) \circ \mb{Cyl}^{\mb{PF}}(X)
\eeqn
at the flow category $\mb{Cyl}(X)$. Notice that one has constructed AMS lifts $\wh{\mb{Cyl}}^{\mb{FP}}(X)$ and $\wh{\mb{Cyl}}^{\mb{PF}}(X)$ to $\outer \uds{\bf SKur}_{\rm rig}^{\rm NC}$. The concatenation can be constructed in a straightforward way. Then we extend the concatenation in the two directions. One direction is the gluing, which has no essential new features for this case. The other direction is to insert gradient flow segments. This is similar to the cases of $\mb{P}$ and $\mb{\Delta}^{\mb{PP}}$. We omit the details.

\newpage

\part{EQUIVARIANT FLOER THEORY AND STEENROD OPERATIONS}

\section*{Outline of Part 4}

Note that we would like to switch to the cohomological convention. The Floer equation will remain the same. However, it will be the equation for the ascending gradient flow equation of the action functional (which differs from the one used previously by a negative sign). Hence the differential increases the action functional. The grading is also changed to the cohomological one so that the pair-of-pants product is degree-preserving. Moreover, we make use of the coherent orientations without explicitly mentioning them.

\section{Morse Theory on Classifying Spaces}

The Borel $\zp$-equivariant cohomology of a point with $\fp$-coefficient is the ring 
\beqn
H_{\zp}^{\rm Borel} ({\rm pt}; \fp) \cong H^*( B\zp; \fp) = \left\{ \begin{array}{ll} \fp [\tu, \uptheta]/  \langle \uptheta^2 = 0 \rangle,\ &\ \p > 2; \\
{\bf F}_2 [ \uptheta], &\ \p = 2.  \end{array}\right.
\eeqn
where $\tu$ has degree $2$ while $\uptheta$ has degree $1$. One usually considers the completion by allowing formal power series rather than polynomials in $\tu$ and $\uptheta$. We denote the completion by 
\beqn
\wt{\bf K}_{\p, 0}:= \left\{ \begin{array}{ll} \fp\psu\langle \uptheta \rangle,\ &\ \p > 2,\\
{\bf F}_2 \pst,\ &\ \p = 2.
\end{array}\right.
\eeqn
By inverting $\tu$ for $p \geq 3$ and inverting $\uptheta$ for $p=2$, one get the Laurent version
\beqn
\tkp:= \left\{ \begin{array}{ll} \fp\lsu \langle \uptheta \rangle,\ &\ \p>2,\\
{\bf F}_2 \lst,\ &\ \p = 2.\end{array}\right.
\eeqn
This ring is not the cohomology of any natural topological space. However, it can be realized as the cohomology of a certain flow category which will be convenient for establishing the Tate equivariant theories. We will also use the notation
\beqn
\kp := \fp\lsu
\eeqn
for $\p > 2$, dropping $\uptheta$, in later discussions.

\subsection{An equivariant Morse flow category}

We provide a concrete Morse model for equivariant cohomology of a point. We will use this Morse model throughout all the equivariant constructions.

\subsubsection{$\p=2$ case}

Consider the ``real $S^\infty$''
\beqn
S^\infty_{\mb R}:= \Big\{  \eta = (\eta_0,
 \eta_1, \ldots )\ |\ \eta_l \in {\mb R},\ \eta_l = 0 \text{ for }\forall l \gg 0,\ \sum_{l=0}^\infty \eta_l^2 = 1 \Big\}
\eeqn
which is a contractible infinite CW complex. The group ${\mb Z}_2$ acts freely by changing $\eta_l$ to $-\eta_l$ with quotient being $\mb{RP}^\infty$, which is a model for the classifying space of ${\mb Z}_2$. As $S_{\mb R}^\infty$ is the colimit of finite-dimensional spheres, we can define notion of smoothness in the naive way. For example, a function on $S_{\mb R}^\infty$ is smooth if its restriction to each truncation $S^k \subset S_{\mb R}^\infty$ is smooth. 

Consider the function 
\beqn
f^{\mb E}: S_{\mb R}^\infty \to {\mb R},\ f^{\mb E}(\eta) = \sum_{l=0}^\infty l \eta_l^2.
\eeqn
Its trunctions to each $S^k$ is a Morse function. The critical point set is 
\beqn
{\rm crit} f^{\mb E} = \{ \eta = \pm w_l \ |\ \eta_j = \pm \delta_{jl},\ l = 0, 1, \ldots \}.
\eeqn
It is easy to see that the standard metric $g^{\mb E}$ on $S^\infty$ makes it Morse--Smale. Define the grading by 
\beqn
|\pm w_l| = l.
\eeqn

The Tate equivariant theory suggests to consider a formal ``negative extension'' of $S_{\mb R}^\infty$. Define
\beqn
\wh S_{\mb R}^\infty = \Big\{  \eta = (\eta_l)_{l \in {\mb Z}}\ |\ \eta_l \in {\mb R},\ \eta_l = 0 \text{ for }\forall |l| \gg 0,\ \sum_{l\in {\mb Z}} \eta_l^2 = 1 \Big\}.
\eeqn
This is no longer a CW complex; however, all structures should be understood in the sense of limits/colimits. For each pair of nonnegative integers $m, n$, there is a finite truncation
\beqn
S^{m,n}_{\mb R}:= \Big\{ \eta \in \wh S_{\mb R}^\infty\ |\ l < -m\ {\rm or}\ l > n \Longrightarrow \eta_l = 0\Big\}.
\eeqn
By abuse of notation, we can extend the Morse function to
\beqn
f^{\mb E}: \wh S_{\mb R}^\infty \to {\mb R},\ f^{\mb E}(\eta) = \sum_{l=-\infty}^\infty l \eta_l^2.
\eeqn
Then with respect to the standard Riemannian metric, the pair is Morse--Smale on any such finite truncation. 

One can formally write down an equivariant Morse flow category ${\mb E}$ as follows. The set of objects ${\rm Ob} {\mb E}$ is the set of critical points of $f^{\mb E}$ on $\wh S^\infty_{\mb R}$. Between $x, y \in {\rm Ob}{\mb E}$, one has the morphism space being the compactified moduli space $M_{xy}^{\mb E}$ of gradient flow lines from $x$ to $y$. Then this is a flow category enriched in $\uds{\bf SMan}$. This flow category satisfies the unobstructedness condition (Definition \ref{defn_unobstructed}) and local finiteness condition (Definition \ref{defn_local_finite_flow}) with respect to the function $f^{\mb E}$. 

Moreover, the flow category is obviously equivariant with respect to a natural Novikov group ${\mb Z}$ action. Here on ${\mb Z}$ we define both the energy and the index function by 
\beqn
\lambda(n) = \mu(n) = n.
\eeqn
The associated Novikov field is
\beqn
\Lambda_{{\mb F}_2}^{\mb Z}:= {\mb F}_2 \lst = \tilde {\bf K}_2.
\eeqn
The ${\mb Z}$-action on ${\mb E}$ is induced by the ${\mb Z}$-action on $\wh S^\infty_{\mb R}$ 
\beqn
n \cdot (\eta_l)_{l\in {\mb Z}}:= (\eta_{l-n})_{l\in {\mb Z}}
\eeqn
which satisfies
\begin{align*}
& f^{\mb E}(n \cdot \eta) = f^{\mb E}(\eta) + \lambda(n),\ &\ |n \cdot \eta| = |\eta| +  \mu(n).
\end{align*}
In characteristic $2$, one does not need orientation. Then by the abstract construction of Floer complexes (see Subsection \ref{subsection_abstract_chain_complex} and Proposition \ref{prop_equivariant_differential}), one obtains a differential on
\beqn
C({\mb E}):= CM(f^\ep) =  \Big\{ \sum_{i=1}^\infty a_i w_i \ |\ a_i \in {\mb F}_2,\ w_i \in {\rm crit} f^{\mb E},\ \lim_{i \to \infty} f^{\ep}(w_i) = +\infty \Big\}
\eeqn
which is a $ \wt{\bf K}_2$-module.

One can compute the cohomology easily. For between any two critical points of $f^{\mb E}$ of adjacent indices, there is exactly one gradient trajectory between them. Hence in each degree $l$, $C^l({\mb E})$ is 2-dimensional while the cohomology is trivial. 

However, this trivial cohomology is not what we need. Notice that the function $f^{\mb E}$ (as well as the metric) is invariant under the natural ${\mb Z}_2$-action on $\wh{S}_{\mb R}^\infty$. One can consider the ${\mb Z}_2$-invariant part 
\beqn
C({\mb E})_{{\mb Z}_2} \subset C({\mb E})
\eeqn
which is 1-dimensional in each degree. The differential is trivial in characteristic $2$, hence the cohomology is also 1-dimensional. Therefore, one has 
\beqn
H ( C({\mb E})_{{\mb Z}_2}) \cong \Lambda_{{\mb F}_2}^{{\mb Z}} \cong \wt{\bf K}_2
\eeqn
as $\wt{\bf K}_2$-modules. Later we will show that the ring structure is also identical to $\wt{\bf K}_2$.

\subsubsection{$\p > 2$ case}

Consider the complex version of the infinite-dimensional sphere. Denote
\beqn
S_{\mb C}^\infty = \Big\{ \eta =  (\eta_0, \eta_1, \ldots)\ |\  \eta_l \in {\mb C},\ \eta_l = 0 \text{ for } \forall l \gg 0,\ \sum_{l=0}^\infty |\eta_l|^2 = 1 \Big\}.
\eeqn
This is a well-defined CW complex and is the universal $S^1$-bundle over $\mb{CP}^\infty$, a  model for the classifying space of $U(1)$. By regarding $\zp$ as a subgroup of $S^1$, $S_{\mb C}^\infty$ is also the universal $\zp$-bundle over $B\zp$. 

Now we define the formal negative extension
\beqn
\wh S_{\mb C}^\infty:= \Big\{ \eta = (\eta_l)_{l \in {\mb Z}}\ |\ \eta_l \in {\mb C},\ \eta_l = 0 \text{ for } \forall |l| \gg 0,\ \sum_{l \in {\mb Z}} |\eta_l|^2 = 1 \Big\}.
\eeqn
The $\zp$-action on $S_{\mb C}^\infty$ extends to $\wh S_{\mb C}^\infty$. Moreover, there is a free ${\mb Z}$-action on $\wh S_{\mb C}^\infty$ by
\beq\label{shift_map}
n\cdot (\eta_l) \mapsto (\eta_{l-n}).
\eeq

Define an $\zp$-invariant function $f^{{\mb E}}: \wh S_{\mb C}^\infty \to {\mb R}$ by 
\beq\label{fep}
f^{{\mb E}} (\eta) =  \left(  \sum_{l \in {\mb Z}} l |\eta_l|^2 \right) + \epsilon \sum_{ l \in {\mb Z}} {\rm Re}(\eta_l^{\p} ),
\eeq
where $\epsilon$ is a sufficiently small positive number. It is a smooth function in the sense that its restriction to each finite truncation is a smooth function. It also respects the ${\mb Z}$-action \eqref{shift_map} in the following way:
\beqn
f^{{\mb E}} (n \cdot \eta) = f^{{\mb E}} (\eta) + n.
\eeqn

If we set $\epsilon = 0$, $f^{{\mb E}}$ is a Morse--Bott function with critical submanifolds being the circles
\beqn
S_l^1 = \{ \eta \in S^\infty\ |\ |\eta_l| = 1\}.
\eeqn
After turning on $\epsilon$, the circle $S_l^1$ splits into nondegenerate critical points
\beqn
w_{2l}^0, \ldots, w_{2l}^{\p-1},\ w_{2l+1}^0, \ldots, w_{2l+1}^{\p-1}
\eeqn
where the subindices indicate the number of positive eigenvalues of the Hessian. One can explicitly write down the coordinates of critical points. For the ones with even degrees, we have
\beqn
w_{2l}^0 = (0, \ldots, 0, 1, 0, \ldots),\ w_{2l}^s = \tau^s w_{2l}^0,\ s = 1, \ldots, \p-1,
\eeqn
where the only nonzero coordinate of $w_{2l}^0$ is the $l$-th one and $\tau = e^{\frac{2\pi {\bf i}}{\p}}$ is the generator of $\zp$. For the ones with odd degrees, we have
\beqn
w_{2l+1}^0 = (0, \ldots, 0, e^{\frac{\pi {\bf i}}{\p}}, 0, \ldots, 0),\ w_{2l+1}^s = \tau^s w_{2l+1}^0,\ s = 1, \ldots, \p-1.
\eeqn
If we define the degree of $w_k^m$ be $k$, then the action by $n \in {\mb Z}$ in \eqref{shift_map} shifts the grading by $2n$.

\begin{lemma}\label{lemma191}
There exists a Riemannian metric $g^\ep$ on $\wh S_{\mb C}^\infty$ satisfying the following conditions.
\begin{enumerate}

\item $g^\ep$ is $\zp$-invariant and invariant under the shift map \eqref{shift_map}.

\item $(f^\ep, g^\ep)$ is Morse--Smale (when restricted to any finite truncation of $\wh S_{\mb C}^\infty$). 

\item If we consider the Morse complex associated to $(f^\ep, g^\ep)$ (restricted to any sufficiently large finite truncation of $\wh S_{\mb C}^\infty$), then the Morse differential, in terms of generators, is
\begin{align*}
&\ d(w_{2l}^i) = w_{2l+1}^i - w_{2l+1}^{i+1},\ &\ d(w_{2l-1}^i) = w_{2l}^0 + \cdots + w_{2l}^{\p-1}.
\end{align*}
\end{enumerate}
\end{lemma}

\begin{proof}
See \cite[Section 4]{shelukhin-zhao} which treated the case on $S_{\mb C}^\infty$. One can then use the ${\mb Z}$-action to extend the metric to the negative part of $\wh S_{\mb C}^\infty$. The expression of the differential can also be found in \cite[Section 2]{Seidel_Wilkins_2022}, and we refer the reader to \cite[Appendix A]{BSWX} for a comparison between the Morse model and the cellular model used in \cite{Seidel_Wilkins_2022}.
\end{proof}

Therefore, one can define a flow category $\ep$ whose objects are critical points of $f^\ep$ and whose morphism spaces are the moduli space of (broken) trajectories between critical points. One can make $\ep$ a ${\mb Z} \rtimes \zp$-equivariant Novikov flow category. Indeed, define
\begin{align*}
&\ \lambda: {\mb Z} \to {\mb R},\ \lambda(a) = a,\ &\ \mu: {\mb Z} \to {\mb Z},\ \mu(a) = 2a.
\end{align*}
They are $\zp$-invariant. Hence the triple $(\zp \times {\mb Z}, \lambda, \mu)$ is a Novikov group (Definition \ref{defn_Novikov_group}). Next, one can make the Morse flow category $\ep$ associated to $(f^\ep, g^\ep)$ a $\zp \times {\mb Z}$-equivariant Novikov flow category. Define the action and grading by
\begin{align*}
&\ {\mc A}(p) = f^\ep(p),\ &\ |w_k^s| = k.
\end{align*}
Then one can see 
\begin{align*}
&\ {\mc A}(ap) = f^{\ep}(p) + \lambda(a),\ &\ |ap| = |p| + \mu(a).
\end{align*}
Then, by the general construction of complexes in Section \ref{section_algebraic}, there is an associated ${\mb Z}$-graded complex
\beqn
C({\mb E}):= CM(f^\ep) = \Big\{ \sum_{i=1}^\infty a_i w_i\ |\ a_i \in \fp,\ w_i \in {\rm crit} f^\ep,\ \lim_{i \to \infty} {\mc A} (w_i) = +\infty \Big\},
\eeqn
which is a graded module over $\Lambda_{\fp}^{\mb Z} \cong \kp$.

One can compute the chain complex easily from Lemma \ref{lemma191}. One can see that the total cohomology of $C({\mb E})$ is trivial. If we consider the $\zp$-invariant part, then in degree $k$, one has 
\beqn
C^k (\ep)_{\zp} = \fp \langle w_k^0 + \cdots + w_k^{\p-1} \rangle
\eeqn
and the differential is trivial in characteristic $\p$. Hence the cohomology of the $\zp$-fixed part is identified with
\beqn
H( C({\mb E})_{\zp}) \cong \tkp
\eeqn
as $\kp$-modules. We will identify below its multiplicative structure using Morse flow trees.

\subsubsection{Equivariant cohomology of free orbits}

We consider a simple case of equivariant localization. Let $A$ be a finite set equipped with a $\zp$-action. Then one has a trivial product $A \times {\mb E}$ as a flow category whose object set is $A \times {\rm Ob}{\mb E}$. Then one can form the complex
\beqn
C(A \times {\mb E}):= \bigoplus_{\alpha \in A} C({\mb E})
\eeqn
which has a $\zp$-action. 

\begin{lemma}\label{lemma_localization}
Let $A_{\zp} \subset A$ be the subset of $\zp$-fixed points. Then 
\beqn
H(C(A \times {\mb E})_{\zp}) \cong \tkp^{A_{\zp}}.
\eeqn
\end{lemma}

\begin{proof}
One only needs to consider the cases when $A$ contains a single element and when $A$ is a $\zp$-torsor. In the former case, the cohomology is $\tkp$. In the latter case, one can explicitly calculate the differential and show that the cohomology is trivial.
\end{proof}

\subsection{Ring structure using Morse model}

We define the Morse-theoretic multiplicative structure on $H_{\zp} ({\rm pt}) \cong \tkp$. In both $\p= 2$ and $\p>2$ cases, we denote by $S^\infty$ the real or complex infinite-dimensional sphere and $\wh S^\infty$ its negative extension.

Let $H_{\zp}^{\rm Borel}({\rm pt})$ be the Borel equivariant cohomology in $\fp$-coefficients, which can be obtained from the above Morse model by truncating at the nonnegative levels, i.e., via a Morse flow category on $S^\infty \subset \wh S^\infty$. We will use Morse flow trees in $S^\infty$ to define a ring structure on $H_{\zp}^{\rm Borel}({\rm pt})$ and then extend to the Tate version.

Let $T_{2, 1}$ denote the tree which is the disjoint union of two incoming edges $T_1\cong T_2 \cong (-\infty, 0]$ and an outgoing edge $T_\infty \cong [0, +\infty)$. Let $s \in T_{2,1}$ denote the coordinate on each edge; when there is no ambiguity, $s$ corresponds to a real number. To achieve transversality, we choose a $\zp$-invariant family of vector fields $Y_s$ on $S^\infty$ parametrized by $s \in T_{2,1}$
\begin{enumerate}

\item $Y (s) = 0$ when $|s|\gg 0$. 

\item For each finite truncation $S^n \subset S^\infty$, $Y|_{S^n}$ is tangent to $S^n \subset S^\infty$ and $\nabla Y$ is sufficiently small in the normal direction of $S^n$.
\end{enumerate}
Let $Y_i$ denote the restriction of $Y$ to $T_i \subset T_{2,1}$. A {\bf $Y$-perturbed flow tree} is a triple $\eta = (\eta_1, \eta_2, \eta_\infty): T_{2,1} \to S^\infty$ which solves the equation
\beqn
\eta'(s) = \nabla f^\ep (\eta(s)) + Y_s(\eta(s))
\eeqn
subject to the incidence relation
\beqn
\eta_1(0) = \eta_2(0) = \eta_\infty(0).
\eeqn
Each flow tree converges to a triple of critical points of $g$ at $\infty$ of the edges.

Given a triple of critical points $w_1, w_2, w_\infty \in {\rm crit}(f^\ep) \cap S^\infty$, let $S^n \subset S^\infty$ be the smallest truncation which contains all of them. One can see that a flow tree connecting $w_1, w_2, w_\infty$ is contained in $S^n$. We say that a $Y$-perturbed flow tree is regular for $(w_1, w_2, w_\infty)$ if it is regular as a flow tree contained in $S^n$.\footnote{In fact if it is regular in $S^n$, then it is regular in any larger truncation.}
 
For a generic perturbation $Y$ and each triple of critical points $w_1, w_2, w_\infty$, one obtains a regular moduli space of solutions, whose mod $p$ count induces a $\zp$-equivariant map
\beqn
CM( f^\ep|_{S^\infty} ) \otimes CM(f^\ep|_{S^\infty} ) \to CM( f^\ep|_{S^\infty} ).
\eeqn
Hence after restriction to the $\zp$-invariant part, one obtains
\beqn
H_{\zp}^{\rm Borel} ({\rm pt}) \otimes_{\fp} H_{\zp}^{\rm Borel}({\rm pt})  \cong CM(f^\ep|_{S^\infty})^{\zp}\otimes_{\fp} CM(f^\ep|_{S^\infty})^{\zp} \to CM_{\zp} (f^\ep|_{S^\infty})^{\zp}\cong H_{\zp}^{\rm Borel}({\rm pt}).
\eeqn 
The usual Morse theoretic technique can be used to show that the bilinear map is independent of the choices of $Y$, and satisfies associativity. Hence $H_{\zp}^{\rm Borel}({\rm pt})$ is equipped with a structure of a ${\mb Z}$-graded, associative and supercommutative algebra over $\fp$. 

We know for compact manifolds, the ring structure defined by counting Morse flow trees coincides with the ordinary cohomology ring. Applying this fact to any finite truncation $S^n$, one can see the above Morse-theoretic definition of the ring structure coincides with the ordinary one on 
\beqn
H^{\rm Borel}_{\zp}({\rm pt}) = \left\{ \begin{array}{ll} \fp \pst,\ &\ \p = 2,\\
\fp \psu \langle \uptheta \rangle,\ &\ \p > 2 \end{array}\right.
\eeqn
By inverting $\uptheta$ or $\tu$, one obtains a ring structure on the Tate equivariant cohomology of a point, which coincides with the standard one on $\tkp$.

Later we will construct cohomology groups which become modules over $\tkp$. The module structure is always defined via a certain Morse-theoretic construction first as a module over the positive truncation, then extended by inverting the corresponding variable $\uptheta$ or $\tu$.

\section{Equivariant Morse and Floer Cohomology}

The purpose of this section is to discuss both an algebraic and Morse-theoretic approach to equivariant cohomology. The former is used for constructing the quasi-Frobenius map (cf. Equation \eqref{quasi_Frobenius}) and the latter uses the Morse flow category discussed in the previous section to define the $\uptheta$-action on equivariant cohomology.

\subsection{Tate Equivariant Cohomology}

\subsubsection{The algebraic construction}

Fix a prime $\p$. The coefficient ring will be $\fp$ unless otherwise declared. Let $D$ be a $\zp$-equivariant cochain complex with differential $d_D$. Let $\tau: D \to D$ be the action of the generator $\tau \in \zp$. Define cochain complexes of $\kp$-modules
\beqn
C_{\zp} (D):=D \otimes \tkp \cong \left\{ \begin{array}{ll} D \otimes {\bf F}_2 \lst,\ &\ \p = 2,\\
                 D \otimes {\bf F}_\p \lsu \langle \uptheta \rangle,\ &\ \p > 2.\end{array} \right.
\eeqn
where the differential $d_{\zp}$ is defined as follows.

\begin{enumerate}

\item When $\p = 2$, $d_{\zp}$ is $\uptheta$-linear and is determined by 
\beq\label{equivariant_differential_1}
d_{\zp}(x\otimes \uptheta^k) = d_D(x) \otimes \uptheta^k + ( 1 - \tau) (x) \otimes \uptheta^{k+1}.
\eeq

\item When $\p > 2$, $d_{\zp}$ is $\tu$-linear and is determined by 
\beq\label{equivariant_differential_2}
\begin{split}
d_\zp (x\otimes 1) &\ = d_D(x) \otimes 1 + (-1)^{|x|} (1-\tau)(x) \otimes \uptheta,\\
d_\zp (x \otimes \uptheta) &\ = d_D(x) \otimes \uptheta + (-1)^{|x|} ( 1+ \tau + \cdots + \tau^{p-1})(x) \otimes \tu.
\end{split}
\eeq
Here $|x|$ is the cohomological grading of $x$.
\end{enumerate}
The resulting cohomology, called the {\bf (Tate) equivariant cohomology} of $D$, is denoted by 
\beqn
H_{\zp} (D)
\eeqn
which is naturally a $\kp$-module. 
In fact the equivariant cohomology $H_{\zp} (D)$ can be made a module over $\tkp$ in a non-obvious way (unless $\p = 2$). We choose to construct the module structure in the Morse-theoretic approach. We will not use the above algebraic definition formally, but only provide it for algebraic intuition.

\begin{rem}
    If we work with $D \otimes \kp$ equipped with the same differentials \eqref{equivariant_differential_1} and \eqref{equivariant_differential_2}, then we obtain the Borel-equivariant cohomology. We mostly stick with the Tate setting because of intended applications.
\end{rem}

The equivariant complex is functorial. Namely, if $f: D \to D'$ is a $\zp$-equivariant chain map, then it induces a chain map
\beqn
f_{\zp}: C_{\zp}(D) \to C_{\zp}(D').
\eeqn
If $f_0$ and $f_1$ are $\zp$-equivariantly chain homotopic, then there is a corresponding $\kp$-linear chain homotopy from $f_{0, \zp}$ to $f_{1, \zp}$.

There are two special cases of the above construction which will appear in this paper soon. The first special case is when the $\zp$-action on the complex $D$ is trivial. In this case
\beqn
H_{\zp} (D) \cong H (D) \otimes \tkp.
\eeqn
Notice that this is a priori only an isomorphism of $\kp$-modules but not $\tkp$-modules.

The second special case is the $\p$-fold tensor power. Let $D$ be a non-equivariant cochain complex of $\fp$-vector spaces. The tensor power
\beqn
( D^{\otimes \p}, d_D^{\otimes \p})
\eeqn
has the $\zp$-action by permuting the factors as follows: let $\tau\in \zp$ be the generator, then $\tau$ acts on the $\p$-fold tensor product by
\beqn
\tau( x_1 \otimes \cdots \otimes x_\p) = (-1)^{|x_1| ( | x_2| + \cdots + |x_\p|) } x_2 \otimes \cdots \otimes x_\p \otimes x_1.
\eeqn
Then  one has the chain complex
\beqn
C_{\zp} (D^{\otimes \p}).
\eeqn

One can calculate the Tate equivariant cohomology by the so-called quasi-Frobenius map. This map is not induced from any chain-level construction. We first define the {\bf quasi-Frobenius map}
\beq\label{quasi_Frobenius}
qF: H  (D) \to H_{\zp} (D^{\otimes \p})
\eeq
in the case when the differential of $D$ is trivial. In this case, $H_{\zp}(D) \cong H(D) \otimes \tkp$. Define
\beqn
qF(x) = \underbrace{x\otimes \cdots \otimes x}_{\p} \in C_{\zp}(D^{\otimes \p})\ \forall x \in H(D) \subset H_{\zp}(D).
\eeqn
The image is contained in the kernel of $d_{\zp}$ hence induces a map to $H_{\zp}(D^{\otimes \p})$. Moreover, we can see that the cohomological map is linear over $\fp$. 

In general, when the differential of $D$ is nontrivial, there exists a $\zp$-equivariant chain homotopy equivalence
\beqn
(H(D), 0) \sim (D, d_D).
\eeqn
Then we define $qF$ as the composition
\beqn
\xymatrix{ H(D) \ar[r] & H_{\zp}(H(D)^{\otimes \p}) \ar[r] & H_{\zp}(D^{\otimes \p})  }
\eeqn
Moreover, any two $\zp$-equivariant chain homotopy equivalences between $H(D)$ and $D$ are equivariantly homotopic. Hence the resulting map 
\beqn
qF: H(D) \to H_{\zp}(D^{\otimes \p})
\eeqn
is well-defined. Furthermore, the following functoriality property of $qF$ holds.

\begin{lemma}
The quasi-Frobenius map is functorial. Namely, if $f: D_1 \to D_2$ is a cochain map, then there is a commutative diagram
\beq\label{quasi_Frobenius_naturality}
\vcenter{ \xymatrix{ H (D_1) \ar[rr]^-{qF} \ar[d]_f & & H_{\zp} (D_1^{\otimes \p}) \ar[d]^{f_\zp} \\
        H (D_2) \ar[rr]_-{qF} & &  H_{\zp}(D_2^{\otimes \p})}}
\eeq
\end{lemma}
This functoriality implies that the quasi-Frobenius map is an isomorphism with field coefficients. 

\begin{prop}\label{prop:qf-iso}
There is a $\kp$-linear extension of $qF$ which is an isomorphism
\beqn
H(D) \otimes \tkp \to H_{\zp}(D^{\otimes \p}).
\eeqn
\end{prop}

\begin{proof}
By the definition of $qF$ and the above lemma, one only needs to verify this for the case when the differential of $D$ is trivial. The fact that $qF$ induces an isomorphism follows from \cite[Lemma 3.1]{shelukhin-zhao}.
\end{proof}

\subsubsection{Novikov coefficients}

Now we extend the above algebraic constructions to complex of modules over Novikov fields. First recall that when $\Pi_1, \Pi_2$ are Novikov groups, then one has a completed tensor product
\beqn
\Lambda^{\Pi_1}\widehat{\otimes} \Lambda^{\Pi_2} \cong \Lambda^{\Pi_1 \oplus \Pi_2}.
\eeqn
Then if $D_1$ resp. $D_2$ are non-Archimedean normed over $\Lambda^{\Pi_1}$ resp. $\Lambda^{\Pi_2}$, then one has a completed tensor product
\beqn
D_1 \widehat{\otimes} D_2
\eeqn
which is a non-Archimedean normed vector space over $\Lambda^{\Pi_1 \oplus \Pi_2}$.

The above equivariant construction requires the following specialization. Let $(\Pi, \lambda, \mu)$ be a Novikov group, which is typically the one defined by \eqref{Novikov_group}. In the equivariant setting we consider another standard one which is 
\beqn
({\mb Z}, \lambda_\p, \mu_\p)
\eeqn
where $\lambda_\p: {\mb Z} \to {\mb R}$ is the canonical inclusion while
\beqn
\mu_\p(a) = \left\{ \begin{array}{cc} a,\ &\ \p = 2,\\
2a,\ &\ \p > 2 \end{array}\right.
\eeqn

Now let $D$ be a complex of $\Lambda_{\fp}^\Pi$-vector spaces. Then the completed tensor product $D \widehat{\otimes} \tkp$ is a vector space over the Novikov field 
\beqn
\Lambda^{\Pi \oplus {\mb Z}}:= \Lambda^{\Pi \oplus {\mb Z}}_{\fp} = \left\{ \sum_{i=1}^\infty c_i T^{a_i + b_i}\ \left| \ \begin{array}{l} c_i \in \fp,\ a_i \in \Pi,\ b_i \in {\mb Z},\\
\displaystyle \inf_i \omega(a_i) > -\infty,\ \inf_i b_i > -\infty,\ \lim_{i +\infty} \omega(a_i) + b_i = +\infty \end{array} \right.\right\}
\eeqn
We typically write $T^{b_i}$ as $\uptheta^{b_i}$ (when $\p =2$) or $\tu^{b_i}$ (when $\p > 2$). Then this Novikov field consists of formal infinite $\fp$-linear combinations of monomials $T^{a_i} P_i( \uptheta)$ or $T^{a_i} P_i(\tu)$ where $P_i$ is a Laurent series with bounded pole orders while $\omega(a_i) \to +\infty$. 

We can define the equivariant differential on the tensor product $D \widehat{\otimes} \tkp$ using the same formula as before. However, most of the time, we will enlarge the Novikov field via
\beqn
\Lambda_{\fp}^{\Pi \oplus {\mb Z}} \hookrightarrow \Lambda_{\kp}^{\Pi}
\eeqn
where the latter does not require bounded pole orders for the coefficient of $T^{a_i}$. Then define 
\beqn
\left( C_{\zp} (D):= D  \widehat{\otimes} \tkp, d_{\zp} \right)
\eeqn
where the differential is still defined in the same way. 

One often needs to take a more convenient field extension. Consider 
\beqn
\Lambda_{\kp}^\Pi
\eeqn
where $\kp$ is the field of Laurent series in $\uptheta$ (for $\p=2$) or $\tu$ (for $\p > 2$). Elements in $\Lambda_{\kp}^\Pi$ are formal infinite $\fp$-linear combinations of monomials $T^{a_i} P_i$ where $P_i$ is a Laurent series without bound on pole orders, or equivalently, formal infinite $\kp$-linear combinations of $T^{a_i}$ with $\omega(a_i) \to +\infty$. Then one often needs to use
\beqn
C_{\zp} (D) \underset{\Lambda^{\Pi \oplus {\mb Z}}}{\otimes} \Lambda_{\kp}^\Pi.
\eeqn

We can also consider $C_{\zp}(D^{\otimes \p})$ when $D$ is a non-equivariant cochain complex of $\Lambda^\Pi$-modules. We consider the completed tensor power (over $\fp$)
\beqn
\overbrace{ D  \widehat{\underset{\fp}{\otimes}} \cdots \widehat{\underset{\fp}{\otimes}} D }^{\p}
\eeqn
which is a module over 
\beqn
\Lambda^{\Pi^\p}:= \Lambda^{\Pi} \widehat{\underset{\fp}{\otimes}} \cdots \widehat{\underset{\fp}{\otimes}} \Lambda^{\Pi}.
\eeqn
Consider the natural morphism of Novikov monoids (see Definition \ref{defn_Novikov_group})
\beqn
\xymatrix{ \Pi \oplus \cdots \oplus \Pi \ar[r] \ar[d] & {\mb R}^{k} \oplus \cdots \oplus {\mb R}^k \ar[d]\\
    \Pi \ar[r] & {\mb R}^k}
\eeqn
where the vertical arrows are $(a_1, \ldots, a_\p) \mapsto a_1 + \cdots + a_\p$. Then there is an induced ring morphism
\beqn
\Lambda_R^{\Pi^{\p}} \to \Lambda_R^{\Pi}.
\eeqn
Denote
\beqn
D^{\widehat{\otimes} \p}:= \left( \overbrace{ D \widehat{\underset{\fp}{\otimes}} \cdots \widehat{\underset{\fp}{\otimes}} D }^{\p} \right) \underset{ \Lambda^{\Pi^{\p}}}{\otimes} \Lambda^{\Pi}
\eeqn
which is a module over $\Lambda_R^{\mb\Pi}$. Then we have the chain complex
\beqn
C_{\zp} ( D^{\widehat{\otimes} \p} ),
\eeqn
whose differential is defined similarly using \eqref{equivariant_differential_1} and \eqref{equivariant_differential_2}. In this situation, we still have the quasi-Frobenius map well-defined on the cohomology level
\beqn
qF: H^*(D) \to H_{\zp}^*( D^{\widehat{\otimes} \p}).
\eeqn
However it is only linear over $\fp$ but not the Novikov field $\Lambda_{\fp}^\Pi$. Instead, one has Frobenius linearity in the Novikov variable
\beqn
qF( T^{a} x ) = T^{\p a} qF(x). 
\eeqn

\subsubsection{Quantitative results}

One considers a key quantitative property of the quasi-Frobenius map. Let $D$ be a Floer-type complex over $\Lambda_{\fp}^\Pi$ (Definition \ref{defn_Floer_type_complex}). Then the tensor product
\beqn
D^{\widehat{\otimes} \p}
\eeqn
is also a Floer-type complex over $\Lambda_{\fp}^\Pi$. The equivariant complex, after tensoring with $\Lambda_{\kp}^\Pi$, becomes a Floer-type complex over this larger Novikov field
\beqn
C_{\zp}(D^{\widehat{\otimes} \p}) \otimes \Lambda_{\kp}^\Pi.
\eeqn
Then for each $\tau \in {\mb R}$, define the filtered cohomology
\beqn
H_{\zp}^{\leq \tau}:= H( (C_{\zp}(D^{
\widehat{\otimes} \p}) \otimes \Lambda_{\kp}^\Pi)^{\leq \tau}).
\eeqn
Via localization, one can prove the following result.

\begin{lemma}
For each $\tau$, there exists a quasi-Frobenius map
\beqn
qF: H^{\leq \tau}(D) \to H_{\zp}^{\leq \p \tau}(D^{\otimes \p})
\eeqn
satisfying
\begin{enumerate}

    \item For $\tau < \sigma$, the following diagram commutes
    \beqn
    \xymatrix{ H^{\leq \tau}(D) \ar[r]^-{qF} \ar[d]   & H_{\zp}^{\leq \p \tau}(D^{\otimes \p}) \ar[d] \\
              H^{\leq \sigma}(D) \ar[r]_-{qF}  & H_{\zp}^{\leq \p \sigma}( D^{\otimes \p})}
    \eeqn

    \item After tensoring with $\tkp$, for each $\tau$, $qF$ induces an isomorphism of $\Lambda_{\kp, 0}^\Pi$-modules
    \beqn
    qF: H^{\leq \tau}(D) \otimes \tkp \to H_{\zp}^{\leq \p \tau}( D^{\otimes \p}).
    \eeqn
\end{enumerate}
\end{lemma}

\begin{proof}
See \cite{shelukhin-zhao}. 
\end{proof}

\subsection{$\zp$-equivariant Morse complexes}

\subsubsection{Equivariant Morse cohomology for $\zp$-manifolds}

Let $M$ be a compact oriented manifold equipped with a $\zp$-action. Let $f: M \to {\mb R}$ be a $\zp$-invariant Morse function. One would like to define a $\zp$-equivariant Morse cochain complex $CM_{\zp} (f)$. 

First we assume that there exists a $\zp$-invariant Riemannian metric $g$ on $M$ making $(f, g)$ Morse--Smale. Namely, one can achieve equivariant transversality. Then one obtains a cochain complex $CM (f, g)$ which has a $\zp$-action. Hence the previous algebraic construction provides a way to obtain an equivariant Morse complex
\beqn
CM_{\zp} (f):= C_{\zp} ( CM^*(f, g)).
\eeqn

In fact the equivariant complex can also be obtained via Morse model. Consider the product function
\beqn
F: M \times \wh S^\infty \to {\mb R},\ F(x, w) = f(x) + f^{\mb E}(w)
\eeqn
and the product metric, which provides a Morse flow category ${\mb M}\#\ep$. Consider the associated Morse cochain complex
\beqn
CM (F)
\eeqn
which coincides with the tensor product $CM (f) \otimes {\mb E}$. This complex has a $\zp$-action given by the $\zp$-action on $CM (f)$ and the $\zp$-action on ${\mb E}$. Then define
\beqn
CM_{\zp} (f):= (CM (f) \otimes {\mb E})_{\zp}
\eeqn
i.e., the $\zp$-fixed part. 

\begin{lemma}
The complex $CM_{\zp} (f)$ coincides with the algebraically defined complex $C_{\zp} (D)$.
\end{lemma}

\begin{proof}
Define a linear map
\beqn
CM(f) \otimes \tkp \to CM(F)_{\zp}
\eeqn
as follows.
\begin{enumerate}

\item When $\p = 2$, define 
\beqn
x\otimes \uptheta^k \mapsto x \otimes q_k^+ + \tau(x) \otimes q_k^-
\eeqn
where $q_k^\pm$ are the two critical points of $f^\ep$ and $\tau$ is the action on $CM(f)$.

\item When $\p > 2$, define
\beqn
\begin{split}
    x \otimes \tu^k & \mapsto \sum_{l=1}^{p} (\tau^l x) \otimes w_{2k}^l,\\
    x \otimes \tu^k \uptheta & \mapsto \sum_{l=1}^p (\tau^l x) \otimes w_{2k+1}^l
    \end{split}
\eeqn
where $w_i^l$ are the critical points of $f^\ep$. 
\end{enumerate}
One can check that this linear map is a bijective cochain map.
\end{proof}

The two special cases mentioned above also arises in the Morse-theoretic context. First, suppose $\zp$ acts trivially on $M$. Then a Morse--Smale pair $(f, g)$ is equivariantly transverse, resulting in the identification
\beqn
HM_{\zp} (f, g) \cong HM (f, g) \otimes \tkp.
\eeqn
On the other hand, if $(f, g)$ is Morse--Smale on $M$, then the $\p$-fold power
\beqn
(f^{\p}, g^{\oplus \p})
\eeqn
is Morse--Smale on $M^{\p}$. One obtains the cochain complex
\beqn
CM (f^\p) \otimes C({\mb E})
\eeqn
whose $\zp$-invariant has cohomology $HM_{\zp} (f^\p, g^{\otimes \p})$. When $M$ is compact, Morse-theoretic technique proves that the resulting cohomology is independent of the choice of $(f, g)$ and coincides with $H_{\zp} (M^\p)$.

Now assume that one does not have equivariant transversality. Let $M$ be equipped with a $\zp$-action, $f: M \to {\mb R}$ is a $\zp$-invariant Morse function. Still consider the product
\beqn
F: M \times \wh S^\infty \to {\mb R}.
\eeqn
Now the $\zp$-action on $M \times \wh S^\infty$ is free, allowing us to achieve equivariant transversality. There are many concrete choices of perturbations. For example, one can perturb the function $f$ to allow it to depend on $w \in \wh S^\infty$, i.e., a $\zp$-invariant function $(w, x) \mapsto f_v(x)$ on $\wh S^\infty \times M$ such that $f_v(x) = f(x)$ when $(v, x)$ is close to a critical point of $g \times f$. One can further ensure that the set of critical points is unchanged under the perturbation provided that the perturbation is sufficiently small in $C^1$-norm. Then one can consider the coupled flow line equation for $(\xi, \eta): {\mb R} \to M \times \wh S^\infty$:
\begin{align*}
    &\ \eta'(s) = \nabla f_{\zeta (s)} (\eta (s)),\ &\ \zeta'(s) = \nabla g(\zeta (s)),
\end{align*}
which has again an obvious translation symmetry. Define the Morse cochain complex as the (completed) tensor product
\beqn
CM (f) \otimes CM (g)
\eeqn
whose differential counts rigid perturbed flow lines. Define
\beqn
CM_{\zp}(f):= (CM (f) \otimes CM (g))^{\zp}
\eeqn
as a subcomplex, whose chain homotopy type is independent of the perturbation. 

By using perturbed flow trees in $\wh S^\infty$, one can also define a module structure of $HM_{\zp}(f)$ over $\tkp \cong H_{\zp}({\rm pt})$ induced by a chain-level bilinear map
\beqn
CM_{\zp} (f) \otimes_{\fp} CM _{\zp} (g) \to CM_{\zp}(f).
\eeqn
More precisely, choose a perturbation $Y$ over the tree $T_{2,1}$ which can be used to define the ring structure on $\wh S^\infty$. Choose another generic family of functions 
\beqn
f_{s, v}: M \to {\mb R}
\eeqn
parametrized by $s \in {\mb R}$ and $v \in S^\infty$ such that 1) $f_{s, v} = f$ near $s = 0$ and 2) $f_{s, v} = f_v$ for $|s|\gg 0$ and 3) $f_{s, \gamma v} (\gamma x) = f_{s, v}(x)$ for all $\gamma \in \zp$. Consider pairs $(\xi, \eta)$ where $\eta = (\zeta_1, \zeta_2, \zeta_\infty): T_{2,1} \to S^\infty$ is a $Y$-perturbed flow tree and $\xi: {\mb R} \to M $ solves the equation 
\beqn
\xi'(s) = \left\{ \begin{array}{cc}  \nabla f_{s, \zeta_1(s)}(\xi(s)),\ &\ s \leq 0,\ \\
 \nabla f_{s, \zeta_\infty(s)}(\xi(s)),\ &\ s \geq 0 \end{array}  \right.
\eeqn
The coefficients of the equation on $\xi$ is smooth. We then impose critical point asymptotics $\xi(\pm \infty) \in {\rm crit}(f)$. By choosing a generic $f_{s, v}$, one obtains equivariant transversality and hence defines a chain map
\beqn
\big( CM(f) \otimes CM(g) \big) \otimes CM(g) \to CM(f) \otimes CM(g).
\eeqn
Its restriction to the $\zp$-invariant part gives
\beqn
CM_{\zp}(f) \otimes CM_{\zp}(g) \hookrightarrow (CM(f) \otimes CM(g) \otimes CM(g))^{\zp} \to CM_{\zp}(f)
\eeqn
hence a bilinear map
\beqn
HM_{\zp}(f) \otimes_{\fp} \tkp \to HM_{\zp}(f).
\eeqn
One can further show that this map defines a $\tkp$-module structure on $HM_{\zp}(f)$ using a standard TQFT argument.

One can further define a multiplicative structure on $HM_{\zp}(f)$ via the cochain-level map
\beqn
CM_{\zp} (f) \otimes_{\fp} CM_{\zp}(f) \to CM_{\zp}(f)
\eeqn
via gradient flow trees in $M \times \wh S^\infty$, whose associativity can be verified in the standard way. Lastly, because the multiplicative structure respects the $\tkp$-module structure, hence it descends to a map
\beqn
HM_{\zp}(f) \underset{\tkp}{\otimes} HM_{\zp}(f) \to HM_{\zp}(f).
\eeqn

\section{Quantum Steenrod Operations---Proof of Theorem \ref{thma_QST}}\label{section_Steenrod}

The classical Steenrod operation can be viewed as a $\zp$-equivariant extension of the map $x \mapsto x^\p$ in cohomology in characteristic $\p$. Proposed by Fukaya \cite{Fukaya_Morse}, the quantum Steenrod operation deforms the classical one in analogy to the case that the quantum cup product deforms the classical cup product. Under topological restrictions such as monotonicity, the quantum Steenrod operation has been constructed by Seidel \cite{Seidel_pants}, Wilkins \cite{Wilkins_2020}, and Seidel--Wilkins \cite{Seidel_Wilkins_2022} (see also the recent survey by Wilkins \cite{Wilkins_survey}).

The use of the FOP transverse perturbations allows us to define integer, hence $\fp$-counts of pseudoholomorphic curves, leading to a definition of the quantum Steenrod operation on all compact symplectic manifolds. In this section we go over the construction using the Morse model.

\subsection{Geometric setup}

\subsubsection{The flow categories and bimodules associated to quantum Steenrod operation}

Fix a prime $\p$. Let $(X, \omega)$ be a compact symplectic manifold and $J$ be a compatible almost complex structure. We abbreviate the data $(X, \omega, J)$ by $X$. We can reduce the consideration to a fixed $J$ using argument such as continuation maps. Therefore we will not address the dependence on $J$. Recall
\beqn
\Pi = \frac{\pi_2(X)}{{\rm ker} \omega \cap {\rm ker} c_1}.
\eeqn

Let $(f, g)$ be a Morse--Smale pair on $X$ which defines a Morse flow category $\uds {\mb M}$. Let ${\mb M}$ be the trivial product $\uds {\mb M}\times \Pi$, which is naturally a $\Pi$-equivariant flow category. Its objects are denoted by 
\beqn
x = (\uds x, a_x) \in {\rm Ob}\uds{\mb M}\times \Pi.
\eeqn

Let $\mb{E}$ be the ${\mb Z}\times \zp$-equivariant flow category associated to the Morse--Smale pair on $\wh S^\infty$. Let $\uds {\mb M} \# {\mb E}$ be the ``coupling'' of $\uds {\mb M}$ and ${\mb E}$, which is a flow category with object set
\beqn
{\rm Ob}(\uds {\mb M}\# {\mb E}) = {\rm Ob}\uds {\mb M}\times {\mb E}.
\eeqn
The morphism space between two objects is the moduli space of flow lines in $X \times \wh S^\infty$. Denote
\beqn
\mb{M}\# {\mb E} = (\uds{\mb M}\# {\mb E}) \times \Pi.
\eeqn
This flow category is equivariant with respect to the free action by 
\beqn
\mb\Pi_\p = \Pi \times {\mb Z} \times \zp =: {\mb \Pi} \times {\mb Z}_{\p}.
\eeqn

We describe a multimodule ${\mb X}^{\mb{QSt}}$ over 
\beqn
(\underbrace{\mb{M}, \ldots, \mb{M}}_{\p}, {\mb E}; {\mb M}\# {\mb E})
\eeqn
as follows. Consider the Riemann surface $\Sigma^{\mb{QSt}}$ which is identified with $\mb{CP}^1$ with $\p+1$ marked points, where the $\p$ negative markings are at all $\p$-th roots of unity 
\beqn
z_1, \ldots, z_\p \in {\mb C} \subset {\mb C} \cup \{\infty\} = \mb{CP}^1
\eeqn
and the positive marking is $z_\infty = \infty$.

Consider pairs
\beqn
(u, v)
\eeqn
where $u: \Sigma^{\mb{QSt}} \to X$ is a $J$-holomorphic sphere in $X$ with these $\p +1$ special markings and $v: {\mb R} \to S^\infty$ is a Morse trajectory of $\ep$ in $S^\infty$.  For each tuple of objects 
\beqn
(x_1, \ldots, x_\p, w) \in {\rm Ob}{\mb M}\times \cdots \times {\rm Ob}{\mb M}\times {\rm Ob}{\mb E},\ p' = (x', w') \in {\rm Ob}{\mb M}\# {\mb E}
\eeqn
let
\beqn
\mathring M_{x_1 \cdots x_\p w;p'}^{\mb{QSt}}
\eeqn
be the space of such pairs $(u, v)$ such that 
\beqn
\int_{S^2} u^* \omega =  a_{x_1} + \cdots + a_{x_\p} - a_{x'}
\eeqn
satisfying the constraints
\beqn
u(z_i) \in W^{\rm us}_{\uds x_i},\ i = 1, \ldots, \p,\ u(z_\infty) \in W^{\rm st}_{{\uds x}'},\ {\rm and}\ \lim_{s \to -\infty} v(s) = w,\ \lim_{s \to +\infty} v(s) = w'.
\eeqn
Let $M_{x_1 \cdots x_\p w; p'}^{\mb{QSt}}$ be the compactification of $\mathring M_{x_1 \cdots x_\p w; p'}^{\mb{QSt}}$ compactified in the following way: at markings, one takes the ordinary Gromov compactification; the Morse breakings at negative ends are allowed independently at different edges while the breakings on the positive ends are required to be ``synchronized.''

\begin{prop}
The spaces $M_{x_1 \cdots x_\p w; p'}^{\mb{QSt}}$ define a Novikov multimodule, denoted by ${\mb X}^{\mb{QSt}}$, over $({\mb M}, \ldots, {\mb M}, \ep; {\mb M}^\zp)$, which is equivariant with respect to the morphism of Novikov groups
\beq\label{eqn201}
\begin{tikzcd}
\overbrace{\Pi \times \cdots \times \Pi}^{\p} \times {\mb Z} \arrow[rr] \arrow[d, "{(A_1, \ldots, A_{\p}, a) \mapsto (A_1 + \cdots + A_{\p}, a)}"'] & & {\mb R}^{\p} \times {\mb R} \arrow[d, "{(t_1, \ldots, t_{\p}, s) \mapsto (t_1 + \cdots + t_{\p}, s)}"] \\
\Pi \times {\mb Z} \arrow[rr] & & {\mb R} \times {\mb R}
\end{tikzcd}
\eeq
Moreover, as a weak bimodule over $({\mb M}\times \cdots \times {\mb M} \times \ep; \mb{M}^\zp)$, it is $\zp$-equivariant. 
\end{prop}

\begin{proof}
    This follows from the usual description of compactifications of moduli spaces of broken Morse trajectories and our choice of Novikov group actions. Note that the sphere bubbles in the Gromov compactification do not introduce new strata.
\end{proof}

Now we state the technical theorem leading to the definition of the quantum Steenrod operation. Notice that the involved flow categories are already enriched in the category of stratified smooth manifolds, in particular, the category $\uds{\bf dOrb}_{\rm rig}^{\rm NC}$. Fix an outercollaring width.

\begin{thm}[The AMS construction for quantum Steenrod operation]\label{thm_QST_AMS}
Fix an outercollaring width.
\begin{enumerate}

\item For any integral symplectic form $\Omega$ tamed by $J$, there exists a class of lifts (called AMS lifts) of the outercollaring $\outer {\mb X}^{\mb{QSt}}$ to $\outer \uds{\bf dOrb}_{\rm rig}^{\rm NC}$, as a multimodule over $(\outer {\mb M}, \cdots, \outer{\mb M}, \outer \ep; \outer {\mb M}^\ep)$, which is equivariant with respect to the morphism \eqref{eqn201} and which is $\zp$-equivariant as a weak bimodule. 

\item Any two AMS lifts of $\mb{X}^{\mb{QSt}}$ are equivariantly homotopic. Namely, if $\tilde {\mb X}_0^{\mb{QSt}}$ and $\tilde {\mb X}_1^{\mb{QSt}}$ are both AMS lifts of $\outer \mb{X}^{\mb{QSt}}$, then there exists a $\mb\Pi_\p$-equivariant lift of the trivial homotopy from $\outer {\mb X}^{\mb{QSt}}$ to $\outer \uds{\bf dOrb}_{\rm rig}^{\rm NC}$.
\end{enumerate}
\end{thm}

The proof of this theorem is postponed to the end of this section. We first derive the consequence, i.e., the definition of the quantum Steenrod operation. 

\begin{rem}
One can also prove the corresponding result for the comparison of using different $J$'s. We omit the detailed statement.
\end{rem}

\subsection{Steenrod operations}

As the flow categories involved are already transverse, they themselves can be regarded as their FOP transverse perturbations. Then by Theorem \ref{thma_FOP}, one can find a $(\Pi \times {\mb Z})\rtimes \zp$-equivariant FOP transverse perturbation of $\tilde {\mb X}^{\mb{QSt}}$, whose isotropy free locus $\mathring{\mb X}^{\mb{QSt}}$ is a multimodule enriched in $\pman$. Notice that equivariant feature of Theorem \ref{thma_FOP} guarantees that all existing symmetries can be respected by the FOP transverse perturbation. 

Then by Proposition \ref{prop_chain_map} and Proposition \ref{prop_Novikov_linear_chain_map}, one obtains a chain map
\beqn
\wt\Phi^{\mb{QSt}}: = \left( \left( \overbrace{ C ({\mb M}) \widehat{ \underset{\fp}{\otimes} }\cdots \widehat{ \underset{\fp}{\otimes} } C ({\mb M})}^{\p} \right)   \widehat{ \underset{\fp}{\otimes} } C (\ep)  \right) \underset{\Lambda^{\Pi^{\p} \times {\mb Z}}}{\otimes} \Lambda^{\mb\Pi} \to  C ({\mb M} \# \ep).
\eeqn
The restriction to the $\zp$-invariant part is denoted by $\Phi^{\mb{QSt}}$.

\begin{prop}
Up to chain homotopy, the map $\wt\Phi^{\mb{QSt}}$ does not depend on choices.
\end{prop}

\begin{proof}
This follows from Part (2) of Theorem \ref{thm_QST_AMS} and the homotopy case of Theorem \ref{thma_FOP}.
\end{proof}

It reduces to a linear map between homology, which, in $\fp$-coefficients, reads
\beqn
\Phi^{\mb{QSt}}: H_{\zp} ( C (\mb{M})^{\otimes \p})  \to H ( \mb{M}^{\zp}).
\eeqn
One can see that 
\beqn
H_{\zp} ( C^*(\mb{M})^{\otimes \p}) \cong H_{\zp} ( CM(f, g)^{\otimes \p}) \otimes \Lambda_{\fp}^\Pi.
\eeqn

Then by pre-composing with the quasi-Frobenius map \eqref{quasi_Frobenius} one obtains a map in homology
\beqn
\xymatrix{ HM (f, g) \ar[r]^-{qF} & H_{\zp} ( CM (f, g)^{\otimes \p}) \ar[r] & H_{\zp} ( C ({\mb M})^{\otimes \p}) \ar[r]^-{\Phi^{\mb{QSt}}} & H ( C({\mb M} \# {\mb E}))}.
\eeqn
Notice that the quasi-Frobenius map is compatible with continuation maps. To see that the composition is independent of the choice of the Morse--Smale pair $(f, g)$, we can concatenate the multimodule $\mb{X}^{\mb{QSt}}$ with Morse bimodules interpolating different choices of the Morse--Smale pairs, therefore we obtain the following statement using standard ideas.

\begin{lemma}
Let ${\mb M}'$ be another Morse flow category on $X$ associated to a different Morse--Smale pair $(f', g')$. Then the following diagram commutes.
\beqn
\xymatrix{ HM^*(f, g) \ar[r]  \ar[d] &    H_{\zp}^*(C^*(\mb{M})^{\otimes \p}) \ar[rr] \ar[d] & &  H^*(C^*(\mb{M}\# {\mb E})) \ar[d] \\
HM^*(f', g') \ar[r]  & H_{\zp}^*(C^*({\mb M}')^{\otimes \p}) \ar[rr]  & &  H^*( C^*(\mb{M}' \# {\mb E}))}
\eeqn
Here vertical arrows are isomorphisms induced from continuation maps.
\end{lemma}

Notice that 
\beqn
H^*(C^*({\mb M}\# {\mb E})) \cong H^*(X; \Lambda_{\fp}^{\mb\Pi} \langle \uptheta \rangle ).
 \eeqn
Then the above lemma shows that, after tensoring with $\Lambda_{\tkp}^\Pi$, one obtains a canonical map
\beqn
QSt_{\p}: H^*(X; \fp) \to H^*(X; \Lambda_{\tkp}^\Pi).
\eeqn
By construction, the chain map $\wt\Phi^{\mb{QSt}}$ does not decrease the $\tu$-grading, hence it indeed defines a map for each degree $k$:
\beqn
QSt_{\p}: H^k( X; {\mb F}_2) \to \left\{ \begin{array}{ll} H^{2k}(X; \Lambda_{{\bf F}_2\pst }^\Pi),\ &\ \p = 2,\\
H^{\p k}(X; \Lambda_{\fp \psu}^{\Pi} \langle \uptheta \rangle),\ &\ \p > 2.  \end{array}\right.
\eeqn
This is called the {\bf quantum Steenrod operation}. The proof of the claim that the classical part of $QSt_{\p}$ (which comes from constant holomorphic spheres) coincides with the classical Steenrod operation is the same as the semi-positive case considered by Seidel--Wilkins \cite{Seidel_Wilkins_2022}. This finishes the proof of part (1) of Theorem \ref{thma_QST}.

\subsection{AMS construction}

We would like to construct $\zp$-equivariant AMS lifts which is also equivariant with respect to the monoid morphism \eqref{eqn201}. It suffices to construct the lift for the corresponding moduli space of $\p+1$-marked curves, from which the flow multimodule is constructed by taking fiber products with moduli spaces of Morse trajectories.

\subsubsection{The domain moduli}

We first consider the moduli space of domain curves. For each $d\geq 0$, consider
\beqn
\ov{\mc M}_{\Sigma^{\mb{QSt}}}(\mb{CP}^d, d)
\eeqn
of holomorphic maps from the special marked curve $\Sigma_\p^{\mb{St}}$ to $\mb{CP}^d$ of degree $d$. Let
\beqn
B_d^{\mb{QSt}} \subset \ov{\mc M}{}_{\Sigma^{\mb{QSt}}}(\mb{CP}^d, d)
\eeqn
be the open subset whose images are not contained in any hyperplane. Then $B_d^{\mb{QSt}}$ is a ${\mc G}_d^{\mb{QSt}}:= PGL(d+1)$-invariant open subset. Moreover, there is an action by $\zp$ which permutes the incoming marked points. Then one has the object
\beqn
{\mc C}_d^{\mb{QSt}} = (\zp \times {\mc G}_d^{\mb{QSt}}, B_d^{\mb{QSt}}, C_d^{\mb{QSt}})\in {\rm Ob} \uds{\bf Curve}.
\eeqn 
Since this is not stratified, we do not need to consider outercollaring or rigidifications at this moment. We can equip ${\mc C}_d^{\mb{QSt}}$ a $\zp\times G_d^{\mb{QSt}}$-equivariant stable complex structure, hence an object of $\uds{\bf Curve}^{\mb{C}}$, using the complex analytic structure. 

\subsubsection{The thickenings}

Choose an integral symplectic form $\Omega$ which is tamed by $J$. Consider the space of maps and framed maps
\begin{align*}
&\ {\rm Map}_d(C_d^{\mb{QSt}}, X),\ &\ {\rm Map}_d^{\rm fr}(C_d^{\mb{QSt}}, X)
\end{align*}
whose elements are $(\phi, u)$ resp. $(\phi, u, F)$ where $\phi \in B_d^{\mb{QSt}}$, $u: C_\phi \to X$ is a smooth map satisfying $\int u^* \Omega = d$, and $F$ is a framing associated to $u$ and the integral symplectic form $\Omega$. 
One can choose a $\zp$-equivariant group reduction scheme
\beqn
\lambda_d^{\mb{QSt}}: {\rm Map}_d^{\rm fr}(C_d^{\mb{QSt}}, X) \to Q_d^{\mb{QSt}}:= T_1 ( {\mc G}_d^{\mb{QSt}}/ G_d^{\mb{QSt}}),
\eeqn
for example, using the formula \eqref{frame_pairing}. On the other hand, choose a thickening datum 
\beqn
\iota_d^{\mb{QSt}}: 
W_d^{\mb{QSt}} \to \Gamma_c( \mathring C_d^{\mb{QSt}}\times X, \Omega^{0,1}_{\mathring C_d^{\mb{QSt}} / B_d^{\mb{QSt}}} \otimes TX)
\eeqn
where $W_d^{\mb{QSt}}$ is a unitary representation of $G_d^{\mb{QSt}} \cong PU(d+1)$ such that $\iota_d^{\mb{QSt}}$ is equivariant with respect to both the $G_d$-action and the $\zp$-action. 

Then one has the thickened moduli space
\beqn
V_d^{\mb{QSt}}:= \left\{ (\phi, u, F, e)\ \left| \begin{array}{c} (\phi, u, F) \in {\rm Map}_d^{\rm fr}(C_d^{\mb{QSt}}, X),\ e \in W_d^{\mb{QSt}},\\
\ov\partial u + \iota_d^{\mb{QSt}} (e) = 0,\ \phi = \phi_F \in B_d^{\mb{QSt}}  \end{array} \right.\right\}
\eeqn
The obstruction bundle 
\beqn
E_d^{\mb{QSt}} \to V_d^{\mb{QSt}}
\eeqn
is the trivial bundle with fiber $W_d^{\mb{QSt}} \oplus Q_d^{\mb{QSt}}$ and the Kuranishi section is 
\beqn
S_d^{\mb{QSt}} (\phi, u, F, e) = (e, \lambda_d^{\mb{QSt}} (\phi, u, F)).
\eeqn
The obtained topological Kuranishi space is denoted by 
\beqn
K_d^{\mb{QSt}}:= (G_d^{\mb{QSt}}, V_d^{\mb{QSt}}, E_d^{\mb{QSt}}, S_d^{\mb{QSt}})
\eeqn
which also has a $\zp$-action commuting with $G_d^{\mb{QSt}}$.
Notice that there is a $\zp$-equivariant evaluation map
\beqn
\ev: V_d^{\mb{QSt}} \to X^{\p + 1}
\eeqn
at the markings, where $\zp$ acts on the first $\p$ copies of $X$ by permutation.

One still have the relative smoothness of the forgetful map 
\beqn
V_d^{\mb{QSt}} \to B_d^{\mb{QSt}}.
\eeqn
One can then also follow the approach of Section \ref{section_NC_structure} (which is much simplified here) to equip the vertical tangent bundle an NC structure. See \cite[Section 6.2]{Bai_Xu_2022} for related setups in Gromov--Witten theory.

Since we are working with a single space without boundary or corners, one can apply the smoothing strategy of \cite{AMS} to obtain a stable $G_d \times \zp$-smoothing of $V_d^{\mb{QSt}}$; the stabilization factor needed for the smoothing can be absorbed into the thickening space $W_d^{\mb{QSt}}$. Moreover, one can make the evaluation map $\ev$ a smooth submersion at the expense of introducing further stabilization. 

Then eventually one can use the stable complex structure of $B_d^{\mb{QSt}}$ and the NC structure on vertical tangent bundles to construct a $\zp\times G_d^{\mb{QSt}}$-equivariant NC structure on $V_d^{\mb{QSt}}$. Then one obtains an object of $\uds{\bf SKur}^{\rm NC}$ written as
\beqn
K_d^{\mb{QSt}} = (\zp\times G_d^{\mb{QSt}}, V_d^{\mb{QSt}}, E_d^{\mb{QSt}}, S_d^{\mb{QSt}})
\eeqn
with an equivariant smooth evaluation map
\beqn
\ev: V_d^{\mb{QSt}} \to X^{\p + 1}.
\eeqn

\subsubsection{The Kuranishi and derived orbifold multimodule}

Now we will insert Morse constraints and couple with the diagonal bimodule of $\ep$. For critical points $\uds x_1, \ldots, \uds x_\p, \uds x'$ of $f^{\mb M}$, $w, w'$ of $f^{\ep}$, consider 
\beqn
\ev^{-1}( W_{\uds x_1}^{\rm us} \times \cdots \times W_{\uds x_\p}^{\rm us} \times W_{x'}^{\rm st}) \times \mathring M_{w;w'}^{{\mb\Delta}^\ep}
\eeqn
where $W_{\uds x_i}^{\rm us}$ resp. $W_{\uds x'}^{\rm st}$ is the uncompactified unstable resp. stable manifold and $\mathring M_{w;w'}^{\mb{\Delta}^\ep}$ is the uncompactified moduli space of parametrized trajectories in $\wh S^\infty$ connecting $w$ and $w'$. Then we compactify by allowing breaking as follows: at negative ends, the breakings at different gradient rays/lines are independent; at the positive end, the breaking is regarded as breaking in the product space $X \times S^\infty$, therefore ``synchronized." Let the compactified space be
\beqn
V_d^{\mb{QSt}} (\uds x_1 \cdots \uds x_\p w; (\uds x',w')).
\eeqn
Then one obtains a stratified smooth Kuranishi space
\beqn
K_d^{\mb{QSt}} ( \uds x_1 \cdots \uds x_\p w; (\uds x', w'))
\eeqn
whose symmetry group is reduced to $G_d^{\mb{QSt}}$.

We define a lift of the Steenrod multimodule to $\uds {\bf SKur}^{\rm NC}$, denoted by $\hat{\mb X}^{\mb{QSt}}$. Given objects $x_i = (a_i, \uds x_i)$ where $a_i \in \Pi$ and $\uds x_i \in {\rm crit} f$, $i = 1, \ldots, \p$ and $x' = (a', \uds x')$, denote $d = \Omega(a' - a_1- \cdots - a_\p)$ and 
\beqn
K_{x_1 \cdots x_{\p} w; (x', w')}^{\mb{QSt}}:= K_d^{\mb{QSt}} (\uds x_1\cdots \uds x_\p w; (x', w')) \in \uds{\bf SKur}^{\rm NC}.
\eeqn
The multimodule structural maps are obtained in a very simple way as we only need to consider broken Morse trajectories. Hence one obtains a multimodule $\hat {\mb X}^{\mb{QSt}}$ over $({\mb M}, \ldots, {\mb M}, \ep; {\mb M}\# {\mb E})$. One is readily to verify that it is a lift of ${\mb X}^{\mb{QSt}}$ and satisfies the expected equivariance property. 

Finally one takes group quotient by $G_d^{\mb{QSt}}$ and obtains a lift to $\uds{\bf dOrb}^{\rm NC}$, denoted by $\tilde{\mb X}^{\mb{QSt}}$. One can do the outercollaring as the last step, providing the desired AMS construction claimed in Theorem \ref{thm_QST_AMS}. 

Lastly, if we make a different set of choices, then the trivial homotopy from ${\mb X}^{\mb{QSt}}$ to itself can be lifted to $\outer \uds{\bf dOrb}_{\rm rig}^{\rm NC}$ in the equivariant way, between the AMS lifts of two different choices. This finishes the proof of Theorem \ref{thm_QST_AMS}.

\subsection{Comparing to Seidel--Wilkins' definition in semi-positive case}

When $(X, \omega)$ is semi-positive, Seidel--Wilkins \cite{Seidel_Wilkins_2022} defined their quantum Steenrod operations using the classical approach. We prove that our definition agrees with theirs.

\begin{thm}\label{thm_steenrod_comparison}
When $(X, \omega)$ is semi-positive, the quantum Steenrod operation defined above agrees with the definition of Seidel--Wilkins.
\end{thm}

\begin{rem}
    Strictly speaking, we are referring to the operations denoted by $QSt$ instead of $Q\Sigma_b$ in \cite{Seidel_Wilkins_2022}, where the latter is featured more in recent studies of arithmetic properties of quantum connections on semi-positive symplectic manifolds (cf. \cite{chen2024quantumsteenrod, chen2024exponential}, \cite{jae, Lee23b}, \cite{bai-lee, bai2026}). We can extend our $QSt_\p$ to define analogs of $Q\Sigma_b$ using our framework and prove that it reduces to the Seidel--Wilkins operation when $(X, \omega)$ is semi-positive. Moreover, such an extension is covariantly constant under the quantum connection defined using FOP perturbations. We will not discuss such extensions in this manuscript.
\end{rem}

\subsubsection{The construction of Seidel--Wilkins}

The quantum Steenrod operation defined in \cite{Seidel_Wilkins_2022} is based on the classical transversality via adding $\zp$-equivariant inhomogeneous terms parametrized by $S^\infty$. We extend the perturbation to $\wh S^\infty$. Consider inhomogeneous terms 
\beqn
\nu_{S^2, \eta, z, x}^{\rm eq}: T_z S^2 \to T_x X\ {\rm for\ } (\eta, z, x) \in \wh S^\infty \times S^2 \times X
\eeqn
which satisfies
\beqn
\nu_{S^2, \sigma (\eta), z, x}^{\rm eq} = \nu_{S^2, \eta, \sigma(z), x}^{\rm eq} \circ D \sigma: T_z S^2 \to T_x M\ \ \forall \sigma \in \zp.
\eeqn
Then one can consider the moduli space of maps $u: S^2 \to X$ satisfying the perturbed Cauchy--Riemann equation. More precisely, consider pairs $(\eta, u)$ where $\eta \in \wh S^\infty$ and $u: S^2 \to X$ such that
\beqn
\ov\partial_J u(z) + \nu_{S^2, \eta, z, u(z)}^{\rm eq} = 0.
\eeqn

We modify Seidel--Wilkins' formulation slightly using the Morse model. For a comprehensive comparison with the original Seidel--Wilkins' setup, which is based on finding explicit cycles representing generators if $H^*(B\zp; {\mb F}_p)$, we refer the reader to \cite[Appendix A]{BSWX}. Recall that one has identified critical points of $f^\ep$. For each pair of critical points $w, w'\in {\rm crit} f^\ep$, denote
\beqn
\mathring \Delta_{w; w'}:= \Big\{ \eta: {\mb R} \to \wh S^\infty\ |\ \eta'(s) = \nabla f^\ep (\eta (s)) = 0,\ \lim_{s \to -\infty}  \eta(s) = w,\ \lim_{s \to +\infty} \eta (s) = w' \Big\}.
\eeqn
Then elements of $\mathring \Delta_{w;w'}$ is determined by its evaluation at $s = 0$. On the other hand, one has
\beqn
\wh S^\infty = \bigsqcup_{w, w'} \mathring \Delta_{w; w'}.
\eeqn

Then given critical points $\uds x_1, \ldots, \uds x_{\p}, \uds x'$ of $f: X \to {\mb R}$ and $a \in \Pi$, consider the moduli space
\beqn
\mathring {\mc M}_{\uds x_1 \cdots \uds x_{\p}; \uds x'}^{\rm SW}(a; \mathring \Delta_{w; w'})
\eeqn
which consists of pairs $(u, \eta)$ where $u: S^2 \to X$ is a smooth map and $\eta \in \mathring \Delta_{w; w'}$ satisfying
\beqn
\ov\partial_J u(z) + \nu_{S^2, \eta, z, u(z)}^{\rm eq} = 0
\eeqn
and 
\beqn
u(z_i) \in W^{\rm us}_{\uds x_i},\ u(z_\infty) \in W^{\rm st}_{\uds x'}.
\eeqn
Seidel--Wilkins proved that for a generic inhomogeneous term $\nu_{S^2}^{\rm eq}$, whenever the expected dimension of the above moduli space is zero or one, the defining perturbed $J$-holomorphic map equation is transverse. 

\begin{lemma}\label{lemma:seidel-wilkins}(cf. \cite[Lemma 4.1]{Seidel_Wilkins_2022}) For a generic $\omega$-compatible almost complex structure $J$ and a ${\mb Z}\times \zp$-equivariant inhomogeneous term $\nu_{S^2, \eta, z, x}^{\rm eq}$, any moduli space $\mathring {\mc M}_{\uds x_1 \cdots \uds x_\p; \uds x'}^{\rm SW}(a; \mathring \Delta_{w; w'})$ is transverse provided that its expected dimension is at most $1$. Moreover,
\begin{enumerate}
    \item When the expected dimension is zero, $\mathring {\mc M}_{\uds x_1 \cdots \uds x_\p; \uds x'}^{\rm SW}(a; \mathring \Delta_{w; w'})$ consists of finitely many points.

    \item When the expected dimension is one, $\mathring {\mc M}_{\uds x_1\cdots \uds x_\p; \uds x'}^{\rm SW}(a; \mathring \Delta_{w; w'})$ admits a natural compactification which is a compact 1-dimensional manifold with boundary. 
\end{enumerate}
\end{lemma}

\begin{proof}
The only difference between our case and the case of \cite{Seidel_Wilkins_2022} is that there is no requirement for ${\mb Z}$-symmetry in \cite{Seidel_Wilkins_2022}. However, this can be achieved as the ${\mb Z}$-action is free. 
\end{proof}

Then the chain-level quantum Steenrod operation can be defined using the count of zero-dimensional moduli spaces. One can use the standard argument to show that up to chain homotopy the chain-level map is well-defined. 

\subsubsection{The AMS lifts and FOP perturbations with an inhomogeneous term}

Now we define a perturbation of the quantum Steenrod multimodule ${\mb X}^{\mb{QSt}}$ by using the inhomogeneous term. Notice that each $\eta \in \wh S^\infty$ corresponds uniquely to a smooth parametrized Morse flow line $y_w(s)$ in $S^\infty$. Then for objects $x_1, \ldots, x_\p, w; (x', w')$, consider the moduli space
\beqn
\mathring M_{x_1 \cdots x_\p w; (x', w')}^{\mb{X}_1^{\mb{QSt}}} = \left\{ (u, \eta)\ \left|\ \begin{array}{c}  u: \Sigma^{\mb{QSt}} \to X,\ \eta \in W^{\rm us}_{w} \cap W^{\rm st}_{w'} \subset \wh S^\infty,\\
\ov\partial_J u + \nu = 0,\\
u(z_i) \in W^{\rm us}_{\uds x_i},\ u(z_\infty) \in W^{\rm st}_{\uds x'}
\end{array}  \right. \right\}.
\eeqn
It can be compactified in a natural way, which is denoted by 
\beqn
M_{x_1 \cdots x_\p w; (x', w')}^{\mb{X}_1^{\mb{QSt}}}.
\eeqn
These moduli spaces form a multimodule ${\mb X}_1^{\mb{QSt}}$ over $({\mb M}, \ldots, \mb{M}, \ep; {\mb M} \# \ep)$, which is a variant of the $\nu = 0$ case considered previously, which is denoted by $\mb{X}_0^{\mb{QSt}}$. After outercollaring, using the same method as proving Theorem \ref{thm_QST_AMS}, one can construct an AMS lift $\hat{\mb X}_1^{\mb{QSt}}$ enriched in $\outer \uds{\bf SKur}_{\rm rig}^{\rm NC}$. The inhomogeneous term $\nu$ can be dealt with by considering the graph $\mathrm{graph}(u): \Sigma^{\mb{QSt}} \to \Sigma^{\mb{QSt}} \times X$, where the product $\Sigma^{\mb{QSt}} \times X$ is equipped with the sheared almost complex structure (e.g. \cite[Lemma 6.18]{AMS}) so that we can go back to the standard situation. The resulting lift in $\outer \uds{\bf dOrb}^{\rm NC}_{\rm rig}$ is denoted by $\tilde{\mb X}_1^{\mb{QSt}}$.

\begin{lemma}
When the virtual dimension is zero, the original section is FOP transverse.
\end{lemma}

\begin{proof}
By the monotonicity condition and the classical transversality condition, when the virtual dimension is zero, the moduli space $\mathring M_{x_1 \cdots x_\p w; (x', w')}^{\mb{X}_1^{\mb{QSt}}}$ agrees with its compactification and there are no orbifold points in the moduli space. The Kuranishi section on the thickening is transverse to the zero section because the Cauchy--Riemann equation is transversely cut out. Then by the {\bf (Classical Transversality)} property of FOP transverse perturbations, the original section is FOP transverse. 
\end{proof}

\begin{prop}
There exists an FOP transverse perturbation of $\tilde {\mb X}_1^{\mb{QSt}}$ such that for moduli spaces of dimension at most $0$, the perturbation agrees with the original Kuranishi section. 
\end{prop}

\begin{proof}
This follows from the same proof of Lemma \ref{lemma_semipositive_FOP}.
\end{proof}

Furthermore, there is a clear homotopy $\mb{H}^{\mb{QSt}}$ from $\mb{X}_0^{\mb{QSt}}$ to $\mb{X}_1^{\mb{QSt}}$ by adding an additional parameter $t \in [0, 1]$ and changing $\nu$ to $t \nu$ in the equation. One can further construct an AMS lift $\hat {\mb H}^{\mb{QSt}}$ as a homotopy from $\hat{\mb X}_0^{\mb{QSt}}$ to $\hat{\mb X}_1^{\mb{QSt}}$. Moreover, one can extend the FOP transverse perturbations $\tilde {\mb X}_0^{\mb{QSt}}$ and $\tilde {\mb X}_1^{\mb{QSt}}$ to an FOP transverse perturbation on $\hat {\mb H}^{\mb{QSt}}$, denoted by $\tilde {\mb H}^{\mb{QSt}}$, in an equivariant way. This leads to a chain homotopy between the chain maps defined with or without the inhomogeneous term. This concludes the proof of Theorem \ref{thm_steenrod_comparison}.

\section{The Equivariant Floer theory for iterations---Proof of Theorem \ref{thma_equivariant_Floer}}

In this section we construct the equivariant Floer theory, including equivariant Floer cohomology, equivariant continuation maps, and equivariant PSS/SSP maps, induced from a prime iteration of a Hamiltonian. Let $\p$ be a prime. 

\subsection{Equivariant Floer flow category as a coupling}

\subsubsection{Coupling of Floer and Morse flow categories}

We define a notion called {\bf coupling} of flow categories. This notion cannot be defined on the abstract level but rely on the geometric setup. 

Let ${\mb F}$ be the Floer flow category on $X$ associated to a pair $(H, J)$. Let ${\mb M}$ be the Morse flow category on another manifold $Y$ associated to a Morse--Smale pair $(f, g)$. Define ${\mb F} \# {\mb M}$ be a flow category as follows. The set of objects is 
\beqn
{\rm Ob} ( {\mb F} \# {\mb M} ) = {\rm Ob} {\mb F} \times {\rm Ob}{\mb M}
\eeqn
whose elements are denoted by $p = ( p^{\mb F}, p^{\mb M})$. For the morphism space $M_{pq}^{{\mb F}\# {\mb M}}$, consider the moduli space of pairs $(u, v)$ where $u: {\mb R}\times S^1 \to X$ is a solution to the Floer equation connection $p^{\mb F}$ and $q^{\mb F}$ and $v: {\mb R} \to Y$ is a Morse trajectory connection $p^{\mb M}$ and $q^{\mb M}$; two pairs $(u, v)$ and $(u', v')$ are equivalent if there exists $s_0 \in {\mb R}$ such that 
\begin{align*}
    &\ u(\cdot, \cdot) = u'(\cdot + s_0, \cdot),\ &\ v(\cdot) = v'(\cdot + s_0).
\end{align*}
One then compactifies this moduli spaces by allowing ``synchronized'' breakings (as well as sphere bubblings).

Now we describe the equivariant Floer flow category as a special coupling. 

\begin{defn}\label{defn_equi_flow_category}
Let $H$ be a 1-periodic Hamiltonian on $(X, \omega)$ such that the $\p$-th iteration $H^{\natural \p}$ is nondegenerate. Let $J$ be a fixed $\omega$-compatible almost complex structure. Define the {\bf $\zp$-equivariant Floer flow category} to be the coupling
\beqn
{\mb F}^{\zp}:= {\mb F} \# {\mb E}
\eeqn
where ${\mb E}$ is the Morse flow category on $\wh S^\infty$.
\end{defn}

On ${\mb F}^\zp$ there is an action
\beqn
{\mc A} = ({\mc A}_1, {\mc A}_2): {\rm Ob}{\mb F}^\zp \to {\mb R}^2
\eeqn
where
\begin{align*}
&\ {\mc A}_1(p^{\mb F}, p^\ep) = {\mc A}_{H^{\natural p}}(p^{\mb F}),\ &\ {\mc A}_2(p^{\mb F}, p^\ep) = f^\ep( p^\ep).
\end{align*}
The grading is given by 
\beqn
{\rm deg} (p^{\mb F}, p^\ep):= {\rm deg}(p^{\mb F}) + {\rm deg}(p^\ep)\in {\mb Z}.
\eeqn
Moreover, there is a free action on ${\rm Ob} {\mb F}^\zp$ the twisted Novikov group
\beqn
{\mb\Pi}_\p:= \Pi \times {\mb Z} \times \zp = {\mb \Pi} \times {\mb Z}_\p.
\eeqn

For different Hamiltonians one can also define an equivariant continuation bimodule. Let $H_i^{\natural \p}$ be pairs as above for $i = 0, 1$. Choose a Floer data on the infinite cylinder connecting $H_0$ and $H_1$, denoted by $\sigma$. Let $\sigma^{\natural \p}$ be the $\p$-th iteration as a Floer daum on the infinite cylinder connecting $H_0^{\natural \p}$ and $H_1^{\natural \p}$. By coupling the bimodule induced from $\sigma^{\natural \p}$ and the diagonal bimodule ${\mb \Delta}^{\ep}$, one obtains a $\mb\Pi_\p$-equivariant  Novikov bimodule over $({\mb F}_0^{\zp}; {\mb F}_1^{\zp})$, denoted by ${\mb B}^{\zp}_{01}$. 

\begin{prop}
$\mb{F}^\zp$ is a $\mb\Pi_\p$-equivariant Novikov flow category enriched in $\uds{\bf Top}$ and $\mb{B}_{01}^\zp$ is a $\mb\Pi_\p$-equivariant Novikov flow bimodule over $(\mb{F}_0^\zp; \mb{F}_1^\zp)$.
\end{prop}

Now we state the technical theorem regarding equivariant AMS lifts of ${\mb F}^\zp$ as well as $\mb{B}_{01}^{\zp}$.

\begin{thm}[Equivariant AMS lifts]\label{thm_equivariant_lift} \hfill

\begin{enumerate}

\item Let ${\mb F}^{\zp}$ be the $\zp$-equivariant Floer flow category defined in Definition \ref{defn_equi_flow_category}. Fix the outercollaring width and an integral action on ${\mb F}$. There is a $\mb\Pi_\p$-equivariant lift $\tilde{\mb F}^{\zp}$ of the outercollaring of ${\mb F}^{\zp}$ to the category $\outer \uds{\bf dOrb}^{\rm NC}_{\rm rig}$, called the {\bf equivariant AMS lift} of $\mb{F}^\zp$ (subject to the outercollaring width and the integral action).

\item Choose integral actions on ${\mb F}_0$ and ${\mb F}_1$ (not necessarily compatible with respect to the interpolating Floer data). Given equivariant AMS lifts of $\tilde {\mb F}_0^{\zp}$ resp. $\tilde {\mb F}_1^{\zp}$ of ${\mb F}_0^{\zp}$ resp. ${\mb F}_1^\zp$ subject to the outercollaring width and the chosen integral actions, there exists a $\mb\Pi_\p$-equivariant lift of the outercollaring of $\mb{B}_{01}^{\zp}$ to $\outer \uds{\bf dOrb}^{\rm NC}_{\rm rig}$ as a $\mb\Pi_\p$-equivariant bimodule over $(\tilde {\mb F}_0^{\zp}; \tilde {\mb F}_1^{\zp})$.
\end{enumerate}
\end{thm}

\begin{proof}
See next Subsection.
\end{proof}

Now assuming this theorem, one can define the equivariant Floer chain complex. Since the $\mb\Pi_\p$-action is free, by Theorem \ref{thma_FOP}, one can find a $\mb\Pi_\p$-equivariant FOP perturbation on $\tilde{\mb F}^\zp$, denoted by $\mathring {\mb F}^{\zp}$. Therefore, by taking the isotropy free part of the perturbed zero locus, one obtains a $\mb\Pi_\p$-equivariant ${\mb Z}$-graded Novikov flow category enriched in $\pman$. Notice that there is still the $\zp$-action which preserves all structures. Then applying the abstract construction in Section \ref{section6}, one obtains a ${\mb Z}$-graded cochain complex 
\beqn
CF^*(H^{\natural \p}, J^{\natural \p}, \Xi^{\natural \p}):= C^*(\mathring {\mb F}^\zp)
\eeqn 
which is a graded $\Lambda_{\mb Z}^{\mb\Pi}$-module with a $\zp$-action. The $\zp$-fixed part, denoted by 
\beqn
CF_{\zp}(H^{\natural \p}),
\eeqn
is called the (Tate) {\bf equivariant Floer complex} of $H^{\natural \p}$.

Furthermore, suppose one has two such equivariant complexes 
\begin{align*}
&\ CF^*(H_1^{\natural \p}),\ &\ CF^*(H_2^{\natural \p}),
\end{align*}
one can extend, by Theorem \ref{thma_FOP}, the $\mb\Pi_\p$-equivariant FOP perturbation over $\tilde {\mb B}_{01}^{\zp}$, which induces a $\zp$-equivariant chain map
\beqn
\Phi_{01}^{\zp}: CF^*(H_1^{\natural \p}) \to CF^*(H_2^{\natural \p}).
\eeqn
Then, upon taking $\zp$-fixed part, one has the continuation maps
\beqn
\Phi_{01}: CF_{\zp}(H_1^{\natural \p}) \to CF_{\zp}(H_2^{\natural \p}).
\eeqn

One can proceed in the same way as proving Theorem \ref{thma_Floer_complex} and Theorem \ref{thma_continuation_map}. Therefore, one obtains a well-defined chain homotopy class of cochain complexes of $\Lambda_{\mb Z}^{\mb\Pi}$-modules, hence a well-defined symplectic invariant
\beqn
HF_{\zp}(X)
\eeqn
which is a module over $\Lambda_{\fp}^{\mb\Pi}$. 

\subsection{AMS construction for coupling of Floer and Morse flow categories}\label{subsection_AMS_equivariant_Floer}

In this subsection we carry out the AMS construction for the coupled flow category. 

\subsubsection{The domain flow category}

We modify the previous construction of AMS lifts for non-equivariant Floer theory. We do not consider the monotone flow category of domains as we have to couple it with the Morse flow category. Instead, one describes a flow category $\mb{Dom}^{\zp}$ enriched in $\outer \uds{\bf Curve}_{\rm rig}^{\mb C}$ with objects the same as those of ${\mb F}^{\zp}$.

First we choose an integral action ${\mc A}_{{\mb F}}^\Omega: {\rm Ob}{\mb F} \to {\mb Z}$ for the Floer flow category associated to $H^{\natural \p}$. Notice that from the construction in the proof of Proposition \ref{prop163}, the integral action can be chosen to be $\zp$-invariant. Then for each pair of objects $p = (p^{\mb F}, p^{\mb E})$, $q = (q^{\mb F}, q^{\mb E})$, denote
\beqn
d_{pq}:= {\mc A}^\Omega(q^{{\mb F}}) - {\mc A}^\Omega(p^{\mb F}) \in {\mb Z}.
\eeqn
We describe an object
\beqn
{\mc C}_{pq}  = ({\mc G}_{pq}, B_{pq}, C_{pq})\in {\rm Ob}\uds{\bf Curve}
\eeqn
as follows. Define ${\mc G}_{xy} = {\mc G}_{d_{pq}}$ as considered before. Consider the  moduli space of pairs
\beqn
\phi: {\mb R} \times S^1 \to \mb{CP}^{d_{pq}},\ v: {\mb R} \to \wh S^\infty
\eeqn
where $\phi$ is holomorphic of degree $d$ and $v$ is a flow line connection $p^{\mb E}$ and $q^{\mb E}$. Two pairs $(u, v)$, $(u', v')$ are equivalent if they differ by a common translation in ${\mb R}$. Then we compactify this space in the usual way, requiring the positive constraint
\beqn
\phi (z_+) = [0,\ldots, 0, 1] \in \mb{CP}^{d_{pq}}
\eeqn
and requiring the image of $\phi$ is not contained in any hyperplane. The space we obtain is denoted by $B_{pq}$, whose dimension is 
\beqn
{\rm dim}_{\mb R} B_{pq} = 2d_{pq} (d_{pq} +1) + {\rm index} ( q^\ep ) - {\rm index} ( p^\ep ) - 1.
\eeqn
Let $C_{pq} \to B_{pq}$ be the universal curve whose fiber at $\phi \in B_{pq}$ is the corresponding rational curve in $\mb{CP}^d$. This is a ${\mc G}_{pq}$-equivariant family, hence one obtains an object ${\mc C}_{xy}$ in $\uds{\bf Curve}$.

Notice that the space $B_{pq}$ is stratified by the following poset
\beqn
A_{pq} = \left\{ \Big( (d_0, \ldots, d_l), p^\ep r_1^\ep \cdots r_l^\ep q^\ep \Big)\ \left|  \begin{array}{c} \ d_i \geq 0,\ \displaystyle \sum_i d_i = d,\\
p^\ep = r_0^\ep \leq r_1^\ep \leq \cdots \leq r_l^\ep \leq r_{l+1}^\ep =  q^\ep,\\
d_i = 0\Longrightarrow  r_i^\ep \neq r_{i+1}^\ep\      \end{array} \right. \right\}.
\eeqn

Notice that there is a free $\zp$-action on the object ${\mc C}_{pq}$, given by the domain rotation and the action on $\wh S^\infty$: for $\gamma \in \zp$ and $\phi \in B_{xy}$ represented by $u: \Sigma \to \mb{CP}^d$ and $v: {\mb R} \to \wh S^\infty$, define 
\beqn
\gamma \cdot \phi = [u \circ \gamma^{-1}, \gamma\cdot v]
\eeqn
where $u\circ \gamma^{-1}$ is the composition with the domain rotation and $\gamma \cdot v$ is induced by the $\zp$-action on $\wh S^\infty$. This action commutes with the action of ${\mc G}_{pq}$ for the following reason: the ${\mc G}_{pq}$-action comes from the geometry of the target space $\mb{CP}^d$ while $\zp$ acts by reparametrizing the maps. On the other hand, there is also an action by ${\mb Z}$ by the corresponding action on $\wh S^\infty$. 

Using the construction of Section \ref{section_AMS_domains}, one can define the following morphisms
\beqn
\zeta_{prq}: {\mc C}_{pr} \times {\mc C}_{rq} \to {\mc C}_{pq}
\eeqn
in the category $\uds{\bf Curve}$. From the explicit construction of the domain curves, this morphism is $\mb\Pi_\p$-equivariant. Hence one obtains a flow category enriched in $\uds{\bf Curve}$, denoted by $\mb{Dom}^{\mb{F}^\zp}$.

\begin{prop}\label{prop245}
Any outercollaring of $\mb{Dom}^{\mb{F}^{\zp}}$ admits a $\mb\Pi_\p$-equivariant lift to $\outer \uds{\bf Curve}_{\rm rig}^{\mb C}$.
\end{prop}

\begin{proof}
The rigidification can be constructed in the same way as in Section \ref{section_AMS_domains}. The freeness of the $\mb\Pi_\p$-action allows us to achieve the equivariance easily. 
\end{proof}

\subsubsection{Thickenings and group reductions}

We follow the construction of Section \ref{section_Floer_lift}. We can write down the flow category 
\beqn
\outer \mb{Map}_{\mb{F}^{\zp}}^{\mb{fr}}
\eeqn
which gives each pair of objects $p, q$ an object 
\beqn
\outer {\rm Map}_{pq}^{\rm fr}
\eeqn
of framed maps. This flow category is canonically induced from previous steps and still carries the same symmetry. 

\begin{lemma}
There exist $\mb\Pi_\p$-equivariant group reduction schemes.
\end{lemma}

\begin{proof}
This can be constructed in the same way as in Proposition \ref{prop_group_reduction}. The freeness of the action allows us to maintain the symmetry. 
\end{proof}

Hence one obtains a compatible system of maps
\beqn
\lambda_{pq}: {\rm Map}_{pq}^{\rm fr} \to Q_{pq}.
\eeqn

Next we choose thickening data. 

\begin{lemma}
There exists a $\mb\Pi_\p$-equivariant and transverse thickening datum on $\outer \mb{Dom}^{\mb{F}^{\zp}}$.
\end{lemma}

\begin{proof}
The same argument as in Proposition \ref{prop_perturbation} applies. The freeness of the $\mb\Pi_\p$-action allows us to maintain the symmetry. 
\end{proof}

Hence, we obtain a compatible system of maps
\beqn
\nu_{pq}: W_{pq} \to \Gamma_c ( \mathring C_{pq} \times X, \Lambda_{\mathring C_{pq}/ B_{pq}}^{0,1}\otimes TX). 
\eeqn

Then similar to Definition \ref{defn_thickened_moduli}, one can write down the thickened moduli spaces 
\beqn
K_{pq} = (G_{xy}, V_{xy}, E_{xy}, S_{xy}) \in {\rm Ob}\outer \uds{\bf Kur}_{\rm rig}
\eeqn
and obtain a lift of the outercollaring of ${\mb F}^{\zp}$ to $\outer \uds{\bf Kur}_{\rm rig}$. 

\subsubsection{Stable smoothing and NC structure}

Next, the stable smoothing procedure works in a similar way. Indeed, the projection map $V_{xy} \to \outer B_{xy}$ has a {\rm rel}-$C^{\infty, 1}$-structure, allowing us to inductively perform stable smoothing. We omit the details. After taking quotient, one obtains a lift to $\outer \uds{\bf dOrb}_{\rm rig}$. The construction of stable complex structure is also the same. Hence one obtains
\beqn
\tilde{\mb F}^\zp
\eeqn
enriched in $\outer \uds{\bf dOrb}^{\rm NC}_{\rm rig}$. This finishes the proof of (1) of Theorem \ref{thm_equivariant_lift}. 

\subsubsection{Equivariant continuation maps}

To prove (2) of Theorem \ref{thm_equivariant_lift}, i.e., the lift of the equivariant continuation bimodule, we also modify the previous non-equivariant case. Let ${\mb F}_0^\zp$ and $\mb{F}_1^\zp$ be the two equivariant flow categories associated to the (iterations of) Hamiltonian $H_0, H_1$. Let $\mb{B}_{01}^\zp$ be the equivariant continuation bimodule associated to the $\p$-th iteration of a flow data on the cylinder interpolating between $H_0^{\natural \p}$ and $H_1^{\natural \p}$. 

We first define a collection of objects
\beqn
{\mc C}_{x_0; x_1} = ({\mc G}_{x_0; x_1}, B_{x_0; x_1}, C_{x_0; x_1}) \in {\rm Ob} \uds{\bf Curve}.
\eeqn
Using the same integral modification of the symplectic action $\Omega$, for $x_0 = (x_0^{\mb F}, x_0^\ep)$ and $x_1 = (x_1^{\mb F}, x_1^\ep)$, consider solutions $(u, v)$ where $u: {\mb R} \times S^1 \to \mb{CP}^d$, $v: {\mb R} \to S^\infty$, $u$ is holomorphic of degree $d = d_{x_0} - d_{x_1}$, and $v$ is a flow line connecting $x_0^\ep$ to $x_1^\ep$. We impose the negative constraint
\beqn
u(z_-) = [1, 0, \ldots, 0] \in \mb{CP}^d.
\eeqn
We compactify this space by allowing sphere bubbles and breakings at both positive and negative infinity. Moreover, breakings on either side are treated in the ``synchronized'' fashion. In this compactification, we also remove those solutions whose images in $\mb{CP}^d$ lie in some hyperplane. Then one obtains the moduli space $B_{x_0; x_1}$ over which there is the universal curve $C_{x_0; x_1}\to B_{x_0; x_1}$. Notice that the space $B_{x_0; x_1}$ is stratified by the poset
\beqn
A_{x_0; x_1}:= \left\{\Big( (d_{-r} \cdots d_{-1}; d_0; d_1\cdots d_s), w_{-r}^\ep \cdots w_0^\ep; w_1^\ep \cdots  w_{s+1}^\ep \Big) \left| \begin{array}{c}  d_i \geq 0,\ \displaystyle \sum_{i=-r}^s d_i = d_{x_0} - d_{x_1} \\
x_0^\ep = w_{-r}^\ep,\ x_1^\ep = w_{s+1}^\ep,\\
d_i = 0,\ i \neq 0, \Longrightarrow w_i^\ep \neq w_{i+1}^\ep
\end{array} \right. \right\}.
\eeqn
Let $A_{x_0; x_1}^\Omega \subset A_{x_0; x_1}$ be the open subset of tuples such that every $d_i$ is equal to the $\Omega$-area of some sphere class. Redefine ${\mc C}_{x_0; x_1}$ to be the restriction to this open subset. 

We can also do the outercollaring and rigidification construction in the equivariant way. We have the following analogue of Proposition \ref{prop245}. 

\begin{prop}
The collection of objects ${\mc C}_{x_0; x_1}$ forms a bimodule $\mb{Dom}^{{\mb B}_{01}^{\zp}}$ over $(\mb{Dom}^{{\mb F}_0^{\zp}}; \mb{Dom}^{{\mb F}^{\zp}_1})$, whose outercollaring admits an enrichment in $\outer \uds{\bf Curve}_{\rm rig}$ which is $\mb\Pi_\p$-equivariant.
\end{prop}

Following Section \ref{section_Floer_lift}, we can then construct a bimodule 
\beqn
\outer \mb{Map}_{\mb{B}^{\zp}_{01}}^{\mb{fr}}
\eeqn
over $(\outer \mb{Map}_{\mb{F}_0^{\zp}}^{\mb{fr}}, \outer \mb{Map}_{\mb{F}^{\zp}_1}^{\mb{fr}})$ enriched in $\outer \uds{\bf Map}^{\rm fr}$. Using the freeness of the $\mb\Pi_\p$-action, we learn that:

\begin{lemma}
    There exist $\mb\Pi_\p$-equivariant group reduction schemes lifting $\mb{Map}_{\mb{B}^{\zp}_{01}}^{\mb{fr}}$ to a flow bimodule enriched in $\outer \uds{{\bf map}}{}^{\rm fr} $.
\end{lemma}

Furthermore, following Proposition \ref{prop_perturbation}, the freeness of the $\mb\Pi_\p$-action guarantees the following:

\begin{lemma}
    There exists a $\mb\Pi_\p$-equivariant and transverse thickening datum on $\outer \mb{Dom}^{{\mb B}_{01}^{\zp}}$, i.e., a transverse lift to $\uds{\bf Thick}_{\rm rig}$.
\end{lemma}

Therefore, the above statements ensure that $\mb{B}_{01}^\zp$ admits a lift over $\outer \uds{\bf Kur}{}_{\rm rig}$. The constructions of normal complex structures and coherent smoothings follow from the recurring discussions of this paper. Therefore, by taking the quotient, we get
\beqn
\tilde{\mb{B}}_{01}^\zp
\eeqn
enriched in $\outer \uds{\bf dOrb}^{\rm NC}_{\rm rig}$ which is a $\mb\Pi_\p$-equivariant bimodule over $(\tilde {\mb F}_0^{\zp}; \tilde {\mb F}_1^{\zp})$. This finishes the proof of Theorem \ref{thm_equivariant_lift} (2).

\subsection{Equivariant PSS and SSP maps}

One uses spiked disks used to define the ordinary PSS and SSP maps. Here we would like to define equivariant analogues, which are maps
\begin{align*}
&\ \wh\Phi_{\zp}^{\mb{MF}}: H^*(X; \Lambda^{{\mb\Pi}}\langle \uptheta \rangle) \to HF_\zp^*(X),\ &\ \Phi_{\zp}^{\mb{FM}}: HF_{\zp}^*(X) \to H^*(X; \Lambda^{{\mb\Pi}} \langle \uptheta \rangle).
\end{align*}

To define these maps, we need to describe certain bimodules, denoted by $\mb{B}^{\mb{MF}, \zp}$ and $\mb{B}^{\mb{FM}, \zp}$ respectively. Consider the Morse-to-Floer map first. Choose a Morse flow category $\uds{\mb M}$ on the symplectic manifold $(X, \omega)$ associated to a Morse--Smale pair $(f^{\mb M}, g^{\mb M})$ and let ${\mb M}$ be the trivial product of $\uds{\mb M}$ with $\Pi$. Consider the coupled flow category 
\beqn
{\mb M}^\zp:=\mb{M} \# \ep.
\eeqn
The equivariant Morse-to-Floer bimodule is a bimodule $\mb{B}^{\mb{MF}, \zp}$ over $(\mb{M}^\zp; \mb{F}^\zp)$. 

We consider spiked disks equipped with a $\zp$-invariant Floer data. Regard $\Sigma^{\mb{MF}} \cong {\mb C}$ as a surface with cylindrical end which has a natural $\zp$-action which rotate the cylindrical coordinates by multiples of $\frac{2\pi}{\p}$. One can choose the Floer datum $\sigma^{\mb{MF}}$ on $\Sigma^{\mb{MF}}$ which is $\zp$-invariant, such that near $\infty$ it is equal to $H^{\natural \p} dt$. 

Now we describe the moduli spaces. Choose a pair of objects $x = (x^{\mb M}, x^{\ep}) \in {\rm Ob} \mb{M}^\zp$ and $y = (y^{\mb F}, y^\ep) \in {\rm Ob}{\mb F}^\zp$. Consider the set of solutions $(\gamma, u, v)$ where $\gamma: (-\infty, 0] \to X$, $u: \Sigma^{\mb{MF}} \to X$, and $v: {\mb R} \to S^\infty$ are smooth maps solving
\beqn
\gamma'(s) + \nabla f^{\mb M}(\gamma(s)) = 0,\ \ov\partial_{\sigma^{\mb{MF}}} u = 0,\ v'(s) + \nabla f^\ep(v(s)) = 0,\ \gamma(0) = u(0) \in X
\eeqn
and subject to the following constraint
\beqn
\lim_{s \to - \infty} \gamma(s) = x^{\mb M},\ \lim_{s \to -\infty} v(s) = x^\ep,\ \lim_{s \to +\infty} v(s) = y^\ep, \lim_{s \to +\infty} u(s,t) = y^{\mb F}.
\eeqn
We compactify the space of such solutions $(\gamma, u, v)$ by allowing synchronized breakings of Morse trajectories in $X \times S^\infty$ at the negative end and synchronized breakings of Floer trajectories in $X$ coupled with Morse trajectories in $S^\infty$ at the positive end. The compactified space is denoted by 
\beqn
M_{x; y}^{\mb{MF}, \zp}.
\eeqn
This space is stratified by the poset
\beqn
A_{x; y}^{\mb{MF}, \zp}:= \Big\{ (x_{-r} \cdots x_{-1}; y_1 \cdots y_s)\ |\ x_{-r} < \cdots < x_{-1},\ y_1 < \cdots < y_s \Big\}.
\eeqn
There are also natural structural maps
\begin{align*}
&\ \iota_{x_1x_0; y}^{\mb{MF}, \zp}: M_{x_1 x_0}^{\mb{M}^\zp} \times M_{x_0; y}^{\mb{MF}, \zp} \to M_{x_1; y}^{\mb{MF}, \zp},\ \ &\ \iota_{x; y_0y_1}^{\mb{MF}, \zp}: M_{x; y_0}^{\mb{MF},\zp}\times M_{y_0 y_1}^{\mb{F}^\zp} \to M_{x; y_1}^{\mb{MF}, \zp}
\end{align*}
satisfying obvious associativity conditions. Hence one obtains a bimodule $\mb{B}^{\mb{MF}, \zp}$ over $(\mb{M}^\zp; \mb{F}^\zp)$. 

Fix an outercollaring width and let the outercollaring of the involved flow categories be
\begin{align*}
&\ \outer \mb{M}^\zp,\ &\ \outer \mb{F}^\zp.
\end{align*}
The source flow category is already enriched in both $\outer \uds{\bf SKur}_{\rm rig}^{\rm NC}$ and $\outer \uds{\bf dOrb}_{\rm rig}^{\rm NC}$. The target flow category $\outer \mb{F}^\zp$ admits an AMS lift $\tilde{\mb F}^\zp$ coming from a lift $\hat{\mb F}^\zp$ in $\outer\uds{\bf SKur}_{\rm rig}^{\rm NC}$ (see Theorem \ref{thm_equivariant_lift}). Here we will establish the lift of the bimodule $\outer {\mb B}^{\mb{MF}, \zp}$ to $\outer \uds{\bf dOrb}_{\rm rig}^{\rm NC}$.

\begin{prop}\label{prop_equivariant_PSS_lift}
Given an equivariant AMS lift $\tilde{\mb F}^\zp$ of $\mb{F}^\zp$, there exists a class of lifts of the outercollaring of ${\mb B}^{\mb{MF}, \zp}$ to $\outer \uds{\bf dOrb}_{\rm rig}^{\rm NC}$, denoted by $\tilde{\mb B}^{\mb{MF}, \zp}$, called an equivariant AMS lift, which is a $\mb\Pi_\p$-equivariant Novikov flow category over $(\outer {\mb M} \# {\mb E}; \tilde {\mb F}^\zp)$. Moreover, two different equivariant AMS lifts are homotopic via a $\mb\Pi_\p$-equivariant Novikov homotopy. 
\end{prop}

The proof of this proposition is given momentarily.

Now using FOP perturbations one can define the equivariant PSS map on the chain level. Choose an FOP perturbation $\mathring{\mb F}^\zp$ of $\tilde{\mb F}^\zp$, one has the associated cochain complex
\beqn
C^*(\mathring {\mb F}^{\zp})
\eeqn
whose $\zp$-invariant part is the equivariant Floer cochain complex
\beqn
CF_{\zp}^*( H^{\natural \p}).
\eeqn
As ${\mb M}\# {\mb E}$ can be viewed as its own FOP perturbation, by the multimodule case of Theorem \ref{thma_FOP}, one can extend the FOP perturbation $\mathring {\mb F}^{\zp}$ to an FOP perturbation $\mathring {\mb B}^{\mb{MF}, \zp}$ on $\tilde {\mb B}^{\mb{MF}}$. Then one obtains a cochain map
\beqn
\wt\Phi^{\mb{MF}, \zp}: C^*( \mb{M}\# {\mb E}) \to C^* (\mathring {\mb F}^\zp)
\eeqn
which is $\Lambda_{\fp}^{{\mb\Pi}}$-linear and $\zp$-equivariant. 

\begin{lemma}
The $\zp$-equivariant chain homotopy class of $\wt\Phi^{\mb{MF}, \zp}$ is independent of choices. Moreover, if ${\mb M}_1, {\mb M}_2$ are two Morse flow category associated to two Morse--Smale pair $(f_1, g_1)$ and $(f_2, g_2)$ on $X$ and ${\mb F}_1, {\mb F}_2$ are Floer flow categories associated to $H_1^{\natural \p}$ and $H_2^{\natural \p}$, then the following diagram commutes
\beqn
\xymatrix{ C^*({\mb M}_1 \# {\mb E}) \ar[rr] \ar[d]  & &  C^*(\mathring {\mb F}_1^{\zp}) \ar[d] \\
          C^*({\mb M}_2 \# {\mb E}) \ar[rr] & & C^*(\mathring {\mb F}_2^\zp) }
\eeqn
up to $\zp$-equivariant chain homotopy. Here the horizontal arrows are the chain-level PSS maps and the vertical arrows are induced from Morse or Floer continuation maps.
\end{lemma}

Therefore, after restricting to the $\zp$-invariant part, one obtains a canonical $\Lambda_{\fp}^{{\mb\Pi}}$-linear map
\beqn
\Phi^{\mb{MF}, \zp}: H^*(X; \Lambda_{\fp}^{{\mb\Pi}} \langle \uptheta \rangle) \to HF_{\zp}^*(X).
\eeqn

In a similar way, one can construct a chain-level Floer-to-Morse (SSP) map, which induces a well-defined $\Lambda_{\fp}^{\mb\Pi}$-linear map
\beqn
\Phi^{\mb{FM}, \zp}: HF_{\zp}^*(X) \to H^*(X; \Lambda_{\fp}^{\mb\Pi}\langle \uptheta \rangle).
\eeqn

Just as in the nonequivariant case, we note that the maps $\Phi^{\mb{MF}, \zp}$ and $\Phi^{\mb{FM}, \zp}$ are isomorphisms.

\begin{thm}
The maps $\Phi^{\mb{MF}, \zp}$ and $\Phi^{\mb{FM}, \zp}$ are isomorphisms on cohomology.
\end{thm}

\begin{proof}
On the chain level, we can turn on a exhaustive, complete, and decreasing filtration on both $CF_{\zp}^*( H^{\natural \p})$ and $C^*( \mb{M}\# {\mb E})^{\mb{Z}_\p}$ by tracking the degree in $\tu$ and $\uptheta$. Then both the equivariant diffferentials and equivariant PSS/SSP maps respect the filtration. On the other hand, geometrically, the constructions in Section \ref{subsection_AMS_equivariant_Floer} can be performed such that the induced differentials and chain maps on the associated graded agrees with the ones coming from the ordinary Morse cohomology, the ordinary Hamiltonian Floer cohomology, and the ordinary PSS/SSP maps. Therefore, based on Theorem \ref{thma_PSS_isomorphism}, we conclude the theorem.
\end{proof}

\subsection{AMS construction for equivariant PSS and SSP maps}

We provide detailed construction for the PSS map. The case for SSP map is very similar. 

Similar to the non-equivariant case, we first choose an integral action ${\mc A}^\Omega: {\rm Ob}{\mb F} \to {\mb Z}$ which is compatible with respect to the bimodule ${\mb B}^{\mb{MF}, \zp}$ (or more precisely, compatible with the Floer data). 

\subsubsection{The domain bimodule for the AMS construction}\label{subsubsec:equiv-AMS}

We first describe the analogue of the space of domains, i.e., a bimodule enriched in $\uds{\bf Curve}$. Recall that the AMS construction for $\mb{F}^\zp$ (see Subsection \ref{subsection_AMS_equivariant_Floer}) provides a flow category $\mb{Dom}^{\mb{F}^\zp}$ enriched in $\outer \uds{\bf Curve}_{\rm rig}^{\mb C}$ whose objects are the same as those of $\mb{F}^\zp$. However, the flow category $\mb{M}^\zp$ does not yet have such a corresponding flow category enriched in $\uds{\bf Curve}$. We define such a flow category as follows. For each pair $x_0, x_1 \in {\rm Ob}{\mb M}^\zp$, define
\beqn
{\mc C}_{x_0 x_1}^{\mb{M}^\zp} = (G, B_{x_0 x_1}, C)
\eeqn
where $G$ is the trivial group, $B_{x_0x_1} = M_{x_0 x_1}^{\mb{M}^\zp}$ is the moduli space of flow lines between $x_0$ and $x_1$, and $C = \emptyset$, regarded as a family of empty curves over $B_{x_0 x_1}$. After outercollaring, we denote this rather trivial flow category enriched in $\outer \uds{\bf Curve}$ by 
\beqn
\mb{Dom}^{\mb{M}^\zp}.
\eeqn
As the concatenation of Morse flow lines does not involve stabilization and the Lie groups do not change, this flow category is automatically enriched in the rigidified version $\outer \uds{\bf Curve}_{\rm rig}$. However we do not consider stable complex structure. 

Now we describe a bimodule over $(\mb{Dom}^{\mb{M}^\zp} ; \mb{Dom}^{\mb{F}^\zp})$, denoted by $\mb{Dom}^{\mb{MF}, \zp}$. Fix a pair of objects
\begin{align*}
&\ x = (x^{\mb M}, x^\ep) \in {\rm Ob} {\mb M}^\zp,\ &\ y = (y^{\mb F}, y^\ep) \in {\rm Ob}{\mb F}^\zp.
\end{align*}
Notice that $x$ and $y$ determines an integer $d_{x;y}$ which is the difference of the symplectic area of $x^{\mb M}$ and the integrally modified symplectic action of $y^{\mb F}$.

Consider the the space of triples $(\gamma, \phi, \eta)$ where $\gamma: (-\infty, 0] \to X$, $\phi: {\mb R}\times S^1 \to \mb{CP}^{d_{x;y}}$, and $\eta: {\mb R} \to \wh S^\infty$ are smooth maps, solving the following equation
\beqn
\gamma'(s)  + \nabla f^{\mb M}(\gamma(s)) = 0,\ \ov\partial \phi = 0,\ \eta'(s) + \nabla f^\ep( \eta (s)) = 0
\eeqn
satisfying the following asymptotic constraints
\beqn
\lim_{s \to -\infty} \gamma(s) = x^{\mb M},\ \lim_{s \to -\infty} \eta(s) = x^\ep,\ \lim_{ s\to +\infty} \eta(s) = y^\ep
\eeqn
and the constraint
\beqn
\phi (+\infty) = [0, \ldots, 0, 1] \in \mb{CP}^{d_{x; y}}.
\eeqn
Notice that there is no requirement for $\gamma(0)$ and $u(-\infty)$. We compactify this space by allowing sychrnoized breakings in $X \times \wh S^\infty$ at the negative end and synchronized breakings of holomorphic cylinders in $\mb{CP}^{d_{x;y}}$ and Morse trajectories in $S^\infty$. Next, we remove those elements whose images in $\mb{CP}^{d_{x;y}}$ lie in any hyperplane. Denote this space by 
\beqn
B_{x; y}.
\eeqn
Then this is a smooth stratified manifold with a ${\mc G}_{d_{x;y}}$-action. On top of it there is the universal curve
\beqn
C_{x; y} \to B_{x;y}.
\eeqn
Then the triple 
\beqn
{\mc C}_{x;y}:= ( {\mc G}_{x;y}, B_{x;y}, C_{x;y})\ {\rm where}\ {\mc G}_{x;y} = {\mc G}_{d_{x;y}}
\eeqn
is an object of $\uds{\bf Curve}$. We can also define structural maps
\beqn
{\mc C}_{x_0 x_1}^{{\mb M}^\zp} \times {\mc C}_{x_1; y}^{\mb{MF}, \zp} \to {\mc C}_{x_0; y}^{\mb{MF}, \zp},\ {\mc C}_{x; y_0}^{\mb{MF},\zp} \times {\mc C}_{y_0 y_1}^{{\mb F}^\zp} \to {\mc C}_{x; y_1}^{\mb{MF}, \zp}.
\eeqn
This provide a flow bimodule  $\mb{Dom}^{\mb{MF}, \zp}$ over $(\mb{Dom}^{{\mb M}^\zp}; \mb{Dom}^{\mb{F}^\zp})$. One can check that this bimodule is $\mb\Pi_\p$-equivariant.

\subsubsection{Thickening}

One can consider the space of framed maps ${\rm Map}_{x; y}^{\rm fr}(X)$ whose elements are $(\phi, u, F)$ where $\phi\in B_{x; y}$, $u: C_\phi \to X$ is an admissible map, and $F$ is a frame determined by $u$ and the integral symplectic form $\Omega$. Recall that $\phi$ carries the information of a flow ray in $X$ and we do not require that the evaluation of $u$ at $-\infty$ coincides with the value of the flow ray at $0$. One chooses an equivariant group reduction scheme 
\beqn
\lambda_{x; y} : \outer {\rm Map}_{x;y}^{\rm fr}(X) \to Q_{x;y}:= Q_{d_{x;y}}
\eeqn
which extends the group reduction scheme chosen for the AMS construction for $\mb{F}^\zp$. On the other hand, we also choose thickening data
\beqn
\nu_{x;y}: 
W_{x;y} \to \Gamma_c( \mathring C_{x;y} \times X, \Lambda_{\mathring C_{x,y}/ B_{x;y}}^{0,1} \otimes TX)
\eeqn
extending the thickening data chosen for the AMS construction for $\hat {\mb F}^{\zp}$. Define the {\bf thickening}
\beqn
\check V_{x;y} := \left\{ (\phi, u, F, e)\ \left|\  \begin{array}{c} (\phi, u, F) \in \outer {\rm Map}_{x; y}^{\rm fr}(X),\ e \in W_{x;y}^{\mb{MF}, \zp} \\
\ov\partial u + \nu_{x;y} (e) = 0,\ \phi = \phi_F 
\end{array}\right. \right\}.
\eeqn
The obstruction bundle $E_{x;y} $ is the trivial bundle with fiber $W_{x;y} \oplus Q_{x;y}$ and the Kuranishi section is
\beqn
S_{x;y} (\phi, u, F, e) = (e, \lambda_{x;y}(\phi, u, F)).
\eeqn

We summarize the construction so far in the following proposition.

\begin{prop}\label{prop2110}
The collection of Kuranishi spaces 
\beqn
\check K_{x;y}:= (G_{x;y}, \check V_{x;y}, E_{x;y}, S_{x;y})
\eeqn
together with the structural maps
\beqn
\outer M_{x_0 x_1}^{\mb{M}^\zp} \times K_{x_1; y} \to K_{x_0; y}
\eeqn
and
\beqn
K_{x; y_0} \times K_{y_0 y_1}^{\mb{F}^\zp} \to K_{x; y_1}
\eeqn
define a bimodule $\check {\mb X}^{\mb{MF}, \zp}$ over $(\outer \mb{M}^\zp; \hat{\mb F}^\zp)$ enriched in $\uds{\bf S^{\rm rel} Kur}_{\rm rig}$.
\end{prop}

Notice that $\check{\mb X}^{\mb{MF}, \zp}$ is not a lift of ${\mb X}^{\mb{MF}, \zp}$ as we did not impose the matching condition in $X$.

We then construct the vertical NC structure. Similar to the case for $\hat {\mb F}^{\zp}$, the construction can be carried out in a $\zp$-equivariant fashion. Hence one obtains a bimodule enriched in $\outer \uds{\bf S^{\rm rel} Kur}_{\rm rig}^{\rm NC}$.

We still apply the same smoothing strategy to obtain a stable smoothing of $\check {\mb M}^{\mb{MF}, \zp}$. To simplify notations, we still denote by the same symbol. Moreover, we may require that the natural map 
\beqn
\ev: \check V_{x; y}^{\mb{MF}, \zp} \to X \times X
\eeqn
is a submersion, hence transverse to the diagonal. Define
\beqn
V_{x;y}^{\mb{MF}, \zp}:= \ev^{-1}(\Delta_X)
\eeqn
and define $K_{x;y}$ be the restriction of $\check K_{x;y}$ to this closed submanifold. Finally, incoporating the stable complex structure on the domain bimodule, one obtains an NC structure. We summarize the construction as follows.

\begin{prop}
The collection of Kuranishi spaces $K_{x;y}$ together with the restriction of the structural morphisms of Proposition \ref{prop2110} define a lift of $\mb{X}^{\mb{MF}, \zp}$ to $\outer \uds{\bf SKur}_{\rm rig}^{\rm NC}$, denoted by $\hat{\mb X}^{\mb{MF}, \zp}$, which is a bimodule over $(\outer {\mb M}\# {\mb E}; \hat {\mb F}^\zp)$. 
\end{prop}

Finally, by taking the quotient by the Lie group actions, we obtain ${\mb B}^{\mb{MF}, \zp}$. One can further follow the machinery of this paper to address the invariance property. This concludes Proposition \ref{prop_equivariant_PSS_lift}.

\section{The Equivariant Pair-of-Pants---Proof of Theorem \ref{thma_Floer_Steenrod}}\label{section_Floer_Steenrod}

The equivariant pair-of-pants product can be viewed as a Floer-theoretic analogue of the quantum Steenrod operation. It is defined via a multimodule associated to a Floer domain having a $\zp$-symmetry. The main results regarding such a construction is stated in Theorem \ref{thma_Floer_Steenrod}. Given the recurring nature of our arguments, we highlight the distinctive features.

\subsection{The equivariant pair-of-pants multimodule}

Let $(X, \omega, J)$ be a compact almost K\"ahler manifold. Again, we do not discuss the dependence on $J$, which could be done without new difficulty. Let ${\mb F}$ be the Floer flow category of a Hamiltonian $H$ and let ${\mb F}^\zp$ be the coupling of the Floer flow category of the $\p$-th iteration $H^{\natural \p}$ with the Morse flow category $\ep$ on $\wh S^\infty$.\footnote{In general one could allow the $\p$-fold iteration of a different Hamiltonian.} We describe a multimodule $\mb{X}^{\mb{FSt}}$ over $({\mb F}, \cdots, {\mb F}, \ep; {\mb F}^\zp)$ where ${\mb F}$ appears $\p$ times in the source.

\subsubsection{The Floer domain and the equation}

Let $z_1, \ldots, z_{\p}\in  {\mb C} \subset \mb{CP}^1$ be all the $\p$-th roots of unity. Let 
\beqn
\Sigma:= \mb{CP}^1 \setminus \{z_1, \ldots, z_{\p}, z_\infty = \infty \}
\eeqn
be the punctured sphere, which has a holomorphic $\zp$-action permuting the finite punctures. Choose a $\zp$-invariant Floer datum $\sigma^{\zp}$ on $\Sigma$ whose restriction on the $\p$ negative ends are $H dt$ and whose restriction on the positive end is the $\p$-th iteration $H^{\natural \p} dt$. In fact, one can choose different $H$ for negative and positive ends as in the setup of Theorem \ref{thma_Floer_Steenrod}, but we spell out the argument for $H_- = H_+ = H$ to simplify the notations. Solutions to the Floer equation
\beqn
\ov\partial_{\sigma^{\zp}} u = 0,\ u: \Sigma \to X
\eeqn
converge to periodic orbits of $H$ at negative ends and periodic orbits of $H^{\natural \p}$ at the positive end. 

\subsubsection{The moduli spaces and the multimodule}

We define the morphism spaces of the multimodule ${\mb X}^\zp$. Choose a collection of objects
\beqn
(x_1^{\mb F}, \ldots x_\p^{\mb F}, x^\ep) \in \underbrace{{\rm Ob} {\mb F} \times \cdots \times {\rm Ob}{\mb F}}_{\p} \times {\rm Ob} {\mb E}
\eeqn
and an object $y = (y^{\mb{F}^{\natural \p}}, y^\ep) \in {\rm Ob} {\mb F}^\zp$, where $y^{\mb{F}^{\natural \p}}$ is the flow category associated with the Hamiltonian $H^{\natural \p}$, we describe a moduli space as follows. Consider pairs
\beqn
u: \Sigma \to X,\ v: {\mb R} \to \wh S^\infty
\eeqn
where $u$ is a solution to $\ov\partial_{\sigma^{\zp}} u = 0$ converging to $x_i^{\mb F}$ at the $i$-th negative cylindrical ends, $i = 1, \ldots, \p$, to $y^{\mb{F}^{\natural \p}}$ at the positive cylindrical end, and where $v$ is a Morse flow line of $f^\ep$ connecting $x^\ep$ and $y^\ep$ at negative and positive ends respectively. We compactify the space of such solutions in the following way: while allowing sphere bubbling, the breakings at negative ends are treated as independent breaking while the breakings at positive ends are synchronized. Let the compactified space be $M_{x_1^{\mb F}\cdots x_\p^{\mb F} x^\ep; y}^{{\mb X}^{\mb{FSt}}}$.

There are obvious inclusions given by breakings of the domains. For $i = 1, \ldots, \p$  and $x_i^{{\mb F}} \leq w_i^{{\mb F}}$ in the $i$-th source flow category, there is a map 
\beqn
M_{x_i^{\mb F} w_i^{\mb F}}^{{\mb F}} \times M_{x_1^{\mb F} \cdots x_{i-1}^{\mb F} w_i^{\mb F} x_{i+1}^{\mb F} \cdots x_\p^{\mb F} x^{\ep}; y}^{\mb{X}^{\mb{FSt}}} \to M_{x_1^{\mb F} \cdots x_\p^{\mb F} x^\ep; y}^{\mb{X}^{\mb{FSt}}};
\eeqn
for $x^\ep \leq w^\ep$ in $\ep$, a map
\beqn
M_{x^\ep w^\ep}^\ep \times M_{x_1^{\mb F} \cdots x_\p^{\mb{F}} w^\ep; y}^{{\mb X}^{\mb{FSt}}} \to M_{x_1^{\mb F} \cdots x_\p^{\mb F} x^\ep; y}^{\mb{X}^{\mf{FSt}}};
\eeqn
for $w \leq y$ in $\mb{F}^{\zp}$, a map
\beqn
M_{x_1^{\mb F} \cdots x_\p^{\mb F} x^\ep; w}^{\mb{X}^{\mb{FSt}}} \times M_{wy}^{\mb{F}^\zp} \to M_{x_1^{\mb F}\cdots x_\p^{\mb F} x^\ep; y}^{\mb{X}^{\mb{FSt}}}.
\eeqn
These maps satisfy the conditions required for a multimodule. Hence one obtains a multimodule ${\mb X}^{\mb{FSt}}$ over $(\mb{F}, \ldots, {\mb F}, \ep; {\mb F}^\zp)$ enriched in $\uds{\bf Top}$.

\subsubsection{The symmetry}

Define the Novikov monoid
\beqn
\mb\Pi^{(\p)}:= \underbrace{\Pi \times \cdots \times \Pi}_{\p} \times {\mb Z}.
\eeqn
There is a $\zp$-action on $\mb\Pi^{(\p)}$ by isomorphisms of Novikov monoids given by permuting the first $\p$ factors in the above product. Then we can form the semi-direct product $\mb\Pi^{(\p)}\rtimes \zp$. 

We can view the multimodule ${\mb X}^{\mb{FSt}}$ as a weak bimodule over $({\mb F}\times \cdots \times {\mb F} \times \ep; {\mb F}^\zp)$. On the source weak flow category, it has the action 
\beqn
{\mc A}(x_1^{\mb F}, \ldots, x^{\mb F}_\p, x^{\ep}) = ({\mc A}_H(x_1^{\mb F}), \ldots, {\mc A}_H(x_\p^{\mb F}), {\mc A}_{\ep}(x^\ep)) \in {\mb R}^{\p + 1}.
\eeqn

\begin{prop}
The collection of moduli spaces $M_{x_1^{\mb F} \cdots x_\p^{\mb F} x^\ep; y}^{{\mb X}^{\mb{FSt}}}$ define a Novikov multimodule over $({\mb F}, \cdots, {\mb F}, \ep; {\mb F}^\zp)$ which is equivariant with respect to the natural morphism
\beqn
\rho: \mb\Pi^{(\p)} \to \mb\Pi,\ (a_1, \ldots, a_\p, l) \mapsto (a_1 + \cdots + a_\p, l)
\eeqn
of Novikov monoids. Moreover, as a weak bimodule over $({\mb F}\times \cdots \times {\mb F} \times \ep; {\mb F}^\zp)$, it is equivariant with respect to $\rho$.
\end{prop}

\begin{proof}
Straightforward.
\end{proof}

\subsubsection{Regularization and the equivariant pair-of-pants product}

Fix an outercollaring width and a compatible choice of integral actions. Recall that the outercollarings of ${\mb F}$ and ${\mb F}^\zp$ both admit AMS lifts to $\outer \uds{\bf dOrb}{}_{\rm rig}^{\rm NC}$.

\begin{prop}
There exists a lift of the outercollaring of ${\mb X}^{\mb{FSt}}$, denoted by $\tilde {\mb X}^{\mb{FSt}}$, to $\outer \uds{\bf dOrb}^{\rm NC}_{\rm rig}$, which is a $\rho$-equivariant Novikov multimodule over $(\tilde {\mb F}, \ldots, \tilde {\mb F}, \outer \ep; \tilde {\mb F}^\zp)$.
\end{prop}

\begin{proof}
The construction is analogous to the AMS construction for a general multimodule induced from a smooth Floer domain, except that the moduli spaces curves now has the parametrization factor from the moduli spaces of Morse trajectories in $\wh S^{\infty}$ as in Section \ref{subsubsec:equiv-AMS}. One can carry out the construction almost verbatim; the $\zp$-equivariance condition can be achieved because the action is free.  
\end{proof}

Now let $\mathring {\mb F}$ be an FOP transverse perturbation on $\tilde{\mb F}$ and $\mathring {\mb F}^\zp$ be a $\zp$-equivariant FOP transverse perturbation on $\tilde{\mb F}^\zp$. Let $CF(H)$ resp. $C(\mathring {\mb F}^\zp)$ be the associated complex over $\Lambda_{\fp}^\Pi$ resp. over $\Lambda_{\fp}^{{\mb\Pi}}$; the $\zp$-invariant part of $C(\mathring {\mb F}^\zp)$ is the equivariant complex $CF_{\zp}(H^{\natural \p})$. By the multimodule case of Theorem \ref{thma_FOP}, one can extend these FOP transverse perturbations to a $\zp$-equivariant FOP transverse perturbation $\mathring {\mb X}^{\mb{FSt}}$ on $\tilde {\mb X}^{\mb{FSt}}$. Then by Proposition \ref{prop_chain_map}, it induces a $\fp$-linear chain map 
\beqn
\wt P: \overbrace{  CF(H) \widehat{\underset{\fp}{\otimes}} \cdots \widehat{\underset{\fp}{\otimes}} CF(H)  }^{\p} \widehat{\underset{\fp}{\otimes}} C ( \ep ) \to C(\mathring {\mb F}^\zp).
\eeqn
The $\rho$-equivariance condition on the multimodule $\mathring{\mb X}^{\mb{FSt}}$ implies that $\Phi^{\rm FSt}$ induces a map
\beqn
\wt P: \left( \overbrace{ CF(H) \widehat{\underset{\fp}{\otimes}} \cdots \widehat{\underset{\fp}{\otimes}} CF(H)  }^{\p} \widehat{\underset{\fp}{\otimes}} CM(f^\ep) \right) \underset{ \Lambda_{\fp}^{{\mb \Pi}^{(\p)}}}{\otimes} \Lambda_{\fp}^{\mb\Pi}
\to C(\mathring {\mb F}^\zp)
\eeqn
which is $\zp$-equivariant and $\Lambda_{\fp}^{\mb\Pi}$-linear.

Restricting to the $\zp$-invariant part, one obtains the chain-level equivariant pair-of-pants product
\beqn
P: C_{\zp}( CF(H)^{\otimes \p})  \to CF_{\zp}(H^{\natural \p})
\eeqn
hence a map on the cohomology level (denoted by the same notation)
\beqn
P: H_{\zp}(CF(H)^{\otimes \p}) \to HF_{\zp} (H^{\natural \p}).
\eeqn

\subsection{Proof of Theorem \ref{thma_Floer_Steenrod}}
As above, we write down the proof for $H_- = H_+ = H$ to simplify the notations.

\subsubsection{Independence of choices}

\begin{lemma}
$\wt P$ is well-defined up to $\zp$-equivariant chain homotopies. Therefore, restricting to the ${\mb Z}_\p$-invariant part, $P$ is well-defined up to chain homotopies.
\end{lemma}

\begin{proof}
Two concrete chain level constructions can be connected by a homotopy: we can carry out the construction of multimodule homotopies in the ${\mb Z}_\p$-equivariant setting using families of curves with parametrization having the additional factor from the moduli spaces of Morse trajectories in $\wh S^{\infty}$, and the ${\mb Z}_\p$-equivariance can be achieved due to the freeness.
\end{proof}

Therefore, following the definition of $\fst$ in Theorem \ref{thma_Floer_Steenrod}, we see that the Floer-theoretic Steenrod operation does not depend on any choices after passing to cohomology. This finishes proving Theorem \ref{thma_Floer_Steenrod} (1) and (3).

\subsubsection{Isomorphism on cohomology}
Theorem \ref{thma_Floer_Steenrod} (2) follows from a stronger assertion in the quantitative setting, Corollary \ref{cor_eq_pants_filtered}, which is proven in the next Section. The idea is that instead of directly constructing a chain-homotopy inverse of $P$ by dualizing the equivariant pair-of-pants product, we turn on the action filtration on the (equivariant) Floer chain complexes, which the operation $P$ respects, so that the induced operation on the associated graded reduces to the equivariant pants product on local Floer cohomology. Then, by appealing to the known isomorphism property in \cite{seidel-pants, shelukhin-zhao} in the aspherical setting, which, in particular, applies to local Floer theory, we conclude that $P$ is an isomorphism on cohomology.  

\begin{rem}
    Using Proposition \ref{prop:qf-iso}, we further see that $\fst$ is an isomorphism. This is the generalization of the localization results in \cite{seidel-pants, shelukhin-zhao} for general compact symplectic manifolds in Hamiltonian Floer theory. We expect that the method in this paper can also be used for proving an analogous result for the fixed point Floer cohomology associated with a general symplectic automorphism, though one needs a slight different regularization procedure.
\end{rem}

\subsubsection{The multimodule homotopy}

Finally, we prove the commutativity of Diagram \eqref{qst_fst_commute} in Theorem \ref{thma_Floer_Steenrod} (4). The core of the proof is the following chain-level assertion. As we only need the cohomological level conclusion, one fixes a Hamiltonian $H = H_- = H_+$ such that $H^{\natural \p}$ is nondegenerate. Moreover, fix a Morse--Smale pair $(f, g)$ on $X$ with associated Morse flow category $\uds {\mb M}$ (whose cohomology is isomorphic to $H(X; \Lambda_{\mb Z}^\Pi)$). 

\begin{lemma}\label{lemma_qst_fst}
The following diagram commutes up to $\zp$-equivariant homotopy.
\beqn
\begin{tikzcd}[column sep=large]
{\overbrace{C({\mb M}) \otimes \cdots \otimes C({\mb M})}^{\p} \otimes CM(f^\ep)} \arrow[rr, "\wt\Phi^{\mb{QSt}}"] \arrow[d, "{(\Phi^{\mb{MF}})^{\otimes \p}}"'] & & {C(\mb{M}) \otimes CM(f^\ep)} \\
{\underbrace{CF(H) \otimes \cdots \otimes CF(H)}_{\p} \otimes CM(f^\ep)} \arrow[rr, "\wt P"'] & & {CF(H^{\natural \p}) \otimes CM(f^\ep)} \arrow[u, "\Phi_{\zp}^{\mb{FM}}"']
\end{tikzcd}
\eeqn
\end{lemma}

Assuming this lemma, the restriction of the commutative diagram to the $\zp$-invariant part reads
\beqn
\vcenter{ \xymatrix{   C_{\zp}( C({\mb M})^{\otimes \p})   \ar[rr] \ar[d] & &   C({\mb M}) \otimes \tkp \\ 
    C_{\zp}(CF(H)^{\otimes \p}) \ar[rr] & &  CF_{\zp}(H^{\natural \p}) \ar[u] }  }.
\eeqn
The objects in this diagram are all modules over $\Lambda_{\kp}^\Pi$. As the quasi-Frobenius map respects equivariant chain homotopies, after passing to cohomology, one obtains the commutative diagram
\beq\label{qst_fst_diagram_2}
\vcenter{\xymatrix{ H(X; \Lambda_{\fp}^\Pi)  \ar[rr]^-{qF} \ar[d]_{\Phi^{\mb{MF}}} & & H_{\zp}(C({\mb M})^{\otimes \p}) \ar[rr]^{\Phi^{\mb{QSt}}} \ar[d]   &  &    H(X; \Lambda_{\tkp}^\Pi)  \\
            HF(X) \ar[rr]_-{qF} & &    H_{\zp} (CF(H)^{\otimes \p}) \ar[rr]_{P}   & & HF_{\zp}(X) \ar[u]_{\Phi_{\zp}^{\mb{FM}}}     }  }.
\eeq
The outer layer of this diagram is just \eqref{qst_fst_commute}. This would finish the proof of Theorem \ref{thma_Floer_Steenrod}.

Now we prove Lemma \ref{lemma_qst_fst}. The proof is based on a multimodule homotopy between the following two multimodules. The first one is the Steenrod multimodule $\mb{X}^{\mb{QSt}}$ over $(\mb{M}, \ldots, \mb{M}, \ep; \mb{M}\# \ep)$, which will lead to the quantum Steenrod operation $\qst$. The second one is the concatenation of three multimodules: first, there is the non-equivariant Morse-to-Floer bimodule $\mb{B}^{\mb{MF}}$. Consider the product of $\p$ copies of this bimodule with a copy of the diagonal bimodule $\mb{\Delta}^{\ep}$
\beqn
(\mb{B}^{\mb{MF}})^\p \times \mb{\Delta}^{\ep}:= \underbrace{\mb{B}^{\mb{MF}}\times \cdots \times \mb{B}^{\mb{MF}}}_{\p} \times \mb{\Delta}^{\ep \ep},
\eeqn
which is a multimodule over $(\mb{M}, \ldots, {\mb M}, \ep; \mb{F}, \ldots, \mb{F}, \ep)$. It can be composed with the multimodule ${\mb X}^{\mb{FSt}}$, followed by the equivariant SSP bimodule ${\mb B}_{\zp}^{\mb{FM}}$.

Before we move on, we specify the symmetries of the flow categories and multimodules. The flow categories ${\mb M}$ and ${\mb F}$ carry the action by the group $\Pi$ and $\ep$ carries the action by ${\mb Z}$. The flow category $\mb{M}\# \ep$ carries the action by (the trivially twisted) $\Pi \times \mb{Z} \times \zp$. The multimodule $\mb{X}^{\mb{QSt}}$ is equivariant with respect to the $\zp$-equivariant morphism of Novikov monoids
\beqn
\rho^{\mb{St}}: \Pi \times \cdots \times \Pi \times {\mb Z} \to \Pi \times {\mb Z}.
\eeqn
Moreover, as a weak bimodule (see Definition \ref{defn_weak_multimodule}), $\mb{X}^{\mb{QSt}}$ is $\zp$-equivariant. It is easy to see that the same is true for the triple concatenation. 

Now fix a common outercollaring width. Also choose an integral action ${\mc A}_{\mb F}^\Omega$ on ${\mb F}$ and an integral action ${\mc A}_{{\mb F}^{\natural \p}}^\Omega$ which are compatible with respect to involved Floer data. 

Suppose $\tilde {\mb F}$ is an AMS lift of ${\mb F}$ subject to the outercollaring width and the action ${\mc A}_{\mb F}^\Omega$. Let $\tilde{\mb F}_\zp$ be an equivariant AMS lift of the outercollaring of ${\mb F}_\zp$ subject to the outercollaring width and the action ${\mc A}_{{\mb F}^{\natural \p}}^\Omega$. Let $\tilde{\mb B}^{\mb{MF}}$ is an AMS lift of the outercollaring of $\mb{B}^{\mb{MF}}$ as a bimodule over $(\outer {\mb M}, \tilde{\mb F})$. Let $\tilde{\mb X}^{\mb{FSt}}$ be an equivariant AMS lift of the outercollaring of ${\mb X}^{\mb{FSt}}$ extending $\tilde{\mb F}$ and $\tilde{\mb F}_{\zp}$. Let $\tilde{\mb B}_{\zp}^{\mb{FM}}$ be an equivariant AMS lift of the outercollaring of ${\mb B}_{\zp}^{\mb{FM}}$. 

\begin{lemma}\label{lemma_qst_fst_lift}
There exists a $\rho^{\mb{St}}$-equivariant homotopy $\tilde{\mb H}$ from $\tilde{\mb X}^{\mb{QSt}}$ to 
\beqn
\Big( (\tilde{\mb B}^{\mb{MF}})^\p \times \outer\mb{\Delta}^{\ep \ep} \Big) \circ \tilde{\mb X}^{\mb{FSt}} \circ \tilde{\mb{B}}_{\zp}^{\mb{FM}}
\eeqn
enriched in the category $\outer\uds{\bf dOrb}_{\rm rig}^{\rm NC}$.
\end{lemma}

Assuming this lemma, we can conclude the assertion of Lemma \ref{lemma_qst_fst}. Since all actions are free, one can choose equivariant FOP transverse perturbations, firstly on the flow categories $\tilde{\mb{F}}$ and $\tilde{\mb{F}}^\zp$, then on the multimodules, and then on the homotopy $\tilde{\mb H}$, thanks to Theorem \ref{thma_FOP}. Then the corresponding homotopy $\mathring{\mb H}$ enriched in $\pman$ induces a chain homotopy $\Psi^{\mathring{\mb H}}$ between the two chain maps, the chain map corresponding to the quantum Steenrod operation and the composition of the three chain maps. This means the diagram \eqref{qst_fst_diagram_2} commutes up to homotopy.

\subsubsection{Proof of Lemma \ref{lemma_qst_fst_lift}}

We first describe the homotopy in the category $\uds{\bf Top}$ before the outercollaring, from $\mb{X}^{\mb{QSt}}$ to the concatenation
\beqn
\mb{Y}:=\left( ( \mb{B}^{\mb{MF}})^\p \times \mb{\Delta}^{\ep\ep} \right) \circ \mb{X}^{\mb{FSt}} \circ \mb{B}_{\zp}^{\mb{FM}}.
\eeqn
Consider the sphere $\Sigma_\p^{\mb{St}} \cong  \mb{CP}^1$ equipped with a 1-parameter family of Floer data $\sigma_t$  parametrized by $t \in [0, 1]$ where, when $t = 0$, $\sigma_t \equiv 0$ while as $t \to 1$, the long necks around each negative markings have the Hamiltonian connection $H dt$ which is $\zp$-equivariant while the long neck around the positive marking has the Hamiltonian $H^{\natural\p} dt$. Then for any $t$, $\sigma_t$ is invariant under the $\zp$-action on the domain. 

\begin{figure}[h]
    \centering
    \includegraphics[width=1\linewidth]{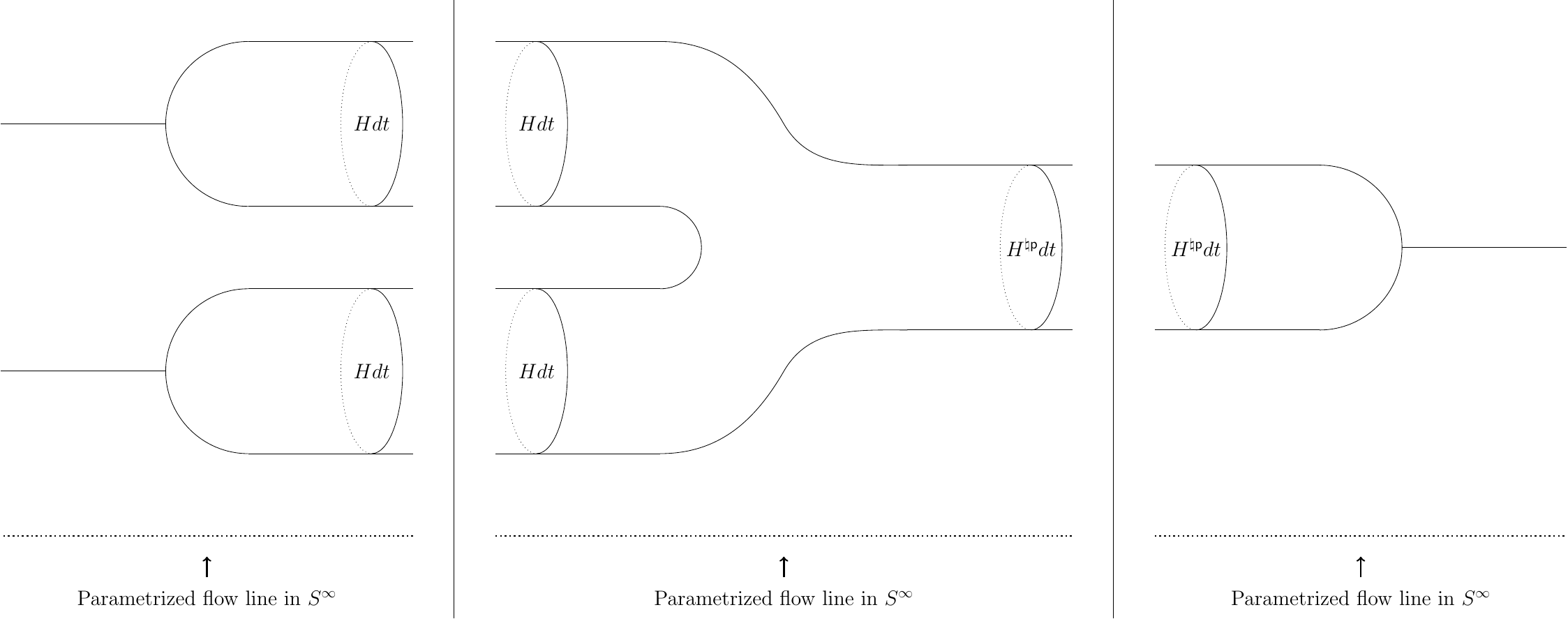}
    \caption{The concatenation of three multimodules. Notice that the gluing on the left have $\p +1$ gluing parameters and the gluing on the right has only one gluing parameter.}
    \label{fig_steenrod_homotopy}
\end{figure}

Now we can describe the moduli spaces responsible for the homotopy. Choose
\beqn
x_1^{\mb M}, \ldots, x_\p^{\mb M} \in {\rm Ob}{\mb M},\ x^\ep \in {\rm Ob}\ep,\ y = (y^{\mb M}, y^\ep) \in {\rm Ob} {\mb M}^\ep.
\eeqn
Consider the moduli space of triples $(t, u, v)$ where $t \in [0, 1)$, $u: \Sigma_\p^{\mb{St}} \to X$ and $v: {\mb R} \to \wh S^\infty$ are smooth maps solving
\beqn
(\nabla^{\sigma_t})^{0,1} u = 0,\ v'(s) + \nabla f^\ep (v(s)) = 0
\eeqn
and 
\beqn
u(z_i) \in W^{\rm us}_{x_i^{\mb M}},\ i = 1, \ldots, \p,\ u(z_\infty) \in W^{\rm st}_{y^{\mb M}}, \lim_{s \to -\infty} v(s) = x^{\ep},\ \lim_{s \to +\infty} v(s) = y^\ep.
\eeqn
Let the space of solutions be $\mathring M_{x_1^{\mb M}\cdots x_\p^{\mb M} x^\ep; y}^{{\mb H}}$. We may regard the constraints of $u$ at $z_i$ as adding gradient rays of $f$ in $X$. Then let
\beqn
M_{x_1^{\mb M}\cdots x_\p^{\mb M} x^\ep; y}^{{\mb H}}
\eeqn
be its natural compactification which allows sphere bubbling, Morse trajectory breaking, and breakings along Floer trajectories as $t \to 1$. Moreover, when Morse trajectories break in the positive end, the trajectories broken off are regarded as flow lines in the product $X \times \wh S^\infty$. The following lemma can be verified easily.

\begin{lemma}
The collection of spaces $M_{x_1^{\mb M}\cdots x_\p^{\mb M} x^\ep; y}^{{\mb H}}$ forms a homotopy from ${\mb X}^{\mb{QSt}}$ to the concatenation ${\mb Y}$. Moreover, this is a $\rho^{\mb{St}}$-equivariant Novikov homotopy and $\zp$-equivariant as a weak homotopy.
\end{lemma}

\subsubsection{The AMS construction}

Now we describe the AMS lift of the outercollaring of the homotopy ${\mb H}$, given AMS lifts of involved Floer flow categories and multimodules. 

Recall the AMS lift $\tilde {\mb F}$ of the Floer flow category ${\mb F}$ contains a domain monotone flow category $\outer {\mb Dom}$ enriched in $\outer \uds{\bf Curve}_{\rm rig}^{\mb{C}}$. The AMS lift $\tilde {\mb F}_{\zp}$ also contains a domain monotone flow category $\outer \mb{Dom}_{\zp}$ enriched in $\outer \uds{\bf Curve}_{\rm rig}^{\mb{C}}$ which also has a $\zp$-action corresponding to the domain rotation. 

The lift $\tilde {\mb X}^{\mb{QSt}}$ contains the domain monotone multimodule $\outer \mb{Dom}^{\mb{QSt}}$. The lift $\tilde {\mb X}^{\mb{FSt}}$ has a similar monotone multimodule $\outer \mb{Dom}^{\mb{FSt}}$. There are also the monotone bimodules $\outer \mb{Dom}^{\mb{MF}}$ and $\outer \mb{Dom}_{\zp}^{\mb{FM}}$. Then similar to previous cases of obtaining monotone homotopies involving concatenations and gluing, one can define a monotone homotopy from $\outer \mb{Dom}^{\mb{QSt}}$ to the concatenation
\beqn
(\outer \mb{Dom}^{\mb{MF}})^{\p} \circ \outer \mb{Dom}^{\mb{FSt}} \circ \outer \mb{Dom}_\zp^{\mb{FM}}.
\eeqn
This monotone homotopy, denoted by $\outer \mb{Dom}^{\mb H}$, is enriched in $\uds{\bf Curve}_{\rm rig}^{\mb{C}}$ and is $\zp$-equivariant. 

The remaining construction will be a tedious repetition of existing constructions. We only give a sketch. For a tuple of objects $x_1^{\mb M}, \ldots, x_\p^{\mb M}, x^{\mb E}; y$, by using the chosen integral actions, one can pull back the monotone homotopy to obtain objects
\beqn
{\mc C}_{x_1^{\mb{M}} \cdots x_\p^{\mb M} x^\ep; y}^{\mb H}\in {\rm Ob}\outer \uds{\bf Curve}_{\rm rig}^{\mb C}
\eeqn
which describe configurations containing not only surfaces but also flow lines/rays. The choices we need for the AMS construction can all be extended from existing ones, including the group reduction and thickening data. All the extensions can be made in a $\zp$-equivariant way. 

After these choices, one obtains a $\zp$-equivariant lift of the outercollaring of ${\mb H}$ to $\outer \uds{\bf S^{\rm rel} Kur}_{\rm rig}$. There is no essential difference from previous cases of constructing vertical NC structures in a $\zp$-equivariant way. The stable smoothings can also be made equivariantly, first without requiring the matching condition between flow rays and maps defined on surfaces. The submersive condition implies that the incidence condition cuts down smooth Kuranishi spaces. By incorporating the domain stable complex structures, one obtains the lift of ${\mb H}$ to $\outer \uds{\bf SKur}_{\rm rig}^{\rm NC}$. By passing to the orbifold quotient, we obtain the final construction.

\section{Quantitative Results about Equivariant Floer Theory}\label{section_equivariant_quantitative}

In this section we discuss the quantitative features of the equivariant Floer theory and prove Theorem \ref{thm_pants_filtered}. 

\subsection{Filtration on the $\p$-fold tensor product}

Let $H$ be a nondegenerate Hamiltonian on $X$. Let $CF(H)$ be a Floer complex associated to $H$, obtained by a particular AMS construction and FOP perturbation. The source of the equivariant pair-of-pants product is the complex
\beqn
C_{\zp}(CF(H)^{\otimes \p})
\eeqn
which is the $\zp$-invariant part of 
\beqn
CF(H)^{\otimes \p} \otimes CM(f^\ep).
\eeqn
By our convention, this tensor product is a module over the Novikov field
\beqn
\Lambda_{\fp}^{\Pi\oplus \cdots \oplus \Pi \oplus {\mb Z}}
\eeqn
(where one has $\p$ copies of $\Pi$ in the exponent). Here within this section, without introducing further notations, we always change coefficients to the Novikov field $\Lambda_{\kp}^\Pi$ using the summation homomorphism $\Pi\oplus \cdots \oplus \Pi \to \Pi$. Then there is a filtration on 
\beqn
CF(H)^{\otimes \p} \otimes CM(f^\ep)
\eeqn
corresponding to the sum of the symplectic actions of the $\p$ tensor factors such that the differential preserves this filtration. The filtration restricts to the $\zp$-invariant part. Moreover, after changing coefficients to $\Lambda_{\kp}^{\Pi}$, one obtains a Floer-type complex over this Novikov field
\beqn
\left( CF(H)^{\otimes \p} \otimes CM(f^\ep) \right) \otimes \Lambda_{\kp}^{\Pi}
\eeqn
as well as its $\zp$-invariant part. 

In principle one should hope to establish a chain-level inverse of the equivariant pair-of-pants product via the coproduct construction (cf. \cite{shelukhin-zhao}). However, the corresponding AMS construction is quite involved. Since our main purpose of the quantitative aspects is to compare barcodes, instead of this potential direct approach, we use more algebraic treatment to bypass the coproduct construction (although we will make use of the local coproduct construction in \cite{shelukhin-zhao}, see the proof of Lemma \ref{lemma2711}). 

\subsubsection{Homological perturbation}

In view of localization, the cohomology of $C_{\zp}(CF(H)^{\otimes \p})$ should only come from generators of the form $x \otimes \cdots \otimes x$ for $x \in {\rm Ob}{\mb F}$, i.e., $\zp$-fixed points of $\p$-tuples of $H$-orbits. We use the homological perturbation argument to remove the effects of the other free $\zp$-orbits. We would like to construct a complex
\beqn
C_{\zp}^{\rm red}(CF(H)^{\otimes \p})
\eeqn
which is filtered chain homotopy equivalent to $C_{\zp}(CF(H)^{\otimes \p})$ and which is only generated by the $\zp$-fixed points.

Let ${\bf x} = (x_1, \ldots, x_\p)$ denote a $\p$-tuple of capped $H$-orbits. We write 
\beqn
CF(H)^{\otimes \p} \otimes CM(f^\ep) = \bigoplus_{{\bf x} \in \tilde {\mc O}(H)^{\p}} \fp \langle {\bf x} \rangle \otimes CM(f^\ep).
\eeqn
Here the infinite direct sum indicates we take a completion. This is a decomposition of $\kp$-modules. Let $[{\bf x}]$ denote the $\zp$-orbit of ${\bf x}$. Then there is a $\zp$-invariant summand
\beqn
\bigoplus_{{\bf x} \in [{\bf x}]} \fp \langle {\bf x} \rangle \otimes CM(f^\ep).
\eeqn
Its $\zp$-fixed part is denoted by 
\beqn
C_{\zp}(CF(H)^{\otimes \p}, [{\bf x}]) \subset C_{\zp}(CF(H)^{\otimes \p}).
\eeqn

\begin{lemma}
The differential $d_{\zp}^{\otimes \p}$ on $C_{\zp}(CF(H)^{\otimes \p})$ admits a decomposition 
\beqn
d_{\zp}^{\otimes \p} = \delta_{\zp}^{\otimes \p} + D_{\zp}^{\otimes \p}
\eeqn
which satisfies the following conditions.
\begin{enumerate}

\item $\delta_{\zp}^{\otimes \p}$ preserves the filtration and for each tuple ${\bf x}$, 
\beqn
\delta_{\zp}^{\otimes \p} \Big( C_{\zp}(CF(H)^{\otimes \p}, [{\bf x}]) \Big) \subset C_{\zp} ( CF(H)^{\otimes \p}, [{\bf x}]).
\eeqn

\item The restriction of $\delta_{\zp}^{\otimes \p}$ to each $C_{\zp}(CF(H)^{\otimes \p}, [{\bf x}])$, denoted by $\delta_{[{\bf x}]}$, is a differential. 

\item If ${\bf x} = (x, \ldots, x)$ is a $\zp$-fixed point, then $\delta_{[{\bf x}]}$ is trivial. If ${\bf x}$ is not a fixed point, then the cohomology of $\delta_{[{\bf x}]}$ is trivial.

\item There exists a constant $\epsilon_H>0$ such that $D_{\zp}^{\otimes \p}$ increases the filtration by at least $\epsilon_H$.
\end{enumerate}
\end{lemma}

\begin{proof}
The differential on the tensor product $CF(H)^{\otimes \p} \otimes CM(f^\ep)$, denoted by $\tilde d_{\zp}^{\otimes\p}$,  is the tensor product differential. Let $d_H: CF(H) \to CF(H)$ be the Floer differential and $d_{\ep}$ be the Morse differential on $CM(f^\ep)$. Then the obvious $\zp$-equivariant decomposition
\beqn
\begin{split}
\tilde d_{\zp}^{\otimes \p} = &\ \tilde \delta_{\zp}^{\otimes \p}  + \tilde D_{\zp}^{\otimes \p} \\
 := &\ \ast {\rm Id}_{CF(H)} \otimes \cdots \otimes {\rm Id}_{CF(H)} \otimes d_{\ep}  + d_H^{\otimes \p} \otimes {\rm Id}_{CM(f^\ep)}.
\end{split}
\eeqn
Here $\ast = \pm 1$ depends on the degrees of the variables. The restriction to the $\zp$-invariant part provides the desired decomposition.
\end{proof}

Now we would like to apply the homological perturbation argument to construct a differential on the chain group
\beq\label{reduced_complex}
C_{\zp}^{\rm red}(CF(H)^{\otimes \p}):= \bigoplus_{ x\in \tilde {\mc O}(H)} C_{\zp}(CF(H)^{\otimes \p}, x\otimes \cdots \otimes x),
\eeq
which is the cohomology of $C_{\zp}(CF(H)^{\otimes \p})$ under the differential $\delta_{\zp}^{\otimes \p}$ defined above. Notice that there is a natural linear inclusion
\beqn
\sigma: C_{\zp}^{\rm red}(CF(H)^{\otimes \p}) \to C_{\zp}(CF(H)^{\otimes \p})
\eeqn
and a projection
\beqn
\pi: C_{\zp}(CF(H)^{\otimes \p}) \to C_{\zp}^{\rm red}(CF(H)^{\otimes \p}).
\eeqn

\begin{lemma}\label{lemma_homological_perturbation_1}
There exists a differential $D_{\zp}^{\rm red}$ on $C_{\zp}^{\rm red}(CF(H)^{\otimes \p})$, together with maps
\begin{align*}
&\ \sigma': C_{\zp}^{\rm red}(CF(H)^{\otimes \p}) \to C_{\zp}(CF(H)^{\otimes \p}),\ &\ \pi': C_{\zp}(CF(H)^{\otimes \p}) \to C_{\zp}^{\rm red}( CF(H)^{\otimes \p})
\end{align*}
satisfying the following conditions.
\begin{enumerate}

\item $D_{\zp}^{\rm red}$ strictly increases the energy filtration.

\item $\sigma'$ and $\pi'$ strictly increase the energy filtration.

\item $\tilde \sigma =\sigma +\sigma'$ and $\tilde \pi = \pi + \pi'$ are filtered chain homotopy inverses to each other.
\end{enumerate}
\end{lemma}

\begin{proof}
For each $[{\bf x}]$ which is not a fixed point, since the cohomology of $\delta_{[{\bf x}]}$ is trivial, one can find a chain homotopy between identity and zero, denoted by $\Theta_{[{\bf x}]}$ satisfying
\beqn
{\rm Id} = \delta_{[{\bf x}]} \Theta_{[{\bf x}]} + \Theta_{[{\bf x}]} \delta_{[{\bf x}]}.
\eeqn
When $[{\bf x}]$ is a $\zp$-fxed point, define $\Theta_{[{\bf x}]} = 0$. One can choose $\Theta_{[{\bf x}]}$ in a $\Pi$-equivariant way. Then define
\beqn
\Theta:= \bigoplus_{[{\bf x}]\in \tilde {\mc O}(H)^{\p}/ \zp} \Theta_{[{\bf x}]}.
\eeqn
Then $\Theta$ preserves the filtration. Then following the explicit formula of \cite{Huebschman_Kadeishvili} and \cite{Markl_perturbation}), define
\beqn
\tilde D_{\zp}^{\otimes \p}:= \pi \left( \sum_{l=0}^\infty (D_{\zp}^{\otimes \p} \circ \Theta )^l \right) D_{\zp}^{\otimes \p} \sigma = \pi D_{\zp}^{\otimes \p} \left( \sum_{l=0}^\infty (\Theta D_{\zp} )^l \right) \sigma.
\eeqn
Its restriction to the $\zp$-invariant part is defined to be $D_{\zp}^{\rm red}$. The terms $\sigma'$ and $\pi'$ are also defined by explicit formulas (see, e.g. \cite[Appendix B]{BSWX}), and we omit the details.
\end{proof}

\subsection{Filtration on equivariant Floer cochain complex}

Let $H$ be a Hamiltonian on $X$ so that an equivariant Floer cochain complex $CF_{\zp} (H^{\natural \p})$ can be defined (see Theorem \ref{thma_equivariant_Floer}). Recall that to define this chain group, we first consider the chain group
\beqn
C ( {\mb F}^{\natural \p} \# {\mb E} ):= \left\{ \sum_{i=1}^\infty c_i p_i \otimes w_i\ \left|\ \begin{array}{l} c_i \in \fp,\ p_i \in \tilde {\mc O}(H^{\natural \p}),\ w_i \in {\rm crit} f^\ep,\\
\displaystyle \inf_i {\mc A}(p_i) >-\infty,\ \inf_i f^\ep(w_i) > -\infty,\\
\displaystyle \lim_{i \to \infty} {\mc A}(p_i) + f^\ep(w_i) = +\infty 
\end{array} \right. \right\}
\eeqn
It is a module with a $\zp$-action, and it is defined over the Novikov field $\Lambda_{\fp}^{\mb\Pi}$. The equivariant complex $CF_{\zp}(H^{\natural \p})$ is (roughly) its $\zp$-fixed part.  

\begin{lemma}
The $\zp$-invariant part $C(F^{\natural \p}\# {\mb E})_{\zp}$ is a free module over $\Lambda_{\fp}^{\mb{\Pi}}$. Its rank is $\# {\mc O}(H^{\natural \p})$ when $\p = 2$ and $2 \# {\mc O}(H^{\natural \p})$ when $\p > 2$.
\end{lemma}

\begin{proof}
Choose for each contractible 1-periodic orbit a capping. Let $p_1, \ldots, p_m$ be these capped 1-periodic orbit. Let $\tau \in \zp$ be the generator. When $\p = 2$, define a map
\beqn
F: \Lambda_{{\bf F}_2}^{\mb\Pi} \langle p_1, \ldots, p_m \rangle \to C ({\mb F}^{\natural \p} \# {\mb E})_{{\mb Z}_2}
\eeqn
by 
\beqn
F(p_\alpha) = p_\alpha \otimes w_0^+ + \tau(p_\alpha) \otimes w_0^- 
\eeqn
where $w_0^\pm$ are the two critical points of $f^\ep$ of index $0$. One can see that this is an isomorphism of free $\Lambda_{\fp}^{\mb{\Pi}}$-modules. On the other hand, when $\p>2$, define
\beqn
F: \Lambda_{\fp}^{\mb\Pi}\langle p_1, \ldots, p_m, p_1 \otimes \uptheta, \ldots, p_m \otimes \uptheta \rangle \to C ({\mb F}^{\natural \p} \# {\mb E})_{\zp}
\eeqn
by sending 
\begin{align*}
&\ F(p_\alpha) = \sum_{j=1}^\p \tau^j(p_\alpha) \otimes w_0^j,\ &\ F(p_\alpha \otimes \uptheta) = \sum_{j=1}^\p \tau^j (p_\alpha) \otimes w_1^j,
\end{align*}
where $w_0^j$, $w_1^j$ are the critical points of $f^\ep$ of index $0$ and $1$. One can also verify that $F$ is a linear isomorphism.
\end{proof}

Note the difference between the Novikov fields $\Lambda_{\fp}^{\mb\Pi}$ and $\Lambda_{\kp}^{\Pi}$, where the latter does not require the pole order of the variable $\tu$ or $\uptheta$ is bounded. Then we can write the chain group after tensoring with this larger Novikov field by
\begin{multline*}
CF_{\zp}(H^{\natural \p}; \Lambda_{\kp}^\Pi)= C({\mb F}^{\natural \p} \# \ep)_{\zp}  \otimes \Lambda_{\kp}^\Pi\\
\cong \Big\{  \sum_{i=1}^\infty \tilde c_i p_i \ |\ \tilde c_i \in \tkp,\ p_i \in  \tilde {\mc O}(H^{\natural \p}),\ \lim_{i \to \infty} {\mc A}_{H^{\natural \p}}(p_i) = +\infty \Big\}.
\end{multline*}

To directly use results from \cite{Usher_Zhang_2016}, change coefficients to $\Lambda_{\kp}^{|\Pi|}$, which becomes only ${\mb Z}_2$-graded. Recall from \eqref{eqn:r-value} the meaning of the notation $|\Pi|$.

\begin{lemma}
The complex 
\beqn
CF_{\zp}(H^{\natural \p}; \Lambda_{\kp}^{|\Pi|}):= CF_{\zp}(H^{\natural \p}; \Lambda_{\kp}^\Pi) \otimes \Lambda_{\kp}^{|\Pi|}
\eeqn
is a Floer-type complex (Definition \ref{defn_Floer_type_complex}) over the Novikov field $\Lambda_{\kp}^{|\Pi|}$.
\end{lemma}

\begin{proof}
To see this, we first need to prove that $CF_{\zp}(H^{\natural  \p}; \Lambda_{\kp}^{|\Pi|})$ is orthogonalizable. Indeed, it is the $\zp$-invariant subspace of the cochain complex $C( {\mb F}^{\natural \p} \# {\mb E}) \otimes \Lambda_{\kp}^{|\Pi|}$, which is orthogonalizable. Then by \cite[Corollary 2.17]{Usher_Zhang_2016}, the subspace $CF_{\zp}(H^{\natural \p}; \Lambda_{\kp}^{|\Pi|})$ is also orthogonalizable. On the other hand, the differential preserves the energy filtration, hence $CF_{\zp}(H^{\natural \p}; \Lambda_{\kp}^{|\Pi|})$ is a Floer-type complex.
\end{proof}

As a corollary, by the basic construction of \cite{Usher_Zhang_2016}, the complex $CF_{\zp}(H^{\natural \p}; \Lambda_{\kp}^{|\Pi|})$ has  well-defined verbose and concise barcodes. 

\begin{prop}
The verbose barcode only depends on $H^{\natural \p}$.
\end{prop}

\begin{proof}
By Theorem \ref{thm_Usher_Zhang_barcodes}, one only needs to show that the filtered isomorphism type only depends on $H^{\natural \p}$. Indeed, given two constructions of $C({\mb F}^{\natural \p} \# {\mb E})$ as a chain complex of $\Lambda_{\fp}^{\mb\Pi}$-modules with $\zp$-action corresponding to two AMS lifts
\begin{align*}
&\ \tilde {\mb F}_0^{\zp}, &\ \tilde {\mb F}_1^{\zp},
\end{align*}
by Theorem \ref{thm_equivariant_lift}, one can construct a lift of the outercollaring of the diagonal bimodule to $\outer \uds{\bf dOrb}_{\rm rig}^{\rm NC}$ as a bimodule over $(\tilde {\mb F}_0^{\zp}; \tilde {\mb F}_1^{\zp})$. Then by choosing FOP perturbations, one obtains a $\zp$-equivariant chain map
\beqn
\Phi: C(\tilde{\mb F}_0^\zp) \to C(\tilde{\mb F}_1^\zp).
\eeqn
Similar to many other cases, this chain map has a leading order term being the identity map on the same chain group. Hence $\Phi$ is invertible. Moreover, $\Phi$ preserves the energy filtration because all the lifts cover the constant homotopy of $H^{\natural \p}$. Therefore, after changing coefficients to $\Lambda_{\kp}^{|\Pi|}$, the chain map $\Phi$ is a filtered isomorphism.
\end{proof}

\subsection{Local equivariant cohomology and homological perturbation}

We would like to compare the barcodes between the ordinary Floer complex of a Hamiltonian and the equivariant Floer complexes of its $\p$-th iteration using the chain-level equivariant pair-of-pants product. In the presence of simple $\p$-periodic orbits, we do not expect to obtain a filtered isomorphism. On the other hand, simple $\p$-periodic orbits, which form free $\zp$-orbits, are not expected to contribute in the Tate equivariant cohomology. We use local equivariant cohomology and the homological perturbation method to remove simple $\p$-periodic orbits from a reduced version of the equivariant Floer complex.

\subsubsection{Local equivariant Floer theory}

Now consider the $\p$-th iteration $H^{\natural \p}$. Recall that the chain group $CF_{\zp}(H^{\natural \p})$ is the $\zp$-invariant part of the module
\beqn
\bigoplus_{ \ov{y} \in {\mc O}(H^{\natural \p})} \Lambda_{\fp}^\Pi \langle \ov{y} \rangle \otimes CM(f^\ep). 
\eeqn
Notice that $\ov{y}$ is either an iterated orbit (fixed point) or a simple $\p$-periodic orbit (free $\zp$-orbit). Then we can write the above chain group as
\beqn
\left(  \bigoplus_{ \uds{x} \in {\mc O}(H)} \Lambda_{\fp}^\Pi \langle \uds x^{(\p)} \rangle \otimes CM(f^\ep) \right) \wh \oplus \left(  \bigoplus_{\uds y \in {\mc O}^{\rm simple}(H^{\natural \p})} \Lambda_{\fp}^\Pi \langle \uds y \rangle \otimes CM(f^\ep)\right).
\eeqn
Notice that $\zp$ acts trivially on the set of iterated orbits and acts freely on the set of simple orbits. Let $[\uds y]$ be the $\zp$-orbit of $\uds y$. Then one reorganizes the decomposition as 
\beqn
\begin{split}
CF_{\zp}(H^{\natural \p}) \otimes \Lambda_{\kp}^\Pi \cong &\ \bigoplus_{\uds x \in {\mc O}(H)} \Lambda_{\tkp}^\Pi \langle \uds x^{(\p)} \rangle \wh \oplus \left( \bigoplus_{\uds y \in {\mc O}^{\rm simple}(H^{\natural \p})} \Lambda_{\fp}^\Pi \langle \uds y \rangle \otimes CM(f^\ep) \right)_{\zp}\\
\cong &\ \bigoplus_{[\uds y]\in {\mc O}(H^{\natural \p})/ \zp} CF_{\zp, {\rm loc}}( H^{\natural \p}, [\uds y] ),
\end{split}
\eeqn
where each $CF_{\zp, {\rm loc}}(H^{\natural \p}, [\uds y])$ is a submodule coming from each $\zp$-orbit. Notice that for each $[y] \in {\mc O}(H)$, there is a $\kp$-summand
\beqn
CF_{\zp, {\rm loc}}(H^{\natural}, [y]) \subset CF_{\zp, {\rm loc}}(H^{\natural \p}, [\uds y]).
\eeqn

The following lemma basically says that one can decompose the equivariant differential into a local part and a positive-energy part, where the local part respects the above decomposition.  

\begin{lemma}
The equivariant differential
\beqn
d_{\zp}: CF_{\zp}(H^{\natural \p}; \Lambda_{\kp}^{\Pi}) \to CF_{\zp}(H^{\natural \p}; \Lambda_{\kp}^\Pi)
\eeqn
can be written as 
\beqn
d_{\zp} = d_{\zp, {\rm loc}} + D_{\zp}
\eeqn
which satisfy the following conditions.

\begin{enumerate}

\item For each $[y] \in \tilde {\mc O}(H^{\natural \p})/\zp$, 
\beqn
d_{\zp, {\rm loc}} \big( CF_{\zp, {\rm loc}}(H^{\natural \p}, [y])  \big) \subset CF_{\zp, {\rm loc}}(H^{\natural \p}, [y]).
\eeqn
Therefore $d_{\zp, {\rm loc}}$ preserves the filtration. Denote the restriction of $d_{\zp, {\rm loc}}$ to the summand $CF_{\zp, {\rm loc}}(H^{\natural \p}, [y])$ by $d_{\zp, [y]}$.

\item For each $[y]$, $d_{\zp, [y]}^2 = 0$.

\item There exists $\delta_{H^{\natural \p}}>0$ such that $D_{\zp}$ increases the filtration by at least $\delta_{H^{\natural \p}}$. 

\item If $y = x^{(\p)}$ for $x \in \tilde {\mc O}(H)$, $d_{\zp, [y]}$ vanishes and denote the cohomology by 
\beqn
HF_{\zp, {\rm loc}}(H^{\natural \p}, x^{(\p)}) \cong \tkp.
\eeqn

\item If $y$ is a simple $\p$-periodic orbit, then the cohomology of $d_{\zp, [y]}$ is trivial.
\end{enumerate}
\end{lemma}

\begin{proof}
The decomposition $d_{\zp} = d_{\zp, {\rm loc}} + D_{\zp}$ and the first three conditions are due to a basic compactness argument, usually termed under cross-energy lemma, see e.g. \cite[Lemma B.1]{BSWX}. In our nondegenerate setting, we know that $d_{\zp, {\rm loc}}$ comes from differentials in $CM(f^\ep)$ and $D_{\zp}$ comes from the contribution from nontrivial Floer cylinders in $X$.

If $y = x^{(\p)}$, the summand $CF_{\zp, {\rm loc}}(H^{\natural\p}, x^{(\p)}) \cong CM(f^\ep)_{\zp}$ and the differential is trivial on $CM(f^\ep)_{\zp}$. Hence (4) is verified.

If $y$ is simple, then we have the following identification. Let $\tau \in \zp$ be the generator and let $w_l^\pm$ resp. $w_l^0, \ldots, w_l^{\p-1}$ be the $\p$ critical points of $f^\ep$ in degree $l$ in $\p=2$ resp. $\p>2$ case. 
\begin{multline*}
CF_{\zp, {\rm loc}}(H^{\natural \p}, [y]) = \bigoplus_{\tau^k \in \zp} \left( \fp \langle \tau^k y \rangle \oplus CM(f^\ep) \right)_{\zp} \\
\cong \left\{ \begin{array}{ll}  {\bf K}_2  \Big\langle y \otimes w_0^- + \tau(y) \otimes w_0^+, y \otimes w_0^+ + \tau(y) \otimes w_0^-  \Big\rangle  ,\ &\ \p = 2,\\[0.2cm]
\kp \Big\langle \left. \sum_{k=0}^{\p-1} \tau^k(y) \otimes w_0^{k+s},\ \sum_{k=0}^{\p-1} \tau^k(y) \otimes w_1^{k+s}\right|\ s = 0, 1, \ldots, \p-1 \Big\rangle &\ \p>2.
\end{array}\right.
\end{multline*}
One can see from the differential of $CM(f^\ep)$ that the cohomology is trivial.
\end{proof}

Then we would like to apply the homological perturbation argument.

\begin{lemma}\label{lemma_reduced_complex}
There exist a differential $D_{\zp}^{\rm red}$ on the ``reduced'' equivariant complex
\beqn
CF_{\zp}^{\rm red}(H^{\natural \p}):= \bigoplus_{x \in \tilde {\mc O}(H)} HF_{\zp, {\rm loc}}(H^{\natural \p}, x^{(\p)}) 
\eeqn
making it a Floer-type complex, filtered chain homotopy equivalences
\begin{align*}
&\ \tilde \sigma: CF_{\zp}^{\rm red}(H^{\natural \p}) \to CF_{\zp}(H^{\natural \p}),\ &\ \tilde \pi: CF_{\zp}(H^{\natural \p}) \to CF_{\zp}^{\rm red}(H^{\natural \p}),
\end{align*}
and a filtered chain homotopy
\beqn
\tilde \Theta: CF_{\zp}(H^{\natural\p}) \to CF_{\zp}(H^{\natural \p}) 
\eeqn
such that 
\begin{align*}
&\ \tilde \pi \circ \tilde \sigma = {\rm Id},\ &\ \tilde \sigma \circ \tilde \pi = {\rm Id} + d_{\zp} \circ \tilde \Theta + \tilde\Theta \circ d_{\zp}.
\end{align*}
\end{lemma}

\begin{proof}
This lemma can be proved in the same way as Lemma \ref{lemma_homological_perturbation_1}. We omit the details.
\end{proof}

\subsection{Equivariant pair-of-pants product}

We consider the quantitative property of the equivariant pair-of-pants product. Let $H$ be as above. First, the complex $CF (H)$ is a Floer-type complex over the Novikov field $\Lambda_{\fp}^{\Pi}$-modules with a filtration ${\mc A}_H$. Then consider the completed tensor product
\beqn
\left( CF(H)^{\otimes \p} \otimes CM(f^\ep) \right) \otimes \Lambda_{\kp}^{\Pi}
\eeqn
which has the energy filtration ${\mc A}_H^{\otimes \p}$. It is a nonarchimedean normed vector space over $\Lambda_{\kp}^{\Pi}$ and is orthogonalizable. Then its $\zp$-invariant subspace
\beqn
C_{\zp}( CF(H)^{\otimes \p}):= \Big( \Big(  CF(H)^{\otimes \p} \otimes CM(f^\ep)\Big)\otimes \Lambda_{\kp}^{\Pi}  \Big)_{\zp}
\eeqn
is orthogonalizable. Moreover, the differential clearly preserves the filtration. Hence one has a Floer-type complex with associated barcode. 

Now consider the chain-level equivariant pair-of-pants product, which was originally defined as 
\beqn
P: \Big( CF(H)^{\otimes \p}\otimes CM(f^\ep) \otimes \Lambda_{\kp}^\Pi \Big)_{\zp} \to CF_{\zp}(H^{\natural \p})
\eeqn
as a $\Lambda_{\kp}^{\Pi}$-linear chain map. 

\begin{prop}
There exists a chain-level construction of $P$ which preserves the nonarchimedean norm, i.e., 
\beqn
{\mc A}(P (x)) \geq {\mc A}_H^{\otimes \p} (x).
\eeqn
\end{prop}

\begin{proof}
Consider the branched covering $\Sigma^{\mb{FSt}} \to {\mb R}\times S^1$, which pulls back the trivial Floer data $H dt$ on the cylinder to a Hamiltonian connection on $\sigma_H^{\mb{FSt}}$ on $\Sigma^{\mb{FSt}}$. Consider the multimodule $\mb{X}^{\mb{FSt}}$. From the energy identity for Floer equation, one can see that the multimodule has a threshold zero (Definition \ref{defn_threshold}) with respect to the symplectic action. Then by Proposition \ref{prop_spectral_continuity}, the chain map preserves the filtration.
\end{proof}

\begin{thm}\label{thm_equivariant_pants_isomorphism}
There exists a filtered chain isomorphism
\beqn
P^{\rm red}: C_{\zp}^{\rm red}(CF(H)^{\otimes \p}) \to CF_{\zp}^{\rm red}(H^{\natural \p})
\eeqn
such that the following diagram commutes up to filtered homotopy, where the vertical arrows are provided by the homological perturbation lemma.
\beqn
\xymatrix{ C_{\zp}^{\rm red}(CF(H)^{\otimes \p})    \ar[r]^{P^{\rm red}} \ar[d]   &    CF_{\zp}^{\rm red}(H^{\natural \p}) \ar[d] \\
          C_{\zp}(CF(H)^{\otimes \p}) \ar[r]_P   &   CF_{\zp}(H^{\natural \p})    }
\eeqn
\end{thm}

The proof is provided momentarily. As the reduced versions of the involved complexes are filtered homotopy equivalent to the original versions, by Theorem \ref{thm_Usher_Zhang_barcodes}, one has the following immediate corollary.

\begin{cor}\label{cor_eq_pants_filtered}
$P$ is a filtered homotopy equivalence. Therefore, the concise barcodes of $C_{\zp}(CF(H)^{\otimes \p})$ and $CF_{\zp}(H^{\natural \p})$ coincide.
\end{cor}

To prove Theorem \ref{thm_equivariant_pants_isomorphism}, we proceed by building the local inverses.

\begin{lemma}\label{lemma2711}
We can write the equivariant pair-of-pants product $P$ as
\beqn
P = P_{\rm loc} + P_{\rm big}
\eeqn
which satisfy the following conditions.

\begin{enumerate}

\item For each $x \in \tilde {\mc O}(H)$, denote $x^{\otimes p} = x\otimes \cdots \otimes x$. Then one has 
\beqn
P_{\rm loc}( C_{\zp, {\rm loc}}(CF(H)^{\otimes \p}, x^{\otimes \p})) \subset CF_{\zp, {\rm loc}}(H^{\natural \p}, x^{(\p)}).
\eeqn
Denote the restriction by 
\beqn
P_x: C_{\zp, {\rm loc}}( CF(H)^{\otimes \p}, x^{\otimes \p}) \to CF_{\zp, {\rm loc}}(H^{\natural \p}, x^{(\p)}).
\eeqn

\item There exists $\delta_H''>0$ such that $P_{\rm big}$ increases the filtration by at least $\delta_H''$.

\item $P_x$ is independent of choices and agrees with the local equivariant pair-of-pants product constructed by Seidel \cite{seidel-pants} (for $\p = 2$) and Shelukhin--Zhao \cite{shelukhin-zhao} (for $\p > 2$).

\item $P_x$ is invertible.
\end{enumerate}
\end{lemma}

\begin{proof}
The decomposition and the first two items are based on the crossing energy argument (cf. \cite[Lemma B.1]{BSWX}). To see that $P_{\rm loc}$ is well-defined, suppose $P'$ is another construction which has a similar decomposition $P' = P_{\rm loc}' + P_{\rm big}'$. Then $P$ and $P'$ are homotopic via a filtered chain homotopy $K$. So 
\beqn
P - P' = K d_{\zp} + d_{\zp} K.
\eeqn
Since $K$ also admits the decomposition into the local component and ``big" component in view of the cross-energy lemma, by passing to the local piece, one has $P_{\rm loc} = P_{\rm loc}'$ on cohomology. 

Furthermore, since the local argument does not involve sphere bubbling, the equivalence between our FOP counts and the counts defined using classical framework as used in \cite{Seidel_pants} and \cite{shelukhin-zhao} can be made using the same argument as proving Part (2) of Theorem \ref{thma_QST}.

Lastly, the invertibility of $P_x$ is proved in \cite{Seidel_pants} by a local computation for $\p = 2$, and in \cite{shelukhin-zhao} by constructing the local equivariant coproduct and prove that the local coproduct is the inverse of $P_x$ for $\p > 2$.
\end{proof}

\begin{proof}[Proof of Theorem \ref{thm_equivariant_pants_isomorphism}]
Now consider the composition
\beqn
P^{\rm red}:= \pi \circ P \circ \sigma: C_{\zp}^{\rm red}(CF(H)^{\otimes \p}) \to CF_{\zp}^{\rm red}(H^{\natural \p})
\eeqn
which is a filtered chain map. Notice that there is an obvious leading order term, denoted by $P_{\rm loc}^{\rm red}$, which is the direct sum of $P_x$, while the difference
\beqn
P_{\rm big}^{\rm red}:= P^{\rm red} - P_{\rm loc}^{\rm red}
\eeqn
strictly increases the energy filtration. Hence $P^{\rm red}$ is invertible. Therefore, after tensoring with the correct Novikov field, $P^{\rm red}$ is a filtered isomorphism of Floer-type complexes.
\end{proof}

\subsection{Proof of Theorem \ref{thm_pants_filtered}}\label{subsec:thm-M}

The first item, which says that the equivariant Floer cochain complex $CF_{\zp}(H^{\natural \p})$ is a Floer-type complex over $\Lambda_{\kp}^\Pi$ follows directly from the construction.

The second item is exactly Corollary \ref{cor_eq_pants_filtered}.

For the third item, by Corollary \ref{cor_eq_pants_filtered}, the barcodes of the two complexes $C_{\zp}(CF(H)^{\otimes \p})$ and $CF_{\zp}(H^{\natural \p})$ are identical. On the other hand, there exists an algebraically defined quasi-Frobenius map (see \eqref{quasi_Frobenius})
\beqn
qF: HF(H) \to H_{\zp}(CF(H)^{\otimes \p}).
\eeqn
Although $qF$ is not induced from a cochain map, it still respects the filtration in the following way. For each $\tau \in {\mb R}$, one has corresponding $\fp$-linear map
\beqn
qF: HF^{\leq \tau}(H) \to H_{\zp}^{\leq \p \tau} (CF(H)^{\otimes \p}).
\eeqn

In the current setting, it is hard to obtain the relationship between barcodes, but only the bar-length spectrum. Then, following the argument in \cite[Section 6.2]{shelukhin-22}, one can define a version of Floer complex over the universal Novikov ring
\beqn
\Lambda_0^{\rm univ}
\eeqn
(with base ring $\fp$), denoted by $CF(H; \Lambda_0^{\rm univ})$. Then by the general structural results of modules over $\Lambda_0^{\rm univ}$, the bar-lengths $\beta_i(H)$ correspond to torsion parts of $HF(H; \Lambda_0^{\rm univ})$ of the forms $\Lambda_0^{\rm univ}/ T^{\beta_i(H)} \Lambda_0^{\rm univ}$. Then one can prove that the quasi-Frobenius map is an isomorphism of $\Lambda_0^{\rm univ}$-modules after the rescaling by $\p$. The result in Theorem \ref{thm_pants_filtered} (3) follows from Theorem \ref{thm_pants_filtered} (2).

\section{Quantitative Theory for Noncontractible Hamiltonian Orbits}\label{section_noncontractible}

In this section we provide constructions for the Floer complex for noncontractible orbits (in the homological sense) and the corresponding quantitative results involving a version of equivariant pair-of-pants product. As mentioned before, the main motivation for developing this variant is to generalize the result of Sugimoto \cite{Sugimoto_2026} to full generality, which will be done in a separate paper with Shelukhin and Wilkins \cite{BSWX}. Here in this section, we will prove Theorem \ref{thm_Floer_noncontractible} and Theorem \ref{thm_noncontractible_pants}. The only notable technical difference is that the capping of a non-contractible orbit is defined relative to a reference homology class; the rest of the construction is exactly the same as the contractible case.

\subsection{The Floer flow category and bimodules}

Let $\gamma$ be a smooth loop in $X$ which represents a nontrivial class in $H_1(X; {\mb Z})_{\rm free}$. Denote
\beqn
\Gamma = H_2(X; {\mb Z})/ {\rm ker} \omega.
\eeqn
For a nondegenerate Hamiltonian $H$, we describe a ${\mb Z}_2$-graded $\Gamma$-equivariant flow category ${\mb F}_\gamma$ of orbits of $H$ which are homologous to $\gamma$. Recall terms defined in Definition \ref{defn_gamma_capping}. Define
\beqn
{\rm Ob} {\mb F}_\gamma:= \tilde {\mc O}_\gamma(H).
\eeqn
To define morphisms of ${\mb F}_\gamma$, choose a 1-periodic family of $\omega$-compatible almost complex structures $J$ and considering moduli spaces of Floer equations and their compactifications. Then the morphism set between two $\gamma$-capped orbits is the moduli space of stable Floer trajectories connecting them, where the Floer cylinders naturally concatenate the $\gamma$-cappings. One then obtains the flow category $\mb{F}_\gamma$ enriched in $\uds{\bf Top}$ where objects are naturally ${\mb Z}_2$-graded. The action functional 
\beqn
{\mc A}_{H, \gamma}: {\rm Ob}{\mb F}_\gamma \to {\mb R}
\eeqn
makes ${\mb F}_\gamma$ a locally finite flow category (Definition \ref{defn_local_finite_flow}). The $\Gamma$-action on objects makes ${\mb F}_\gamma$ a $\Gamma$-equivariant flow category (Definition \ref{defn_equivariant_flow_category}). 

Notice that if we choose a different loop $\gamma'$ homologous to $\gamma$, then one essentially get the same flow category except that the actions ${\mc A}_{H, \gamma}$ and ${\mc A}_{H, \gamma'}$ differ by a constant.

If we choose a different pair $(H', J')$ and consider the same construction (for a fixed loop $\gamma$), then one obtains a similar Floer flow category ${\mb F}_\gamma'$. By choosing a Floer datum on the infinite cylinder interpolating between $(H, J)$ and $(H', J')$ and considering moduli spaces of Floer equations associated to the Floer datum, one obtains a bimodule $\mb{B}_\gamma$ over $({\mb F}_\gamma; {\mb F}_\gamma')$, which is naturally $\Gamma$-equivariant (Definition \ref{defn_equivariant_bimodule}). This is completely analogous to the usual contractible case.

\subsection{The complex and continuation map}

Since the Floer homology of the noncontractible case is (expected to be) trivial, in fact, the Floer complex is null-homotopic, the only interesting theory is the quantitative theory. Consider
\beqn
CF_\gamma(H):= \Big\{ \sum_{i=1}^\infty a_i p_i\ |\ a_i \in {\mb Z},\ p_i \in {\rm Ob}{\mb F}_\gamma,\ \lim_{i \to \infty} {\mc A}_{H,\gamma} (p_i) = -\infty \Big\}
\eeqn
which is a module over $\Lambda_{{\mb Z}}^
\Gamma$. 

\begin{proof}[Proof of Theorem \ref{thm_Floer_noncontractible}]
The AMS construction for the flow category ${\mb F}_\gamma$ is almost identical to the contractible case. There are only two minor places which we would like to comment. The first one is the existence of integral actions. The proof of Proposition \ref{prop163} works equally well for the noncontractible case. Then one can go through the AMS construction for the flow category. The change of Novikov group from $\Pi$ to $\Gamma$ does not affect our argument to make each step of the construction equivariant. Then upon making choices, one obtains a lift of the outercollaring of ${\mb F}_\gamma$ to $\outer \uds{\bf dOrb}_{\rm rig}^{\rm NC}$. Choosing FOP perturbations allows us to construct the differential, hence obtaining a ${\mb Z}_2$-graded chain complex
\beqn
CF_\gamma(H, J, \Xi)
\eeqn
(linear over $\Lambda_{\mb Z}^\Gamma$). 

After changing coefficients to $\Lambda_{\bf K}^\Gamma$, one obtains a Floer-type complex because the differential increases the filtration. 

To see that the filtered isomorphism class only depends on $H$, we compare two different almost complex structures $J$, $J'$ and two different sets of choices $\Xi$, $\Xi'$, resulting in AMS lifts $\tilde {\mb F}_\gamma$ and $\tilde {\mb F}_\gamma'$. Choose an interpolation between $J$ and $J'$ and consider the corresponding bimodule ${\mb B}_\gamma$ over $({\mb F}_\gamma; {\mb F}_\gamma')$ via continuation map moduli spaces. One apply the AMS construction for ${\mb B}_\gamma$, ending with a chain map 
\beqn
\Phi^{{\mb B}_\gamma}: CF_\gamma(H, J, \Xi) \to CF_\gamma(H, J', \Xi').
\eeqn
Notice that the two chain groups are canonically identified since the Hamiltonian $H$ is fixed. The leading order term of $\Phi^{{\mb B}_\gamma}$ is then the identity map on the chain group $CF_\gamma(H)$. Hence after tensoring with $\Lambda_{{\bf K}}^\Gamma$, the chain map is a filtered isomorphism.
\end{proof}

Therefore, if we switch to field coefficients, for each field ${\bf K}$, we obtain a well-defined verbose barcode ${\mc B}_\gamma(H)$. It still depends on the concrete choice of the loop $\gamma$. We shorten the notation for the chain complex by 
\beqn
CF_\gamma(H; \Lambda_{\bf K}^\Gamma).
\eeqn

\begin{cor}
For two homologous loops $\gamma$ and $\gamma'$, the bar-length spectrum of $CF_\gamma(H; \Lambda_{\bf K}^{\Gamma})$ coincides with that of $CF_{\gamma'}(H; \Lambda_{\bf K}^\Gamma)$.
\end{cor}

\begin{proof}
Let $\gamma'$ be another loop which is homologous to $\gamma$. Choosing a concrete cobordism between them. Then one obtains a $\Gamma$-equivariant isomorphism of flow categories ${\mb F}_\gamma \cong {\mb F}_{\gamma'}$ which shifts the action by a constant. An AMS construction and FOP perturbation for ${\mb F}_\gamma$ can be transported to ${\mb F}_{\gamma'}$. Hence the bar-length spectrum does not depend on the representating loop $\gamma$.
\end{proof}

\begin{rem}
It is standard to argue that the barcode ${\mc B}_\gamma(H)$ varies Lipschitz continuously with $H$, where on barcodes we use the bottleneck distance. Hence the barcode ${\mc B}_\gamma(H)$ also extends to degenerate Hamiltonians.
\end{rem}

\subsection{Equivariant Floer complex}\label{subsec:equiv-noncontrat}

Now we can repeat the equivariant Floer construction to the noncontractible case. Let $\p$ be a prime and suppose $H^{\natural \p}$ is nondegenerate, which may have noncontractible orbits homologous to $\gamma$. By coupling with the Morse flow category $\ep$ on $\wh S^\infty$, one obtains a flow category (see Definition \ref{defn_equi_flow_category})
\beqn
{\mb F}_{\gamma}^{\zp}
\eeqn
which has a free action by the twisted Novikov group $(\Gamma \times {\mb Z}) \rtimes \zp$. Notice that although the reference loop $\gamma$ is not invariant under the $\zp$-action via loop rotation, the symplectic action is still well-defined. 

Using exactly the same construction as the contractible case, one can construct an equivariant lift of the outercollaring of ${\mb F}_\gamma^{\zp}$ to $\outer \uds{\bf dOrb}_{\rm rig}^{\rm NC}$. After choose an equivariant FOP perturbation, denoted by $\mathring {\mb F}_\gamma^\zp$, one obtains a complex with a $\zp$-action, denoted by 
\beqn
C(\mathring {\mb F}_\gamma^\zp).
\eeqn
Its $\zp$-invariant part provides a complex
\beqn
CF_{\zp}(H^{\natural\p}, \gamma)
\eeqn
which is a module over $\Lambda_{\fp}^{\Gamma \oplus {\mb Z}}$. By tensoring with the larger Novikov field $\Lambda_{\kp}^{\Gamma}$ one obtains a Floer-type complex
\beqn
CF_{\zp}(H^{\natural \p}, \gamma; \Lambda_{\kp}^\Gamma).
\eeqn
By using the equivariant continuation map, one proves that the filtered isomorphism class is well-defined. Hence one obtains a barcode
\beqn
{\mc B}_{\zp} (H^{\natural \p}, \gamma)
\eeqn
which only depends on $H$. Again, choosing a homologous loop $\gamma'$ will shift the barcodes but does not change the bar-length spectrum.

\subsection{Equivariant pair-of-pants product}

\begin{figure}[h]
    \centering
    \includegraphics[width=1\linewidth]{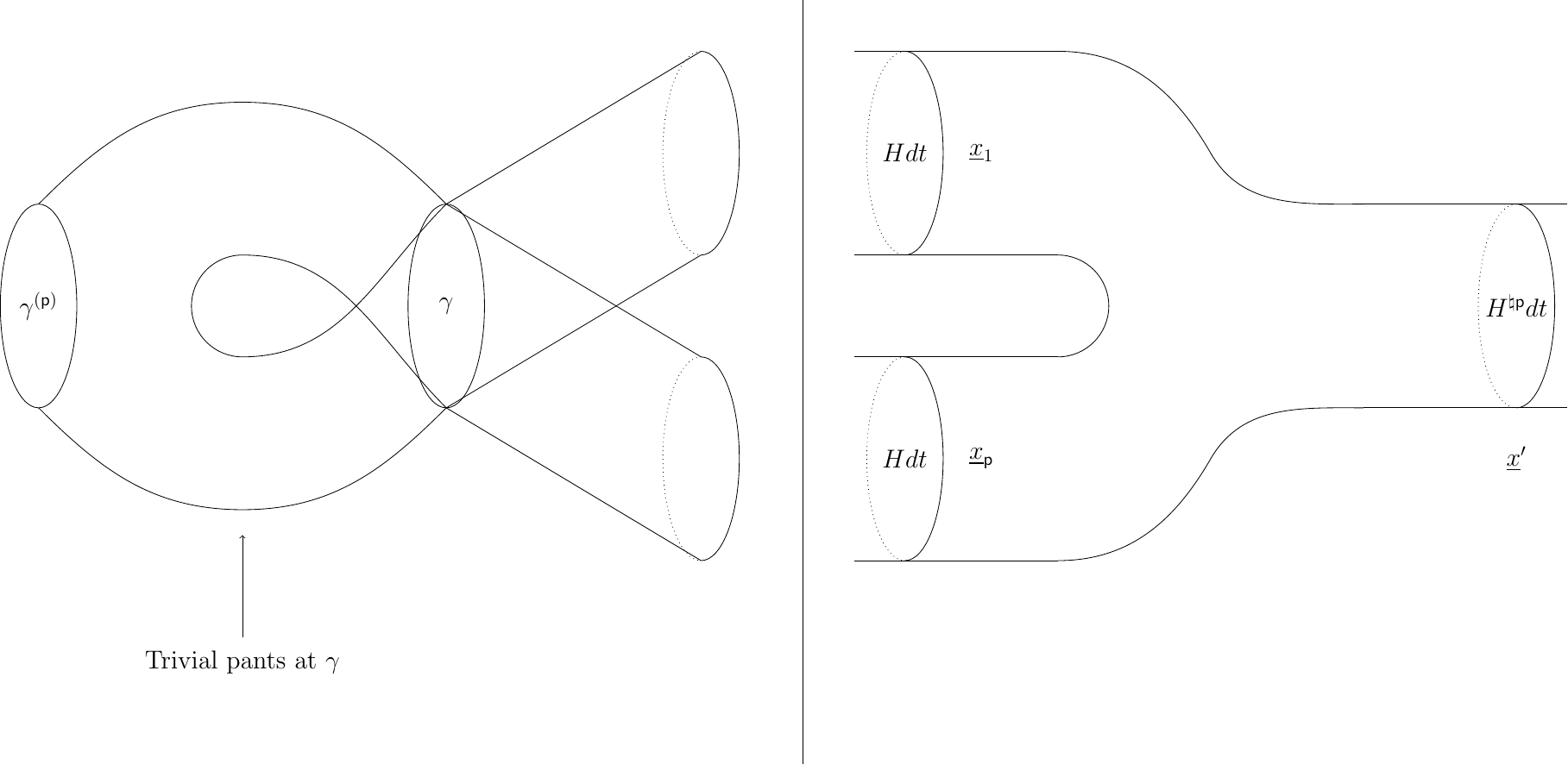}
    \caption{Configurations of equivariant pair-of-pants product for noncontractible orbits}
    \label{Figure_noncontractible}
\end{figure}

One can construct the noncontractible version of the equivariant pair-of-pants product. Let $H$ be a Hamiltonian such that $H^{\natural \p}$ is nondegenerate. Let $\gamma^{(\p)}$ be the $\p$-th iteration of $\gamma$. Basically, the same setup as the contractible case provides a multimodule
\beqn
{\mb X}_\gamma^{\mb{FSt}}
\eeqn
over $({\mb F}_\gamma, \ldots, {\mb F}_\gamma, {\mb E}; {\mb F}_{\gamma^{(\p)}}^{\zp})$. More precisely, let $x_1, \ldots, x_\p$ be objects of ${\mb F}_\gamma$ whose underlying loops $\uds x_1, \ldots, \uds x_\p$ are homologous to $\gamma$. Let $\uds x'$ be a Hamiltonian loop of $H^{\natural \p}$ homologous to $\gamma^{(\p)}$. Consider $p$-legged pair-of-pants forming a cobordism (up to equality of symplectic area) between $\uds x_1 \sqcup \cdots \sqcup \uds x_\p$ represented by a map $u: \Sigma^{\mb{FSt}} \to X$. Then the $\gamma$-cappings $x_1, \ldots, x_\p$ together with the relative homology class of $u$ induces a $\gamma^{(\p)}$-capping of $\uds x'$ as illustrated in Figure \ref{Figure_noncontractible}. 

The AMS construction and FOP perturbation provide a chain map
\beqn
\underbrace{CF_\gamma(H) \otimes \cdots \otimes \cdots CF_\gamma(H)}_{\p} \otimes CM(f^\ep) \to C(\mathring {\mb F}_{\gamma^{(\p)}}^\zp).
\eeqn
Restricting to the $\zp$-invariant part one obtains
\beqn
P_\gamma: C_{\zp}( CF_\gamma(H)^{\otimes \p}) \to CF_{\zp}(H^{\natural \p}, \gamma^{(\p)}).
\eeqn

\begin{prop}\label{prop:pants-noncontract}
$P_\gamma$ is a filtered chain homotopy equivalence of Floer-type complexes.
\end{prop}

\begin{proof}
The argument is similar to the contractible case. One can first define the reduced versions
\begin{align*}
&\ C_{\zp}^{\rm red}(CF_\gamma(H)^{\otimes \p}),\ &\ CF_{\zp}^{\rm red}(H^{\natural \p}, \gamma^{(\p)})
\end{align*}
(see \eqref{reduced_complex} and Lemma \ref{lemma_reduced_complex}) which are filtered chain homotopy equivalent to the unreduced versions via the homological perturbation construction. Moreover, for each orbit $x \in \tilde {\mc O}_\gamma(H)$, one has a local equivariant pair-of-pants product $P_x$ which can be proved to be invertible. Moreover, the direct sum of $P_x$ is the leading order term of $P_\gamma$. Then $P_\gamma$ is transferred to a filtered isomorphism
\beqn
P_\gamma^{\rm red}: C_{\zp}^{\rm red}(CF_\gamma(H)^{\otimes \p}) \to CF_{\zp}(H^{\natural \p}, \gamma^{(\p)}). \qedhere
\eeqn
\end{proof}

\begin{proof}[Proof of Theorem \ref{thm_noncontractible_pants}]
    The isomorphism class of the Floer-type complex is disucssed in Section \ref{subsec:equiv-noncontrat}. Proposition \ref{prop:pants-noncontract} shows Theorem \ref{thm_noncontractible_pants} (2) holds. Finally, Theorem \ref{thm_noncontractible_pants} (3) follows from the same argument as in Section \ref{subsec:thm-M} with the above ingredients.
\end{proof}

\bibliography{reference}

\bibliographystyle{amsalpha}

\end{document}